\documentclass[12pt]{article}
\usepackage{full page, 
amssymb, amscd, graphicx, 
}
    \title{{\bf Logarithmic tensor product theory for generalized modules for
a conformal vertex algebra}}
    \author{Yi-Zhi Huang, James Lepowsky and Lin Zhang}
    \date{}
    \begin{document}
    \bibliographystyle{alpha}
    \maketitle

    \newtheorem{rema}{Remark}[section]
    \newtheorem{propo}[rema]{Proposition}
    \newtheorem{theo}[rema]{Theorem}
   \newtheorem{defi}[rema]{Definition}
    \newtheorem{lemma}[rema]{Lemma}
    \newtheorem{corol}[rema]{Corollary}
     \newtheorem{exam}[rema]{Example}
\newtheorem{assum}[rema]{Assumption}
     \newtheorem{nota}[rema]{Notation}
        \newcommand{\ba}{\begin{array}}
        \newcommand{\ea}{\end{array}}
        \newcommand{\be}{\begin{equation}}
        \newcommand{\ee}{\end{equation}}
        \newcommand{\bea}{\begin{eqnarray}}
        \newcommand{\eea}{\end{eqnarray}}
        \newcommand{\nno}{\nonumber}
        \newcommand{\nn}{\nonumber\\}
        \newcommand{\lbar}{\bigg\vert}
        \newcommand{\p}{\partial}
        \newcommand{\dps}{\displaystyle}
        \newcommand{\bra}{\langle}
        \newcommand{\ket}{\rangle}
 \newcommand{\res}{\mbox{\rm Res}}
\newcommand{\wt}{\mbox{\rm wt}\;}
\newcommand{\swt}{\mbox{\scriptsize\rm wt}\;}
 \newcommand{\pf}{{\it Proof}\hspace{2ex}}
 \newcommand{\epf}{\hspace{2em}$\square$}
 \newcommand{\epfv}{\hspace{1em}$\square$\vspace{1em}}
        \newcommand{\ob}{{\rm ob}\,}
        \renewcommand{\hom}{{\rm Hom}}
\newcommand{\C}{\mathbb{C}}
\newcommand{\R}{\mathbb{R}}
\newcommand{\Z}{\mathbb{Z}}
\newcommand{\N}{\mathbb{N}}
\newcommand{\A}{\mathcal{A}}
\newcommand{\Y}{\mathcal{Y}}
\newcommand{\comp}{\mathrm{COMP}}
\newcommand{\lgr}{\mathrm{LGR}}

\newcommand{\dlt}[3]{#1 ^{-1}\delta \bigg( \frac{#2 #3 }{#1 }\bigg) }

\newcommand{\dlti}[3]{#1 \delta \bigg( \frac{#2 #3 }{#1 ^{-1}}\bigg) }

 \makeatletter
\newlength{\@pxlwd} \newlength{\@rulewd} \newlength{\@pxlht}
\catcode`.=\active \catcode`B=\active \catcode`:=\active \catcode`|=\active
\def\sprite#1(#2,#3)[#4,#5]{
   \edef\@sprbox{\expandafter\@cdr\string#1\@nil @box}
   \expandafter\newsavebox\csname\@sprbox\endcsname
   \edef#1{\expandafter\usebox\csname\@sprbox\endcsname}
   \expandafter\setbox\csname\@sprbox\endcsname =\hbox\bgroup
   \vbox\bgroup
  \catcode`.=\active\catcode`B=\active\catcode`:=\active\catcode`|=\active
      \@pxlwd=#4 \divide\@pxlwd by #3 \@rulewd=\@pxlwd
      \@pxlht=#5 \divide\@pxlht by #2
      \def .{\hskip \@pxlwd \ignorespaces}
      \def B{\@ifnextchar B{\advance\@rulewd by \@pxlwd}{\vrule
         height \@pxlht width \@rulewd depth 0 pt \@rulewd=\@pxlwd}}
      \def :{\hbox\bgroup\vrule height \@pxlht width 0pt depth
0pt\ignorespaces}
      \def |{\vrule height \@pxlht width 0pt depth 0pt\egroup
         \prevdepth= -1000 pt}
   }
\def\endsprite{\egroup\egroup}
\catcode`.=12 \catcode`B=11 \catcode`:=12 \catcode`|=12\relax
\makeatother

\def\hboxtr{\FormOfHboxtr} 
\sprite{\FormOfHboxtr}(25,25)[0.5 em, 1.2 ex] 

:BBBBBBBBBBBBBBBBBBBBBBBBB |
:BB......................B |
:B.B.....................B |
:B..B....................B |
:B...B...................B |
:B....B..................B |
:B.....B.................B |
:B......B................B |
:B.......B...............B |
:B........B..............B |
:B.........B.............B |
:B..........B............B |
:B...........B...........B |
:B............B..........B |
:B.............B.........B |
:B..............B........B |
:B...............B.......B |
:B................B......B |
:B.................B.....B |
:B..................B....B |
:B...................B...B |
:B....................B..B |
:B.....................B.B |
:B......................BB |
:BBBBBBBBBBBBBBBBBBBBBBBBB |

\endsprite

\def\shboxtr{\FormOfShboxtr} 
\sprite{\FormOfShboxtr}(25,25)[0.3 em, 0.72 ex] 

:BBBBBBBBBBBBBBBBBBBBBBBBB |
:BB......................B |
:B.B.....................B |
:B..B....................B |
:B...B...................B |
:B....B..................B |
:B.....B.................B |
:B......B................B |
:B.......B...............B |
:B........B..............B |
:B.........B.............B |
:B..........B............B |
:B...........B...........B |
:B............B..........B |
:B.............B.........B |
:B..............B........B |
:B...............B.......B |
:B................B......B |
:B.................B.....B |
:B..................B....B |
:B...................B...B |
:B....................B..B |
:B.....................B.B |
:B......................BB |
:BBBBBBBBBBBBBBBBBBBBBBBBB |

\endsprite

\vspace{2em}

\begin{abstract}
We generalize the tensor product theory for modules for a vertex
operator algebra previously developed in a series of papers by the
first two authors to suitable module categories for a ``conformal
vertex algebra'' or even more generally, for a ``M\"obius vertex
algebra.''  We do not require the module categories to be semisimple,
and we accommodate modules with generalized weight spaces.  As in the
earlier series of papers, our tensor product functors depend on a
complex variable, but in the present generality, the logarithm of the
complex variable is required; the general representation theory of
vertex operator algebras requires logarithmic structure.  The first
part of this work is devoted to the study of logarithmic intertwining
operators and their role in the construction of the tensor product
functors.  The remainder of this work is devoted to the construction
of the appropriate natural associativity isomorphisms between triple
tensor product functors, to the proof of their fundamental properties,
and to the construction of the resulting braided tensor category
structure.  This work includes the complete proofs in the present
generality and can be read independently of the earlier series of
papers.
\end{abstract}

\vspace{2em}

\tableofcontents
\vspace{1em}
\noindent{\large \bf References}\hfill 

\newpage

\renewcommand{\theequation}{\thesection.\arabic{equation}}
\renewcommand{\therema}{\thesection.\arabic{rema}}
\setcounter{equation}{0}
\setcounter{rema}{0}

\section{Introduction}
In a series of papers (\cite{tensorAnnounce}, \cite{tensorK},
\cite{tensor1}, \cite{tensor2}, \cite{tensor3}, \cite{tensor4}), the
first two authors have developed a tensor product theory for modules
for a vertex operator algebra under suitable conditions.  A structure
called ``vertex tensor category structure,'' which is much richer than
tensor category structure, has thereby been established for many
important categories of modules for classes of vertex operator
algebras, since the conditions needed for invoking the general theory
have been verified for these categories.  The most important such
families of examples of this theory are listed in Subsection 1.1
below.  In the present work, which has been announced in \cite{HLZ},
we generalize this tensor product theory to a larger family of module
categories, for a ``conformal vertex algebra,'' or even more
generally, for a ``M\"obius vertex algebra,'' under suitably relaxed
conditions.  A conformal vertex algebra is just a vertex algebra in
the sense of Borcherds \cite{B} equipped with a conformal vector
satisfying the usual axioms; a M\"obius vertex algebra is a variant of
a ``quasi-vertex operator algebra'' as in \cite{FHL}.  Central
features of the present work are that we do not require the modules in
our categories to be completely reducible and that we accommodate
modules with generalized weight spaces.

As in the earlier series of papers, our tensor product functors depend
on a complex variable, but in the present generality, the logarithm of
the complex variable is required.  The first part of this work is
devoted to the study of logarithmic intertwining operators and their
role in the construction of the tensor product functors.  The
remainder of this work is devoted to the construction of the
appropriate natural associativity isomorphisms between triple tensor
product functors, to the proof of their fundamental properties, and to
the construction of the resulting braided tensor category structure.
This leads to vertex tensor category structure for further important
families of examples, or, in the M\"obius case, to ``M\"obius vertex
tensor category'' structure.

We emphasize that we develop our representation theory (tensor
category theory) in a very general setting; the vertex (operator)
algebras that we consider are very general, and the ``modules'' that
we consider are very general.  We call them ``generalized modules'';
they are not assumed completely reducible.  Many extremely important
(and well-understood) vertex operator algebras have semisimple module
categories, but in fact, now that the theory of vertex operator
algebras and of their representations is as highly developed as it has
come to be, it is in fact possible, and very fruitful, to work in the
greater generality.  Focusing mainly on the representation theory of
those vertex operator algebras for which every module is completely
reducible would be just as restrictive as focusing, classically, on
the representation theory of semisimple Lie algebras as opposed to the
representation theory of Lie algebras in general.  In addition, once
we consider suitably general vertex (operator) algebras, it is
unnatural to focus on only those modules that are completely
reducible.  As we explain below, such a general representation theory
of vertex (operator) algebras requires logarithmic structure.

A general representation theory of vertex operator algebras is crucial
for a range of applications, and we expect that it will be a
foundation for future developments.  One example is that the original
formulation of the uniqueness conjecture \cite{FLM2} for the moonshine
module vertex operator algebra $V^{\natural}$ (again see \cite{FLM2})
requires (general) vertex operator algebras whose modules might not be
completely reducible.  Another example is that this general theory is
playing a deep role in the (mathematical) construction of conformal
field theories (cf. \cite{HPNAS}, \cite{Hconference},
\cite{HVerlindeconjecture}, \cite{rigidity}, \cite{LPNAS}), which in
turn correspond to the perturbative part of string theory.  Just as
the classical (general) representation theory of groups, or of Lie
groups, or of Lie algebras, is not about any particular group or Lie
group or Lie algebra (although one of its central goals is certainly
to understand the representations of particular structures), the
general representation theory of suitably general vertex operator
algebras is ``background independent,'' in the terminology of string
theory.  In addition, the general representation theory of vertex
(operator) algebras can be thought of as a ``symmetry'' theory, where
vertex (operator) algebras play a role analogous to that of groups or
of Lie algebras in classical theories; deep and well-known analogies
between the notion of vertex operator algebra and the classical notion
of, for example, Lie algebra are discussed in several places,
including \cite{FLM2}, \cite{FHL} and \cite{LL}.

The present work includes the complete proofs in the present
generality and can be read independently of the earlier series of
papers of the first two authors constructing tensor categories.  Our
treatment is based on the theory of vertex operator algebras and their
modules as developed in \cite{FLM2}, \cite{FHL}, \cite{DL} and
\cite{LL}.

\subsection{Tensor product theory for finitely reductive 
vertex operator algebras}

The main families for which the conditions needed for invoking the
first two authors' general tensor product theory have been verified,
thus yielding vertex tensor category structure on these module
categories, are the module categories for the following classes of
vertex operator algebras (or, in the last case, vertex operator
superalgebras):

\begin{enumerate}
\item The vertex operator algebras $V_L$ associated with positive
definite even lattices $L$; see \cite{B}, \cite{FLM2} for these vertex
operator algebras and see \cite{D1}, \cite{DL} for the conditions
needed for invoking the general tensor product theory.

\item The vertex operator algebras $L(k,0)$ associated with affine Lie
algebras and positive integers $k$; see \cite{FZ} for these vertex
operator algebras and \cite{FZ}, \cite{HLaffine} for the conditions.

\item The ``minimal series'' of vertex operator algebras associated
with the Virasoro algebra; see \cite{FZ} for these vertex operator
algebras and \cite{W}, \cite{H3} for the conditions.

\item Frenkel, Lepowsky and Meurman's moonshine module $V^{\natural}$;
see \cite{FLM1}, \cite{B}, \cite{FLM2} for this vertex operator
algebra and \cite{D2} for the conditions.

\item The fixed point vertex operator subalgebra of $V^{\natural}$
under the standard involution; see \cite{FLM1}, \cite{FLM2} for this
vertex operator algebra and \cite{D2}, \cite{H4} for the conditions.

\item The ``minimal series'' of vertex operator superalgebras
(suitably generalized vertex operator algebras) associated with the
Neveu-Schwarz superalgebra and also the ``unitary series'' of vertex
operator superalgebras associated with the $N=2$ superconformal
algebra; see \cite{KW} and \cite{A2} for the corresponding $N=1$ and
$N=2$ vertex operator superalgebras, respectively, and \cite{A1},
\cite{A3}, \cite{HM1}, \cite{HM2} for the conditions.
\end{enumerate}

In addition, vertex tensor category structure has also been
established for the module categories for certain vertex operator
algebras built {}from the vertex operator algebras just mentioned, such
as tensor products of such algebras; this is carried out in certain of
the papers listed above.

For all of the six classes of vertex operator algebras (or
superalgebras) listed above, each of the algebras is ``rational'' in
the specific sense of Huang-Lepowsky's work on tensor product theory.
This particular ``rationality'' property is easily proved to be a
sufficient condition for insuring that the tensor product modules
exist; see for instance \cite{tensor1}.  Unfortunately, the phrase
``rational vertex operator algebra'' also has several other distinct
meanings in the literature.  Thus we find it convenient at this time
to assign a new term, ``finite reductivity,'' to our particular notion
of ``rationality'': We say that a vertex operator algebra (or
superalgebra) $V$ is {\it finitely reductive} if
\begin{enumerate}
\item Every $V$-module is completely reducible.
\item There are only finitely many irreducible $V$-modules (up to
equivalence).
\item All the fusion rules (the dimensions of the spaces of
intertwining operators among triples of modules) for $V$ are finite.
\end{enumerate}
We choose the term ``finitely reductive'' because we think of the term
``reductive'' as describing the complete reducibility---the first of
the conditions (that is, the algebra ``(completely) reduces'' every
module); the other two conditions are finiteness conditions.

The vertex-algebraic study of tensor category structure on module
categories for certain vertex algebras was stimulated by the work of
Moore and Seiberg \cite{MS}, in which, in the study of what they
termed ``rational'' conformal field theory, they constructed a tensor
category structure based on the assumption of the existence of a
suitable operator product expansion for ``chiral vertex operators''
(which correspond to intertwining operators in vertex algebra theory).
Earlier, in \cite{BPZ}, Belavin, Polyakov, and Zamolodchikov had
already formalized the relation between the (nonmeromorphic) operator
product expansion, chiral correlation functions and representation
theory, for the Virasoro algebra in particular, and Knizhnik and
Zamolodchikov \cite{KZ} had established fundamental relations between
conformal field theory and the representation theory of affine Lie
algebras.  As we have discussed in the introductory material in
\cite{tensorK}, \cite{tensor1} and \cite{HLaffine}, such study of
conformal field theory is deeply connected with the vertex-algebraic
construction and study of tensor categories, and also with other
mathematical approaches to the construction of tensor categories in
the spirit of conformal field theory.  Concerning the latter
approaches, we would like to mention that the method used by Kazhdan
and Lusztig, especially in their construction of the associativity
isomorphisms, in their breakthrough work in \cite{KL1}--\cite{KL5}, is
related to the algebro-geometric formulation and study of
conformal-field-theoretic structures in the influential works of
Tsuchiya-Ueno-Yamada \cite{TUY}, Drinfeld \cite{Dr} and
Beilinson-Feigin-Mazur \cite{BFM}.  See also the important work of
Deligne \cite{De}, Finkelberg (\cite{F1} \cite{F2}), Bakalov-Kirillov
\cite{BK} and Nagatomo-Tsuchiya \cite{NT} on the construction of
tensor categories in the spirit of this approach to conformal field
theory.

\subsection{Logarithmic tensor product theory}

The semisimplicity of the module categories mentioned in the examples
above is related to another property of these modules, namely, that
each module is a direct sum of its ``weight spaces,'' which are the
eigenspaces of a special operator $L(0)$ coming {}from the Virasoro
algebra action on the module.  But there are important situations in
which module categories are not semisimple and in which modules are
not direct sums of their weight spaces.  Notably, for the vertex
operator algebras $L(k,0)$ associated with affine Lie algebras, when
the sum of $k$ and the dual Coxeter number of the corresponding Lie
algebra is not a nonnegative rational number, the vertex operator
algebra $L(k,0)$ is not finitely reductive, and, working with Lie
algebra theory rather than with vertex operator algebra theory,
Kazhdan and Lusztig constructed a natural braided tensor category
structure on a certain category of modules of level $k$ for the affine
Lie algebra (\cite{KL1}, \cite{KL2}, \cite{KL3}, \cite{KL4},
\cite{KL5}).  This work of Kazhdan-Lusztig in fact motivated the first
two authors to develop an analogous theory for vertex operator
algebras rather than for affine Lie algebras, as was explained in
detail in the introductory material in \cite{tensorAnnounce},
\cite{tensorK}, \cite{tensor1}, \cite{tensor2}, and \cite{HLaffine}.
However, this general theory, in its original form, did not apply to
Kazhdan-Lusztig's context, because the vertex-operator-algebra modules
considered in \cite{tensorAnnounce}, \cite{tensorK}, \cite{tensor1},
\cite{tensor2}, \cite{tensor3}, \cite{tensor4} are assumed to be the
direct sums of their weight spaces (with respect to $L(0)$), and the
non-semisimple modules considered by Kazhdan-Lusztig fail in general
to be the direct sums of their weight spaces.  Although their setup,
based on Lie theory, and ours, based on vertex operator algebra
theory, are very different (as was discussed in the introductory
material in our earlier papers), we expected to be able to recover
(and further extend) their results through our vertex operator
algebraic approach, which is very general, as we discussed above.
This motivated us, in the present work, to generalize the work of the
first two authors by considering modules with generalized weight
spaces, and especially, intertwining operators associated with such
generalized kinds of modules.  As we discuss below, this required us
to use logarithmic intertwining operators and logarithmic formal
calculus, and we have been able to construct braided tensor category
structure, and even vertex tensor category structure, on important
module categories that are not semisimple.  Using the present theory,
the third author (\cite{Z1}, \cite{Z2}) has indeed recovered the
braided tensor category structure of Kazhdan-Lusztig, and has also
extended it to vertex tensor category structure.

{}From the viewpoint of the general representation theory of vertex
operator algebras, it would be unnatural to study only semisimple
modules or only $L(0)$-semisimple modules; focusing only on such
modules would be analogous to focusing only on semisimple modules for
general (nonsemisimple) finite-dimensional Lie algebras.  And as we
have pointed out, working in this generality leads to logarithmic
structure; the general representation theory of vertex operator
algebras requires logarithmic structure.

Logarithmic structure in conformal field theory was in fact first
introduced by physicists to describe disorder phenomena \cite{Gu}.  A
lot of progress has been made on this subject.  We refer the interested
reader to the review articles \cite{Ga2}, \cite{Fl2}, \cite{RT} and
\cite{Fu}, and references therein.  One particularly interesting class
of logarithmic conformal field theories is the class associated to the
triplet $\mathcal{W}$-algebras of central charge
$1-6\frac{(p-1)^{2}}{p}$, which were introduced by Kausch \cite{K1} and
studied extensively both in physics and in mathematics by Flohr
\cite{Fl1} \cite{Fl2}, Gaberdiel-Kausch \cite{GK1} \cite{GK2}, Kausch
\cite{K2}, Fuchs-Hwang-Semikhatov-Tipunin \cite{FHST}, Abe \cite{A},
Feigin-Ga{\u\i}nutdinov-Semikhatov-Tipunin \cite{FGST1} \cite{FGST2}
\cite{FGST3}, Carqueville-Flohr \cite{CF}, Flohr-Gaberdiel \cite{FG},
Fuchs \cite{Fu}, Adamovi\'{c}-Milas \cite{AM1} \cite{AM2},
Flohr-Grabow-Koehn \cite{FGK} and Flohr-Knuth \cite{FK}.  A family of
$N=1$ triplet vertex operator superalgebras has been constructed and
studied recently by Adamovi\'{c} and Milas in \cite{AM3}. The paper
\cite{FHST} initiated a study of a possible generalization of the
Verlinde conjecture for rational conformal field theories to these
theories (see also the recent work \cite{FG}). The paper \cite{Fu}
assumed the existence of braided tensor category structures on the
categories of modules for the vertex operator algebras considered;
together with \cite{H13}, the present work gives a construction of these
structures. The paper \cite{CF} used the results in the present work as
announced in \cite{HLZ}. 

Here is how such logarithmic structure also arises naturally in the
representation theory of vertex operator algebras: In the construction
of intertwining operator algebras, the first author proved (see
\cite{diff-eqn}) that if modules for the vertex operator algebra
satisfy a certain cofiniteness condition, then products of the usual
intertwining operators satisfy certain systems of differential
equations with regular singular points. In addition, it was proved in
\cite{diff-eqn} that if the vertex operator algebra satisfies certain
finite reductivity conditions, then the analytic extensions of
products of the usual intertwining operators have no logarithmic
terms.  In the case when the vertex operator algebra satisfies only
the cofiniteness condition but not the finite reductivity conditions,
the products of intertwining operators still satisfy systems of
differential equations with regular singular points.  But in this
case, the analytic extensions of such products of intertwining
operators might have logarithmic terms. This means that if we want to
generalize the results in \cite{tensorAnnounce},
\cite{tensorK}--\cite{tensor3}, \cite{tensor4} and \cite{diff-eqn} to
the case in which the finite reductivity properties are not always
satisfied, we have to consider intertwining operators involving
logarithmic terms.

Logarithmic structure also appears naturally in modular invariance
results for vertex operator algebras and in the genus-one parts of
conformal field theories.  For a vertex operator algebra $V$
satisfying certain finiteness and reductivity conditions, Zhu proved
in \cite{Zhu} a modular invariance result for $q$-traces of products
of vertex operators associated to $V$-modules.  Zhu's result was
generalized to the case involving twisted vertex operators by Dong, Li
and Mason in \cite{DLM} and to the case of $q$-traces of products of
one intertwining operator and arbitrarily many vertex operators by
Miyamoto in \cite{M1}.  In \cite{M2}, Miyamoto generalized Zhu's
modular invariance result to a modular invariance result involving the
logarithm of $q$ for vertex operator algebras not necessarily
satisfying the reductivity condition.  In \cite{Hmodular}, for vertex
operator algebras satisfying certain finiteness and reductivity
conditions, by overcoming the difficulties one encounters if one tries
to generalize Zhu's method, the first author was able to prove the
modular invariance for $q$-traces of products and iterates of more
than one intertwining operator, using certain differential equations
and duality properties for intertwining operators.  If the vertex
operator algebra satisfies only Zhu's cofiniteness condition but not
the reductivity condition, the $q$-traces of products and iterates of
intertwining operators still satisfy the same differential equations,
but now they involve logarithms of all the variables.  To generalize
the general Verlinde conjecture proved in \cite{HVerlindeconjecture}
and the modular tensor category structure on the category of
$V$-modules obtained in \cite{rigidity}, one will need such general
logarithmic modular invariance.  See \cite{FHST} for research in this
direction.

In \cite{Mi}, Milas introduced and studied what he called
``logarithmic modules'' and ``logarithmic intertwining operators.''
See also \cite{Mi2}.
Roughly speaking, logarithmic modules are weak modules for a vertex
operator algebra that are direct sums of generalized eigenspaces for
the operator $L(0)$.  We will call such weak modules ``generalized
modules'' in this work.  Logarithmic intertwining operators are
operators that depend not only on powers of a (formal or complex)
variable $x$, but also on its logarithm $\log x$.

The special features of the logarithm function make the logarithmic
theory very subtle and interesting.  Although we show that all the
main theorems in the original tensor product theory developed by the
first two authors still hold in the logarithmic theory, many of the
proofs involve certain new techniques and have surprising connections
with certain combinatorial identities.

As we mentioned above, one important application of our generalization
is to the category ${\cal O}_\kappa$ of certain modules for an affine
Lie algebra studied by Kazhdan and Lusztig in their series of papers
\cite{KL1}--\cite{KL5}. It has been shown in \cite{Z1} and \cite{Z2}
by the third author that, among other things, all the conditions
needed to apply our theory to this module category are satisfied.  As
a result, it is proved in \cite{Z1} and \cite{Z2} that there is a
natural vertex tensor category structure on this module category,
giving in particular a new construction, in the context of general
vertex operator algebra theory, of the braided tensor category
structure on ${\cal O}_\kappa$.  The methods used in
\cite{KL1}--\cite{KL5} were very different.

The triplet $\mathcal{W}$-algebras belong to a different class of
vertex operator algebras, satisfying certain finiteness, boundedness
and reality conditions.  In this case, it has been shown in \cite{H13}
by the first author that all the conditions needed to apply the theory
carried out in the present work to the category of grading-restricted
modules for the vertex operator algebra are also satisfied.  Thus, by
the results obtained in this work, there is a natural vertex tensor
category structure on this category.

In addition to these logarithmic issues, another way in which the
present work generalizes the earlier tensor product theory for module
categories for a vertex operator algebra is that we now allow the
algebras to be somewhat more general than vertex operator algebras, in
order, for example, to accommodate module categories for the vertex
algebras $V_L$ where $L$ is a nondegenerate even lattice that is not
necessarily positive definite (cf.\ \cite{B}, \cite{DL}); see
\cite{Z1}.

What we accomplish in this work, then, is the following: We generalize
essentially all the results in \cite{tensor1}, \cite{tensor2},
\cite{tensor3} and \cite{tensor4} {}from the category of modules for a
vertex operator algebra to categories of suitably generalized modules
for a conformal vertex algebra or a M\"obius vertex algebra equipped
with an additional suitable grading by an abelian group.  The algebras
that we consider include not only vertex operator algebras but also
such vertex algebras as $V_L$ where $L$ is a nondegenerate even
lattice, and the modules that we consider are not required to be the
direct sums of their weight spaces but instead are required only to be
the (direct) sums of their ``generalized weight spaces,'' in a
suitable sense.  In particular, in this work we carry out, in the
present greater generality, the construction theory for the
``$P(z)$-tensor product'' functor originally done in \cite{tensor1},
\cite{tensor2} and \cite{tensor3} and the associativity theory for
this functor---the construction of the natural associativity
isomorphisms between suitable ``triple tensor products'' and the proof
of their important properties, including the isomorphism
property---originally done in \cite{tensor4}. This leads, as in
\cite{tensor5}, to the proof of the coherence properties for vertex
tensor categories, and in the M\"obius case, the coherence properties
for M\"obius vertex tensor categories.

The general structure of much of this work essentially follows that of
\cite{tensor1}, \cite{tensor2}, \cite{tensor3} and \cite{tensor4}.
However, the results here are much stronger and more general than in
these earlier works, and in addition, many of the results here have no
counterparts in those works.  Moreover, many ideas, formulations and
proofs in this work are necessarily quite different {}from those in the
earlier papers, and we have chosen to give some proofs that are new
even in the finitely reductive case studied in the earlier papers.

Some of the new ingredients that we are introducing into the theory
here are (as we shall explain in detail): an analysis of logarithmic
intertwining operators, including ``logarithmic formal calculus''; a
notion of ``$P(z_1,z_2)$-intertwining map'' and a study of its
properties; new ``compatibility conditions''; a generalization of the
result that the homogeneous components of the products and iterates of
intertwining maps span the appropriate tensor product modules; results
strengthening the relation between products and iterates of
intertwining maps; and a generalized sufficient condition for the
applicability of the theory, a condition that can be applied in the
case of our suitably generalized modules.

The contents of the sections of this work are as follows: In the rest
of this Introduction we compare classical tensor product theory for
Lie algebra modules with tensor product theory for vertex operator
algebra modules.  A crucial difference between the two theories is
that in the vertex operator algebra setting, the theory depends on an
``extra parameter'' $z$, which must be understood as a (nonzero)
complex variable rather than as a formal variable (although one needs,
and indeed we very heavily use, an extensive ``calculus of formal
variables'' in order to develop the theory).  In Section 2 we recall
some basic concepts in the theory of vertex (operator) algebras.  We
use the treatments of \cite{FLM2}, \cite{FHL}, \cite{DL} and
\cite{LL}; in particular, it is crucial in this tensor product theory
to use the formal-calculus point of view as developed in these works.
Readers can consult these works for further detail.  We also set up
the notation that will be used in this work, and we describe the main
category of (generalized) modules that we will consider. In Section 3
we introduce the notion of logarithmic intertwining operator as in
\cite{Mi} and present a detailed study of some of its properties.  In
Section 4 and 5 we present the definitions and constructions of
$P(z)$- and $Q(z)$-tensor products, generalizing those in
\cite{tensor1}, \cite{tensor2} and \cite{tensor3}. Some of the proofs
of results in Section 5 are given in Section 6. In Section 7 the
convergence condition introduced in \cite{tensor4} for constructing
the associativity isomorphism is given in the present context.  The
new notion of $P(z_1,z_2)$-intertwining map, generalizing the
corresponding concept in \cite{tensor4}, is introduced in Section 8.
This will play a crucial role in the construction of the associativity
isomorphisms.  In Section 9 we prove some properties that are
satisfied by vectors in the dual space of the vector-space tensor
product of three modules that arise {}from products and {}from iterates of
intertwining maps.  This enables us to define two crucial subspaces of
this dual space, by means of suitable compatibility and local grading
restriction conditions.  By relating these two subspaces, we construct
the associativity isomorphism in Section 10.  In Section 11, we
generalize a certain sufficient condition for the existence of
associativity isomorphisms in \cite{tensor4}, and we prove the
relevant conditions using differential equations.  In Section 12, we
establish the coherence of our tensor category.  

\subsection{The Lie algebra case}\label{LA}

It is heuristically useful to start by considering the tensor product
theory for modules for a Lie algebra in a somewhat unusual way---a way
that motivates our approach for the case of vertex algebras.

In the theory of tensor products for modules for a Lie algebra, the
tensor product of two modules is defined, or rather, constructed, as
the vector-space tensor product of the two modules, equipped with a
Lie algebra module action given by the familiar diagonal action of the
Lie algebra.  In the vertex algebra case, however, the vector-space
tensor product of two modules for a vertex algebra is {\it not} the
correct underlying vector space for the tensor product of the
vertex-algebra modules.  In this subsection we therefore consider
another approach to the tensor product theory for modules for a Lie
algebra---an approach, based on ``intertwining maps,'' that will show
how the theory proceeds in the vertex algebra case.  Then, in the next
subsection, we shall lay out the corresponding ``road map'' for the
tensor product theory in the vertex algebra case, which we then carry
out in the body of this work.

We first recall the following elementary but crucial background about
tensor product vector spaces: Given vector spaces $W_1$ and $W_2$, the
corresponding tensor product structure consists of a vector space $W_1
\otimes W_2$ equipped with a bilinear map
$$W_1 \times W_2 \longrightarrow W_1 \otimes W_2,$$
denoted
$$(w_{(1)},w_{(2)}) \mapsto w_{(1)}\otimes w_{(2)}$$
for $w_{(1)}\in W_1$ and $w_{(2)}\in W_2$, such that for any vector
space $W_3$ and any bilinear map
$$B:W_1 \times W_2 \longrightarrow W_3,$$
there is a unique linear map
$$L:W_1 \otimes W_2 \longrightarrow W_3$$
such that
$$B(w_{(1)},w_{(2)}) = L(w_{(1)}\otimes w_{(2)})$$
for $w_{(i)}\in W_i$, $i=1, 2$.  This universal property characterizes the
tensor product structure $W_1 \otimes W_2$, equipped with its bilinear
map $\cdot \otimes \cdot$, up to unique isomorphism.  In addition, the
tensor product structure in fact exists.

As was illustrated in \cite{tensorK}, and as is well known, the notion
of tensor product of modules for a Lie algebra can be formulated in
terms of what can be called ``intertwining maps'': Let $W_1$, $W_2$,
$W_3$ be modules for a fixed Lie algebra $V$.  (We are calling our Lie
algebra $V$ because we shall be calling our vertex algebra $V$, and we
would like to emphasize the analogies between the two theories.)  An
{\it intertwining map of type ${W_3 \choose {W_1 W_2}}$} is a linear
map $I: W_1 \otimes W_2 \longrightarrow W_3$ (or equivalently, {}from
what we have just recalled, a bilinear map $W_1 \times W_2
\longrightarrow W_3$) such that
\begin{equation}\label{intwmap}
\pi_3 (v)I(w_{(1)}\otimes w_{(2)})=I(\pi_1 (v)w_{(1)}\otimes
w_{(2)})+I(w_{(1)}\otimes \pi_2 (v)w_{(2)})
\end{equation}
for $v\in V$ and $w_{(i)}\in W_i$, $i=1,2$, where $\pi_1, \pi_2,
\pi_3$ are the module actions of $V$ on $W_1$, $W_2$ and $W_3$,
respectively.  (Clearly, such an intertwining map is the same as a
module map {}from $W_1\otimes W_2$, equipped with the tensor product
module structure, to $W_3$, but we are now temporarily ``forgetting''
what the tensor product module is.)

A {\it tensor product of the $V$-modules $W_1$ and $W_2$} is then a
pair $(W_0,I_0)$, where $W_0$ is a $V$-module and $I_0$ is an
intertwining map of type ${W_0 \choose {W_1 W_2}}$ (which, again,
could be viewed as a suitable bilinear map $W_1 \times W_2
\longrightarrow W_0$), such that for any pair $(W,I)$ with $W$ a
$V$-module and $I$ an intertwining map of type ${W \choose {W_1
W_2}}$, there is a unique module homomorphism $\eta:
W_0\longrightarrow W$ such that $I=\eta \circ I_0$. This universal
property of course characterizes $(W_0, I_0)$ up to canonical
isomorphism.  Moreover, it is obvious that the tensor product in fact
exists, and may be constructed as the vector-space tensor product
$W_1\otimes W_2$ equipped with the diagonal action of the Lie algebra,
together with the identity map {}from $W_1\otimes W_2$ to itself (or
equivalently, the canonical bilinear map $W_1 \times W_2
\longrightarrow W_1 \otimes W_2$).  We shall denote the tensor product
$(W_0,I_0)$ of $W_1$ and $W_2$ by $(W_1 \boxtimes W_2, \boxtimes)$,
where it is understood that the image of $w_{(1)}\otimes w_{(2)}$
under our canonical intertwining map $\boxtimes$ is $w_{(1)}\boxtimes
w_{(2)}$.  Thus $W_1 \boxtimes W_2=W_1\otimes W_2$, and under our
identifications, $\boxtimes = 1_{W_1\otimes W_2}$.

\begin{rema}
{\rm This classical explicit construction of course shows that the
tensor product functor exists for the category of modules for a Lie
algebra.  For vertex algebras, it will be relatively straightforward
to {\it define} the appropriate tensor product functor(s) (see
\cite{tensorK}, \cite{tensor1}, \cite{tensor2}, \cite{tensor3}), but
it will be a nontrivial matter to {\it construct} this functor (or
more precisely, these functors) and thereby prove that the
(appropriate) tensor product of modules for a (suitable) vertex
algebra exists.  The reason why we have formulated the notion of
tensor product module for a Lie algebra in the way that we just did is
that this formulation motivates the correct notion of tensor product
functor(s) in the vertex algebra case.}
\end{rema}

\begin{rema}
{\rm Using this explicit construction of the tensor product functor
and our notation $w_{(1)}\boxtimes w_{(2)}$ for the tensor product of
elements, the standard natural associativity isomorphisms among tensor
products of triples of Lie algebra modules are expressed as follows:
Since $w_{(1)}\boxtimes w_{(2)} = w_{(1)}\otimes w_{(2)}$, we have
\begin{eqnarray}\
(w_{(1)}\boxtimes w_{(2)})\boxtimes w_{(3)}&=&
(w_{(1)}\otimes w_{(2)})\otimes w_{(3)},\nno\\
w_{(1)}\boxtimes(w_{(2)}\boxtimes w_{(3)})&=&
w_{(1)}\otimes(w_{(2)}\otimes w_{(3)})\nno
\end{eqnarray}
for $w_{(i)}\in W_i$, $i=1,2,3$, and so
the canonical identification between
$w_{(1)}\otimes(w_{(2)}\otimes w_{(3)})$ and $(w_{(1)}\otimes
w_{(2)})\otimes w_{(3)}$ gives the standard natural isomorphism
\begin{eqnarray}
(W_1\boxtimes W_2)\boxtimes W_3 &\to & W_1\boxtimes (W_2\boxtimes W_3)\nno\\
(w_{(1)}\boxtimes w_{(2)})\boxtimes w_{(3)} &\mapsto &
w_{(1)}\boxtimes(w_{(2)}\boxtimes w_{(3)}).\label{elemap}
\end{eqnarray}
This collection of natural associativity isomorphisms of course
satisfies the classical coherence conditions for associativity
isomorphisms among multiple nested tensor product modules---the
conditions that say that in nested tensor products involving any
number of tensor factors, the placement of parentheses (as in
(\ref{elemap}), the case of three tensor factors) is immaterial; we
shall discuss coherence conditions in detail later.  Now, as was
discovered in \cite{tensor4}, it turns out that maps analogous to
(\ref{elemap}) can also be constructed in the vertex algebra case,
giving natural associativity isomorphisms among triples of modules for
a (suitable) vertex operator algebra.  However, in the vertex algebra
case, the elements ``$w_{(1)}\boxtimes w_{(2)}$,'' which indeed exist
(under suitable conditions) and are constructed in the theory, lie in
a suitable ``completion'' of the tensor product module rather than in
the module itself.  Correspondingly, it is a nontrivial matter to
construct the triple-tensor-product elements on the two sides of
(\ref{elemap}); in fact, one needs to prove certain convergence, under
suitable additional conditions.  Even after the triple-tensor-product
elements are constructed (in suitable completions of the
triple-tensor-product modules), it is a delicate matter to construct
the appropriate natural associativity maps, analogous to
(\ref{elemap}), to prove that they are well defined, and to prove that
they are isomorphisms.  In the present work, we shall generalize these
matters (in a self-contained way) {}from the context of \cite{tensor4}
to a more general one.  In the rest of this subsection, for triples of
modules for a Lie algebra, we shall now describe a construction of the
natural associativity isomorphisms that will seem roundabout and
indirect, but this is the method of construction of these isomorphisms
that will give us the correct ``road map'' for the corresponding
construction (and theorems) in the vertex algebra case, as in
\cite{tensor1}, \cite{tensor2}, \cite{tensor3} and \cite{tensor4}.}
\end{rema}

A significant feature of the constructions in the earlier works (and
in the present work) is that the tensor product of modules $W_1$ and
$W_2$ for a vertex operator algebra $V$ is the contragredient module
of a certain $V$-module that is typically a {\it proper} subspace of
$(W_1\otimes W_2)^*$, the dual space of the vector-space tensor
product of $W_1$ and $W_2$.  In particular, our treatment, which
follows, of the Lie algebra case will use contragredient modules, and
we will therefore restrict our attention to {\it finite-dimensional}
modules for our Lie algebra.  It will be important that the
double-contragredient module of a Lie algebra module is naturally
isomorphic to the original module.  We shall sometimes denote the
contragredient module of a $V$-module $W$ by $W'$, so that $W''=W$.
(We recall that for a module $W$ for a Lie algebra $V$, the
corresponding contragredient module $W'$ consists of the dual vector
space $W^*$ equipped with the action of $V$ given by: $(v \cdot
w^*)(w) = - w^*(v \cdot w)$ for $v \in V$, $w^* \in W^*$, $w \in W$.)

Let us, then, now restrict our attention to finite-dimensional modules
for our Lie algebra $V$.  The dual space $(W_1\otimes W_2)^*$ carries
the structure of the classical contragredient module of the tensor
product module.  Given any intertwining map of type ${W_3 \choose {W_1
W_2}}$, using the natural linear isomorphism
\begin{equation}\label{corpd2}
\hom (W_1\otimes W_2, W_3)\tilde{\longrightarrow} \hom (W^*_3,
(W_1\otimes W_2)^*)
\end{equation}
we have a corresponding linear map in $\hom (W^*_3, (W_1\otimes
W_2)^*)$, and this must be a map of $V$-modules.  In the vertex
algebra case, given $V$-modules $W_1$ and $W_2$, it turns out that
with a suitable analogous setup, the union in the vector space
$(W_1\otimes W_2)^*$ of the ranges of all such $V$-module maps, as
$W_3$ and the intertwining map vary (and with $W^*_3$ replaced by the
contragredient module $W'_3$), is the correct candidate for the
contragredient module of the tensor product module $W_1\boxtimes W_2$.
Of course, in the Lie algebra situation, this union is $(W_1\otimes
W_2)^*$ itself (since we are allowed to take $W_3=W_1\otimes W_2$ and
the intertwining map to be the canonical map), but in the vertex
algebra case, this union is typically much smaller than $(W_1\otimes
W_2)^*$.  In the vertex algebra case, we will use the notation $W_1
\hboxtr \, W_2$ to designate this union, and if the tensor product
module $W_1\boxtimes W_2$ in fact exists, then
\begin{eqnarray}
W_1 \boxtimes W_2 &=& (W_1 \,\hboxtr \; W_2)',\label{LAhbox1}\\
W_1 \,\hboxtr \; W_2 &=& (W_1 \boxtimes W_2)'.\label{LAhbox2}
\end{eqnarray}
Thus in the Lie algebra case we will write
\begin{eqnarray}
W_1 \,\hboxtr \; W_2 = (W_1 \otimes W_2)^*,
\end{eqnarray}
and (\ref{LAhbox1}) and (\ref{LAhbox2}) hold.  (In the Lie algebra
case we prefer to write $(W_1 \otimes W_2)^*$ rather than $(W_1
\otimes W_2)'$, because in the vertex algebra case, $W_1 \otimes W_2$
is typically not a $V$-module, and so we will not be allowed to write
$(W_1 \otimes W_2)'$ in the vertex algebra case.)

The subspace $W_1 \,\hboxtr \; W_2$ of $(W_1\otimes W_2)^*$ was in
fact further described in the following terms in \cite{tensor1} and
\cite{tensor3}, in the vertex algebra case: For any map in $\hom
(W'_3, (W_1\otimes W_2)^*)$ corresponding to an intertwining map
according to (\ref{corpd2}), the image of any $w'_{(3)}\in W'_3$ under
this map satisfies certain subtle conditions, called the
``compatibility condition'' and the ``local grading restriction
condition''; these conditions are not ``visible'' in the Lie algebra
case.  These conditions in fact precisely describe the proper subspace
$W_1 \hboxtr \, W_2$ of $(W_1\otimes W_2)^*$.  We will discuss such
conditions further in Subsection 1.4 and in the body of this work.  As we
shall explain, this idea of describing elements in certain dual spaces
was also used in constructing the natural associativity isomorphisms
between triples of modules for a vertex operator algebra in
\cite{tensor4}.

In order to give the reader a guide to the vertex algebra case, we now
describe the analogue for the Lie algebra case of this construction of
the associativity isomorphisms.  To construct the associativity
isomorphism {}from $(W_1 \boxtimes W_2)\boxtimes W_3$ to $W_1\boxtimes
(W_2\boxtimes W_3)$, it is equivalent (by duality) to give a suitable
isomorphism {}from $W_1\hboxtr\, (W_2\boxtimes W_3)$ to $(W_1 \boxtimes
W_2)\, \hboxtr\, W_3$ (recall (\ref{LAhbox1}), (\ref{LAhbox2})).

Rather than directly constructing an isomorphism between these two
$V$-modules, it turns out that we want to embed both of them,
separately, into the single space $(W_1\otimes W_2 \otimes W_3)^*$.
Note that $(W_1\otimes W_2 \otimes W_3)^*$
is naturally a $V$-module, via the contragredient of the
diagonal action, that is,
\begin{eqnarray}\label{actiononW1W2W2*}
(\pi(v)\lambda)(w_{(1)}\otimes w_{(2)}\otimes w_{(3)})&=&
-\lambda(\pi_1(v)w_{(1)}\otimes w_{(2)}\otimes w_{(3)})\nno\\
& -&\lambda(w_{(1)}\otimes \pi_2(v)w_{(2)}\otimes w_{(3)})\nno\\
& -&\lambda(w_{(1)}\otimes w_{(2)}\otimes \pi_3(v)w_{(3)}),
\end{eqnarray}
for $v\in V$ and $w_{(i)}\in W_i$, $i=1,2,3$, where $\pi_1, \pi_2,
\pi_3$ are the module actions of $V$ on $W_1$, $W_2$ and $W_3$,
respectively.  A concept related to this is the notion of
{\it intertwining map {}from $W_1 \otimes W_2 \otimes W_3$ to a module
$W_4$}, a natural analogue of (\ref{intwmap}), defined to be a linear map
\begin{equation}\label{intwmap3}
F: W_1\otimes W_2 \otimes W_3 \longrightarrow W_4
\end{equation}
such that
\begin{eqnarray}\label{intwmapfor3} \pi_4(v)F(w_{(1)}\otimes w_{(2)}
\otimes w_{(3)})&=&F(\pi_1(v)w_{(1)} \otimes w_{(2)} \otimes
w_{(3)})\nno\\ &+&F(w_{(1)} \otimes \pi_2(v)w_{(2)} \otimes
w_{(3)})\nno\\ &+&F(w_{(1)} \otimes w_{(2)} \otimes
\pi_{(3)}(v)w_3),\label{11i}
\end{eqnarray}
with the obvious notation.  The
relation between (\ref{actiononW1W2W2*}) and (\ref{intwmapfor3})
comes directly {}from the natural linear isomorphism
\begin{equation}
\hom (W_1\otimes W_2\otimes W_3, W_4) \tilde{\longrightarrow} \hom
(W^*_4, (W_1\otimes W_2\otimes W_3)^*);
\end{equation}
given $F$, we have
\begin{eqnarray}
W^*_4 &\longrightarrow &(W_1\otimes W_2\otimes
W_3)^*\nno\\
\nu&\mapsto &\nu\circ F.
\end{eqnarray}
Under this natural linear isomorphism, the intertwining maps
correspond precisely to the $V$-module maps {}from $W^*_4$ to
$(W_1\otimes W_2\otimes W_3)^*$.  In the situation for vertex
algebras, as was the case for tensor products of two rather than three
modules, there are analogues of all of the notions and comments
discussed in this paragraph {\it except that we will not put
$V$-module structure onto the vector space} $W_1\otimes W_2\otimes
W_3$; as we have emphasized, we will instead base the theory on
intertwining maps.

Two important ways of constructing maps of the type (\ref{intwmap3})
are as follows: For modules $W_1$, $W_2$, $W_3$, $W_4$, $M_1$ and
intertwining maps $I_1$ and $I_2$ of types ${W_4 \choose {W_1 M_1}}$
and ${M_1 \choose {W_2 W_3}}$, respectively, by definition the
composition $I_1\circ (1_{W_1}\otimes I_2)$ is an intertwining map
{}from $W_1\otimes W_2 \otimes W_3$ to $W_4$. Similarly, for
intertwining maps $I^1$, $I^2$ of types ${W_4 \choose {M_2 W_3}}$ and
${M_2 \choose {W_1 W_2}}$, respectively, the composition $I^1\circ
(I^2\otimes 1_{W_3})$ is an intertwining map {}from $W_1\otimes W_2
\otimes W_3$ to $W_4$. Hence we have two $V$-module homomorphisms
\begin{eqnarray}
W^*_4 &\longrightarrow &(W_1\otimes W_2\otimes
W_3)^*\nno\\
\nu&\mapsto &\nu\circ F_1,
\label{injint1}
\end{eqnarray}
where $F_1$ is the intertwining map $I_1\circ (1_{W_1}\otimes I_2)$;
and
\begin{eqnarray}
W^*_4 &\longrightarrow & (W_1\otimes W_2\otimes
W_3)^*\nno\\
\nu&\mapsto &\nu\circ F_2,
\label{injint2}
\end{eqnarray}
where $F_2$ is the intertwining map $I^1\circ (I^2\circ 1_{W_3})$.

The special cases in which the modules $W_4$ are two iterated tensor
product modules and the ``intermediate'' modules $M_1$ and $M_2$ are
two tensor product modules are particularly interesting: When
$W_4=W_1\boxtimes (W_2\boxtimes W_3)$ and $M_1=W_2\boxtimes W_3$, and
$I_1$ and $I_2$ are the corresponding canonical intertwining maps,
(\ref{injint1}) gives the natural $V$-module homomorphism
\begin{eqnarray}
W_1\,\hboxtr \,(W_2 \boxtimes W_3)&\longrightarrow &(W_1\otimes
W_2\otimes W_3)^*\nno\\
\nu&\mapsto &(w_{(1)}\otimes w_{(2)}\otimes w_{(3)}\mapsto \nu
(w_{(1)}\boxtimes (w_{(2)}\boxtimes w_{(3)})));\nno\\
\label{inj1}
\end{eqnarray}
when $W_4=(W_1\boxtimes W_2)\boxtimes W_3$ and $M_2=W_1\boxtimes W_2$,
and $I^1$ and $I^2$ are the corresponding canonical intertwining maps,
(\ref{injint2}) gives the natural $V$-module homomorphism
\begin{eqnarray}
(W_1 \boxtimes W_2)\,\hboxtr \,W_3&\longrightarrow & (W_1\otimes
W_2\otimes W_3)^*\nno\\
\nu&\mapsto &(w_{(1)}\otimes w_{(2)}\otimes w_{(3)}\mapsto \nu
((w_{(1)}\boxtimes w_{(2)})\boxtimes w_{(3)})).\nno\\
\label{inj2}
\end{eqnarray}

Clearly, in our Lie algebra case, both of the maps (\ref{inj1}) and
(\ref{inj2}) are isomorphisms, since they both in fact amount to the
identity map on $(W_1\otimes W_2\otimes W_3)^*$.  However, in the
vertex algebra case the analogues of these two maps are only injective
homomorphisms, and typically not isomorphisms.  (Recall the analogous
situation, mentioned above, for double rather than triple tensor
products.)  These two maps enable us to identify both $W_1\,\hboxtr
\,(W_2 \boxtimes W_3)$ and $(W_1 \boxtimes W_2)\,\hboxtr \,W_3$ with
subspaces of $(W_1\otimes W_2\otimes W_3)^*$.  In the vertex algebra
case we will have certain ``compatibility conditions'' and ``local
grading restriction conditions'' on elements of $(W_1\otimes
W_2\otimes W_3)^*$ to describe each of the two subspaces.  In either
the Lie algebra or the vertex algebra case, the construction of our
desired natural associativity isomorphism between the two modules
$(W_1 \boxtimes W_2)\boxtimes W_3$ and $W_1\boxtimes(W_2 \boxtimes
W_3)$ follows {}from showing that the ranges of homomorphisms
(\ref{inj1}) and (\ref{inj2}) are equal to each other, which is of
course obvious in the Lie algebra case since both (\ref{inj1}) and
(\ref{inj2}) are isomorphisms to $(W_1\otimes W_2\otimes W_3)^*$.  It
turns out that, under this associativity isomorphism, (\ref{elemap})
holds in both the Lie algebra case and the vertex algebra case; in the
Lie algebra case, this is obvious because all the maps are the
``tautological'' ones.

Now we give the reader a preview of how, in the vertex algebra case,
these compatibility and local grading restriction conditions on
elements of $(W_1\otimes W_2\otimes W_3)^*$ will arise.  As we have
mentioned, in the Lie algebra case, an intertwining map {}from
$W_1\otimes W_2 \otimes W_3$ to $W_4$ corresponds to a module map {}from
$W^*_4$ to $(W_1\otimes W_2\otimes W_3)^*$.  As was discussed in
\cite{tensor4}, for the vertex operator algebra analogue, the image of
any $w'_{(4)}\in W'_4$ under such an analogous map satisfies certain
``compatibility'' and ``local grading restriction'' conditions, and so
these conditions must be satisfied by those elements of $(W_1\otimes
W_2\otimes W_3)^*$ lying in the ranges of the vertex-operator-algebra
analogues of either of the maps (\ref{inj1}) and (\ref{inj2}) (or the
maps (\ref{injint1}) and (\ref{injint2})).

Besides these two conditions, satisfied by the elements of the ranges
of the maps of both types (\ref{inj1}) and (\ref{inj2}), the elements
of the ranges of the analogues of the homomorphisms (\ref{inj1}) and
(\ref{inj2}) have their own separate properties.  First note that any
$\lambda\in (W_1\otimes W_2\otimes W_3)^*$ induces the two maps
\begin{eqnarray}\label{mu1}
\mu^{(1)}_\lambda: W_1 &\to & (W_2\otimes W_3)^*\nno\\
 w_{(1)}&\mapsto & \lambda(w_{(1)} \otimes \cdot \otimes \cdot)
\end{eqnarray}
and
\begin{eqnarray}\label{mu2}
\mu^{(2)}_\lambda: W_3 &\to & (W_1\otimes W_2)^*\nno\\
w_{(3)}& \mapsto &\lambda(\cdot \otimes \cdot \otimes w_{(3)}).
\end{eqnarray}
In the vertex operator algebra analogue \cite{tensor4}, if $\lambda$
lies in the range of (\ref{inj1}), then it must satisfy the condition
that the elements $\mu^{(1)}_\lambda(w_{(1)})$ all lie in a suitable
completion of the subspace $W_2\hboxtr\, W_3$ of $(W_2\otimes W_3)^*$,
and if $\lambda$ lies in the range of (\ref{inj2}), then it must
satisfy the condition that the elements $\mu^{(2)}_\lambda(w_{(3)})$
all lie in a suitable completion of the subspace $W_1\hboxtr\, W_2$ of
$(W_1\otimes W_2)^*$.  (Of course in the Lie algebra case, these
statements are tautological.)  In \cite{tensor4}, these important
conditions, that $\mu^{(1)}_\lambda(W_1)$ lies in a suitable
completion of $W_2\hboxtr\, W_3$ and that $\mu^{(2)}_\lambda(W_3)$
lies in a suitable completion of $W_1\hboxtr\, W_2$, are understood as
``local grading restriction conditions'' with respect to the two
different ways of composing intertwining maps.

In the construction of our desired natural associativity isomorphism,
since we want the ranges of (\ref{inj1}) and (\ref{inj2}) to be the
same submodule of $(W_1\otimes W_2\otimes W_3)^*$, the ranges of both
(\ref{inj1}) and (\ref{inj2}) should satisfy both of these conditions.
This amounts to a certain ``extension condition'' in the vertex
algebra case.  When all these conditions are satisfied, it can in fact
be proved \cite{tensor4} that the associativity isomorphism does
indeed exist and that in addition, the ``associativity of intertwining
maps'' holds; that is, the ``product'' of two suitable intertwining
maps can be written, in a certain sense, as the ``iterate'' of two
suitable intertwining maps, and conversely.  This equality of products
with iterates, highly nontrivial in the vertex algebra case, amounts
in the Lie algebra case to the easy statement that in the notation
above, any intertwining map of the form $I_1\circ (1_{W_1}\otimes
I_2)$ can also be written as an intertwining map of the form $I^1\circ
(I^2\otimes 1_{W_3})$, for a suitable ``intermediate module'' $M_2$
and suitable intertwining maps $I^1$ and $I^2$, and conversely.  The
reason why this statement is easy in the Lie algebra case is that in
fact {\it any} intertwining map $F$ of the type (\ref{intwmap3}) can
be ``factored'' in either of these two ways; for example, to write $F$
in the form $I_1\circ (1_{W_1}\otimes I_2)$, take $M_1$ to be
$W_2\otimes W_3$, $I_2$ to be the canonical (identity) map and $I_1$
to be $F$ itself (with the appropriate identifications having been
made).

We are now ready to discuss the vertex algebra case.

\subsection{The vertex algebra case}

In this subsection, which should be carefully compared with the
previous one, we shall lay out our ``road map'' of the constructions of
the tensor product functors and the associativity isomorphisms for a
suitable class of vertex algebras, generalizing, and also following
the ideas of, the corresponding theory developed in \cite{tensor1},
\cite{tensor2}, \cite{tensor3} and \cite{tensor4} for vertex operator
algebras.  Without yet specifying the precise class of vertex algebras
that we shall be using in the body of this work, except to say that
our vertex algebras will be $\mathbb{Z}$-graded and our modules will be
$\mathbb{C}$-graded, we now discuss the vertex algebra case.

In this case, the concept of intertwining map involves the moduli
space of Riemann spheres with one negatively oriented puncture and two
positively oriented punctures and with local coordinates around each
puncture; the details of the geometric structures needed in this
theory are presented in \cite{H1}.  For each element of this moduli
space there is a notion of intertwining map adapted to the particular
element.  Let $z$ be a nonzero complex number and let $P(z)$ be the
Riemann sphere $\hat{\mathbb C}$ with one negatively oriented puncture
at $\infty$ and two positively oriented punctures at $z$ and $0$, with
local coordinates $1/w$, $w-z$ and $w$ at these three punctures,
respectively.

Let $V$ be a vertex algebra (on which appropriate assumptions,
including the existence of a suitable ${\mathbb Z}$-grading, will be made
later), and let $Y(\cdot,x)$ be the vertex operator map defining the
algebra structure (see Section 2 below for a brief summary of basic
notions and notation, including the formal delta function).  Let
$W_1$, $W_2$ and $W_3$ be modules for $V$, and let $Y_1(\cdot,x)$,
$Y_2(\cdot,x)$ and $Y_3(\cdot,x)$ be the corresponding vertex operator
maps.  (The cases in which some of the $W_i$ are $V$ itself, and some
of the $Y_i$ are, correspondingly, $Y$, are important, but the most
interesting cases are those where all three modules are different {}from
$V$.)  A ``$P(z)$-intertwining map of type ${W_3 \choose {W_1 W_2}}$''
is a linear map
\begin{equation}
I: W_1 \otimes W_2 \longrightarrow \overline{W}_3,
\end{equation}
where $\overline{W}_3$ is a certain completion of $W_3$, related to
its ${\mathbb C}$-grading, such that
\begin{eqnarray}\label{im-jacobi}
\lefteqn{x_0^{-1}\delta\left(\frac{ x_1-z}{x_0}\right)
Y_3(v, x_1)I(w_{(1)}\otimes w_{(2)})}\nno\\
&&=z^{-1}\delta\left(\frac{x_1-x_0}{z}\right)
I(Y_1(v, x_0)w_{(1)}\otimes w_{(2)})\nno\\
&&\hspace{2em}+x_0^{-1}\delta\left(\frac{z-x_1}{-x_0}\right)
I(w_{(1)}\otimes Y_2(v, x_1)w_{(2)})
\end{eqnarray}
for $v\in V$, $w_{(1)}\in W_1$, $w_{(2)}\in W_2$, where $x_0$,
$x_1$ and $x_2$ are commuting independent formal variables.  This
notion is motivated in detail in \cite{tensorK}, \cite{tensor1} and
\cite{tensor3}; we shall recall the motivation below.

\begin{rema}\label{formalandcomplexvariables}
{\rm In this theory, it is crucial to distinguish between formal
variables and complex variables. Thus we shall use the following
notational convention: {\it Throughout this work, unless we specify
otherwise, the symbols $x$, $x_0$, $x_1$, $x_2$, $\dots$ , $y$, $y_0$,
$y_1$, $y_2$, $\dots$ will denote commuting independent formal
variables, and by contrast, the symbols $z$, $z_0$, $z_1$, $z_2$,
$\dots$ will denote complex numbers in specified domains, not formal
variables.}}
\end{rema}

\begin{rema}\label{im-io}{\rm
Recall {}from \cite{FHL} the definition of the notion of intertwining
operator ${\cal Y}(\cdot, x)$ in the theory of vertex (operator)
algebras.  Given $(W_1,Y_1)$, $(W_2,Y_2)$ and $(W_3,Y_3)$ as above, an
intertwining operator of type ${W_3 \choose {W_1 W_2}}$ can be viewed
as a certain type of linear map ${\cal Y}(\cdot,x)\cdot$ {}from $W_1
\otimes W_2$ to the vector space of formal series in $x$ of the form
$\sum_{n \in {\mathbb C}} w(n) x^n$, where the coefficients $w(n)$ lie in
$W_3$, and where we are allowing arbitrary complex powers of $x$,
suitably ``truncated {}from below'' in this sum.  The main property of an
intertwining operator is the following ``Jacobi identity'':
\begin{eqnarray}\label{io-jacobi}
\lefteqn{\dps x^{-1}_0\delta \left( {x_1-x_2\over x_0}\right)
Y_3(v,x_1){\cal Y}(w_{(1)},x_2)w_{(2)}}\nno\\
&&\hspace{2em}- x^{-1}_0\delta \left( {x_2-x_1\over -x_0}\right)
{\cal Y}(w_{(1)},x_2)Y_2(v,x_1)w_{(2)}\nno \\
&&{\dps = x^{-1}_2\delta \left( {x_1-x_0\over x_2}\right)
{\cal Y}(Y_1(v,x_0)w_{(1)},x_2)
w_{(2)}}
\end{eqnarray}
for $v\in V$, $w_{(1)}\in W_1$ and $w_{(2)}\in W_2$.  (When all three
modules $W_i$ are $V$ itself and all four operators $Y_i$ and ${\cal
Y}$ are $Y$ itself, (\ref{io-jacobi}) becomes the usual Jacobi
identity in the definition of the notion of vertex algebra.  When
$W_1$ is $V$, $W_2=W_3$ and ${\cal Y}=Y_2=Y_3$, (\ref{io-jacobi})
becomes the usual Jacobi identity in the definition of the notion of
$V$-module.)  The point is that by ``substituting $z$ for $x_2$'' in
(\ref{io-jacobi}), we obtain (\ref{im-jacobi}), where we make the
identification
\begin{eqnarray}\label{intwmap=intwopatz}
I(w_{(1)} \otimes w_{(2)}) = {\cal Y}(w_{(1)},z)w_{(2)};
\end{eqnarray}
the resulting complex powers of the complex number $z$ are made
precise by a choice of branch of the $\log$ function.  The nonzero
complex number $z$ in the notion of $P(z)$-intertwining map thus
``comes {}from'' the substitution of $z$ for $x_2$ in the Jacobi
identity in the definition of the notion of intertwining operator.  In
fact, this correspondence (given a choice of branch of $\log$)
actually defines an isomorphism between the space of
$P(z)$-intertwining maps and the space of intertwining operators of
the same type (\cite{tensor1}, \cite{tensor3}); this will be discussed
below.}
\end{rema}

There is a natural linear injection
\begin{equation}\label{homva}
\hom (W_1\otimes W_2, \overline{W}_3)\longrightarrow
\hom (W'_3, (W_1\otimes W_2)^*),
\end{equation}
where here and below we denote by $W'$ the (suitably defined)
contragredient module of a $V$-module $W$; we have $W''=W$.  Under
this injection, a map $I\in \hom (W_1\otimes W_2, \overline{W}_3)$
amounts to a map $I':
W'_3\longrightarrow (W_1\otimes W_2)^*$:
\begin{equation}\label{I'}
w'_{(3)} \mapsto \langle w'_{(3)}, I(\cdot\otimes \cdot)\rangle,
\end{equation}
where $\langle \cdot,\cdot \rangle$ denotes the natural pairing
between the contragredient of a module and its completion. If $I$ is a
$P(z)$-intertwining map, then as in the Lie algebra case (see above),
where such a map is a module map, the map (\ref{I'}) intertwines two
natural $V$-actions on $W'_3$ and $(W_1\otimes W_2)^*$. We will see
that in the present (vertex algebra) case, $(W_1\otimes W_2)^*$ is
typically not a $V$-module.  The images of all the elements $w'_{(3)}\in
W'_3$ under this map satisfy certain conditions, called the
``$P(z)$-compatibility condition'' and the ``$P(z)$-local grading
restriction condition,'' as formulated in \cite{tensor1} and
\cite{tensor3}; we shall discuss these below.

Given a category of $V$-modules and two modules $W_1$ and $W_2$ in
this category, as in the Lie algebra case, the ``$P(z)$-tensor product
of $W_1$ and $W_2$'' is then defined to be a pair $(W_0,I_0)$, where
$W_0$ is a module in the category and $I_0$ is a $P(z)$-intertwining
map of type ${W_0 \choose {W_1 W_2}}$, such that for any pair $(W,I)$
with $W$ a module in the category and $I$ a $P(z)$-intertwining map of
type ${W \choose {W_1 W_2}}$, there is a unique morphism $\eta:
W_0\longrightarrow W$ such that $I=\bar\eta \circ I_0$; here and
throughout this work we denote by $\bar\chi$ the linear map naturally
extending a suitable linear map $\chi$ {}from a graded space to its
appropriate completion. This universal property characterizes $(W_0,
I_0)$ up to canonical isomorphism, {\it if it exists}.  We will denote
the $P(z)$-tensor product of $W_1$ and $W_2$, if it exists, by
$(W_1\boxtimes_{P(z)} W_2, \boxtimes_{P(z)})$, and we will denote the
image of $w_{(1)}\otimes w_{(2)}$ under $\boxtimes_{P(z)}$ by
$w_{(1)}\boxtimes_{P(z)} w_{(2)}$, which is an element of
$\overline{W_1\boxtimes_{P(z)} W_2}$, not of $W_1\boxtimes_{P(z)}
W_2$.

{}From this definition and the natural map (\ref{homva}), we will see
that if the $P(z)$-tensor product of $W_1$ and $W_2$ exists, then its
contragredient module can be realized as the union of ranges of all
maps of the form (\ref{I'}) as $W'_3$ and $I$ vary.  Even if the
$P(z)$-tensor product of $W_1$ and $W_2$ does not exist, we denote
this union (which is always a subspace stable under a natural action
of $V$) by $W_1\hboxtr_{P(z)} W_2$.  If the tensor product does exist,
then
\begin{eqnarray}
W_1 \boxtimes_{P(z)} W_2 &=& (W_1 \hboxtr_{P(z)}
W_2)',\label{vertexhbox1}\\
W_1 \hboxtr_{P(z)} W_2 &=& (W_1 \boxtimes_{P(z)}
W_2)';\label{vertexhbox2}
\end{eqnarray}
examining (\ref{vertexhbox1}) will show the reader why the notation
$\hboxtr$ was chosen in the earlier papers ($\boxtimes = \hboxtr
\,'$!).  Several critical facts about $W_1\hboxtr_{P(z)} W_2$ were
proved in \cite{tensor1}, \cite{tensor2} and \cite{tensor3}, notably,
$W_1\hboxtr_{P(z)} W_2$ is equal to the subspace of $(W_1\otimes
W_2)^*$ consisting of all the elements satisfying the
$P(z)$-compatibility condition and the $P(z)$-local grading
restriction condition, and in particular, this subspace is $V$-stable;
and the condition that $W_1\hboxtr_{P(z)} W_2$ is a module is
equivalent to the existence of the $P(z)$-tensor product
$W_1\boxtimes_{P(z)} W_2$.  All these facts will be proved below.

In order to construct vertex tensor category structure, we need to
construct appropriate natural associativity isomorphisms.  Assuming
the existence of the relevant tensor products, we in fact need to
construct an appropriate natural isomorphism {}from $(W_1
\boxtimes_{P(z_1-z_2)} W_2)\boxtimes_{P(z_2)} W_3$ to
$W_1\boxtimes_{P(z_1)} (W_2\boxtimes_{P(z_2)} W_3)$ for complex
numbers $z_1$, $z_2$ satisfying $|z_1|>|z_2|>|z_1-z_2|>0$.  Note that
we are using two distinct nonzero complex numbers, and that certain
inequalities hold.  This situation corresponds to the fact that a
Riemann sphere with one negatively oriented puncture and three
positively oriented punctures can be seen in two different ways as the
``product'' of two Riemann spheres each with one negatively oriented
puncture and two positively oriented punctures; the detailed geometric
motivation is presented in \cite{H1}, \cite{tensorK} and
\cite{tensor4}.

To construct this natural isomorphism, we first consider compositions
of certain intertwining maps.  As we have mentioned, a
$P(z)$-intertwining map $I$ of type ${W_3 \choose {W_1 W_2}}$ maps
into $\overline{W}_3$ rather than $W_3$.  Thus the existence of
compositions of suitable intertwining maps always entails certain
convergence.  In particular, the existence of the composition
$w_{(1)}\boxtimes_{P(z_1)} (w_{(2)}\boxtimes_{P(z_2)} w_{(3)})$ when
$|z_1|>|z_2|>0$ and the existence of the composition
$(w_{(1)}\boxtimes_{P(z_1-z_2)} w_{(2)})\boxtimes_{P(z_2)} w_{(3)}$
when $|z_2|>|z_1-z_2|>0$, for general elements $w_{(i)}$ of $W_i$,
$i=1,2,3$, requires the proof of certain convergence conditions.
These conditions will be discussed in detail below.

Let us now assume these convergence conditions and let $z_1$, $z_2$
satisfy $|z_1|>|z_2|>|z_1-z_2|>0$. To construct the desired
associativity isomorphism {}from $(W_1 \boxtimes_{P(z_1-z_2)}
W_2)\boxtimes_{P(z_2)} W_3$ to $W_1\boxtimes_{P(z_1)}
(W_2\boxtimes_{P(z_2)} W_3)$, it is equivalent (by duality) to give a
suitable natural isomorphism {}from $W_1\hboxtr_{P(z_1)}
(W_2\boxtimes_{P(z_2)} W_3)$ to $(W_1 \boxtimes_{P(z_1-z_2)}
W_2)\hboxtr_{P(z_2)} W_3$. As we mentioned in the previous subsection,
instead of constructing this isomorphism directly, we shall embed both
of these spaces, separately, into the single space $(W_1\otimes W_2
\otimes W_3)^*$.

We will see that $(W_1\otimes W_2 \otimes W_3)^*$ carries a natural
$V$-action analogous to the contragredient of the diagonal action in
the Lie algebra case (recall the similar action of $V$ on $(W_1\otimes
W_2)^*$ mentioned above).  Also, for four $V$-modules $W_1$, $W_2$,
$W_3$ and $W_4$, we have a canonical notion of ``$P(z_1,
z_2)$-intertwining map {}from $W_1 \otimes W_2 \otimes W_3$ to
$\overline{W}_4$'' given by a vertex-algebraic analogue of
(\ref{11i}); for this notion, we need only that $z_1$ and $z_2$ are
nonzero and distinct. The relation between these two concepts comes
{}from the natural linear injection
\begin{eqnarray}
\hom (W_1\otimes W_2\otimes W_3, \overline{W}_4) &\longrightarrow &
\hom (W'_4, (W_1\otimes W_2\otimes W_3)^*)\nno\\ F&\mapsto & F',
\end{eqnarray}
where $F': W'_4\longrightarrow (W_1\otimes W_2\otimes W_3)^*$ is given
by
\begin{equation}\label{F'}
\nu\mapsto \nu\circ F,
\end{equation}
which is indeed well defined.  Under this natural map, the
$P(x_1,z_2)$-intertwining maps correspond precisely to the maps {}from
$W'_4$ to $(W_1\otimes W_2\otimes W_3)^*$ that intertwine the two
natural $V$-actions on $W'_4$ and $(W_1\otimes W_2\otimes W_3)^*$.

Now for modules $W_1$, $W_2$, $W_3$, $W_4$, $M_1$, and a
$P(z_1)$-intertwining map $I_1$ and a $P(z_2)$-intertwining map $I_2$
of types ${W_4 \choose {W_1 M_1}}$ and ${M_1 \choose {W_2 W_3}}$,
respectively, where $|z_1|>|z_2|>0$, it turns out that by definition
the composition $\bar I_1\circ (1_{W_1}\otimes I_2)$ is a $P(z_1,
z_2)$-intertwining map.  Similarly, for a $P(z_2)$-intertwining map
$I^1$ and a $P(z_1-z_2)$-intertwining map $I^2$ of types ${W_4 \choose
{M_2 W_3}}$ and ${M_2 \choose {W_1 W_2}}$, respectively, where
$|z_2|>|z_1-z_2|>0$, the composition $\bar I^1\circ (I^2\otimes
1_{W_3})$ is a $P(z_1, z_2)$-intertwining map. Hence we have two
maps intertwining the $V$-actions:
\begin{eqnarray}
W'_4 &\longrightarrow &(W_1\otimes W_2\otimes
W_3)^*\nno\\
\nu&\mapsto &\nu\circ F_1,
\label{iiv1}
\end{eqnarray}
where $F_1$ is the intertwining map $\bar I_1\circ (1_{W_1}\otimes
I_2)$, and
\begin{eqnarray}
W'_4 &\longrightarrow & (W_1\otimes W_2\otimes
W_3)^*\nno\\
\nu&\mapsto &\nu\circ F_2,
\label{iiv2}
\end{eqnarray}
where $F_2$ is the intertwining map $\bar I^1\circ (I^2\circ 1_{W_3})$.

It is important to note that we can express these compositions $\bar
I_1\circ (1_{W_1}\otimes I_2)$ and $\bar I^1\circ (I^2\otimes 1_{W_3})$
in terms of intertwining operators, as discussed in Remark
\ref{im-io}.  Let ${\cal Y}_1$, ${\cal Y}_2$, ${\cal Y}^1$ and ${\cal
Y}^2$ be the intertwining operators corresponding to $I_1$, $I_2$,
$I^1$ and $I^2$, respectively.  Then the compositions $\bar I_1\circ
(1_{W_1}\otimes I_2)$ and $\bar I^1\circ (I^2\otimes 1_{W_3})$
correspond to the ``product'' ${\cal Y}_1(\cdot, x_1){\cal Y}_2(\cdot,
x_2)\cdot$ and ``iterate'' ${\cal Y}^1({\cal Y}^2(\cdot, x_0)\cdot,
x_2)\cdot$ of intertwining operators, respectively, and we make the
``substitutions'' (in the sense of Remark \ref{im-io}) $x_1 \mapsto
z_1$, $x_2 \mapsto z_2$ and $x_0 \mapsto z_1-z_2$ in order to express
the two compositions of intertwining maps as the ``product'' ${\cal
Y}_1(\cdot, z_1){\cal Y}_2(\cdot, z_2)\cdot$ and ``iterate'' ${\cal
Y}^1({\cal Y}^2(\cdot, z_1-z_2)\cdot, z_2)\cdot$ of intertwining maps,
respectively.  (These products and iterates involve a branch of the
$\log$ function and also certain convergence.)

Just as in the Lie algebra case, the special cases in which the
modules $W_4$ are two iterated tensor product modules and the
``intermediate'' modules $M_1$ and $M_2$ are two tensor product
modules are particularly interesting: When $W_4=W_1\boxtimes_{P(z_1)}
(W_2\boxtimes_{P(z_2)} W_3)$ and $M_1=W_2\boxtimes_{P(z_2)} W_3$, and
$I_1$ and $I_2$ are the corresponding canonical intertwining maps,
(\ref{iiv1}) gives the natural $V$-homomorphism
\begin{eqnarray}
W_1\hboxtr_{P(z_1)}(W_2 \boxtimes_{P(z_2)} W_3)&\longrightarrow &
(W_1\otimes W_2\otimes W_3)^*\nno\\
\nu&\mapsto &(w_{(1)}\otimes w_{(2)}\otimes w_{(3)}\mapsto \nno\\
&&\nu(w_{(1)}\boxtimes_{P(z_1)} (w_{(2)}\boxtimes_{P(z_2)}
w_{(3)})))\nno;\\
\label{injva1}
\end{eqnarray}
when $W_4=(W_1\boxtimes_{P(z_1-z_2)} W_2)\boxtimes_{P(z_2)} W_3$ and
$M_2=W_1\boxtimes_{P(z_1-z_2)}W_2$, and $I^1$ and $I^2$ are the
corresponding canonical intertwining maps, (\ref{iiv2}) gives the
natural $V$-homomorphism
\begin{eqnarray}
(W_1 \boxtimes_{P(z_1-z_2)} W_2)\hboxtr_{P(z_2)} W_3&\longrightarrow &
(W_1\otimes W_2\otimes
W_3)^*\nno\\
\nu&\mapsto &(w_{(1)}\otimes w_{(2)}\otimes w_{(3)}\mapsto \nno\\
&&\nu((w_{(1)}\boxtimes_{P(z_1-z_2)} w_{(2)})\boxtimes_{P(z_2)}
w_{(3)})).\nno\\
\label{injva2}
\end{eqnarray}

It turns out that both of these maps are injections \cite{tensor4}, as
we prove below, so that we are embedding the spaces
$W_1\hboxtr_{P(z_1)}(W_2 \boxtimes_{P(z_2)} W_3)$ and $(W_1
\boxtimes_{P(z_1-z_2)} W_2)\hboxtr_{P(z_2)} W_3$ into the space
$(W_1\otimes W_2 \otimes W_3)^*$. Following the ideas in
\cite{tensor4}, we shall give a precise description of the ranges of
these two maps, and under suitable conditions, prove that the two
ranges are the same; this will establish the associativity
isomorphism.

More precisely, as in \cite{tensor4}, we prove that for any $P(z_1,
z_2)$-intertwining map $F$, the image of any $\nu\in W'_4$ under $F'$
(recall (\ref{F'})) satisfies certain conditions that we call the
``$P(z_1, z_2)$-compatibility condition'' and the ``$P(z_1,
z_2)$-local grading restriction condition.''  Hence, as special cases,
the elements of $(W_1\otimes W_2 \otimes W_3)^*$ in the ranges of
either of the maps (\ref{iiv1}) or (\ref{iiv2}), and in particular, of
(\ref{injva1}) or (\ref{injva2}), satisfy these conditions.

In addition, any $\lambda\in (W_1\otimes W_2\otimes W_3)^*$ induces
two maps $\mu^{(1)}_\lambda$ and $\mu^{(2)}_\lambda$ as in (\ref{mu1})
and (\ref{mu2}).  We will see that any element $\lambda$ of the range
of (\ref{iiv1}), and in particular, of (\ref{injva1}), must satisfy
the condition that the elements $\mu^{(1)}_\lambda(w_{(1)})$ all lie
in a suitable completion of the subspace $W_2\hboxtr_{P(z_2)} W_3$ of
$(W_2\otimes W_3)^*$, and any element $\lambda$ of the range of
(\ref{iiv2}), and in particular, of (\ref{injva2}), must satisfy the
condition that the elements $\mu^{(2)}_\lambda(w_{(3)})$ all lie in a
suitable completion of the subspace $W_1\hboxtr_{P(z_1-z_2)} W_2$ of
$(W_1\otimes W_2)^*$.  These conditions will be called the
``$P^{(1)}(z)$-local grading restriction condition'' and the
``$P^{(2)}(z)$-local grading restriction condition,'' respectively.

It turns out that the construction of the desired natural
associativity isomorphism follows {}from showing that the ranges of both
of (\ref{injva1}) and (\ref{injva2}) satisfy both of these conditions.
This amounts to a certain ``extension condition'' on our module
category.  When this extension condition and a suitable convergence
condition are satisfied, as in \cite{tensor4} we show below that the
desired associativity isomorphisms do exist, and that in addition, the
associativity of intertwining maps holds.  That is, let $z_1$ and
$z_2$ be complex numbers satisfying the inequalities
$|z_1|>|z_2|>|z_1-z_2|>0$.  Then for any $P(z_1)$-intertwining map
$I_1$ and $P(z_2)$-intertwining map $I_2$ of types ${W_4\choose {W_1\,
M_1}}$ and ${M_1\choose {W_2\, W_3}}$, respectively, there is a
suitable module $M_2$, and a $P(z_2)$-intertwining map $I^1$ and a
$P(z_1-z_2)$-intertwining map $I^2$ of types ${W_4\choose{M_2\, W_3}}$
and ${M_2\choose{W_1\, W_2}}$, respectively, such that
\begin{equation}\label{iiii}
\langle w'_{(4)}, \bar I_1(w_{(1)}\otimes I_2(w_{(2)}\otimes w_{(3)}))
\rangle =
\langle w'_{(4)}, \bar I^1(I^2(w_{(1)}\otimes w_{(2)})\otimes w_{(3)})
\rangle
\end{equation}
for $w_{(1)}\in W_1, w_{(2)}\in W_2$, $w_{(3)} \in W_3$ and
$w'_{(4)}\in W'_4$; and conversely, given $I^1$ and $I^2$ as
indicated, there exist a suitable module $M_1$ and maps $I_1$ and
$I_2$ with the indicated properties.  In terms of intertwining
operators (recall the comments above), the equality (\ref{iiii}) reads
\begin{eqnarray}\label{yyyy}
\lefteqn{\langle w'_{(4)}, {\cal Y}_1(w_{(1)}, x_1){\cal
Y}_2(w_{(2)},x_2)
w_{(3)}\rangle|_{x_1=z_1,\; x_2=z_2}}\nno\\
&&=\langle w'_{(4)},{\cal Y}^1({\cal Y}^2(w_{(1)}, x_0)w_{(2)},
x_2)w_{(3)})\rangle |_{x_0=z_1-z_2,\; x_2=z_2},
\end{eqnarray}
where ${\cal Y}_1$, ${\cal Y}_2$, ${\cal Y}^1$ and ${\cal Y}^2$ are
the intertwining operators corresponding to $I_1$, $I_2$, $I^1$ and
$I^2$, respectively.  (As we have been mentioning, the substitution of
complex numbers for formal variables involves a branch of the $\log$
function and also certain convergence.)  In this sense, the
associativity asserts that the ``product'' of two suitable
intertwining maps can be written as the ``iterate'' of two suitable
intertwining maps, and conversely.

{}From this construction of the natural associativity isomorphisms we
will see, by analogy with (\ref{elemap}), that
$(w_{(1)}\boxtimes_{P(z_1-z_2)} w_{(2)})\boxtimes_{P(z_2)} w_{(3)}$ is
mapped naturally to $w_{(1)}\boxtimes_{P(z_1)}
(w_{(2)}\boxtimes_{P(z_2)} w_{(3)})$ under the natural extension of
the corresponding associativity isomorphism (these elements in general
lying in the algebraic completions of the corresponding tensor product
modules). In fact, this property
\begin{equation}
(w_{(1)}\boxtimes_{P(z_1-z_2)} w_{(2)})\boxtimes_{P(z_2)} w_{(3)}
\mapsto
w_{(1)}\boxtimes_{P(z_1)}
(w_{(2)}\boxtimes_{P(z_2)} w_{(3)})
\end{equation}
for $w_{(1)}\in W_{1}$, $w_{(2)}\in W_{2}$ and $w_{(3)}\in W_{3}$
characterizes
the associativity isomorphism
\begin{equation}
(W_{1}\boxtimes_{P(z_1-z_2)} W_{2})\boxtimes_{P(z_2)} W_{3}\to
W_{1}\boxtimes_{P(z_1)}
(W_{2}\boxtimes_{P(z_2)} W_{3})
\end{equation}
(cf. (\ref{elemap})).  The coherence property of the associativity
isomorphisms will follow {}from this fact.

\begin{rema}{\rm
Note that equation (\ref{yyyy}) can be written as
\begin{equation}\label{yyyy2}
{\cal Y}_1(w_{(1)}, z_1){\cal Y}_2(w_{(2)},z_2)=
{\cal Y}^1({\cal Y}^2(w_{(1)}, z_1-z_2)w_{(2)},
z_2),
\end{equation}
with the appearance of the complex numbers being understood as
substitutions in the sense mentioned above, and with the ``generic''
vectors $w_{(3)}$ and $w'_{(4)}$ being implicit.  This (rigorous)
equation amounts to the ``operator product expansion'' in the physics
literature on conformal field theory; indeed, in our language, if we
expand the right-hand side of (\ref{yyyy2}) in powers of $z_1-z_2$, we
find that a product of intertwining maps is expressed as an expansion
in powers of $z_1-z_2$, with coefficients that are again intertwining
maps, of the form ${\cal Y}^1(w,z_2)$.  When all three modules are the
vertex algebra itself, and all the intertwining operators are the
canonical vertex operator $Y(\cdot,x)$ itself, this ``operator product
expansion'' follows easily {}from the Jacobi identity.  But for
intertwining operators in general, it is hard to prove the operator
product expansion, that is, to prove the assertions involving
(\ref{iiii}) and (\ref{yyyy}) above.}
\end{rema}

\begin{rema}{\rm
The constructions of the tensor product modules and of the
associativity isomorphisms previewed above for suitably general vertex
algebras follow those in \cite{tensor1}, \cite{tensor2},
\cite{tensor3} and \cite{tensor4}.  Alternative constructions are
certainly possible.  For example, an alternative construction of the
tensor product modules was given in \cite{Li}.  However, no matter
what construction is used for the tensor product modules of suitably
general vertex algebras, one cannot avoid constructing structures and
proving results equivalent to what is carried out in this work.  The
constructions in this work of the tensor product functors and of the
natural associativity isomorphisms are crucial in the deeper part of
the theory of vertex tensor categories. }
\end{rema}

\begin{rema}{\rm

A braided tensor category structure on certain module categories for
affine Lie algebras, and more generally, on certain module categories
for ``chiral algebras'' associated with ``rational conformal field
theories,'' was discovered by Moore and Seiberg \cite{MS} in their
important study of conformal field theory.  However, they constructed
this structure based on the assumption of the existence of a suitable
tensor product functor (including a tensor product module) and also
the assumption of the existence of a suitable operator product
expansion for chiral vertex operators, which is essentially equivalent
to assuming the associativity of intertwining maps, as we have
expressed it above.  As we have discussed, the desired tensor product
modules were constructed under suitable conditions in the series of
papers \cite{tensor1}, \cite{tensor2} and \cite{tensor3}, and in
\cite{tensor4} the appropriate natural associativity isomorphisms
among tensor products of triples of modules were constructed and it
was shown that this is equivalent to the desired associativity of
intertwining maps (and thus the existence of a suitable operator
product expansion).  The results in these papers will now be
generalized in this work. }

\end{rema}

\paragraph{Acknowledgments}
We would like to thank Sasha Kirillov Jr., Masahiko Miyamoto, Kiyokazu
Nagatomo and Bin Zhang for their early interest in our work on
logarithmic tensor product theory, and Bin Zhang in particular for
inviting L.Z. to lecture on this ongoing work at Stony Brook.
Y.-Z.H. would also like to thank Paul Fendley for raising his interest
in logarithmic conformal field theory.  We are grateful to the
participants in our courses and seminars for their insightful comments
on this work, especially Geoffrey Buhl, William Cook, Thomas Robinson
and Brian Vancil.  The authors gratefully acknowledge partial support
{}from NSF grants DMS-0070800 and DMS-0401302.

\newpage

\setcounter{equation}{0}
\setcounter{rema}{0}

\section{The setting}

In this section we define and discuss the basic structures and
introduce some notation that will be used in this work. More
specifically, we first introduce the notions of ``conformal vertex
algebra'' and ``M\"obius vertex algebra.''  A conformal vertex algebra
is just a vertex algebra equipped with a conformal vector satisfying
the usual axioms; a M\"obius vertex algebra is a variant of a
``quasi-vertex operator algebra'' as in \cite{FHL}, with the
difference that the two grading restriction conditions in the
definition of vertex operator algebra are not required. We then define
the notion of module for each of these types of vertex algebra.
Relaxing the $L(0)$-semisimplicity in the definition of module we
obtain the notion of ``generalized module.''  Finally, we notice that
in order to have a contragredient functor on the module category under
consideration, we need to impose a stronger grading condition.  This
leads to the notions of ``strong gradedness'' of M\"obius vertex
algebras and their generalized modules. In this work we are mainly
interested in certain full subcategories of the category of strongly
graded generalized modules for certain strongly graded M\"obius vertex
algebras.  Throughout the work we shall assume some familiarity with
the material in \cite{FLM2}, \cite{FHL}, \cite{DL} and \cite{LL}.

Throughout, we shall use the notation ${\mathbb N}$ for the nonnegative
integers and $\Z_{+}$ for the positive integers.

We shall continue to use the notational convention concerning formal
variables and complex variables given in Remark
\ref{formalandcomplexvariables}.  Recall {}from \cite{FLM2},
\cite{FHL} or \cite{LL} that the ``formal delta function'' is defined
as the formal Laurent series
\[
\delta(x)=\sum_{n\in {\mathbb Z}} x^n.
\]
We will consistently use the {\em binomial expansion convention}: For
any complex number $\lambda$, $(x+y)^\lambda$ is to be expanded in
nonnegative integral powers of the second variable, i.e.,
\[
(x+y)^\lambda=\sum_{n\in {\mathbb N}} {\lambda \choose n} x^{\lambda -n}
y^n.
\]
Here $x$ or $y$ might be something other than a formal
variable (or a nonzero complex multiple of a formal variable); for
instance, $x$ or $y$ (but not both!) might be a nonzero complex
number, or $x$ or $y$ might be some more complicated object.  The use
of the binomial expansion convention will be clear in context.

Objects like $\delta(x)$ and $(x+y)^\lambda$ lie in spaces of formal
series.  Some of the spaces that we will use are, with $W$ a vector
space (over $\C$) and $x$ a formal variable:
\[
W[x]=\biggl\{\sum_{n\in \mathbb{N}}a_{n}x^{n}| a_{n}\in W,\;
\mbox{all but finitely many} \; a_{n} = 0\biggr\},
\]
\[
W[x,x^{-1}]=\biggl\{\sum_{n\in \mathbb{Z}}a_{n}x^{n}| a_{n}\in W,\;
\mbox{all but finitely many} \; a_{n} = 0\biggr\},
\]
\[
W[[x]]=\biggl\{\sum_{n\in \mathbb{N}}a_{n}x^{n}| a_{n}\in W\;
\mbox{(with possibly infinitely many} \; a_{n} \; \mbox{not} \; 0)\biggr\},
\]
\[
W((x))=\biggl\{\sum_{n\in \mathbb{Z}}a_{n}x^{n}| a_{n}\in W,\;
a_{n} = 0 \; \mbox{for sufficiently small} \; n \biggr\},
\]
\[
W[[x,x^{-1}]]=\biggl\{\sum_{n\in \mathbb{Z}}a_{n}x^{n}| a_{n}\in W\;
\mbox{(with possibly infinitely many} \; a_{n} \; \mbox{not} \; 0)\biggr\}.
\]
We will also need
\begin{equation}\label{formalserieswithcomplexpowers}
W\{ x\}=\biggl\{\sum_{n\in \mathbb{C}}a_{n}x^{n}| a_{n}\in W \;
\mbox{for} \; n\in {\mathbb C}\biggr\}
\end{equation}
as in \cite{FLM2}; here the powers of the formal variable are complex,
and the coefficients may all be nonzero.  We will also use analogues
of these spaces involving two or more formal variables.

The following formal version of Taylor's theorem is easily verified by
direct expansion (see Proposition 8.3.1 of \cite{FLM2}):  For $f(x)
\in W\{ x\}$,
\begin{equation}\label{formalTaylortheorem}
e^{y \frac{d}{dx}}f(x) = f(x + y),
\end{equation}
where the exponential denotes the formal exponential series, and where
we are using the binomial expansion convention on the right-hand side.
It is important to note that this formula holds for arbitrary formal
series $f(x)$ with complex powers of $x$, where $f(x)$ need not be an
expansion in any sense of an analytic function (again, see Proposition
8.3.1 of \cite{FLM2}).

The formal delta function $\delta(x)$ has the following simple and
fundamental property: For any $f(x)\in W[x, x^{-1}]$,
\begin{equation}
f(x)\delta(x)=f(1)\delta(x).
\end{equation}
(Here we are taking the liberty of writing complex numbers to the
right of vectors in $W$.)  This is proved immediately by observing its
truth for $f(x)=x^n$ and then using linearity.  This property has many
important variants; in general, whenever an expression is multiplied
by the formal delta function, we may formally set the argument
appearing in the delta function equal to 1, provided that the relevant
algebraic expressions make sense.  For example, for any
$$X(x_{1},
x_{2})\in (\mbox{End }W)[[x_{1}, x_{1}^{-1}, x_{2}, x_{2}^{-1}]]$$
such that
\begin{equation}\label{limx1approachesx2}
\lim_{x_{1}\to x_{2}}X(x_{1}, x_{2})=X(x_{1},
x_{2})\lbar_{x_{1}=x_{2}}
\end{equation}
exists, we have
\begin{equation}\label{Xx1x2=Xx2x2}
X(x_{1}, x_{2})\delta\left(\frac{x_{1}}{x_{2}}\right)=X(x_{2}, x_{2})
\delta\left(\frac{x_{1}}{x_{2}}\right).
\end{equation}
The existence of the ``algebraic limit'' defined in
(\ref{limx1approachesx2}) means that for an arbitrary vector $w\in W$,
the coefficient of each power of $x_{2}$ in the formal expansion
$X(x_{1}, x_{2})w\lbar_{x_{1}=x_{2}}$ is a finite sum.  In general,
the existence of such ``algebraic limits,'' and also such products of
formal sums, always means that the coefficient of each monomial in the
relevant formal variables gives a finite sum.  Often, proving the
existence of the relevant algebraic limits (or products) is a much
more subtle matter than computing such limits (or products), just as
in analysis.  (In this work, we will typically use ``substitution
notation'' like $\lbar_{x_{1}=x_{2}}$ or $X(x_{2},x_{2})$ rather than
the formal limit notation on the left-hand side of
(\ref{limx1approachesx2}).)  Below, we will give a more sophisticated
analogue of the delta-function substitution principle
(\ref{Xx1x2=Xx2x2}), an analogue that we will need in this work.

This analogue, and in fact, many fundamental principles of vertex
operator algebra theory, are based on certain delta-function
expressions of the following type, involving three (commuting and
independent, as usual) formal variables:
\[
x_{0}^{-1}\delta\left(\frac{x_{1}-x_{2}}{x_{0}}\right)=\sum_{n\in
{\mathbb Z}}
\frac{(x_{1}-x_{2})^{n}}{x_{0}^{n+1}}=\sum_{m\in {\mathbb N},\; n\in
{\mathbb Z}}
(-1)^{m}{{n}\choose {m}} x_{0}^{-n-1}x_{1}^{n-m}x_{2}^{m};
\]
here the binomial expansion convention is of course being used.

The following important identities involving such three-variable
delta-function expressions are easily proved (see \cite{FLM2} or
\cite{LL}, where extensive motivation for these formulas is also
given):
\begin{equation}\label{2termdeltarelation}
x_{2}^{-1}\left(\frac{x_{1}-x_{0}}{x_{2}}\right)=
x_{1}^{-1}\delta\left(\frac{x_{2}+x_{0}}{x_{1}}\right),
\end{equation}
\begin{equation}\label{3termdeltarelation}
x_{0}^{-1}\delta\left(\frac{x_{1}-x_{2}}{x_{0}}\right)-
x_{0}^{-1}\delta\left(\frac{x_{2}-x_{1}}{-x_{0}}\right)=
x_{2}^{-1}\delta\left(\frac{x_{1}-x_{0}}{x_{2}}\right).
\end{equation}
Note that the three terms in (\ref{3termdeltarelation}) involve
nonnegative integral powers of $x_2$, $x_1$ and $x_0$, respectively.
In particular, the two terms on the left-hand side of
(\ref{3termdeltarelation}) are unequal formal Laurent series in three
variables, even though they might appear equal at first glance.  We
shall use these two identities extensively.

\begin{rema}\label{deltafunctionsubstitutionremark}{\rm
Here is the useful analogue, mentioned above, of the delta-function
substitution principle (\ref{Xx1x2=Xx2x2}):  Let
\begin{equation}
f(x_1,x_2,y) \in (\mbox{End }W)[[x_{1}, x_{1}^{-1}, x_{2},
x_{2}^{-1},y,y^{-1}]]
\end{equation}
be such that
\begin{equation}
\lim_{x_{1}\to x_{2}}f(x_1,x_2,y) \; \mbox{exists}
\end{equation}
and such that for any $w \in W$,
\begin{equation}
f(x_1,x_2,y)w \in W[[x_1,x_{1}^{-1},x_2,x_{2}^{-1}]]((y)).
\end{equation}
Then
\begin{equation}\label{deltafunctionsubstitutionformula}
x_{1}^{-1}\delta\left(\frac{x_{2}-y}{x_{1}}\right)f(x_1,x_2,y)=
x_{1}^{-1}\delta\left(\frac{x_{2}-y}{x_{1}}\right)f(x_2-y,x_2,y).
\end{equation}
For this principle, see Remark 2.3.25 of \cite{LL}, where the proof is
also presented.}
\end{rema}

The following formal residue notation will be useful: For
\[
f(x) = \sum_{n\in \mathbb{C}}a_{n}x^{n} \in W\{x\}
\]
(note that the powers of $x$ need not be integral),
\[
{\res_x}f(x) = a_{-1}.
\]
For instance, for the expression in (\ref{2termdeltarelation}),
\begin{equation}
\res_{x_2}x_{2}^{-1}\left(\frac{x_{1}-x_{0}}{x_{2}}\right)=1.
\end{equation}

For a vector space $W$, we will denote its vector space dual by $W^*$
($=\hom_{\mathbb C}(W,{\mathbb C})$), and we will use the notation
$\langle\cdot,\cdot\rangle_W$, or $\langle\cdot,\cdot\rangle$ if the
underlying space $W$ is clear, for the canonical pairing between $W^*$
and $W$.

We will use the following version of the notion of ``conformal vertex
algebra'': A conformal vertex algebra is a vertex algebra (in the
sense of Borcherds \cite{B}; see \cite{LL}) equipped with a ${\mathbb
Z}$-grading and with a conformal vector satisfying the usual
compatibility conditions.  Specifically:

\begin{defi}\label{cva}
{\rm A {\it conformal vertex algebra} is a ${\mathbb Z}$-graded vector
space
\begin{equation}\label{Vgrading}
V=\coprod_{n\in{\mathbb Z}} V_{(n)}
\end{equation}
(for $v\in V_{(n)}$, we say the {\it weight} of $v$ is $n$ and we write
$\mbox{wt}\, v=n$)
equipped with a linear map $V\otimes V\to V[[x,
x^{-1}]]$, or equivalently,
\begin{eqnarray}\label{YforV}
V&\to&({\rm End}\; V)[[x, x^{-1}]] \nno\\
v&\mapsto& Y(v,
x)={\displaystyle \sum_{n\in{\mathbb Z}}}v_{n}x^{-n-1} \;\;(
\mbox{where }v_{n}\in{\rm End} \;V),
\end{eqnarray}
$Y(v, x)$ denoting the {\it vertex operator associated with} $v$, and
equipped also with two distinguished vectors ${\bf 1}\in V_{(0)}$ (the
{\it vacuum vector}) and $\omega\in V_{(2)}$ (the {\it conformal
vector}), satisfying the following conditions for $u,v \in V$:
the {\it lower truncation condition}:
\begin{equation}\label{ltc}
u_{n}v=0\;\;\mbox{ for }n\mbox{ sufficiently large}
\end{equation}
(or equivalently, $Y(u, x)v\in V((x))$); the {\it vacuum property}:
\begin{equation}\label{1left}
Y({\bf 1}, x)=1_V;
\end{equation}
the {\it creation property}:
\begin{equation}\label{1right}
Y(v, x){\bf 1} \in V[[x]]\;\;\mbox{ and }\;\lim_{x\rightarrow
0}Y(v, x){\bf 1}=v
\end{equation}
(that is, $Y(v, x){\bf 1}$ involves only nonnegative integral powers
of $x$ and the constant term is $v$); the {\it Jacobi identity} (the
main axiom):
\begin{eqnarray}
&x_0^{-1}\delta
\bigg({\displaystyle\frac{x_1-x_2}{x_0}}\bigg)Y(u, x_1)Y(v,
x_2)-x_0^{-1} \delta
\bigg({\displaystyle\frac{x_2-x_1}{-x_0}}\bigg)Y(v, x_2)Y(u,
x_1)&\nno \\ &=x_2^{-1} \delta
\bigg({\displaystyle\frac{x_1-x_0}{x_2}}\bigg)Y(Y(u, x_0)v,
x_2)&\label{Jacobi}
\end{eqnarray}
(note that when each expression in (\ref{Jacobi}) is applied to any
element of $V$, the coefficient of each monomial in the formal
variables is a finite sum; on the right-hand side, the notation
$Y(\cdot, x_2)$ is understood to be extended in the obvious way to
$V[[x_0, x^{-1}_0]]$); the {\em Virasoro algebra relations}:
\begin{equation}\label{vir1}
[L(m), L(n)]=(m-n)L(m+n)+{\displaystyle\frac{1}{12}}
(m^3-m)\delta_{n+m,0}c
\end{equation}
for $m, n \in {\mathbb Z}$, where
\begin{equation}\label{vir2}
L(n)=\omega _{n+1}\;\; \mbox{ for } \;  n\in{\mathbb Z}, \;\;
\mbox{i.e.},\;\;Y(\omega, x)=\sum_{n\in{\mathbb Z}}L(n)x^{-n-2},
\end{equation}
\begin{equation}\label{vir3}
c\in {\mathbb C}
\end{equation}
(the {\it central charge} or {\it rank} of $V$);
\begin{equation}\label{L-1derivativeproperty}
{\displaystyle \frac{d}{dx}}Y(v,
x)=Y(L(-1)v, x)
\end{equation}
(the {\it  $L(-1)$-derivative property}); and
\begin{equation}\label{L0gradingproperty}
L(0)v=nv=(\wt v)v \;\; \mbox{ for }\; n\in {\mathbb Z}\; \mbox{ and }\;
v\in V_{(n)}.
\end{equation}
}
\end{defi}

This completes the definition of the notion of conformal vertex
algebra.  We will denote such a conformal vertex algebra by $(V,Y,{\bf
1},\omega)$ or simply by $V$.

The only difference between the definition of conformal vertex algebra
and the definition of {\it vertex operator algebra} (in the sense of
\cite{FLM2} and \cite{FHL}) is that a vertex operator algebra $V$
also satisfies the two {\it grading restriction conditions}
\begin{equation}\label{gr1}
V_{(n)}=0 \;\; \mbox{ for }n\mbox{ sufficiently negative},
\end{equation}
and
\begin{equation}\label{gr2}
{\rm dim} \; V_{(n)} < \infty \;\;\mbox{ for }n \in {\mathbb Z}.
\end{equation}
(As we mentioned above, a {\it vertex algebra} is the same thing as a
conformal vertex algebra but without the assumptions of a grading or a
conformal vector, or, of course, the $L(n)$'s.)

\begin{rema}\label{va>cva}{\rm
Of course, not every vertex algebra is conformal. For example, it is
well known \cite{B} that any commutative associative algebra $A$ with
unit $1$, together with a derivation $D:A\to A$ can be equipped with a
vertex algebra structure, by:
\[
Y(\cdot,x)\cdot : A\times A\to A[[x]],\;\;Y(a,x)b=(e^{xD}a)b,
\]
and ${\bf 1}=1$.  In particular, $u_n=0$ for any $u\in A$ and $n\geq
0$.  If $\omega$ is a conformal vector for such a vertex algebra, then
for any $u\in A$, $Du=u_{-2}{\bf 1}=L(-1)u$ {}from (\ref{1right}) and
(\ref{L-1derivativeproperty}), so $D=L(-1)=\omega_0$, which equals $0$
because $\omega=L(0)\omega/2=\omega_1\omega/2=0$.  Thus a vertex
algebra constructed {}from a commutative associative algebra with
nonzero derivation in this way cannot be conformal.  }
\end{rema}

\begin{rema}\label{motivate-Mobius}{\rm
The theory of vertex tensor categories inherently uses the whole
moduli space of spheres with two positively oriented punctures and one
negatively oriented puncture (and in fact, more generally, with
arbitrary numbers of positively oriented punctures and one negatively
oriented puncture) equipped with general (analytic) local coordinates
vanishing at the punctures.  Because of the analytic local
coordinates, our constructions require certain conditions on the
Virasoro algebra operators.  However, recalling the definition of the moduli
space elements $P(z)$ {}from Subsection 1.4, we point out that if we
restrict our attention to elements of the moduli space of only the
type $P(z)$, then the relevant operations of sewing and subsequently
decomposing Riemann spheres continue to yield spheres of the same
type, and rather than general conformal transformations around the
punctures, only M\"obius (projective) transformations around the
punctures are needed.  This makes it possible to develop the essential
structure of our tensor product theory by working entirely with
spheres of this special type; the general vertex tensor category
theory then follows {}from the structure thus developed.  This is why,
in the present work, we are focusing on the theory of $P(z)$-tensor
products.  Correspondingly, it turns out that it is very natural for
us to consider, along with the notion of conformal vertex algebra
(Definition \ref{cva}), a weaker notion of vertex algebra
involving only the three-dimensional subalgebra of the Virasoro
algebra corresponding to the group of M\"obius transformations.  That
is, instead of requiring an action of the whole Virasoro algebra, we
use only the action of the Lie algebra ${\mathfrak s}{\mathfrak l}(2)$
generated by $L(-1)$, $L(0)$ and $L(1)$.  Thus we get a notion
essentially identical to the notion of ``quasi-vertex operator
algebra'' in \cite{FHL}; the reason for focusing on this notion here
is the same as the reason why it was considered in \cite{FHL}.  Here
we designate this notion by the term ``M\"obius vertex algebra''; the
only difference between the definition of M\"obius vertex algebra and
the definition of quasi-vertex operator algebra \cite{FHL} is that a
quasi-vertex operator algebra $V$ also satisfies the two grading
restriction conditions (\ref{gr1}) and (\ref{gr2}).}
\end{rema}

Thus we formulate:

\begin{defi}\label{mobdef}{\rm
The notion of {\it M\"obius vertex algebra} is defined in the same way
as that of conformal vertex algebra except that in addition to the
data and axioms concerning $V$, $Y$ and ${\bf 1}$ (through
(\ref{Jacobi}) in Definition \ref{cva}), we assume (in place of
the existence of the conformal vector $\omega$ and the Virasoro
algebra conditions (\ref{vir1}), (\ref{vir2}) and (\ref{vir3})) the
following: We have a representation $\rho$ of ${\mathfrak s}{\mathfrak l}(2)$
on $V$ given by
\begin{equation}\label{Lrho}
L(j)=\rho(L_j),\;\;j=0,\pm 1,
\end{equation}
where $\{L_{-1}, L_0, L_1\}$ is a basis of ${\mathfrak s}{\mathfrak l}(2)$
with Lie brackets
\begin{equation}\label{L_*}
[L_0, L_{-1}]=L_{-1},\;\;[L_0,L_1]=-L_1,{\rm\;\;and\;\;}
[L_{-1},L_1]=-2L_0,
\end{equation}
and the following conditions hold for $v \in V$:
\begin{equation}\label{sl2-1}
[L(-1), Y(v,x)]=Y(L(-1)v,x),
\end{equation}
\begin{equation}\label{sl2-2}
[L(0), Y(v,x)]=Y(L(0)v,x)+xY(L(-1)v,x),
\end{equation}
\begin{equation}\label{sl2-3}
[L(1),Y(v,x)]=Y(L(1)v,x)+2xY(L(0)v,x)+x^2Y(L(-1)v,x),
\end{equation}
and also, (\ref{L-1derivativeproperty}) and
(\ref{L0gradingproperty}).  Of course, (\ref{sl2-1})--(\ref{sl2-3})
can be written as
\begin{eqnarray}\label{sl2-all}
[L(j),Y(v,x)]&=&\sum_{k=0}^{j+1}{j+1\choose k}x^kY(L(j-k)v,x)\nno\\
&=&\sum_{k=0}^{j+1}{j+1\choose k}x^{j+1-k}Y(L(k-1)v,x)
\end{eqnarray}
for $j=0,\pm 1$.}
\end{defi}

We will denote such a M\"obius vertex algebra by $(V,Y,{\bf 1}, \rho)$
or simply by $V$. Note that there is no notion of central charge (or
rank) for a M\"obius vertex algebra.  Also, a conformal vertex algebra
can certainly be viewed as a M\"obius vertex algebra in the obvious
way.  (Of course, a conformal vertex algebra could have other ${\mathfrak
s}{\mathfrak l}(2)$-structures making it a M\"obius vertex algebra in a
different way.)

\begin{rema}{\rm
By (\ref{Lrho}) and (\ref{L_*}) we have $[L(0), L(j)]=-jL(j)$ for
$j=0,\pm 1$. Hence
\begin{equation}\label{degL(j)}
L(j)V_{(n)}\subset V_{(n-j)},\;\; \mbox{ for }\;j=0,\pm1.
\end{equation}
Moreover, {}from (\ref{sl2-1}), (\ref{sl2-2}) and (\ref{sl2-3}) with
$v={\bf 1}$ we get, by (\ref{1left}) and (\ref{1right}),
\[
L(j){\bf 1}=0\;\; \mbox{ for } \; j=0,\pm 1.
\]
}
\end{rema}

\begin{rema}{\rm
Not every M\"obius vertex algebra is conformal. As an example, take
the commutative associative algebra ${\mathbb C}[t]$ with derivation
$D=-d/dt$, and form a vertex algebra as in Remark \ref{va>cva}. By
Remark \ref{va>cva}, this vertex algebra is not conformal. However,
define linear operators
\[
L(-1)=D,\quad L(0)=tD,\quad L(1)=t^2D
\]
on ${\mathbb C}[t]$. Then it is straightforward to verify that ${\mathbb
C}[t]$ becomes a M\"obius vertex algebra with these operators giving a
representation of ${\mathfrak s}{\mathfrak l}(2)$ having the desired
properties and with the ${\mathbb Z}$-grading (by nonpositive integers)
given by the eigenspace decomposition with respect to $L(0)$.}
\end{rema}

\begin{rema}{\rm
It is also easy to see that not every vertex algebra is M\"obius. For
example, take the two-dimensional commutative associative algebra
$A={\mathbb C}1\oplus{\mathbb C}a$ with $1$ as identity and $a^2=0$. The
linear operator $D$ defined by $D(1)=0$, $D(a)=a$ is a nonzero
derivation of $A$. Hence $A$ has a vertex algebra structure by Remark
\ref{va>cva}. Now if it is a module for ${\mathfrak s}{\mathfrak l}(2)$ as in
Definition \ref{mobdef}, since $A$ is two-dimensional and $L(0)1=0$,
$L(0)$ must act as $0$. But then $D=L(-1)=[L(0),L(-1)]=0$, a
contradiction.  }
\end{rema}

A module for a conformal vertex algebra $V$ is a module for $V$ viewed
as a vertex algebra such that the conformal element acts in the same
way as in the definition of vertex operator algebra. More precisely:

\begin{defi}\label{cvamodule}
{\rm  Given a conformal vertex algebra
$(V,Y,{\bf 1},\omega)$,  a {\it module} for $V$ is a ${\mathbb C}$-graded
vector space
\begin{equation}\label{Wgrading}
W=\coprod_{n\in{\mathbb C}} W_{(n)}
\end{equation}
(graded by {\it weights}) equipped with a linear map
$V\otimes W \rightarrow W[[x,x^{-1}]]$, or equivalently,
\begin{eqnarray}\label{YforW}
V &\rightarrow & (\mbox{End}\ W)[[x,x^{-1}]] \nno\\
v&\mapsto & Y(v,x) =\sum_{n\in {\mathbb Z}}v_nx^{-n-1}\;\;\;
(\mbox{where}\;\; v_n \in \mbox{End}\ W)
\end{eqnarray}
(note that the sum is over ${\mathbb Z}$, not ${\mathbb C}$), $Y(v,x)$
denoting the {\it vertex operator on $W$ associated with $v$}, such
that all the defining properties of a conformal vertex algebra that
make sense hold.  That is, the following conditions are satisfied: the
lower truncation condition: for $v \in V$ and $w \in W$,
\begin{equation}\label{ltc-w}
v_nw = 0 \;\;\mbox{ for }\;n\;\mbox{ sufficiently large}
\end{equation}
(or equivalently, $Y(v, x)w\in W((x))$); the vacuum property:
\begin{equation}\label{m-1left}
Y(\mbox{\bf 1},x) = 1_W;
\end{equation}
the Jacobi identity for vertex operators on $W$: for $u, v \in V$,
\begin{eqnarray}\label{m-Jacobi}
&{\dps x^{-1}_0\delta \bigg( {x_1-x_2\over x_0}\bigg)
Y(u,x_1)Y(v,x_2) - x^{-1}_0\delta \bigg( {x_2-x_1\over -x_0}\bigg)
Y(v,x_2)Y(u,x_1)}&\nno\\
&{\dps = x^{-1}_2\delta \bigg( {x_1-x_0\over x_2}\bigg)
Y(Y(u,x_0)v,x_2)}
\end{eqnarray}
(note that on the right-hand side, $Y(u,x_0)$ is the operator on $V$
associated with $u$); the Virasoro algebra relations on $W$ with
scalar $c$ equal to the central charge of $V$:
\begin{equation}\label{m-vir1}
[L(m), L(n)]=(m-n)L(m+n)+{\displaystyle\frac{1}{12}}
(m^3-m)\delta_{n+m,0}c
\end{equation}
for $m,n \in {\mathbb Z}$, where
\begin{equation}\label{m-vir2}
L(n)=\omega _{n+1}\;\; \mbox{ for }n\in{\mathbb Z}, \;\;{\rm
i.e.},\;\;Y(\omega, x)=\sum_{n\in{\mathbb Z}}L(n)x^{-n-2};
\end{equation}
\begin{equation}\label{L-1}
\displaystyle \frac{d}{dx}Y(v, x)=Y(L(-1)v, x)
\end{equation}
(the $L(-1)$-derivative property); and
\begin{equation}\label{wl0}
(L(0)-n)w=0\;\;\mbox{ for }\;n\in {\mathbb C}\;\mbox{ and }\;w\in W_{(n)}.
\end{equation}
}
\end{defi}

This completes the definition of the notion of module for a conformal
vertex algebra.

\begin{rema}\label{virrelationsformodule}{\rm
The Virasoro algebra relations (\ref{m-vir1}) for a module action follow
{}from the corresponding relations (\ref{vir1}) for $V$ together with the
Jacobi identities (\ref{Jacobi}) and (\ref{m-Jacobi}) and the
$L(-1)$-derivative properties (\ref{L-1derivativeproperty}) and
(\ref{L-1}), as we recall {}from (for example) \cite{FHL} or \cite{LL}.}
\end{rema}

We also have:

\begin{defi}\label{moduleMobius}{\rm
The notion of {\it module} for a M\"obius vertex algebra is defined in
the same way as that of module for a conformal vertex algebra except
that in addition to the data and axioms concerning $W$ and $Y$
(through (\ref{m-Jacobi}) in Definition \ref{cvamodule}), we assume
(in place of the Virasoro algebra conditions (\ref{m-vir1}) and
(\ref{m-vir2})) a representation $\rho$ of ${\mathfrak s}{\mathfrak l}(2)$ on
$W$ given by (\ref{Lrho}) and the conditions (\ref{sl2-1}),
(\ref{sl2-2}) and (\ref{sl2-3}), for operators acting on $W$, and
also, (\ref{L-1}) and (\ref{wl0}).  }
\end{defi}

In addition to modules, we have the following notion of {\em
generalized module} (or {\em logarithmic module}, as in, for example,
\cite{Mi}):

\begin{defi}\label{definitionofgeneralizedmodule}{\rm
A {\it generalized module} for a conformal (respectively, M\"obius)
vertex algebra is defined in the same way as a module for a conformal
(respectively, M\"obius) vertex algebra except that in the grading
(\ref{Wgrading}), each space $W_{(n)}$ is replaced by $W_{[n]}$, where
$W_{[n]}$ is the generalized $L(0)$-eigenspace corresponding to the
(generalized) eigenvalue $n\in {\mathbb C}$; that is, (\ref{Wgrading})
and (\ref{wl0}) in the definition are replaced by
\begin{equation}\label{Wgeneralizedgrading}
W=\coprod_{n\in{\mathbb C}} W_{[n]}
\end{equation}
and
\begin{equation}\label{gerwt}
\mbox{for }\;n\in {\mathbb C}\;\mbox{ and }\;w\in W_{[n]},\;
(L(0)-n)^mw=0 \;\mbox{ for }\;m\in {\mathbb N}\;\mbox{ sufficiently
large},
\end{equation}
respectively.  For $w \in W_{[n]}$, we still write $\wt w = n$ for the
(generalized) weight of $w$.
}
\end{defi}

We will denote such a module or generalized module just defined by
$(W,Y)$, or sometimes by $(W,Y_W)$ or simply by $W$. We will use the
notation
\begin{equation}\label{pi_n}
\pi_n: W\to W_{[n]}
\end{equation}
for the projection {}from $W$ to its subspace of (generalized) weight
$n$, and for its natural extensions to spaces of formal series with
coefficients in $W$. In either the conformal or M\"obius case, a
module is of course a generalized module.

\begin{rema}\label{generalizedeigenspacedecomp}
{\rm For any vector space $U$ on which an operator, say, $L(0)$, acts
in such a way that
\begin{equation}\label{U=directsum}
U=\coprod_{n\in{\mathbb C}} U_{[n]}
\end{equation}
where for $n\in {\mathbb C}$,
\[
U_{[n]} = \{ u \in U | (L(0)-n)^m u=0 \;\mbox{ for }\;m\in {\mathbb
N}\;\mbox{ sufficiently large} \},
\]
we shall typically use the same projection notation
\begin{equation}
\pi_n: U\to U_{[n]}
\end{equation}
as in (\ref{pi_n}).  If instead of (\ref{U=directsum}) we have only
\[
U=\sum_{n\in{\mathbb C}} U_{[n]},
\]
then in fact this sum is indeed direct, and for any $L(0)$-stable
subspace $T$ of $U$, we have
\[
T=\coprod_{n\in{\mathbb C}} T_{[n]}
\]
(as with ordinary rather than generalized eigenspaces).}
\end{rema}

\begin{rema}\label{modulesaremodules}{\rm
A module for a conformal vertex algebra $V$ is obviously again a
module for $V$ viewed as a M\"obius vertex algebra, and conversely, a
module for $V$ viewed as a M\"obius vertex algebra is a module for $V$
viewed as a conformal vertex algebra, by Remark
\ref{virrelationsformodule}.  Similarly, the generalized modules for a
conformal vertex algebra $V$ are exactly the generalized modules for
$V$ viewed as a M\"obius vertex algebra.}
\end{rema}

\begin{rema}{\rm
A conformal or M\"obius vertex algebra is a module for itself (and 
in particular, a generalized module for itself).}
\end{rema}

\begin{rema}{\rm
In either the conformal or M\"obius vertex algebra case, we have the
obvious notions of $V$-{\em module homomorphism}, {\em submodule},
{\em quotient module}, and so on; in particular, homomorphisms are
understood to be grading-preserving.  We sometimes write the vector
space of (generalized-) module maps (homomorphisms) $W_1 \to W_2$ for
(generalized) $V$-modules $W_1$ and $W_2$ as $\hom_{V}(W_1,W_2)$.  }
\end{rema}

\begin{rema}{\rm
We have chosen the name ``generalized module'' here because the vector
space underlying the module is graded by generalized
eigenvalues. (This notion is different {}from the notion of
``generalized module'' used in \cite{tensor1}. A generalized module
for a vertex operator algebra $V$ as defined in, for example,
Definition 2.11 of \cite{tensor1} is precisely a module for $V$ viewed
as a conformal vertex algebra.)  }
\end{rema}

We will use the following notion of (formal algebraic) completion of a
generalized module:

\begin{defi}\label{Wbardef}{\rm
Let $W=\coprod_{n\in{\mathbb C}}W_{[n]}$ be a generalized module for a
M\"obius (or conformal) vertex algebra. We denote by $\overline{W}$
the (formal) completion of $W$ with respect to the ${\mathbb
C}$-grading, that is,
\begin{equation}\label{Wbar}
\overline{W}=\prod _{n\in {\mathbb C}} W_{[n]}.
\end{equation}
We will use the same notation $\overline{U}$ for any ${\mathbb
C}$-graded subspace $U$ of $W$.  We will continue to use the notation
$\pi_n$ for the projection {}from $\overline{W}$ to $W_{[n]}$.  We will
also continue to use the notation $\langle\cdot,\cdot\rangle_W$, or
$\langle\cdot,\cdot\rangle$ if the underlying space is clear, for the
canonical pairing between the subspace $\coprod _{n\in {\mathbb C}}
(W_{[n]})^*$ of $W^*$, and $\overline{W}$.  We are of course viewing
$(W_{[n]})^*$ as embedded in $W^*$ in the natural way, that is, for
$w^*\in (W_{[n]})^*$,
\begin{equation}\label{Wnstar}
\langle w^*, w\rangle_W=\langle w^*, w_n\rangle_{W_{[n]}}
\end{equation}
for any $w=\sum_{m\in {\mathbb C}} w_m$ (finite sum) in $W$, where
$w_m\in W_{[m]}$.}
\end{defi}

The following weight formula holds for generalized modules,
generalizing the corresponding formula in the module case
(cf.\ \cite{Mi}):

\begin{propo}\label{gweight}
Let $W$ be a generalized module for a M\"obius (or conformal) vertex
algebra $V$. Let both $v\in V$ and $w\in W$ be homogeneous. Then
\begin{eqnarray}
&\wt(v_n w)=\wt v +\wt w-n-1 \;\; \mbox{ for any }\; n\in {\mathbb Z},
&\label{set:wtvn}\\
&\wt(L(j)w)=\wt w-j \;\;\mbox{ for }\; j=0,\pm 1.&\label{set:wtsl2}
\end{eqnarray}
\end{propo}
\pf Applying the $L(-1)$-derivative property (\ref{L-1}) to formula
(\ref{sl2-2}), with the operators acting on $W$, and extracting the
coefficient of $x^{-n-1}$, we obtain:
\begin{equation}\label{[L(0),v_n]}
[L(0), v_n]=(L(0)\,v)_n+(-n-1)v_n.
\end{equation}
This can be written as
\[
(L(0)-(\wt v-n-1))v_n=v_n L(0),
\]
and so we have
\[
(L(0)-(\wt v+m-n-1))v_n=v_n (L(0)-m)
\]
for any $m\in {\mathbb C}$. Applying this repeatedly we get
\[
(L(0)-(\wt v+m-n-1))^t v_n=v_n (L(0)-m)^t
\]
for any $t\in {\mathbb N}$, $m\in {\mathbb C}$, and (\ref{set:wtvn})
follows.

For (\ref{set:wtsl2}), since as operators acting on $W$ we have
\begin{equation}\label{set:0j}
[L(0),L(j)]=-jL(j)
\end{equation}
for $j=0,\pm 1$, we get $(L(0)+j)L(j)=L(j)L(0)$ so that
\[
(L(0)-m+j)L(j)=L(j)(L(0)-m)
\]
for any $m\in {\mathbb C}$. Thus
\[
(L(0)-m+j)^tL(j)=L(j)(L(0)-m)^t
\]
for any $t\in {\mathbb N}$, $m\in {\mathbb C}$, and (\ref{set:wtsl2})
follows. \epf

\begin{rema}\label{congruent}{\rm
{}From Proposition \ref{gweight} we see that a generalized $V$-module
$W$ decomposes into submodules corresponding to the congruence classes
of its weights modulo $\Z$: For $\mu \in \C/\Z$, let
\begin{equation}
W_{[\mu]} = \coprod_{\bar n = \mu} W_{[n]},
\end{equation}
where $\bar n$ denotes the equivalence class of $n \in \C$ in
$\C/\Z$.  Then
\begin{equation}
W = \coprod_{\mu \in \C/\Z} W_{[\mu]}
\end{equation}
and each $W_{[\mu]}$ is a $V$-submodule of $W$.  Thus if a generalized
module $W$ is indecomposable (in particular, if it is irreducible),
then all complex numbers $n$ for which $W_{[n]}\neq 0$ are congruent
modulo ${\mathbb Z}$ to each other. }
\end{rema}

\begin{rema}\label{set:L(0)s}{\rm
Let $W$ be a generalized module for a M\"obius (or conformal) vertex
algebra $V$.  We consider the ``semisimple part'' $L(0)_s\in
\mbox{End}\ W$ of the operator $L(0)$:
\[
L(0)_sw=nw \;\; \mbox{ for } \; w\in W_{[n]},\; n \in {\mathbb C}.
\]
Then on $W$ we have
\begin{eqnarray}
&{}[L(0)_s, v_n]=[L(0), v_n]\;\;\mbox{ for all }\;v\in V\;\mbox{ and }
\; n\in {\mathbb Z};&\label{L0s,vn}\\
&{}[L(0)_s, L(j)]=[L(0), L(j)]\;\;\mbox{ for }\; j=0,\pm 1.
&\label{L0s,Lj}
\end{eqnarray}
Indeed, for homogeneous elements $v\in V$ and $w\in W$,
(\ref{set:wtvn}) and (\ref{[L(0),v_n]}) imply that
\begin{eqnarray*}
[L(0)_s, v_n]w&=&L(0)_s(v_nw)-v_n(L(0)_sw)\\
&=&(\wt v +\wt w-n-1)v_nw-(\wt w)v_nw\\
&=&(\wt v)v_nw+(-n-1)v_nw\\
&=&(L(0)v)_nw+(-n-1)v_nw\\
&=&[L(0), v_n]w.
\end{eqnarray*}
Similarly, for any homogeneous element $w\in W$ and $j=0,\pm 1$,
(\ref{set:wtsl2}) and (\ref{set:0j}) imply that
\begin{eqnarray*}
[L(0)_s, L(j)]w&=&L(0)_s(L(j)w)-L(j)(L(0)_sw)\\
&=&(\wt w-j)L(j)w-(\wt w)L(j)w\\
&=&-jL(j)w\\
&=&[L(0), L(j)]w.
\end{eqnarray*}
Thus the ``locally nilpotent part'' $L(0)-L(0)_s$ of $L(0)$ commutes
with the action of $V$ and of ${\mathfrak s}{\mathfrak l}(2)$ on $W$.  In
other words, $L(0)-L(0)_s$ is a $V$-homomorphism {}from $W$ to itself.}
\end{rema}

Now suppose that $L(1)$ acts locally nilpotently on a M\"obius (or
conformal) vertex algebra $V$, that is, for
any $v\in V$, there is $m\in {\mathbb N}$ such that $L(1)^mv=0$.  Then
generalizing formula (3.20) in \cite{tensor1} (the case of ordinary
modules for a vertex operator algebra), we define the {\it opposite
vertex operator} on a generalized $V$-module $(W,Y_W)$ associated to
$v\in V$ by
\begin{equation}\label{yo}
Y^o_W(v,x)=Y_W(e^{xL(1)}(-x^{-2})^{L(0)}v,x^{-1}),
\end{equation}
that is, for $k\in {\mathbb Z}$ and $v\in V_{(k)}$,
\begin{eqnarray}\label{yo1}
Y^o_W(v,x)&=&\sum_{n \in {\mathbb Z}} v^o_n x^{-n-1}\nno\\
&=&\sum_{n\in {\mathbb Z}}\bigg((-1)^k\sum_{m\in {\mathbb N}}
\frac{1}{m!}(L(1)^mv)_{-n-m-2+2k}\bigg)x^{-n-1},
\end{eqnarray}
as in \cite{tensor1}. (In the present work, we are replacing the 
symbol $*$ used in \cite{tensor1} for opposite vertex operators
by the symbol $o$; see also Subsection 5.1 below.)
Here we are defining the component operators
\begin{equation}\label{v^o_n}
v^o_n=(-1)^k\sum_{m\in {\mathbb N}}\frac{1}{m!}(L(1)^mv)_{-n-m-2+2k}
\end{equation}
for $v\in V_{(k)}$ and $n,k\in {\mathbb Z}$. Note that the $L(1)$-local
nilpotence ensures well-definedness here.  Clearly, $v \mapsto
Y^o_W(v,x)$ is a linear map $V \to (\mbox{End} \ W)[[x,x^{-1}]]$ such
that $V \otimes W \to W((x^{-1}))$ ($v \otimes w \mapsto
Y^o_W(v,x)w)$.

By (\ref{v^o_n}), (\ref{degL(j)}) and (\ref{set:wtvn}), we see that
for $n,k\in {\mathbb Z}$ and $v\in V_{(k)}$, the operator $v^o_n$ is of
generalized weight $n+1-k\;(=n+1-\wt v)$, in the sense that
\begin{equation}\label{v^o-deg}
v^o_nW_{[m]}\subset W_{[m+n+1-k]}\;\;\mbox{ for any }\;m\;\in{\mathbb C}.
\end{equation}

As mentioned in \cite{tensor1} (see (3.23) in \cite{tensor1}), the
proof of Jacobi identity in Theorem 5.2.1 of \cite{FHL} proves the
following {\it opposite Jacobi identity} for $Y^o_{W}$ in the case
where $V$ is a vertex operator algebra and $W$ is a $V$-module:
\begin{eqnarray}\label{op-jac-id}
\lefteqn{\dps x_{0}^{-1}\delta\bigg(\frac{x_{1}-x_{2}}{x_{0}}\bigg)
Y_{W}^o(v, x_{2})Y^o_{W}(u, x_{1})}\nno\\
&&\hspace{2em}-x_{0}^{-1}\delta\bigg(\frac{x_{2}-x_{1}}{-x_{0}}\bigg)
Y_{W}^o(u, x_{1})Y^o_{W}(v, x_{2})\nno\\
&&{\dps =x_{2}^{-1}\delta\bigg(\frac{x_{1}-x_{0}}{x_{2}}\bigg)
Y_{W}^o(Y(u, x_{0})v, x_{2})}
\end{eqnarray}
for $u,v \in V$, and taking ${\rm Res}_{x_{0}}$ gives us the {\it
opposite commutator formula}.  Similarly, the proof of the
$L(-1)$-derivative property in Theorem 5.2.1 of \cite{FHL} proves the
following {\it $L(-1)$-derivative property} for $Y^o_{W}$ in the same
case:
\begin{equation}\label{yo-l-1}
\frac{d}{dx} Y^o_{W}(v,x) = Y^o_{W}(L(-1)v,x).
\end{equation}
The same proofs carry over and prove the opposite Jacobi identity and
the $L(-1)$-derivative property
for $Y^o_{W}$  in the present case, where $V$ is a M\"obius (or
conformal) vertex algebra with $L(1)$ acting locally nilpotently and
where $W$ is a generalized $V$-module.  In the case in which $V$ is a
conformal vertex algebra, we have
\begin{equation}\label{Yoppositeomega}
Y^o_{W}(\omega,x) = Y_{W}(x^{-4}\omega,x^{-1})
=\sum_{n\in {\mathbb Z}}L(n)x^{n-2}
\end{equation}
since $L(1)\omega = 0$.

For opposite vertex operators, we have the following analogues of
(\ref{sl2-1})--(\ref{sl2-all}) in the M\"obius case:

\begin{lemma}\label{sl2opposite}
For $v \in V$,
\begin{equation}\label{sl2opp-1}
[Y^o_W(v,x),L(1)]=Y^o_W(L(-1)v,x),
\end{equation}
\begin{equation}\label{sl2opp-2}
[Y^o_W(v,x),L(0)]=Y^o_W(L(0)v,x)+xY^o_W(L(-1)v,x),
\end{equation}
\begin{equation}\label{sl2opp-3}
[Y^o_W(v,x),L(-1)]=Y^o_W(L(1)v,x)+2xY^o_W(L(0)v,x)+x^2Y^o_W(L(-1)v,x).
\end{equation}
Equivalently,
\begin{eqnarray}\label{sl2opp-all}
[Y^o_W(v,x),L(-j)]&=&\sum_{k=0}^{j+1}{j+1\choose k}x^k
Y^o_W(L(j-k)v,x)\nno\\
&=&\sum_{k=0}^{j+1}{j+1\choose k}x^{j+1-k}
Y^o_W(L(k-1)v,x)
\end{eqnarray}
for $j=0,\pm 1$.
\end{lemma}
\pf
For $j=0, \pm 1$, by definition and (\ref{sl2-all}) we have
\begin{eqnarray}\label{sl2opp-all-1}
[Y^o_W(v,x),L(j)]&=&-[L(j), Y_W(e^{xL(1)}(-x^{-2})^{L(0)}v,x^{-1})]\nn
&=&-\sum_{k=0}^{j+1}{j+1\choose k}x^{-k}Y_{W}(L(j-k)
e^{xL(1)}(-x^{-2})^{L(0)}v,x^{-1}).
\end{eqnarray}
By (5.2.14) in \cite{FHL} and the fact that
\begin{equation}\label{xL(0)L(j)}
x^{L(0)}L(j)x^{-L(0)}=x^{-j}L(j)
\end{equation}
(easily proved by applying to a homogeneous vector),
\begin{eqnarray}\label{sl2opp-all-2}
\lefteqn{L(-1)e^{xL(1)}(-x^{-2})^{L(0)}}\nn
&&=
e^{xL(1)}L(-1)(-x^{-2})^{L(0)}-2xe^{xL(1)}L(0)(-x^{-2})^{L(0)}
+x^{2}e^{xL(1)}L(1)(-x^{-2})^{L(0)}\nn
&&=-x^{2}e^{xL(1)}(-x^{-2})^{L(0)}L(-1)-2xe^{xL(1)}(-x^{-2})^{L(0)}L(0)
-e^{xL(1)}(-x^{-2})^{L(0)}L(1)\nn
&&=-e^{xL(1)}(-x^{-2})^{L(0)}(x^{2}L(-1)+2xL(0)+L(1)).
\end{eqnarray}
We also have
\begin{eqnarray}\label{sl2opp-all-4}
L(1)e^{xL(1)}(-x^{-2})^{L(0)}&=&
e^{xL(1)}L(1)(-x^{-2})^{L(0)}\nn
&=&-x^{-2}e^{xL(1)}(-x^{-2})^{L(0)}L(1).
\end{eqnarray}
By (\ref{sl2opp-all-2}), (\ref{sl2opp-all-4}),
$L(0)=\frac{1}{2}[L(1), L(-1)]$ and $[L(1), L(0)]=L(1)$, we have
\begin{eqnarray}\label{sl2opp-all-3}
\lefteqn{L(0)e^{xL(1)}(-x^{-2})^{L(0)}}\nn
&&=\frac{1}{2}L(1)L(-1)e^{xL(1)}(-x^{-2})^{L(0)}
-\frac{1}{2}L(-1)L(1)e^{xL(1)}(-x^{-2})^{L(0)}\nn
&&=-\frac{1}{2}L(1)e^{xL(1)}(-x^{-2})^{L(0)}(x^{2}L(-1)+2xL(0)+L(1))\nn
&&\quad +\frac{1}{2}x^{-2}L(-1)e^{xL(1)}(-x^{-2})^{L(0)}L(1)\nn
&&=\frac{1}{2}x^{-2}e^{xL(1)}(-x^{-2})^{L(0)}L(1)(x^{2}L(-1)+2xL(0)+L(1))\nn
&&\quad
-\frac{1}{2}x^{-2}e^{xL(1)}(-x^{-2})^{L(0)}(x^{2}L(-1)+2xL(0)+L(1))L(1)\nn
&&=e^{xL(1)}(-x^{-2})^{L(0)}L(0)+x^{-1}e^{xL(1)}(-x^{-2})^{L(0)}L(1)\nn
&&=e^{xL(1)}(-x^{-2})^{L(0)}(L(0)+x^{-1}L(1)).
\end{eqnarray}
Thus we obtain
\begin{eqnarray*}
\lefteqn{[Y^o_W(v,x),L(1)]}\nn
&&=-\sum_{k=0}^{2}{2\choose k}x^{-k}Y_{W}(L(1-k)
e^{xL(1)}(-x^{-2})^{L(0)}v,x^{-1})\nn
&&=-Y_{W}(L(1)e^{xL(1)}(-x^{-2})^{L(0)}v,x^{-1})
-2x^{-1}Y_{W}(L(0)
e^{xL(1)}(-x^{-2})^{L(0)}v,x^{-1})\nn
&&\quad -x^{-2}Y_{W}(L(-1)
e^{xL(1)}(-x^{-2})^{L(0)}v,x^{-1})\nn
&&=x^{-2}Y_{W}(e^{xL(1)}(-x^{-2})^{L(0)}L(1)v,x^{-1})
\nn
&&\quad -2x^{-1}Y_{W}(
e^{xL(1)}(-x^{-2})^{L(0)}(L(0)+x^{-1}L(1))v,x^{-1})
\nn
&&\quad +x^{-2}Y_{W}(
e^{xL(1)}(-x^{-2})^{L(0)}(x^{2}L(-1)+2xL(0)+L(1))v,x^{-1})\nn
&&=Y_{W}(
e^{xL(1)}(-x^{-2})^{L(0)}L(-1)v,x^{-1})\nn
&&=Y_{W}^{o}(L(-1)v,x),
\end{eqnarray*}
\begin{eqnarray*}
\lefteqn{[Y^o_W(v,x),L(0)]}\nn
&&=-\sum_{k=0}^{1}{1\choose k}x^{-k}Y_{W}(L(-k)
e^{xL(1)}(-x^{-2})^{L(0)}v,x^{-1})\nn
&&=-Y_{W}(L(0)
e^{xL(1)}(-x^{-2})^{L(0)}v,x^{-1})
-x^{-1}Y_{W}(L(-1)
e^{xL(1)}(-x^{-2})^{L(0)}v,x^{-1})\nn
&&=-Y_{W}(e^{xL(1)}(-x^{-2})^{L(0)}(L(0)+x^{-1}L(1))v,x^{-1})\nn
&&\quad
+x^{-1}Y_{W}(
e^{xL(1)}(-x^{-2})^{L(0)}(x^{2}L(-1)+2xL(0)+L(1))v,x^{-1})\nn
&&=Y_{W}(
e^{xL(1)}(-x^{-2})^{L(0)}(xL(-1)+L(0))v,x^{-1})\nn
&&=Y_{W}^{o}(L(0)v,x)+xY_{W}^{o}(L(-1)v,x)
\end{eqnarray*}
and
\begin{eqnarray*}
[Y^o_W(v,x),L(-1)]&=&-Y_{W}(L(-1)
e^{xL(1)}(-x^{-2})^{L(0)}v,x^{-1})\nn
&=&Y_{W}(
e^{xL(1)}(-x^{-2})^{L(0)}(x^{2}L(-1)+2xL(0)+L(1))v,x^{-1})\nn
&=&Y_{W}^{o}(L(1)v,x)+2xY_{W}^{o}(L(0)v,x)+x^{2}Y_{W}^{o}(L(-1)v,x),
\end{eqnarray*}
proving the lemma.
\epfv

As in Section 5.2 of \cite{FHL}, we can define a $V$-action on $W^*$
as follows:
\begin{equation}\label{y'}
\langle Y'(v,x)w',w\rangle = \langle w', Y^o_W(v,x)w\rangle
\end{equation}
for $v\in V$, $w'\in W^*$ and $w\in W$; the correspondence $v\mapsto
Y'(v,x)$ is a linear map {}from $V$ to $({\rm End}\,W^*)[[x,x^{-1}]]$.
Writing
\[
Y'(v,x)=\sum_{n\in {\mathbb Z}} v_n x^{-n-1}
\]
($v_n\in {\rm End}\,W^*)$, we have
\begin{equation}\label{v'vo}
\langle v_n w', w\rangle = \langle w', v^o_n w\rangle
\end{equation}
for $v\in V$, $w'\in W^*$ and $w\in W$.  (Actually, in \cite{FHL} this
$V$-action was defined on a space smaller than $W^*$, but this
definition holds without change on all of $W^*$.)  In the case in
which $V$ is a conformal vertex algebra we define the operators
$L'(n)\;(n\in {\mathbb Z})$ by
\[
Y'(\omega,x)=\sum_{n\in {\mathbb Z}}L'(n)x^{-n-2};
\]
then, by extracting the coefficient of $x^{-n-2}$ in (\ref{y'}) with
$v=\omega$ and using the fact that $L(1)\omega=0$ we have
\begin{equation}\label{L'(n)}
\langle L'(n)w',w\rangle=\langle w',L(-n)w\rangle\;\;\mbox{ for }\;n\in
{\mathbb Z}
\end{equation}
(see (\ref{Yoppositeomega})), as in Section 5.2 of \cite{FHL}.
In the case where $V$ is only a M\"obius vertex algebra, we define operators
$L'(-1)$, $L'(0)$ and $L'(1)$ on $W^*$ by formula (\ref{L'(n)}) for $n=0,
\pm 1$. It follows {}from (\ref{set:wtsl2}) that
\begin{equation}\label{L'(n)2}
L'(j)(W_{[m]})^*\subset (W_{[m-j]})^*
\end{equation}
for $m\in {\mathbb C}$ and $j=0,\pm 1$.  By combining (\ref{v'vo}) with
(\ref{v^o-deg}) we get
\begin{equation}\label{stable0}
v_n(W_{[m]})^*\subset (W_{[m+k-n-1]})^*
\end{equation}
for any $n,k\in {\mathbb Z}$, $v\in V_{(k)}$ and $m\in {\mathbb C}$.

We have just seen that the $L(1)$-local nilpotence
condition enables us to define a natural vertex operator action on the
vector space dual of a generalized module for a M\"obius (or
conformal) vertex algebra. This condition is satisfied by all
vertex operator algebras, due to (\ref{degL(j)}) and the grading
restriction condition (\ref{gr1}). However, the functor $W\mapsto W^*$
is certainly
not involutive, and $W^*$ is not in general a generalized module. In this work
we will need certain module categories equipped with an involutive
``contragredient functor'' $W\mapsto W'$ which generalizes the
contragredient functor for the category of modules for vertex operator
algebras. For this purpose, we introduce the following:

\begin{defi}\label{def:dgv}
{\rm Let $A$ be an abelian group.  A M\"obius (or conformal) vertex
algebra
\[
V=\coprod_{n\in {\mathbb Z}} V_{(n)}
\]
is said to be {\em strongly graded with respect to $A$} (or {\em
strongly $A$-graded}, or just {\em strongly graded} if the abelian
group $A$ is understood) if there is a second gradation on $V$, by $A$,
\[
V=\coprod _{\alpha \in A} V^{(\alpha)},
\]
such that the following conditions are satisfied: the two gradations
are compatible, that is,
\[
V^{(\alpha)}=\coprod_{n\in {\mathbb Z}} V^{(\alpha)}_{(n)} \;\;
(\mbox{where}\;V^{(\alpha)}_{(n)}=V_{(n)}\bigcap V^{(\alpha)})\;
\mbox{ for any }\;\alpha \in A;
\]
for any $\alpha,\beta\in A$ and $n\in {\mathbb Z}$,
\begin{eqnarray}
&V^{(\alpha)}_{(n)}=0\;\;\mbox{ for }\;n\;\mbox{ sufficiently
negative};&\label{dua:ltc}\\
&\dim V^{(\alpha)}_{(n)} <\infty;&\label{dua:fin}\\
&{\bf 1}\in V^{(0)}_{(0)};&\\
&v_l V^{(\beta)} \subset V^{(\alpha+\beta)}\;\;
\mbox{ for any }\;v\in V^{(\alpha)},\;l\in {\mathbb Z};&\label{v_l-A}
\end{eqnarray}
and
\begin{equation}\label{L(n)-A}
L(j)V^{(\alpha)} \subset V^{(\alpha)}\;\;\mbox{ for }\;j=0,\pm 1.
\end{equation}
If $V$ is in fact a conformal vertex algebra, we in addition require
that
\begin{equation}\label{omega0}
\omega\in V^{(0)}_{(2)},
\end{equation}
so that for all $j\in {\mathbb Z}$, (\ref{L(n)-A}) follows {}from
(\ref{v_l-A}).  }
\end{defi}

\begin{rema}\label{rm1}{\rm
Note that the notion of conformal vertex algebra strongly
graded with respect to the trivial group is exactly the notion of
vertex operator algebra. Also note that (\ref{degL(j)}),
(\ref{dua:ltc}) and (\ref{L(n)-A}) imply the local nilpotence of
$L(1)$ acting on $V$, and hence we have the construction and
properties of opposite
vertex operators on a generalized module for a strongly graded
M\"obius (or conformal) vertex algebra. }
\end{rema}

For (generalized) modules for a strongly graded algebra we will also
have a second grading by an abelian group, and it is natural to allow
this group to be larger than the second grading group $A$ for the
algebra.  (Note that this already occurs for the {\em first} grading
group, which is ${\mathbb Z}$ for algebras and ${\mathbb C}$ for
(generalized) modules.)

\begin{defi}\label{def:dgw}{\rm
Let $A$ be an abelian group and $V$ a strongly $A$-graded M\"obius (or
conformal)
vertex algebra. Let $\tilde A$ be an abelian group containing $A$ as a
subgroup. A $V$-module (respectively, generalized $V$-module)
\[
W=\coprod_{n\in{\mathbb C}} W_{(n)} \;\;\;(\mbox{respectively, }\;
W=\coprod_{n\in{\mathbb C}} W_{[n]})
\]
is said to be {\em strongly graded with respect to $\tilde A$} (or
{\em strongly $\tilde A$-graded}, or just {\em strongly graded}) if the
abelian group $\tilde A$ is understood) if there is a second gradation
on $W$, by $\tilde A$,
\begin{equation}\label{2ndgrd}
W=\coprod _{\beta \in \tilde A} W^{(\beta)},
\end{equation}
such that the following conditions are satisfied: the two gradations
are compatible, that is, for any $\beta \in \tilde A$,
\[
W^{(\beta)}=\coprod_{n\in {\mathbb C}} W^{(\beta)}_{(n)} \;\;(\mbox{where }\;
W^{(\beta)}_{(n)}=W_{(n)}\bigcap W^{(\beta)})
\]
\[
(\mbox{respectively, }\;
W^{(\beta)}=\coprod_{n\in {\mathbb C}} W^{(\beta)}_{[n]} \;\;(\mbox{where }\;
W^{(\beta)}_{[n]}=W_{[n]}\bigcap W^{(\beta)});
\]
for any $\alpha\in A$, $\beta\in \tilde A$ and $n\in {\mathbb C}$,
\begin{eqnarray}
&W^{(\beta)}_{(n+k)}=0 \;\; (\mbox{respectively, } \; W^{(\beta)}_{[n+k]}=0)
\;\;\mbox{ for }\;k\in {\mathbb Z}\;\mbox{ sufficiently
negative};&\label{set:dmltc}\\
&\dim W^{(\beta)}_{(n)} <\infty \;\; (\mbox{respectively, } \;
\dim W^{(\beta)}_{[n]} <\infty);&\label{set:dmfin}\\
&v_l W^{(\beta)} \subset W^{(\alpha+\beta)}\;\;\mbox{ for any }\;v\in
V^{(\alpha)},\;l\in {\mathbb Z};&\label{m-v_l-A}
\end{eqnarray}
and
\begin{equation}\label{m-L(n)-A}
L(j)W^{(\beta)} \subset W^{(\beta)}\;\;\mbox{ for }\;j=0,\pm 1.
\end{equation}
(Note that if $V$ is in fact a conformal vertex algebra, then for all
$j\in {\mathbb Z}$, (\ref{m-L(n)-A}) follows {}from (\ref{omega0}) and
(\ref{m-v_l-A}).)  }
\end{defi}

\begin{rema}\label{v-str-module}{\rm
A strongly  $A$-graded conformal or M\"{o}bius vertex algebra is a 
strongly  $A$-graded module for itself (and in particular, a 
strongly  $A$-graded generalized module for itself).}
\end{rema}

\begin{rema}\label{moduleswiththetrivialgroup}{\rm
Let $V$ be a vertex operator algebra, viewed (equivalently) as a
conformal vertex algebra strongly graded with respect to the trivial
group (recall Remark \ref{rm1}).  Then the $V$-modules that are
strongly graded with respect to the trivial group (in the sense of
Definition \ref{def:dgw}) are exactly the ($\C$-graded) modules for
$V$ as a vertex operator algebra, with the grading restrictions as
follows: For $n \in \C$,
\begin{equation}\label{Wn+k=0}
W_{(n+k)}=0 \;\; \mbox { for }\;k\in {\mathbb Z}\;\mbox{ sufficiently
negative}
\end{equation}
and
\begin{equation}\label{dimWnfinite}
\dim W_{(n)} <\infty.
\end{equation}
Also, the generalized $V$-modules that are strongly graded with
respect to the trivial group are exactly the generalized $V$-modules
(in the sense of Definition \ref{definitionofgeneralizedmodule}) such
that for $n \in \C$,
\begin{equation}
W_{[n+k]}=0 \;\; \mbox { for }\;k\in {\mathbb Z}\;\mbox{ sufficiently
negative}
\end{equation}
and
\begin{equation}
\dim W_{[n]} <\infty.
\end{equation}}
\end{rema}

\begin{rema}\label{homsaregradingpreserving}{\rm
In the strongly graded case, algebra and module homomorphisms are of
course understood to preserve the grading by $A$ or $\tilde A$.}
\end{rema}

\begin{exam}{\rm
An important source of examples of strongly graded conformal vertex
algebras and modules comes {}from the vertex algebras and modules
associated with even lattices.  Let $L$ be an even lattice, i.e., a
finite-rank free abelian group equipped with a nondegenerate symmetric
bilinear form $\langle\cdot,\cdot\rangle$, not necessarily positive
definite, such that $\langle a,a\rangle\in 2{\mathbb Z}$ for all $a\in
L$.  Then there is a natural structure of conformal vertex algebra on
a certain vector space $V_L$; see \cite{B} and Chapter 8 of
\cite{FLM2}.  If the form $\langle\cdot,\cdot\rangle$ on $L$ is also
positive definite, then $V_L$ is a vertex operator algebra (that is,
the grading restrictions hold).  If $L$ is not necessarily positive
definite, then $V_L$ is equipped with a natural second grading given
by $L$ itself, making $V_L$ a strongly $L$-graded conformal vertex
algebra in the sense of Definition \ref{def:dgv}.  Any (rational)
sublattice $M$ of the ``dual lattice'' $L^{\circ}$ of $L$ containing
$L$ gives rise to a strongly $M$-graded module for the strongly
$L$-graded conformal vertex algebra (see Chapter 8 of \cite{FLM2};
cf.\ \cite{LL}).}
\end{exam}

\begin{rema}{\rm
As mentioned in Remark \ref{rm1}, strong gradedness for a M\"obius (or
conformal) vertex algebra $V$ implies the
local nilpotence of $L(1)$ acting on $V$. In fact,
strong gradedness implies much more that will be important for us:
{}From (\ref{dua:ltc}),
(\ref{dua:fin}), (\ref{v_l-A}) and (\ref{L(n)-A}) (and (\ref{omega0})
in the conformal vertex algebra case), it is clear that strong gradedness
for $V$ implies the following {\em local grading restriction
condition} on $V$ (see \cite{H-codes}):
\begin{enumerate}
\item[(i)] for any $m>0$ and $v_{(1)}, \dots, v_{(m)}\in V$, there
exists $r\in {\mathbb Z}$ such that the coefficient of each monomial in
$x_1, \dots, x_{m-1}$ in the formal series $Y(v_{(1)}, x_1)\cdots
Y(v_{(m-1)}, x_{m-1})v_{(m)}$ lies in $\coprod_{n>r}V_{(n)}$;
\item[(ii)] in the conformal vertex algebra case: for any element of
the conformal vertex algebra $V$ homogeneous with respect to the
weight grading, the Virasoro-algebra submodule $M=\coprod_{n\in {\mathbb
Z}}M_{(n)}$ (where $M_{(n)}=M\bigcap V_{(n)}$) of $V$ generated by this
element satisfies the following grading restriction conditions:
$M_{(n)}=0$ when $n$ is sufficiently negative and $\dim
M_{(n)}<\infty$ for $n\in {\mathbb Z}$
\end{enumerate}
or
\begin{enumerate}
\item[(ii$'$)] in the M\"obius vertex algebra case: for any element of
the M\"obius vertex algebra $V$ homogeneous with respect to the weight
grading, the ${\mathfrak s}{\mathfrak l}(2)$-submodule $M=\coprod_{n\in {\mathbb
Z}}M_{(n)}$ (where $M_{(n)}=M\bigcap V_{(n)}$) of $V$ generated by this
element satisfies the following grading restriction conditions:
$M_{(n)}=0$ when $n$ is sufficiently negative and $\dim
M_{(n)}<\infty$ for $n\in {\mathbb Z}$.
\end{enumerate}
As was pointed out in \cite{H-codes}, Condition (i) above was first stated
in \cite{DL} (see formula (9.39), Proposition 9.17 and Theorem 12.33
in \cite{DL}) for generalized vertex algebras and abelian intertwining
algebras (certain generalizations of vertex algebras); it guarantees
the convergence, rationality and commutativity properties of the
matrix coefficients of products of more than two vertex operators.
Conditions (i) and (ii) (or (ii$'$)) together ensure that all the
essential results involving the Virasoro operators and the geometry of
vertex operator algebras in \cite{H1} still hold for these algebras.}
\end{rema}

\begin{rema}{\rm
Similarly, {}from (\ref{set:dmltc}), (\ref{set:dmfin}), (\ref{m-v_l-A})
and (\ref{m-L(n)-A}) (and (\ref{omega0}) in the conformal vertex
algebra case), it is clear that strong gradedness for (generalized)
modules implies the following {\it local grading restriction
condition} for a (generalized) module $W$ for a strongly graded
M\"obius (or conformal) vertex algebra $V$:
\begin{enumerate}
\item[(i)] for any $m>0$, $v_{(1)}, \dots, v_{(m-1)}\in V$, $n\in{\mathbb
C}$ and $w\in W_{[n]}$, there exists $r\in {\mathbb Z}$ such that the
coefficient of each monomial in $x_1, \dots, x_{m-1}$ in the formal
series $Y(v_{(1)}, x_1)\cdots Y(v_{(m-1)}, x_{m-1})w$ lies in
$\coprod_{k>r}W_{[n+k]}$;
\item[(ii)] in the conformal vertex algebra case: for any $w\in
W_{[n]}$ ($n\in{\mathbb C}$), the Virasoro-algebra submodule
$M=\coprod_{k\in {\mathbb Z}}M_{[n+k]}$ (where $M_{[n+k]}=M\bigcap
W_{[n+k]}$) of $W$ generated by $w$ satisfies the following grading
restriction conditions: $M_{[n+k]}=0$ when $k$ is sufficiently
negative and $\dim M_{[n+k]}<\infty$ for $k\in {\mathbb Z}$
\end{enumerate}
or
\begin{enumerate}
\item[(ii$'$)] in the M\"obius vertex algebra case: for any $w\in
W_{[n]}$ ($n\in {\mathbb C}$), the ${\mathfrak s}{\mathfrak l}(2)$-submodule
$M=\coprod_{k\in {\mathbb Z}}M_{[n+k]}$ (where $M_{[n+k]}=M\bigcap
W_{[n+k]}$) of $W$ generated by $w$ satisfies the following grading
restriction conditions: $M_{[n+k]}=0$ when $k$ is sufficiently
negative and $\dim M_{[n+k]}<\infty$ for $k\in {\mathbb Z}$.
\end{enumerate}
Note that in the case of ordinary (as opposed to generalized) modules,
all the generalized weight spaces such as $W_{[n]}$ mentioned here are
ordinary weight spaces $W_{(n)}$.}
\end{rema}

With the strong gradedness condition on a (generalized) module, we can
now define the corresponding notion of contragredient module. First we
give:

\begin{defi}\label{defofWprime}{\rm
Let $W=\coprod_{\beta \in \tilde A,\;n\in{\mathbb C}}W^{(\beta)}_{[n]}$
be a strongly $\tilde A$-graded generalized module for a strongly
$A$-graded M\"obius (or conformal) vertex algebra.  For each $\beta\in
\tilde A$ and $n \in {\mathbb C}$, let us identify
$(W^{(\beta)}_{[n]})^*$ with the subspace of $W^*$ consisting of the
linear functionals on $W$ vanishing on each $W^{(\gamma)}_{[m]}$ with
$\gamma \neq \beta$ or $m \neq n$ (cf.\ (\ref{Wnstar})).  We define
$W'$ to be the $(\tilde A\times{\mathbb C})$-graded vector subspace of
$W^*$ given by
\begin{equation}
W'=\coprod_{\beta \in \tilde A,\; n\in{\mathbb C}}
(W')^{(\beta)}_{[n]}, \;\;\mbox{ where }\;
(W')^{(\beta)}_{[n]}=(W^{(-\beta)}_{[n]})^*;
\end{equation}
we also use the notations
\begin{equation}\label{W'beta}
(W')^{(\beta)}=\coprod_{n\in{\mathbb C}}(W^{(-\beta)}_{[n]})^*
\subset(W^{(-\beta)})^* \subset W^*
\end{equation}
(where $(W^{(\beta)})^*$ consists of the linear functionals on $W$
vanishing on all $W^{(\gamma)}$ with $\gamma \neq \beta$) and
\begin{equation}
(W')_{[n]}=\coprod_{\beta\in\tilde A}(W^{(-\beta)}_{[n]})^*
\subset(W_{[n]})^* \subset W^*
\end{equation}
for the homogeneous subspaces of $W'$ with respect to the $\tilde A$-
and ${\mathbb C}$-grading, respectively (The reason for the minus
signs here will become clear below.)  We will still use the notation
$\langle\cdot,\cdot\rangle_W$, or $\langle\cdot,\cdot\rangle$ when the
underlying space is clear, for the canonical pairing between $W'$ and
$\overline{W}\subset \prod_{\beta \in \tilde A,\; n\in{\mathbb
C}}W^{(\beta)}_{[n]}$ (recall (\ref{Wbar})).  }
\end{defi}

\begin{rema}{\rm
In the case of ordinary rather than generalized modules, Definition
\ref{defofWprime} still applies, and all of the generalized weight
subspaces $W_{[n]}$ of $W$ are ordinary weight spaces $W_{(n)}$.  In
this case, we can write $(W')_{(n)}$ rather than $(W')_{[n]}$ for the
corresponding subspace of $W'$.}
\end{rema}

Let $W$ be a strongly graded (generalized) module for a strongly
graded M\"obius (or conformal) vertex algebra $V$. Recall that we have
the action (\ref{y'}) of $V$ on $W^*$ and that (\ref{stable0}) holds.
Furthermore, (\ref{v^o_n}), (\ref{v'vo}) and (\ref{m-v_l-A}) imply
for any $n,k\in {\mathbb Z}$, $\alpha\in A$, $\beta\in \tilde A$, $v\in
V^{(\alpha)}_{(k)}$ and $m\in {\mathbb C}$,
\begin{equation}\label{shift}
v_n((W')^{(\beta)}_{[m]})=v_n((W^{(-\beta)}_{[m]})^*)\subset
(W^{(-\alpha-\beta)}_{[m+k-n-1]})^* =(W')^{(\alpha+\beta)}_{[m+k-n-1]}.
\end{equation}
Thus $v_n$ preserves $W'$ for $v \in V$, $n \in {\mathbb Z}$.  Similarly
(in the M\"obius case), (\ref{L'(n)}), (\ref{L'(n)2}) and
(\ref{m-L(n)-A}) imply that $W'$ is stable under the operators
$L'(-1)$, $L'(0)$ and $L'(1)$, and in fact
\[
L'(j)(W')^{(\beta)}_{[n]}\subset (W')^{(\beta)}_{[n-j]}
\]
for any $j=0,\pm 1$, $\beta\in \tilde A$ and $n\in {\mathbb C}$.  In
the case or ordinary rather than generalized modules, the symbols
$(W')^{(\beta)}_{[n]}$, etc., can be replaced by
$(W')^{(\beta)}_{(n)}$, etc.

For any fixed $\beta\in \tilde A$ and $n\in {\mathbb C}$, by
(\ref{gerwt}) and the finite-dimensionality (\ref{set:dmfin}) of
$W^{(-\beta)}_{[n]}$, there exists $N\in {\mathbb N}$
such that $(L(0)-n)^NW^{(-\beta)}_{[n]}=0$. But then for any $w'\in
(W')^{(\beta)}_{[n]}$,
\begin{equation}\label{L(0)N}
\langle (L'(0)-n)^Nw',w\rangle=\langle w',(L(0)-n)^Nw\rangle=0
\end{equation}
for all $w\in W$. Thus $(L'(0)-n)^Nw'=0$. So (\ref{gerwt}) holds with
$W$ replaced by $W'$.  In the case of ordinary modules, we of course
take $N=1$.

By (\ref{set:dmltc}) and (\ref{shift}) we have the lower truncation
condition for the action $Y'$ of $V$ on $W'$:
\begin{equation}\label{truncationforY'}
\mbox{For any }\;v\in V\;\mbox{ and }\;w'\in W',\; v_nw'=0\;\;
\mbox{ for }\;n\;\mbox{ sufficiently large}.
\end{equation}
As a consequence, the Jacobi identity can now be formulated on
$W'$. In fact, by the above, and using the same proofs as those of
Theorems 5.2.1 and 5.3.1 in \cite{FHL}, together with Lemma
\ref{sl2opposite}, we obtain:

\begin{theo}\label{set:W'}
Let $\tilde A$ be an abelian group containing $A$ as a subgroup and
$V$ a strongly $A$-graded M\"obius (or conformal) vertex algebra. Let
$(W,Y)$ be a strongly $\tilde A$-graded $V$-module (respectively,
generalized $V$-module). Then the pair $(W',Y')$ carries a strongly
$\tilde A$-graded $V$-module (respectively, generalized $V$-module)
structure, and $(W'',Y'')=(W,Y)$. \epf
\end{theo}

\begin{defi}{\rm
The pair $(W',Y')$ in Theorem \ref{set:W'} will be called the {\em
contragredient module} of $(W,Y)$.}
\end{defi}

Let $W_1$ and $W_2$ be strongly $\tilde A$-graded (generalized)
$V$-modules and let $f:W_1\to W_2$ be a module homomorphism (which is
of course understood to preserve both the ${\mathbb C}$-grading and the
$\tilde A$-grading, and to preserve the action of ${\mathfrak s}{\mathfrak
l}(2)$ in the M\"obius case). Then by (\ref{v'vo}) and (\ref{L'(n)}),
the linear map $f':W'_2\to W'_1$ given by
\begin{equation}\label{fprime}
\langle f'(w'_{(2)}), w_{(1)}\rangle=\langle
w'_{(2)},f(w_{(1)})\rangle
\end{equation}
for any $w_{(1)}\in W_1$ and $w'_{(2)}\in W'_2$ is well defined and is
clearly a module homomorphism {}from $W'_2$ to $W'_1$.

\begin{nota}\label{MGM}
{\rm In this work we will be especially interested in the case where
$V$ is strongly graded, and we will be focusing on the category of all
strongly graded modules, for which we will use the notation ${\cal
M}_{sg}$, or the category of all strongly graded generalized modules,
which we will call ${\cal GM}_{sg}$.  {}From the above we see that in
the strongly graded case we have contravariant functors
\[
(\cdot)': (W,Y)\mapsto (W',Y'),
\]
the {\it contragredient functors}, {}from ${\cal M}_{sg}$ to itself and
{}from ${\cal GM}_{sg}$ to itself.  We also know that $V$ itself is an
object of ${\cal M}_{sg}$ (and thus of ${\cal GM}_{sg}$ as well);
recall Remark \ref{v-str-module}. Our main objects of study will be
certain full subcategories ${\cal C}$ of ${\cal M}_{sg}$ or ${\cal
GM}_{sg}$ that are closed under the contragredient functor and such
that $V\in \ob {\cal C}$.}
\end{nota}

\begin{rema}{\rm
In order to formulate certain results in this work, even in the case
when our M\"obius or conformal vertex algebra $V$ is strongly graded
we will in fact sometimes use the category whose objects are {\it all}
the modules for $V$ and whose morphisms are all the $V$-module
homomorphisms, and also the category of {\it all} the generalized
modules for $V$.  (If $V$ is conformal, then the category of all the
$V$-modules is the same whether $V$ is viewed as either conformal or
M\"obius, by Remark \ref{modulesaremodules}, and similarly for the
category of all the generalized $V$-modules.)  Note that in view of
Remark \ref{homsaregradingpreserving}, the categories ${\cal M}_{sg}$
and ${\cal GM}_{sg}$ are not full subcategories of these categories of
{\it all} modules and generalized modules.}
\end{rema}

We now recall {}from \cite{FLM2}, \cite{FHL}, \cite{DL} and \cite{LL}
the well-known principles that vertex operator algebras (which are
exactly conformal vertex algebras strongly graded with respect to the
trivial group; recall Remark \ref{rm1}) and their modules have
important ``rationality,'' ``commutativity'' and ``associativity''
properties, and that these properties can in fact be used as axioms
replacing the Jacobi identity in the definition of the notion of
vertex operator algebra.  (These principles in fact generalize to all
vertex algebras, as in \cite{LL}.)

In the propositions below, ${\C}[x_{1}, x_{2}]_{S}$ is the ring of
formal rational functions obtained by inverting (localizing with
respect to) the products of (zero or more) elements of the set $S$ of
nonzero homogeneous linear polynomials in $x_{1}$ and $x_{2}$. Also,
$\iota_{12}$ (which might also be written as $\iota_{x_{1}x_{2}}$) is
the operation of expanding an element of ${\C}[x_{1}, x_{2}]_{S}$,
that is, a polynomial in $x_{1}$ and $x_{2}$ divided by a product of
homogeneous linear polynomials in $x_{1}$ and $x_{2}$, as a formal
series containing at most finitely many negative powers of $x_{2}$
(using binomial expansions for negative powers of linear polynomials
involving both $x_{1}$ and $x_{2}$); similarly for $\iota_{21}$ and so
on. (The distinction between rational functions and formal Laurent
series is crucial.)

Let $V$ be a vertex operator algebra.  For $W$ a ($\C$-graded)
$V$-module (including possibly $V$ itself), the space $W'$ is just the
``restricted dual space''
\begin{equation}
W'=\coprod W_{(n)}^{*}.
\end{equation}

\begin{propo}\label{rationalityandcommutativity}
We have:
\begin{enumerate}
\item[(a)] (rationality of products) For $v$, $v_{1}$, $v_{2}\in V$
and $v'\in V'$, the formal series
\begin{equation}\label{v'Yv1v2v}
\left\langle v', Y(v_{1}, x_{1})Y(v_{2}, x_{2})v\right\rangle,
\end{equation}
which involves only finitely
many negative powers of $x_{2}$ and only finitely many positive powers
of $x_{1}$, lies in the image of the map $\iota_{12}$:
\begin{equation}\label{Yv1v2}
\left\langle v', Y(v_{1}, x_{1})Y(v_{2}, x_{2})v\right\rangle
=\iota_{12}f(x_{1}, x_{2}),
\end{equation}
where the (uniquely determined) element $f\in {\C}[x_{1},
x_{2}]_{S}$ is of the form
\begin{equation}
f(x_{1}, x_{2})={\displaystyle \frac{g(x_{1},
x_{2})}{x_{1}^{r}x_{2}^{s}(x_{1}-x_{2})^{t}}}
\end{equation}
for some $g\in {\C}[x_{1}, x_{2}]$ and $r, s, t\in {\Z}$.
\item[(b)] (commutativity) We also have 
\begin{equation}\label{Yv2v1}
\left\langle v', Y(v_{2}, x_{2})Y(v_{1}, x_{1})v\right\rangle
=\iota_{21}f(x_{1}, x_{2}).
\end{equation}
\end{enumerate}
\end{propo}

\begin{propo}\label{rationalityofiterates}
We have:
\begin{enumerate}
\item[(a)] (rationality of iterates) For $v$, $v_{1}$, $v_{2}\in V$
and $v'\in V'$, the formal series
\begin{equation}\label{v'YYv1v2v}
\left\langle v', Y(Y(v_{1}, x_{0})v_{2}, x_{2})v\right\rangle,
\end{equation}
which involves only finitely many negative powers of $x_{0}$ and only
finitely many positive powers of $x_{2}$, lies in the image of the map
$\iota_{20}$:
\begin{equation}
\left\langle v', Y(Y(v_{1}, x_{0})v_{2},
x_{2})v\right\rangle=\iota_{20}h(x_{0}, x_{2}),
\end{equation}
where the (uniquely determined) element $h\in {\C}[x_{0},
x_{2}]_{S}$ is of the form
\begin{equation}
h(x_{0}, x_{2})={\displaystyle \frac{k(x_{0},
x_{2})}{x_{0}^{r}x_{2}^{s}(x_{0}+x_{2})^{t}}}
\end{equation}
for some $k\in {\C}[x_{0}, x_{2}]$ and $r, s, t\in {\Z}$.
\item[(b)] The formal series $\left\langle v', Y(v_{1},
x_{0}+x_{2})Y(v_{2}, x_{2})v\right\rangle,$ which involves only
finitely many negative powers of $x_{2}$ and only finitely many
positive powers of $x_{0}$, lies in the image of $\iota_{02}$, and in
fact
\begin{equation}
\left\langle v', Y(v_{1}, x_{0}+x_{2})Y(v_{2},
x_{2})v\right\rangle=\iota_{02}h(x_{0}, x_{2}).
\end{equation}
\end{enumerate}
\end{propo}

\begin{propo}\label{associativity} 
(associativity) We have the following equality of formal rational
functions:
\begin{equation}
\iota_{12}^{-1}\left\langle v', Y(v_{1}, x_{1})Y(v_{2}, x_{2})v\right\rangle
=(\iota_{20}^{-1}\left\langle v', Y(Y(v_{1}, x_{0})v_{2},
x_{2})v\right\rangle)\lbar_{x_{0}=x_{1}-x_{2}},
\end{equation}
that is,
\[
f(x_1,x_2)=h(x_1-x_2,x_2).
\]
\end{propo}

\begin{propo}\label{commandassocequivtoJacobi} 
In the presence of the other axioms for the notion of vertex operator
algebra, the Jacobi identity follows {from} the rationality of
products and iterates, and commutativity and associativity.  In
particular, in the definition of vertex operator algebra, the Jacobi
identity may be replaced by these properties.
\end{propo}

The rationality, commutativity and associativity properties
immediately imply the following result, in which the formal variables
$x_1$ and $x_2$ are specialized to nonzero complex numbers in suitable
domains:

\begin{corol}\label{dualitywithcovergence}
The formal series obtained by specializing $x_1$ and $x_2$ to
(nonzero) complex numbers $z_1$ and $z_2$, respectively, in
(\ref{v'Yv1v2v}) converges to a rational function of $z_1$ and $z_2$
in the domain
\begin{equation}
|z_1| > |z_2| > 0       
\end{equation}
and the analogous formal series obtained by specializing $x_1$ and
$x_2$ to $z_1$ and $z_2$, respectively, in (\ref{Yv2v1}) converges to
the same rational function of $z_1$ and $z_2$ in the (disjoint) domain
\begin{equation}
|z_2| > |z_1| > 0. 
\end{equation}
Moreover, the formal series obtained by specializing $x_0$ and $x_2$
to $z_1-z_2$ and $z_2$, respectively, in (\ref{v'YYv1v2v}) converges
to this same rational function of $z_1$ and $z_2$ in the domain
\begin{equation}
|z_2| > |z_1-z_2| > 0. 
\end{equation}
In particular, in the common domain
\begin{equation}
|z_1| > |z_2| > |z_1-z_2| > 0,
\end{equation}
we have the equality
\begin{equation}\label{associativitywithz1,z2}
\left\langle v', Y(v_{1}, z_{1})Y(v_{2}, z_{2})v\right\rangle=
\left\langle v', Y(Y(v_{1}, z_1-z_2)v_{2},z_{2})v\right\rangle
\end{equation}
of rational functions of $z_1$ and $z_2$.
\end{corol}

\begin{rema}{\rm
These last five results also hold for modules for a vertex operator
algebra $V$; in the statements, one replaces the vectors $v$ and $v'$
by elements $w$ and $w'$ of a $V$-module $W$ and its restricted dual
$W'$, respectively, and Proposition \ref{commandassocequivtoJacobi}
becomes: Given a vertex operator algebra $V$, in the presence of the
other axioms for the notion of $V$-module, the Jacobi identity follows
{from} the rationality of products and iterates, and commutativity and
associativity.  In particular, in the definition of $V$-module, the
Jacobi identity may be replaced by these properties.  }
\end{rema}

For either vertex operator algebras or modules, it is sometimes
convenient to express the equalities of rational functions in
Corollary \ref{dualitywithcovergence} informally as follows:
\begin{equation}\label{commutativityasoperatorvaluedratfns}
Y(v_{1}, z_{1})Y(v_{2}, z_{2}) \sim Y(v_{2}, z_{2})Y(v_{1}, z_{1})
\end{equation}
and
\begin{equation}\label{associativityasoperatorvaluedratfns}
Y(v_{1}, z_{1})Y(v_{2}, z_{2}) \sim Y(Y(v_{1}, z_1-z_2)v_{2},z_{2}),
\end{equation}
meaning that these expressions, defined in the domains indicated in
Corollary \ref{dualitywithcovergence} when the ``matrix coefficients''
of these expressions are taken as in this corollary, agree as
operator-valued rational functions, up to analytic continuation.

\begin{rema}\label{OPE}{\rm
Formulas (\ref{commutativityasoperatorvaluedratfns}) and
(\ref{associativityasoperatorvaluedratfns}) (or more precisely,
(\ref{associativitywithz1,z2})), express the meromorphic, or
single-valued, version of ``duality,'' in the language of conformal
field theory.  Formulas (\ref{associativityasoperatorvaluedratfns})
(and (\ref{associativitywithz1,z2})) express the existence and
associativity of the single-valued, or meromorphic, operator product
expansion.  This is the statement that the product of two (vertex)
operators can be expanded as a (suitable, convergent) infinite sum of
vertex operators, and that this sum can be expressed in the form of an
iterate of vertex operators, parametrized by the complex numbers
$z_1-z_2$ and $z_2$, in the format indicated; the infinite sum comes
{}from expanding $Y(v_{1}, z_1-z_2)v_{2},z_{2})$ in powers of $z_1-z_2$.
A central goal of this work is to generalize
(\ref{commutativityasoperatorvaluedratfns}) and
(\ref{associativityasoperatorvaluedratfns}), or more precisely,
(\ref{associativitywithz1,z2}), to logarithmic intertwining operators
in place of the operators $Y(\cdot,z)$.  This will give the existence
and also the associativity of the general, nonmeromorphic operator
product expansion.  This was done in the non-logarithmic setting in
\cite{tensor1}--\cite{tensor3} and \cite{tensor4}.  In the next
section, we shall develop the concept of logarithmic intertwining
operator.
}
\end{rema}

\newpage

\setcounter{equation}{0} \setcounter{rema}{0}

\section{Logarithmic intertwining operators}

In this section we study the notion of ``logarithmic intertwining
operator'' introduced in \cite{Mi}.  For this, we will need to discuss
spaces of formal series in powers of both $x$ and ``$\log x$'', a new
formal variable, with coefficients in certain vector spaces. We
establish certain properties, some of them quite subtle, of the formal
derivative operator $d/dx$ acting on such spaces. Then, following
\cite{Mi} with a slight variant (see Remark \ref{log:compM}), we
introduce the notion of logarithmic intertwining operator. These are
the appropriate replacements of ordinary intertwining operators when
$L(0)$-semisimplicity is relaxed.  In the strongly graded setting, it
will be natural to consider the associated ``grading-compatible''
logarithmic intertwining operators.  We work out some important
principles and formulas concerning logarithmic intertwining operators,
certain of which turn out to be the same as in the ordinary
intertwining operator case.  Some of these results require proofs that
are quite delicate.

Recall the notation ${\cal W}\{x\}$
(\ref{formalserieswithcomplexpowers}) for the space of formal series
in a formal variable $x$ with coefficients in a vector space ${\cal
W}$, with arbitrary complex powers of $x$.

{}From now on we will sometimes need and use new independent
(commuting) formal variables called $\log x$, $\log y$, $\log x_1$,
$\log x_2, \dots,$ etc. We will work with formal series in such formal
variables together with the ``usual'' formal variables $x$, $y$,
$x_1$, $x_2, \dots,$ etc., with coefficients in certain vector spaces,
and the powers of the monomials in {\it all} the variables can be
arbitrary complex numbers.  (Later we will restrict our attention to
only nonnegative integral powers of the ``$\log$'' variables.)

Given a formal variable $x$, we will use the notation $\frac d{dx}$ to
denote the linear map on ${\cal W}\{x,\log x\}$, for any vector space
${\cal W}$ not involving $x$, defined (and indeed well defined) by the
(expected) formula
\begin{eqnarray}\label{ddxdef}
\lefteqn{\frac d{dx}\bigg(\sum_{m,n\in {\mathbb C}}w_{n,m}x^n(\log
x)^m\bigg)}\nno\\
&&\;\;\;=\sum_{m,n\in{\mathbb C}}((n+1)w_{n+1,m}+ (m+1)w_{n+1,m+1})x^{n}(\log
x)^{m}\nno\\
&&\bigg(=\sum_{m,n\in{\mathbb C}}nw_{n,m}x^{n-1}(\log x)^m+ \sum_{m,n\in
{\mathbb C}}mw_{n,m}x^{n-1}(\log x)^{m-1}\bigg)
\end{eqnarray}
where $w_{n,m}\in {\cal W}$ for all $m, n\in {\mathbb C}$. We will
also use the same notation for the restriction of $\frac d{dx}$ to any
subspace of ${\cal W}\{x,\log x\}$ which is closed under $\frac
d{dx}$, e.g., ${\cal W}\{x\}[[\log x]]$ or
${\mathbb C}[x,x^{-1},\log x]$. Clearly, $\frac
d{dx}$ acting on ${\cal W}\{x\}$ coincides with the usual formal
derivative.

\begin{rema}\label{ddxchk}{\rm
Let $f$, $g$ and $f_i$, $i$ in some index set $I$, all be formal
series of the form
\begin{equation}\label{log:f}
\sum_{m,n\in {\mathbb C}}w_{n,m}x^n(\log x)^m\in {\cal W}\{x,\log x\},\;\;
w_{n,m}\in {\cal W}.
\end{equation}
One checks the following straightforwardly:
Suppose that the sum of $f_i$, $i\in I$, exists (in the obvious
sense). Then the sum of the $\frac d{dx}f_i$,
$i\in I$, also exists and is equal to $\frac d{dx}\sum_{i\in I}f_i$.
More generally, for any $T=p(x)\frac{d}{dx}$, $p(x)\in {\mathbb
C}[x,x^{-1}]$, the sum of $Tf_i$, $i\in I$, exists and is equal to
$T\sum_{i\in I}f_i$. Thus the sum of $e^{yT}f_i$, $i\in I$, exists and is
equal to $e^{yT}\sum_{i\in I}f_i$ ($e^{yT}$ being the
formal exponential series, as usual).  Suppose that ${\cal W}$ is an
(associative) algebra or that the coefficients of either $f$ or $g$
are complex numbers. If the product of $f$ and $g$
exists, then the product of $\frac d{dx}f$ and $g$ and the product of
$f$ and $\frac d{dx}g$ both exist, and $\frac d{dx}(fg)=(\frac
d{dx}f)g+f(\frac d{dx}g)$. Furthermore, for any $T$ as before, the
product of $Tf$ and $g$ and the product of $f$ and $Tg$ both exist,
and $T(fg)=(Tf)g+f(Tg)$. In addition, the product of $e^{yT}f$ and
$e^{yT}g$ exists and is equal to $e^{yT}(fg)$, just as in formulas
(8.2.6)--(8.2.10) of \cite{FLM2}.  The point here, of
course, is just the formal derivation property of $\frac d{dx}$,
except that sums and products of expressions do not exist in general.
}
\end{rema}

\begin{rema}{\rm
Note that the ``equality'' $x=e^{\log x}$ does not hold, since the
left-hand side is a formal variable, while the right-hand side
is a formal series in another formal variable. In fact, this formula
should not be assumed, since, for example,
the formal delta function $\delta(x)=\sum_{n\in {\mathbb Z}}x^n$ would not exist
in the sense of formal calculus, if $x$ were allowed to be replaced by
the formal series $e^{\log x}$. By contrast, note that the equality
\begin{equation}\label{log:logex}
\log e^x=x
\end{equation}
does indeed hold. This is because the formal series $e^x$ is of the
form $1+X$ where $X$ involves only positive integral powers of
$x$ and in (\ref{log:logex}), ``$\log$'' refers to the usual formal
logarithmic series
\begin{equation}\label{log:usual}
\log(1+X)=\sum_{i\geq 1} \frac{(-1)^{i-1}}i X^i,
\end{equation}
{\em not} to the ``$\log$'' of a formal variable. We will use the
symbol ``$\log$'' in both ways, and the meaning will be clear in
context.  }
\end{rema}

We will typically use notations of the form $f(x)$, instead of
$f(x,\log x)$, to denote elements of ${\cal W}\{x,\log x\}$ for some
vector space ${\cal W}$ as above. For this reason, we need to
interpret carefully the meaning of symbols such as $f(x+y)$, or more
generally, symbols obtained by replacing $x$ in $f(x)$ by something
other than just a single formal variable (since $\log x$ is a formal
variable and not the image of some operator acting on $x$). For the
three main types of cases that will be encountered in this work, we
use the following notational conventions; the existence of the
expressions will be justified in Remark \ref{log:exist}:

\begin{nota}
{\rm For formal variables $x$ and $y$, and $f(x)$ of the form
(\ref{log:f}), we define
\begin{eqnarray}\label{log:not1}
f(x+y)&=&\sum_{m,n\in {\mathbb C}}w_{n,m}(x+y)^n\bigg(\log x+
\log\bigg(1+\frac yx\bigg)\bigg)^m\nno\\
&=&\sum_{m,n\in {\mathbb C}}w_{n,m}(x+y)^n\bigg(\log x+
\sum_{i\geq 1} \frac{(-1)^{i-1}}i \bigg(\frac yx\bigg)^i\bigg)^m
\end{eqnarray}
(recall (\ref{log:usual})); in the right-hand side, $(\log x+\sum_{i\geq
1}\frac{(-1)^{i-1}}i
(\frac yx)^i)^m$, according to the binomial expansion convention, is to be
expanded in nonnegative integral powers of the second summand
$\sum_{i\geq 1} \frac{(-1)^{i-1}}i (\frac yx)^i$, so the right-hand
side of (\ref{log:not1}) is equal to
\begin{equation}\label{log:1-tmp}
\sum_{m,n\in {\mathbb C}}w_{n,m}(x+y)^n\sum_{j\in {\mathbb N}}{m\choose j}
(\log x)^{m-j}\bigg(\sum_{i\geq 1} \frac{(-1)^{i-1}}i \bigg(\frac
yx\bigg)^i\bigg)^j
\end{equation}
when expanded one step further. Also define
\begin{equation}\label{log:not2}
f(xe^y)=\sum_{m,n\in {\mathbb C}}w_{n,m}x^ne^{ny}(\log x+y)^m
\end{equation}
and
\begin{equation}\label{log:not3}
f(xy)=\sum_{m,n\in {\mathbb C}}w_{n,m}x^ny^n(\log x+\log y)^m,
\end{equation}
where the binomial expansion convention is again of course being
used.}
\end{nota}

\begin{rema}\label{log:exist}{\rm
The existence of the right-hand side of (\ref{log:not1}), or
(\ref{log:1-tmp}), can be seen by writing $(x+y)^n$ as $x^n(1+\frac
yx)^n$ and observing that
\[
\bigg(\sum_{i\geq 1} \frac{(-1)^{i-1}}i \bigg(\frac yx\bigg)^i\bigg)^j
\in \bigg(\frac yx\bigg)^j{\mathbb C}\bigg[\bigg[\frac yx\bigg]\bigg].
\]
The existence of the right-hand sides of (\ref{log:not2}) and of
(\ref{log:not3}) is clear. Furthermore, both $f(x+y)$ and $f(xe^y)$
lie in ${\cal W}\{x,\log x\}[[y]]$, while $f(xy)$ lies in ${\cal
W}\{xy,\log x\}[[\log y]]$. One might expect that $f(x+y)$ can be
written as $e^{y\frac{d}{dx}}f(x)$, and $f(xe^y)$ as
$e^{yx\frac{d}{dx}}f(x)$ (cf.\ Section 8.3 of \cite{FLM2}), but these
formulas must be verified (see Theorem \ref{log:ids} below).
}
\end{rema}

\begin{rema}\label{subchk}{\rm
It is clear that when there is no $\log x$ involved in $f(x)$, the
expression $f(x+y)$ (respectively, $f(xe^y)$, $f(xy)$) coincides with
the usual formal operation of substitution of $x+y$ (respectively,
$xe^y$, $xy$) for $x$ in $f(x)$. In general, it is straightforward to check
that if the sum of $f_i(x)$, $i\in I$, exists and is equal to $f(x)$,
then the sum of $f_i(x+y)$, $i\in I$ (respectively, $f_i(xe^y), i\in I$,
$f_i(xy), i\in I$), also exists and is equal to $f(x+y)$ (respectively,
$f(xe^y)$, $f(xy)$). Also, suppose that ${\cal W}$ is an (associative)
algebra or that the coefficients of either $f$ or $g$ are complex numbers.
If the product of $f(x)$ and $g(x)$ exists, then
the product of $f(x+y)$ and $g(x+y)$ (respectively, $f(xe^y)$ and
$g(xe^y)$, $f(xy)$ and $g(xy)$) also exists and is equal to
$(fg)(x+y)$ (respectively, $(fg)(xe^y)$, $(fg)(xy)$).  }
\end{rema}

Note that by (\ref{log:not1}),
\begin{equation}\label{logx+y}
\log(x+y)=\log x+\sum_{i\geq 1} \frac{(-1)^{i-1}}i \bigg(\frac
yx\bigg)^i=e^{y\frac d{dx}}\log x.
\end{equation}
The next result includes a generalization of this to arbitrary
elements of ${\cal W}\{x,\log x\}$.  Formula (\ref{log:ck1}) is a
formal ``Taylor theorem'' for logarithmic formal series.  In the case
of non-logarithmic formal series, this principle is used extensively
in vertex operator algebra theory; recall (\ref{formalTaylortheorem})
above and see Proposition 8.3.1 of \cite{FLM2} for the proof in the
generality of formal series with arbitrary complex powers of the
formal variable.  In the logarithmic case below, a much more elaborate
proof is required than in the non-logarithmic case.  The other formula
in Theorem \ref{log:ids}, formula (\ref{log:ck2}), is easier to prove.
It too is important (in the non-logarithmic case) in vertex operator
algebra theory; again see Proposition 8.3.1 of \cite{FLM2}.

\begin{theo}\label{log:ids}
For $f(x)$ as in (\ref{log:f}), we have
\begin{equation}\label{log:ck1}
e^{y\frac d{dx}}f(x)=f(x+y)
\end{equation}
(``Taylor's theorem'' for logarithmic formal series) and
\begin{equation}\label{log:ck2}
e^{yx\frac d{dx}}f(x)=f(xe^y).
\end{equation}
\end{theo}
\pf By Remarks \ref{ddxchk} and \ref{subchk} (or \ref{log:exist}), we
need only prove these equalities for $f(x)=x^n$ and $f(x)=(\log x)^m$,
$m,n\in {\mathbb C}$.  The case $f(x)=x^n$ easily follows {}from the
direct expansion of the two sides of (\ref{log:ck1}) and of
(\ref{log:ck2}) (see Proposition 8.3.1 in \cite{FLM2}). Now assume
that $f(x)=(\log x)^m$, $m\in {\mathbb C}$.

Formula (\ref{log:ck2}) is easier, so we prove it first. By
(\ref{ddxdef}) we have
\[
x\frac d{dx}(\log x)^m=m(\log x)^{m-1},
\]
so that for $k\in {\mathbb N}$,
\[
\bigg(x\frac d{dx}\bigg)^k(\log x)^m=m(m-1)\cdots(m-k+1)(\log
x)^{m-k}=k!{m\choose k}(\log x)^{m-k}.
\]
Thus
\[
e^{yx\frac d{dx}}(\log x)^m=\sum_{k\in {\mathbb N}}
\frac{y^k}{k!}\bigg(x\frac{d}{dx}\bigg)^k(\log x)^m=\sum_{k\in {\mathbb
N}} {m\choose k}y^k (\log x)^{m-k}=(\log x+y)^m,
\]
as we want.

For (\ref{log:ck1}), we shall give two proofs --- an analytic proof and an
algebraic proof. First, consider the analytic function $(\log
z)^m=e^{m\log\log z}$ over, say, $|z-3|<1$ in the complex plane.
In this proof we take the branch of $\log z$ so that
\begin{equation}\label{log:br1}
-\pi< {\rm Im }(\log z)\leq\pi.
\end{equation}
Then by analyticity, for any $z$ in
this domain, when $|z_1|$ is small enough the Taylor series expansion
$e^{z_1\frac{d}{dz}}(\log z)^m$ converges absolutely to
$(\log(z+z_1))^m$. That is,
\begin{equation}\label{log:ana1}
(\log(z+z_1))^m=e^{z_1\frac{d}{dz}}(\log z)^m
=e^{\frac{z_1}{z}z\frac{d}{dz}}(\log z)^m.
\end{equation}
Observe that as a formal series, the right-hand side of
(\ref{log:ana1}) is in the space $(\log z)^m{\mathbb C}[(\log
z)^{-1}][[z_1/z]]$.

On the other hand, by the choice of domain and the branch of $\log$ we
have
\[
\log(z+z_1)=\log z+\log(1+z_1/z)
\]
and
\[
|\log z|>\log 2>|\log(1+z_1/z)|
\]
when $|z_1|$ is small enough. So when $|z_1|$ is small enough we have
\begin{eqnarray}\label{log:ana2}
(\log(z+z_1))^m&=&(\log z+\log(1+z_1/z))^m\nno\\
&=& \sum_{j\in {\mathbb N}}{m\choose j} (\log z)^{m-j}\bigg(\sum_{i\geq
1} \frac{(-1)^{i-1}}i \bigg(\frac{z_1}z\bigg)^i\bigg)^j.
\end{eqnarray}
Since as formal series, the right-hand sides of (\ref{log:ana1}) and
(\ref{log:ana2}) are both in the space $(\log z)^m{\mathbb C}[(\log
z)^{-1}][[z_1/z]]$, and both converge to the same analytic function
$(\log(z+z_1))^m$ in the above domain, by setting $z_1=0$ in these two
functions and their derivatives with respect to $z_1$ we see that
their corresponding coefficients of powers of $z_1/z$ and further, of
all monomials in $\log z$ and $z_1/z$ must be the same. Hence we can
replace $z$ and $z_1$ by formal variables $x$ and $y$, respectively,
and obtain (\ref{log:ck1}) for $f(x)=(\log x)^m$.

An algebraic proof of (\ref{log:ck1}) (for $(\log x)^m$) can be given
as follows: Since
\[
\frac d{dx}(\log x)^m=mx^{-1}(\log x)^{m-1}
\]
and higher derivatives involve derivatives of products of powers of
$x$ and powers of $\log x$, let us first compute $(d/dx)^k(x^n(\log
x)^m)$ directly for all $m,n\in {\mathbb C}$ and $k\in{\mathbb N}$. Define
linear maps $T_0$ and $T_1$ on ${\mathbb C}\{x,\log x\}$ by setting
\[
T_0x^n(\log x)^m=nx^{n-1}(\log x)^m\;\;\mbox{ and }\;\;
T_1x^n(\log x)^m=mx^{n-1}(\log x)^{m-1},
\]
respectively, and extending to all of ${\mathbb C}\{x,\log x\}$ by formal
linearity. Then the formula
\[
\frac d{dx}x^n(\log x)^m=nx^{n-1}(\log x)^m+mx^{n-1}(\log x)^{m-1}
\]
(extended to ${\mathbb C}\{x,\log x\}$) can be written as
\[
\frac d{dx}=T_0+T_1
\]
on ${\mathbb C}\{x,\log x\}$. So for $k\geq 1$, on ${\mathbb C}\{x,\log x\}$,
\begin{eqnarray*}
&&\bigg(\frac d{dx}\bigg)^k=\sum_{(i_1,\dots,i_k)\in\{0,1\}^k}
T_{i_1}\cdots T_{i_k}\\
&=&\sum_{j=0}^{k-1}\,\sum_{0\leq t_1<t_2<\cdots<t_{k-j}<k}
T_1^{k-t_{k-j}-1}T_0T_1^{t_{k-j}-t_{k-j-1}-1}T_0\cdots
T_0T_1^{t_2-t_1-1}T_0T_1^{t_1}+T_1^k,
\end{eqnarray*}
where $j$ gives the number of $T_1$'s in the product $T_{i_1}\cdots
T_{i_k}$; there are $k-j$ $T_0$'s, which are in the following
positions, reading {}from the right: $t_1+1$, $t_2+1, \dots,$
$t_{k-j}+1$. That is, for $k\geq 1$,
\begin{eqnarray*}
\lefteqn{\bigg(\frac d{dx}\bigg)^kx^n(\log x)^m=
\sum_{j=0}^km(m-1)\cdots(m-j+1)\cdot}\\
&&\cdot\bigg(\sum_{0\leq t_1<t_2<\cdots<t_{k-j}<k}
(n-t_1)(n-t_2)\cdots(n-t_{k-j})\bigg) x^{n-k}(\log x)^{m-j},
\end{eqnarray*}
where it is understood that if $j=k$, then the latter sum (in
parentheses) is $1$. In this formula, setting $n=0$, multiplying by
$y^k/k!$, and then summing over $k\in {\mathbb N}$, we get (noting that
$t_1=0$ contributes $0$)
\begin{eqnarray}\label{log:alg1}
\lefteqn{e^{y\frac d{dx}}(\log x)^m=\sum_{k\in {\mathbb N}}\bigg(\frac
yx\bigg)^k
\sum_{j=0}^k{m\choose j}(\log x)^{m-j}\frac{j!}{k!}\cdot}\nno\\
&&\cdot\bigg(\sum_{0<t_1<t_2<\cdots<t_{k-j}<k}
(-t_1)(-t_2)\cdots(-t_{k-j})\bigg).
\end{eqnarray}
So (\ref{log:ck1}) for $f(x)=(\log x)^m$ is equivalent to equating the
right-hand side of (\ref{log:alg1}) to
\begin{eqnarray}\label{log:alg2}
\lefteqn{\bigg(\log x+\sum_{i\geq 1} \frac{(-1)^{i-1}}i \bigg(\frac
  yx\bigg)^i\bigg)^m=}\nno\\
&&=\sum_{j\in {\mathbb N}}{m\choose j}(\log x)^{m-j} \bigg(\sum_{i\geq
1}\frac{(-1)^{i-1}}i \bigg(\frac yx\bigg)^i\bigg)^j\nno\\
&&=\sum_{j\in {\mathbb N}}{m\choose j}(\log x)^{m-j}\sum_{k\geq j}
\Bigg(\sum_{
\mbox{
\tiny
$\begin{array}{c}i_1+\cdots+i_j=k\\1\leq
i_1,\dots,i_j\leq k\end{array}$
}
} \frac{(-1)^{k-j}}{i_1i_2\cdots
i_j}\Bigg)\bigg(\frac yx\bigg)^k\nno\\
&&=\sum_{k\in {\mathbb N}}\bigg(\frac yx\bigg)^k\sum_{j=0}^k{m\choose
j}(\log x)^{m-j}\Bigg(\sum_{
\mbox{
\tiny
$\begin{array}{c}i_1+\cdots+i_j=k\\1
\leq i_1,\dots,i_j\leq k\end{array}$
}
} \frac{(-1)^{k-j}}{i_1i_2\cdots
i_j}\Bigg).
\end{eqnarray}
Comparing the right-hand sides of (\ref{log:alg1}) and
(\ref{log:alg2}) we see that it is equivalent to proving the
combinatorial identity
\begin{eqnarray}\label{log:comb}
\frac{j!}{k!}\sum_{0<t_1<t_2<\cdots<t_{k-j}<k} t_1t_2\cdots
t_{k-j}= \sum_{
\mbox{
\tiny
$\begin{array}{c}i_1+\cdots+i_j=k\\1\leq
i_1,\dots,i_j\leq k\end{array}$
}
} \frac{1}{i_1i_2\cdots i_j}
\end{eqnarray}
for all $k\in {\mathbb N}$ and $j=0,\dots,k$.
Note that there is no $m$ involved here.  But for $m$ a
nonnegative integer, (\ref{log:ck1}) for $f(x)=(\log x)^m$ follows
{}from
\[
e^{y\frac d{dx}}(\log x)^m=(e^{y\frac d{dx}}\log x)^m
\]
(recall Remark \ref{ddxchk}) and
\[
e^{y\frac d{dx}}\log x=\log x+\sum_{i\geq 1}\frac{(-1)^{i-1}}i \bigg(\frac
yx\bigg)^i
\]
(recall (\ref{logx+y})).  Thus the expressions in (\ref{log:alg1}) and
(\ref{log:alg2}) are equal for any such $m$.  Equating coefficients
and choosing $m\geq j$ gives us (\ref{log:comb}).  Therefore
(\ref{log:ck1}) also holds for any $m\in{\mathbb C}$.  \epf

\begin{rema}\label{log:ids-rm}
{\rm Note that in both the left- and right-hand sides of
(\ref{log:ck1}) or (\ref{log:ck2}), $y$ can be replaced by a suitable
formal series, for example, by an element of $y\C[x][[y]]$, and
Theorem \ref{log:ids} in fact still holds if $y$ is replaced by an
arbitrary formal series in $y\C[x][[y]]$.  We will exploit this
later.}
\end{rema}

\begin{rema}{\rm
Here is an amusing sidelight: When we were writing up the proof above,
one of us (L.Z.) happened to pick up the then-current issue of the American
Mathematical Monthly and happened to notice the following problem
{}from the Problems and Solutions section, proposed by D. Lubell
\cite{Lu}:
\begin{quote}
Let $N$ and $j$ be positive integers, and let $S=\{(w_1,\dots, w_j)\in
{\mathbb Z}_+^j\,|\,0<w_1+\cdots+w_j\leq N\}$ and $T=\{(w_1,\cdots,w_j)\in
{\mathbb Z}_+^j\,|\linebreak w_1,\dots,w_j \mbox{ are distinct and bounded by }N\}.$
Show that
$$\sum_S\frac 1{w_1\cdots w_j}=\sum_T\frac 1{w_1\cdots w_j}.$$
\end{quote}
But this follows immediately {}from (\ref{log:comb}) (which is in fact a
refinement), since the left-hand side of (\ref{log:comb}) is equal to
\[
j!\sum_{1\leq w_1<w_2<\cdots<w_{j-1}\leq k-1} \frac
1{w_1w_2\cdots w_{j-1}k}=\sum_{T_k} \frac 1{w_1w_2\cdots w_j}
\]
where
\[
T_k=\{(w_1,\dots,w_j)\in\{1,2,\dots,k\}^j\,|\,w_i\mbox{ distinct,
with maximum exactly }k\},
\]
the right-hand side is
\[
\sum_{S_k} \frac 1{w_1w_2\cdots w_j}
\]
where
\[
S_k=\{(w_1,\dots,w_j)\in\{1,2,\dots,k\}^j\,|\,
w_1+\cdots+w_j=k\},
\]
and one has $S=\coprod_{k=1}^N S_k$ and $T=\coprod_{k=1}^N T_k$. }
\end{rema}

When we define the notion of logarithmic intertwining operator below,
we will impose a condition requiring certain formal series to lie in
spaces of the type ${\cal W}[\log x]\{x\}$ (so that for each power of
$x$, possibly complex, we have a {\it polynomial} in $\log x$), partly
because such results as the following (which is expected) will indeed
hold in our formal setup when the powers of the formal variables are
restricted in this way (cf.\ Remark \ref{log:[[]]} below).

\begin{lemma}\label{log:de}
Let $a\in {\mathbb C}$ and $m\in {\mathbb Z}_+$. If $f(x)\in {\cal W}[\log
x]\{x\}$ (${\cal W}$ any vector space not involving $x$ or $\log x$)
satisfies the formal differential equation
\begin{equation}\label{de:(xdx-a)^m}
\bigg(x\frac{d}{dx}-a\bigg)^mf(x)=0,
\end{equation}
then $f(x)\in {\cal W}x^a\oplus{\cal W}x^a\log x \oplus\cdots\oplus{\cal
W}x^a(\log x)^{m-1}$, and furthermore, if $m$ is the smallest integer
so that (\ref{de:(xdx-a)^m}) is satisfied, then the coefficient of
$x^a(\log x)^{m-1}$ in $f(x)$ is nonzero.
\end{lemma}
\pf For any $f(x)=\sum_{n,k}w_{n,k}x^n(\log x)^k\in {\cal W}\{x,\log
x\}$,
\begin{eqnarray*}
x\frac d{dx}f(x)&=&\sum_{n,k}nw_{n,k}x^n(\log x)^k+\sum_{n,k}kw_{n,k}
x^n(\log x)^{k-1}\\
&=&\sum_{n,k}(nw_{n,k}+(k+1)w_{n,k+1})x^n(\log x)^k.
\end{eqnarray*}
Thus for any $a\in {\mathbb C}$,
\begin{equation}\label{de:act1}
\bigg(x\frac d{dx}-a\bigg)f(x)=\sum_{n,k}((n-a)w_{n,k}+(k+1)w_{n,k+1})
x^n(\log x)^k.
\end{equation}

Now suppose that $f(x)$ lies in ${\cal W}[\log x]\{x\}$. Let us prove
the assertion of the lemma by induction on $m$.

If $m=1$, by (\ref{de:act1}) we see that $(x\frac{d}{dx}-a)f(x)=0$
means that
\begin{equation}\label{de:m=1}
(n-a)w_{n,k}+(k+1)w_{n,k+1}=0\;\;\mbox{ for any }n\in{\mathbb C},\,
k\in{\mathbb Z}.
\end{equation}
Fix $n$. If $w_{n,k}\neq 0$ for some $k$, let $k_n$ be the smallest
nonnegative integer such that $w_{n,k}=0$ for any $k>k_n$ (such a
$k_n$ exists because $f(x)\in {\cal W}[\log x]\{x\}$). Then
\[
(n-a)w_{n,k_n}=-(k_n+1)w_{n,k_n+1}=0.
\]
But $w_{n,k_n}\neq 0$ by the choice of $k_n$, so we must have
$n=a$. Now (\ref{de:m=1}) becomes $(k+1)w_{a,k+1}=0$ for any $k\in
{\mathbb Z}$, so that $w_{a,k}=0$ unless $k=0$. Thus $f(x)=w_{a,0}x^a$.
If in addition $m=1$ is the smallest integer such that
(\ref{de:(xdx-a)^m}) holds, then $f(x)\neq 0$. So $w_{a,0}$, the
coefficient of $x^a$, is not zero.

Suppose the statement is true for $m$. Then for the case $m+1$, since
\begin{equation}\label{de:m+1}
0=\bigg(x\frac d{dx}-a\bigg)^{m+1}f(x)=\bigg(x\frac
d{dx}-a\bigg)^m\bigg(x\frac d{dx}-a\bigg)f(x)
\end{equation}
implies that
\begin{equation}\label{de:wbar}
\bigg(x\frac d{dx}-a\bigg)f(x)=\bar w_0x^a+\bar w_1x^a\log x+\cdots +
\bar w_{m-1}x^a(\log x)^{m-1}
\end{equation}
for some $\bar w_0, \bar w_1,\cdots, \bar w_{m-1}\in {\cal W}$, by
(\ref{de:act1}) we get
\begin{eqnarray*}
(n-a)w_{n,j}+(j+1)w_{n,j+1}=0&\;\;&\mbox{for any }n\neq a\mbox{ and
    any }j\in{\mathbb Z}\\
(j+1)w_{a,j+1}=\bar w_j&\;&\mbox{for any }j\in\{0,1,\dots,m-1\}\\
(j+1)w_{a,j+1}=0&\;&\mbox{for any }j\notin\{0,1,\dots,m-1\}
\end{eqnarray*}
By the same argument as above we get $w_{n,j}=0$ for any $n\neq a$ and
any $j$. So
\[
f(x)=w_{a,0}x^a+\bar w_0x^a\log x+\frac{\bar w_1}{2}x^a(\log x)^2+\cdots
+\frac{\bar w_{m-1}}mx^a(\log x)^m,
\]
as we want.  If in addition $m+1$ is the smallest integer so that
(\ref{de:m+1}) is satisfied, then by the induction assumption, $\bar
w_{m-1}$ in (\ref{de:wbar}) is not zero. So the coefficient in $f(x)$
of $x^a(\log x)^m$, $\bar w_{m-1}/m$, is not zero, as we want.  \epf

\begin{rema}\label{log:[[]]}{\rm
Note that there are solutions of the equation (\ref{de:(xdx-a)^m})
outside ${\cal W}[\log x]\{x\}$, for example, $f(x)=wx^be^{(a-b)\log
x}\in x^b{\cal W}[[\log x]]$ for any complex number $b\neq a$ and any
$0\neq w\in {\cal W}$.}
\end{rema}

Following \cite{Mi}, with a slight generalization (see Remark
\ref{log:compM}), we now introduce the notion of logarithmic
intertwining operator, together with the notion of
``grading-compatible logarithmic intertwining operator,'' adapted to
the strongly graded case.  We will later see that the axioms in these
definitions correspond exactly to those in the notion of certain
``intertwining maps'' (see Definition \ref{im:imdef} below).

\begin{defi}\label{log:def}{\rm
Let $(W_1,Y_1)$, $(W_2,Y_2)$ and $(W_3,Y_3)$ be generalized modules
for a M\"obius (or conformal) vertex algebra $V$. A {\em logarithmic
intertwining operator of type ${W_3\choose W_1\,W_2}$} is a linear map
\begin{equation}\label{log:map0}
{\cal Y}(\cdot, x)\cdot: W_1\otimes W_2\to W_3[\log x]\{x\},
\end{equation}
or equivalently,
\begin{equation}\label{log:map}
w_{(1)}\otimes w_{(2)}\mapsto{\cal Y}(w_{(1)},x)w_{(2)}=\sum_{n\in {\mathbb
C}}\sum_{k\in {\mathbb N}}{w_{(1)}}_{n;\,k}^{\cal Y}w_{(2)}x^{-n-1}(\log
x)^k\in W_3[\log x]\{x\}
\end{equation}
for all $w_{(1)}\in W_1$ and $w_{(2)}\in W_2$, such that the following
conditions are satisfied: the {\em lower truncation condition}: for
any $w_{(1)}\in W_1$, $w_{(2)}\in W_2$ and $n\in {\mathbb C}$,
\begin{equation}\label{log:ltc}
{w_{(1)}}_{n+m;\,k}^{\cal Y}w_{(2)}=0\;\;\mbox{ for }\;m\in {\mathbb N}
\;\mbox{ sufficiently large,\, independently of}\;k;
\end{equation}
the {\em Jacobi identity}:
\begin{eqnarray}\label{log:jacobi}
\lefteqn{\dps x^{-1}_0\delta \bigg( {x_1-x_2\over x_0}\bigg)
Y_3(v,x_1){\cal Y}(w_{(1)},x_2)w_{(2)}}\nno\\
&&\hspace{2em}- x^{-1}_0\delta \bigg( {x_2-x_1\over -x_0}\bigg)
{\cal Y}(w_{(1)},x_2)Y_2(v,x_1)w_{(2)}\nno \\
&&{\dps = x^{-1}_2\delta \bigg( {x_1-x_0\over x_2}\bigg) {\cal
Y}(Y_1(v,x_0)w_{(1)},x_2) w_{(2)}}
\end{eqnarray}
for $v\in V$, $w_{(1)}\in W_1$ and $w_{(2)}\in W_2$ (note that the
first term on the left-hand side is meaningful because of
(\ref{log:ltc})); the {\em $L(-1)$-derivative property:} for any
$w_{(1)}\in W_1$,
\begin{equation}\label{log:L(-1)dev}
{\cal Y}(L(-1)w_{(1)},x)=\frac d{dx}{\cal Y}(w_{(1)},x);
\end{equation}
and the {\em ${\mathfrak s}{\mathfrak l}(2)$-bracket relations:} for any
$w_{(1)}\in W_1$,
\begin{equation}\label{log:L(j)b}
{}[L(j), {\cal Y}(w_{(1)},x)]=\sum_{i=0}^{j+1}{j+1\choose i}x^i{\cal
Y}(L(j-i)w_{(1)},x)
\end{equation}
for $j=-1, 0$ and $1$ (note that if $V$ is in fact a conformal vertex
algebra, this follows automatically {}from the Jacobi identity
(\ref{log:jacobi}) by setting $v=\omega$ and then taking
$\res_{x_0}\res_{x_1}x_1^{j+1}$).  }
\end{defi}

\begin{rema}{\rm
We will sometimes write the Jacobi identity (\ref{log:jacobi}) as
\begin{eqnarray}
\lefteqn{\dps x^{-1}_0\delta \bigg( {x_1-x_2\over x_0}\bigg)
Y(v,x_1){\cal Y}(w_{(1)},x_2)w_{(2)}}\nno\\
&&\hspace{2em}- x^{-1}_0\delta \bigg( {x_2-x_1\over -x_0}\bigg)
{\cal Y}(w_{(1)},x_2)Y(v,x_1)w_{(2)}\nno \\
&&{\dps = x^{-1}_2\delta \bigg( {x_1-x_0\over x_2}\bigg) {\cal
Y}(Y(v,x_0)w_{(1)},x_2) w_{(2)}}
\end{eqnarray}
(dropping the subscripts on the module actions) for brevity.}
\end{rema}

\begin{rema}\label{ordinaryandlogintwops}{\rm
The ordinary intertwining operators (as in, for example,
\cite{tensor1}) among triples of modules for a vertex operator algebra
are exactly the logarithmic intertwining operators that do not involve
the formal variable $\log x$, except for our present relaxation of the
lower truncation condition.  The lower truncation condition that we
use here can be equivalently stated as: For any $w_{(1)}\in W_1$,
$w_{(2)}\in W_2$ and $n\in{\mathbb C}$, there is no nonzero term
involving $x^{n-m}$ appearing in ${\cal Y}(w_{(1)},x)w_{(2)}$ when
$m\in{\mathbb N}$ is large enough.  In \cite{tensor1}, the lower
truncation condition in the definition of the notion of intertwining
operator states: For any $w_{(1)}\in W_1$ and $w_{(2)}\in W_2$,
\[
(w_{(1)})_nw_{(2)}=0\;\mbox{for}\;n\;\mbox{whose real part is
sufficiently large}.
\]
This is slightly stronger than the lower truncation condition that we
use here, even if no $\log x$ is involved, when the powers of $x$ in
${\cal Y}(w_{(1)},x)w_{(2)}$ belong to infinitely many different
congruence classes modulo ${\mathbb Z}$. }
\end{rema}

\begin{rema}\label{g-mod-as-l-int}{\rm
Given a generalized module $(W, Y_{W})$ for a M\"{o}bius (or conformal) 
vertex algebra $V$, the vertex operator map $Y_{W}$ itself is clearly 
a logarithmic intertwining operator of type ${W\choose VW}$; in fact, it 
does not involve $\log x$ and its powers of $x$ are all integers. In 
particular, taking $(W, Y_{W})$ to be $(V, Y)$ itself, we have that the 
vertex operator map $Y$ is a logarithmic intertwining operator of type
${V\choose VV}$ not involving $\log x$ and having only integral powers 
of $x$.}
\end{rema}

The logarithmic intertwining operators of a fixed type ${W_3\choose
W_1\,W_2}$ form a vector space.

\begin{defi}\label{gradingcompatintwop}{\rm
In the setting of Definition \ref{log:def}, suppose in addition that
$V$ and $W_1$, $W_2$ and $W_3$ are strongly graded (recall Definitions
\ref{def:dgv} and \ref{def:dgw}).  A logarithmic intertwining operator
${\cal Y}$ as in Definition \ref{log:def} is a {\em grading-compatible
logarithmic intertwining operator} if for $\beta, \gamma \in \tilde A$
(recall Definition \ref{def:dgw}) and $w_{(1)} \in W_{1}^{(\beta)}$,
$w_{(2)} \in W_{2}^{(\gamma)}$, $n \in \C$ and $k \in \mathbb N$, we
have
\begin{equation}\label{gradingcompatcondn}
{w_{(1)}}_{n;\,k}^{\cal Y}w_{(2)} \in W_{3}^{(\beta + \gamma)}.
\end{equation}}
\end{defi}

\begin{rema}{\rm
The term ``grading-compatible'' in Definition
\ref{gradingcompatintwop} refers to the $\tilde A$-gradings; any
logarithmic intertwining operator is compatible with the $\C$-gradings
of $W_1$, $W_2$ and $W_3$, in view of Proposition \ref{log:logwt}(b)
below.}
\end{rema}

\begin{rema}\label{str-graded-g-mod-as-l-int}
{\rm Given a strongly graded generalized module $(W, Y_{W})$ for a
strongly graded M\"{o}bius (or conformal) vertex algebra $V$, the
vertex operator map $Y_{W}$ is a grading-compatible logrithmic
intertwining operator of type ${W\choose VW}$ not involving $\log x$
and having only integral powers of $x$.  Taking $(W, Y_{W})$ in
particular to be $(V, Y)$ itself, we have that the vertex operator map
$Y$ is a grading-compatible logarithmic intertwining operator of type
${V\choose VV}$ not involving $\log x$ and having only integral powers
of $x$.}
\end{rema}

In the strongly graded context (the main context for our tensor
product theory), we will use the following notation and terminology,
traditionally used in the setting of ordinary intertwining operators,
as in \cite{FHL}:

\begin{defi}\label{fusionrule}{\rm
In the setting of Definition \ref{gradingcompatintwop}, the
grading-compatible logarithmic intertwining operators of a fixed type
${W_3\choose W_1\, W_2}$ form a vector space, which we denote by
${\cal V}^{W_3}_{W_1\,W_2}$.  We call the dimension of ${\cal
V}^{W_3}_{W_1\,W_2}$ the {\it fusion rule} for $W_1$, $W_2$ and $W_3$
and denote it by $N^{W_3}_{W_1\,W_2}$.}
\end{defi}

\begin{rema}{\rm
In the strongly graded context, suppose that $W_1$, $W_2$ and $W_3$ in
Definition \ref{log:def} are expressed as finite direct sums of
submodules.  Then the space ${\cal V}^{W_3}_{W_1\,W_2}$ can be
naturally expressed as the corresponding (finite) direct sum of the
spaces of (grading-compatible) logarithmic intertwining operators
among the direct summands, and the fusion rule $N^{W_3}_{W_1\,W_2}$ is
thus the sum of the fusion rules for the direct summands.
}
\end{rema}

\begin{rema}{\rm
As we shall point out in Remark \ref{log:ordi} below, it turns out
that the notion of fusion rule in Definition \ref{fusionrule} agrees
with the traditional notion, in the case of a vertex operator algebra
and ordinary modules.  The justification of this assertion uses Parts
(b) and (c), or alternatively, Part (a), of the next proposition.
Part (a), whose proof uses Lemma \ref{log:de}, shows how logarithmic
intertwining operators yield expansions involving only finitely many
powers of $\log x$.  Part (b) is a generalization of formula
(\ref{set:wtvn}).
}
\end{rema}

\begin{propo}\label{log:logwt}
Let $W_1$, $W_2$, $W_3$ be generalized modules for a M\"obius (or
conformal) vertex algebra $V$, and let ${\cal Y}(\cdot, x)\cdot$ be a
logarithmic intertwining operator of type ${W_3\choose W_1\,W_2}$.
Let $w_{(1)}$ and $w_{(2)}$ be homogeneous elements of $W_1$ and $W_2$
of generalized weights $n_1$ and $n_2 \in {\mathbb C}$, respectively, and
let $k_1$ and $k_2$ be positive integers such that
$(L(0)-n_1)^{k_1}w_{(1)}=0$ and $(L(0)-n_2)^{k_2}w_{(2)}=0$.  Then we
have:

(a) (\cite{Mi}) For any $w'_{(3)}\in W_3^*$, $n_3\in {\mathbb C}$ and
$k_3\in {\mathbb Z}_+$ such that $(L'(0)-n_3)^{k_3}w'_{(3)}=0$,
\begin{eqnarray}\label{log:k}
\lefteqn{\langle w'_{(3)}, {\cal Y}(w_{(1)}, x)w_{(2)}\rangle}\nno\\
&&\in {\mathbb C}x^{n_3-n_1-n_2}\oplus{\mathbb C}x^{n_3-n_1-n_2}\log
x\oplus\cdots\oplus {\mathbb C}x^{n_3-n_1-n_2}(\log
x)^{k_1+k_2+k_3-3}.\nno\\
\end{eqnarray}

(b) For any $n\in {\mathbb C}$ and $k\in {\mathbb N}$, ${w_{(1)}}^{\cal
Y}_{n;\,k}w_{(2)}\in W_3$ is homogeneous of generalized weight
$n_1+n_2-n-1$.

(c) Fix $n\in {\mathbb C}$ and $k\in {\mathbb N}$. For each $i,j\in {\mathbb
N}$, let $m_{ij}$ be a nonnegative integer such that
\[
(L(0)-n_1-n_2+n+1)^{m_{ij}}(((L(0)-n_1)^iw_{(1)})^{\cal
Y}_{n;\,k}(L(0)-n_2)^jw_{(2)})=0.
\]
Then for all $t\geq \max\{m_{ij}\,|\,0\leq i<k_1,\; 0\leq
j<k_2\}+k_1+k_2-2$,
\[
{w_{(1)}}^{\cal Y}_{n;\,k+t}w_{(2)}=0.
\]
\end{propo}

We will need the following lemma in the proof:

\begin{lemma}\label{log:lemma}
Let $W_1$, $W_2$, $W_3$ be generalized modules for a M\"obius (or
conformal) vertex algebra $V$. Let
\begin{eqnarray}
{\cal Y}(\cdot,x)\cdot: W_1\otimes
W_2 & \to & W_3\{x,\log x\} \nonumber\\
w_{(1)}\otimes w_{(2)} & \mapsto & {\cal Y}(w_{(1)},x)w_{(2)}=\sum_{n,k\in
{\mathbb C}}{w_{(1)}}_{n;\,k}^{\cal Y}w_{(2)}x^{-n-1}(\log x)^k
\nonumber\\
& & \label{intertwopinlemma}
\end{eqnarray}
be a linear map that satisfies the L(-1)-derivative property
(\ref{log:L(-1)dev}) and the $L(0)$-bracket relation, that is,
(\ref{log:L(j)b}) with $j=0$. Then for any $a,b,c\in{\mathbb C}$,
$t\in{\mathbb N}$, $w_{(1)}\in W_1$ and $w_{(2)}\in W_2$,
\begin{eqnarray}\label{log:ty}
\lefteqn{(L(0)-c)^t{\cal Y}(w_{(1)},x)w_{(2)}=\sum_{i,j,l\in {\mathbb
N},\;i+j+l=t}\frac{t!}{i!j!l!}\cdot}\nno\\
&&\cdot\left(x\frac d{dx}-c+a+b\right)^l{\cal Y}((L(0)-a)^iw_{(1)},x)
(L(0)-b)^jw_{(2)}.\nno\\
\end{eqnarray}
Also, for any $a,b,n,k\in{\mathbb C}$, $t\in{\mathbb N}$,
$w_{(1)}\in W_1$ and $w_{(2)}\in W_2$, we have
\begin{eqnarray}\label{log:t00}
\lefteqn{(L(0)-a-b+n+1)^t({w_{(1)}}_{n;\,k}^{\cal Y}w_{(2)})}\nno\\
&&=t!\sum_{i,j,l\geq 0,\;i+j+l=t}{k+l\choose l}
\bigg(\frac{(L(0)-a)^i}{i!}w_{(1)}\bigg)^{\cal
Y}_{n;\,k+l}\bigg(\frac{(L(0)-b)^j}{j!}w_{(2)}\bigg);\nno\\
\end{eqnarray}
in generating function form, this gives
\begin{eqnarray}\label{log:e^L(0)}
\lefteqn{e^{y(L(0)-a-b+n+1)}({w_{(1)}}_{n;\,k}^{\cal
Y}w_{(2)})}\nno\\ 
&&=\sum_{l\in{\mathbb N}}{k+l\choose
l}(e^{y(L(0)-a)}w_{(1)})_{n;\,k+l}^{\cal Y}(e^{y(L(0)-b)}w_{(2)})y^l.
\end{eqnarray}
\end{lemma}

\pf {}From (\ref{log:L(-1)dev}) and (\ref{log:L(j)b}) with $j=0$ we have
\begin{eqnarray}\label{log:pf1}
\lefteqn{L(0){\cal Y}(w_{(1)},x)w_{(2)}={\cal Y}(w_{(1)},x)
L(0)w_{(2)}}\nno\\
&&+x\frac d{dx}{\cal Y}(w_{(1)},x)w_{(2)}+
{\cal Y}(L(0)w_{(1)},x)w_{(2)}.
\end{eqnarray}
Hence
\begin{eqnarray*}
\lefteqn{(L(0)-c){\cal Y}(w_{(1)},x)w_{(2)}={\cal Y}(w_{(1)},x)
(L(0)-b)w_{(2)}}\\
&&+\left(x\frac d{dx}-c+a+b\right){\cal Y}(w_{(1)},x)w_{(2)}+ {\cal
Y}((L(0)-a)w_{(1)},x)w_{(2)}
\end{eqnarray*}
for any complex numbers $a$, $b$ and $c$. In view of the fact that the
actions of $L(0)$ and $d/dx$ commute with each other, this implies
(\ref{log:ty}) essentially because of the expansion formula for powers
of a sum of commuting operators, that is, for any commuting operators
$T_1, \dots, T_s$ and $t\in {\mathbb N}$,
\begin{equation}\label{log:expand}
(T_1+\cdots+T_s)^t=\sum_{i_1,\dots,i_s\in{\mathbb N},\; i_1+\cdots+i_s=t}
\frac{t!}{i_1!\cdots i_s!}T_1^{i_1}\cdots T_s^{i_s}.
\end{equation}

On the other hand, by taking coefficient of $x^{-n-1}(\log x)^k$ in
(\ref{log:pf1}) we get
\begin{eqnarray*}
\lefteqn{L(0){w_{(1)}}_{n;\,k}^{\cal Y}w_{(2)}={w_{(1)}}_{n;\,k}^{\cal
Y}L(0)w_{(2)}+(-n-1){w_{(1)}}_{n;\,k}^{\cal Y}w_{(2)}}\nno\\
&&+(k+1){w_{(1)}}_{n;\,k+1}^{\cal
Y}w_{(2)}+(L(0)w_{(1)})_{n;\,k}^{\cal Y}w_{(2)}.
\end{eqnarray*}
So for any $a,b,n,k\in {\mathbb C}$,
\begin{eqnarray}\label{log:t=1}
\lefteqn{(L(0)-a-b+n+1)({w_{(1)}}_{n;\,k}^{\cal Y}w_{(2)})=
((L(0)-a)w_{(1)})_{n;\,k}^{\cal Y}w_{(2)}}\nno\\
&&\hspace{2cm}+{w_{(1)}}_{n;\,k}^{\cal
Y}(L(0)-b)w_{(2)}+(k+1){w_{(1)}}_{n;\,k+1}^{\cal Y}w_{(2)}.
\end{eqnarray}
For $p,q\in {\mathbb N}$ and $n,k\in {\mathbb C}$, let us write
\begin{equation}\label{log:Tprq}
T_{p,k,q}=((L(0)-a)^pw_{(1)})^{\cal Y}_{n;\,k}((L(0)-b)^qw_{(2)}).
\end{equation}
Then {}from (\ref{log:t=1}) we see that for any $p,q\in {\mathbb N}$ and
$a,b,n,k\in {\mathbb C}$,
\begin{equation}\label{log:L(0)^1}
(L(0)-a-b+n+1)T_{p,k,q}=T_{p+1,k,q}+(k+1)T_{p,k+1,q}+T_{p,k,q+1}.
\end{equation}
Hence by (\ref{log:expand}) we have
\[
(L(0)-a-b+n+1)^tT_{p,k,q}=t!\sum_{i,j,l\geq 0,\;i+j+l=t}\frac
{(k+1)(k+2)\cdots(k+l)}{i!j!l!}T_{p+i,k+l,q+j}
\]
for any $a,b,n,k\in {\mathbb C}$ and $p,q\in {\mathbb N}$. In particular, by
setting $p=q=0$ we get (\ref{log:t00}), and (\ref{log:e^L(0)}) follows
easily {}from (\ref{log:t00}) by multiplying by $y^t/t!$ and then
summing over $t\in{\mathbb N}$.  \epfv

{\it Proof of Proposition \ref{log:logwt}}\hspace{2ex} (a): Under the
assumptions of the proposition, let us show that
\begin{equation}\label{log:ode}
\left\langle w'_{(3)},\left(x\frac d{dx}-n_3+n_1+n_2\right)^
{k_3+k_1+k_2-2}{\cal Y}(w_{(1)},x)w_{(2)}\right\rangle=0
\end{equation}
by induction on $k_1+k_2$.

For $k_1=k_2=1$, {}from (\ref{log:ty}) with $a=n_1$, $b=n_2$, $c=n_3$
and $t=k_3$ we have
\begin{eqnarray*}
0&=&\langle (L'(0)-n_3)^{k_3}w'_{(3)},{\cal
Y}(w_{(1)},x)w_{(2)}\rangle\\
&=&\langle w'_{(3)},(L(0)-n_3)^{k_3}{\cal
Y}(w_{(1)},x)w_{(2)}\rangle\\
&=&\left\langle w'_{(3)},\left(x\frac d{dx}-n_3+n_1+n_2\right)^{k_3}{\cal
Y}(w_{(1)},x)w_{(2)}\right\rangle,
\end{eqnarray*}
which is (\ref{log:ode}) in the case $k_1=k_2=1$.

Suppose that (\ref{log:ode}) is true for all the cases with smaller
$k_1+k_2$. Then {}from (\ref{log:ty}) with $a=n_1$, $b=n_2$, $c=n_3$ and
$t=k_3+k_1+k_2-2$ we have
\begin{eqnarray*}
0&=&\langle (L'(0)-n_3)^{k_3+k_1+k_2-2}w'_{(3)},{\cal
Y}(w_{(1)},x)w_{(2)}\rangle\\
&=&\langle w'_{(3)},(L(0)-n_3)^{k_3+k_1+k_2-2}{\cal
Y}(w_{(1)},x)w_{(2)}\rangle\\
&=&\Biggl\langle w'_{(3)},\sum_{i,j,k\in {\mathbb
N},\; i+j+k=k_3+k_1+k_2-2}\frac{(k_3+k_1+k_2-2)!}{i!j!k!}\cdot\\
&&\quad\cdot\left(x\frac d{dx}-n_3+n_1+n_2\right)^k{\cal Y}
((L(0)-n_1)^iw_{(1)},x) (L(0)-n_2)^jw_{(2)}\Biggr\rangle\\
&=&\left\langle w'_{(3)},\left(x\frac d{dx}-n_3+n_1+n_2\right)^
{k_3+k_1+k_2-2}{\cal Y}(w_{(1)},x)w_{(2)}\right\rangle,
\end{eqnarray*}
where the last equality uses the induction assumption for the pair of
elements $(L(0)-n_1)^iw_{(1)}$ and $(L(0)-n_2)^jw_{(2)}$ for all
$(i,j)\neq (0,0)$. So (\ref{log:ode}) is established, that is, we have
the formal differential equation
\[
\left(x\frac d{dx}-n_3+n_1+n_2\right)^{k_3+k_1+k_2-2}\langle w'_{(3)},
{\cal Y}(w_{(1)},x)w_{(2)}\rangle=0.
\]
This implies (a) by Lemma \ref{log:de}.

(b): This follows {}from (\ref{log:t00}) with $a=n_1$, $b=n_2$ and the
fact that for any $\bar w_{(1)}\in W_1$, $\bar w_{(2)}\in W_2$ and
$\bar n\in {\mathbb C}$, there exists $K\in{\mathbb N}$ so that $(\bar
w_{(1)})^{\cal Y}_{\bar n;\bar k}\bar w_{(2)}=0$ for all $\bar k>K$,
due to (\ref{log:map0}).

(c): Let us prove (c) by induction on $k_1+k_2$ again. For
$k_1=k_2=1$, (\ref{log:t00}) with $a=n_1$ and $b=n_2$ gives
\[
(L(0)-n_1-n_2+n+1)^t({w_{(1)}}_{n;\,k}^{\cal
Y}w_{(2)})=((k+t)!/k!){w_{(1)}}_{n;\,k+t}^{\cal Y}w_{(2)},
\]
that is,
\[
{w_{(1)}}_{n;\,k+t}^{\cal Y}w_{(2)}= (k!/(k+t)!)
(L(0)-n_1-n_2+n+1)^t({w_{(1)}}_{n;\,k}^{\cal Y}w_{(2)}).
\]
So for $t\geq m_{00}$, ${w_{(1)}}_{n;\,k+t}^{\cal Y}w_{(2)}=0$,
proving the statement in case $k_1=k_2=1$.

Suppose that the statement (c) is true for all smaller $k_1+k_2$.
Then for $(i,j)\neq (0,0)$,
\[
\bigg(\frac{(L(0)-n_1)^i}{i!}w_{(1)}\bigg)^{\cal
Y}_{n;\,k+l}\bigg(\frac{(L(0)-n_2)^j}{j!}w_{(2)}\bigg)=0
\]
when $l\geq \max\{m_{i'j'}\,|\,i\leq i'<k_1,\; j\leq
j'<k_2\}+(k_1-i)+(k_2-j)-2$, and in particular, when $l\geq
\max\{m_{i'j'}\,|\,0\leq i'<k_1,\; 0\leq j'<k_2\}+k_1+k_2-i-j-2$. But
then for all $t\geq\max\{m_{ij}\,|\,0\leq i<k_1,\; 0\leq
j<k_2\}+k_1+k_2-2$, (\ref{log:t00}) gives
$0=((k+t)!/k!){w_{(1)}}_{n;\,k+t}^{\cal Y}w_{(2)}$, proving what we
need.  \epfv

The following corollary is immediate {}from Proposition \ref{log:logwt}(b):

\begin{corol}\label{powerscongruentmodZ}
Let $V$ be a M\"obius (or conformal) vertex algebra and let $W_1$,
$W_2$ and $W_3$ be generalized $V$-modules whose weights are all
congruent modulo ${\mathbb Z}$ to complex numbers $h_1$, $h_2$ and $h_3$,
respectively.  (For example, $W_1$, $W_2$ and $W_3$ might be
indecomposable; recall Remark \ref{congruent}.)  Let ${\cal Y}(\cdot,
x)\cdot$ be a logarithmic intertwining operator of type ${W_3\choose
W_1\,W_2}$.  Then all powers of $x$ in ${\cal Y}(\cdot, x)\cdot$ are
congruent modulo ${\mathbb Z}$ to $h_3-h_1-h_2$. \epf
\end{corol}

\begin{rema}\label{log:ordi}{\rm
Let $W_1$, $W_2$ and $W_3$ be (ordinary) modules for a M\"obius (or
conformal) vertex algebra $V$. Then any logarithmic intertwining
operator of type ${W_3\choose W_1\,W_2}$ is just an ordinary
intertwining operator of this type, i.e., it does not involve $\log
x$. This clearly follows {}from Proposition \ref{log:logwt}(b) and (c),
where $k_1$ and $k_2$ are chosen to be $1$, $k$ is chosen to be $0$,
and $m_{00}$ is chosen to be $1$.  It also follows, alternatively,
{}from Proposition \ref{log:logwt}(a).  As a result, for $V$ a vertex
operator algebra (viewed as a conformal vertex algebra strongly graded
with respect to the trivial group; recall Remark \ref{rm1}) and $W_1$,
$W_2$ and $W_3$ $V$-modules in the sense of Remark
\ref{moduleswiththetrivialgroup}, the notion of fusion rule defined in
this work (recall Definition \ref{fusionrule}) coincides with the
notion of fusion rule defined in, for example, \cite{tensor1} (except
for the minor issue of the truncation condition for an intertwining
operator, discussed in Remark \ref{ordinaryandlogintwops}).}
\end{rema}

\begin{rema}\label{log:compM}{\rm
Our definition of logarithmic intertwining operator is identical to
that in \cite{Mi} (in case $V$ is a vertex operator algebra) except
that in \cite{Mi}, a logarithmic intertwining operator ${\cal
Y}$ of type ${W_3\choose W_1\,W_2}$ is required to be a
linear map $W_1\to \hom(W_2,W_3)\{x\}[\log x]$, instead of as in
(\ref{log:map0}), and the lower truncation condition (\ref{log:ltc})
is replaced by: For any $w_{(1)}\in W_1$, $w_{(2)}\in W_2$ and $k\in
{\mathbb N}$,
\[
{w_{(1)}}^{\cal Y}_{n;k}w_{(2)}=0\;\;\mbox{for}\;n\;\mbox{whose real
part is sufficiently large}.
\]
Given generalized $V$-modules $W_1$, $W_2$ and $W_3$, suppose that for
each $i=1,2,3$, there exists some $K_i\in {\mathbb Z}_+$ such that
$(L(0)-L(0)_s)^{K_i}W_i=0$ (this is satisfied by many interesting
examples and is assumed in \cite{Mi} for all generalized modules under
consideration). Then for any logarithmic intertwining operator ${\cal
Y}$, any homogeneous elements $w_{(1)}\in {W_1}_{[n_1]}$, $w_{(2)}\in
{W_2}_{[n_2]}$, $n_1, n_2\in{\mathbb C}$, any $n\in{\mathbb C}$ and
any $k\in{\mathbb N}$, all the $m_{ij}$'s in Proposition
\ref{log:logwt}(c) can be chosen to be no greater than $K_3$, while
$k_1$ and $k_2$ can be chosen to be no greater than $K_1$ and $K_2$,
respectively.  Proposition \ref{log:logwt}(c) thus
implies that the largest power of $\log x$ that is involved in ${\cal
Y}$ is no greater than $K_1+K_2+K_3-3$. In particular, ${\cal Y}$ maps
$W_1$ to $\hom(W_2,W_3)\{x\}[\log x]$, and in fact, we even have that
$K_1+K_2+K_3-3$ is a global bound on the powers of $\log x$,
independently of $w_{(1)}\in W_1$, so that
\begin{equation}
{\cal Y}(\cdot,x)\cdot \in \hom(W_1 \otimes W_2,W_3)\{x\}[\log x].
\end{equation}
}
\end{rema}

\begin{rema}{\rm
Given a logarithmic intertwining operator ${\cal Y}$ as in
(\ref{log:map}), set
\[
{\cal Y}^{(k)}(w_{(1)},x)w_{(2)}=\sum_{n\in {\mathbb
C}}{w_{(1)}}_{n;\,k}^{\cal Y}w_{(2)}x^{-n-1}
\]
for $k\in {\mathbb N}$, $w_{(1)}\in W_1$ and $w_{(2)}\in W_2$, so that
\[
{\cal Y}(w_{(1)},x)w_{(2)}=\sum_{k\in {\mathbb N}}{\cal
  Y}^{(k)}(w_{(1)},x)w_{(2)}(\log x)^k.
\]
By taking the coefficients of the powers of $\log x_2$ and $\log x$ in
(\ref{log:jacobi}) and (\ref{log:L(j)b}), respectively, we see that
each ${\cal Y}^{(k)}$ satisfies the Jacobi identity and the ${\mathfrak
s}{\mathfrak l}(2)$-bracket relations. On the other hand, taking the
coefficients of the powers of $\log x$ in (\ref{log:L(-1)dev}) gives
\begin{equation}\label{log:L(-1)comp}
{\cal Y}^{(k)}(L(-1)w_{(1)},x)=\frac d{dx}{\cal Y}^{(k)}(w_{(1)},x)
+\frac{k+1}x{\cal Y}^{(k+1)}(w_{(1)},x)
\end{equation}
for any $k\in {\mathbb N}$ and $w_{(1)}\in W_1$. So ${\cal Y}^{(k)}$
does not in general satisfy the $L(-1)$-derivative property.  (If
${\cal Y}^{(k+1)}=0$, then ${\cal Y}^{(k)}$ of course does satisfy the
$L(-1)$-derivative property and so is an (ordinary) intertwining
operator; this certainly happens for $k=0$ in the context of Remark
\ref{log:ordi} and for $k=K_1+K_2+K_3-3$ in the context of Remark
\ref{log:compM}.)  However, in the following we will see that suitable
formal linear combinations of certain modifications of ${\cal
Y}^{(k)}$ (depending on $t\in{\mathbb N}$; see below) form a sequence
of logarithmic intertwining operators.  }
\end{rema}

\begin{rema}\label{log:mu}{\rm
Given a logarithmic intertwining operator ${\cal Y}$, let us write
\begin{eqnarray*}
{\cal Y}(w_{(1)},x)w_{(2)}&=&\sum_{n\in {\mathbb C}}\sum_{k\in {\mathbb N}}
{w_{(1)}}^{\cal Y}_{n;\,k}w_{(2)}x^{-n-1}(\log x)^k\\
&=&\sum_{\mu\in {\mathbb C}/{\mathbb Z}}\sum_{\bar n=\mu}\sum_{k\in {\mathbb
N}}{w_{(1)}}^{\cal Y}_{n;\,k}w_{(2)}x^{-n-1}(\log x)^k
\end{eqnarray*}
for any $w_{(1)}\in W_1$ and $w_{(2)}\in W_2$, where $\bar n$ denotes
the equivalence class of $n$ in ${\mathbb C}/{\mathbb Z}$. By extracting
summands corresponding to the same congruence class modulo ${\mathbb Z}$
of the powers of $x$ in (\ref{log:jacobi}), (\ref{log:L(-1)dev}) and
(\ref{log:L(j)b}) we see that for each $\mu\in {\mathbb C}/{\mathbb Z}$,
\begin{equation}\label{log:c+n}
w_{(1)}\otimes w_{(2)}\mapsto {{\cal Y}^\mu}(w_{(1)},x)w_{(2)}= 
\sum_{\bar n=\mu}\sum_{k\in {\mathbb N}}
{w_{(1)}}^{\cal Y}_{n;\,k}w_{(2)}x^{-n-1}(\log x)^k
\end{equation}
still defines a logarithmic intertwining operator.  In the strongly
graded case, if ${\cal Y}$ is grading-compatible, then so is the
operator ${\cal Y}^\mu$ in (\ref{log:c+n}).  Conversely, suppose that
we are given a family of logarithmic intertwining operators $\{{\cal
Y}^\mu | \mu\in{\mathbb C}/{\mathbb Z}\}$ parametrized by
$\mu\in{\mathbb C}/{\mathbb Z}$ such that the powers of $x$ in ${\cal
Y}^\mu$ are restricted as in (\ref{log:c+n}). Then the formal sum
$\sum_{\mu\in{\mathbb C}/{\mathbb Z}}{\cal Y}^\mu$ is well defined and
is a logarithmic intertwining operator.  In the strongly graded case,
if each ${\cal Y}^\mu$ is grading-compatible, then so is this sum.}
\end{rema}

In the setting of Definition \ref{log:def}, for any integer $p$, set
\begin{equation}\label{substitutionofe2piipx}
{\cal Y} (w_{(1)},e^{2\pi ip}x)w_{(2)} = {\cal
Y}(w_{(1)},y)w_{(2)}\lbar_{y^n=e^{2\pi ipn}x^n,\; (\log y)^k=({2\pi
ip} +\log x)^k,\;n\in{\mathbb C},\;k\in{\mathbb N}}.
\end{equation}
This is in fact a well-defined element of $W_3[\log x]\{x\}$.  Note
that this element certainly depends on $p$, not just on $e^{2\pi ip}$
($=1$).  This substitution, which can be thought of as ``$x \mapsto
e^{2\pi ip}x$,'' will be considered in a more general form in
(\ref{log:subs}) below.

\begin{rema}\label{formalinvariance}{\rm
It is clear that in Definition \ref{log:def}, for any integer $p$, all
the axioms are formally invariant under the substitution $x \mapsto
e^{2\pi ip}x$ given by (\ref{substitutionofe2piipx}).  That is, if we
apply this substitution to each axiom, the axiom keeps the same form,
with the operator ${\cal Y}(\cdot,x)$ replaced by ${\cal
Y}(\cdot,e^{2\pi ip}x)$.  For example, for the Jacobi identity
(\ref{log:jacobi}), we perform the substitution $x_{2} \mapsto e^{2\pi
ip}x_{2}$; the formal delta-functions remain unchanged because they
involve only integral powers of $x_2$ and no logarithms.  It follows
that ${\cal Y}(\cdot,e^{2\pi ip}x)$ is again a logarithmic
intertwining operator.  }
\end{rema}

{}From Remark \ref{formalinvariance}, for any $\mu\in {\mathbb
C}/{\mathbb Z}$ and logarithmic intertwining operator ${\cal Y}^\mu$
as in (\ref{log:c+n}), the linear map defined by
\[
w_{(1)}\otimes w_{(2)}\mapsto {{\cal Y}^\mu}(w_{(1)},e^{2\pi ip}x)
w_{(2)}= \sum_{\bar n=\mu}\sum_{k\in {\mathbb N}}
{w_{(1)}}_{n;\,k}^{\cal Y}w_{(2)}e^{2\pi ip(-n-1)}x^{-n-1}(\log x+
2\pi ip)^k
\]
is also a logarithmic intertwining operator.  In the strongly graded
case, if the operator in (\ref{log:c+n}) is grading-compatible, then
so is this one.  The right-hand side above can be written as
\begin{eqnarray*}
\lefteqn{e^{-2\pi ip\mu}\sum_{\bar n=\mu}\sum_{k\in {\mathbb
N}}{w_{(1)}}_{n;\,k}^{\cal Y}w_{(2)}x^{-n-1}\sum_{t\in {\mathbb
N}}{k\choose t}(\log x)^{k-t}(2\pi ip)^t}\\
&&=e^{-2\pi ip\mu}\sum_{t\in {\mathbb N}}(2\pi ip)^t\sum_{k\in {\mathbb
N}}{k+t\choose t}\sum_{\bar n=\mu}{w_{(1)}}_{n;\,k+t}^{\cal
Y}w_{(2)}x^{-n-1}(\log x)^k
\end{eqnarray*}
(the coefficient of each power of $x$ being a finite sum over $t$ and
$k$).  We now have:

\begin{propo}
Let $W_1$, $W_2$, $W_3$ be generalized modules for a M\"obius (or
conformal) vertex algebra $V$, and let ${\cal Y}(\cdot, x)\cdot$ be a
logarithmic intertwining operator of type ${W_3\choose W_1\,W_2}$.
For $\mu\in {\mathbb C}/{\mathbb Z}$ and $t \in {\mathbb N}$, define ${\cal
X}^{\mu}_t:W_1\otimes W_2\to W_3[\log x]\{x\}$ by:
\[
{\cal X}^{\mu}_t: w_{(1)}\otimes w_{(2)}\mapsto \sum_{k\in {\mathbb N}}
{k+t\choose t}\sum_{\bar n=\mu}{w_{(1)}}_{n; \,k+t}^{\cal Y}
w_{(2)}x^{-n-1}(\log x)^k.
\]
Then each ${\cal X}^{\mu}_t$ is a logarithmic intertwining operator of
type ${W_3\choose W_1\,W_2}$. In particular, the operator ${\cal X}_t$
defined by
\begin{equation}\label{newio}
{\cal X}_t: w_{(1)}\otimes w_{(2)}\mapsto \sum_{k\in {\mathbb N}}
{k+t\choose t}\sum_{n\in{\mathbb C}}{w_{(1)}}_{n; \,k+t}^{\cal Y}
w_{(2)}x^{-n-1}(\log x)^k
\end{equation}
is a logarithmic intertwining operator of the same type.  In the
strongly graded case, if ${\cal Y}$ is grading-compatible, then so are
${\cal X}^{\mu}_t$ and ${\cal X}_t$.
\end{propo}
\pf {}From the above we see that for any integer $p$, $e^{-2\pi
ip\mu}\sum_{t\in {\mathbb N}}(2\pi ip)^t{\cal X}^{\mu}_t$, and hence
\begin{equation}\label{sumx}
\sum_{t\in {\mathbb N}}(2\pi ip)^t{\cal X}^{\mu}_t,
\end{equation}
is a logarithmic intertwining operator. Let us now prove that ${\cal
X}^{\mu}_m$ is a logarithmic intertwining operator by induction on
$m$. The $m=0$ case follows immediately {}from setting $p=0$ in
(\ref{sumx}). Suppose that ${\cal X}^{\mu}_0,\dots,{\cal
X}^{\mu}_{m-1}$ are all logarithmic intertwining operators. Then for
any integer $p$,
\[
\sum_{t\geq m}(2\pi ip)^t{\cal X}^{\mu}_t=\sum_{t\in {\mathbb N}}(2\pi
ip)^t{\cal X}^{\mu}_t-\sum_{t=0}^{m-1}(2\pi ip)^t{\cal X}^{\mu}_t
\]
is also a logarithmic intertwining operator. Dividing this by $(2\pi
ip)^m$ and then setting $p=0$ we see that ${\cal X}^{\mu}_m$ is also a
logarithmic intertwining operator.  The second assertion follows {}from
Remark \ref{log:mu}, and the last assertion is clear. \epf

\begin{rema}\label{log:fcf}{\rm
Let $W_i$, $W^i$, $i=1,2,3$, be generalized modules for a M\"obius (or
conformal) vertex algebra $V$. If ${\cal Y}(\cdot,x)\cdot$ is a
logarithmic intertwining operator of type ${W_3\choose W_1\,W_2}$ and
$\sigma_1: W^1\to W_1$, $\sigma_2: W^2\to W_2$ and $\sigma_3: W_3\to
W^3$ are $V$-module homomorphisms, then it is easy to see that
$\sigma_3{\cal Y}(\sigma_1\cdot,x)\sigma_2\cdot$ is a logarithmic
intertwining operator of type ${W^3\choose W^1\,W^2}$.  In the
strongly graded case, if ${\cal Y}$ is grading-compatible, then so is
$\sigma_3{\cal Y}(\sigma_1\cdot,x)\sigma_2\cdot$ (recall {}from Remark
\ref{homsaregradingpreserving} that each $\sigma_j$ preserves the
$\tilde A$-grading).  That is, in categorical language, with ${\cal
C}$ a full subcategory of the category of either ${\cal M}$ (the
category of $V$-modules; recall Notation \ref{MGM}), ${\cal GM}$ (the
category of generalized $V$-modules), ${\cal M}_{sg}$ (the category of
strongly graded $V$-modules) or ${\cal GM}_{sg}$ (the category of
strongly graded generalized $V$-modules), the correspondence {}from
${\cal C}\times{\cal C}\times{\cal C}$ to the category ${\bf Vect}$ of
vector spaces given by $(W_1, W_2, W_3)\mapsto {\cal
V}^{W_3}_{W_1\,W_2}$ is functorial in the third slot and cofunctorial
in the first two slots.  }
\end{rema}

\begin{rema}\label{log:newiorm}{\rm
Now recall {}from Remark \ref{set:L(0)s} that $L(0)-L(0)_s$ commutes
with the actions of both $V$ and ${\mathfrak s}{\mathfrak l}(2)$. So
$(L(0)-L(0)_s)^i$ is a $V$-module homomorphism {}from a generalized
module to itself for any $i\in {\mathbb N}$, and this remains true in
the strongly graded case.  Hence by Remark \ref{log:fcf}, given any
logarithmic intertwining operator ${\cal Y}(\cdot,x)\cdot$ as in
(\ref{log:map0}) and any $i,j,k\in {\mathbb N}$,
\begin{equation}\label{newio'}
(L(0)-L(0)_s)^k{\cal Y}((L(0)-L(0)_s)^i\cdot,x) (L(0)-L(0)_s)^j\cdot
\end{equation}
is again a logarithmic intertwining operator, and in the
strongly graded case, if ${\cal Y}$ is grading-compatible, so is this
operator.  In the next remark we will see that the logarithmic
intertwining operators (\ref{newio}) are just linear combinations of
these.  }
\end{rema}

\begin{rema}{\rm
Let $W_1$, $W_2$, $W_3$ and ${\cal Y}$ be as above and let
$w_{(1)}\in{W_1}_{[n_1]}$ and $w_{(2)}\in{W_2}_{[n_2]}$ for some
complex numbers $n_1$ and $n_2$. Fixing $n\in{\mathbb C}$ and using the notation
$T_{p,k,q}$ in (\ref{log:Tprq}) with $a=n_1$ and $b=n_2$, we rewrite
formula (\ref{log:L(0)^1}) in the proof of Lemma \ref{log:lemma} as
\[
(k+1)T_{p,k+1,q}=(L(0)-n_1-n_2+n+1)T_{p,k,q}-T_{p+1,k,q}-T_{p,k,q+1}
\]
for any $p,q,k\in {\mathbb N}$.  {}From this, by (\ref{log:expand}) we see
that for any $t\in {\mathbb N}$,
\[
{k+t\choose t}T_{p,k+t,q}=\sum_{i,j,l\in {\mathbb N},\;i+j+l=t}
\frac{1}{i!j!l!}(-1)^{i+j}(L(0)-n_1-n_2+n+1)^l T_{p+i,k,q+j}.
\]
Setting $p=q=0$ we get
\begin{eqnarray}\label{log:r+t=?}
\lefteqn{{k+t\choose t}{w_{(1)}}^{\cal Y}_{n;\,k+t}w_{(2)}=
\sum_{i,j,l\in {\mathbb N},\;i+j+l=t}\frac{1}{i!j!l!}(-1)^{i+j}\cdot}\nno\\
&&\cdot(L(0)-n_1-n_2+n+1)^l ((L(0)-n_1)^i{w_{(1)}})^{\cal Y}_{n;\,k}
(L(0)-n_2)^jw_{(2)})\nno\\
\end{eqnarray}
for any $t\in {\mathbb N}$.  (Note that this formula gives an
alternate proof of Proposition \ref{log:logwt}(c).)  Now multiplying
by $x^{-n-1}(\log x)^k$ and then summing over $n\in{\mathbb C}$ and
over $k\in{\mathbb N}$ we see that for every $t\in{\mathbb N}$ the
intertwining operator ${\cal X}_t$ in (\ref{newio}) is a linear
combination of intertwining operators of the form (\ref{newio'}).
}
\end{rema}

In preparation for generalizing basic results {}from \cite{FHL} on
intertwining operators to the logarithmic case, we need to generalize
more of the basic tools.  We now define operators ``$x^{\pm L(0)}$''
for generalized modules, in the natural way:

\begin{defi}{\rm
Let $W$ be a generalized module for a M\"obius (or conformal) vertex
algebra. We define $x^{\pm L(0)}: W\to W\{x\}[\log x]\subset W[\log
x]\{x\}$ as follows: For any $w\in W_{[n]}$ ($n\in{\mathbb C}$), define
\begin{equation}\label{log:xpmL}
x^{\pm L(0)}w=x^{\pm n}e^{\pm\log x(L(0)-n)}w
\end{equation}
(note that the local nilpotence of $L(0)-n$ on $W_{[n]}$ insures that
the formal exponential series terminates) and then extend linearly to
all $w\in W$.  (Of course, we could also write
\begin{equation}
x^{\pm L(0)}=x^{\pm L(0)_s}e^{\pm\log x(L(0)-L(0)_s)},
\end{equation}
using the notation $L(0)_s$.)  We also define operators $x^{\pm
L'(0)}$ on $W^*$ by the condition that for all $w'\in W^*$ and $w\in
W$,
\begin{equation}\label{log:xpmL'}
\langle x^{\pm L'(0)}w',w\rangle=\langle w',x^{\pm L(0)}w\rangle\;
(\in {\mathbb C}\{x\}[\log x]),
\end{equation}
so that $x^{\pm L'(0)}: W^*\to W^*\{x\}[[\log x]]$.  }
\end{defi}

\begin{rema}{\rm
Note that these definitions are (of course) compatible with the usual
definitions if $W$ is just an (ordinary) module.  In formula
(\ref{log:xpmL}), $x^{\pm L(0)}$ is defined in a naturally factored
form reminiscent of the factorization invoked in Remark \ref{log:[[]]}
(providing counterexamples there); the symbol $x^{\pm (L(0)-n)}$ is
given meaning by its replacement by $e^{\pm\log x(L(0)-n)}$.  Also
note that in case both $V$ and $W$ are strongly graded, the definition
of $x^{\pm L'(0)}$ given by (\ref{log:xpmL'}), when applied to the
subspace $W'$ of $W^*$, coincides with the definition of $x^{\pm
L'(0)}$ given by (\ref{log:xpmL}) induced {}from the contragredient
module action of $V$ on $W'$. (Recall Theorem \ref{set:W'}.) }
\end{rema}

\begin{rema}{\rm
Note that for $w\in W_{[n]}$, by definition we have
\begin{equation}\label{log:x^L(0)}
x^{\pm L(0)}w=x^{\pm n}\sum_{i\in {\mathbb N}}\frac{(L(0)-n)^iw}{i!}(\pm
\log x)^i\in x^{\pm n}W_{[n]}[\log x].
\end{equation}
It is also handy to have that for any $w\in W$,
\begin{equation}\label{log:inv}
x^{L(0)}x^{-L(0)}w=w=x^{-L(0)}x^{L(0)}w,
\end{equation}
which is clear {}from definition.  Later we will also need the
formula
\begin{equation}\label{log:dx^}
\frac{d}{dx}x^{\pm L(0)}w=\pm x^{-1}x^{\pm L(0)}L(0)w
\end{equation}
for any $w\in W$, i.e.,
\begin{equation}
x\frac{d}{dx}x^{\pm L(0)}w=\pm x^{\pm L(0)}L(0)w,
\end{equation}
or equivalently,
\begin{equation}
\left(x\frac{d}{dx}\mp L(0)\right)x^{\pm L(0)}w=0
\end{equation}
(cf. Lemma \ref{log:de} and Remark \ref{log:[[]]}).  This can be
proved by directly checking that for $w$ homogeneous of generalized
weight $n$,
\[
\frac{d}{dx}e^{\pm\log x(L(0)-n)}w=\pm x^{-1}e^{\pm \log x(L(0)-n)}
(L(0)-n)w,
\]
and hence for such $w$,
\begin{eqnarray*}
\frac{d}{dx}x^{\pm L(0)}w&=&\frac{d}{dx}(x^{\pm n}e^{\pm\log x(L(0)-n)}w)\\
&=&\pm nx^{\pm n-1}e^{\pm\log x(L(0)-n)}w\pm x^{\pm n-1}e^{\pm \log x
(L(0)-n)}(L(0)-n)w\\
&=&\pm x^{\pm n-1}e^{\pm\log x(L(0)-n)}L(0)w=\pm x^{-1}x^{\pm L(0)}L(0)w.
\end{eqnarray*}
}
\end{rema}

In the statement (and proof) of the next result, we shall use
expressions of the type
\begin{eqnarray*}
(1-x)^{L(0)}&=&\sum_{k\in \mathbb{N}}{L(0)\choose k}(-x)^{k}\nn
&=&\sum_{k\in \mathbb{N}}\frac{L(0)(L(0)-1)\cdots
(L(0)-k+1)}{k!}(-x)^{k},
\end{eqnarray*}
which also equals
\[
e^{L(0)\log (1-x)},
\]
as well as expressions involving
\[
(x(1-yx)^{-1})^{n}=\sum_{k\in \mathbb{N}}{-n\choose k}x^{n}(-yx)^{k}
\]
for $n \in \C$ and
\begin{eqnarray*}
\log (x(1-yx)^{-1})&=&\log x+\log (1-yx)^{-1}\nn
&=&\log x+\sum_{k\ge 1}\frac{1}{k}(yx)^{k}.
\end{eqnarray*}

We can now state and prove generalizations to logarithmic intertwining
operators of three standard formulas for (ordinary) intertwining
operators, namely, formulas (5.4.21), (5.4.22) and (5.4.23) of
\cite{FHL}.  The result is (see also \cite{Mi} for Parts (a) and (b)):

\begin{propo}
Let ${\cal Y}$ be a logarithmic intertwining operator of type
${W_3\choose W_1\,W_2}$ and let $w\in W_1$. Then
\begin{description}
\item{(a)}
\begin{equation}\label{log:p1}
e^{yL(-1)}{\cal Y}(w,x)e^{-yL(-1)}={\cal Y}(e^{yL(-1)}w,x)={\cal Y}(w,x+y)
\end{equation}
(recall (\ref{log:not1}))
\item{(b)}
\begin{equation}\label{log:p2}
y^{L(0)}{\cal Y}(w,x)y^{-L(0)}={\cal Y}(y^{L(0)}w,xy)
\end{equation}
(recall (\ref{log:not3}))
\item{(c)}
\begin{equation}\label{log:p3}
e^{yL(1)}{\cal Y}(w,x)e^{-yL(1)}={\cal
Y}(e^{y(1-yx)L(1)}(1-yx)^{-2L(0)}w,x(1-yx)^{-1}).
\end{equation}
\end{description}
\end{propo}
\pf {}From (\ref{log:L(j)b}) with $j=-1$ we see that for any $w\in W_1$,
\[
L(-1){\cal Y}(w,x)={\cal Y}(L(-1)w,x)+{\cal Y}(w,x)L(-1).
\]
This implies
\[
\frac{y^n(L(-1))^n}{n!}{\cal Y}(w,x)=\sum_{i,j\in{\mathbb N},\;i+j=n}{\cal
Y} \bigg(\frac {y^i(L(-1))^i}{i!}w,x\bigg)\frac{y^j(L(-1))^j}{j!}
\]
for any $n\in{\mathbb N}$, where $y$ is a new formal
variable. Summing over $n\in{\mathbb N}$ we see that for any $w\in W_1$,
\[
e^{yL(-1)}{\cal Y}(w,x)={\cal Y}(e^{yL(-1)}w,x)e^{yL(-1)},
\]
and hence
\begin{equation}\label{log:p1-1}
e^{yL(-1)}{\cal Y}(w,x)e^{-yL(-1)}={\cal Y}(e^{yL(-1)}w,x)=e^{y\frac
d{dx}}{\cal Y}(w,x)={\cal Y}(w,x+y),
\end{equation}
where in the last equality we have used (\ref{log:ck1}).  Note that
all the expressions in (\ref{log:p1-1}) remain well defined if we
replace $y$ by any element of $y\C[x][[y]]$.  Thus by Remark
\ref{log:ids-rm}, (\ref{log:p1-1}) still holds if we replace $y$ by
any such element.

For (b), note that for homogeneous $w_{(1)}\in {W_1}$ and $w_{(2)}\in
{W_2}$, by (\ref{log:e^L(0)}) with the formal variable $y$ replaced by
the formal variable $\log y$, we get (recalling Proposition
\ref{log:logwt}(b) and (\ref{log:xpmL}))
\begin{eqnarray*}
\lefteqn{y^{L(0)}({w_{(1)}}_{n;\,k}^{\cal Y}w_{(2)})=}\\
&&\sum_{l\in{\mathbb N}}{k+l\choose k}(y^{L(0)}w_{(1)})_{n;\,k+l}^{\cal
Y}(y^{L(0)}w_{(2)})y^{-n-1}(\log y)^l
\end{eqnarray*}
for any $n\in{\mathbb C}$ and $k\in{\mathbb N}$.
Multiplying this by $x^{-n-1}(\log x)^k$, summing over $n\in
{\mathbb C}$ and $k\in {\mathbb N}$ and using (\ref{log:not3}) we get
\[
y^{L(0)}{\cal Y}(w_{(1)},x)w_{(2)}
={\cal Y}(y^{L(0)}w_{(1)},xy)y^{L(0)}w_{(2)}.
\]
Formula (\ref{log:p2}) then follows {}from (\ref{log:inv}).

Finally we prove (c).  {}From (\ref{log:L(j)b}) with $j=1$ we see that
for any $w\in W_1$,
\[
L(1){\cal Y}(w,x)={\cal Y}((L(1)+2xL(0)+x^2L(-1))w,x)+
{\cal Y}(w,x)L(1).
\]
This implies that
\[
e^{yL(1)}{\cal Y}(w,x)={\cal Y}(e^{y(L(1)+2xL(0)+x^2L(-1))}w,x)
e^{yL(1)},
\]
or
\[
e^{yL(1)}{\cal Y}(w,x)e^{-yL(1)}={\cal Y}
(e^{y(L(1)+2xL(0)+x^2L(-1))}w,x).
\]
Using the identity
\[
e^{y(L(1)+2xL(0)+x^2L(-1))}=e^{yx^2(1-yx)^{-1}L(-1)}
e^{y(1-yx)L(1)}(1-yx)^{-2L(0)},
\]
whose proof is exactly the same as the proof of formula (5.2.41) of
\cite{FHL}, we obtain
\begin{equation}\label{log:p4}
e^{yL(1)}{\cal Y}(w,x)e^{-yL(1)}={\cal
Y}(e^{yx^2(1-yx)^{-1}L(-1)}e^{y(1-yx)L(1)}(1-yx)^{-2L(0)}w,x).
\end{equation}
But by (\ref{log:p1-1}) with $y$ replaced by $yx^2(1-yx)^{-1}$, 
the right-hand side of (\ref{log:p4}) 
is equal to 
\begin{equation}\label{log:p4-r}
{\cal Y}(e^{y(1-yx)L(1)}(1-yx)^{-2L(0)}w,x+yx^2(1-yx)^{-1}).
\end{equation}
Since 
\[
(x+yx^2(1-yx)^{-1})^{n}=(x(1-yx)^{-1})^{n}
\]
for $n \in \C$ and
\[
\log (x+yx^2(1-yx)^{-1})=\log (x(1-yx)^{-1}),
\]
(\ref{log:p4-r}) is equal to the right-hand side of 
(\ref{log:p3}), proving (\ref{log:p3}). 
\epf

\begin{rema}{\rm
The following formula, also a generalization of the corresponding
formula in the ordinary case (see (\ref{xL(0)L(j)})), will be needed:
For $j=-1, 0, 1$,
\begin{equation}\label{log:xLx^}
x^{L(0)}L(j)x^{-L(0)}=x^{-j}L(j).
\end{equation}
To prove this, we first observe that for any $m\in{\mathbb C}$,
$[L(0)-m,L(j)]=-jL(j)$ implies that
\[
e^{\log x(L(0)-m)}L(j)e^{-\log x(L(0)-m)}=e^{-j\log x}L(j).
\]
Hence, for a generalized module element $w$ homogeneous of generalized
weight $n$,
\begin{eqnarray*}
x^{L(0)}L(j)w&=&x^{n-j}e^{\log x(L(0)-n+j)}L(j)w\\
&=&x^{n-j}e^{-j\log x}L(j)e^{\log x(L(0)-n+j)}w\\
&=&x^{n-j}L(j)e^{\log x(L(0)-n)}w\\
&=&x^{-j}L(j)x^{L(0)}w,
\end{eqnarray*}
and (\ref{log:xLx^}) then follows immediately {}from (\ref{log:inv}).
}
\end{rema}

\begin{rema}{\rm {}From (\ref{log:xLx^}) we see that
\begin{equation}\label{xe^Lx}
x^{L(0)}e^{yL(j)}x^{-L(0)}=e^{yx^{-j}L(j)}.
\end{equation}
}
\end{rema}

For an ordinary module $W$ for a vertex operator algebra and any
$a\in{\mathbb C}$, the operator $e^{aL(0)}$ on $W$ is defined
by
\begin{equation}\label{eaL0ordinary}
e^{aL(0)}w=e^{ah}w
\end{equation}
for any homogeneous $w\in W_{(h)}$, $h\in{\mathbb C}$ and then by linear
extension to any $w\in W$.  More generally, for a generalized module
$W$ for a M\"obius (or conformal) vertex algebra and any $a\in{\mathbb
C}$, we define the operator $e^{aL(0)}$ on $W$ by
\begin{equation}\label{eaL0}
e^{aL(0)}w=e^{ah}e^{a(L(0)-h)}w
\end{equation}
for any homogeneous $w\in W_{[h]}$, $h\in{\mathbb C}$ and then by
linear extension to all $w\in W$. (Note that for a formal variable
$x$, we already have $e^{xL(0)}w=e^{hx}e^{x(L(0)-h)}w$.)  {}From the
definition,
\begin{equation}
e^{aL(0)}e^{-aL(0)}w=w.
\end{equation}
Recalling Remark
\ref{set:L(0)s} for the notation $L(0)_s$, we see that
\begin{equation}
e^{aL(0)}=e^{aL(0)_s}e^{a(L(0)-L(0)_s)}\;\;\mbox{on}\;W,
\end{equation}
where the exponential series $e^{a(L(0)-L(0)_s)}$ terminates on each
element of $W$.

\begin{rema}\label{analyticallyconvervent}{\rm
The operators defined in (\ref{eaL0ordinary}) and (\ref{eaL0}) can be
alternatively defined or viewed as the (analytically) convergent sums
of the indicated exponential series of operators; these operators act
on the (finite-dimensional) subspaces of $W$ generated by the repeated
action of $L(0)$ on homogeneous vectors $w \in W$.
}
\end{rema}

\begin{rema}\label{exponentialaVhom}{\rm
The operator $e^{a(L(0)-L(0)_s)}$ on $W$ is a $V$-homomorphism, in
view of Remark \ref{set:L(0)s} (cf.\ Remark \ref{log:newiorm}).  Let
$r$ be an integer.  Then $e^{2\pi irL(0)_s}$ is also a
$V$-homomorphism, by Remark \ref{congruent}.  Thus for $r \in \Z$,
$e^{2\pi irL(0)}$ is a $V$-homomorphism.  In the strongly graded
case, all of these $V$-homomorphisms are grading-preserving.}
\end{rema}

\begin{rema}{\rm
We now recall some identities about the action of ${\mathfrak s}{\mathfrak l}
(2)$ on any of its modules. For convenience we put them in the
following form:
\begin{equation}\label{log:SL2-1}
e^{xL(-1)}\left(\begin{array}{c}L(-1)\\L(0)\\L(1)\end{array}\right)
e^{-xL(-1)}=\left(\begin{array}{ccc}1&0&0\\-x&1&0\\x^2&-2x&1\end{array}
\right)\left(\begin{array}{c}L(-1)\\L(0)\\L(1)\end{array}\right)
\end{equation}
\begin{equation}\label{log:SL2-2}
e^{xL(0)}\left(\begin{array}{c}L(-1)\\L(0)\\L(1)\end{array}\right)
e^{-xL(0)}=\left(\begin{array}{ccc}e^x&0&0\\0&1&0\\0&0&e^{-x}\end{array}
\right)\left(\begin{array}{c}L(-1)\\L(0)\\L(1)\end{array}\right)
\end{equation}
\begin{equation}\label{log:SL2-3}
e^{xL(1)}\left(\begin{array}{c}L(-1)\\L(0)\\L(1)\end{array}\right)
e^{-xL(1)}=\left(\begin{array}{ccc}1&2x&x^2\\0&1&x\\0&0&1\end{array}
\right)\left(\begin{array}{c}L(-1)\\L(0)\\L(1)\end{array}\right).
\end{equation}
Formula (\ref{log:SL2-2}) follows {}from (5.2.12) and (5.2.13) in
\cite{FHL}; Formula (\ref{log:SL2-3}) follows {}from (5.2.14) in
\cite{FHL} and $[L(1),L(-1)]=2L(0)$; and formula (\ref{log:SL2-1})
follows {}from (\ref{log:SL2-3}) and the fact that
\[
L(-1)\mapsto L(1),\; L(0)\mapsto-L(0),\; L(1)\mapsto L(-1)
\]
is a Lie algebra automorphism of ${\mathfrak s}{\mathfrak l}(2)$.
}
\end{rema}

\begin{rema}\label{log:Lj2rema}{\rm
It is convenient to note that the ${\mathfrak s}{\mathfrak l}(2)$-bracket
relations (\ref{log:L(j)b}) are equivalent to
\begin{equation}\label{log:L(j)b2}
{\cal Y}(L(j)w_{(1)},x)=\sum_{i=0}^{j+1}{j+1\choose i}(-x)^i
[L(j-i),{\cal Y}(w_{(1)},x)]
\end{equation}
for $j=-1$, $0$ and $1$. This can be easily checked by writing
(\ref{log:L(j)b}) as
\[
\left(\begin{array}{c}[L(-1),{\cal Y}(w_{(1)},x)]\\{}[L(0),{\cal Y}
(w_{(1)},x)]\\{}[L(1),{\cal Y}(w_{(1)},x)]\end{array}\right)
=\left(\begin{array}{ccc}1&0&0\\x&1&0\\x^2&2x&1\end{array}\right)
\left(\begin{array}{c}{\cal Y}(L(-1)w_{(1)},x)\\{\cal Y}
(L(0)w_{(1)},x)\\{\cal Y}(L(1)w_{(1)},x)\end{array}\right)
\]
and then multiplying this by the inverse of the invertible matrix on
the right-hand side, obtained by replacing $x$ by $-x$.  (Of course,
in the case where $V$ is conformal, this equivalence is already
encoded in the symmetry of the Jacobi identity.)  }
\end{rema}

We have already defined the natural process of multiplying the formal
variable in a logarithmic intertwining operator by $e^{2\pi ip}$ for
$p \in \Z$ (recall (\ref{substitutionofe2piipx})), and this process
yields another logarithmic intertwining operator (recall Remark
(\ref{formalinvariance})).  It is natural to generalize this
substitution process to that of multiplying the formal variable in a
logarithmic intertwining operator by the exponential $e^{\zeta}$ of
any complex number $\zeta$.  As in the special case $\zeta = 2\pi ip$,
the process will depend on $\zeta$, not just on $e^{\zeta}$, but we
will still find it convenient to use the shorthand symbol $e^{\zeta}$
in our notation for the process.  In Section 7 of \cite{tensor2}, we
introduced this procedure in the case of ordinary (nonlogarithmic)
intertwining operators, and we now carry it out in the general
logarithmic case.  We are about to use this substitution mostly for
$\zeta = (2r+1)\pi i$, $r \in \Z$.

Let $(W_1,Y_1)$, $(W_2,Y_2)$ and $(W_3,Y_3)$ be generalized modules
for a M\"obius (or conformal) vertex algebra $V$.  Let ${\cal Y}$ be a
logarithmic intertwining operator of type ${W_3\choose W_1\, W_2}$.
For any complex number $\zeta$ and any $w_{(1)}\in W_1$ and
$w_{(2)}\in W_2$, set
\begin{equation}\label{log:subs}
{\cal Y} (w_{(1)},e^\zeta x)w_{(2)} = {\cal
Y}(w_{(1)},y)w_{(2)}\lbar_{y^n=e^{\zeta n}x^n,\; (\log y)^k=(\zeta+\log
x)^k,\;n\in{\mathbb C},\;k\in{\mathbb N}},
\end{equation}
a well-defined element of $W_3[\log x]\{x\}$.  Note that this element
indeed depends on $\zeta$, not just on $e^\zeta$.

\begin{rema}{\rm
In Section 4 below we will take the further step of specializing the
formal variable $x$ to $1$ (or equivalently, the formal variable $y$
to $e^{\zeta}$) in (\ref{log:subs}); that is, we will consider ${\cal
Y} (w_{(1)},e^\zeta)w_{(2)}$.
}
\end{rema}

Given any $r\in {\mathbb Z}$, we define
\[
\Omega_r({\cal Y}): W_2\otimes W_1\to W_3[\log x]\{x\}
\]
by the formula
\begin{equation}\label{Omega_r}
\Omega_{r}({\cal Y})(w_{(2)},x)w_{(1)} = e^{xL(-1)} {\cal
Y}(w_{(1)},e^{ (2r+1)\pi i}x)w_{(2)}
\end{equation}
for $w_{(1)}\in W_{1}$ and $w_{(2)}\in W_{2}$.  This expression is
indeed well defined because of the truncation condition
(\ref{log:ltc}) (recall Remark \ref{ordinaryandlogintwops}).  The
following result generalizes Proposition 7.1 of \cite{tensor2} (for
the ordinary intertwining operator case) and has essentially the same
proof as that proposition, which in turn generalized Proposition 5.4.7
of \cite{FHL}:

\begin{propo}\label{log:omega}
The operator $\Omega_r({\cal Y})$ is a logarithmic intertwining
operator of type ${W_3\choose W_2\, W_1}$. Moreover,
\begin{equation}\label{log:or}
\Omega_{-r-1}(\Omega_r({\cal Y}))=\Omega_r(\Omega_{-r-1}({\cal Y}))
={\cal Y}.
\end{equation}
In the strongly graded case, if ${\cal Y}$ is grading-compatible, then
so is $\Omega_r({\cal Y})$, and in particular, the correspondence
${\cal Y}\mapsto \Omega_r({\cal Y})$ defines a linear isomorphism {}from
${\cal V}_{W_1\,W_2}^{W_3}$ to ${\cal V}_{W_2\,W_1}^{W_3}$, and we
have
\[
 N_{W_1\,W_2}^{W_3}=N_{W_2\,W_1}^{W_3}.
\]
\end{propo}
\pf The lower truncation condition (\ref{log:ltc}) is clear. {From}
the Jacobi identity (\ref{log:jacobi}) for ${\cal Y}$,
\begin{eqnarray}
\lefteqn{\dps x^{-1}_0\delta \left( {x_1-y\over x_0}\right)
Y_3(v,x_1){\cal Y}(w_{(1)},y)w_{(2)}}\nno\\
&&\hspace{2em}- x^{-1}_0\delta \left( {y-x_1\over -x_0}\right)
{\cal Y}(w_{(1)},y)Y_2(v,x_1)w_{(2)}\nno \\
&&{\dps = y^{-1}\delta \left( {x_1-x_0\over y}\right)
{\cal Y}(Y_1(v,x_0)w_{(1)},y)
w_{(2)}},
\end{eqnarray}
with $v \in V$, $w_{(1)} \in W_1$ and $w_{(2)} \in W_2$, we obtain
\begin{eqnarray}\label{75}
\dps x^{-1}_0\delta \left( {x_1-y\over x_0}\right)
e^{-yL(-1)}Y_3(v,x_1){\cal Y}(w_{(1)},y)w_{(2)}
\lbar_{y^n=e^{(2r+1)\pi in}x_2^n,\; (\log y)^k=((2r+1)\pi i+\log
x_2)^k,\;n\in{\mathbb C},\;k\in{\mathbb N}}
\hspace{1.5em}\nno\\
- x^{-1}_0\delta \left( {y-x_1\over -x_0}\right)
e^{-yL(-1)}{\cal Y}(w_{(1)},y)Y_2(v,x_1)w_{(2)}
\lbar_{y^n=e^{(2r+1)\pi in}x_2^n,\; (\log y)^k=((2r+1)\pi i+\log
x_2)^k,\;n\in{\mathbb C},\;k\in{\mathbb N}}\nno \\
{\dps = y^{-1}\delta \left( {x_1-x_0\over y}\right)
e^{-yL(-1)}{\cal Y}(Y_1(v,x_0)w_{(1)},y)
w_{(2)}\lbar_{y^n=e^{(2r+1)\pi in}x_2^n,\; (\log y)^k=((2r+1)\pi i+\log
x_2)^k,\;n\in{\mathbb C},\;k\in{\mathbb N}}}.
\end{eqnarray}
The first term of the left-hand side of (\ref{75}) is equal to
\begin{eqnarray*}
&{\dps x^{-1}_0\delta \left( {x_1-y\over x_0}\right)
Y_3(v,x_1-y)e^{-yL(-1)}{\cal Y}(w_{(1)},y)w_{(2)}
\lbar_{y^n=e^{(2r+1)\pi in}x_2^n,\; (\log y)^k=((2r+1)\pi i+\log
x_2)^k,\;n\in{\mathbb C},\;k\in{\mathbb N}}}&\nno\\
&{\dps =x^{-1}_0\delta \left( {x_1+x_{2}\over x_0}\right)
Y_3(v,x_0)\Omega_{r}({\cal Y})(w_{(2)},x_{2})w_{(1)},}&
\end{eqnarray*}
the second term  is equal to
$$
- x^{-1}_0\delta \left( {-x_{2}-x_1\over -x_0}\right)
\Omega_{r}({\cal Y})(Y_2(v,x_1)w_{(2)},x_{2})w_{(1)}
$$
and the right-hand side of (\ref{75}) is equal to
$$
-x_{2}^{-1}\delta \left( {x_1-x_0\over -x_{2}}\right)
\Omega_{r}({\cal Y})(w_{(2)},x_{2})Y_1(v,x_0)w_{(1)}.
$$
Substituting  into (\ref{75}) we obtain
\begin{eqnarray}
\lefteqn{x^{-1}_0\delta \left( {x_1+x_{2}\over x_0}\right)
Y_3(v,x_0)\Omega_{r}({\cal Y})(w_{(2)},x_{2})w_{(1)}}\nno\\
&&\hspace{2em}- x^{-1}_0\delta \left( {x_{2}+x_1\over x_0}\right)
\Omega_{r}({\cal Y})(Y_2(v,x_1)w_{(2)},x_{2})w_{(1)}\nno\\
&&=-x_{2}^{-1}\delta \left( {x_1-x_0\over -x_{2}}\right)
\Omega_{r}({\cal Y})(w_{(2)},x_{2})Y_1(v,x_0)w_{(1)},
\end{eqnarray}
which is equivalent to
\begin{eqnarray}
\lefteqn{x^{-1}_1\delta \left( {x_0-x_{2}\over x_1}\right)
Y_3(v,x_0)\Omega_{r}({\cal Y})(w_{(2)},x_{2})w_{(1)}}\nno\\
&&\hspace{2em}-x_{1}^{-1}\delta \left( {x_2-x_0\over -x_{1}}\right)
\Omega_{r}({\cal Y})(w_{(2)},x_{2})Y_1(v,x_0)w_{(1)}\nno\\
&&= x^{-1}_2\delta \left( {x_{0}-x_1\over x_2}\right)
\Omega_{r}({\cal Y})(Y_2(v,x_1)w_{(2)},x_{2})w_{(1)}
\end{eqnarray}
(recall (\ref{2termdeltarelation})).  This in turn is the Jacobi
identity for $\Omega_{r}({\cal Y})$ (with the roles of $x_0$ and $x_1$
reversed in (\ref{log:jacobi})).

To prove the $L(-1)$-derivative property (\ref{log:L(-1)dev}) for
$\Omega_{r}({\cal Y})$, first note that {}from (\ref{log:subs}) and the
$L(-1)$-derivative property for ${\cal Y}$,
\begin{eqnarray}
\frac{d}{dx}{\cal Y} (w_{(1)},e^{\zeta} x)w_{(2)}
&=& e^{\zeta}\left(\frac{d}{dy} {\cal
Y}(w_{(1)},y)w_{(2)}\right)\lbar_{y^n=e^{\zeta n}x^n,\; (\log
y)^k=(\zeta+\log x)^k,\;n\in{\mathbb C},\;k\in{\mathbb N}}\nno\\
&=& e^{\zeta}{\cal Y} (L(-1)w_{(1)},e^\zeta x)w_{(2)},
\end{eqnarray}
and in particular,
\begin{equation}
\frac{d}{dx}{\cal Y} (w_{(1)},e^{(2r+1)\pi i} x)w_{(2)}
= -{\cal Y} (L(-1)w_{(1)},e^{(2r+1)\pi i} x)w_{(2)}.
\end{equation}
Thus, by using formula (\ref{log:L(j)b}) with $j=-1$ we have
\begin{eqnarray}
\lefteqn{\frac{d}{dx}\Omega_{r}({\cal Y})(w_{(2)},x)w_{(1)}=
\frac{d}{dx}e^{xL(-1)}{\cal Y}(w_{(1)},e^{ (2r+1)\pi i}x)w_{(2)}}\nno\\
&&=e^{xL(-1)}L(-1){\cal Y}(w_{(1)},e^{ (2r+1)\pi i}x)w_{(2)}+
e^{xL(-1)}\frac{d}{dx}{\cal Y}(w_{(1)},e^{ (2r+1)\pi i}x)w_{(2)}\nno\\
&&=e^{xL(-1)}L(-1){\cal Y}(w_{(1)},e^{ (2r+1)\pi i}x)w_{(2)}-
e^{xL(-1)}{\cal Y}(L(-1)w_{(1)},e^{ (2r+1)\pi i}x)w_{(2)}\nno\\
&&=e^{xL(-1)}{\cal Y}(w_{(1)},e^{ (2r+1)\pi i}x)L(-1)w_{(2)}\nno\\
&&=\Omega_{r}({\cal Y})(L(-1)w_{(2)},x)w_{(1)},
\end{eqnarray}
as desired.

In the case that $V$ is M\"obius, we prove the ${\mathfrak
s}{\mathfrak l}(2)$-bracket relations (\ref{log:L(j)b}) for
$\Omega_r({\cal Y})$. By using the ${\mathfrak s}{\mathfrak
l}(2)$-bracket relations for ${\cal Y}$ and the relations
\[
e^{xL(-1)}L(j)e^{-xL(-1)}=\sum_{i=0}^{j+1}{j+1\choose i}(-x)^iL(j-i)
\]
for $j=-1, 0$ and $-1$ (see (\ref{log:SL2-1})), we have
\begin{eqnarray*}
\lefteqn{\Omega_r({\cal Y})(L(j)w_{(2)},x)w_{(1)}=
e^{xL(-1)}{\cal Y}(w_{(1)},e^{(2r+1)\pi i}x)L(j)w_{(2)}=}\\
&&=e^{xL(-1)}L(j){\cal Y}(w_{(1)},e^{(2r+1)\pi i}x)w_{(2)}\\
&&\hspace{4em}-e^{xL(-1)}\sum_{i=0}^{j+1}{j+1\choose i}(-x)^i{\cal Y}
(L(j-i)w_{(1)},e^{(2r+1)\pi i}x)w_{(2)}\\
&&=\sum_{i=0}^{j+1}{j+1\choose i}(-x)^iL(j-i)e^{xL(-1)}{\cal
Y}(w_{(1)},e^{(2r+1)\pi i}x)w_{(2)}\\
&&\hspace{4em}-e^{xL(-1)}\sum_{i=0}^{j+1}{j+1\choose i}(-x)^i{\cal Y}
(L(j-i)w_{(1)},e^{(2r+1)\pi i}x)w_{(2)}\\
&&=\sum_{i=0}^{j+1}{j+1\choose i}(-x)^i(L(j-i)\Omega_r({\cal
Y})(w_{(2)},x)-\Omega_r({\cal Y})(w_{(2)},x)L(j-i))w_{(1)},
\end{eqnarray*}
which is the alternative form (\ref{log:L(j)b2}) of the ${\mathfrak
s}{\mathfrak l}(2)$-bracket relations for $\Omega_r({\cal Y})$.

The identity (\ref{log:or}) is clear {from} the definitions of
$\Omega_{r}({\cal Y})$ and $\Omega_{-r-1}({\cal Y})$, and the
remaining assertions are clear. \epf

\begin{rema}\label{Ys1s2s3}{\rm
For each triple $s_1, s_2, s_3\in{\mathbb Z}$, the
logarithmic intertwining operator ${\cal Y}$ gives rise to a
logarithmic intertwining
operator ${\cal Y}_{[s_1,s_2, s_3]}$ of the same type, defined by
\[
{\cal Y}_{[s_1,s_2, s_3]}(w_{(1)},x)=e^{2\pi is_1 L(0)}{\cal
Y}(e^{2\pi is_2 L(0)} w_{(1)},x)e^{2\pi is_3 L(0)}
\]
for $w_{(1)}\in W_1$, by Remarks \ref{log:fcf} and
\ref{exponentialaVhom}.  In the strongly graded case, if ${\cal Y}$ is
grading-compatible, so is ${\cal Y}_{[s_1,s_2, s_3]}$.  Clearly,
\[
{\cal Y}_{[0,0,0]}={\cal Y}
\]
and for $r_1,r_2,r_3,s_1,s_2,s_3\in{\mathbb Z}$,
\[
({\cal Y}_{[r_1,r_2, r_3]})_{[s_1,s_2, s_3]}={\cal Y}_{[r_1+s_1,
r_2+s_2, r_3+s_3]}.
\]
For any $a\in{\mathbb C}$, we have the formula
\begin{equation}\label{710}
e^{aL(0)}{\cal Y}(w_{(1)},x)e^{-aL(0)}={\cal Y}
(e^{aL(0)}w_{(1)},e^ax)
\end{equation}
(cf.\ (\ref{log:p2})).  This is proved by imitating the proof of
(\ref{log:p2}), replacing $y^{L(0)}$ by $e^{aL(0)}$, $y$ by $e^a$ and
$\log y$ by $a$ in that proof, using (\ref{eaL0}) in place of
(\ref{log:xpmL}) and keeping in mind formula (\ref{log:subs}).  (When
(\ref{log:e^L(0)}) is used in this proof, for homogeneous elements
$w_{(1)}$ and $w_{(2)}$, the exponential series all terminate, as does
the sum over $l \in {\mathbb{N}}$.)  {}From this, we see that
(\ref{log:or}) generalizes to
\[
\Omega_s(\Omega_r({\cal Y}))={\cal Y}_{[r+s+1,-(r+s+1),-(r+s+1)]}
={\cal Y}(\cdot,e^{2\pi i (r+s+1)}\cdot)
\]
for all $r,s \in \Z$.
}
\end{rema}

In case $V$, $W_1$, $W_2$ and $W_3$ are strongly graded, which we now
assume, we have the concept of ``$r$-contragredient operator'' as
follows (in the ordinary intertwining operator case this was
introduced in \cite{tensor2}): Given a grading-compatible logarithmic
intertwining operator ${\cal Y}$ of type ${W_3\choose W_1\,W_2}$ and
an integer $r$, we define the {\em $r$-contragredient operator of
${\cal Y}$} to be the linear map
\begin{eqnarray*}
W_1\otimes W'_3&\to&W'_2\{x\}[[\log x]]\\
w_{(1)}\otimes w'_{(3)}&\mapsto&A_r({\cal Y})(w_{(1)},x)w'_{(3)}
\end{eqnarray*}
given by
\begin{eqnarray}\label{log:Ardef}
\lefteqn{\langle A_r({\cal Y})(w_{(1)},x)w'_{(3)},w_{(2)}
\rangle_{W_2}=}\nno\\
&&=\langle w'_{(3)},{\cal Y}(e^{xL(1)}e^{(2r+1)\pi iL(0)}(x^{-L(0)})^2
w_{(1)},x^{-1})w_{(2)}\rangle_{W_3},
\end{eqnarray}
for any $w_{(1)}\in W_1$, $w_{(2)}\in W_2$ and $w'_{(3)}\in W'_3$,
where we use the notation
\[
f(x^{-1})=\sum_{m\in{\mathbb N},\,n\in{\mathbb C}}w_{n,m} x^{-n}(-\log x)^m
\]
for any $f(x)=\sum_{m\in{\mathbb N},\,n\in {\mathbb C}}w_{n,m}
x^n(\log x)^m\in{\cal W}\{x\}[[\log x]]$, ${\cal W}$ any vector space
(not involving $x$).  Note that for the case $W_1=V$ and $W_2=W_3=W$,
the operator $A_r({\cal Y})$ agrees with the contragredient vertex
operator ${\cal Y}'$ (recall (\ref{yo}) and (\ref{y'})) for any
$r\in{\mathbb Z}$.

We have the following result generalizing Proposition 7.3 in
\cite{tensor2} for ordinary intertwining operators, and having
essentially the same proof (and also generalizing Theorem 5.5.1 and
Proposition 5.5.2 of \cite{FHL}):

\begin{propo}\label{log:A}
The $r$-contragredient operator $A_r({\cal Y})$ of a
grading-compatible logarithmic intertwining operator ${\cal Y}$ of
type ${W_3\choose W_1\,W_2}$ is a grading-compatible logarithmic
intertwining operator of type ${W'_2\choose W_1\,W'_3}$. Moreover,
\begin{equation}\label{log:ar}
A_{-r-1}(A_r({\cal Y}))=A_r(A_{-r-1}({\cal Y}))={\cal Y}.
\end{equation}
In particular, the correspondence ${\cal Y}\mapsto A_r({\cal Y})$
defines a linear isomorphism {}from ${\cal V}_{W_1\,W_2}^{W_3}$ to
${\cal V}_{W_1\,W'_3}^{W'_2}$, and we have
\[
N_{W_1\,W_2}^{W_3}=N_{W_1\,W'_3}^{W'_2}.
\]
\end{propo}
\pf First we need to show that for $w_{(1)}\in W_1$ and $w'_{(3)}\in
W'_3$,
\begin{equation}\label{finitelymanypowersoflogx}
A_r({\cal Y})(w_{(1)},x)w'_{(3)} \in W'_2[\log x]\{x\},
\end{equation}
that is, for each power of $x$ there are only finitely many powers of
$\log x$.  This and the lower truncation condition (\ref{log:ltc}) as
well as the grading-compatibility condition (\ref{gradingcompatcondn})
follow {}from a variant of the argument proving the lower truncation
condition for contragredient vertex operators (recall
(\ref{truncationforY'})):

Fix elements $w_{(1)}\in W_1$ and $w'_{(3)}\in W'_3$ homogeneous with
respect to the double gradings of $W_1$ and $W'_3$ (it will suffice to
prove the desired assertions for such elements), and in fact take
\[
w_{(1)}\in W_1^{(\beta)} \;\;\mbox{and}\;\; w'_{(3)}\in
(W'_3)^{(\gamma)},
\]
where $\beta$ and $\gamma$ are elements of the abelian group $\tilde
A$, in the notation of Definition \ref{def:dgw}, and fix $n \in \C$.
The right-hand side of (\ref{log:Ardef}) is a (finite) sum of terms of
the form
\begin{equation}\label{w3Yw2}
\langle w'_{(3)},{\cal Y}(w,x^{-1})w_{(2)}\rangle x^p (\log x)^q
\end{equation}
where $w \in W_1$ is doubly homogeneous and in fact $w \in
W_1^{(\beta)}$ (by (\ref{m-L(n)-A})), and where $p \in \C$ and $q \in
{\mathbb N}$.  (The pairing $\langle \cdot,\cdot \rangle$ is between
$W'_3$ and $W_3$.)  Let $w_{(2)}^{(-\beta-\gamma)}$ be the component
of (the arbitrary element) $w_{(2)}$ in $W_2^{(-\beta-\gamma)}$, with
respect to the $\tilde A$-grading.  Then (\ref{w3Yw2}) equals
\begin{equation}\label{w3Yw2betagamma}
\langle w'_{(3)},{\cal Y}(w,x^{-1})w_{(2)}^{-\beta-\gamma}\rangle x^p
(\log x)^q
\end{equation}
because of the grading-compatibility condition
(\ref{gradingcompatcondn}) together with (\ref{W'beta}).  (This is why
we need our logarithmic intertwining operators to be
grading-compatible.)  This shows in particular that
\begin{equation}
{w_{(1)}}_{N;\,K}^{\cal Y}w'_{(3)} \in (W'_{2})^{(\beta + \gamma)}
\end{equation}
for $N \in \C$ and $K \in {\mathbb N}$, so that
(\ref{gradingcompatcondn}) holds for $A_r({\cal Y})$.  Let us write
(\ref{w3Yw2betagamma}) as
\begin{equation}\label{w3Yw2betagammaexpanded}
\sum_{l\in{\C}}\sum_{k\in {\mathbb N}}\langle w'_{(3)},
w_{l;\,k}^{\cal Y}w_{(2)}^{-\beta-\gamma}\rangle x^{l+1+p}(\log
x)^{k+q} =\sum_{m\in{\C}}\sum_{k\in {\mathbb N}}\langle w'_{(3)},
w_{n-p-1-m;\,k}^{\cal Y}w_{(2)}^{-\beta-\gamma}\rangle x^{n-m}(\log
x)^{k+q}
\end{equation}
(recall that we have fixed $n \in \C$).  But each term $\langle
w'_{(3)}, w_{n-p-1-m;\,k}^{\cal Y}w_{(2)}^{-\beta-\gamma}\rangle$ in
(\ref{w3Yw2betagammaexpanded}) can be replaced by
\begin{equation}
\langle w'_{(3)}, w_{n-p-1-m;\,k}^{\cal Y}u^{[m]}\rangle,
\end{equation}
where $u^{[m]} \in W_3^{(-\beta-\gamma)}$ is the component of
$w_{(2)}^{-\beta-\gamma}$, with respect to the generalized-weight
grading, of (generalized) weight
\[
\wt u^{[m]} = \wt w'_{(3)} - \wt w + n - p - m,
\]
by Proposition \ref{log:logwt}(b).  To see that the coefficient of
$x^n$ in (\ref{w3Yw2betagammaexpanded}) involves only finitely many
powers of $\log x$, independently of the element $w_{(2)}$, we take
$m=0$ in (\ref{w3Yw2betagammaexpanded}) and we observe that the
possible elements $u^{[0]}$ range through the space
\[
(W_3)^{(-\beta -\gamma)}_{[\swt w'_{(3)}-\swt w + n - p]},
\]
which is finite dimensional by the grading restriction condition
(\ref{set:dmfin}).  This proves (\ref{finitelymanypowersoflogx}).  To
prove the lower truncation condition (\ref{log:ltc}), what we must
show is that for sufficiently large $m \in {\mathbb N}$, the
coefficient of $x^{n-m}$ in (\ref{w3Yw2betagammaexpanded}) is $0$
(independently of $w_{(2)}$).  But by the grading restriction
condition (\ref{set:dmltc}),
\[
(W_3)^{(-\beta -\gamma)}_{[\swt w'_{(3)}-\swt w + n - p - m]} = 0
\;\;\mbox{ for }\;m\in {\mathbb N}\;\mbox{ sufficiently large.}
\]
Hence the coefficient of $x^{n-m}$ in (\ref{w3Yw2betagammaexpanded})
is zero for $m\in {\mathbb N}$ sufficiently large, as desired, proving
the lower truncation condition.

For the Jacobi identity, we need to show that
\begin{eqnarray}\label{716}
\lefteqn{\left\langle x^{-1}_0\delta \left( {x_1-x_2\over x_0}\right)
Y_{2}'(v,x_1)A_{r}({\cal
Y})(w_{(1)},x_2)w'_{(3)},w_{(2)}\right\rangle_{W_2}}
\nno \\
&&\;\;\;\;- \left\langle x^{-1}_0\delta \left( {x_2-x_1\over -x_0}\right)
A_{r}({\cal
Y})(w_{(1)},x_2)Y_{3}'(v,x_1)w'_{(3)},w_{(2)}\right\rangle_{W_2}
\nno\\
&&=\left\langle x^{-1}_2\delta \left( {x_1-x_0\over x_2}\right)
A_{r}({\cal Y})(Y_{1}(v,x_0)w_{(1)},x_2)w'_{(3)},w_{(2)}\right
\rangle_{W_2}.
\end{eqnarray}
By the definitions (\ref{y'}) and (\ref{log:Ardef}) we have
\begin{eqnarray}\label{717}
\lefteqn{\langle Y_{2}'(v,x_1)A_{r}({\cal Y})
(w_{(1)},x_2)w'_{(3)},w_{(2)}\rangle_{W_2}}\nno  \\
&&= \langle w'_{(3)},{\cal Y}(e^{x_2L(1)}e^{(2r+1)\pi iL(0)}
(x_2^{-L(0)})^2 w_{(1)},x^{-1}_2)\cdot \nno\\
&&\hspace{8em}\cdot Y_{2}(e^{x_1L(1)}(-x^{-2}_1)^{L(0)}v,x^{-1}_1)w_{(2)}
\rangle_{W_3},
\end{eqnarray}
\begin{eqnarray}\label{718}
\lefteqn{\langle A_{r}({\cal Y})(w_{(1)},x_2)
Y_{3}'(v,x_1)w'_{(3)},w_{(2)}\rangle_{W_2}}\nno  \\
&&=\langle w'_{(3)},Y_{3}(e^{x_1L(1)}(-x^{-2}_1)^{L(0)}v,x^{-1}_1)\cdot \nno\\
&&\hspace{6em}\cdot {\cal Y}(e^{x_2L(1)}e^{(2r+1)\pi iL(0)}
(x_2^{-L(0)})^2 w_{(1)},x^{-1}_2)w_{(2)}\rangle_{W_3},
\end{eqnarray}
\begin{eqnarray}\label{719}
\lefteqn{\langle A_{r}({\cal Y})
(Y_{1}(v,x_0)w_{(1)},x_2)w'_{(3)},w_{(2)}\rangle_{W_2}}\nno \\
&&= \langle w'_{(3)},{\cal Y}(e^{x_2L(1)}e^{(2r+1)\pi iL(0)}
(x_2^{-L(0)})^2
Y_{1}(v,x_0)w_{(1)},x^{-1}_2)w_{(2)}\rangle_{W_3}.
\end{eqnarray}
{From} the Jacobi identity for  ${\cal Y}$  we have
\begin{eqnarray}\label{720}
\lefteqn{\Biggl\langle w'_{(3)},\left( {-x_0\over x_1x_2}\right) ^{-1}\delta
 \left( {x^{-1}_1-x^{-1}_2\over -x_0/x_1x_2} \right)
Y_{3}(e^{x_1L(1)}(-x^{-2}_1)^{L(0)}v,x^{-1}_1)\cdot}\nno \\
&&\;\;\;\;\;\;\;\;\;\;\cdot {\cal Y}(e^{x_2L(1)}e^{(2r+1)\pi iL(0)}
(x_2^{-L(0)})^2 w_{(1)},x^{-1}_2)w_{(2)}\Biggr\rangle_{W_3}\nno \\
&&- \Biggl\langle w'_{(3)},\left( {-x_0\over x_1x_2}\right) ^{-1}\delta
\left( {x^{-1}_2-x^{-1}_1\over x_0/x_1x_2}\right)
{\cal Y}(e^{x_2L(1)}e^{(2r+1)\pi iL(0)}
(x_2^{-L(0)})^2 w_{(1)},x^{-1}_2)\cdot \nno\\
&&\;\;\;\;\;\;\;\;\;\;\cdot
 Y_{2}(e^{x_1L(1)}(-x^{-2}_1)^{L(0)}v,x^{-1}_1)w_{(2)}\Biggr
\rangle_{W_3} \nno\\
&&= \Biggl\langle w'_{(3)},(x^{-1}_2)^{-1}
\delta \left( {x^{-1}_1+x_0/x_1x_2\over x^{-1}_2}\right)
{\cal Y}(Y_{1}(e^{x_1L(1)}(-x^{-2}_1)^{L(0)}v,
-x_0/x_1x_2)\cdot \nno\\
&&\;\;\;\;\;\;\;\;\;\;\cdot e^{x_2L(1)}e^{(2r+1)\pi iL(0)}
(x_2^{-L(0)})^2 w_{(1)},x^{-1}_2)w_{(2)}\Biggr\rangle_{W_3},
\end{eqnarray}
or equivalently,
\begin{eqnarray}\label{721}
\lefteqn{- \Biggl\langle w'_{(3)},x^{-1}_0\delta \left( {x_2-x_1\over -x_0}\right)
Y_{3}(e^{x_1L(1)}(-x^{-2}_1)^{L(0)}v,x^{-1}_1)\cdot }\nno\\
&&\;\;\;\;\;\;\;\;\;\;\cdot {\cal Y}(e^{x_2L(1)}e^{(2r+1)\pi iL(0)}
(x_2^{-L(0)})^2 w_{(1)},
x^{-1}_2)w_{(2)}\Biggr\rangle_{W_3} \nno\\
&&\;\;\;\;+ \Biggl\langle w'_{(3)},x^{-1}_0
\delta \left( {x_1-x_2\over x_0}\right)
{\cal Y}(e^{x_2L(1)}e^{(2r+1)\pi
iL(0)}(x_2^{-L(0)})^2 w_{(1)},x^{-1}_2)\cdot \nno\\
&&\;\;\;\;\;\;\;\;\;\;\cdot Y_{2}(e^{x_1L(1)}(-x^{-2}_1)^{L(0)}v,
x^{-1}_1)w_{(2)}\Biggr\rangle_{W_3} \nno\\
&&= \Biggl\langle w'_{(3)},x^{-1}_1\delta
\left( {x_2+x_0\over x_{1}}\right) {\cal Y}(Y_{1}
(e^{x_1L(1)}(-x^{-2}_1)^{L(0)}v,-x_0/x_1x_2)\cdot\nno\\
&&\;\;\;\;\;\;\;\;\;\;\cdot e^{x_2L(1)}e^{(2r+1)\pi iL(0)}
(x_2^{-L(0)})^2 w_{(1)},x^{-1}_2)w_{(2)}\Biggr\rangle_{W_3}.
\end{eqnarray}
Substituting (\ref{717}), (\ref{718}) and (\ref{719}) into (\ref{716})
and then comparing with (\ref{721}), we see
that the proof of (\ref{716}) is reduced to the proof of the formula
\begin{eqnarray}
\lefteqn{x^{-1}_1\delta \left( {x_2+x_0\over x_1}\right)
{\cal Y}(e^{x_2L(1)}e^{(2r+1)\pi iL(0)}
(x_2^{-L(0)})^2 Y_{1}(v,x_0)w_{(1)},x^{-1}_2)}\nno \\
&&= x^{-1}_1\delta \left( {x_2+x_0\over x_1}\right)
 {\cal Y}(Y_{1}(e^{x_1L(1)}(-x^{-2}_1)^{L(0)}v,-x_0/x_1x_2)\cdot \nno\\
&&\;\;\;\;\;\;\;\;\;\;\cdot e^{x_2L(1)}e^{(2r+1)\pi iL(0)}
(x_2^{-L(0)})^2 w_{(1)},x^{-1}_2),
\end{eqnarray}
or of
\begin{eqnarray}
\lefteqn{{\cal Y}(e^{x_2L(1)}e^{(2r+1)\pi iL(0)}
(x_2^{-L(0)})^2 Y_{1}(v,x_0)w_{(1)},x^{-1}_2)}\nno \\
&&= {\cal Y}(Y_{1}
(e^{(x_2+x_0)L(1)}(-(x_2+x_0)^{-2})^{L(0)}v,-x_0/(x_2+x_0)x_2)\cdot\nno \\
&&\;\;\;\;\;\;\;\;\;\;\cdot e^{x_2L(1)}e^{(2r+1)\pi iL(0)}
(x_2^{-L(0)})^2 w_{(1)},x^{-1}_2).
\end{eqnarray}
We see that we need only prove
\begin{eqnarray}
\lefteqn{e^{x_2L(1)}e^{(2r+1)\pi iL(0)}
(x_2^{-L(0)})^2 Y_{1}(v,x_0)}\nno\\
&&=Y_{1}(e^{(x_2+x_0)L(1)}(-(x_2+x_0)^{-2})^{L(0)}v,
-x_0/(x_2+x_0)x_2)\cdot\nno\\
&&\;\;\;\;\;\;\;\;\;\;\cdot e^{x_2L(1)}e^{(2r+1)\pi iL(0)}
(x_2^{-L(0)})^2
\end{eqnarray}
or equivalently, the conjugation formula
\begin{eqnarray}\label{725}
\lefteqn{e^{x_{2}L(1)}e^{(2r+1)\pi iL(0)}
(x_2^{-L(0)})^2 Y_{1}(v,x_0)(x_2^{L(0)})^2 e^{-(2r+1)\pi iL(0)}
e^{-x_{2}L(1)}}\nno \\
&&= Y_1(e^{(x_{2}+x_0)L(1)}(-(x_{2}+x_0)^{-2})^{L(0)}v,-x_0/(x_{2}+x_0)x_{2})
\end{eqnarray}
for $v\in V$, acting on the module $W_{1}$. But
formula (\ref{725}) follows {from} (\ref{log:p2}), (\ref{log:p3}) and the
formula
\begin{equation}
e^{(2r+1)\pi iL(0)}Y_{1}(v,x)e^{-(2r+1)\pi iL(0)}=Y_{1}((-1)^{L(0)}v, -x),
\end{equation}
which is a special case of (\ref{710}).
This establishes the Jacobi identity.

The $L(-1)$-derivative property follows {}from the same argument used
in the proof of Theorem 5.5.1 of \cite{FHL}: We have (omitting the
subscript $W_3$ on the pairings after a certain point)
\begin{eqnarray}\label{log:ArL(-1)}
\lefteqn{\left\langle \frac{d}{dx}A_r({\cal Y})(w_{(1)},x)w'_{(3)},w_{(2)}
\right\rangle_{W_2}=\frac{d}{dx}\langle A_r({\cal Y})(w_{(1)},x)w'_{(3)},w_{(2)}
\rangle_{W_2}}\nno\\
&&=\frac{d}{dx}\langle w'_{(3)},{\cal Y}(e^{xL(1)}e^{(2r+1)\pi iL(0)}
(x^{-L(0)})^2 w_{(1)},x^{-1})w_{(2)}\rangle_{W_3}\nno\\
&&=\langle w'_{(3)},\frac{d}{dx}{\cal Y}(e^{xL(1)}e^{(2r+1)\pi iL(0)}
(x^{-L(0)})^2 w_{(1)},x^{-1})w_{(2)}\rangle\nno\\
&&=\langle w'_{(3)},{\cal Y}(\frac{d}{dx}(e^{xL(1)}e^{(2r+1)\pi iL(0)}
(x^{-L(0)})^2 w_{(1)}),x^{-1})w_{(2)}\rangle\nno\\
&&\quad+\langle w'_{(3)},\frac{d}{dx}{\cal Y}(w,x^{-1})|_{w=
e^{xL(1)}e^{(2r+1)\pi iL(0)} (x^{-L(0)})^2 w_{(1)}}
w_{(2)}\rangle\nno\\
&&=\langle w'_{(3)},{\cal Y}(e^{xL(1)}L(1)e^{(2r+1)\pi iL(0)}
(x^{-L(0)})^2 w_{(1)},x^{-1})w_{(2)}\rangle\nno\\
&&\quad+\langle w'_{(3)},{\cal Y}(e^{xL(1)}(-2L(0)x^{-1})
e^{(2r+1)\pi iL(0)}(x^{-L(0)})^2w_{(1)},x^{-1})w_{(2)}\rangle\nno\\
&&\quad+\langle w'_{(3)},\frac{d}{d{x^{-1}}}{\cal Y}(w,x^{-1})|_{w=
e^{xL(1)}e^{(2r+1)\pi iL(0)} (x^{-L(0)})^2 w_{(1)}}
w_{(2)}\rangle(-x^{-2})\nno\\
&&=\langle w'_{(3)},{\cal Y}(e^{xL(1)}(2xL(0)-x^2L(1))e^{(2r+1)\pi iL(0)}
(x^{-L(0)})^2(-x^{-2})w_{(1)},x^{-1})w_{(2)}\rangle\nno\\
&&\quad+\langle w'_{(3)},{\cal Y}(L(-1)e^{xL(1)}e^{(2r+1)\pi iL(0)}
(x^{-L(0)})^2 w_{(1)},x^{-1})w_{(2)}\rangle(-x^{-2})\nno\\
&&=\langle w'_{(3)},{\cal Y}(e^{xL(1)}(2xL(0)-x^2L(1))e^{(2r+1)\pi iL(0)}
(x^{-L(0)})^2(-x^{-2})w_{(1)},x^{-1})w_{(2)}\rangle\nno\\
&&\quad+\langle w'_{(3)},{\cal Y}(L(-1)e^{xL(1)}e^{(2r+1)\pi iL(0)}
(x^{-L(0)})^2 (-x^{-2})w_{(1)},x^{-1})w_{(2)}\rangle
\end{eqnarray}
Now by (\ref{log:SL2-3}), with $x$ replaced by $-x$, we have
\[
L(-1)e^{xL(1)}=e^{xL(1)}(L(-1)-2xL(0)+x^2L(1)).
\]
Using this together with (\ref{log:xLx^}) and (\ref{log:SL2-2}) (with
$x$ specialized to $-(2r+1)\pi i$, and the convergence of the
exponential series invoked; recall Remark
\ref{analyticallyconvervent}), we see that the right-hand side of
(\ref{log:ArL(-1)}) equals
\begin{eqnarray*}
\lefteqn{\langle w'_{(3)},{\cal Y}(e^{xL(1)}L(-1)e^{(2r+1)\pi iL(0)}
(x^{-L(0)})^2(-x^{-2})w_{(1)},x^{-1})w_{(2)}\rangle}\\
&&=\langle w'_{(3)},{\cal Y}(e^{xL(1)}e^{(2r+1)\pi iL(0)}(-L(-1))
(x^{-L(0)})^2(-x^{-2})w_{(1)},x^{-1})w_{(2)}\rangle\\
&&=\langle w'_{(3)},{\cal Y}(e^{xL(1)}e^{(2r+1)\pi iL(0)}
(x^{-L(0)})^2L(-1)w_{(1)},x^{-1})w_{(2)}\rangle\\
&&=\langle A_r({\cal Y})(L(-1)w_{(1)},x)w'_{(3)},w_{(2)}\rangle_{W_2},
\end{eqnarray*}
as desired.

We now show that, in case $V$ is M\"obius, the ${\mathfrak
s}{\mathfrak l}(2)$-bracket relations (\ref{log:L(j)b}) hold for
$A_r({\cal Y})$.  For these, we first see that, for $j=-1,0,1$, by
using (\ref{L'(n)}), (\ref{log:Ardef}) and the ${\mathfrak
s}{\mathfrak l}(2)$-bracket relations (\ref{log:L(j)b}) for ${\cal Y}$
we have
\begin{eqnarray}\label{log:tmp1}
\lefteqn{\langle L'(j)A_r({\cal Y})(w_{(1)},x)w'_{(3)},w_{(2)}
\rangle_{W_2}}\nno\\
&&=\langle A_r({\cal Y})(w_{(1)},x)w'_{(3)},L(-j)w_{(2)}
\rangle_{W_2}\nno\\
&&=\langle w'_{(3)},{\cal Y}(e^{xL(1)}e^{(2r+1)\pi iL(0)}(x^{-L(0)})^2
w_{(1)},x^{-1})L(-j)w_{(2)}\rangle_{W_3}\nno\\
&&=\langle w'_{(3)},L(-j){\cal Y}(e^{xL(1)}e^{(2r+1)\pi iL(0)}
(x^{-L(0)})^2w_{(1)},x^{-1})w_{(2)}\rangle_{W_3}\nno\\
&&\hspace{1em}-\Biggl\langle w'_{(3)},\sum_{i=0}^{-j+1}{-j+1\choose i}x^{-i}
\cdot\nno\\
&&\hspace{8em}\cdot{\cal Y}(L(-j-i)e^{xL(1)}e^{(2r+1)\pi iL(0)}
(x^{-L(0)})^2w_{(1)},x^{-1})w_{(2)}\Biggr\rangle_{W_3}.
\end{eqnarray}
Now {}from (\ref{log:SL2-3}), (\ref{log:SL2-2}) (with $x$ specialized to
$-(2r+1)\pi i$) and (\ref{log:xLx^}), one computes that
\begin{eqnarray*}
\lefteqn{(x^{L(0)})^2e^{-(2r+1)\pi iL(0)}e^{-xL(1)}
\left(\begin{array}{c}L(-1)\\L(0)\\L(1)\end{array}\right)
e^{xL(1)}e^{(2r+1)\pi iL(0)}(x^{-L(0)})^2}\\
&&\hspace{5cm}=\left(\begin{array}{ccc}-x^2&-2x&-1\\0&1&x^{-1}
\\0&0&-x^{-2}\end{array}\right)
\left(\begin{array}{c}L(-1)\\L(0)\\L(1)\end{array}\right)
\end{eqnarray*}
on $W_1$, which implies that
\begin{eqnarray*}
\lefteqn{\sum_{i=0}^{-j+1}{-j+1\choose i}x^{-i}L(-j-i)e^{xL(1)}
e^{(2r+1)\pi iL(0)}(x^{-L(0)})^2}\\
&&\hspace{2em}=-\sum_{i=0}^{j+1}{j+1\choose i}x^ie^{xL(1)}e^{(2r+1)
\pi iL(0)}(x^{-L(0)})^2L(j-i)
\end{eqnarray*}
on $W_1$.  Hence the right-hand side of
(\ref{log:tmp1}) is equal to
\begin{eqnarray*}
&&\lefteqn{\langle L'(j)w'_{(3)},{\cal Y}(e^{xL(1)}e^{(2r+1)\pi iL(0)}
(x^{-L(0)})^2w_{(1)},x^{-1})w_{(2)}\rangle_{W_3}}\\
&&\hspace{1em}+\sum_{i=0}^{j+1}{j+1\choose i}x^i\langle w'_{(3)},
{\cal Y}(e^{xL(1)}e^{(2r+1)\pi iL(0)}(x^{-L(0)})^2
L(j-i)w_{(1)},x^{-1})w_{(2)} \rangle_{W_2}\\
&&=\langle A_r({\cal Y})(w_{(1)},x)L'(j)w'_{(3)},w_{(2)}
\rangle_{W_2}\\
&&\hspace{1em}+\sum_{i=0}^{j+1}{j+1\choose i}x^i\langle A_r({\cal Y})
(L(j-i)w_{(1)},x)w'_{(3)},w_{(2)} \rangle_{W_2},
\end{eqnarray*}
and the ${\mathfrak s}{\mathfrak l}(2)$-bracket relations for
$A_r({\cal Y})$ are proved.  (Note that this argument essentially
generalizes the proof of Lemma \ref{sl2opposite}.)  We have finished
proving that $A_r({\cal Y})$ is a grading-compatible logarithmic
intertwining operator.

Finally, for the relation (\ref{log:ar}), we of course identify
$W''_2$ with $W_2$ and $W''_3$ with $W_3$, according to Theorem
\ref{set:W'}.  Let us view ${\cal Y}$ as a grading-compatible
logarithmic intertwining operator of type ${W'_2\choose W_1\,W'_3}$,
so that $A_r({\cal Y})$ is such an operator of type ${W_3\choose
W_1\,W_2}$.  We have
\begin{eqnarray*}
\lefteqn{\langle A_{-r-1}A_r({\cal Y})(w_{(1)},x)w'_{(3)},w_{(2)}
\rangle_{W_2}}\\
&&=\langle w'_{(3)},A_r({\cal Y})(e^{xL(1)}e^{(-2r-1)\pi iL(0)}
(x^{-L(0)})^2 w_{(1)},x^{-1})w_{(2)}\rangle_{W_3}\\
&&=\langle {\cal Y}(e^{x^{-1}L(1)}e^{(2r+1)\pi iL(0)}(x^{L(0)})^2
e^{xL(1)}e^{(-2r-1)\pi iL(0)}(x^{-L(0)})^2w_{(1)},x)w'_{(3)},w_{(2)}
\rangle_{W_2}\\
&&=\langle {\cal Y}(w_{(1)},x)w'_{(3)},w_{(2)}\rangle_{W_2},
\end{eqnarray*}
where the last equality is due to the relation
\begin{equation}\label{conjrelation}
e^{(2r+1)\pi iL(0)}(x^{L(0)})^2 e^{xL(1)}(x^{-L(0)})^2 e^{-(2r+1)\pi
iL(0)} = e^{-x^{-1}L(1)}
\end{equation}
on $W_{1}$, whose proof is similar to that of formula (5.3.1) of
\cite{FHL}.  Namely, (\ref{conjrelation}) follows {}from the relation
\begin{equation}\label{xto-1/x}
e^{(2r+1)\pi iL(0)}(x^{L(0)})^2 xL(1) (x^{-L(0)})^2 e^{-(2r+1)\pi
iL(0)} = -x^{-1}L(1),
\end{equation}
which realizes the transformation $x \mapsto -\frac{1}{x}$, and 
(\ref{xto-1/x}) follows {}from (\ref{log:xLx^}) together with 
(\ref{log:SL2-2}) specialized as above.  \epf

\begin{rema}{\rm
The last argument in the proof shows that for any $r,s\in{\mathbb Z}$,
formula (\ref{log:ar}) generalizes to:
\[
A_s(A_r({\cal Y}))={\cal Y}_{[0,r+s+1,0]}
\]
(recall Remark \ref{Ys1s2s3}).
}
\end{rema}

With $V$, $W_1$, $W_2$ and $W_3$ strongly graded, set
\begin{equation}
N_{W_1W_2W_3}=N_{W_1\,W_2}^{W'_3}.
\end{equation}
Then Proposition \ref{log:omega} gives 
\[
N_{W_1W_2W_3}=N_{W_2W_1W_3}
\]
and Proposition \ref{log:A} gives
\[
N_{W_1W_2W_3}=N_{W_1W_3W_2}.
\]
Thus for any permutation $\sigma$ of $(1,2,3)$,
\begin{equation}
N_{W_1W_2W_3}=N_{W_{\sigma(1)}W_{\sigma(2)}W_{\sigma(3)}}.
\end{equation}

It is clear {}from Proposition \ref{log:logwt}(b) that in the nontrivial
logarithmic intertwining operator case, taking projections of ${\cal
Y}(w_{(1)},x)w_{(2)}$ to (generalized) weight subspaces is not enough
to recover its coefficients of $x^n(\log x)^k$ for $n\in{\mathbb C}$
and $k\in{\mathbb N}$, in contrast with the (ordinary) intertwining
operator case (cf.\ \cite{tensor1}, the paragraph containing formula
(4.17)).  However, taking projections of certain related intertwining
operators does indeed suffice for this purpose:

\begin{propo}\label{log:proj}
Let $W_1$, $W_2$, $W_3$ be generalized modules for a M\"obius (or
conformal) vertex algebra $V$ and let ${\cal Y}$ be a logarithmic
intertwining operator of type ${W_3\choose W_1\,W_2}$.  Let
$w_{(1)}\in W_1$ and $w_{(2)}\in W_2$ be homogeneous of generalized
weights $n_1$ and $n_2$, respectively. Then for any $n\in {\mathbb C}$
and any $r\in {\mathbb N}$, ${w_{(1)}}^{\cal Y}_{n;\,r}w_{(2)}$ can be
written as a certain linear combination of products of the component of 
weight $n_1+n_2-n-1$ of
\[
(L(0)-n_1-n_2+n+1)^l{\cal Y}((L(0)-n_1)^iw_{(1)},x)(L(0)-n_2)^jw_{(2)}
\]
for certain $i,j,l\in{\mathbb N}$ with monomials of the form
$x^{n+1}(\log x)^m$ for certain $m\in{\mathbb N}$.
\end{propo}
\pf Multiplying (\ref{log:r+t=?}) by $x^{-n-1}(\log x)^k$ and summing
over $k\in {\mathbb N}$ (a finite sum by definition) we have that for
any $t\in {\mathbb N}$,
\begin{eqnarray}\label{log:proj1}
\lefteqn{\sum_{k\in {\mathbb N}}{k+t\choose t}{w_{(1)}}^{\cal
Y}_{n;\,k+t}w_{(2)}x^{-n-1}(\log x)^k \;
\bigg(=x^{-n-1}\sum_{k\in {\mathbb N}}{k\choose t}
{w_{(1)}}^{\cal Y}_{n;\,k}w_{(2)}(\log x)^{k-t}\bigg)}\nno\\
&&=\sum_{i,j,l\in {\mathbb N}, \, i+j+l=t}\frac{1}{i!j!l!}(-1)^{i+j}
\sum_{k\in {\mathbb N}}(L(0)-n_1-n_2+n+1)^l\cdot\nno\\
&&\quad\cdot((L(0)-n_1)^i{w_{(1)}})^{\cal Y}_{n;\,k}
(L(0)-n_2)^jw_{(2)})x^{-n-1}(\log x)^k\nno\\
&&=\sum_{i,j,l\in {\mathbb N}, \, i+j+l=t}\frac{1}{i!j!l!}(-1)^{i+j}
\pi_{n_1+n_2-n-1}((L(0)-n_1-n_2+n+1)^l\cdot\nno\\
&&\quad\cdot{\cal Y}((L(0)-n_1)^i{w_{(1)}},x) (L(0)-n_2)^jw_{(2)}).
\end{eqnarray}
Let $K$ be a positive integer such that ${w_{(1)}}^{\cal Y}_{n;\,k'}
w_{(2)}=0$ for all $k'\geq K$, Denote the right-hand side of
(\ref{log:proj1}) by $\pi(t,w_{(1)},w_{(2)},x,\log x)$.  Then by
putting the identities (\ref{log:proj1}) for $t=0,1,\dots,K-1$ together in
matrix form we have
\begin{equation}\label{log:projmat}
x^{-n-1}A\left(\begin{array}{c}{w_{(1)}}^{\cal Y}_{n;\,0}w_{(2)}\\
{w_{(1)}}^{\cal Y}_{n;\,1}w_{(2)}\\\vdots\\
{w_{(1)}}^{\cal Y}_{n;\,K-1}w_{(2)}\end{array}\right)=
\left(\begin{array}{c}\pi(0,w_{(1)},w_{(2)},x,\log x)\\
\pi(1,w_{(1)},w_{(2)},x,\log x)\\\vdots\\
\pi(K-1,w_{(1)},w_{(2)},x,\log x)\end{array}\right)
\end{equation}
where $A$ is the $K\times K$ matrix whose $(i,j)$-entry is equal to
$\displaystyle{j-1\choose i-1}(\log x)^{j-i}$.  Letting $P_K$ be the
triangular matrix whose $(i,j)$-entry is $\displaystyle{j-1\choose
i-1}$ (an upper triangular ``Pascal matrix''), we have
\[
A={\rm diag}(1,(\log x)^{-1},\dots,(\log x)^{-(K-1)})\cdot P_K\cdot
{\rm diag}(1,\log x,\dots,(\log x)^{K-1}).
\]
Its inverse is
\[
A^{-1}={\rm diag}(1,(\log x)^{-1},\dots,(\log x)^{-(K-1)})\cdot
P_K^{-1}\cdot{\rm diag}(1,\log x,\dots,(\log x)^{K-1})
\]
and the $(i,j)$-entry of $P_K^{-1}$ is $\displaystyle(-1)^{i+j}
{j-1\choose i-1}$.  Now multiplying the left-hand side of
(\ref{log:projmat}) by $x^{n+1}A^{-1}$ we obtain
\[
\left(\begin{array}{c}{w_{(1)}}^{\cal Y}_{n;\,0}w_{(2)}\\
{w_{(1)}}^{\cal Y}_{n;\,1}w_{(2)}\\\vdots\\
{w_{(1)}}^{\cal Y}_{n;\,K-1}w_{(2)}\end{array}\right)=x^{n+1}A^{-1}
\left(\begin{array}{c}\pi(0,w_{(1)},w_{(2)},x,\log x)\\
\pi(1,w_{(1)},w_{(2)},x,\log x)\\\vdots\\
\pi(K-1,w_{(1)},w_{(2)},x,\log x)\end{array}\right)
\]
or explicitly,
\begin{equation}\label{log:last}
(w_{(1)})^{\cal Y}_{n;\,r}w_{(2)}=x^{n+1}\sum_{t=r}^{K-1}(-1)^{r+t}
{t\choose r}(\log x)^{t-r}\pi(t,w_{(1)},w_{(2)},x,\log x)
\end{equation}
for $r=0,1,\dots,K-1$. (In particular, all $x$'s and $\log x$'s cancel
out in the right-hand side of (\ref{log:last}).)  \epfv

\newpage

\setcounter{equation}{0}
\setcounter{rema}{0}

\section{The notions of $P(z)$- and $Q(z)$-tensor product}

We now generalize to the setting of the present work the notions of
$P(z)$- and $Q(z)$-tensor product of modules introduced in
\cite{tensor1}, \cite{tensor2} and \cite{tensor3}.  We introduce the
notions of $P(z)$- and $Q(z)$-intertwining map among strongly
$\tilde{A}$-graded generalized modules for a strongly $A$-graded
M\"{o}bius or conformal vertex algebra $V$ and establish the
relationship between such intertwining maps and grading-compatible
logarithmic intertwining operators.  We define the $P(z)$- and
$Q(z)$-tensor products of two strongly $\tilde{A}$-graded generalized
$V$-modules using these intertwining maps and natural universal
properties. As examples, for a strongly $\tilde{A}$-graded generalized
module $W$, we construct and describe the $P(z)$-tensor products of
$V$ and $W$ and also of $W$ and $V$; the underlying strongly
$\tilde{A}$-graded generalized modules of the tensor product
structures are $W$ itself, in both of these cases.  In the case in
which $V$ is a finitely reductive vertex operator algebra (recall the
Introduction), we construct and describe the $P(z)$- and $Q(z)$-tensor
products of arbitrary $V$-modules, and we use this structure to
motivate the construction of associativity isomorphisms that we will
carry out in later sections.

In view of the results in Sections 2 and 3 involving contragredient
modules, it is natural for us to work in the strongly-graded setting
{}from now on:

\begin{assum}\label{assum}
Throughout this section and the remainder of this work, we shall
assume the following, unless other assumptions are explicitly made:
$A$ is an abelian group and $\tilde{A}$ is an abelian group containing
$A$ as a subgroup; $V$ is a strongly $A$-graded M\"{o}bius or
conformal vertex algebra; all $V$-modules and generalized $V$-modules
considered are strongly $\tilde{A}$-graded; and all intertwining
operators and logarithmic intertwining operators considered are
grading-compatible.  (Recall Definitions \ref{def:dgv}, \ref{def:dgw},
\ref{log:def} and \ref{gradingcompatintwop}.)  Also, in this section,
$z$ will be a fixed nonzero complex number.
\end{assum}

\subsection{The notion of $P(z)$-tensor product}

We first generalize the notion of $P(z)$-intertwining map given in
Section 4 of \cite{tensor1}; our $P(z)$-intertwining maps will
automatically be grading-compatible by definition.  We use the
notations given in Definition \ref{Wbardef}.  The main part of the
following definition, the Jacobi identity (\ref{im:def}), was
previewed in the Introduction (formula (\ref{im-jacobi})).  It should
be compared with the corresponding formula (\ref{intwmap}) in the Lie
algebra setting, and with the Jacobi identity (\ref{log:jacobi}) in
the definition of the notion of logarithmic intertwining operator;
note that the formal variable $x_2$ in that Jacobi identity is
specialized here to the nonzero complex number $z$.  Also, the
${\mathfrak s}{\mathfrak l}(2)$-bracket relations (\ref{im:Lj}) should
be compared with the corresponding relations (\ref{log:L(j)b}).  There
is no $L(-1)$-derivative formula for intertwining maps; as we shall
see, the $P(z)$-intertwining maps are obtained {}from logarithmic
intertwining operators by a process of specialization of the formal
variable to the complex variable $z$.

\begin{defi}\label{im:imdef}{\rm
Let $(W_1,Y_1)$, $(W_2,Y_2)$ and $(W_3,Y_3)$ be generalized
$V$-modules.  A {\it $P(z)$-intertwining map of type ${W_3\choose
W_1\,W_2}$} is a linear map
\begin{equation}\label{PzintwmapI}
I: W_1\otimes W_2 \to \overline{W}_3
\end{equation}
(recall {}from Definition \ref{Wbardef} that $\overline{W}_3$ is the
formal completion of $W_3$ with respect to the ${\mathbb C}$-grading)
such that the following conditions are satisfied: the {\it grading
compatibility condition}: for $\beta, \gamma\in \tilde{A}$ and
$w_{(1)}\in W_{1}^{(\beta)}$, $w_{(2)}\in W_{2}^{(\gamma)}$,
\begin{equation}\label{grad-comp}
I(w_{(1)}\otimes w_{(2)})\in \overline{W_{3}^{(\beta+\gamma)}};
\end{equation}
the
{\em lower truncation condition}: for any elements
$w_{(1)}\in W_1$, $w_{(2)}\in W_2$, and any $n\in {\mathbb C}$,
\begin{equation}\label{im:ltc}
\pi_{n-m}I(w_{(1)}\otimes w_{(2)})=0\;\;\mbox{ for }\;m\in {\mathbb N}
\;\mbox{ sufficiently large}
\end{equation}
(which follows {}from (\ref{grad-comp}), in view of the
grading restriction condition (\ref{set:dmltc}); recall the notation
$\pi_n$ {}from Definition \ref{Wbardef}); the {\em Jacobi identity}:
\begin{eqnarray}\label{im:def}
\lefteqn{x_0^{-1}\delta\bigg(\frac{ x_1-z}{x_0}\bigg)
Y_3(v, x_1)I(w_{(1)}\otimes w_{(2)})}\nno\\
&&=z^{-1}\delta\bigg(\frac{x_1-x_0}{z}\bigg)
I(Y_1(v, x_0)w_{(1)}\otimes w_{(2)})\nno\\
&&\hspace{2em}+x_0^{-1}\delta\bigg(\frac{z-x_1}{-x_0}\bigg)
I(w_{(1)}\otimes Y_2(v, x_1)w_{(2)})
\end{eqnarray}
for $v\in V$, $w_{(1)}\in W_1$ and $w_{(2)}\in W_2$ (note that all the
expressions in the right-hand side of (\ref{im:def}) are well defined,
and that the left-hand side of (\ref{im:def}) is meaningful because
any infinite linear combination of $v_n$ ($n\in{\mathbb Z}$) of the
form $\sum_{n<N}a_nv_n$ ($a_n\in {\mathbb C}$) acts in a well-defined
way on any $I(w_{(1)}\otimes w_{(2)})$, in view of (\ref{im:ltc}));
and the {\em ${\mathfrak s}{\mathfrak l}(2)$-bracket relations}: for any
$w_{(1)}\in W_1$ and $w_{(2)}\in W_2$,
\begin{equation}\label{im:Lj}
L(j)I(w_{(1)}\otimes w_{(2)})=I(w_{(1)}\otimes L(j)w_{(2)})+
\sum_{i=0}^{j+1}{j+1\choose i}z^iI((L(j-i)w_{(1)})\otimes w_{(2)})
\end{equation}
for $j=-1, 0$ and $1$ (note that if $V$ is in fact a conformal vertex
algebra, this follows automatically {}from (\ref{im:def}) by setting
$v=\omega$ and taking $\res_{x_0}\res_{x_1}x_1^{j+1}$). The vector
space of $P(z)$-intertwining maps of type ${W_3}\choose {W_1W_2}$ is
denoted by ${\cal M}[P(z)]^{W_3}_{W_1W_2}$, or simply by ${\cal
M}^{W_3}_{W_1W_2}$ if there is no ambiguity. }
\end{defi}

\begin{rema}\label{P(z)geometry}
{\rm As we mentioned in the Introduction, $P(z)$ is the Riemann sphere
$\hat{\mathbb C}$ with one negatively oriented puncture at $\infty$
and two ordered positively oriented punctures at $z$ and $0$, with
local coordinates $1/w$, $w-z$ and $w$, respectively, vanishing at
these three punctures.  The geometry underlying the notion of
$P(z)$-intertwining map and the notions of $P(z)$-product and
$P(z)$-tensor product (see below) is determined by $P(z)$.}
\end{rema}

\begin{rema}{\rm
In the case of $\C$-graded ordinary modules for a vertex operator
algebra, where the grading restriction condition (\ref{Wn+k=0}) for a
module $W$ is replaced by the (more restrictive) condition
\begin{equation}
W_{(n)}=0 \;\; \mbox { for }\;n\in {\C}\;\mbox{ with sufficiently
negative real part}
\end{equation}
as in \cite{tensor1} (and where, in our context, the abelian groups
$A$ and $\tilde{A}$ are trivial), the notion of $P(z)$-intertwining
map above agrees with the earlier one introduced in \cite{tensor1}; in
this case, the conditions (\ref{grad-comp}) and (\ref{im:ltc}) are
automatic.}
\end{rema}

\begin{rema}{\rm
As in Remark \ref{log:Lj2rema}, it is clear that the ${\mathfrak s}{\mathfrak
l}(2)$-bracket relations (\ref{im:Lj}) can equivalently be written as
\begin{eqnarray}\label{im:Lj2}
I(L(j)w_{(1)}\otimes w_{(2)})&=& \sum_{i=0}^{j+1}{j+1\choose
i}(-z)^i L(j-i)I(w_{(1)}\otimes w_{(2)})\nno\\
&&-\sum_{i=0}^{j+1}{j+1\choose i}(-z)^iI(w_{(1)}\otimes
L(j-i)w_{(2)})
\end{eqnarray}
for $w_{(1)}\in W_1$, $w_{(2)}\in W_2$ and $j=-1,0$ and $1$.
}
\end{rema}

Following \cite{tensor1} we will choose the branch of $\log z$ so that
\begin{equation}\label{branch1}
0\leq {\rm Im }(\log z) < 2 \pi
\end{equation}
(despite the fact that we happened to have used a different branch in
(\ref{log:br1}) in the proof of Theorem \ref{log:ids}).  We will also
use the notation
\begin{equation}\label{branch2}
l_p(z)=\log z+2\pi ip, \ p\in {\mathbb Z},
\end{equation}
as in \cite{tensor1}, for arbitrary values of the $\log$ function. For
a formal expression $f(x)$ as in (\ref{log:f}), but involving only
nonnegative integral powers of $\log x$, and $\zeta\in {\mathbb C}$,
whenever
\begin{equation}\label{log:fsub}
f(x)\lbar_{x^n=e^{\zeta n},\;(\log x)^m=\zeta^m,\;n\in{\mathbb C},\;
m\in \N}
\end{equation}
exists algebraically, we will write (\ref{log:fsub}) simply as
$f(x)\lbar_{x=e^{\zeta}}$ or $f(e^\zeta)$, and we will call this
``substituting $e^\zeta$ for $x$ in $f(x)$,'' even though, in general,
it depends on $\zeta$, not just on $e^\zeta$.  (See also
(\ref{log:subs}).) In addition, for a fixed integer $p$, we will
sometimes write
\begin{equation}\label{im:f(z)}
f(x)\lbar_{x=z}\;\;\mbox{or}\;\;f(z)
\end{equation}
instead of $f(x)\lbar_{x=e^{l_p(z)}}$ or $f(e^{l_p(z)})$.  We will
sometimes say that ``$f(e^{\zeta})$ exists'' or that ``$f(z)$
exists.''

\begin{rema}{\rm
A very important example of an $f(z)$ existing in this sense occurs
when $f(x)={\cal Y}(w_{(1)},x)w_{(2)}$ ($\in W_3[\log x]\{x\}$) for
$w_{(1)}\in W_1$, $w_{(2)}\in W_2$ and a logarithmic intertwining
operator ${\cal Y}$ of type ${W_3\choose W_1\,W_2}$, in the notation
of Definition \ref{log:def}; note that (\ref{log:fsub}) exists (as an
element of $\overline{W_{3}}$) in this case because of Proposition
\ref{log:logwt}(b).  Note also that in particular, ${\cal
Y}(w_{(1)},e^{\zeta})$ (or ${\cal Y}(w_{(1)},z)$) exists as a linear
map {}from $W_2$ to $\overline{W_{3}}$, and that ${\cal
Y}(\cdot,z)\cdot$ exists as a linear map
\begin{eqnarray}
W_1\otimes W_2 &\rightarrow & \overline{W}_3 \nonumber\\ 
w_{(1)}\otimes w_{(2)} &\mapsto & {\cal Y}(w_{(1)},z)w_{(2)}.
\end{eqnarray}
}
\end{rema}

Now we use these considerations to construct correspondences between
(grading-compatible) logarithmic intertwining operators and
$P(z)$-intertwining maps.  Fix an integer $p$. Let ${\cal Y}$ be a
logarithmic intertwining operator of type ${W_3\choose
W_1\,W_2}$.  Then we have a linear map
\begin{equation}
I_{{\cal Y},p}: W_1\otimes W_2\to \overline{W}_3
\end{equation}
defined by
\begin{equation}\label{log:IYp}
I_{{\cal Y},p}(w_{(1)}\otimes w_{(2)})={\cal
Y}(w_{(1)},e^{l_p(z)})w_{(2)}
\end{equation}
for all $w_{(1)}\in W_1$ and $w_{(2)}\in W_2$.  The
grading-compatibility condition (\ref{gradingcompatcondn}) yields the
grading-compatibility condition (\ref{grad-comp}) for $I_{{\cal
Y},p}$, and (\ref{im:ltc}) follows.  By substituting $e^{l_p(z)}$ for
$x_2$ in (\ref{log:jacobi}) and for $x$ in (\ref{log:L(j)b}), we see
that $I_{{\cal Y},p}$ satisfies the Jacobi identity (\ref{im:def}) and
the ${\mathfrak s}{\mathfrak l}(2)$-bracket relations (\ref{im:Lj}).
Hence $I_{{\cal Y},p}$ is a $P(z)$-intertwining map.

On the other hand, we note that (\ref{log:p2}) is equivalent to
\begin{equation}\label{log:4.14}
\langle y^{L'(0)}w'_{(3)}, {\cal Y}(y^{-L(0)}w_{(1)},
x)y^{-L(0)}w_{(2)}\rangle_{W_3} =\langle w'_{(3)}, {\cal
Y}(w_{(1)}, xy)w_{(2)}\rangle_{W_3}
\end{equation}
for all $w_{(1)}\in W_1$, $w_{(2)}\in W_2$ and $w'_{(3)}\in W'_3$,
where we are using the pairing between the contragredient module
$W'_3$ and $W_3$ or $\overline{W}_3$ (recall Definition
\ref{defofWprime}, Theorem \ref{set:W'}, (\ref{L'(n)}),
(\ref{truncationforY'}) and (\ref{log:x^L(0)})).  Substituting
$e^{l_{p}(z)}$ for $x$ and then $e^{-l_{p}(z)}x$ for $y$, we obtain
\begin{eqnarray*}
&\langle y^{L'(0)}x^{L'(0)}w'_{(3)}, {\cal Y}
(y^{-L(0)}x^{-L(0)}w_{(1)}, e^{l_{p}(z)})
y^{-L(0)}x^{-L(0)}w_{(2)}\rangle_{W_3}\lbar_{y=e^{-l_{p}(z)}}&\\
&=\langle w'_{(3)}, {\cal Y}(w_{(1)}, x)w_{(2)}\rangle_{W_3},&
\end{eqnarray*}
or equivalently, using the notation (\ref{log:IYp}),
\begin{eqnarray*}
&\langle w'_{(3)}, y^{L(0)}x^{L(0)}I_{{\cal Y}, p}
(y^{-L(0)}x^{-L(0)}w_{(1)}\otimes
y^{-L(0)}x^{-L(0)}w_{(2)})\rangle_{W_3}\lbar_{y=e^{-l_{p}(z)}}&\\
&=\langle w'_{(3)}, {\cal Y}(w_{(1)}, x)w_{(2)}\rangle_{W_3}.&
\end{eqnarray*}
Thus we have recovered ${\cal Y}$ {}from $I_{{\cal Y},p}$.

This motivates the following definition: Given a $P(z)$-intertwining
map $I$ and an integer $p$, we define a linear map
\begin{equation}\label{YIp}
{\cal Y}_{I,p}:W_1\otimes W_2\to W_3[\log x]\{x\}
\end{equation}
by
\begin{eqnarray}\label{recover}
\lefteqn{{\cal Y}_{I,p}(w_{(1)}, x)w_{(2)}}\nno\\
&&=y^{L(0)}x^{L(0)}I(y^{-L(0)}x^{-L(0)}w_{(1)}\otimes
y^{-L(0)}x^{-L(0)}w_{(2)})\lbar_{y=e^{-l_{p}(z)}}
\end{eqnarray}
for any $w_{(1)}\in W_1$ and $w_{(2)}\in W_2$ (this is well defined
and indeed maps to $W_3[\log x]\{x\}$, in view of (\ref{log:x^L(0)})).
We will also use the notation ${w_{(1)}}_{n;k}^{I,p}w_{(2)}\in W_3$
defined by
\begin{equation}\label{wInkw}
{\cal Y}_{I,p}(w_{(1)}, x)w_{(2)}=\sum_{n\in{\mathbb
C}}\sum_{k\in {\mathbb N}} {w_{(1)}}_{n;k}^{I,p}w_{(2)}
x^{-n-1}(\log x)^k.
\end{equation}
Observe that since the operator $x^{\pm L(0)}$ always increases the
power of $x$ in an expression homogeneous of generalized weight $n$ by
$\pm n$, we see {}from (\ref{recover}) that
\begin{equation}
{w_{(1)}}_{n;k}^{I,p}w_{(2)}\in(W_3)_{[n_1+n_2-n-1]}
\end{equation}
for $w_{(1)}\in (W_1)_{[n_1]}$ and $w_{(2)}\in (W_2)_{[n_2]}$.
Moreover, for $I=I_{{\cal Y},p}$, we have ${\cal Y}_{I,p}={\cal Y}$,
and for ${\cal Y}={\cal Y}_{I,p}$, we have $I_{{\cal Y},p}=I$.

We can now prove the following proposition generalizing Proposition
12.2 in \cite{tensor3}.

\begin{propo}\label{im:correspond}
For $p\in {\mathbb Z}$, the correspondence ${\cal Y}\mapsto I_{{\cal
Y}, p}$ is a linear isomorphism {}from the space ${\cal
V}^{W_3}_{W_1W_2}$ of (grading-compatible) logarithmic intertwining
operators of type ${W_3\choose W_1\,W_2}$ to the space ${\cal
M}^{W_3}_{W_1W_2}$ of $P(z)$-intertwining maps of the same type. Its
inverse map is given by $I\mapsto {\cal Y}_{I,p}$.
\end{propo}
\pf We need only show that for any $P(z)$-intertwining map $I$ of type
${W_3\choose W_1\,W_2}$, ${\cal Y}_{I,p}$ is a logarithmic
intertwining operator of the same type. The lower truncation condition 
(\ref{im:ltc}) implies that
the lower truncation condition (\ref{log:ltc}) for logarithmic
intertwining operator holds for ${\cal Y}_{I,p}$. Let us now prove the
Jacobi identity for ${\cal Y}_{I,p}$.

Changing the formal variables $x_0$ and $x_1$ to
$x_0e^{l_p(z)}x_2^{-1}$ and $x_1e^{l_p(z)}x_2^{-1}$, respectively, in
the Jacobi identity (\ref{im:def}) for $I$, and then changing $v$ to
$y^{-L(0)}x_2^{-L(0)}v\lbar_{y=e^{-l_{p}(z)}}$ we obtain (noting that
at first, $e^{l_{p}(z)}$ could be written simply as $z$ because only
integral powers occur)
\begin{eqnarray*}
\lefteqn{x^{-1}_0\delta\left(\frac{x_1-x_2}{x_0}\right)
Y_{3}(y^{-L(0)}x_2^{-L(0)}v,x_1y^{-1}x_2^{-1}) I(w_{(1)}\otimes
w_{(2)})\lbar_{y=e^{-l_{p}(z)}}}\nno\\
&&=x_2^{-1}\delta\left(\frac{x_1-x_0}{x_2}\right)
I(Y_1(y^{-L(0)}x_2^{-L(0)}v, x_0y^{-1}x_2^{-1})w_{(1)}
\otimes w_{(2)})\lbar_{y=e^{-l_{p}(z)}}\nno\\
&&\quad +x_0^{-1}\delta\left(\frac{x_2-x_1}{-x_0}\right)
I(w_{(1)}\otimes Y_2(y^{-L(0)}x_2^{-L(0)}v, x_1
y^{-1}x_2^{-1})w_{(2)})\lbar_{y=e^{-l_{p}(z)}}.
\end{eqnarray*}
Using the formula
\[
Y_{3}(y^{-L(0)}x_2^{-L(0)}v,x_1y^{-1}x_2^{-1})=
y^{-L(0)}x_2^{-L(0)}Y_{3}(v,x_1)y^{L(0)}x_2^{L(0)},
\]
which holds on the generalized module $W_3$, by (\ref{log:p2}), and
the similar formulas for $Y_1$ and $Y_2$, we get
\begin{eqnarray*}
\lefteqn{x^{-1}_0\delta\left(\frac{x_1-x_2}{x_0}\right)
y^{-L(0)}x_2^{-L(0)}Y_{3}(v,x_1)y^{L(0)}x_2^{L(0)}I(w_{(1)}\otimes
w_{(2)})\lbar_{y=e^{-l_{p}(z)}}}\nno\\
&&=x_2^{-1}\delta\left(\frac{x_1-x_0}{x_2}\right)
I(y^{-L(0)}x_2^{-L(0)}Y_1(v, x_0)y^{L(0)}x_2^{L(0)}w_{(1)}
\otimes w_{(2)})\lbar_{y=e^{-l_{p}(z)}}\nno\\
&&\quad +x_0^{-1}\delta\left(\frac{x_2-x_1}{-x_0}\right)
I(w_{(1)}\otimes y^{-L(0)}x_2^{-L(0)}Y_2(v, x_1)
y^{L(0)}x_2^{L(0)}w_{(2)})\lbar_{y=e^{-l_{p}(z)}}.
\end{eqnarray*}
Replacing $w_{(1)}$ by
$y^{-L(0)}x_2^{-L(0)}w_{(1)}\lbar_{y=e^{-l_{p}(z)}}$ and $w_{(2)}$ by
$y^{-L(0)}x_2^{-L(0)}w_{(2)}\lbar_{y=e^{-l_{p}(z)}}$, and then
applying $y^{L(0)}x_2^{L(0)}\lbar_{y=e^{-l_{p}(z)}}$ to the whole
equation, we obtain
\begin{eqnarray*}
\lefteqn{x^{-1}_0\delta\left(\frac{x_1-x_2}{x_0}\right)
Y_{3}(v,x_1)y^{L(0)}x_2^{L(0)}\cdot}\nno\\
&&\hspace{2em}\cdot I(y^{-L(0)}x_2^{-L(0)}w_{(1)}\otimes
y^{-L(0)}x_2^{-L(0)}w_{(2)})\lbar_{y=e^{-l_{p}(z)}}\nno\\
&&=x_2^{-1}\delta\left(\frac{x_1-x_0}{x_2}\right)
y^{L(0)}x_2^{L(0)}\cdot\nno\\
&&\hspace{2em}\cdot I(y^{-L(0)}x_2^{-L(0)}Y_1(v, x_0)w_{(1)}\otimes
y^{-L(0)}x_2^{-L(0)}w_{(2)})\lbar_{y=e^{-l_{p}(z)}}\nno\\
&&\quad +x_0^{-1}\delta\left(\frac{x_2-x_1}{-x_0}\right)
y^{L(0)}x_2^{L(0)}\cdot\nno\\
&&\hspace{2em}\cdot I(y^{-L(0)}x_2^{-L(0)}w_{(1)}\otimes
y^{-L(0)}x_2^{-L(0)}Y_2(v,x_1) w_{(2)})\lbar_{y=e^{-l_{p}(z)}}.
\end{eqnarray*}
But using (\ref{recover}), we can write this as
\begin{eqnarray*}
\lefteqn{x^{-1}_0\delta\left(\frac{x_1-x_2}{x_0}\right)
Y_{3}(v, x_1){\cal Y}_{I, p}(w_{(1)}, x_2)w_{(2)}}\nno\\
&&=x_2^{-1}\delta\left(\frac{x_1-x_0}{x_2}\right) {\cal Y}_{I,
p}( Y_1(v, x_0) w_{(1)}, x_2)w_{(2)}\nno\\
&&\hspace{2em}+x_0^{-1}\delta\left(\frac{x_2-x_1}{-x_0}\right)
{\cal Y}_{I, p}(w_{(1)}, x_2) Y_2(v, x_1) w_{(2)}.
\end{eqnarray*}
That is, the Jacobi identity for ${\cal Y}_{I, p}$ holds.

Similar procedures show that the ${\mathfrak s}{\mathfrak l}(2)$-bracket
relations for $I$ imply the ${\mathfrak s}{\mathfrak l}(2)$-bracket
relations for ${\cal Y}_{I, p}$, as follows: Let $j$ be $-1$,
$0$ or $1$. By multiplying (\ref{im:Lj}) by $(yx)^j$ and using
(\ref{log:xLx^}) we obtain
\begin{eqnarray*}
\lefteqn{(yx)^{-L(0)}L(j)(yx)^{L(0)}I(w_{(1)}\otimes w_{(2)})}\nn
&&=I(w_{(1)}\otimes (yx)^{-L(0)}L(j)(yx)^{L(0)}w_{(2)})\nn
&&\quad +\sum_{i=0}^{j+1} {j+1\choose i}z^i(yx)^i
I(((yx)^{-L(0)}L(j-i)(yx)^{L(0)}w_{(1)})\otimes w_{(2)})
\end{eqnarray*}
Replacing $w_{(1)}$ by $(yx)^{-L(0)}w_{(1)}$ and $w_{(2)}$ by
$(yx)^{-L(0)}w_{(2)}$, and then applying $(yx)^{L(0)}$ to the whole
equation, we obtain
\begin{eqnarray*}
\lefteqn{L(j)(yx)^{L(0)}I((yx)^{-L(0)}w_{(1)}\otimes
(yx)^{-L(0)}w_{(2)})}\\
&&\hspace{3em}=(yx)^{L(0)}I((yx)^{-L(0)}w_{(1)}\otimes
(yx)^{-L(0)}L(j)w_{(2)})\\
&&\hspace{4em}+\sum_{i=0}^{j+1} {j+1\choose i}z^i(yx)^i
(yx)^{L(0)}I(((yx)^{-L(0)}L(j-i)w_{(1)})\otimes(yx)^{-L(0)}w_{(2)}).
\end{eqnarray*}
Evaluating at $y=e^{-l_p(z)}$ and using (\ref{recover}) we see that
this gives exactly the ${\mathfrak s}{\mathfrak l}(2)$-bracket relations
(\ref{log:L(j)b}) for ${\cal Y}_{I, p}$.

Finally, we prove the $L(-1)$-derivative property for ${\cal
Y}_{I,p}$. This follows {}from (\ref{recover}), (\ref{log:dx^}), and the
${\mathfrak s}{\mathfrak l}(2)$-bracket relation with $j=0$ for ${\cal
Y}_{I, p}$, namely,
\[
[L(0),{\cal Y}_{I, p}(w_{(1)},x)]={\cal Y}_{I, p}(L(0)w_{(1)},x)+
x{\cal Y}_{I, p}(L(-1)w_{(1)},x),
\]
as follows:
\begin{eqnarray*}
\lefteqn{\frac{d}{dx}{\cal Y}_{I, p}(w_{(1)},x)w_{(2)}}\\
&&=\frac{d}{dx}e^{-l_p(z)L(0)}x^{L(0)}I(e^{l_p(z)L(0)}x^{-L(0)}w_{(1)}
\otimes e^{l_p(z)L(0)}x^{-L(0)}w_{(2)})\\
&&=e^{-l_p(z)L(0)}x^{-1}x^{L(0)}L(0)I(e^{l_p(z)L(0)}x^{-L(0)}w_{(1)}
\otimes e^{l_p(z)L(0)}x^{-L(0)}w_{(2)})\\
&&\hspace{1em}-e^{-l_p(z)L(0)}x^{L(0)}I(e^{l_p(z)L(0)}x^{-1}x^{-L(0)}
L(0)w_{(1)}\otimes e^{l_p(z)L(0)}x^{-L(0)}w_{(2)})\\
&&\hspace{1em}-e^{-l_p(z)L(0)}x^{L(0)}I(e^{l_p(z)L(0)}x^{-L(0)}w_{(1)}
\otimes e^{l_p(z)L(0)}x^{-1}x^{-L(0)}L(0)w_{(2)})\\
&&=x^{-1}L(0){\cal Y}_{I, p}(w_{(1)},x)w_{(2)}-x^{-1}{\cal Y}_{I, p}
(w_{(1)},x)L(0)w_{(2)}\\
&&\hspace{1em}-x^{-1}{\cal Y}_{I, p}(L(0)w_{(1)},x)w_{(2)}\\
&&={\cal Y}_{I, p}(L(-1)w_{(1)},x)w_{(2)}. \hspace{16em}\square
\end{eqnarray*}

\vspace{.1em}

\begin{rema}\label{mod-sub}{\rm
Given a generalized $V$-module $(W, Y_{W})$, recall {}from Remark
\ref{str-graded-g-mod-as-l-int} that $Y_{W}$ is a logarithmic
intertwining operator of type ${W\choose VW}$ not involving $\log x$
and having only integral powers of $x$. Then the substitution
$x\mapsto z$ in (\ref{log:IYp}) is very simple; it is independent of
$p$ and $Y_{W}(\cdot, z)\cdot$ entails only the substitutions
$x^{n}\mapsto z^{n}$ for $n\in \Z$. As a special case, we can take
$(W, Y_{W})$ to be $(V, Y)$ itself.}
\end{rema}

\begin{rema}{\rm
Let $I$ be a $P(z)$-intertwining map of type ${W_3\choose W_1\,W_2}$
and let $p,p' \in {\mathbb Z}$.  {}From (\ref{recover}), we see that the
logarithmic intertwining operators ${\cal Y}_{I, p}$ and ${\cal Y}_{I,
p'}$ of this same type differ as follows:
\begin{eqnarray}\label{YIp'YIp}
\lefteqn{{\cal Y}_{I,p'}(w_{(1)}, x)w_{(2)}}\nno\\
&&=e^{2\pi i(p-p')L(0)}{\cal Y}_{I,p}(e^{2\pi i(p'-p)L(0)}w_{(1)},
x)e^{2\pi i(p'-p)L(0)}w_{(2)}
\end{eqnarray}
for $w_{(1)}\in W_1$ and $w_{(2)}\in W_2$.  Using the notation in
Remark \ref{Ys1s2s3}, we thus have
\begin{eqnarray}
{\cal Y}_{I,p'}&=&({\cal Y}_{I,p})_{[p-p',p'-p,p'-p]}\nno\\
&=&{\cal Y}_{I,p}(\cdot,e^{2\pi i(p'-p)} \cdot).
\end{eqnarray}
}
\end{rema}

\begin{rema}\label{II1}{\rm
Let $I$ be a $P(z)$-intertwining map of type ${W_3}\choose
{W_1W_2}$. Then {}from the correspondence between $P(z)$-intertwining
maps and logarithmic intertwining operators in Proposition
\ref{im:correspond}, we see that for any nonzero complex number $z_1$,
the linear map $I_1$ defined by
\begin{equation}\label{log:zz_1}
I_1(w_{(1)}\otimes w_{(2)})=\sum_{n\in{\mathbb C}}\sum_{k\in {\mathbb N}}
{w_{(1)}}_{n;k}^{I,p}w_{(2)} e^{l_p(z_1)(-n-1)}(l_p(z_1))^k
\end{equation}
for $w_{(1)}\in W_1$ and $w_{(2)}\in W_2$ (recall (\ref{wInkw})) is a
$P(z_1)$-intertwining map of the same type.  In this sense,
${w_{(1)}}_{n;k}^{I,p}w_{(2)}$ is independent of $z$. We will hence
sometimes write $I(w_{(1)}\otimes w_{(2)})$ as
\begin{eqnarray}\label{imz}
I(w_{(1)},z)w_{(2)},
\end{eqnarray}
indicating that $z$ can be replaced by any nonzero complex number.
However, for a general intertwining map associated to a sphere with
punctures not necessarily of type $P(z)$, the corresponding element
${w_{(1)}}_{n;k}^{I,p}w_{(2)}$ will in general be different.  }
\end{rema}

We now proceed to the definition of the $P(z)$-tensor product.  As in
\cite{tensor1}, this will be a suitably universal ``$P(z)$-product.''
We generalize these notions {}from \cite{tensor1} using the notations
${\cal M}_{sg}$ and ${\cal GM}_{sg}$ (the categories of strongly
graded $V$-modules and generalized $V$-modules, respectively; recall
Notation \ref{MGM}) as follows:

\begin{defi}\label{pz-product}{\rm
Let ${\cal C}_1$ be either of the categories ${\cal M}_{sg}$ or ${\cal
GM}_{sg}$.  For $W_1, W_2\in \ob{\cal C}_1$, a {\em $P(z)$-product of
$W_1$ and $W_2$} is an object $(W_3,Y_3)$ of ${\cal C}_1$ equipped
with a $P(z)$-intertwining map $I_3$ of type ${W_3\choose
W_1\,W_2}$. We denote it by $(W_3,Y_3;I_3)$ or simply by
$(W_3;I_3)$. Let $(W_4,Y_4;I_4)$ be another $P(z)$-product of $W_1$
and $W_2$. A {\em morphism} {}from $(W_3,Y_3;I_3)$ to $(W_4,Y_4;I_4)$ is
a module map $\eta$ {}from $W_3$ to $W_4$ such that the
diagram
\begin{center}
\begin{picture}(100,60)
\put(-5,0){$\overline W_3$}
\put(13,4){\vector(1,0){104}}
\put(119,0){$\overline W_4$}
\put(41,50){$W_1\otimes W_2$}
\put(61,45){\vector(-3,-2){50}}
\put(68,45){\vector(3,-2){50}}
\put(65,8){$\bar\eta$}
\put(20,27){$I_3$}
\put(98,27){$I_4$}
\end{picture}
\end{center}
commutes, that is,
\begin{equation}
I_4=\bar\eta\circ I_3,
\end{equation}
where
\begin{equation}
\bar\eta : \overline{W}_3 \to \overline{W}_4
\end{equation}
is the natural map uniquely extending $\eta$.  (Note that $\bar\eta$
exists because $\eta$ preserves ${\mathbb C}$-gradings; we shall use
the notation $\bar\eta$ for any such map $\eta$.)
}
\end{defi}

\begin{rema}{\rm
In this setting, let $\eta$ be a morphism {}from $(W_3,Y_3;I_3)$ to
$(W_4,Y_4;I_4)$.  We know {}from (\ref{YIp})--(\ref{wInkw}) that for $p
\in {\mathbb Z}$, the coefficients ${w_{(1)}}_{n;k}^{I_{3},p}w_{(2)}$
and ${w_{(1)}}_{n;k}^{I_{4},p}w_{(2)}$ in the formal expansion
(\ref{wInkw}) of ${\cal Y}_{I_{3},p}(w_{(1)}, x)w_{(2)}$ and ${\cal
Y}_{I_{4},p}(w_{(1)}, x)w_{(2)}$, respectively, are determined by
$I_3$ and $I_4$, and that
\begin{equation}
\eta ({w_{(1)}}_{n;k}^{I_{3},p}w_{(2)}) =
{w_{(1)}}_{n;k}^{I_{4},p}w_{(2)},
\end{equation}
as we see by applying $\bar\eta$ to (\ref{recover}).
}
\end{rema}

The notion of $P(z)$-tensor product is now defined by means of a
universal property as follows:

\begin{defi}\label{pz-tp}{\rm
Let ${\cal C}$ be a full subcategory of either ${\cal M}_{sg}$ or
${\cal GM}_{sg}$.  For $W_1, W_2\in \ob{\cal C}$, a {\em $P(z)$-tensor
product of $W_1$ and $W_2$ in ${\cal C}$} is a $P(z)$-product $(W_0,
Y_0; I_0)$ with $W_0\in{\rm ob\,}{\cal C}$ such that for any
$P(z)$-product $(W,Y;I)$ with $W\in{\rm ob\,}{\cal C}$, there is a
unique morphism {}from $(W_0, Y_0; I_0)$ to $(W,Y;I)$.  Clearly, a
$P(z)$-tensor product of $W_1$ and $W_2$ in ${\cal C}$, if it exists,
is unique up to unique isomorphism. In this case we will denote it by
\[
(W_1\boxtimes_{P(z)} W_2, Y_{P(z)}; \boxtimes_{P(z)})
\]
and call the object $(W_1\boxtimes_{P(z)} W_2, Y_{P(z)})$ the {\em
$P(z)$-tensor product module of $W_1$ and $W_2$ in ${\cal C}$}.  We
will skip the phrase ``in ${\cal C}$'' if the category ${\cal C}$
under consideration is clear in context.  }
\end{defi}

\begin{rema}{\rm
Consider the functor {}from ${\cal C}$ to the category ${\bf Set}$
defined by assigning to $W\in {\rm ob\,}{\cal C}$ the set ${\cal
M}^W_{W_1\,W_2}$ of all $P(z)$-intertwining maps of type ${W\choose
W_1\,W_2}$.  Then if the $P(z)$-tensor product of $W_1$ and $W_2$
exists, it is just the universal element for this functor, and this
functor is representable, represented by the $P(z)$-tensor product.
(Recall that given a functor $f$ {}from a category ${\cal K}$ to ${\bf
Set}$, a universal element for $f$, if it exists, is a pair $(X,x)$
where $X\in{\rm ob\,}{\cal K}$ and $x\in f(X)$ such that for any pair
$(Y,y)$ with $Y\in {\rm ob\,}{\cal K}$ and $y\in f(Y)$, there is a
unique morphism $\sigma: X\to Y$ such that $f(\sigma)(x)=y$; in this
case, $f$ is represented by $X$.)  }
\end{rema}

Definition \ref{pz-tp} and Proposition \ref{im:correspond} immediately
give the following result relating the module maps {}from a
$P(z)$-tensor product module with the $P(z)$-intertwining maps and the
logarithmic intertwining operators:

\begin{propo}\label{pz-iso}
Suppose that $W_1\boxtimes_{P(z)}W_2$ exists. We have a natural
isomorphism
\begin{eqnarray}\label{isofromhomstointwmaps}
\hom_{V}(W_1\boxtimes_{P(z)}W_2, W_3)&\stackrel{\sim}{\to}&
{\cal M}^{W_3}_{W_1W_2}\nno\\
\eta&\mapsto& \overline{\eta}\circ
\boxtimes_{P(z)}
\end{eqnarray}
and for $p\in {\mathbb Z}$, a natural isomorphism
\begin{eqnarray}
\hom_{V}(W_1\boxtimes_{P(z)} W_2, W_3)&
\stackrel{\sim}{\rightarrow}& {\cal V}^{W_3}_{W_1W_2}\nno\\
\eta&\mapsto & {\cal Y}_{\eta, p}
\end{eqnarray}
where ${\cal Y}_{\eta, p}={\cal Y}_{I, p}$ with
$I=\overline{\eta}\circ \boxtimes_{P(z)}$.\epf
\end{propo}

Suppose that the $P(z)$-tensor product $(W_1\boxtimes_{P(z)} W_2,
Y_{P(z)}; \boxtimes_{P(z)})$ of $W_1$ and $W_2$ exists.  We will
sometimes denote the action of the canonical $P(z)$-intertwining map
\begin{equation}\label{actionofboxtensormap}
w_{(1)} \otimes w_{(2)}\mapsto\boxtimes_{P(z)}(w_{(1)} \otimes
w_{(2)})=\boxtimes_{P(z)}(w_{(1)}, z)w_{(2)} \in 
\overline{W_1\boxtimes_{P(z)}W_2}
\end{equation}
(recall (\ref{imz})) on elements simply by $w_{(1)}\boxtimes_{P(z)}
w_{(2)}$:
\begin{equation}\label{boxtensorofelements}
w_{(1)}\boxtimes_{P(z)} w_{(2)}=\boxtimes_{P(z)}(w_{(1)} \otimes
w_{(2)})=\boxtimes_{P(z)}(w_{(1)}, z)w_{(2)}.
\end{equation}

\begin{rema} {\rm
We emphasize that the element $w_{(1)}\boxtimes_{P(z)} w_{(2)}$
defined here is an element of the formal completion
$\overline{W_1\boxtimes_{P(z)}W_2}$, and {\it not} (in general) of the
module $W_1\boxtimes_{P(z)}W_2$ itself.  This is different {}from the
classical case for modules for a Lie algebra (recall Section
\ref{LA}), where the tensor product of elements of two modules is an
element of the tensor product module.}
\end{rema}

\begin{rema} {\rm
Note that under the natural isomorphism (\ref{isofromhomstointwmaps})
for the case $W_{3} = W_1\boxtimes_{P(z)}W_2$, the identity map {}from
$W_1\boxtimes_{P(z)}W_2$ to itself corresponds to the canonical
intertwining map $\boxtimes_{P(z)}$.  Furthermore, for $p \in {\mathbb
Z}$, the $P(z)$-tensor product of $W_1$ and $W_2$ gives rise to a
logarithmic intertwining operator ${\cal Y}_{\boxtimes_{P(z)}, p}$ of
type ${W_1\boxtimes_{P(z)}W_2\choose W_1\,W_2}$, according to formula
(\ref{recover}).  If $p$ is changed to $p' \in {\mathbb Z}$, this
logarithmic intertwining operator changes according to
(\ref{YIp'YIp}).  Note that the $P(z)$-intertwining map
$\boxtimes_{P(z)}$ is canonical and depends only on $z$, while a
corresponding logarithmic intertwining operator is not; it depends on
$p \in {\mathbb Z}$.
}
\end{rema}

\begin{rema} {\rm
Sometimes it will be convenient, as in the next proposition, to use
the particular isomorphism associated with $p=0$ (in Proposition
\ref{im:correspond}) between the spaces of $P(z)$-intertwining maps
and of logarithmic intertwining operators of the same type.  In this
case, we shall sometimes simplify the notation by dropping the $p$
($=0$) in the notation ${w_{(1)}}_{n;k}^{I,0}w_{(2)}$ (recall
\ref{wInkw})):
\begin{equation}
{w_{(1)}}_{n;k}^{I}w_{(2)} = {w_{(1)}}_{n;k}^{I,0}w_{(2)}.
\end{equation}
}
\end{rema}

\begin{propo}\label{4.19}
Suppose that the $P(z)$-tensor product $(W_1\boxtimes_{P(z)} W_2,
Y_{P(z)}; \boxtimes_{P(z)})$ of $W_1$ and $W_2$ in ${\cal C}$ exists.
Then for any complex number $z_1\neq 0$, the $P(z_1)$-tensor product
of $W_1$ and $W_2$ in ${\cal C}$ also exists, and is given by
$(W_1\boxtimes_{P(z)} W_2, Y_{P(z)}; \boxtimes_{P(z_1)})$, where the
$P(z_1)$-intertwining map $\boxtimes_{P(z_1)}$ is defined by
\begin{equation}\label{tpzz_1}
\boxtimes_{P(z_1)}(w_{(1)}\otimes w_{(2)})=\sum_{n\in{\mathbb
C}}\sum_{k\in {\mathbb N}} {w_{(1)}}_{n;k}^{\boxtimes_{P(z)}}w_{(2)}
e^{\log z_1 (-n-1)}(\log z_1)^k
\end{equation}
for $w_{(1)}\in W_1$ and $w_{(2)}\in W_2$.
\end{propo}
\pf By Remark \ref{II1}, (\ref{tpzz_1}) indeed defines a
$P(z_1)$-product.  Given any $P(z_1)$-product $(W_3, Y_3; I_1)$ of
$W_1$ and $W_2$, let $I$ be the $P(z)$-product related to $I_1$ by
formula (\ref{log:zz_1}) with $I_1$, $I$ and $z_1$ in (\ref{log:zz_1})
replaced by $I$, $I_1$ and $z$, respectively.  Then {}from the
definition of $P(z)$-tensor product, there is a unique morphism $\eta$
{}from $(W_1\boxtimes_{P(z)} W_2, Y_{P(z)}; \boxtimes_{P(z)})$ to $(W_3,
Y_3; I)$.  Thus by (\ref{boxtensorofelements}) and (\ref{tpzz_1}) we
see that $\eta$ is also a morphism {}from the $P(z_1)$-product
$(W_1\boxtimes_{P(z)} W_2, Y_{P(z)}; \boxtimes_{P(z_1)})$ to $(W_3,
Y_3; I_1)$. The uniqueness of such a morphism follows similarly {}from
the uniqueness of a morphism {}from $(W_1\boxtimes_{P(z)} W_2, Y_{P(z)};
\boxtimes_{P(z)})$ to $(W_3, Y_3; I)$.  Hence $(W_1\boxtimes_{P(z)}
W_2, Y_{P(z)}; \boxtimes_{P(z_1)})$ is the $P(z_1)$-tensor product of
$W_1$ and $W_2$.  \epfv

\begin{rema}\label{intwmapdependsongeomdata} {\rm
In general, it will turn out that the existence of tensor product, and
the tensor product module itself, do not depend on the geometric data.
It is the intertwining map {}from the two modules to the completion of
their tensor product that encodes the geometric information.
}
\end{rema}

Generalizing Lemma 4.9 of \cite{tensor4}, we have:

\begin{propo}\label{span}
The module $W_1\boxtimes_{P(z)}W_2$ (if it exists) is spanned (as a
vector space) by the (generalized-) weight components of the elements
of $\overline{W_1\boxtimes_{P(z)}W_2}$ of the form
$w_{(1)}\boxtimes_{P(z)} w_{(2)}$, for all $w_{(1)}\in W_1$ and
$w_{(2)}\in W_2$.
\end{propo}

\pf Denote by $W_0$ the vector subspace of $W_1\boxtimes_{P(z)}W_2$
spanned by all the weight components of all the elements of
$\overline{W_1\boxtimes_{P(z)} W_2}$ of the form
$w_{(1)}\boxtimes_{P(z)}w_{(2)}$ for $w_{(1)}\in W_1$ and $w_{(2)}\in
W_2$.  For a homogeneous vector $v\in V$ and arbitrary elements
$w_{(1)}\in W_1$ and $w_{(2)}\in W_2$, equating the
$x_0^{-1}x_1^{-m-1}$ coefficients of the Jacobi identity
(\ref{im:def}) gives
\begin{equation}\label{elm}
v_m({w_{(1)}}\boxtimes_{P(z)} w_{(2)})={w_{(1)}}\boxtimes_{P(z)} (v_m
w_{(2)})+\sum_{i\in {\mathbb N}}{m\choose i}z^{m-i}(v_i
w_{(1)})\boxtimes_{P(z)} w_{(2)}
\end{equation}
for all $m\in {\mathbb Z}$.  Note that the summation in the right-hand
side of (\ref{elm}) is always finite. Hence by taking arbitrary weight
components of (\ref{elm}) we see that $W_0$ is closed under the action
of $V$.  In case $V$ is M\"obius, a similar argument, using
(\ref{im:Lj}), shows that $W_0$ is stable under the action of
${\mathfrak s}{\mathfrak l}(2)$.  It is clear that $W_0$ is ${\mathbb
C}$-graded and $\tilde{A}$-graded.  Thus $W_0$ is a submodule of
$W_1\boxtimes_{P(z)} W_2$.  Now consider the quotient module
\[
W=(W_1\boxtimes_{P(z)} W_2)/W_0
\]
and let $\pi_{W}$ be the canonical map {}from $W_1\boxtimes_{P(z)}
W_2$ to $W$.  By the definition of $W_0$, we have
\[
\overline{\pi}_{W}\circ \boxtimes_{P(z)}=0,
\]
using the notation (\ref{actionofboxtensormap}).  The universal
property of the $P(z)$-tensor product then demands that $\pi_W=0$,
i.e., that $W_0=W_1\boxtimes_{P(z)} W_2$.  \epfv

It is clear {}from Definition \ref{pz-tp} that the tensor product
operation distributes over direct sums in the following sense:

\begin{propo}\label{tensorproductdistributes}
For $U_1, \dots, U_{k}$, $W_1, \dots, W_{l}\in \ob{\cal C}$,
suppose that each $U_{i}\boxtimes_{P(z)}W_{j}$ exists. Then
$(\coprod_{i}U_{i})\boxtimes_{P(z)}(\coprod_{j}W_{j})$ exists and
there is a natural isomorphism
\[
\biggl(\coprod_{i}U_{i}\biggr)\boxtimes_{P(z)}
\biggl(\coprod_{j}W_{j}\biggr)\stackrel{\sim}
{\rightarrow} \coprod_{i,j}U_{i}\boxtimes_{P(z)}W_{j}.\hspace{4em}\square
\]
\end{propo}

\begin{rema}\label{bifunctor} {\rm
It is of course natural to view the $P(z)$-tensor product as a
bifunctor: Suppose that ${\cal C}$ is a full subcategory of either
${\cal M}_{sg}$ or ${\cal GM}_{sg}$ (recall Notation \ref{MGM}) such
that for all $W_{1}, W_{2}\in \ob{\cal C}$, the $P(z)$-tensor product
of $W_{1}$ and $W_{2}$ exists in ${\cal C}$.  Then $\boxtimes_{P(z)}$
provides a (bi)functor
\begin{equation}
\boxtimes_{P(z)}:{\cal C} \times  {\cal C} \rightarrow {\cal C}
\end{equation}
as follows: For $W_{1}, W_{2}\in \ob{\cal C}$,
\begin{equation}
\boxtimes_{P(z)}(W_{1}, W_{2}) = W_1\boxtimes_{P(z)} W_2 \in \ob{\cal C}
\end{equation}
and for $V$-module maps
\begin{equation}
\sigma_1 : W_1 \rightarrow W_3,
\end{equation}
\begin{equation}
\sigma_2 : W_2 \rightarrow W_4
\end{equation}
with $W_{3}, W_{4}\in \ob{\cal C}$, we have the $V$-module map,
denoted
\begin{equation}
\boxtimes_{P(z)}(\sigma_{1}, \sigma_{2}) = \sigma_1\boxtimes_{P(z)}
\sigma_2,
\end{equation}
{}from $W_1\boxtimes_{P(z)} W_2$ to $W_3\boxtimes_{P(z)} W_4$, defined
by the universal property of the $P(z)$-tensor product
$W_1\boxtimes_{P(z)} W_2$ and the fact that the composition of
$\boxtimes_{P(z)}$ with $\sigma_1 \otimes \sigma_2$ is a
$P(z)$-intertwining map
\begin{equation}
\boxtimes_{P(z)} \circ (\sigma_{1} \otimes \sigma_{2}):W_1\otimes W_2
\rightarrow \overline{W_3\boxtimes_{P(z)} W_4}.
\end{equation}
Note that it is the effect of this bifunctor on morphisms (rather than
on objects) that exhibits the role of the geometric data.
}
\end{rema}

We now discuss the simplest examples of $P(z)$-tensor products---those
in which one or both of $W_{1}$ or
$W_{2}$ is $V$ itself (viewed as a (generalized) $V$-module); 
we suppose here that $V\in \ob \mathcal{C}$. Since the discussion of the 
case in which both
$W_{1}$ and $W_{2}$ are $V$ turns out to be no simpler than the case 
in which $W_{1}=V$,
we shall discuss only the two more general 
cases $W_{1}=V$ and $W_{2}=V$. 

\begin{exam}\label{expl-vw}
{\rm Let $(W, Y_{W})$ be an object of $\mathcal{C}$. 
The vertex operator map $Y_{W}$ gives a $P(z)$-intertwining 
map 
\[
I_{Y_{W}, p}=Y_{W}(\cdot, z)\cdot: V\otimes W\to \overline{W}
\]
for any fixed $p\in \Z$ (recall Proposition \ref{im:correspond} and Remark 
\ref{mod-sub}). We claim that 
$(W, Y_{W}; Y_{W}(\cdot, z)\cdot)$ is the $P(z)$-tensor product of $V$ and 
$W$ in $\mathcal{C}$. In fact, let $(W_{3}, Y_{3}; I)$ be a $P(z)$-product 
of $V$ and $W$ in $\mathcal{C}$ and suppose that there exists 
a module map $\eta: W\to W_{3}$ such that 
\begin{equation}\label{v-tensor-w-1}
\overline{\eta}\circ (Y_{W}(\cdot, z)\cdot)=I.
\end{equation}
Then for $w\in W$, we must have 
\begin{eqnarray}\label{v-tensor-w-2}
\eta(w)&=&\eta(Y_{W}(\mathbf{1}, z)w)\nn
&=&(\overline{\eta}\circ (Y_{W}(\cdot, z)\cdot))(\mathbf{1}\otimes w)\nn
&=&I(\mathbf{1}\otimes w),
\end{eqnarray}
so that $\eta$ is unique if it exists.
We now define $\eta: W\to \overline{W}_{3}$ using (\ref{v-tensor-w-2}). 
We shall show that $\eta(W)\subset W_{3}$ and that $\eta$ has the 
desired properties.
Since 
$I$ is a $P(z)$-intertwining map of type ${W_{3}\choose VW}$, it corresponds 
to a logarithmic intertwining 
operator $\mathcal{Y}=\mathcal{Y}_{I, p}$ of the same type, according to 
Proposition \ref{im:correspond}. Since $L(-1)\mathbf{1}=0$, 
we have 
\[
\frac{d}{dx}\mathcal{Y}(\mathbf{1}, x)=\mathcal{Y}(L(-1)\mathbf{1}, x)=0.
\]
Thus $\mathcal{Y}(\mathbf{1}, x)$ 
is simply the constant map $\mathbf{1}_{-1; 0}^{\mathcal{Y}}: W\to W_{3}$ 
(using the notation (\ref{log:map})), and this map preserves (generalized) 
weights, by Proposition \ref{log:logwt}(b).  By
Proposition 
\ref{im:correspond}, $I=I_{\mathcal{Y}, p}$, so that
\begin{eqnarray*}
\eta(w)&=&I(\mathbf{1}\otimes w)\nn
&=&I_{\mathcal{Y}, p}(\mathbf{1}\otimes w)\nn
&=&\mathbf{1}_{-1; 0}^{\mathcal{Y}}w
\end{eqnarray*}
for $w\in W$.
So $\eta=\mathbf{1}_{-1; 0}^{\mathcal{Y}}$ is a linear map {}from 
$W$ to $W_{3}$ preserving (generalized) weights. Using the 
Jacobi identity for the $P(z)$-intertwining map $I$ and the 
fact that $Y(u, x_{0})\mathbf{1}\in V[[x_{0}]]$ for $u\in V$, we obtain
\begin{eqnarray*}
\eta(Y_{W}(u, x)w)&=&I(\mathbf{1} \otimes Y_{W}(u, x)w)\nn
&=&Y_{3}(u, x)I(\mathbf{1} \otimes w)
-\res_{x_{0}}z^{-1}\delta\left(\frac{x-x_{0}}{z}\right)
I(Y(u, x_{0})\mathbf{1} \otimes w)\nn
&=&Y_{3}(u, x)I(\mathbf{1} \otimes w)\nn
&=&Y_{3}(u, x)\eta(w)
\end{eqnarray*}
for $u\in V$ and $w\in W$, proving that
$\eta$ is a module map. 
For $w\in W$, 
\begin{eqnarray}\label{v-tensor-w-3}
(\overline{\eta}\circ (Y_{W}(\cdot, z)\cdot))(\mathbf{1}\otimes w)&=&
\overline{\eta}(Y_{W}(\mathbf{1}, z)w)\nn
&=&\eta(w)\nn
&=&I(\mathbf{1}\otimes w).
\end{eqnarray}
Using the Jacobi identity for $P(z)$-intertwining maps, we obtain
\begin{eqnarray}\label{int-recurrence-rel}
\lefteqn{I(Y(u, x_{0})v\otimes w)}\nn
&&=\res_{x}x_{0}^{-1}\delta\left(\frac{x-z}{x_{0}}\right)
Y_3(u, x)I(v\otimes w)-\res_{x}x_{0}^{-1}\delta\left(\frac{z-x}{-x_{0}}\right)
I(v\otimes Y_W(u, x)w) \;\;\;\;\;\;\;
\end{eqnarray}
for $u, v\in V$ and $w\in W$. 
Since $\eta$ is a module map and $Y_{W}(\cdot, z)\cdot$ is a $P(z)$-intertwining 
map of type ${W\choose VW}$, $\overline{\eta}\circ Y_{W}(\cdot, z)\cdot$
is a $P(z)$-intertwining map of  type ${W_{3}\choose VW}$. In particular, 
(\ref{int-recurrence-rel}) holds when we replace $I$ by 
$\overline{\eta}\circ Y_{W}(\cdot, z)\cdot$. Using
(\ref{int-recurrence-rel}) for $v=\mathbf{1}$ together with (\ref{v-tensor-w-3}), 
we obtain
\[
(\overline{\eta}\circ (Y_{W}(\cdot, z)\cdot))(u\otimes w)=I(u\otimes w)
\]
for $u\in V$ and $w\in W$, proving (\ref{v-tensor-w-1}), as desired.
Thus $(W, Y_{W}; Y_{W}(\cdot, z)\cdot)$ is the $P(z)$-tensor product
of $V$ and $W$ in $\mathcal{C}$. }
\end{exam}

\begin{exam}
{\rm Let $(W, Y_{W})$ be an object of $\mathcal{C}$. In order to
construct the $P(z)$-tensor product $W\boxtimes_{P(z)} V$, recall {}from
(\ref{Omega_r}) and Proposition \ref{log:omega} that
$\Omega_{p}(Y_{W})$ is a logarithmic intertwining operator of type
${W\choose WV}$. It involves only integral powers of the formal
variable and no logarithms, and it is independent of $p$.  In fact,
\[
\Omega_{p}(Y_{W})(w, x)v=e^{xL(-1)}Y_{W}(v, -x)w
\]
for $v\in V$ and $w\in W$. For $q\in \Z$, 
\[
I_{\Omega_{p}(Y_{W}), q}
=\Omega_{p}(Y_{W})(\cdot, z)\cdot: W\otimes 
V\to \overline{W}
\]
is a $P(z)$-intertwining map of the same type and is independent of
$q$.  We claim that $(W, Y_{W}; \Omega_{p}(Y_{W})(\cdot, z)\cdot)$ is
the $P(z)$-tensor product of $W$ and $V$ in $\mathcal{C}$. In fact,
let $(W_{3}, Y_{3}; I)$ be a $P(z)$-product of $W$ and $V$ in
$\mathcal{C}$ and suppose that there exists a module map $\eta: W\to
W_{3}$ such that
\begin{equation}\label{w-tensor-v-1}
\overline{\eta}\circ \Omega_{p}(Y_{W})(\cdot, z)\cdot=I.
\end{equation}
For $w\in W$, we must have 
\begin{eqnarray}\label{w-tensor-v-3}
\eta(w)&=&\eta(Y_{W}(\mathbf{1}, -z)w)\nn
&=&e^{-zL(-1)}\overline{\eta}(e^{zL(-1)}Y_{W}(\mathbf{1}, -z)w)\nn
&=&e^{-zL(-1)}\overline{\eta}(\Omega_{p}(Y_{W})(w, z)\mathbf{1}))\nn
&=&e^{-zL(-1)}(\overline{\eta}\circ (\Omega_{p}(Y_{W})(\cdot, z)\cdot))
(w \otimes \mathbf{1})\nn
&=&e^{-zL(-1)}I(w\otimes \mathbf{1}),
\end{eqnarray}
and so $\eta$ is unique if it exists.
We now define $\eta: W\to \overline{W}_{3}$ by (\ref{w-tensor-v-3}).
Consider the logarithmic intertwining operator $\mathcal{Y}=
\mathcal{Y}_{I, q}$ that corresponds to $I$ by Proposition  
\ref{im:correspond}. 
Using Proposition \ref{im:correspond},
(\ref{branch1})--(\ref{log:fsub}),
(\ref{log:subs})
and the equality
\begin{eqnarray*}
l_{q}(-z)&=&\log |-z|+i(\arg (-z)+2\pi q)\nn
&=&\left\{\begin{array}{ll}
\log |z|+i(\arg z+\pi+2\pi q),&0\le \arg z< \pi\\
\log |z|+i(\arg z-\pi+2\pi q),&\pi\le \arg z< 2\pi
\end{array}\right.\nn
&=&\left\{\begin{array}{ll}
l_{q}(z)+\pi i,&0\le \arg z< \pi\\
l_{q}(z)-\pi i,&\pi\le \arg z< 2\pi,
\end{array}\right.\nn
\end{eqnarray*}
we have
\begin{eqnarray*}
e^{-zL(-1)}I(w\otimes \mathbf{1})
&=&e^{-zL(-1)}\mathcal{Y}(w, e^{l_{q}(z)})\mathbf{1}\nn
&=&e^{-xL(-1)}\mathcal{Y}(w, x)\mathbf{1}|_{x^{n}=e^{nl_{q}(z)},\; 
(\log x)^{m}=(l_{q}(z))^{m},\; n\in \C, \;m\in \N}\nn
&=&e^{yL(-1)}\mathcal{Y}(w, e^{\pm \pi i}y)\mathbf{1}
|_{y^{n}=e^{nl_{q}(-z)},\; (\log y)^{m}=(l_{q}(-z))^{m},\;
n\in \C, \; m\in \N},
\end{eqnarray*}
where $e^{\pm \pi i}$ is $e^{-\pi i}$ when $0\le \arg z<\pi$ and is 
$e^{\pi i}$ when $\pi\le \arg z<2 \pi$. Then by (\ref{Omega_r}),
we see that 
$\eta(w)=e^{-zL(-1)}I(w\otimes \mathbf{1})$ is equal to 
$\Omega_{-1}(\mathcal{Y})(\mathbf{1}, e^{l_{q}(-z)})w$ 
when $0\le \arg z<\pi$ and 
is equal to $\Omega_{0}(\mathcal{Y})(\mathbf{1}, e^{l_{q}(-z)})w$ when 
$\pi\le \arg z<2 \pi$. By Proposition \ref{log:omega}, 
$\Omega_{-1}(\mathcal{Y})$ and $\Omega_{0}(\mathcal{Y})$ are
logarithmic intertwining operators of type ${W_{3}\choose VW}$. 
As in Example \ref{expl-vw}, we see that 
$\Omega_{-1}(\mathcal{Y})(\mathbf{1}, y)$ 
and $\Omega_{0}(\mathcal{Y})(\mathbf{1}, y)$ are equal to 
$\mathbf{1}_{-1, 0}^{\Omega_{-1}(\mathcal{Y})}$ and
$\mathbf{1}_{-1, 0}^{\Omega_{0}(\mathcal{Y})}$, respectively, and 
these maps preserve (generalized) weights. Therefore $\eta$ 
is a linear map {}from $W$ to $W_{3}$ preserving (generalized) weights.
Using the Jacobi identity for the $P(z)$-intertwining map $I$ and the 
fact that $Y(u, x_{1})\mathbf{1}\in V[[x_{1}]]$, we have
\begin{eqnarray*}
\eta(Y_{W}(u, x_{0})w)&=&e^{-zL(-1)}I(Y_{W}(u, x_{0})w\otimes \mathbf{1})\nn
&=&\res_{x_{1}}x_{0}^{-1}\delta\left(\frac{x_{1}-z}{x_{0}}\right)
e^{-zL(-1)}Y_{3}(u, x_{1})I(w\otimes \mathbf{1})\nn
&& -\res_{x_{1}}x_{0}^{-1}\delta\left(\frac{z-x_{1}}{-x_{0}}\right)
e^{-zL(-1)}I(w\otimes Y(u, x_{1})\mathbf{1})\nn
&=&e^{-zL(-1)}Y_{3}(u, x_{0}+z)I(w\otimes \mathbf{1})\nn
&=&Y_{3}(u, x_{0})e^{-zL(-1)}I(w\otimes \mathbf{1})\nn
&=&Y_{3}(u, x_{0})\eta(w)
\end{eqnarray*}
for $u\in V$ and $w\in W$, proving that $\eta$ is a module map.
For $w\in W$, 
\begin{eqnarray}\label{w-tensor-v-4}
(\overline{\eta}\circ (\Omega_{p}(Y_{W})(\cdot, z)\cdot))(w \otimes \mathbf{1})
&=&\overline{\eta}(e^{zL(-1)}Y_{W}(\mathbf{1}, -z)w)\nn
&=&e^{zL(-1)}\eta(w)\nn
&=&e^{zL(-1)}e^{-zL(-1)}I(w\otimes \mathbf{1})\nn
&=&I(w\otimes \mathbf{1}).
\end{eqnarray}
Since both $\overline{\eta}\circ (\Omega_{p}(Y_{W})(\cdot, z)\cdot)$ and $I$ 
are $P(z)$-intertwining maps of type ${W_{3}\choose WV}$,
using the Jacobi identity for $P(z)$-intertwining operators
and (\ref{w-tensor-v-4}) (cf. Example \ref{expl-vw}), we have
\[
(\overline{\eta}\circ (\Omega_{p}(Y_{W})(\cdot, z)\cdot))
(w\otimes v)=I(w\otimes v)
\]
for $v\in V$ and $w\in W$, proving (\ref{w-tensor-v-1}). 
Thus $(W, Y_{W}; \Omega_{p}(Y_{W})(\cdot, z)\cdot)$  is the 
$P(z)$-tensor product of $W$ and 
$V$ in $\mathcal{C}$.}
\end{exam}

We discussed the important special class of finitely reductive vertex
operator algebras in the Introduction.  In case $V$ is a finitely
reductive vertex operator algebra, the $P(z)$-tensor product always
exists, as we are about to establish (following \cite{tensor1} and
\cite{tensor3}).  As in the Introduction, the definition of finite
reductivity is:

\begin{defi}\label{finitelyreductive}{\rm
A vertex operator algebra $V$ is {\it finitely reductive} if
\begin{enumerate}
\item Every $V$-module is completely reducible.
\item There are only finitely many irreducible $V$-modules (up to
equivalence).
\item All the fusion rules (the dimensions of the spaces of
intertwining operators among triples of modules) for $V$ are finite.
\end{enumerate}
}
\end{defi}

\begin{rema} {\rm
In this case, every $V$-module is of course a {\it finite} direct sum
of irreducible modules.  Also, the third condition holds if the
finiteness of the fusion rules among triples of only {\it irreducible}
modules is assumed.
}
\end{rema}

\begin{rema} {\rm
We are of course taking the notion of $V$-module so that the grading
restriction conditions are the ones described in Remark
\ref{moduleswiththetrivialgroup}, formulas (\ref{Wn+k=0}) and
(\ref{dimWnfinite}); in particular, $V$-modules are understood to be
$\C$-graded.  Recall {}from Remark \ref{congruent} that for an
irreducible module, all its weights are congruent to one another
modulo $\Z$.  Thus for an irreducible module, our grading-truncation
condition (\ref{Wn+k=0}) amounts exactly to the condition that the
real parts of the weights are bounded {}from below.  In
\cite{tensor1}--\cite{tensor3}, boundedness of the real parts of the
weights {}from below was our grading-truncation condition in the
definition of the notion of module for a vertex operator algebra.
Thus the first two conditions in the notion of finite reductivity are
the same whether we use the current grading restriction conditions in
the definition of the notion of module or the corresponding conditions
in \cite{tensor1}--\cite{tensor3}.  As for intertwining operators,
recall {}from Remark \ref{ordinaryandlogintwops} and Corollary
\ref{powerscongruentmodZ} that when the first two conditions are
satisfied, the notion of (ordinary, non-logarithmic) intertwining
operator here coincides with that in \cite{tensor1} because the
truncation conditions agree.  Also, in this setting, by Remark
\ref{log:ordi}, the logarithmic and ordinary intertwining operators
are the same, and so the spaces of intertwining operators ${\cal
V}^{W_3}_{W_1\,W_2}$ and fusion rules $N^{W_3}_{W_1\,W_2}$ in
Definition \ref{fusionrule} have the same meanings as in
\cite{tensor1}.  Thus the notion of finite reductivity for a vertex
operator algebra is the same whether we use the current grading
restriction and truncation conditions in the definitions of the
notions of module and of intertwining operator or the corresponding
conditions in \cite{tensor1}--\cite{tensor3}.  In particular, finite
reductivity of $V$ according to Definition \ref{finitelyreductive} is
equivalent to the corresponding notion, ``rationality'' (recall the
Introduction) in \cite{tensor1}--\cite{tensor3}.  }
\end{rema}

\begin{rema} {\rm
For a vertex operator algebra $V$ (in particular, a finitely reductive
one), the category ${\cal M}$ of $V$-modules coincides with the
category ${\cal M}_{sg}$ of strongly graded $V$-modules; recall
Notation \ref{MGM}.
}
\end{rema}

For the rest of Section 4.1, let us assume that $V$ is a finitely
reductive vertex operator algebra.  We shall now show that
$P(z)$-tensor products always exist in the category ${\cal M}$ ($=
{\cal M}_{sg}$) of $V$-modules, in the sense of Definition
\ref{pz-tp}.

Consider $V$-modules $W_{1}$, $W_{2}$ and $W_{3}$.  We know that
\begin{eqnarray}
N^{W_3}_{W_1\,W_2} = \dim {\cal V}^{W_3}_{W_1\,W_2} < \infty
\end{eqnarray}
and {}from Proposition \ref{im:correspond}, we also have
\begin{eqnarray}
N^{W_3}_{W_1\,W_2} = \dim {\cal M}[P(z)]_{W_{1}W_{2}}^{W_{3}} = \dim
{\cal M}_{W_{1}W_{2}}^{W_{3}} < \infty
\end{eqnarray}
(recall Definition \ref{im:imdef}).

The natural evaluation map
\begin{eqnarray}
W_{1}\otimes W_{2}\otimes {\cal M}^{W_{3}}_{W_{1}W_{2}}&\to& \overline{W}_{3}
\nno\\
w_{(1)}\otimes w_{(2)}\otimes I&\mapsto& I(w_{(1)}\otimes w_{(2)})
\end{eqnarray}
gives a natural map
\begin{equation}
{\cal F}[P(z)]^{W_{3}}_{W_{1}W_{2}}: W_{1}\otimes W_{2}\to 
\mbox{\rm Hom}({\cal M}_{W_{1}W_{2}}^{W_{3}}, \overline{W}_{3})=
({\cal M}^{W_{3}}_{W_{1}W_{2}})^{*}\otimes \overline{W}_{3}.
\end{equation}
Since $\dim {\cal M}_{W_{1}W_{2}}^{W_{3}} < \infty$, $({\cal
M}^{W_{3}}_{W_{1}W_{2}})^{*}\otimes W_{3}$ is a $V$-module (with
finite-dimensional weight spaces) in the obvious way, and the map
${\cal F}[P(z)]^{W_{3}}_{W_{1}W_{2}}$ is clearly a $P(z)$-intertwining
map, where we make the identification
\begin{equation}
({\cal M}^{W_{3}}_{W_{1}W_{2}})^{*}\otimes \overline{W}_{3}
=\overline{({\cal M}^{W_{3}}_{W_{1}W_{2}})^{*}\otimes W_{3}}.
\end{equation}
This gives us a natural $P(z)$-product for the category ${\cal M} =
{\cal M}_{sg}$ (recall Definition \ref{pz-product}).  Moreover, we
have a natural linear injection
\begin{eqnarray}
i: {\cal M}^{W_{3}}_{W_{1}W_{2}}&\to &
\mbox{\rm Hom}_{V}(({\cal M}^{W_{3}}_{W_{1}W_{2}})^{*}\otimes W_{3}, W_{3})\nno\\
I&\mapsto &(f\otimes w_{(3)}\mapsto f(I)w_{(3)})
\end{eqnarray}
which is an isomorphism if $W_{3}$ is irreducible, since in this
case, 
\[
\mbox{\rm Hom}_{V}(W_{3}, W_{3})\simeq {\C}
\]
(see \cite{FHL}, Remark 4.7.1).  On the other hand, the natural map
\begin{eqnarray}
h:\mbox{\rm Hom}_{V}(({\cal M}^{W_{3}}_{W_{1}W_{2}})^{*}
\otimes W_{3}, W_{3})&\to &
{\cal M}^{W_{3}}_{W_{1}W_{2}}\nno\\
\eta&\mapsto &\overline{\eta}\circ {\cal F}[P(z)]^{W_{3}}_{W_{1}W_{2}}
\end{eqnarray}
given by composition clearly satisfies the condition that
\begin{equation}\label{hiI=I}
h(i(I))=I,
\end{equation}
so that if $W_{3}$ is irreducible, the maps $h$ and $i$ are mutually inverse 
isomorphisms and we have the universal property that for any $I\in 
{\cal M}^{W_{3}}_{W_{1}W_{2}}$, there exists a unique $\eta$ such that 
\begin{equation}\label{I=etabarF}
I=\overline{\eta}\circ {\cal I}^{W_{3}}_{W_{1}W_{2}}
\end{equation}
(cf. Definition \ref{pz-tp}).

Using this, we can now show, in the next result, that $P(z)$-tensor
products always exist for the category of modules for a finitely
reductive vertex operator algebra, and we shall in fact exhibit the
$P(z)$-tensor product.  Note that there is no need to assume that
$W_{1}$ and $W_{2}$ are irreducible in the formulation or proof, but
by Proposition \ref{tensorproductdistributes}, the case in which
$W_{1}$ and $W_{2}$ are irreducible is in fact sufficient, and the
tensor product operation is canonically described using only the
spaces of intertwining maps among triples of {\it irreducible}
modules.

\begin{propo}\label{construcofPztensorprod-finredcase}
Let $V$ be a finitely reductive vertex operator algebra and let
$W_{1}$ and $W_{2}$ be $V$-modules.  Then
$(W_{1}\boxtimes_{P(z)}W_{2}, Y_{P(z)}; \boxtimes_{P(z)})$ exists, and
in fact
\begin{equation}\label{Pztensorprodfinitelyredcase}
W_{1}\boxtimes_{P(z)}W_{2}=\coprod_{i=1}^{k}
({\cal M}^{M_{i}}_{W_{1}W_{2}})^{*}\otimes M_{i},
\end{equation}
where $\{ M_{1}, \dots, M_{k}\}$ is a set of representatives of the
equivalence classes of irreducible $V$-modules, and the right-hand
side of (\ref{Pztensorprodfinitelyredcase}) is equipped with the
$V$-module and $P(z)$-product structure indicated above. That is,
\begin{equation}
\boxtimes_{P(z)}=\sum_{i=1}^{k}{\cal F}[P(z)]^{M_{i}}_{W_{1}W_{2}}.
\end{equation}
\end{propo}
\pf
{From} the comments above and the definitions, it is clear that we have a 
$P(z)$-product. Let $(W_{3}, Y_{3}; I)$ be any $P(z)$-product. Then 
$W_{3}=\coprod_{j}U_{j}$ where $j$ ranges through a finite set and each
$U_{j}$ is irreducible. Let $\pi_{j}: W_{3}\to U_{j}$ denote the $j$-th 
projection. A module map $\eta:\coprod_{i=1}^{k}
({\cal M}^{M_{i}}_{W_{1}W_{2}})^{*}\otimes M_{i}\to W_{3}$ amounts to 
module maps
$$\eta_{ij}:  ({\cal M}^{M_{i}}_{W_{1}W_{2}})^{*}\otimes M_{i}\to U_{j}$$
for each $i$ and $j$ such that $U_{j}\simeq M_{i}$, and 
$I=\overline{\eta}\circ \boxtimes_{P(z)}$ if and only if
\[
\overline{\pi}_{j}\circ I=\overline{\eta}_{ij}\circ 
{\cal F}^{M_{i}}_{W_{1}W_{2}}
\]
for each $i$ and $j$, the bars having the obvious meaning. But 
$\overline{\pi}_{j}\circ I$ is a $P(z)$-intertwining map of type 
${U_{j}}\choose {W_{1}W_{2}}$, and so 
$\overline{\iota}\circ \overline{\pi}_{j}\circ I\in 
{\cal M}^{M_{i}}_{W_{1}W_{2}}$,
where $\iota: U_{j}\stackrel{\sim}{\to}M_{i}$ is a fixed isomorphism. 
Denote this map by $\tau$. Thus what we finally want is a unique module map
\[
\theta: ({\cal M}^{M_{i}}_{W_{1}W_{2}})^{*}\otimes M_{i}\to M_{i}
\]
such that 
\[
\tau=\overline{\theta}\circ {\cal F}[P(z)]^{M_{i}}_{W_{1}W_{2}}.
\]
But we in fact have such a unique $\theta$, by
(\ref{hiI=I})--(\ref{I=etabarF}).  \epfv

\begin{rema} {\rm
By combining Proposition \ref{construcofPztensorprod-finredcase} with
Proposition \ref{im:correspond}, we can express $W_{1}\boxtimes_{P(z)}
W_{2}$ in terms of ${\cal V}^{M_{i}}_{W_{1}W_{2}}$ in place of ${\cal
M}^{M_{i}}_{W_{1}W_{2}}$.
}
\end{rema}

\begin{rema} {\rm
If we know the fusion rules among triples of irreducible $V$-modules,
then {}from Proposition \ref{construcofPztensorprod-finredcase} we know
all the $P(z)$-tensor product modules, up to equivalence; that is, we
know the multiplicity of each irreducible $V$-module in each
$P(z)$-tensor product module.  But recall that the $P(z)$-tensor
product structure of $W_{1}\boxtimes_{P(z)} W_{2}$ involves much more
than just the $V$-module structure.
}
\end{rema}

As we discussed in the Introduction, the main theme of this work is to
construct natural ``associativity'' isomorphisms between triple tensor
products of the shape $W_{1}\boxtimes (W_{2}\boxtimes W_{3})$ and
$(W_{1}\boxtimes W_{2})\boxtimes W_{3}$, for (generalized) modules
$W_{1}$, $W_{2}$ and $W_{3}$.  In the finitely reductive case, let
$W_{1}$, $W_{2}$ and $W_{3}$ be $V$-modules.  By Proposition
\ref{construcofPztensorprod-finredcase}, we have, as $V$-modules,
\begin{eqnarray}\label{W1(W2W3)}
\lefteqn{W_{1}\boxtimes_{P(z)}(W_{2}\boxtimes_{P(z)}W_{3})
= W_{1}\boxtimes_{P(z)}\left(\coprod_{i=1}^{k} M_{i} \otimes ({\cal
M}^{M_{i}}_{W_{2}W_{3}})^{*}\right)}
\nno\\
&&
=\coprod_{i=1}^{k} (W_{1}\boxtimes_{P(z)} M_{i}) \otimes ({\cal
M}^{M_{i}}_{W_{2}W_{3}})^{*}
\nno\\
&&
=\coprod_{i=1}^{k} \left(\coprod_{j=1}^{k} ({\cal
M}^{M_{j}}_{W_{1}M_{i}})^{*} \otimes M_{j}\right) \otimes ({\cal
M}^{M_{i}}_{W_{2}W_{3}})^{*}
\nno\\
&&
=\coprod_{j=1}^{k} \left(\coprod_{i=1}^{k} ({\cal
M}^{M_{j}}_{W_{1}M_{i}})^{*} \otimes ({\cal
M}^{M_{i}}_{W_{2}W_{3}})^{*}\right) \otimes M_{j}
\nno\\
&&
=\coprod_{j=1}^{k} \left(\coprod_{i=1}^{k} ({\cal
M}^{M_{j}}_{W_{1}M_{i}} \otimes {\cal
M}^{M_{i}}_{W_{2}W_{3}})^{*}\right) \otimes M_{j}
\end{eqnarray}
and
\begin{eqnarray}\label{(W1W2)W3}
\lefteqn{(W_{1}\boxtimes_{P(z)}W_{2})\boxtimes_{P(z)}W_{3}
=\left(\coprod_{i=1}^{k} M_{i} \otimes ({\cal
M}^{M_{i}}_{W_{1}W_{2}})^{*} \right) \boxtimes_{P(z)} W_{3}}
\nno\\
&&
=\coprod_{i=1}^{k} (M_{i} \boxtimes_{P(z)} W_{3}) \otimes ({\cal
M}^{M_{i}}_{W_{1}W_{2}})^{*}
\nno\\
&&
=\coprod_{i=1}^{k} \left(\coprod_{j=1}^{k} ({\cal
M}^{M_{j}}_{M_{i}W_{3}})^{*} \otimes M_{j}\right) \otimes ({\cal
M}^{M_{i}}_{W_{1}W_{2}})^{*}
\nno\\
&&
=\coprod_{j=1}^{k} \left(\coprod_{i=1}^{k} ({\cal
M}^{M_{j}}_{M_{i}W_{3}})^{*} \otimes ({\cal
M}^{M_{i}}_{W_{1}W_{2}})^{*}\right) \otimes M_{j}
\nno\\
&&
=\coprod_{j=1}^{k} \left(\coprod_{i=1}^{k} ({\cal
M}^{M_{j}}_{M_{i}W_{3}} \otimes {\cal
M}^{M_{i}}_{W_{1}W_{2}})^{*}\right) \otimes M_{j}.
\end{eqnarray}
These two $V$-modules will be equivalent if for each $j=1,\dots,k$,
their $M_j$-multiplicities are the same, that is, if
\begin{equation}\label{fusionrulerelation}
\sum_{i=1}^{k} N^{M_{j}}_{W_{1}M_{i}} N^{M_{i}}_{W_{2}W_{3}}=
\sum_{i=1}^{k} N^{M_{i}}_{W_{1}W_{2}} N^{M_{j}}_{M_{i}W_{3}}.
\end{equation}

{\it However,} knowing only that these two $V$-modules are equivalent
(knowing that $\boxtimes$ is ``associative'' in only a rough sense) is
far {}from enough.  What we need is a natural isomorphism between these
two modules analogous to the natural isomorphism
\begin{equation}\label{calWassociativity}
{\cal W}_1 \otimes ({\cal W}_2 \otimes {\cal W}_3)
\stackrel{\sim}{\longrightarrow} ({\cal W}_1 \otimes {\cal W}_2) \otimes
{\cal W}_3
\end{equation}
of vector spaces ${\cal W}_i$ determined by the natural condition
\begin{equation}\label{wassociativity}
w_{(1)} \otimes (w_{(2)} \otimes w_{(3)}) \mapsto (w_{(1)} \otimes w_{(2)})
\otimes w_{(3)}
\end{equation}
on elements (recall the Introduction).  Suppose that ${\cal W}_1,$
${\cal W}_2$ and ${\cal W}_3$ are finite-dimensional completely
reducible modules for some Lie algebra.  Then we of course have the
analogue of the relation (\ref{fusionrulerelation}).  But knowing the
equality of these multiplicities certainly does not give the natural
isomorphism (\ref{calWassociativity})--(\ref{wassociativity}).

Our intent to construct a natural isomorphism between the spaces
(\ref{W1(W2W3)}) and (\ref{(W1W2)W3}) (under suitable conditions) in
fact provides a guide to what we need to do.  In (\ref{W1(W2W3)}),
each space ${\cal M}^{M_{j}}_{W_{1}M_{i}} \otimes {\cal
M}^{M_{i}}_{W_{2}W_{3}}$ suggests combining an intertwining map ${\cal
Y}_1$ of type ${M_{j}}\choose {W_{1}M_{i}}$ with an intertwining map
${\cal Y}_2$ of type ${M_{i}}\choose {W_{2}W_{3}}$, presumably by
composition:
\begin{equation}\label{Y1zY2z}
{\cal Y}_1 (w_{(1)},z){\cal Y}_2 (w_{(2)},z).
\end{equation}
But this will not work, since this compostion does not exist because
the relevant formal series in $z$ does not converge; we must instead
take
\begin{equation}\label{Y1z1Y2z2}
{\cal Y}_1 (w_{(1)},z_1){\cal Y}_2 (w_{(2)},z_2),
\end{equation}
where the complex numbers $z_1$ and $z_2$ are such that
\[
|z_1|>|z_2|>0,
\]
by analogy with, and generalizing, the situation in Corollary
\ref{dualitywithcovergence}.  The composition (\ref{Y1z1Y2z2}) must be
understood using convergence and ``matrix coefficients,'' again as in
Corollary \ref{dualitywithcovergence}.

Similarly, in (\ref{(W1W2)W3}), each space ${\cal
M}^{M_{j}}_{M_{i}W_{3}} \otimes {\cal M}^{M_{i}}_{W_{1}W_{2}}$
suggests combining an intertwining map ${\cal Y}^1$ of type
${M_{j}}\choose {M_{i}W_{3}}$ with an intertwining map of type ${\cal
Y}^2$ of type ${M_{i}}\choose {W_{1}W_{2}}$:
\[
{\cal Y}^1 ({\cal Y}^2 (w_{(1)},z_1-z_2)w_{(2)},z_2),
\]
a (convergent) iterate of intertwining maps as in 
(\ref{associativitywithz1,z2}), with
\[
|z_2|>|z_1-z_2|>0,
\]
{\it not}
\begin{equation}
{\cal Y}^1 ({\cal Y}^2 (w_{(1)},z),z),
\end{equation}
which fails to converge.

The natural way to construct a natural associativity isomorphism
between (\ref{W1(W2W3)}) and (\ref{(W1W2)W3}) will in fact, then, be
to implement a correspondence of the type
\begin{equation}\label{YY=Y(Y)}
{\cal Y}_1 (w_{(1)},z_1){\cal Y}_2 (w_{(2)},z_2)={\cal Y}^1 ({\cal
Y}^2 (w_{(1)},z_1-z_2)w_{(2)},z_2),
\end{equation}
as we have previewed in the Introduction (formula (\ref{yyyy2})) and
also in (\ref{associativitywithz1,z2}).  Formula (\ref{YY=Y(Y)})
expresses the existence and associativity of the general
nonmeromorphic operator product expansion, as discussed in Remark
\ref{OPE}.  Note that this viewpoint shows that we should not try
directly to construct a natural isomorphism
\begin{equation}
W_1\boxtimes_{P(z)} (W_2\boxtimes_{P(z)} W_3)
\stackrel{\sim}{\longrightarrow} (W_1 \boxtimes_{P(z)}
W_2)\boxtimes_{P(z)} W_3,
\end{equation}
but rather a natural isomorphism
\begin{equation}\label{naturalassociso}
W_1\boxtimes_{P(z_1)} (W_2\boxtimes_{P(z_2)} W_3)
\stackrel{\sim}{\longrightarrow} (W_1 \boxtimes_{P(z_1-z_2)}
W_2)\boxtimes_{P(z_2)} W_3.
\end{equation}
This is what we will actually do in this work, in the general
logarithmic, not-necessarily-finitely-reductive case, under suitable
conditions.  The natural isomorphism (\ref{naturalassociso}) will act
as follows on elements of the completions of the relevant
(generalized) modules:
\begin{equation}
w_{(1)}\boxtimes_{P(z_1)} (w_{(2)}\boxtimes_{P(z_2)} w_{(3)})
\mapsto (w_{(1)} \boxtimes_{P(z_1-z_2)}
w_{(2)})\boxtimes_{P(z_2)} w_{(3)},
\end{equation}
implementing the stategy suggested by the classical natural
isomorphism (\ref{calWassociativity})--(\ref{wassociativity}).  Recall
that we previewed this strategy in the Introduction.

It turns out that in order to carry out this program, including the
construction of equalities of the type (\ref{YY=Y(Y)}) (the existence
and associativity of the nonmeromorphic operator product expansion) in
general, we cannot use the realization of the $P(z)$-tensor product
given in Proposition \ref{construcofPztensorprod-finredcase}, {\it
even when} $V$ {\it is a finitely reductive vertex operator algebra.}
As in \cite{tensor1}--\cite{tensor3} and \cite{tensor4}, what we do
instead is to construct $P(z)$-tensor products in a completely
different way (even in the finitely reductive case), a way that allows
us to also construct the natural associativity isomorphisms.  Section
5 is devoted to this construction of $P(z)$- (and $Q(z)$-)tensor
products.

\subsection{The notion of $Q(z)$-tensor product}

We now generalize the notion of $Q(z)$-tensor product of modules {}from
\cite{tensor1} to the setting of the present work, parallel to what we
did for the $P(z)$-tensor product above.  Here we give only the
results that we will need later.  Other results similar to those for
$P(z)$-tensor products certainly also carry over to the case of
$Q(z)$, for example, the results above on the finitely reductive case,
as were presented in \cite{tensor1}.

\begin{defi}\label{im:qimdef}{\rm
Let $(W_1, Y_{1})$, $(W_2, Y_{2})$ and $(W_3, Y_{3})$ 
be generalized $V$-modules.  A {\it $Q(z)$-intertwining map of type
${W_3\choose W_1\,W_2}$} is a linear map 
\[
I: W_1\otimes W_2 \to
\overline{W}_3
\]
such that the following conditions are satisfied: 
the {\it grading
compatibility condition}: for $\beta, \gamma\in \tilde{A}$ and
$w_{(1)}\in W_{1}^{(\beta)}$, $w_{(2)}\in W_{2}^{(\gamma)}$,
\begin{equation}\label{grad-comp-qz}
I(w_{(1)}\otimes w_{(2)})\in \overline{W_{3}^{(\beta+\gamma)}};
\end{equation}
the
{\em lower truncation condition:} for any elements
$w_{(1)}\in W_1$, $w_{(2)}\in W_2$, and any $n\in {\mathbb C}$,
\begin{equation}\label{imq:ltc}
\pi_{n-m}I(w_{(1)}\otimes w_{(2)})=0\;\;\mbox{ for }\;m\in {\mathbb
N}\;\mbox{ sufficiently large}
\end{equation}
(which follows {}from (\ref{grad-comp-qz}), in view of the
grading restriction condition (\ref{set:dmltc}); cf. (\ref{im:ltc}));
the {\em Jacobi identity}:
\begin{eqnarray}\label{imq:def}
\lefteqn{z^{-1}\delta\left(\frac{x_1-x_0}{z}\right)
Y^o_3(v, x_0)I(w_{(1)}\otimes w_{(2)})}\nonumber\\
&&=x_0^{-1}\delta\left(\frac{x_1-z}{x_0}\right)
I(Y_1^{o}(v, x_1)w_{(1)}\otimes w_{(2)})\nonumber\\
&&\hspace{2em}-x_0^{-1}\delta\left(\frac{z-x_1}{-x_0}\right)
I(w_{(1)}\otimes Y_2(v, x_1)w_{(2)})
\end{eqnarray}
for $v\in V$, $w_{(1)}\in W_1$ and $w_{(2)}\in W_2$ (recall
(\ref{yo}) for the notation $Y^{o}$, and note that the 
left-hand side of (\ref{imq:def}) is
meaningful because any infinite linear combination of $v_n$ of the
form $\sum_{n<N}a_nv_n$ ($a_n\in {\mathbb C}$) acts on any
$I(w_{(1)}\otimes w_{(2)})$, in view of (\ref{imq:ltc})); and the {\em
${\mathfrak s}{\mathfrak l}(2)$-bracket relations}: for any
$w_{(1)}\in W_1$ and $w_{(2)}\in W_2$,
\begin{eqnarray}\label{imq:Lj}
L(-j)I(w_{(1)}\otimes w_{(2)})&=&\sum_{i=0}^{j+1}{j+1\choose
i}(-z)^iI((L(-j+i)w_{(1)})\otimes w_{(2)})\nno\\
&&-\sum_{i=0}^{j+1}{j+1\choose i}(-z)^iI(w_{(1)}\otimes
L(j-i)w_{(2)})
\end{eqnarray}
for $j=-1, 0$ and $1$ (note that if $V$ is in fact a conformal vertex
algebra, this follows automatically {}from (\ref{imq:def}) by setting
$v=\omega$ and taking $\res_{x_1}\res_{x_0}x_0^{j+1}$). The vector
space of $Q(z)$-intertwining maps of type ${W_3\choose W_1\,W_2}$ is
denoted by ${\cal M}[Q(z)]^{W_3}_{W_1W_2}$.}
\end{defi}

\begin{rema}\label{Q(z)geometry}{\rm
As was explained in \cite{tensor1}, the symbol $Q(z)$ represents the
Riemann sphere ${\mathbb C}\cup \{ \infty \}$ with one negatively
oriented puncture at $z$ and two ordered positively oriented punctures
at $\infty$ and $0$, with local coordinates $w-z$, $1/w$ and $w$,
respectively, vanishing at these punctures.  In fact, this structure
is conformally equivalent to the Riemann sphere ${\mathbb C}\cup \{
\infty \}$ with one negatively oriented puncture at $\infty$ and two
ordered positively oriented punctures $1/z$ and $0$, with local
coordinates $z/(zw-1)$, $(zw-1)/z^2w$ and $z^2w/(zw-1)$ vanishing at
$\infty$, $1/z$ and $0$, respectively. }
\end{rema}

\begin{rema}{\rm
In the case of $\C$-graded ordinary modules for a vertex operator
algebra, where the grading restriction condition (\ref{Wn+k=0}) for a
module $W$ is replaced by the (more restrictive) condition
\begin{equation}
W_{(n)}=0 \;\; \mbox { for }\;n\in {\C}\;\mbox{ with sufficiently
negative real part}
\end{equation}
as in \cite{tensor1} (and where, in our context, the abelian groups
$A$ and $\tilde{A}$ are trivial), the notion of $Q(z)$-intertwining
map above agrees with the earlier one introduced in \cite{tensor1}; in
this case, the conditions (\ref{grad-comp-qz}) and (\ref{imq:ltc}) are
automatic.}
\end{rema}

In view of Remarks \ref{P(z)geometry} and \ref{Q(z)geometry}, we can
now give a natural correspondence between $P(z)$- and
$Q(z)$-intertwining maps.  (See the next three results.)  Recall that
since our generalized $V$-modules are strongly graded, we have
contragredient generalized modules of generalized modules.

\begin{propo}\label{qp:qp}
Let $I: W_1\otimes W_2\to \overline{W}_3$ and $J: W'_3\otimes W_2\to
\overline{W'_1}$ be linear maps related to each other by:
\begin{equation}\label{qz:qtop}
\langle w_{(1)},J(w'_{(3)}\otimes w_{(2)})\rangle=
\langle w'_{(3)},I(w_{(1)}\otimes w_{(2)})\rangle
\end{equation}
for any $w_{(1)}\in W_1$, $w_{(2)}\in W_2$ and $w'_{(3)}\in
W'_3$. Then $I$ is a $Q(z)$-intertwining map of type ${W_3\choose
W_1\,W_2}$ if and only if $J$ is a $P(z)$-intertwining map of type
${W'_1\choose W'_3\,W_2}$.
\end{propo}
\pf Suppose that $I$ is a $Q(z)$-intertwining map of type ${W_3\choose
W_1\,W_2}$.  We shall show that $J$ is a $P(z)$-intertwining map of
type ${W'_1\choose W'_3\,W_2}$.

Since $I$ satisfies the grading compatibility condition, it is clear
that $J$ also satisfies this condition.  For the lower truncation
condition for $J$, it suffices to show that for any $w_{(2)}\in
W_2^{(\beta)}$ and $w'_{(3)}\in (W'_3)^{(\gamma)}$, where $\beta, \gamma
\in \tilde{A}$, and any $n\in{\mathbb C}$, $\langle
\pi_{[n-m]}W_1^{(-\beta -\gamma)}, J(w'_{(3)}\otimes
w_{(2)})\rangle=0$ for $m\in{\mathbb N}$ sufficiently large, or that
\begin{equation}\label{qz:Jltrp}
\langle w'_{(3)},I(\pi_{[n-m]}W_1^{(-\beta -\gamma)} \otimes
w_{(2)})\rangle=0\;\;\mbox{ for }\;m\in{\mathbb
N}\;\;\mbox{sufficiently large.}
\end{equation}
But (\ref{qz:Jltrp}) follows immediately {}from (\ref{set:dmltc}).

Now we prove the Jacobi identity for $J$. The Jacobi identity for $I$
gives
\begin{eqnarray}\label{qz:jcba}
\lefteqn{z^{-1}\delta\left(\frac{x_1-x_0}{z}\right)\langle w'_{(3)},
Y^o_3(v, x_0)I(w_{(1)}\otimes w_{(2)})\rangle}\nno\\
&&=x_0^{-1}\delta\left(\frac{x_1-z}{x_0}\right)\langle w'_{(3)},
I(Y_1^{o}(v, x_1)w_{(1)}\otimes w_{(2)})\rangle\nno\\
&&\hspace{2em}-x_0^{-1}\delta\left(\frac{z-x_1}{-x_0} \right)\langle
w'_{(3)},I(w_{(1)}\otimes Y_2(v,x_1)w_{(2)})\rangle
\end{eqnarray}
for any $v\in V$, $w_{(1)}\in W_1$, $w_{(2)}\in W_2$ and $w'_{(3)}\in
W'_3$. By (\ref{y'}) the left-hand side is equal to
\[
z^{-1}\delta\left(\frac{x_1-x_0}{z}\right)\langle Y'_3(v,x_0)w'_{(3)},
I(w_{(1)}\otimes w_{(2)})\rangle
\]
So by (\ref{qz:qtop}), the identity (\ref{qz:jcba}) can be written as
\begin{eqnarray*}
\lefteqn{z^{-1}\delta\left(\frac{x_1-x_0}{z}\right)\langle w_{(1)},
J(Y'_3(v,x_0)w'_{(3)}\otimes w_{(2)})\rangle}\nno\\
&&=x_0^{-1}\delta\left(\frac{x_1-z}{x_0}\right)\langle Y_1^{o}(v,
x_1)w_{(1)}, J(w'_{(3)}\otimes w_{(2)})\rangle\nno\\
&&\hspace{2em}-x_0^{-1}\delta\left(\frac{z-x_1}{-x_0} \right)\langle
w_{(1)},J(w'_{(3)}\otimes Y_2(v,x_1)w_{(2)})\rangle.
\end{eqnarray*}
Applying (\ref{y'}) to the first term of the right-hand side
we see that this can be written as
\begin{eqnarray*}
\lefteqn{z^{-1}\delta\left(\frac{x_1-x_0}{z}\right)\langle w_{(1)},
J(Y'_3(v,x_0)w'_{(3)}\otimes w_{(2)})\rangle}\nno\\
&&=x_0^{-1}\delta\left(\frac{x_1-z}{x_0}\right)\langle w_{(1)},
Y'_1(v, x_1)J(w'_{(3)}\otimes w_{(2)})\rangle\nno\\
&&\hspace{2em}-x_0^{-1}\delta\left(\frac{z-x_1}{-x_0} \right)\langle
w_{(1)},J(w'_{(3)}\otimes Y_2(v,x_1)w_{(2)})\rangle
\end{eqnarray*}
for any $v\in V$, $w_{(1)}\in W_1$, $w_{(2)}\in W_2$ and $w'_{(3)}\in
W'_3$. This is exactly the Jacobi identity for $J$.

The ${\mathfrak s}{\mathfrak l}(2)$-bracket relations can be proved similarly,
as follows: The ${\mathfrak s}{\mathfrak l}(2)$-bracket relations for $I$ give
\begin{eqnarray*}
\langle w'_{(3)},L(-j)I(w_{(1)}\otimes w_{(2)})\rangle&=&
\sum_{i=0}^{j+1}{j+1\choose i}(-z)^i\langle
w'_{(3)},I((L(-j+i)w_{(1)})\otimes w_{(2)})\rangle\nno\\
&&-\sum_{i=0}^{j+1}{j+1\choose i}(-z)^i\langle
w'_{(3)},I(w_{(1)}\otimes L(j-i)w_{(2)})\rangle
\end{eqnarray*}
for any $w_{(1)}\in W_1$, $w_{(2)}\in W_2$, $w'_{(3)}\in W'_3$ and
$j=-1, 0, 1$. Using (\ref{L'(n)}) and then applying (\ref{qz:qtop}) we
get
\begin{eqnarray*}
\langle w_{(1)},J(L'(j)w'_{(3)}\otimes w_{(2)})\rangle&=&
\sum_{i=0}^{j+1}{j+1\choose i}(-z)^i\langle L(-j+i)w_{(1)},
J(w'_{(3)}\otimes w_{(2)})\rangle\nno\\
&&-\sum_{i=0}^{j+1}{j+1\choose i}(-z)^i\langle
w_{(1)},J(w'_{(3)}\otimes L(j-i)w_{(2)})\rangle,
\end{eqnarray*}
or
\begin{eqnarray*}
J(L'(j)w'_{(3)}\otimes w_{(2)})&=& \sum_{i=0}^{j+1}{j+1\choose
i}(-z)^i L(j-i)J(w'_{(3)}\otimes w_{(2)})\nno\\
&&-\sum_{i=0}^{j+1}{j+1\choose i}(-z)^iJ(w'_{(3)}\otimes
L(j-i)w_{(2)}),
\end{eqnarray*}
for $j=-1,0,1$. This is the alternative form (\ref{im:Lj2}) of the
${\mathfrak s} {\mathfrak l}(2)$-bracket relations for $J$. Hence $J$ is a
$P(z)$-intertwining map.

The other direction of the proposition is proved by simply reversing
the order of the arguments. \epfv

Let $W_{1}$, $W_{2}$ and $W_{3}$ be generalized $V$-modules, as above.
We shall call an element $\lambda$ of $(W_{1}\otimes W_{2}\otimes
W_{3})^{*}$ {\it $\tilde{A}$-compatible} if
\[
\lambda((W_{1})^{(\beta)}\otimes (W_{2})^{(\gamma)}\otimes
(W_{3})^{(\delta)})=0
\]
for $\beta, \gamma, \delta\in \tilde{A}$ satisfying
$\beta+\gamma+\delta\ne 0$.  Recall {}from Definitions \ref{Wbardef} and
\ref{defofWprime} that for a generalized $V$-module $W$,
$\overline{W'}$ can be viewed as a (usually proper) subspace of
$W^{*}$.  We shall call a linear map
\[
I: W_{1}\otimes W_{2} \rightarrow W_{3}^{*}
\]
{\it $\tilde{A}$-compatible} if its image lies in $\overline{W_{3}'}$,
that is,
\begin{equation}\label{IAtildecompat}
I:W_{1}\otimes W_{2} \rightarrow \overline{W_{3}'},
\end{equation}
and if $I$ satisfies the usual grading compatibility condition
(\ref{grad-comp}) or (\ref{grad-comp-qz}) for $P(z)$- or
$Q(z)$-intertwining maps.  
Now an element $\lambda$ of $(W_{1}\otimes
W_{2}\otimes W_{3})^{*}$ amounts exactly to a linear map
\[
I_{\lambda}:W_{1}\otimes W_{2} \rightarrow W_{3}^{*}.
\]
If $\lambda$ is $\tilde{A}$-compatible, then for $w_{(1)}\in
W_{1}^{(\beta)}$, $w_{(2)}\in W_{2}^{(\gamma)}$ and $w_{(3)}\in
W_{3}^{(\delta)}$ such that $\delta\ne -(\beta +\gamma)$,
\[
\langle w_{(3)}, I_{\lambda}(w_{(1)}\otimes w_{(2)})\rangle=
\lambda(w_{(1)}\otimes w_{(2)}\otimes w_{(3)})=0,
\]
so that $I_{\lambda}(w_{(1)}\otimes w_{(2)})\in
\overline{(W_{3}')^{(\beta+\gamma)}}$ and $I_{\lambda}$ is
$\tilde{A}$-compatible.  Similarly, if $I_{\lambda}$ is
$\tilde{A}$-compatible, then so is $\lambda$.  Thus we have the
following straightforward result relating $\tilde{A}$-compatibility of
$\lambda$ with that of $I_{\lambda}$:

\begin{lemma}\label{4.36}
The linear functional $\lambda\in (W_{1}\otimes W_{2}\otimes
W_{3})^{*}$ is $\tilde{A}$-compatible if and only if $I_{\lambda}$ is
$\tilde{A}$-compatible.  The map given by $\lambda\mapsto I_{\lambda}$
is the unique linear isomorphism {}from the space of
$\tilde{A}$-compatible elements of $(W_{1}\otimes W_{2}\otimes
W_{3})^{*}$ to the space of $\tilde{A}$-compatible linear maps {}from
$W_{1}\otimes W_{2}$ to $\overline{W_{3}'}$ such that
\[
\langle w_{(3)}, I_{\lambda}(w_{(1)}\otimes w_{(2)})\rangle
=\lambda(w_{(1)}\otimes w_{(2)}\otimes w_{(3)})
\]
for $w_{(1)}\in W_{1}$, $w_{(2)}\in W_{2}$ and $w_{(3)}\in W_{3}$.
Similarly, there are canonical linear isomorphisms {}from the space of
$\tilde{A}$-compatible elements of $(W_{1}\otimes W_{2}\otimes
W_{3})^{*}$ to the space of $\tilde{A}$-compatible linear maps {}from
$W_{1}\otimes W_{3}$ to $\overline{W_{2}'}$ and to the space of
$\tilde{A}$-compatible linear maps {}from $W_{2}\otimes W_{3}$ to
$\overline{W_{1}'}$ satisfying the corresponding conditions.  In
particular, there is a canonical linear isomorphism {}from the space of
$\tilde{A}$-compatible linear maps {}from $W_{1}\otimes W_{2}$ to
$\overline{W_{3}}$ to the space of $\tilde{A}$-compatible linear maps
{}from $W_{3}'\otimes W_{2}$ to $\overline{W_{1}'}$ given by
(\ref{qz:qtop}).  \epf
\end{lemma}

Using this lemma and Proposition \ref{qp:qp}, we have:

\begin{corol}\label{Q(z)P(z)iso}
The formula (\ref{qz:qtop}) gives a canonical linear isomorphism
between the space of $Q(z)$-intertwining maps of type ${W_3\choose
W_1\,W_2}$ and the space of $P(z)$-intertwining maps of type
${W'_1\choose W'_3\,W_2}$.  \epf 
\end{corol}

We can now use Proposition \ref{im:correspond} together with
Proposition \ref{qp:qp} and Corollary \ref{Q(z)P(z)iso} to construct a
correspondence between the logarithmic intertwining operators of type
${W'_1}\choose {W'_3W_2}$ and the $Q(z)$-intertwining maps of type
${W_3}\choose {W_1W_2}$; this generalizes the corresponding result in
the finitely reductive case, with ordinary modules, in \cite{tensor1}.
Fix an integer $p$.  Let ${\cal Y}$ be a logarithmic intertwining
operator of type ${W'_1}\choose {W'_3W_2}$, and use (\ref{log:IYp}) to
define a linear map $I_{{\cal Y}, p}: W'_3\otimes W_2\to
\overline{W_{1}'}$; by Proposition \ref{im:correspond}, this is a
$P(z)$-intertwining map of the same type.  Then use Proposition
\ref{qp:qp} and Corollary \ref{Q(z)P(z)iso} to define a
$Q(z)$-intertwining map $I^{Q(z)}_{{\cal Y}, p}: W_1\otimes W_2\to
\overline{W}_3$ of type ${W_3}\choose {W_1W_2}$ (uniquely) by
\begin{eqnarray}\label{imq:IYp}
\lefteqn{\langle w'_{(3)}, I^{Q(z)}_{{\cal Y}, p}(w_{(1)}\otimes
w_{(2)})\rangle_{W_3} = \langle w_{(1)}, I_{{\cal Y},
p}(w'_{(3)}\otimes w_{(2)})\rangle_{W'_1}}\nno\\
&&=\langle w_{(1)}, {\cal Y}(w'_{(3)},
e^{l_{p}(z)})w_{(2)}\rangle_{W'_1}
\;\;\;\;\;\;\;\;\;\;\;\;\;\;\;\;\;\;\;\;\;\;\;\;\;\;\;\;\;\;\;\;\;\;\;\;\;\;\;
\end{eqnarray}
for all $w_{(1)}\in W_1$, $w_{(2)}\in W_2$, $w'_{(3)}\in W'_3$. (We
are using the symbol $Q(z)$ to distinguish this {}from the $P(z)$ case
above.)  Then the correspondence ${\cal Y} \mapsto I^{Q(z)}_{{\cal Y},
p}$ is an isomorphism {}from ${\cal V}^{W_{1}'}_{W_{3}'W_{2}}$ to ${\cal
M}[Q(z)]^{W_3}_{W_1W_2}$.  {}From Proposition \ref{im:correspond} and
(\ref{recover}), its inverse is given by sending a $Q(z)$-intertwining
map $I$ of type ${W_3}\choose {W_1W_2}$ to the logarithmic
intertwining operator ${\cal Y}^{Q(z)}_{I,p}:W'_3\otimes W_2\to
W'_1[\log x]\{x\}$ defined by
\begin{eqnarray*}
\lefteqn{\langle w_{(1)}, {\cal Y}^{Q(z)}_{I,p}(w'_{(3)},
x)w_{(2)}\rangle_{W'_1}}\\ 
&&=\langle y^{-L'(0)}x^{-L'(0)}w'_{(3)},
I(y^{L(0)}x^{L(0)}w_{(1)}\otimes
y^{-L(0)}x^{-L(0)}w_{(2)})\rangle_{W_3}\lbar_{y=e^{-l_{p}(z)}}
\end{eqnarray*}
for any $w_{(1)}\in W_1$, $w_{(2)}\in W_2$ and $w'_{(3)}\in
W'_3$.  Thus we have:

\begin{propo}\label{Q-cor}
For $p\in {\mathbb Z}$, the correspondence ${\cal Y}\mapsto
I^{Q(z)}_{{\cal Y}, p}$ is a linear isomorphism {}from the space ${\cal
V}^{W'_1}_{W'_3W_2}$ of logarithmic intertwining operators of type
${W'_1}\choose {W'_3\; W_2}$ to the space ${\cal
M}[Q(z)]^{W_3}_{W_1W_2}$ of $Q(z)$-intertwining maps of type
${W_3}\choose {W_1W_2}$. Its inverse is given by $I\mapsto {\cal
Y}^{Q(z)}_{I, p}$. \epf
\end{propo}

We now give the definition of $Q(z)$-tensor product.

\begin{defi}\label{qz-product}{\rm
Let ${\cal C}_{1}$ be either ${\cal M}_{sg}$ or ${\cal GM}_{sg}$.
For $W_1, W_2\in \ob{\cal C}_{1}$, a {\it $Q(z)$-product of $W_1$ and
$W_2$} is an object $(W_3, Y_3)$ of ${\cal C}_{1}$ together with a
$Q(z)$-intertwining map $I_3$ of type ${W_3}\choose {W_1W_2}$. We
denote it by $(W_3, Y_3; I_3)$ or simply by $(W_3, I_3)$.  Let
$(W_4,Y_4; I_4)$ be another $Q(z)$-product of $W_1$ and $W_2$. A {\em
morphism} {}from $(W_3, Y_3; I_3)$ to $(W_4, Y_4; I_4)$ is a module map
$\eta$ {}from $W_3$ to $W_4$ such that the
diagram
\begin{center}
\begin{picture}(100,60)
\put(-5,0){$\overline W_3$}
\put(13,4){\vector(1,0){104}}
\put(119,0){$\overline W_4$}
\put(41,50){$W_1\otimes W_2$}
\put(61,45){\vector(-3,-2){50}}
\put(68,45){\vector(3,-2){50}}
\put(65,8){$\bar\eta$}
\put(20,27){$I_3$}
\put(98,27){$I_4$}
\end{picture}
\end{center}
commutes, that is,
\[
I_4=\overline{\eta}\circ I_3.
\]
where, as before, $\overline{\eta}$ is the natural map {from}
$\overline{W}_3$ to $\overline{W}_4$ uniquely extending $\eta$. }
\end{defi}

\begin{defi}\label{qz-tp}
{\rm Let ${\cal C}$ be a full subcategory of either ${\cal M}_{sg}$ or
${\cal GM}_{sg}$. For $W_1, W_2\in \ob{\cal C}$, a {\em $Q(z)$-tensor
product of $W_1$ and $W_2$ in ${\cal C}$} is a $Q(z)$-product $(W_0,
Y_0; I_0)$ with $W_0\in{\rm ob\,}{\cal C}$ such that for any
$Q(z)$-product $(W,Y;I)$ with $W\in{\rm ob\,}{\cal C}$, there is a
unique morphism {}from $(W_0, Y_0; I_0)$ to $(W,Y;I)$. Clearly, a
$Q(z)$-tensor product of $W_1$ and $W_2$ in ${\cal C}$, if it exists,
is unique up to unique isomorphism.  In this case we will denote it as
$(W_1\boxtimes_{Q(z)} W_2, Y_{Q(z)}; \boxtimes_{Q(z)})$ and call the
object $(W_1\boxtimes_{Q(z)} W_2, Y_{Q(z)})$ the {\em $Q(z)$-tensor
product module of $W_1$ and $W_2$ in ${\cal C}$}. Again we will skip
the phrase ``in ${\cal C}$'' if the category ${\cal C}$ under
consideration is clear in context. }
\end{defi}

The following immediate consequence of Definition \ref{qz-tp}
and Proposition \ref{Q-cor}
relates module maps {}from a
$Q(z)$-tensor product module with $Q(z)$-intertwining maps and
logarithmic intertwining operators:

\begin{propo}
Suppose that $W_1\boxtimes_{Q(z)}W_2$ exists. We have a natural
isomorphism
\begin{eqnarray*}
\hom_{V}(W_1\boxtimes_{Q(z)}W_2, W_3)&\stackrel{\sim}{\to}&
{\cal M}[Q(z)]^{W_3}_{W_1W_2}\\ \eta&\mapsto& \overline{\eta}\circ
\boxtimes_{Q(z)}
\end{eqnarray*}
and for $p\in {\mathbb Z}$, a natural isomorphism
\begin{eqnarray*}
\hom_{V}(W_1\boxtimes_{Q(z)} W_2, W_3)&
\stackrel{\sim}{\rightarrow}& {\cal V}^{W'_1}_{W'_3W_2}\\ \eta&\mapsto
& {\cal Y}^{Q(z)}_{\eta, p}
\end{eqnarray*}
where ${\cal Y}^{Q(z)}_{\eta, p}={\cal Y}^{Q(z)}_{I, p}$ with
$I=\overline{\eta}\circ \boxtimes_{Q(z)}$.\epf
\end{propo}

Suppose that the $Q(z)$-tensor product $(W_1\boxtimes_{Q(z)} W_2,
Y_{Q(z)}; \boxtimes_{Q(z)})$ of $W_1$ and $W_2$ exists.  We will
sometimes denote the action of the canonical $Q(z)$-intertwining map
\begin{equation}\label{q-actionofboxtensormap}
w_{(1)} \otimes w_{(2)}\mapsto\boxtimes_{Q(z)}(w_{(1)} \otimes
w_{(2)})=\boxtimes_{Q(z)}(w_{(1)}, z)w_{(2)} \in 
\overline{W_1\boxtimes_{Q(z)}W_2}
\end{equation}
on elements simply by $w_{(1)}\boxtimes_{Q(z)}
w_{(2)}$:
\begin{equation}\label{q-boxtensorofelements}
w_{(1)}\boxtimes_{Q(z)} w_{(2)}=\boxtimes_{Q(z)}(w_{(1)} \otimes
w_{(2)})=\boxtimes_{Q(z)}(w_{(1)}, z)w_{(2)}.
\end{equation}

Using Propositions \ref{log:omega} and \ref{log:A}, we have the
following result, generalizing Proposition 4.9 and Corollary 4.10 in
\cite{tensor1}:

\begin{propo}\label{b-r}
For any integer $r$, there is a natural isomorphism
\[
B_{r}: {\cal V}^{W_3}_{W_1W_2}\to {\cal V}^{W'_1}_{W'_3W_2}
\]
defined by the condition that for any logarithmic intertwining
operator ${\cal Y}$ in ${\cal V}^{W_3}_{W_1W_2}$ and $w_{(1)}\in W_1$,
$w_{(2)}\in W_2$, $w'_{(3)}\in W'_3$,
\begin{eqnarray}\label{4.31}
\lefteqn{\langle w_{(1)}, B_{r}({\cal Y})(w'_{(3)}, x)
w_{(2)}\rangle_{W'_1}}\nno\\
&&=\langle e^{-x^{-1}L(1)}w'_{(3)}, {\cal Y}(e^{xL(1)}w_{(1)},
x^{-1})e^{-xL(1)}e^{(2r+1)\pi iL(0)}
(x^{-L(0)})^{2}w_{(2)}\rangle_{W_3}.
\end{eqnarray}
\end{propo}
\pf
{From} Proposition \ref{log:omega},  for any integer $r_{1}$ we have an
isomorphism $\Omega_{r_{1}}$ {from} ${\cal
V}^{W_{3}}_{W_{1}W_{2}}$ to ${\cal V}^{W_{3}}_{W_{2}W_{1}}$, and {from}
Proposition \ref{log:A},  for any integer $r_{2}$ we have an
isomorphism $A_{r_{2}}$ {from} ${\cal V}^{W_{3}}_{W_{2}W_{1}}$
to ${\cal V}^{W'_{1}}_{W_{2}W'_{3}}$. By Proposition \ref{log:omega} again, 
for any integer $r_{3}$ there is an  isomorphism, which we again denote 
$\Omega_{r_{3}}$, {from} ${\cal V}^{W'_{1}}_{W_{2}W'_{3}}$ to ${\cal
V}^{W'_{1}}_{W'_{3}W_{2}}$. Thus for any triple $(r_{1}, r_{2},
r_{3})$ of integers, we have an isomorphism $\Omega_{r_{3}}\circ A_{r_{2}}\circ
\Omega_{r_{1}}$ {from} ${\cal V}^{W_{3}}_{W_{1}W_{2}}$ to ${\cal
V}^{W'_{1}}_{W'_{3}W_{2}}$.  Let ${\cal Y}$ be a logarithmic intertwining
operator in ${\cal V}^{W_{3}}_{W_{1}W_{2}}$ and $w_{(1)}$, $w_{(2)}$,
$w'_{(3)}$ elements of $W_{1}$, $W_{2}$, $W'_{3}$, respectively.
{From} the definitions of $\Omega_{r_{1}}$, $A_{r_{2}}$ and
$\Omega_{r_{3}}$, we have
\begin{eqnarray}\label{7.29}
\lefteqn{\langle (\Omega_{r_{3}}\circ A_{r_{2}}\circ 
\Omega_{r_{1}})({\cal Y})(w'_{(3)}, x)w_{(2)}, w_{(1)}\rangle_{W_1}=}\nno\\
&&=\langle e^{xL(-1)}A_{r_{2}}(\Omega_{r_{1}}
({\cal Y}))(w_{(2)}, e^{(2r_{3}+1)\pi i}x)w'_{(3)}, 
w_{(1)}\rangle_{W_1}\nno\\
&&=\langle A_{r_{2}}(\Omega_{r_{1}}
({\cal Y}))(w_{(2)}, e^{(2r_{3}+1)\pi i}x)w'_{(3)}, 
e^{xL(1)}w_{(1)}\rangle_{W_1}\nno\\
&&=\langle w'_{(3)}, \Omega_{r_{1}}({\cal Y})(e^{-xL(1)}e^{(2r_{2}+1)\pi iL(0)}
e^{-2(2r_{3}+1)\pi iL(0)}(x^{-L(0)})^{2}w_{(2)}, \nno\\
&&\hspace{10em}e^{-(2r_{3}+1)\pi i}x^{-1})
e^{xL(1)}w_{(1)}\rangle_{W_3}
\nno\\
&&=\langle w'_{(3)}, e^{-x^{-1}L(-1)}{\cal Y}(e^{xL(1)}w_{(1)}, 
e^{(2r_{1}+1)\pi i}e^{-(2r_{3}+1)\pi i}x^{-1})\cdot \nno\\
&&\hspace{4em}\cdot e^{-xL(1)}e^{(2r_{2}+1)\pi iL(0)}
e^{-2(2r_{3}+1)\pi iL(0)}(x^{-L(0)})^{2}w_{(2)}\rangle_{W_3}
\nno\\
&&=\langle e^{-x^{-1}L(1)}w'_{(3)}, {\cal Y}(e^{xL(1)}w_{(1)}, 
e^{2(r_{1}-r_{3})\pi i}x^{-1})\cdot \nno\\
&&\hspace{6em}\cdot e^{-xL(1)}e^{(2(r_{2}-2r_{3}-1)+1)\pi iL(0)}
(x^{-L(0)})^{2}w_{(2)}\rangle_{W_3}.
\end{eqnarray}
{From} (\ref{7.29}) we see that $\Omega_{r_{3}}\circ A_{r_{2}}\circ
\Omega_{r_{1}}$ depends only on $r_{2}-2r_{3}-1$ and $r_{1}-r_{3}$,
and the operators $\Omega_{r_{3}}\circ A_{r_{2}}\circ \Omega_{r_{1}}$
with different $r_{1}-r_{3}$ but the same $r_{2}-2r_{3}-1$ differ
{from} each other only by automorphisms of ${\cal
V}^{W_{3}}_{W_{1}W_{2}}$ (recall Remarks \ref{log:fcf},
\ref{exponentialaVhom} and \ref{Ys1s2s3}).  Thus for our purpose, we
need only consider those isomorphisms such that $r_{1}-r_{3}=0$.
Given any integer $r$, we choose two integers $r_{2}$ and $r_{3}$ such
that $r=r_{2}-2r_{3}-1$ and we define
\begin{equation}
B_{r}=\Omega_{r_{3}}\circ A_{r_{2}}\circ \Omega_{r_{3}}.
\end{equation}
{From} (\ref{7.29}) we see that $B_{r}$ is independent of the choices of $r_{2}$ and
$r_{3}$ and that 
(\ref{4.31}) holds.
\epfv

Combining the last two results, we obtain:

\begin{corol}
 For any $W_1, W_2, W_3\in \ob{\cal C}$ such that
$W_1\boxtimes_{Q(z)}W_2$ exists and any integers $p$ and $r$, we have
a natural isomorphism
\begin{eqnarray}
\hom_{V}(W_1\boxtimes_{Q(z)} W_2, W_3)&
\stackrel{\sim}{\rightarrow}&
{\cal V}^{W_3}_{W_1W_2}\nno\\
\eta&\mapsto &B^{-1}_{r}({\cal Y}^{Q(z)}_{\eta, p}).\hspace{2em}\square
\end{eqnarray}
\end{corol}

\subsection{$P(z)$-tensor products and $Q(z^{-1})$-tensor products}

Here we prove the following result:

\begin{theo}\label{pz-qz-1}
Let $W_{1}$ and $W_{2}$ be objects of a full subcategory $\mathcal{C}$
of either $\mathcal{M}_{sg}$ or $\mathcal{GM}_{sg}$. Then the
$P(z)$-tensor product of $W_{1}$ and $W_{2}$ exists if and only if the
$Q(z^{-1})$-tensor product of $W_{1}$ and $W_{2}$ exists.
\end{theo}
\pf 
Recalling our choice of branch (\ref{branch1}), let
\[
p=-\frac{\log (z^{-1})+\log z}{2\pi i}.
\]
Then $p$ is an integer and we have
\[
-(\log (z^{-1})+2\pi pi)=\log z,
\]
and 
\[
e^{-n(\log (z^{-1})+2\pi pi)}=e^{n\log z}
\]
for $n\in \C$.

{}From Propositions \ref{im:correspond}, \ref{b-r} and
\ref{Q-cor}, we see that for $W_{1}, W_{2}, W_{3}\in \ob \mathcal{C}$,
there is a linear isomorphism $\mu_{W_{1}W_{2}}^{W_{3}}:
\mathcal{M}[P(z)]_{W_{1}W_{2}}^{W_{3}} \to
\mathcal{M}[Q(z^{-1})]_{W_{1}W_{2}}^{W_{3}}$ defined by
\[
\mu_{W_{1}W_{2}}^{W_{3}}(I)=I^{Q(z^{-1})}_{B_{2p}(\mathcal{Y}_{I,0}),p}
\]
for $I\in \mathcal{M}[P(z)]_{W_{1}W_{2}}^{W_{3}}$.  
By definition, $\mu_{W_{1}W_{2}}^{W_{3}}(I)$ is determined uniquely
by (recalling (\ref{branch1})--(\ref{branch2}))
\begin{eqnarray}\label{mu}
\lefteqn{\langle w_{(3)}', \mu_{W_{1}W_{2}}^{W_{3}}(I)(w_{(1)}
\otimes w_{(2)})\rangle}\nn
&&=\langle w_{(3)}', I^{Q(z^{-1})}_{B_{2p}(\mathcal{Y}_{I, 0}),p}
(w_{(1)}\otimes w_{(2)})\rangle\nn
&&=\langle w_{(1)}, B_{2p}(\mathcal{Y}_{I, 0})(w'_{(3)},
e^{l_p(z^{-1})})w_{(2)}\rangle\nn
&&=\langle w_{(1)}, B_{2p}(\mathcal{Y}_{I, 0})(w'_{(3)},
e^{\log (z^{-1})+2\pi pi})w_{(2)}\rangle\nn
&&=\langle e^{-zL(1)}w'_{(3)}, 
\mathcal{Y}_{I, 0}(e^{z^{-1}L(1)}w_{(1)}, e^{\log z})e^{-z^{-1}L(1)}
e^{(2(2p)+1)i\pi L(0)}e^{-2(\log z^{-1}+2\pi pi))L(0)}w_{(2)}\rangle\nn
&&=\langle e^{-zL(1)}w'_{(3)}, 
\mathcal{Y}_{I, 0}(e^{z^{-1}L(1)}w_{(1)}, z)e^{-z^{-1}L(1)}
e^{i\pi L(0)}e^{-2(\log z^{-1})L(0)}w_{(2)}\rangle\nn
&&=\langle e^{-zL(1)}w'_{(3)}, I((e^{z^{-1}L(1)}w_{(1)})
\otimes (e^{-z^{-1}L(1)}
e^{i\pi L(0)}e^{-2(\log z^{-1})L(0)}w_{(2)}))\rangle\nn
\end{eqnarray}
for $w_{(1)}\in W_{1}$, $w_{(2)}\in W_{2}$ and $w_{(3)}'\in W_{3}'$.
{}From (\ref{mu}), we also see that
for $J\in \mathcal{M}[Q(z^{-1})]_{W_{1}W_{2}}^{W_{3}}$,
$(\mu_{W_{1}W_{2}}^{W_{3}})^{-1}(J)$ is determined uniquely
by 
\begin{eqnarray}\label{mu-1}
\lefteqn{\langle w_{(3)}', (\mu_{W_{1}W_{2}}^{W_{3}})^{-1}(J)(w_{(1)}
\otimes w_{(2)})\rangle}\nn
&&=\langle e^{zL(1)}w'_{(3)}, 
J((e^{-z^{-1}L(1)}w_{(1)})
\otimes (e^{2(\log z^{-1})L(0)}
e^{-i\pi L(0)}e^{z^{-1}L(1)}w_{(2)}))\rangle\nn
\end{eqnarray}
for $w_{(1)}\in W_{1}$, $w_{(2)}\in W_{2}$ and $w_{(3)}'\in W_{3}'$.

Assume that the
$P(z)$-tensor product $(W_{1}\boxtimes_{P(z)} W_{2}, Y_{P(z)};
\boxtimes_{P(z)})$ exists.  Then
\[
\boxtimes_{Q(z^{-1})}=\mu_{W_{1}W_{2}}^{W_{1}\boxtimes_{P(z)}
W_{2}}(\boxtimes_{P(z)})=
I^{Q(z^{-1})}_{B_{2p}(\mathcal{Y}_{\boxtimes_{P(z)}, 0}),p}
\]
is a $Q(z^{-1})$-intertwining
map
of type ${W_{1}\boxtimes_{P(z)}
W_{2}\choose W_{1}W_{2}}$. 
We claim that $(W_{1}\boxtimes_{P(z)} W_{2}, Y_{P(z)};
\boxtimes_{Q(z^{-1})})$ is the $Q(z^{-1})$-tensor product of $W_{1}$ and
$W_{2}$.  

In fact, for any $Q(z^{-1})$-product $(W, Y; I)$ of $W_{1}$ and
$W_{2}$,
\[
(\mu_{W_{1}W_{2}}^{W})^{-1}(I)
=I_{B^{-1}_{2p}(\mathcal{Y}_{I,p}^{Q(z^{-1})}), 0}
\]
is a $P(z)$-intertwining map of type ${W\choose W_{1}W_{2}}$ and thus
$(W, Y; (\mu_{W_{1}W_{2}}^{W})^{-1}(I))$ is a $P(z)$-product of
$W_{1}$ and $W_{2}$. Since $(W_{1}\boxtimes_{P(z)} W_{2}, Y_{P(z)};
\boxtimes_{P(z)})$ is the $P(z)$-tensor product of $W_{1}$ and
$W_{2}$, there is a unique morphism of $P(z)$-products {}from
$(W_{1}\boxtimes_{P(z)} W_{2}, Y_{P(z)}; \boxtimes_{P(z)})$ to $(W, Y;
(\mu_{W_{1}W_{2}}^{W})^{-1}(I))$, that is, there exists a unique
module map 
\[
\eta^{P(z)}: W_{1}\boxtimes_{P(z)} W_{2} \to W
\]
such that
\[
(\mu_{W_{1}W_{2}}^{W})^{-1}(I)=\overline{\eta^{P(z)}}
\circ \boxtimes_{P(z)},
\]
or equivalently,
\begin{eqnarray}\label{pzt-qzt-equiv-1}
I&=&\mu_{W_{1}W_{2}}^{W}(\overline{\eta^{P(z)}}
\circ \boxtimes_{P(z)})\nn
&=&\mu_{W_{1}W_{2}}^{W}(\overline{\eta^{P(z)}}\circ 
(\mu_{W_{1}W_{2}}^{W_{1}\boxtimes_{P(z)}
W_{2}})^{-1}(\mu_{W_{1}W_{2}}^{W_{1}\boxtimes_{P(z)}
W_{2}}(\boxtimes_{P(z)})))\nn
&=&\mu_{W_{1}W_{2}}^{W}(\overline{\eta^{P(z)}}\circ 
(\mu_{W_{1}W_{2}}^{W_{1}\boxtimes_{P(z)}
W_{2}})^{-1}(\boxtimes_{Q(z^{-1})})).
\end{eqnarray}

{}From (\ref{mu}) and (\ref{mu-1}), we see that the right-hand side
of (\ref{pzt-qzt-equiv-1}) is determined uniquely by 
\begin{eqnarray}\label{pzt-qzt-equiv-2}
\lefteqn{\langle w', (\mu_{W_{1}W_{2}}^{W}(\overline{\eta^{P(z)}}\circ 
(\mu_{W_{1}W_{2}}^{W_{1}\boxtimes_{P(z)}
W_{2}})^{-1}(\boxtimes_{Q(z^{-1})})))(w_{(1)}\otimes w_{(2)})\rangle}\nn
&&=\langle e^{-zL(1)}w', (\overline{\eta^{P(z)}}\circ 
(\mu_{W_{1}W_{2}}^{W_{1}\boxtimes_{P(z)}
W_{2}})^{-1}(\boxtimes_{Q(z^{-1})}))\nn
&&\quad\quad\quad\quad\quad\quad\quad((e^{z^{-1}L(1)}w_{(1)})
\otimes (e^{-z^{-1}L(1)}
e^{i\pi L(0)}e^{-2(\log z^{-1})L(0)}w_{(2)}))\rangle\nn
&&=\langle e^{-zL(1)}w', \overline{\eta^{P(z)}}
((\mu_{W_{1}W_{2}}^{W_{1}\boxtimes_{P(z)}
W_{2}})^{-1}(\boxtimes_{Q(z^{-1})})\nn
&&\quad\quad\quad\quad\quad\quad\quad((e^{z^{-1}L(1)}w_{(1)})
\otimes (e^{-z^{-1}L(1)}
e^{i\pi L(0)}e^{-2(\log z^{-1})L(0)}w_{(2)})))\rangle\nn
&&=\langle (\eta^{P(z)})'(e^{-zL(1)}w'), 
(\mu_{W_{1}W_{2}}^{W_{1}\boxtimes_{P(z)}
W_{2}})^{-1}(\boxtimes_{Q(z^{-1})})\nn
&&\quad\quad\quad\quad\quad\quad\quad((e^{z^{-1}L(1)}w_{(1)})
\otimes (e^{-z^{-1}L(1)}
e^{i\pi L(0)}e^{-2(\log z^{-1})L(0)}w_{(2)}))\rangle\nn
&&=\langle (\eta^{P(z)})'(w'), 
\boxtimes_{Q(z^{-1})}(w_{(1)}
\otimes w_{(2)})\rangle\nn
&&=\langle w',  (\overline{\eta^{P(z)}}\circ \boxtimes_{Q(z^{-1})})
(w_{(1)}\otimes w_{(2)})\rangle
\end{eqnarray}
for $w_{(1)}\in W_{1}$, $w_{(2)}\in W_{2}$ and $w'\in W'$. {}From 
(\ref{pzt-qzt-equiv-1}) and (\ref{pzt-qzt-equiv-2}), we
see that 
\begin{equation}\label{pzt-qzt-equiv-3}
I=\overline{\eta^{P(z)}}\circ \boxtimes_{Q(z^{-1})}.
\end{equation}

We also need to show the uniqueness---that any module map $\eta:
W_{1}\boxtimes_{P(z)} W_{2}\to W$ such that $I=\overline{\eta}\circ
\boxtimes_{Q(z^{-1})}$ must be equal to $\eta^{P(z)}$.  For this,
it is sufficient to show that $\eta_1 = 0$, where
\[
\eta_1 = \eta^{P(z)} - \eta,
\]
given that 
\[
\overline{\eta_1}(w_{(1)}\boxtimes_{Q(z^{-1})}w_{(2)})=0
\]
for $w_{(1)}\in W_{1}$ and $w_{(2)}\in W_{2}$.  But for $w'\in
(W_{1}\boxtimes_{P(z)}W_{2})'$
\[
\langle e^{zL(1)} w', \overline{\eta_1}
(w_{(1)}\boxtimes_{Q(z^{-1})}w_{(2)})\rangle=0,
\]
so that
\[
\langle e^{zL(1)} {\eta_1}'(w'), 
w_{(1)}\boxtimes_{Q(z^{-1})}w_{(2)}\rangle
=\langle {\eta_1}' (e^{zL(1)} w'), 
w_{(1)}\boxtimes_{Q(z^{-1})}w_{(2)}\rangle = 0.
\]
{}From the definition of $\boxtimes_{Q(z)}$ and (\ref{mu}), we have
\begin{eqnarray*}
\lefteqn{\langle 
e^{zL(1)} {\eta_1}' (w'), w_{(1)}\boxtimes_{Q(z^{-1})}w_{(2)}\rangle}\nn
&&=\langle {\eta_1}' (w'), 
(e^{z^{-1}L(1)}w_{(1)})\boxtimes_{P(z)}(e^{-z^{-1}L(1)}
e^{i\pi L(0)}e^{-2(\log z^{-1})L(0)}w_{(2)})\rangle, \nn
\end{eqnarray*}
and thus
\begin{equation}\label{span2}
\langle {\eta_1}'(w'), 
(e^{z^{-1}L(1)}w_{(1)})\boxtimes_{P(z)}(e^{-z^{-1}L(1)}
e^{i\pi L(0)}e^{-2(\log z^{-1})L(0)}w_{(2)})\rangle=0.
\end{equation}
Since  $e^{z^{-1}L(1)}$ and $e^{-z^{-1}L(1)}
e^{i\pi L(0)}e^{-2(\log z^{-1})L(0)}$ are invertible operators
on $W_{1}$ and $W_{2}$, (\ref{span2}) for all $w_{(1)}\in W_{1}$, 
$w_{(2)}\in W_{2}$ is equivalent to 
\[
\langle {\eta_1}'(w'), w_{(1)}\boxtimes_{P(z)}w_{(2)}\rangle=0
\]
for all $w_{(1)}\in W_{1}$, $w_{(2)}\in W_{2}$.  Thus by
Proposition \ref{span},
\[
{\eta_1}'(w')=0
\]
for all homogeneous $w'$ and hence for all $w'$, showing that indeed
$\eta_1 = 0$ and proving the uniqueness of $\eta$.  Thus
$(W_{1}\boxtimes_{P(z)} W_{2}, Y_{P(z)}; \boxtimes_{Q(z^{-1})})$ is
the $Q(z^{-1})$-tensor product of $W_{1}$ and $W_{2}$.

Conversely, by essentially reversing these arguments we see that if
the $Q(z^{-1})$-tensor product of $W_{1}$ and $W_{2}$ exists, then so
does the $P(z)$-tensor product.  \epfv

{}From Theorem \ref{pz-qz-1} and Proposition \ref{4.19}, we
immediately obtain:

\begin{corol}\label{pz-qz}
Let $W_{1}$ and $W_{2}$ be objects of a full subcategory $\mathcal{C}$
of either $\mathcal{M}_{sg}$ or $\mathcal{GM}_{sg}$. Then the
$P(z)$-tensor product of $W_{1}$ and $W_{2}$ exists if and only if the
$Q(z)$-tensor product of $W_{1}$ and $W_{2}$ exists. \epf
\end{corol}

\begin{rema}{\rm
{}From the proof we see that as $V$-modules, $W_{1}\boxtimes_{P(z)}
W_{2}$ and $W_{1}\boxtimes_{Q(z^{-1})} W_{2}$ are equivalent, but the
main issue is that the intertwining maps $\boxtimes_{P(z)}$ and
$\boxtimes_{Q(z^{-1})}$, which encode the geometric information, are
very different; as $V$-modules {\it only}, $W_{1}\boxtimes_{P(z)}
W_{2}$ and $W_{1}\boxtimes_{Q(z)} W_{2}$ are equivalent.  Compare this
with Remark \ref{intwmapdependsongeomdata}.}
\end{rema}

\newpage

\setcounter{equation}{0}
\setcounter{rema}{0}

\section{Constructions of $P(z)$- and
$Q(z)$-tensor products}

We now generalize the constructions of $P(z)$- and $Q(z)$-tensor
products in \cite{tensor1}--\cite{tensor3} to the setting of the
present work.  In the earlier work \cite{tensor1}-\cite{tensor3} of
the first two authors, the $Q(z)$-tensor product of two modules was
studied and developed first, in \cite{tensor1} and \cite{tensor2}. The
$P(z)$-tensor product was then studied systematically in
\cite{tensor3}, and many proofs for the $P(z)$ case were given by
using the results established for the $Q(z)$ case in \cite{tensor1}
and \cite{tensor2}, rather than by carrying out the subtle arguments
in the $P(z)$ case itself, arguments that are similar to (but
different {}from) those in the $Q(z)$ case.  In the present section and
the next section, instead of following this approach of
\cite{tensor1}--\cite{tensor3}, we shall construct the $P(z)$-tensor
product and $Q(z)$-tensor product of two modules independently. In
particular, even for the finitely reductive case carried out in
\cite{tensor1}--\cite{tensor3}, some of the present results and proofs
of the main theorems are completely new.  One new result is
Proposition \ref{pz-comm} below, which was not stated or proved (or
needed) in the finitely reductive case in \cite{tensor3}. This is
proved below, by a direct argument in the $P(z)$ setting, rather than
by the use of the $Q(z)$ structure. Theorems \ref{comp=>jcb},
\ref{stable}, \ref{6.1} and \ref{6.2} formulated below will be proved
in the next section.  The proofs of Theorems \ref{comp=>jcb} and
\ref{stable} are new, even in the finitely reductive case.  Recall
Assumption \ref{assum}.

\subsection{Affinizations of vertex algebras and the
opposite-operator map}

Just as in \cite{tensor1}--\cite{tensor3}, 
we shall use the Jacobi identity as a motivation 
to construct tensor products of
(generalized) $V$-modules in a suitable category. To do this, 
we need to study  various
``affinizations'' of a vertex algebra with respect to certain
algebras and vector spaces of formal Laurent series and formal 
rational functions.  The treatment of these matters below is 
very similar to that 
in \cite{tensor1}, but here we must take into account 
the gradings by $A$ and $\tilde{A}$. Here, as in Section 2 above,  
we are replacing the symbol $*$ for the 
``opposite-operator map''  in \cite{tensor1} by $o$. 
In Subsections 5.2 and 5.3 below, we will be using the material 
in this subsection to construct certain actions $\tau_{P(z)}$
and $\tau_{Q(z)}$, in order to construct $P(z)$- and $Q(z)$-tensor 
products.

Let $(W, Y_{W})$ be a generalized $V$-module. We adjoin the
formal variable $t$ to our list of commuting formal variables. This
variable will play a special role.
Consider the vector spaces $$V[t,t^{-1}] = V \otimes {\C}[t,t^{-1}]
\subset V \otimes {\C} ((t)) \subset V\otimes {\C} [[t,t^{-1}]]
\subset V[[t,t^{-1}]] $$ (note carefully the distinction between the
last two, since $V$ is typically infinite-dimensional) and $W \otimes
{\C} \{t\}
\subset W \{t\}$
(recall (\ref{formalserieswithcomplexpowers})). The linear map
\begin{eqnarray}\label{tauW}
\tau_W: V[t,t^{-1}] &\to& \mbox{End} \;W\nno\\
v \otimes t^n &\mapsto&
v_n
\end{eqnarray}
($v\in V$, $n\in {\Z}$) extends canonically to 
\begin{eqnarray}\label{tauw}
\tau_W:& V \otimes {\C}((t)) &\to\; \mbox{End} \;W \nno\\
&v \otimes {\dps \sum_{n > N}} a_nt^n
&\mapsto\; \sum_{n > N} a_nv_n
\end{eqnarray}
 (but not to $V((t))$), in view of (\ref{set:wtvn}) and Assumption 
\ref{assum}. 
It further extends canonically to 
\begin{equation}
\tau_W: (V
\otimes {\C}((t)))[[x,x^{-1}]] \to (\mbox{End} \;W)[[x,x^{-1}]],
\end{equation}
where of course $(V
\otimes {\C}((t)))[[x,x^{-1}]]$ can be viewed as the subspace of
$V[[t,t^{-1},x,x^{-1}]]$ such that the coefficient of each power of
$x$ lies in $V \otimes {\C}((t))$.  

Let $v \in V$ and define the ``generic vertex operator''
\begin{equation}\label{3.4}
Y_t(v,x)  = \sum_{n \in {\Z}}(v \otimes
t^n)x^{-n-1} \in  (V \otimes {\C}[t,t^{-1}])[[x,x^{-1}]].
\end{equation}  
Then 
\begin{eqnarray}\label{3.5}
Y_t(v,x) & = & v\otimes x^{-1} \delta \left(\frac{t}{x}\right)\nno\\
 & = & v \otimes t^{-1} \delta \left(\frac{x}{t}\right)\nno\\
 & \in & V \otimes {\C}
[[t,t^{-1},x,x^{-1}]]\nno\\
 & (\subset & V[[t,t^{-1},x,x^{-1}]]),
\end{eqnarray} 
and the linear map 
\begin{eqnarray}\label{3.6}
V& \to &V\otimes {\C}[[t,t^{-1},x,x^{-1}]]\nno\\
v &\mapsto & Y_t(v,x)
\end{eqnarray}
 is simply the  map  given by tensoring by the
``universal element'' $x^{-1} \delta
\left(\frac{t}{x}\right)$.  We have 
\begin{equation}\label{3.7}
\tau_W(Y_t(v,x)) =
Y_W(v,x).
\end{equation}
 For all $f(x) \in {\C} [[x,x^{-1}]]$,
$f(x)Y_t(v,x)$ is defined and 
\begin{equation}
f(x)Y_t(v,x) = f(t) Y_t(v,x).
\end{equation}
In case $f(x) \in {\C}((x))$, then $\tau_{W}(f(x)Y_{t}(v, x))$ is
also defined, and
\begin{equation}\label{3.9}
f(x)Y_W(v,x)  =  f(x)\tau_W(Y_t(v,x)) = 
\tau_W(f(x)Y_t(v,x)) =  \tau_W(f(t)Y_t(v,x)).
\end{equation}
The expansion coefficients, in powers of $x$, of $Y_t(v,x)$
span $v \otimes {\C}[t,t^{-1}],$ the
$x$-expansion coefficients
of $Y_W(v,x)$ span $\tau_W(v \otimes {\C}[t,t^{-1}])$ and 
for $f(x)\in \C[[x, x^{-1}]]$, the
$x$-expansion coefficients of $f(x)Y_t(v,x)$ span $v \otimes f(t)
{\C} [t,t^{-1}]$. In case $f(x) \in {\C}((x))$, the 
$x$-expansion coefficients of $f(x)Y_W(v,x)$ span $\tau_W(v \otimes
f(t) {\C}[t,t^{-1}])$.
 
Using this viewpoint, we shall examine each of the three terms in
the Jacobi identity (\ref{log:jacobi}) in the definition of logarithmic
intertwining operator. First we consider the formal Laurent series in
$x_{0}$, $x_{1}$, $x_{2}$ and $t$ given by
\begin{eqnarray}\label{3.10}
{\dps x^{-1}_2 \delta\left(\frac{x_1-x_0}{x_2}\right)
Y_t(v,x_0) }&=& {\dps x^{-1}_1\delta \left(\frac{x_2+x_0}{x_1}\right)
Y_t(v,x_0)}\nno\\
&=&{\dps  v \otimes x^{-1}_1\delta\left(\frac{x_2 +
t}{x_1}\right)x^{-1}_0 \delta\left(\frac{t}{x_0}\right)}
\end{eqnarray}
(cf. the right-hand side of (\ref{log:jacobi})). The
expansion coefficients in powers of $x_0$, $x_1$ and $x_2$ of (\ref{3.10}) 
span just the space $v \otimes {\C}[t,t^{-1}]$. However,
the expansion coefficients in $x_0$ and $x_1$ only (but not in $x_{2}$) of 
\begin{eqnarray}\label{3.11}
x^{-1}_1\delta\left(\frac{x_{2}+x_0}{x_1}\right) Y_t(v,x_0)
&  =& v \otimes x^{-
1}_1\delta\left(\frac{x_{2}+t}{x_1}\right)x^{-
1}_0\delta\left(\frac{t}{x_0}\right)\nno\\ 
&  =& v \otimes
\left(\sum_{m \in {\Z}} (x_{2} + t)^m x^{-m-
1}_1\right)\left(\sum_{n \in {\Z}}t^n x^{-n-
1}_0\right)
\end{eqnarray}
 span $$v \otimes \iota_{x_{2},t}{\C}[t,t^{- 1}, x_{2}+t,
(x_{2}+t)^{-1}]\subset v \otimes {\C}[x_{2}, x_{2}^{-1}]((t)),$$
where $\iota_{x_{2}, t}$ is the operation of expanding a formal rational 
function in the indicated algebra as a formal Laurent series
involving only finitely many negative powers of $t$ (cf.
the notation $\iota_{12}$, etc., considered at the end of Section 2). 
We shall use similar
$\iota$-notations below. Specifically, the coefficient of 
$x^{-n-1}_0 x^{-m-1}_1$ $(m,n
\in {\Z})$ in (\ref{3.11}) is $v \otimes (x_{2} + t)^m t^n$.  

We may
specialize $x_2 \mapsto z\in {\C}^{\times}$, and (\ref{3.11}) becomes
\begin{eqnarray}\label{3.12}
z^{-1}\delta\left(\frac{x_{1}-x_0}{z}\right) Y_t(v,x_0)
&=&x^{-1}_1\delta\left(\frac{z+x_0}{x_1}\right) Y_t(v,x_0)\nno\\
&  =& v \otimes x^{-
1}_1\delta\left(\frac{z+t}{x_1}\right)x^{-
1}_0\delta\left(\frac{t}{x_0}\right)\nno\\ 
&  =& v \otimes
\left(\sum_{m \in {\Z}} (z + t)^m x^{-m-
1}_1\right)\left(\sum_{n \in {\Z}}t^n x^{-n-
1}_0\right).
\end{eqnarray}
The coefficient of $x^{-n-1}_0 x^{-m-1}_1$ $(m,n \in {\Z})$ in (\ref{3.12})
is $v \otimes (z + t)^m t^n\in V\otimes {\C}((t))$, and these
coefficients span
\begin{equation}\label{3.13}
v\otimes {\C}[t, t^{-1}, (z+t)^{-1}]\subset v\otimes {\C}((t)).
\end{equation}
Our $Q(z)$-tensor product construction in Subsection 5.3 below
will be based on a certain action of the 
space $V\otimes {\C}[t, t^{-1}, (z+t)^{-1}]$, and the description of 
this space as the span of the coefficients of the expression (\ref{3.12}) (as
$v\in V$ varies) will be very useful.

Now consider 
\begin{eqnarray}\label{3.14}
x^{-1}_0 \delta\left(\frac{-x_2 +
x_1}{x_0}\right) Y_t(v,x_1)& =&  v
\otimes x^{-1}_0 \delta\left(\frac{-x_2 + t}{x_0}\right) x^{-1}_1
\delta \left(\frac{t}{x_1}\right)\nno\\
&=&{\dps    v\otimes \biggr(\sum_{n \in {\Z}} (-x_2 + t)^n x^{-n-
1}_0\biggr)\biggr(\sum_{m \in {\Z}} t^m x^{-m-1}_1\biggr)}
\end{eqnarray}
(cf. the second term on the left-hand side of (\ref{log:jacobi})). 
The expansion coefficients in powers of $x_0$ and $x_1$ (but not $x_2$)
span 
$$v \otimes \iota_{x_{2}, t}{\C}[t,t^{-1},-x_2 + t, (-x_2 +
t)^{-1}],$$ and in fact the coefficient of $x^{-n-1}_0 x^{-m-1}_1$
$(m,n \in {\Z})$ in (\ref{3.14}) is $v \otimes (-x_2 + t)^n t^m$. Again
specializing $x_2 \mapsto z\in {\C}^{\times}$, we obtain
\begin{eqnarray}\label{3.15}
x^{-1}_0 \delta \left(\frac{-z+x_1}{x_0}\right)
Y_t(v,x_1) 
&=& v\otimes x^{-1}_0 \delta \left(\frac{-
z+t}{x_0}\right) x^{-1}_1 \delta \left(\frac{t}{x_1}\right)\nn
& =& v \otimes \biggr(\sum_{n \in {\Z}} (-z+t)^n x^{-n- 1}_0\biggr)
\biggr(\sum_{m \in
{\Z}}t^m x^{-m-1}_1\biggr).
\end{eqnarray}
 The coefficient of $x^{-n-1}_0 x^{- m-1}_1$ $(m,n \in {\Z})$ in
(\ref{3.15}) is $v \otimes (-z+t)^n t^m$, and these coefficients span
\begin{equation}\label{3.16}
v\otimes {\C}[t, t^{-1}, (-z+t)^{-1}]\subset v\otimes {\C}((t)).
\end{equation}

Finally, consider 
\begin{eqnarray}
{\dps x^{-1}_0 \delta \left(\frac{x_1 -
x_2}{x_0}\right) Y_t (v,x_1)}
&=&{\dps v \otimes x^{-1}_0
\delta \left(\frac{t-x_2}{x_0}\right) x^{-1}_1 \delta
\left(\frac{t}{x_1}\right).}
\end{eqnarray}
The coefficient of $x_{0}^{-n-1}x_{1}^{-m-1}$ ($m, n\in {\Z}$) is
$v\otimes (t-x_{2})^{n}t^{m}$, and these expansion cofficients  span
$$v\otimes \iota_{t, x_{2}}{\C}[t, t^{-1}, t-x_{2}, (t-x_{2})^{-1}].$$
If we again specialize $x_2
\mapsto z$, we get
\begin{equation}\label{3.18}
x^{-1}_0\delta\left(\frac{x_1-z}{x_0}\right)Y_t(v,x_1) = v\otimes
x_{0}^{-1}\delta\left(\frac{t-z}{x_{0}}\right)x_{1}^{-1}\delta\left(
\frac{t}{x_{1}}\right),
\end{equation}
whose coefficient of $x^{-n-1}_0 x^{-m-1}_1$  is $v
\otimes (t - z)^n t^m$.  These coefficients span 
\begin{equation}\label{3.19}
v\otimes {\C}[t, t^{-1}, (t-z)^{-1}]\subset v\otimes {\C}((t^{-1}))
\end{equation}
(cf. (\ref{3.13}), (\ref{3.16})).

In the construction of $P(z)$-tensor products in Subsection 5.2, we 
shall also need the following expression, 
which is slightly different {}from what we have analyzed 
above:
\begin{eqnarray}\label{y-t-delta}
{\dps x^{-1}_0 \delta \left(\frac{x_1^{-1} -
x_2}{x_0}\right) Y_t (v, x_1)}
&=&{\dps v \otimes x^{-1}_0
\delta \left(\frac{t^{-1}-x_2}{x_0}\right) x^{-1}_1 \delta
\left(\frac{t}{x_1}\right).}
\end{eqnarray}
The coefficient of $x_{0}^{-n-1}x_{1}^{-m-1}$ ($m, n\in {\Z}$) is
$v\otimes (t^{-1}-x_{2})^{n}t^{m}$, and these expansion cofficients  span
$$v\otimes \iota_{t^{-1}, x_{2}}{\C}[t, t^{-1}, t^{-1}-x_{2}, (t^{-1}-x_{2})^{-1}].$$
If we again specialize $x_2
\mapsto z$, we get
\begin{equation}\label{3.18-1}
x^{-1}_0\delta\left(\frac{x_1^{-1}-z}{x_0}\right)Y_t(v,x_1) = v\otimes
x_{0}^{-1}\delta\left(\frac{t^{-1}-z}{x_{0}}\right)x_{1}^{-1}\delta\left(
\frac{t}{x_{1}}\right),
\end{equation}
whose coefficient of $x^{-n-1}_0 x^{-m-1}_1$  is $v
\otimes (t^{-1} - z)^n t^m$.  These coefficients span 
\begin{equation}\label{3.19-1}
v\otimes {\C}[t, t^{-1}, (t^{-1}-z)^{-1}]
=v\otimes {\C}[t, t^{-1}, (z^{-1}-t)^{-1}]\subset v\otimes {\C}((t)).
\end{equation}
Our $P(z)$-tensor product construction in Subsection 5.2 below
will be based on a certain action of the 
space $V\otimes {\C}[t, t^{-1}, (z^{-1}-t)^{-1}]$.

Later we shall evaluate the identity (\ref{log:jacobi}) 
on the elements of the
contragredient module $W'_{3}$. This will allow us to convert the
expansion (\ref{3.19}) into an expansion in positive powers of $t$. 
It will be useful to examine the
notions of opposite and contragredient vertex operators  more
closely (recall Section 2, in particular, (\ref{yo})). 

We shall interpret the opposite vertex operator map $Y^{o}_{W}$ by means
of an operation on $V\otimes {\C}[[t, t^{-1}]]$ that will convert
vertex operators into their opposites.  We shall write this
``opposite-operator'' map, in various contexts, as ``$o$.''  The
operation $o$ will be an involution.  We proceed as follows: First we
generalize $Y^{o}$ in the following way: Recall that 
by Assumption \ref{assum}, $L(1)$ acts nilpotently on 
any element $v\in V$. In particular, $e^{xL(1)}v$ is a polynomial
in the formal variable $x$. Given any vector space $U$
and any linear map
\begin{eqnarray}
Z(\cdot,x):\ V & \rightarrow &  U[[x,x^{-
1}]]\;\;\biggr(=\prod_{n \in {\Z}} U \otimes x^n\biggr)\nno\\ 
v & \mapsto & Z(v,x)
\end{eqnarray}
{from} $V$ into $U[[x,x^{-1}]]$
(i.e., given any family of linear maps {from} $V$ into the spaces
$U \otimes x^n$), we define $Z^o(\cdot,x):\ V \to U[[x,x^{-1}]]$
by 
\begin{equation}\label{3.21}
Z^o(v,x) = Z(e^{xL(1)}(-x^{-2})^{L(0)}v,x^{-1}),
\end{equation}
 where
we use the obvious linear map $Z(\cdot,x^{-1}):\ V \to U[[x,x^{-
1}]],$ and where we extend $Z(\cdot,x^{-1})$ canonically to a
linear map $Z(\cdot,x^{-1}):\ V[x,x^{-1}] \to U[[x,x^{-1}]].$
Then by formula (5.3.1) in \cite{FHL}  (the proof of Proposition
5.3.1), we have
\begin{eqnarray}
Z^{oo}(v,x) & = & Z^o(e^{xL(1)}(-x^{-
2})^{L(0)}v,x^{-1})\nno\\ 
& = & Z(e^{x^{-1}L(1)}(-
x^2)^{L(0)}e^{xL(1)}(-x^{-2})^{L(0)}v,x)\nno\\
&=& Z(v,x).
\end{eqnarray}
That is,
\begin{equation}
Z^{oo}(\cdot,x) = Z(\cdot,x).
\end{equation}
 Moreover, if
$Z(v,x) \in U((x))$, then $Z^o(v,x) \in U((x^{-1}))$ and vice
versa.  
 
Now we expand $Z(v,x)$ and $Z^o(v,x)$ in components. Write
\begin{equation}
Z(v,z) = \sum_{n \in {\Z}}v_{(n)} x^{-n-1},
\end{equation}
 where for all $n\in {\Z}$,
\begin{eqnarray}
V & \to & U\nno\\
v & \mapsto & v_{(n)}
\end{eqnarray}
 is a linear map depending on $Z(\cdot, x)$ (and in fact, as
$Z(\cdot, x)$ varies, these linear maps are arbitrary).
Also write 
\begin{equation}
Z^o(v,x)
= \sum_{n \in {\Z}} v^o_{(n)}x^{-n-1}
\end{equation}
 where
\begin{eqnarray}
V & \to & U\nno\\
v & \mapsto & v^o_{(n)}
\end{eqnarray}
 is a linear map depending on $Z(\cdot, x)$.  We shall compute $v_{(n)}^{o}$.
First note that
\begin{equation}\label{3.32}
\sum_{n
\in {\Z}}v^o_{(n)}x^{-n-1} = \sum_{n \in {\Z}}(e^{xL(1)}(-
x^{-2})^{L(0)}v)_{(n)}x^{n+1}.
\end{equation}
  For convenience, suppose that $v
\in V_{(h)}$, for $h \in {\Z}$.  Then the right-hand side of (\ref{3.32}) 
is equal to
\begin{eqnarray}
\lefteqn{(-1)^h\sum_{n \in {\Z}}(e^{xL(1)}v)_{(-
n)}x^{-n+1-2h}}\nno\\
&&=  (-1)^h \sum_{n \in {\Z}} \sum_{m \in 
{\N}} \frac{1}{m!}(L(1)^m v)_{(-n)}x^{m-n+1-2h}\nno\\
&&= 
(-1)^h\sum_{m
\in {\N}} \frac{1}{m!} \sum_{n \in {\Z}}(L(1)^mv)_{(-n-m-
2+2h)}x^{-n-1}, 
\end{eqnarray} 
that is, 
\begin{equation}\label{vo}
v^o_{(n)} = (-
1)^h\sum_{m \in {\N}}\frac{1}{m!} (L(1)^mv)_{(-n-m-2+2h)}.
\end{equation}
(Recall that by Assumption \ref{assum}, $L(1)^{m}v=0$ when  $m$ 
is sufficiently large, so that these expressions 
are well defined.)
For $v\in V$ not necessarily homogeneous, $v^o_{(n)}$ is given by the
appropriate sum of such expressions.

Now consider the special case where $U = V \otimes {\C}[t,t^{-1}]$ and
where $Z(\cdot, x)$ is the ``generic'' linear map
\begin{eqnarray}
Y_t(\cdot,x): V & \rightarrow & (V \otimes
{\C}[t,t^{-1}])[[x,x^{-1}]]\nno\\
v & \mapsto & Y_t(v,x) = \sum_{n
\in {\Z}}(v \otimes t^n)x^{-n-1}
\end{eqnarray} 
(recall (\ref{3.4})), i.e.,
\begin{equation}
v_{(n)} = v \otimes t^n.
\end{equation}
Then for $v\in V_{(h)}$, 
\begin{equation}
v^o_{(n)} = (-1)^h\sum_{m \in {\N}}
\frac{1}{m!}((L(1))^mv) \otimes t^{-n-m-2+2h}
\end{equation}
 in this case.
 
This motivates defining an $o$-operation on $V \otimes {\C}[t,t^{-1}]$
as follows: For any $n, h\in {\Z}$ and $v\in V_{(h)}$, define
\begin{equation}\label{3.38}
(v\otimes t^n)^o = (-1)^h\sum_{m \in {\N}} \frac{1}{m!}
(L(1)^mv) \otimes t^{-n-m-2+2h}\in V \otimes {\C}[t,t^{-
1}],
\end{equation} 
and extend by linearity to $V \otimes {\C}[t,t^{-1}]$. 
That is, $(v \otimes t^n)^o = v^o_{(n)}$ for the special case $Z(\cdot, x)
= Y_t(\cdot, x)$ discussed above. (Note that for general $Z$, 
we cannot expect to
be able to define an analogous $o$-operation on $U$.)  Also consider the map
\begin{eqnarray}
Y^o_t(\cdot,x) = (Y_t(\cdot,x))^o: V & \rightarrow & (V \otimes
{\C}[t,t^{- 1}])[[x,x^{-1}]]\nno\\ v & \mapsto & Y^o_t(v,x) = \sum_{n
\in {\Z}}(v \otimes t^n)^ox^{-n-1}.
\end{eqnarray} 
Then for general
$Z(\cdot, x)$ as above, we can define a linear map
\begin{eqnarray}\label{3.40}
\varepsilon_{Z}: \ V \otimes {\C}[t,t^{-1}]
& \rightarrow & U\nno\\
v \otimes t^n & \mapsto & v_{(n)}
\end{eqnarray}
(``evaluation with respect to $Z$''), i.e.,
\begin{equation}
\varepsilon_{Z}:\ Y_t(v,x) \mapsto Z(v,x),
\end{equation}
 and a linear map
\begin{eqnarray}
\varepsilon^o_{Z}:\ V \otimes {\C}[t,t^{-1}]
& \rightarrow & U\nno\\
v \otimes t^n & \mapsto & v^o_{(n)},
\end{eqnarray}
 i.e., 
\begin{equation}
\varepsilon^o_{Z}:\ Y_t(v,x) \mapsto
Z^o(v,x).
\end{equation}
 Then 
\begin{equation}
\varepsilon^o_{Z} = \varepsilon_{Z} \circ o,
\end{equation}
that is,
\begin{equation}
\varepsilon_{Z}(Y^o_t(v,x)) = Z^o(v,x),
\end{equation}
or equivalently, the diagram
\begin{eqnarray}
Y_t(v,x) & \stackrel{\varepsilon_{Z}}{\longmapsto}
& Z(v,x)\nno\\
{o} \bar{\downarrow}\hspace{1.5em} & & \hspace{1.5em}
\bar{\downarrow} \;(Z(\cdot, x) \mapsto Z^o(\cdot, x))\nno\\
Y^o_t(v,x) &
\stackrel{\varepsilon_{Z}}{\longmapsto} & Z^o(v,x) 
\end{eqnarray} 
commutes. Note that the components $v^o_{(n)}$ of $Z^{o}(v, x)$ depend on
all the components $v_{(n)}$ of $Z(v,z)$ (for
arbitrary $v$), whereas the component $(v \otimes t^n)^o$ of
$Y^o_t(v,z)$ can be defined generically and abstractly; $(v
\otimes t^n)^o$ depends linearly on $v \in V$ alone.
 
     Since in general $Z^{oo}(v, x) = Z(v, x)$, we know that 
\begin{equation}
Y^{oo}_t(v, x) =Y_t(v, x)
\end{equation}
as a special case, and in particular (and equivalently),
\begin{equation}
(v \otimes t^n)^{oo} = v \otimes t^n
\end{equation}
for all $v\in V$ and $n\in {\Z}$. Thus $o$ is an involution of $V
\otimes {\C}[t,t^{-1}]$.
 
     Furthermore, the involution $o$ of $V \otimes {\C}[t,t^{-1}]$
extends canonically to a linear map
$$V \otimes {\C}[[t,t^{-1}]] \stackrel{o}{\rightarrow} V \otimes
{\C}[[t,t^{-1}]].$$
In fact, consider the restriction of $o$ to $V=V \otimes t^0$:
\begin{eqnarray}
V &\stackrel{o}{\rightarrow} &V \otimes
{\C}[t,t^{-1}]\nno\\
v &\mapsto& v^o = (-1)^h \sum_{m \in {\N}}
\frac{1}{m!}(L(1)^m v) \otimes t^{-m-2+2h},
\end{eqnarray}
 extended by
linearity {from} $V_{(h)}$ to $V$.  Then for $v\in V$, we may write 
\begin{equation}\label{vo1}
v^{o}=e^{t^{-1}L(1)}(-t^{2})^{L(0)}vt^{-2}.
\end{equation}
Also, for $v\in V$ and $n\in {\Z}$,
\begin{equation}
(v \otimes t^n)^o = v^ot^{-n},
\end{equation}
 and it is clear that $o$ extends to $V\otimes {\C}[[t, t^{-1}]]$: 
For $f(t) \in {\C}[[t,t^{-1}]],$ 
\begin{equation}
(v \otimes
f(t))^o = v^of(t^{-1}).
\end{equation}

To see that $o$ is an involution of this larger space, first note that 
\begin{equation}
v^{oo} = v
\end{equation}
(although $v^{o}\not\in V$ in general). (This could of course
alternatively be proved by direct calculation using formula (\ref{3.38}).)
Also, for $g(t) \in {\C}[t,t^{-1}]$ and $f(t) \in {\C}[[t,t^{-1}]],$
\begin{equation}
(v \otimes
g(t)f(t))^o = v^og(t^{-1})f(t^{-1})= (v \otimes
g(t))^of(t^{-1}).
\end{equation}
 Thus for all $x \in V\otimes {\C}[t,t^{-
1}]$ and $f(t) \in {\C}[[t,t^{-1}]],$ 
\begin{equation}
(xf(t))^o = x^of(t^{-
1}).
\end{equation}
  It follows that 
\begin{eqnarray}
(v \otimes f(t))^{oo} &
= & (v^of(t^{-1}))^o\nno\\ 
& = & v^{oo}f(t)\nno\\
& = & vf(t)\nno\\ 
& = & v
\otimes f(t),
\end{eqnarray}
and we have shown that $o$ is an involution of $V
\otimes {\C}[[t,t^{-1}]]$.
We have  
\begin{equation}
o: V \otimes {\C}((t)) \leftrightarrow V
\otimes {\C}((t^{-1})).
\end{equation}

 Note that
\begin{eqnarray}\label{op-y-t}
Y^o_t(v,x) & = & \sum_{n \in {\Z}}(v \otimes
t^n)^ox^{-n-1}\nno\\
& = & v^o \sum_{n \in {\Z}}t^{-n}x^{-n-1}\nno\\
&= & v^ox^{-1} \delta(tx)\nno\\
& = & v^ot \delta(tx)\nno\\
& \in & V
\otimes {\C}[[t,t^{-1},x,x^{-1}]].
\end{eqnarray}
Thus the map $v \mapsto Y^o_t(v,x)$ is the linear map given by
multiplying $v^o$ by the ``universal element'' $t \delta(tx)$ 
(cf. the comment following (\ref{3.6})). By (\ref{vo1}), we also have 
\begin{eqnarray}\label{op-y-t-2}
Y^o_t(v,x)&=&e^{t^{-1}L(1)}(-t^{2})^{L(0)}vt^{-1}\delta(tx)\nn
&=&e^{xL(1)}(-x^{-2})^{L(0)}v\otimes x\delta(tx).
\end{eqnarray}
For all $f(x) \in {\C}[[x,x^{-1}]],$ $f(x)Y^o_t(v,x)$
is defined and 
\begin{eqnarray}
f(x)Y^o_t(v,x) & = & f(t^{-
1})Y^o_t(v,x)\nno\\
& = & v^of(t^{-1})t\delta(tx).
\end{eqnarray}

Now we return to the starting point---the original special case: $U =
\mbox{End}\; W$ and $Z(\cdot,z) = Y_W(\cdot,z):\ V \to (\mbox{End}\;
W)[[x,x^{-1}]].$ The corresponding map
\begin{eqnarray}
\varepsilon_{Z} =
\varepsilon_{Y_W}:\ V[t,t^{-1}] & \to & \mbox{End}\; W\nno\\
v \otimes
t^n & \mapsto & v_{(n)}
\end{eqnarray} 
(recall (\ref{3.40})) is just the map
$\tau_W: v \otimes t^n
\mapsto v_n$ (recall (\ref{tauW})), i.e.,
 $v_{(n)} = v_n$ in this case.  Recall that this
map extends canonically to $V \otimes {\C}((t))$.  The map
$\varepsilon^o_{Z}$  is 
$\tau_{W}\circ o:\ V \otimes {\C}[t,t^{-1}] \to \mbox{End}\;
W$ and this map extends canonically to $V \otimes {\C}((t^{-
1}))$. In addition to (\ref{3.7}),  we have 
\begin{equation}\label{tauw-yto}
\tau_W(Y^o_t(v,x))  =  Y^o_W(v,x)
\end{equation} 
($v^o_{(n)} =
v^o_n$ in this case; recall (\ref{yo1})).  In case $f(x) \in {\C}((x^{-1})),$
$$f(x)Y^o_W(v,x) = \tau_W(f(x)Y^o_t(v,x))$$
 is defined and is equal to
$\tau_W(f(t^{-1})Y^o_t(v,z))$ (which is also defined).
 
     The $x$-expansion coefficients of $f(x)Y^o_t(v,x)$, for
$f(x) \in {\C}[[x,x^{-1}]]$, span 
\begin{eqnarray}
v^of(t^{-
1}){\C}[t,t^{-1}] & = & (v{\C}[t,t^{-1}])^of(t^{-1})\nno\\
& = &
(vf(t){\C}[t,t^{-1}])^o.
\end{eqnarray}
The $x$-expansion coefficients of $Y^o_W(v,x)$ span
\begin{eqnarray}
\tau_W(v^o{\C}[t,t^{-1}])&=& \tau_W((v \otimes
{\C}[t,t^{-1}])^o)\nno\\ &=& \tau^o_W(v \otimes {\C}[t,t^{- 1}]).
\end{eqnarray}
In case $f(x) \in {\C}((x^{-1}))$, the $x$-expansion coefficients of
$f(x)Y^o_W(v,x)$ span 
\[
\tau_W(v^of(t^{-1}){\C}[t,t^{-1}])=
\tau^o_W(vf(t){\C}[t,t^{-1}]).
\]
  (Cf. the comments after (\ref{3.9}).)

We shall need spaces of the forms $V\otimes \iota_{+} {\C}[t, t^{-1},
(z+t)^{-1}]$ and 
$V\otimes \iota_{-} {\C}[t, t^{-1},
(z+t)^{-1}]$, where we use the notations
\begin{eqnarray}\label{iota+-}
\iota_+: {\C}(t) & \hookrightarrow & {\C}((t)) \subset
{\C}[[t,t^{-1}]]\nno\\ \iota_-:{\C}(t) & \hookrightarrow &
{\C}((t^{-1})) \subset {\C} [[t,t^{-1}]]
\end{eqnarray}
to denote the operations of expanding a rational function of the
formal variable $t$ in
the indicated directions (as in Section 8.1 of
\cite{FLM2}).  
We shall also 
need certain translation operations, as well as the $o$-operation. 
For $a \in \mathbb{C}$, we define the translation
isomorphism
\begin{eqnarray}
T_a:\ \mathbb{C}(t) & \stackrel{\sim}{\rightarrow}
& \mathbb{C}(t)\nno\\f(t) & \mapsto & f(t+a)
\end{eqnarray} 
and (for our use below) we also set
\begin{equation}
T^\pm_a = \iota_\pm \circ T_a \circ
\iota^{-1}_+:\
\iota_+\mathbb{C}(t) \hookrightarrow \mathbb{C}((t^{\pm1})).
\end{equation}
(Note that the domains of these maps consist of certain series expansions
of formal rational functions rather than of 
formal rational functions themselves.)
The following lemma will be needed for our action $\tau_{P(z)}$ in 
Subsection 5.2 below (we shall sometimes write $o(Y_{t}(v, x_{1}))$ for 
$Y_{t}^{o}(v, x_{1})$, etc.):

\begin{lemma}\label{tauP}
Let $z\in \C^{\times}$. Then
\begin{eqnarray}
\lefteqn{o\left(x_0^{-1}\delta\left(\frac{x^{-1}_1-z}{x_0}\right)
Y_{t}(v,x_1)\right)=x_0^{-1}\delta\left(\frac{x^{-1}_1-z}{x_0}\right)
Y^o_{t}(v,x_1),}\label{ztr1}\\
\lefteqn{(\iota_+\circ\iota_-^{-1}\circ o)\left(x_0^{-1}\delta\left(
\frac{x^{-1}_1 -z}{x_0}\right)Y_{t}(v,x_1)\right)=x_0^{-1}\delta\left(
\frac{z-x^{-1}_1}{-x_0}\right) Y^o_{t}(v,x_1),}\label{ztr2}\\
\lefteqn{(\iota_+\circ T_{z}\circ\iota_-^{-1}\circ o)\left(x_0^{-1}
\delta\left(\frac{x^{-1}_1-z}{x_0}\right) Y_{t}(v,x_1)\right)}\nno\\
&&\hspace{4.5em}=z^{-1}\delta \left(\frac{x^{-1}_1-x_0}{z}\right)
Y_{t}(e^{x_{1}L(1)}(-x_{1}^{-2})^{L(0)}v,x_0).\quad\quad\quad\quad\quad
\label{ztr3}
\end{eqnarray}
\end{lemma}
\pf 
Formula (\ref{ztr1}) is immediate {}from the definition 
of the map $o$ (recall (\ref{3.38})). 
By (\ref{ztr1}), (\ref{op-y-t}) and (\ref{Xx1x2=Xx2x2}), we have 
\begin{eqnarray}
\lefteqn{(\iota_+\circ\iota_-^{-1}\circ o)\left(x_0^{-1}\delta\left(
\frac{x^{-1}_1 -z}{x_0}\right)Y_{t}(v,x_1)\right)}\nn
&&=(\iota_+\circ\iota_-^{-1})\left(x_0^{-1}\delta\left(
\frac{x^{-1}_1 -z}{x_0}\right)Y^{o}_{t}(v,x_1)\right)\nn
&&=(\iota_+\circ\iota_-^{-1})\left(x_0^{-1}\delta\left(
\frac{x^{-1}_1 -z}{x_0}\right)v^{o}t\delta(tx_{1})\right)\nn
&&=(\iota_+\circ\iota_-^{-1})\left(x_0^{-1}\delta\left(
\frac{t -z}{x_0}\right)v^{o}t\delta(tx_{1})\right)\nn
&&=x_0^{-1}\delta\left(
\frac{z-t}{-x_0}\right)v^{o}t\delta(tx_{1})\nn
&&=x_0^{-1}\delta\left(
\frac{z-x^{-1}}{-x_0}\right)v^{o}t\delta(tx_{1})\nn
&&=x_0^{-1}\delta\left(
\frac{z-x^{-1}_1}{-x_0}\right) Y^o_{t}(v,x_1),
\end{eqnarray}
proving (\ref{ztr2}).
For (\ref{ztr3}), note that by (\ref{op-y-t-2}), 
the coefficient of $x_0^{-n-1}$ in the right-hand side of (\ref{ztr1}) is
\begin{eqnarray*}
\lefteqn{(x_1^{-1}-z)^n\left(e^{x_1L(1)}(-x_1^{-2})^{L(0)}v\otimes
x_1\delta\left(\frac t{x_1^{-1}}\right)\right)}\\
&&=(t-z)^n\left(e^{x_1L(1)}(-x_1^{-2})^{L(0)}v\otimes
x_1\delta\left(\frac t{x_1^{-1}}\right)\right).
\end{eqnarray*}
Acted by $\iota_+\circ T_z\circ\iota_-^{-1}$, this becomes
\begin{eqnarray*}
\lefteqn{t^n\left(e^{x_1L(1)}(-x_1^{-2})^{L(0)}v\otimes
x_1\delta\left(\frac {z+t}{x_1^{-1}}\right)\right)}\\
&&=z^{-1}\delta\left(\frac {x_1^{-1}-t}{z}\right)\left(e^{x_1L(1)}
(-x_1^{-2})^{L(0)}v\otimes t^n\right),
\end{eqnarray*}
which by (\ref{3.5}) is the  coefficient of $x_0^{-n-1}$ in the right-hand side of
(\ref{ztr3}).  
\epfv

We shall be interested in 
\begin{equation}
T^\pm_{-z}: \iota_+\mathbb{C}[t,t^{-1},(z +
t)^{-1}]
\hookrightarrow \mathbb{C}((t^{\pm1})),
\end{equation}
where $z$ is an arbitrary nonzero complex number, as above.
The images of these two maps are $\iota_{\pm}\mathbb{C}[t,t^{-1},(z-t)^{-1}]$.
 
Extend the maps $T^\pm_{-z}$ to linear isomorphisms 
\begin{equation}\label{Tpm-z}
T^\pm_{-z}:\ V \otimes
\iota_+\mathbb{C}[t,t^{-1},(z+t)^{-1}] \stackrel{\sim}{\rightarrow}
V
\otimes \iota_\pm\mathbb{C}[t,t^{-1},(z-t)^{-1}]
\end{equation}
 given by $1\otimes T^\pm_{-z}$ with $T^\pm_{-z}$ as defined above.
Note that the domain of these two maps is described by 
(\ref{3.12})--(\ref{3.13}),
that the image of the map $T^{+}_{-z}$ is described by 
(\ref{3.15})--(\ref{3.16})
and that the image of the map $T^{-}_{-z}$ is described by
(\ref{3.18})--(\ref{3.19}).

We have the two mutually inverse maps
\begin{eqnarray}
V \otimes \iota_-\mathbb{C}[t,t^{-
1},(z-t)^{-1}] & \stackrel{o}{\rightarrow}& 
V \otimes \iota_+\mathbb{C}[t,t^{-1},(z^{-1}-t)^{-1}]\nno\\
v \otimes f(t) &\mapsto &v^of(t^{-1})  
\end{eqnarray}
 and 
\begin{eqnarray}
V \otimes \iota_+{\C}[t,t^{-1},(z^{-1}-t)^{-1}] &
\stackrel{o}{\rightarrow}& V \otimes
\iota_-\mathbb{C}[t,t^{-1},(z-t)^{-1}]\nno\\ v \otimes f(t) & \mapsto
&v^of(t^{-1}),
\end{eqnarray}
which are both isomorphisms. We form the composition
\begin{equation}\label{To-z}
T^o_{-z} = o
\circ T^-_{-z}
\end{equation}
to obtain another isomorphism
\[
T^o_{-z}: V \otimes \iota_+{\C}[t,t^{-1},(z+t)^{-1}]
\stackrel{\sim}{\rightarrow} V \otimes
\iota_+\mathbb{C}[t,t^{-1},(z^{-1}-t)^{-1}].
\]
The maps $T^{+}_{-z}$ and $T^o_{-z}$ will be the main ingredients of
our action $\tau_{Q(z)}$ (see Subsection 5.3 below). The following
result asserts that $T^{+}_{-z}$, $T^{-}_{-z}$ and $T^o_{-z}$
transform the expression (\ref{3.12}) into (\ref{3.15}), (\ref{3.18})
and the $o$-transform of (\ref{3.18}), respectively:

\begin{lemma}\label{lemma5.2}
We have
\begin{eqnarray}
T_{-z}^{+}\left(z^{-1}\delta\left(\frac{x_{1}-x_{0}}{z}\right)Y_{t}(v,
x_{0})\right)
=x_{0}^{-1}\delta\left(\frac{z-x_{1}}{-x_{0}}\right) Y_{t}(v, x_{1}),
\label{3.71}\\
T_{-z}^{-}\left(z^{-1}\delta\left(\frac{x_{1}-x_{0}}{z}\right)Y_{t}(v,
x_{0})\right)
=x_{0}^{-1}\delta\left(\frac{x_{1}-z}{x_{0}}\right) Y_{t}(v, x_{1}),
\label{3.72}\\
T_{-z}^{o}\left(z^{-1}\delta\left(\frac{x_{1}-x_{0}}{z}\right)Y_{t}(v,
x_{0})\right)
=x_{0}^{-1}\delta\left(\frac{x_{1}-z}{x_{0}}\right) Y^{o}_{t}(v, x_{1}).
\label{3.73}
\end{eqnarray}
\end{lemma}
\pf
We  prove (\ref{3.71}): {From} (\ref{3.12}), the coefficient of 
$x_{0}^{-n-1}x_{1}^{-m-1}$ in the left-hand side of
(\ref{3.71}) is $T_{-z}^{+}(v\otimes (z+t)^{m}t^{n})$. By the definitions,
\begin{eqnarray}
T^{+}_{-z}(v\otimes (z+t)^{m}t^{n})
=v\otimes t^{m}(-(z-t))^{n}.
\end{eqnarray}
On the other hand, the right-hand side of (\ref{3.71}) can be written as
\begin{eqnarray}\label{3.75}
v\otimes
x_{0}^{-1}\delta\left(\frac{z-x_{1}}{-x_{0}}\right)  x_{1}^{-1}
\delta\left(\frac{t}{x_{1}}\right)
=v\otimes x_{0}^{-1}\delta\left(\frac{z-t}{-x_{0}}\right) x_{1}^{-1}
\delta\left(\frac{t}{x_{1}}\right),
\end{eqnarray}
where we have used (\ref{3.5}) and the fundamental property 
(\ref{Xx1x2=Xx2x2}) of the
formal $\delta$-function. The coefficient of
$x_{0}^{-n-1}x_{1}^{-m-1}$ in the right-hand side of (\ref{3.75}) is also
$v\otimes t^{m}(-(z-t))^{n}$,  proving (\ref{3.71}).  Formula
(\ref{3.72}) is proved similarly, and (\ref{3.73}) is obtained 
{from} (\ref{3.72}) by the
application of the map $o$.
\epf

\subsection{Constructions of $P(z)$-tensor products}

We  proceed to the construction of $P(z)$-tensor products.  While
one can certainly consider categories ${\cal C}$ in Remark
\ref{bifunctor} that are not closed under the contragredient functor,
it is most natural to consider such categories ${\cal C}$ that are
indeed closed under this functor (recall Notation \ref{MGM}).  Our
constructions of $P(z)$-tensor products will in fact use the
contragredient functor; the $P(z)$-tensor product of (generalized)
modules $W_1$ and $W_2$ will arise as the contragredient module of a
certain subspace of the vector space dual $(W_1 \otimes W_2)^*$.  We
now present this ``double-dual'' approach to the construction of
$P(z)$-tensor products, generalizing the double-dual approach carried
out in \cite{tensor1}--\cite{tensor3}.  At first, we need not fix any
subcategory ${\cal C}$ of ${\cal M}_{sg}$ or ${\cal GM}_{sg}$.
As usual, we take $z\in \C^{\times}$.

We shall be constructing an action of the space $V\otimes \mathbb{C}[t,
t^{-1}, (z^{-1}-t)^{-1}]$ on the space $(W_1 \otimes W_2)^*$, given
generalized $V$-modules $W_1$ and $W_2$. This action will be based
on the translation operations and on the $o$-operation discussed in the 
preceding subsection.  More
precisely, it is the space $V\otimes \iota_{+} \mathbb{C}[t, t^{-1},
(z^{-1}-t)^{-1}]$ whose action we shall define.

Let $I$ be a $P(z)$-intertwining map of type ${W_3\choose W_1\, W_2}$,
as in Definition \ref{im:imdef}. Consider the contragredient generalized 
$V$-module $(W_{3}', Y_{3}')$, recall the opposite vertex operator 
(\ref{yo}) and formula (\ref{y'}), and recall why the ingredients of 
formula (\ref{im:def})  are well defined. For 
$v\in V$, $w_{(1)}\in W_{1}$, $w_{(2)}\in W_{2}$ and $w'_{(3)}\in W'_3$, 
applying $w'_{(3)}$ to (\ref{im:def}), replacing $x_1$ by
$x_1^{-1}$ in the resulting formula 
and then replacing $v$ by $e^{x_1L(1)}(-x_1^{-2})^{L(0)}v$,
we get:
\begin{eqnarray}\label{im:def'}
\lefteqn{\left\langle x_0^{-1}\delta\bigg(\frac{x^{-1}_1-z}{x_0}\bigg)
Y_{3}'(v, x_1)w'_{(3)}, I(w_{(1)}\otimes w_{(2)})\right\rangle}\nno\\
&&=\left\langle w'_{(3)},z^{-1}\delta\bigg(\frac{x^{-1}_1-x_0}{z}\bigg)
I(Y_1(e^{x_1L(1)}(-x_1^{-2})^{L(0)}v, x_0)w_{(1)}\otimes
w_{(2)})\right\rangle\nno\\
&&\quad +\left\langle w'_{(3)}, x^{-1}_0\delta\bigg(\frac{z-x^{-1}_1}{-x_0}\bigg)
I(w_{(1)}\otimes Y_2^{o}(v, x_1)w_{(2)})\right\rangle.
\end{eqnarray}
We shall use this to motivate our action.

As we discussed in the preceding subsection (see (\ref{3.18-1}) and
(\ref{3.19-1})), in the left-hand side of (\ref{im:def'}), the
coefficients of
\begin{equation}\label{deltaY3'}
x_0^{-1}\delta\bigg(\frac{x^{-1}_1-z}{x_0}\bigg) Y_{3}'(v, x_1)
\end{equation}
in powers of $x_0$ and $x_1$, for all $v\in V$, span
\begin{equation}\label{tausubW3'}
\tau_{W'_3}(V\otimes \iota_{+}{\mathbb C}[t,t^{-1},(z^{-1}-t)^{-1}])
\end{equation}
(recall (\ref{tauw}) and (\ref{3.7})).
Let us now define an action of $V\otimes \iota_{+}{\mathbb
C}[t,t^{-1},(z^{-1}-t)^{-1}]$ on $(W_1\otimes W_2)^*$. 

\begin{defi}\label{deftau}{\rm
Define the linear action $\tau_{P(z)}$ of
\[
V \otimes \iota_{+}{\mathbb C}[t,t^{- 1}, (z^{-1}-t)^{-1}]
\]
on $(W_1 \otimes W_2)^*$ by 
\begin{eqnarray}\label{taudef0}
(\tau_{P(z)}(\xi)\lambda)(w_{(1)}\otimes w_{(2)})&=&\lambda
(\tau_{W_{1}}((\iota_+\circ T_z\circ\iota_-^{-1}\circ o)\xi)w_{(1)}
\otimes w_{(2)})\nno\\
&&+\lambda (w_{(1)}\otimes\tau_{W_{2}}((\iota_+\circ
\iota_-^{-1}\circ o)\xi)w_{(2)})
\end{eqnarray}
for $\xi\in V \otimes \iota_{+}{\mathbb C}[t,t^{- 1},
(z^{-1}-t)^{-1}]$, $\lambda\in (W_1\otimes W_2)^*$, $w_{(1)}\in W_1$
and $w_{(2)}\in W_2$. (The fact that the right-hand side is well
defined follows immediately {}from the generating-function reformulation
of (\ref{taudef0}) given in (\ref{taudef}) below.) Denote by
$Y'_{P(z)}$ the action of $V\otimes{\mathbb C}[t,t^{-1}]$ on
$(W_1\otimes W_2)^*$ thus defined, that is,
\begin{equation}\label{y'-p-z}
Y'_{P(z)}(v,x)=\tau_{P(z)}(Y_t(v,x))
\end{equation}
for $v\in V\otimes \C[t, t^{-1}]$.
}
\end{defi}

By Lemma \ref{tauP}, (\ref{3.7}) and (\ref{tauw-yto}), we see that
formula (\ref{taudef0}) can be written in terms of
generating functions as
\begin{eqnarray}\label{taudef}
\lefteqn{\bigg(\tau_{P(z)}
\bigg(x_0^{-1}\delta\bigg(\frac{x^{-1}_1-z}{x_0}\bigg)
Y_{t}(v, x_1)\bigg)\lambda\bigg)(w_{(1)}\otimes w_{(2)})}\nno\\
&&=z^{-1}\delta\bigg(\frac{x^{-1}_1-x_0}{z}\bigg)
\lambda(Y_1(e^{x_1L(1)}(-x_1^{-2})^{L(0)}v, x_0)w_{(1)}\otimes w_{(2)})
\nno\\
&&\quad +x^{-1}_0\delta\bigg(\frac{z-x^{-1}_1}{-x_0}\bigg)
\lambda(w_{(1)}\otimes Y_2^{o}(v, x_1)w_{(2)})
\end{eqnarray}
for $v\in V$, $\lambda\in (W_1\otimes W_2)^{*}$, $w_{(1)}\in W_1$,
$w_{(2)}\in W_2$; note that by (\ref{3.18-1})--(\ref{3.19-1}), the
expansion coefficients in $x_{0}$ and $x_{1}$ of the left-hand side
span the space of elements in the left-hand side of (\ref{taudef0}). Compare 
formula (\ref{taudef}) with the motivating formula (\ref{im:def'}).
  The generating function form of the action
$Y'_{P(z)}$ can be obtained by taking $\res_{x_0}$ of both sides of
(\ref{taudef}), that is,
\begin{eqnarray}\label{Y'def}
\lefteqn{(Y'_{P(z)}(v,x_1)\lambda)(w_{(1)}\otimes w_{(2)})=
\lambda(w_{(1)}\otimes Y_2^{o}(v, x_1)w_{(2)})}\nno\\
&&\quad+\res_{x_0}z^{-1}\delta\bigg(\frac{x^{-1}_1-x_0}{z}\bigg)
\lambda(Y_1(e^{x_1L(1)}(-x_1^{-2})^{L(0)}v, x_0)w_{(1)}\otimes
w_{(2)}).
\end{eqnarray}

\begin{rema}\label{I-intw}{\rm
Using the actions $\tau_{W'_3}$ and $\tau_{P(z)}$, we can write
(\ref{im:def'}) as
\[
\left(x_0^{-1}\delta\left(\frac{x^{-1}_1-z}{x_0}\right)
Y_{3}'(v, x_1)w'_{(3)}\right)\circ I=
\tau_{P(z)}\left(x_0^{-1}\delta\left(\frac{x^{-1}_1-z}{x_0}\right)
Y_t(v, x_1)\right)(w'_{(3)}\circ I)
\]
or equivalently, as
\[
\left(\tau_{W'_3}\left(x_0^{-1}\delta\left(\frac{x^{-1}_1-z}{x_0}\right)
Y_t(v, x_1)\right)w'_{(3)}\right)\circ I=
\tau_{P(z)}\left(x_0^{-1}\delta\left(\frac{x^{-1}_1-z}{x_0}\right)
Y_t(v, x_1)\right)(w'_{(3)}\circ I).
\]
}
\end{rema}

In the spirit of the discussion related to Lemma \ref{4.36}, we find
it natural to introduce subspaces of $(W_{1}\otimes W_{2})^{*}$
homogeneous with respect to $\tilde{A}$.  Since $W_{1}$ and $W_{2}$
are $\tilde{A}$-graded, $W_{1}\otimes W_{2}$ also has a natural
$\tilde{A}$-grading---the tensor product grading, and we shall write
$(W_{1}\otimes W_{2})^{(\beta)}$ for the homogeneous subspace of
degree $\beta\in \tilde{A}$ of $W_{1}\otimes W_{2}$. For $\beta\in
\tilde{A}$, let $((W_{1}\otimes W_{2})^{*})^{(\beta)}$ be the space
consisting of the elements $\lambda\in (W_{1}\otimes W_{2})^{*}$ such
that $\lambda(\tilde{w})=0$ for $\tilde{w}\in (W_{1}\otimes
W_{2})^{(\gamma)}$ with $\gamma\ne -\beta$. (Of course, the full space
$(W_{1}\otimes W_{2})^{*}$ is not $\tilde{A}$-graded since it is not a
direct sum of subspaces homogeneous with respect to $\tilde{A}$.)

The space 
$V \otimes \iota_{+}{\mathbb C}[t,t^{- 1}, (z^{-1}-t)^{-1}]$
also has an $A$-grading, induced {}from the $A$-grading on $V$:
For $\alpha\in A$, 
\begin{equation}
(V\otimes 
\iota_{+}{\mathbb C}[t,t^{- 1}, (z^{-1}-t)^{-1}])^{(\alpha)}
=V^{(\alpha)}\otimes 
\iota_{+}{\mathbb C}[t,t^{- 1}, (z^{-1}-t)^{-1}].
\end{equation}

Using these gradings, we formulate:

\begin{defi}\label{linearactioncompatible}
{\rm We call a linear action $\tau$ of 
$V \otimes \iota_{+}{\mathbb C}[t,t^{- 1}, (z^{-1}-t)^{-1}]$
on $(W_1 \otimes W_2)^*$
{\it $\tilde{A}$-compatible} if 
for $\alpha\in A$, $\beta\in \tilde{A}$,
$\xi\in 
(V \otimes \iota_{+}{\mathbb C}[t,t^{- 1}, (z^{-1}-t)^{-1}])^{(\alpha)}$
and $\lambda\in ((W_{1}\otimes W_{2})^{*})^{(\beta)}$,
\[
\tau(\xi)\lambda\in ((W_{1}\otimes W_{2})^{*})^{(\alpha+\beta)}.
\]}
\end{defi}

{}From (\ref{taudef0}) or (\ref{taudef}), we have:

\begin{propo}\label{tau-a-comp}
The action $\tau_{P(z)}$ is $\tilde{A}$-compatible. \epf
\end{propo}

\begin{rema}
{\rm Notice that Proposition \ref{tau-a-comp} is analogous to the
condition (\ref{m-v_l-A}) in the definition of the notion of
(generalized) module.  We now proceed to establish several more of the
module-action properties for our action $\tau_{P(z)}$ on $(W_1 \otimes
W_2)^*$, in both the conformal and M\"obius cases.  However, while we
will prove the commutator formula for our action (see Proposition
\ref{pz-comm} below), we will {\em not} be able to prove the Jacobi
identity on an element $\lambda \in (W_1 \otimes W_2)^*$ until we
assume the ``$P(z)$-compatibility condition'' for the element
$\lambda$ (see Theorem \ref{comp=>jcb} below).  We shall be
constructing a certain subspace $W_1\hboxtr_{P(z)} W_2$ of $(W_1
\otimes W_2)^*$ which under suitable conditions will be a generalized
$V$-module and whose contragredient module will be
$W_1\boxtimes_{P(z)} W_2$ (see Remark \ref{motivationofbackslash} and
Proposition \ref{tensor1-13.7}), and we shall use the
$P(z)$-compatibility condition to describe this subspace (see Theorem
\ref{characterizationofbackslash}).  }
\end{rema}

We have the following result generalizing Proposition 13.3 in
\cite{tensor3}:

\begin{propo}\label{id-dev}
The action $Y'_{P(z)}$ has the property
\[
Y'_{P(z)}({\bf 1},x)=1,
\]
where $1$ on the right-hand side is the identity map of $(W_1\otimes
W_2)^*$. It also has the $L(-1)$-derivative property
\[
\frac{d}{dx}Y'_{P(z)}(v,x)=Y'_{P(z)}(L(-1)v,x)
\]
for $v\in V$.
\end{propo}
\pf The first statement follows directly {}from the definition. We
prove the $L(-1)$-derivative property. {}From (\ref{Y'def}), we
obtain, using (\ref{yo-l-1}),
\begin{eqnarray}\label{der-1}
\lefteqn{\left(\frac{d}{dx}Y'_{P(z)}(v, x)\lambda\right)
(w_{(1)}\otimes w_{(2)})}\nno\\
&&=\frac{d}{dx}\lambda(w_{(1)}\otimes Y_{2}^{o}(v, x)w_{(2)})\nn
&&\quad +\frac{d}{dx}\res_{x_{0}}z^{-1}\delta\left(\frac{x^{-1}-x_{0}}{z}\right)
\lambda(Y_{1}(e^{xL(1)}(-x^{-2})^{L(0)}v, x_{0})w_{(1)}\otimes w_{(2)})\nno\\
&&=\lambda\left(w_{(1)}\otimes \frac{d}{dx}Y_{2}^{o}(v, x)w_{(2)}\right)\nn
&&\quad+ \res_{x_{0}}\left(\frac{d}{dx}\left(z^{-1}
\delta\left(\frac{x^{-1}-x_{0}}{z}\right)\right)\right)
\lambda(Y_{1}(e^{xL(1)}(-x^{-2})^{L(0)}v, x_{0})w_{(1)}\otimes w_{(2)})\nno\\
&&\quad +\res_{x_{0}}z^{-1}\delta\left(\frac{x^{-1}-x_{0}}{z}\right)
\frac{d}{dx}\lambda(Y_{1}(e^{xL(1)}(-x^{-2})^{L(0)}v, x_{0})w_{(1)}
\otimes w_{(2)})
\nonumber\\
&&=\lambda(w_{(1)}\otimes Y_{2}^{o}(L(-1)v, x)w_{(2)})\nn
&&\quad+ \res_{x_{0}}\left(\frac{d}{dx}\left(z^{-1}
\delta\left(\frac{x^{-1}-x_{0}}{z}\right)\right)\right)
\lambda(Y_{1}(e^{xL(1)}(-x^{-2})^{L(0)}v, x_{0})w_{(1)}\otimes w_{(2)})\nno\\
&&\quad +\res_{x_{0}}z^{-1}\delta\left(\frac{x^{-1}-x_{0}}{z}\right)
\lambda(Y_{1}(e^{xL(1)}L(1)(-x^{-2})^{L(0)}v, x_{0})w_{(1)}\otimes w_{(2)})
\nonumber\\
&&\quad -2\res_{x_{0}}z^{-1}\delta\left(\frac{x^{-1}-x_{0}}{z}\right)
\lambda(Y_{1}(e^{xL(1)}L(0)x^{-1}(-x^{-2})^{L(0)}v, x_{0})w_{(1)}
\otimes w_{(2)}).
\end{eqnarray}
The second term
on the right-hand side of (\ref{der-1}) is equal to
\begin{eqnarray}\label{der-2}
\lefteqn{-\res_{x_{0}}x^{-2}\left(\frac{d}{dx^{-1}}\left(z^{-1}
\delta\left(\frac{x^{-1}-x_{0}}{z}\right)\right)\right)\cdot}\nno\\
&&\quad \quad\quad\quad \cdot
\lambda(Y_{1}(e^{xL(1)}(-x^{-2})^{L(0)}v, x_{0})w_{(1)}\otimes w_{(2)})\nno\\
&&=\res_{x_{0}}x^{-2}\left(\frac{d}{dx_{0}}\left(z^{-1}
\delta\left(\frac{x^{-1}-x_{0}}{z}\right)\right)\right)\cdot\nno\\
&&\quad \quad\quad\quad \cdot
\lambda(Y_{1}(e^{xL(1)}(-x^{-2})^{L(0)}v, x_{0})w_{(1)}\otimes w_{(2)})\nno\\
&&=-\res_{x_{0}}x^{-2}z^{-1}
\delta\left(\frac{x^{-1}-x_{0}}{z}\right)\cdot\nno\\
&&\quad \quad\quad\quad \cdot
\frac{d}{dx_{0}}\lambda(Y_{1}(e^{xL(1)}(-x^{-2})^{L(0)}v, x_{0})w_{(1)}
\otimes w_{(2)})\nno\\
&&=-\res_{x_{0}}x^{-2}z^{-1}
\delta\left(\frac{x^{-1}-x_{0}}{z}\right)\cdot\nno\\
&&\quad \quad\quad\quad \cdot
\lambda(Y_{1}(L(-1)e^{xL(1)}(-x^{-2})^{L(0)}v, x_{0})w_{(1)}
\otimes w_{(2)}).
\end{eqnarray}
By (\ref{der-2}), (\ref{log:SL2-3}) and (\ref{log:xLx^}),
the right-hand side of (\ref{der-1}) is equal to
\begin{eqnarray*}
\lefteqn{\displaystyle \lambda(w_{(1)}\otimes 
Y_{2}^{o}(L(-1)v, x)w_{(2)})}\nno\\
&&\quad + \res_{x_{0}}z^{-1}\delta\left(\frac{x^{-1}-x_{0}}{z}\right)
\lambda(Y_{1}(e^{xL(1)}(-x^{-2})^{L(0)}L(-1)v, x_{0})w_{(1)}\otimes w_{(2)})
\nonumber\\
&& =(Y'_{P(z)}(L(-1)v, x)\lambda)(w_{(1)}\otimes w_{(2)}),
\end{eqnarray*}
proving the $L(-1)$-derivative property.
\epfv

\begin{propo}\label{pz-comm}
The action $Y'_{P(z)}$ satisfies the commutator formula for vertex
operators, that is, on $(W_1\otimes W_2)^*$,
\begin{eqnarray*}
\lefteqn{[Y'_{P(z)}(v_1,x_1),Y'_{P(z)}(v_2,x_2)]}\\
&&=\res_{x_0}x_2^{-1}\delta\bigg(\frac{x_1-x_0}{x_2}\bigg)
Y'_{P(z)}(Y(v_1,x_0)v_2,x_2)
\end{eqnarray*}
for $v_1,v_2\in V$.
\end{propo}
\pf
In the following proof, the reader should note the
well-definedness of each expression and the justifiability of each use of
a $\delta$-function property.

Let $\lambda
\in (W_{1}\otimes W_{2})^{*}$, $v_{1}, v_{2}\in V$, $w_{(1)}\in W_{1}$
and $w_{(2)}\in W_{2}$. By (\ref{Y'def}),
\bea\label{y-12}
\lefteqn{(Y'_{P(z)}(v_{1}, x_{1})Y'_{P(z)}(v_{2}, x_{2})
\lambda)(w_{(1)}\otimes w_{(2)})}\nno\\
&&=(Y'_{P(z)}(v_{2}, x_{2})\lambda)(w_{(1)}\otimes Y_2^{o}(v_{1}, x_1)w_{(2)})\nno\\
&&\quad+\res_{y_{1}}z^{-1}\delta\bigg(\frac{x^{-1}_1-y_{1}}{z}\bigg)
(Y'_{P(z)}(v_{2}, x_{2})\lambda)(Y_1(e^{x_1L(1)}(-x_1^{-2})^{L(0)}v_{1}, y_{1})w_{(1)}
\otimes w_{(2)})\nn
&&=\lambda(w_{(1)}\otimes Y_2^{o}(v_{2}, x_2)Y_2^{o}(v_{1}, x_1)w_{(2)})\nno\\
&&\quad+\res_{y_{2}}z^{-1}\delta\bigg(\frac{x^{-1}_{2}-y_{2}}{z}\bigg)
\lambda(Y_1(e^{x_{2}L(1)}(-x_{2}^{-2})^{L(0)}v_{2}, y_{2})w_{(1)}\otimes
Y_2^{o}(v_{1}, x_1)w_{(2)})\nn
&&\quad+\res_{y_{1}}z^{-1}\delta\bigg(\frac{x^{-1}_1-y_{1}}{z}\bigg)
\lambda(Y_1(e^{x_1L(1)}(-x_1^{-2})^{L(0)}v_{1}, y_{1})w_{(1)}\otimes
Y_2^{o}(v_{2}, x_2)w_{(2)})\nn
&&\quad+\res_{y_{1}}\res_{y_{2}}z^{-1}\delta\bigg(\frac{x^{-1}_{1}-y_{1}}{z}\bigg)
z^{-1}\delta\bigg(\frac{x^{-1}_2-y_{2}}{z}\bigg)\cdot\nn
&&\quad\quad\quad\quad\cdot\lambda(Y_1(e^{x_2L(1)}(-x_2^{-2})^{L(0)}v_{2}, y_{2})
Y_1(e^{x_1 L(1)}(-x_1^{-2})^{L(0)}v_{1}, y_{1})w_{(1)}\otimes
w_{(2)}).
\eea
Transposing the subscripts $1$ and $2$ of the symbols $v$, $x$ and $y$,
we also have
\bea\label{y-21}
\lefteqn{(Y'_{P(z)}(v_{2}, x_{2})Y'_{P(z)}(v_{1}, x_{1})
\lambda)(w_{(1)}\otimes w_{(2)})}\nno\\
&&=\lambda(w_{(1)}\otimes Y_2^{o}(v_{1}, x_1)Y_2^{o}(v_{2}, x_2)w_{(2)})\nno\\
&&\quad+\res_{y_{1}}z^{-1}\delta\bigg(\frac{x^{-1}_{1}-y_{1}}{z}\bigg)
\lambda(Y_1(e^{x_{1}L(1)}(-x_{1}^{-2})^{L(0)}v_{1}, y_{1})w_{(1)}\otimes
Y_2^{o}(v_{2}, x_2)w_{(2)})\nn
&&\quad+\res_{y_{2}}z^{-1}\delta\bigg(\frac{x^{-1}_2-y_{2}}{z}\bigg)
\lambda(Y_1(e^{x_2L(1)}(-x_2^{-2})^{L(0)}v_{2}, y_{2})w_{(1)}\otimes
Y_2^{o}(v_{1}, x_1)w_{(2)})\nn
&&\quad+\res_{y_{2}}\res_{y_{1}}z^{-1}\delta\bigg(\frac{x^{-1}_{2}-y_{2}}{z}\bigg)
z^{-1}\delta\bigg(\frac{x^{-1}_1-y_{1}}{z}\bigg)\cdot\nn
&&\quad\quad\quad\quad\cdot\lambda(Y_1(e^{x_1L(1)}(-x_1^{-2})^{L(0)}v_{1}, y_{1})
Y_1(e^{x_2 L(1)}(-x_2^{-2})^{L(0)}v_{2}, y_{2})w_{(1)}\otimes
w_{(2)}).
\eea
The equalities (\ref{y-12}) and (\ref{y-21}) give
\bea\label{y-bracket}
\lefteqn{([Y'_{P(z)}(v_{1}, x_{1}), Y'_{P(z)}(v_{2}, x_{2})]
\lambda)(w_{(1)}\otimes w_{(2)})}\nno\\
&&=\lambda(w_{(1)}\otimes [Y_2^{o}(v_{2}, x_2), Y_2^{o}(v_{1}, x_1)]w_{(2)})\nn
&&\quad-\res_{y_{1}}\res_{y_{2}}z^{-1}\delta\bigg(\frac{x^{-1}_{1}-y_{1}}{z}\bigg)
z^{-1}\delta\bigg(\frac{x^{-1}_2-y_{2}}{z}\bigg)\cdot\nn
&&\quad\quad\quad\quad\cdot\lambda([Y_1(e^{x_1 L(1)}(-x_1^{-2})^{L(0)}v_{1}, y_{1}),
Y_1(e^{x_2L(1)}(-x_2^{-2})^{L(0)}v_{2}, y_{2})]w_{(1)}\otimes
w_{(2)})\nn
&&=\res_{x_{0}}x_{2}^{-1}
\delta\left(\frac{x_{1}-x_{0}}{x_{2}}\right)\lambda(w_{(1)}\otimes
Y_2^{o}(Y(v_{1}, x_0)v_{2}, x_2)w_{(2)})\nn
&&\quad-\res_{y_{1}}\res_{y_{2}}\res_{x_{0}}
z^{-1}\delta\bigg(\frac{x^{-1}_{1}-y_{1}}{z}\bigg)
z^{-1}\delta\bigg(\frac{x^{-1}_2-y_{2}}{z}\bigg)
y_{2}^{-1}\delta\left(\frac{y_{1}-x_{0}}{y_{2}}\right)\cdot\nn
&&\quad\quad\quad\quad\cdot\lambda(
Y_{1}(Y(e^{x_1 L(1)}(-x_1^{-2})^{L(0)}v_{1}, x_{0})
e^{x_2L(1)}(-x_2^{-2})^{L(0)}v_{2}, y_{2})w_{(1)}\otimes
w_{(2)})
\eea
(recall (\ref{op-jac-id})).

But we have
\bea\label{delta-idty}
\lefteqn{z^{-1}\delta\bigg(\frac{x^{-1}_{1}-y_{1}}{z}\bigg)
z^{-1}\delta\bigg(\frac{x^{-1}_2-y_{2}}{z}\bigg)
y_{2}^{-1}\delta\left(\frac{y_{1}-x_{0}}{y_{2}}\right)}\nno\\
&&=\left({\displaystyle \sum_{m,n\in {\mathbb Z}}}
\frac{(x_{1}^{-1}-y_{1})^{m}}{z^{m+1}}
\frac{(x_{2}^{-1}-y_{2})^{n}}{z^{n+1}}\right)
y_{2}^{-1}\delta\left(\frac{y_{1}-x_{0}}{y_{2}}\right)\nno\\
&&=\left({\displaystyle \sum_{m,n\in {\mathbb Z}}}(x_{2}^{-1}-y_{2})^{-1}
\left(\frac{x_{1}^{-1}-y_{1}}{x_{2}^{-1}-y_{2}}\right)^{m}
\frac{(x_{2}^{-1}-y_{2})^{m+n+1}}
{z^{m+n+2}} \right)
y_{2}^{-1}\delta\left(\frac{y_{1}-x_{0}}{y_{2}}\right)\nno\\
&&=\left({\displaystyle \sum_{m,k\in {\mathbb Z}}}(x_{2}^{-1}-y_{2})^{-1}
\left(\frac{x_{1}^{-1}-y_{1}}{x_{2}^{-1}-y_{2}}\right)^{m}
z^{-1}\left(\frac{x_{2}^{-1}-y_{2}}{z}\right)^{k}\right)
y_{2}^{-1}\delta\left(\frac{y_{1}-x_{0}}{y_{2}}\right)\nno\\
&&=(x_{2}^{-1}-y_{2})^{-1}\delta\left(\frac{x_{1}^{-1}-y_{1}}{x_{2}^{-1}-y_{2}}\right)
z^{-1}\delta\left(\frac{x_{2}^{-1}-y_{2}}{z}\right)
y_{2}^{-1}\delta\left(\frac{y_{1}-x_{0}}
{y_{2}}\right)\nno\\
&&=x_{2}\delta\left(\frac{x_{1}^{-1}-(y_{1}-y_{2})}{x_{2}^{-1}}\right)
z^{-1}\delta\left(\frac{x_{2}^{-1}-y_{2}}{z}\right)
y_{1}^{-1}\delta\left(\frac{y_{2}+x_{0}}
{y_{1}}\right)\nno\\
&&=x_{2}\delta\left(\frac{x_{1}^{-1}-x_{0}}{x_{2}^{-1}}\right)
z^{-1}\delta\left(\frac{x_{2}^{-1}-y_{2}}{z}\right)y_{1}^{-1}\delta
\left(\frac{y_{2}+x_{0}}{y_{1}}\right).
\eea
By (\ref{log:p2}), (\ref{log:p3}) and (\ref{xe^Lx}), we also have
\bea\label{sl2-idty}
\lefteqn{Y(e^{x_1 L(1)}(-x_1^{-2})^{L(0)}v_{1}, x_{0})
e^{x_2L(1)}(-x_2^{-2})^{L(0)}}\nn
&&=e^{x_2L(1)}Y\left(e^{-x_{2}(1+x_{0}x_{2})L(1)}(1+x_{0}x_{2})^{-2L(0)}
e^{x_1 L(1)}(-x_1^{-2})^{L(0)}v_{1}, \frac{x_{0}}{1+x_{0}x_{2}}\right)
(-x_2^{-2})^{L(0)}\nn
&&=e^{x_2L(1)}(-x_2^{-2})^{L(0)}
Y\Biggl((-x_2^{2})^{L(0)}e^{-x_{2}(1+x_{0}x_{2})L(1)}\cdot\nn
&&\quad\quad\quad\quad\quad\quad\quad\quad\quad\quad\quad\cdot (1+x_{0}x_{2})^{-2L(0)}
e^{x_1 L(1)}(-x_1^{-2})^{L(0)}v_{1}, -\frac{x_{0}x_{2}^{2}}{1+x_{0}x_{2}}\Biggr)
\nn
&&=e^{x_2L(1)}(-x_2^{-2})^{L(0)}
Y\Biggl(e^{-x_{2}(1+x_{0}x_{2})(-x_2^{-2})L(1)}
(-x_2^{2})^{L(0)}\cdot\nn
&&\quad\quad\quad\quad\quad\quad\quad\quad\quad\quad\quad\cdot(1+x_{0}x_{2})^{-2L(0)}
e^{x_1 L(1)}(-x_1^{-2})^{L(0)}v_{1}, -\frac{x_{0}x_{2}^{2}}{1+x_{0}x_{2}}\Biggr)
\nn
&&=e^{x_2L(1)}(-x_2^{-2})^{L(0)}
Y\Biggl(e^{(x_{2}^{-1}+x_{0})L(1)}
(-(x_{2}^{-1}+x_{0})^{2})^{-L(0)}
e^{x_1 L(1)}(-x_1^{-2})^{L(0)}v_{1}, -\frac{x_{0}x_{2}}{x_{2}^{-1}+x_{0}}\Biggr)
\nn
&&=e^{x_2L(1)}(-x_2^{-2})^{L(0)}
Y\Biggl(e^{(x_{2}^{-1}+x_{0})L(1)}e^{-x_1 (x_{2}^{-1}+x_{0})^{2})L(1)}
\cdot\nn
&&\quad\quad\quad\quad\quad\quad\quad\quad\quad\quad\quad\cdot
(-(x_{2}^{-1}+x_{0})^{2})^{-L(0)}(-x_1^{-2})^{L(0)}v_{1},
-\frac{x_{0}x_{2}}{x_{2}^{-1}+x_{0}}\Biggr)
\nn
&&=e^{x_2L(1)}(-x_2^{-2})^{L(0)}
Y\left(e^{(x_{2}^{-1}+x_{0})L(1)}
e^{-x_1 (x_{2}^{-1}+x_{0})^{2}L(1)}
((x_{2}^{-1}+x_{0})x_{1})^{-2L(0)}v_{1},
-\frac{x_{0}x_{2}}{x_{2}^{-1}+x_{0}}\right).\nn
&&
\eea

Using (\ref{delta-idty}), (\ref{sl2-idty}) and the basic
properties of the formal delta function, we see that
(\ref{y-bracket}) becomes
\bea
\lefteqn{([Y'_{P(z)}(v_{1}, x_{1}), Y'_{P(z)}(v_{2}, x_{2})]\lambda)(w_{(1)}
\otimes w_{(2)})}\nno\\
&&=\res_{x_{0}}x_{2}^{-1}
\delta\left(\frac{x_{1}-x_{0}}{x_{2}}\right)\lambda(w_{(1)}\otimes
Y_2^{o}(Y(v_{1}, x_0)v_{2}, x_2)w_{(2)})\nn
&&\quad-\res_{y_{1}}\res_{y_{2}}\res_{x_{0}}
x_{2}\delta\left(\frac{x_{1}^{-1}-x_{0}}{x_{2}^{-1}}\right)
z^{-1}\delta\left(\frac{x_{2}^{-1}-y_{2}}{z}\right)y_{1}^{-1}\delta
\left(\frac{y_{2}+x_{0}}{y_{1}}\right)\cdot\nn
&&\quad\quad\quad\quad\cdot\lambda\Biggl(
Y_{1}\Biggl(e^{x_2L(1)}(-x_2^{-2})^{L(0)}
Y\Biggl(e^{(x_{2}^{-1}+x_{0})L(1)}
e^{-x_1 (x_{2}^{-1}+x_{0})^{2}L(1)}\cdot\nn
&&\quad\quad\quad\quad\quad\quad\quad\quad\quad\quad\quad\cdot
((x_{2}^{-1}+x_{0})x_{1})^{-2L(0)}v_{1},
-\frac{x_{0}x_{2}}{x_{2}^{-1}+x_{0}}\Biggr)v_{2}, y_{2}\Biggr)w_{(1)}\otimes
w_{(2)}\Biggr)\nn
&&=\res_{x_{0}}x_{2}^{-1}
\delta\left(\frac{x_{1}-x_{0}}{x_{2}}\right)\lambda(w_{(1)}\otimes
Y_2^{o}(Y(v_{1}, x_0)v_{2}, x_2)w_{(2)})\nn
&&\quad-\res_{x_{0}}\res_{y_{2}}
x_{2}\delta\left(\frac{x_{1}^{-1}-x_{0}}{x_{2}^{-1}}\right)
z^{-1}\delta\left(\frac{x_{2}^{-1}-y_{2}}{z}\right)\cdot\nn
&&\quad\quad\quad\quad\cdot\lambda(
Y_{1}(e^{x_2L(1)}(-x_2^{-2})^{L(0)}
Y(e^{x_{1}^{-1}L(1)}
e^{-x_1^{-1}L(1)}\cdot\nn
&&\quad\quad\quad\quad\quad\quad\quad\quad\quad\quad\quad\cdot
(x_{1}^{-1}x_{1})^{-2L(0)}v_{1},
-x_{0}x_{1}x_{2})v_{2}, y_{2})w_{(1)}\otimes
w_{(2)})\nn
&&=\res_{x_{0}}x_{2}^{-1}
\delta\left(\frac{x_{1}-x_{0}}{x_{2}}\right)\lambda(w_{(1)}\otimes
Y_2^{o}(Y(v_{1}, x_0)v_{2}, x_2)w_{(2)})\nn
&&\quad-\res_{x_{0}}x_{2}\delta\left(\frac{x_{2}+(-x_{0}x_{1}x_{2})}{x_{1}}\right)
\res_{y_{2}}
z^{-1}\delta\left(\frac{x_{2}^{-1}-y_{2}}{z}\right)\cdot\nn
&&\quad\quad\quad\quad\cdot\lambda(
Y_{1}(e^{x_2L(1)}(-x_2^{-2})^{L(0)}
Y(v_{1},
-x_{0}x_{1}x_{2})v_{2}, y_{2})w_{(1)}\otimes
w_{(2)})\nn
&&=\res_{x_{0}}x_{2}^{-1}
\delta\left(\frac{x_{1}-x_{0}}{x_{2}}\right)\lambda(w_{(1)}\otimes
Y_2^{o}(Y(v_{1}, x_0)v_{2}, x_2)w_{(2)})\nn
&&\quad+\res_{y_{0}}x_{1}^{-1}\delta\left(\frac{x_{2}+y_{0}}{x_{1}}\right)
\res_{y_{2}}
z^{-1}\delta\left(\frac{x_{2}^{-1}-y_{2}}{z}\right)\cdot\nn
&&\quad\quad\quad\quad\cdot\lambda(
Y_{1}(e^{x_2L(1)}(-x_2^{-2})^{L(0)}
Y(v_{1},
y_{0})v_{2}, y_{2})w_{(1)}\otimes
w_{(2)})\nn
&&=\res_{x_{0}}x_{2}^{-1}
\delta\left(\frac{x_{1}-x_{0}}{x_{2}}\right)\lambda(w_{(1)}\otimes
Y_2^{o}(Y(v_{1}, x_0)v_{2}, x_2)w_{(2)})\nn
&&\quad+\res_{x_{0}}x_{2}^{-1}\delta\left(\frac{x_{1}-x_{0}}{x_{2}}\right)
\res_{y_{2}}
z^{-1}\delta\left(\frac{x_{2}^{-1}-y_{2}}{z}\right)\cdot\nn
&&\quad\quad\quad\quad\cdot\lambda(
Y_{1}(e^{x_2L(1)}(-x_2^{-2})^{L(0)}
Y(v_{1},
x_{0})v_{2}, y_{2})w_{(1)}\otimes
w_{(2)})\nn
&&=\res_{x_{0}}x_{2}^{-1}
\delta\left(\frac{x_{1}-x_{0}}{x_{2}}\right)
(Y'_{P(z)}(Y(v_{1}, x_{0})v_{2}, x_{2})\lambda)(w_{(1)}\otimes w_{(2)}).
\eea
Since $\lambda$, $w_{(1)}$ and $w_{(2)}$ are arbitrary,
this  equality gives the commutator formula for $Y'_{P(z)}$.
\epfv

The following observations are analogous to those in Remark 8.1 of
\cite{tensor2} (concerning the case of $Q(z)$ rather than $P(z)$):

\begin{rema}
{\rm The proof of Proposition \ref{pz-comm} suggests the following:  
Using the definitions (\ref{taudef0}) and (\ref{taudef})
as motivation, we define 
a (linear) action
$\sigma_{P(z)}$ of $V \otimes \iota_{+}{\mathbb C}[t,t^{- 1}, 
(z^{-1}-t)^{-1}]$
on the 
vector space $W_{1}\otimes W_{2}$ (as opposed to $(W_{1}\otimes W_{2})^{*}$)
as follows:
\begin{equation}
\sigma_{P(z)}(\xi)(w_{(1)}\otimes w_{(2)})
=
\tau_{W_{1}}((\iota_+\circ T_z\circ\iota_-^{-1}\circ o)\xi)w_{(1)}
\otimes w_{(2)}
+w_{(1)}\otimes\tau_{W_{2}}((\iota_+\circ
\iota_-^{-1}\circ o)\xi)w_{(2)}
\end{equation}
for $\xi\in V \otimes \iota_{+}{\mathbb C}[t,t^{- 1}, 
(z^{-1}-t)^{-1}]$, $w_{(1)}\in
W_{1}$, $w_{(2)}\in W_{2}$, or equivalently,
\begin{eqnarray}\label{sigma-p-z}
\lefteqn{\bigg(\sigma_{P(z)}\bigg(x_0^{-1}
\delta\bigg(\frac{x^{-1}_1-z}{x_0}\bigg)
Y_{t}(v, x_1)\bigg)\bigg)
(w_{(1)}\otimes w_{(2)})}\nno\\
&&=z^{-1}\delta\bigg(\frac{x^{-1}_1-x_0}{z}\bigg)
Y_1(e^{x_1L(1)}(-x_1^{-2})^{L(0)}v, x_0)w_{(1)}\otimes w_{(2)}
\nno\\
&&\quad +x^{-1}_0\delta\bigg(\frac{z-x^{-1}_1}{-x_0}\bigg)
w_{(1)}\otimes Y_2^{o}(v, x_1)w_{(2)}.
\end{eqnarray}
That is, the operators $\sigma_{P(z)}(\xi)$ and $\tau_{P(z)}(\xi)$ are 
mutually adjoint:
\begin{equation}
(\tau_{P(z)}(\xi)\lambda)(w_{(1)}\otimes w_{(2)})
=\lambda(\sigma_{P(z)}(\xi)(w_{(1)}\otimes w_{(2)})).
\end{equation}
While this action on $W_{1}\otimes W_{2}$ is not very useful, it has the 
following three properties:
\begin{equation}\label{sigma-id}
\sigma_{P(z)}(Y_{t}({\bf 1}, x))=1,
\end{equation}
\begin{equation}\label{sigma-dev}
\frac{d}{dx}\sigma_{P(z)}(Y_{t}(v, x))=\sigma_{P(z)}(Y_{t}(L(-1)v, x))
\end{equation}
for $v\in V$,
\begin{eqnarray}\label{sigma-comm}
\lefteqn{[\sigma_{P(z)}(Y_{t}(v_{2}, x_{2})), \sigma_{P(z)}(Y_{t}(v_{1}, 
x_{1}))]}\nno\\
&&=\res_{x_{0}}x_{2}^{-1}\delta\left(\frac{x_{1}-x_{0}}{x_{2}}\right)
\sigma_{P(z)}(Y_{t}(Y(v_{1}, x_{0})v_{2}, x_{2}))
\end{eqnarray}
for $v_{1}, v_{2}\in V$ (the opposite commutator formula). These
follow either {}from the assertions of Propositions \ref{id-dev} and 
\ref{pz-comm} or,
better, {}from the fact that it was actually 
(\ref{sigma-id})--(\ref{sigma-comm}) that the
proofs of these propositions were proving.}
\end{rema}

\begin{rema}
{\rm Taking $\res_{x_{0}}$ of (\ref{sigma-p-z}),
we obtain
\begin{eqnarray}\label{sigma-p-z-1}
\lefteqn{(\sigma_{P(z)}(Y_{t}(v, x_1)))
(w_{(1)}\otimes w_{(2)})}\nno\\
&&=\res_{x_{0}}z^{-1}\delta\bigg(\frac{x^{-1}_1-x_0}{z}\bigg)
Y_1(e^{x_1L(1)}(-x_1^{-2})^{L(0)}v, x_0)w_{(1)}\otimes w_{(2)}
\nno\\
&&\quad +
w_{(1)}\otimes Y_2^{o}(v, x_1)w_{(2)}.
\end{eqnarray}
Substituting first $(-x_1^{-2})^{-L(0)}e^{-x_1L(1)}v$ for $v$ in
(\ref{sigma-p-z-1}) and then $x^{-1}_{1}$ for $x_{1}$ in the same
formula and using (\ref{xe^Lx}), (\ref{3.5}) and (\ref{op-y-t-2}), we
obtain
\begin{eqnarray}\label{sigma-p-z-1.5}
(\sigma_{P(z)}(Y^{o}_{t}(v, x_1)))
(w_{(1)}\otimes w_{(2)})
&=&\res_{x_{0}}z^{-1}\delta\bigg(\frac{x_1-x_0}{z}\bigg)
Y_1(v, x_0)w_{(1)}\otimes w_{(2)}
\nno\\
&& +
w_{(1)}\otimes Y_2(v, x_1)w_{(2)}.
\end{eqnarray}
Using this, we see that $\sigma_{P(z)}$ can actually be viewed as a
map {}from $V[t, t^{-1}]$ to $V((t))\otimes V[t, t^{-1}]$ defined
by
\begin{eqnarray}\label{sigma-p-z-2}
\sigma_{P(z)}(Y^{o}_{t}(v, x_1))
&=&\res_{x_{0}}z^{-1}\delta\bigg(\frac{x_1-x_0}{z}\bigg)
Y_t(v, x_0)\otimes {\mathbf{1}}
\nno\\
&& +
{\mathbf{1}}\otimes Y_t(v, x_1).
\end{eqnarray}
Let $\Delta_{P(z)}=\sigma_{P(z)}\circ o$.  Then (\ref{sigma-p-z-2})
becomes
\begin{eqnarray}\label{sigma-p-z-3}
\Delta_{P(z)}(Y_{t}(v, x_1))
&=&\res_{x_{0}}z^{-1}\delta\bigg(\frac{x_1-x_0}{z}\bigg)
Y_t(v, x_0)\otimes {\mathbf{1}}
\nno\\
&& +
{\mathbf{1}}\otimes Y_t(v, x_1),
\end{eqnarray}
which again can be viewed as a map
\begin{equation}
\Delta_{P(z)}:V[t, t^{-1}] \rightarrow V((t))\otimes V[t, t^{-1}].
\end{equation}
In formula (2.4) of \cite{MS}, Moore and Seiberg introduced a map
$\Delta_{z, 0}$ which in fact corresponds exactly to the map
$\Delta_{P(z)}$ defined by (\ref{sigma-p-z-3}).  They proposed to
define a $V$-module structure (called ``a representation of
$\mathcal{A}$'' in \cite{MS}, where $\mathcal{A}$ corresponds to our
vertex algebra $V$) on $W_{1}\otimes W_{2}$ by using this map, which
can be viewed as a sort of analogue of a coproduct, but they
acknowledged that E. Witten pointed out ``subtleties in this
definition which are related to the fact that [a] representation of
$\mathcal{A}$ obtained this way is not always a highest weight
representation.''  In fact, it is these subtleties that make it
impossible to work with $W_{1}\otimes W_{2}$; in virtually all
interesting cases, $W_{1}\otimes W_{2}$ does not have a natural
(generalized) $V$-module structure.  This is exactly the reason why we
had to use a completely different approarch to construct our
$P(z)$-tensor product of $W_{1}$ and $W_{2}$.}
\end{rema}

When $V$ is in fact a conformal (rather than M\"{o}bius) vertex
algebra, we will write
\begin{equation}\label{13.11}
Y'_{P(z)}(\omega,x)=\sum_{n\in {\mathbb Z}} L'_{P(z)}(n)x^{-n-2}.
\end{equation}
Then {}from the last two propositions we see that the coefficient operators
of $Y'_{P(z)}(\omega, x)$ satisfy the Virasoro algebra commutator
relations, that is,
\[
[L'_{P(z)}(m), L'_{P(z)}(n)]
=(m-n)L'_{P(z)}(m+n)+{\displaystyle\frac1{12}}
(m^3-m)\delta_{m+n,0}c,
\]
with $c\in \C$ the central charge of $V$ (recall Definition 
\ref{cva}).
Moreover, in this case, by setting $v=\omega$ in (\ref{Y'def}) and
taking the coefficient of $x_1^{-j-2}$ for $j=-1, 0, 1$, we find that
\begin{eqnarray}\label{LP'(j)}
\lefteqn{(L'_{P(z)}(j)\lambda)(w_{(1)}\otimes w_{(2)})}\nn
&&=\lambda\Biggl(w_{(1)}\otimes L(-j)w_{(2)}+\Biggl(\sum_{i=0}
^{1-j}{1-j\choose i}z^iL(-j-i)\Biggr)w_{(1)}\otimes w_{(2)}
\Biggr),
\end{eqnarray}
by (\ref{Yoppositeomega}). If $V$ is just a M\"obius vertex algebra, we
define the actions $L'_{P(z)}(j)$ on $(W_1\otimes W_2)^*$ by
(\ref{LP'(j)}) for $j=-1, 0$ and $1$. 

\begin{rema}\label{I-intw2}{\rm
In view of the action $L'_{P(z)}(j)$, the ${\mathfrak s}{\mathfrak
l}(2)$-bracket relations (\ref{im:Lj}) for a $P(z)$-intertwining map
$I$, with notation as in Definition \ref{im:imdef}, can be written as
\begin{equation}\label{I-intw2f}
(L'(j)w'_{(3)})\circ I=L'_{P(z)}(j)(w'_{(3)}\circ I)
\end{equation}
for $w'_{(3)}\in W'_{3}$ and  $j=-1$, $0$ and $1$ (cf. 
(\ref{im:def'}) and Remark
\ref{I-intw}).
}
\end{rema}

\begin{rema}\label{L'jpreservesbetaspace}
{\rm We have 
\[
L'_{P(z)}(j)((W_{1}\otimes W_{2})^{*})^{(\beta)}\subset ((W_{1}\otimes
W_{2})^{*})^{(\beta)}
\]
for $j=-1, 0, 1$ and $\beta\in \tilde{A}$ (cf. Proposition \ref{3.6}).}
\end{rema}

When $V$ is a conformal vertex algebra, {}from the commutator formula
for $Y'_{P(z)}(\omega,x)$, we see that $L'_{P(z)}(-1)$, $L'_{P(z)}(0)$
and $L'_{P(z)}(1)$ realize the actions of $L_{-1}$, $L_0$ and $L_1$ in
${\mathfrak s}{\mathfrak l}(2)$ (recall (\ref{L_*})) on $(W_1\otimes
W_2)^*$.  When $V$ is just a M\"obius vertex algebra, the same
conclusion still holds but a proof is needed.  We state this as a
proposition:

\begin{propo}\label{sl-2}
Let $V$ be a M\"obius vertex algebra and let $W_{1}$ and $W_{2}$ be
generalized $V$-modules.  Then the operators $L'_{P(z)}(-1)$,
$L'_{P(z)}(0)$ and $L'_{P(z)}(1)$ realize the actions of $L_{-1}$,
$L_0$ and $L_1$ in ${\mathfrak s}{\mathfrak l}(2)$ on $(W_1\otimes
W_2)^*$, according to (\ref{L_*}).
\end{propo}
\pf
For $\lambda\in (W_{1}\otimes W_{2})^{*}$, $w_{(1)}\in W_{1}$,
$w_{(2)}\in W_{2}$ and $j, k=-1, 0, 1$, we have 
\begin{eqnarray}\label{kj}
\lefteqn{(L'_{P(z)}(j)L'_{P(z)}(k)\lambda)(w_{(1)}\otimes w_{(2)})}\nn
&&=(L'_{P(z)}(k)\lambda)
\Biggl(w_{(1)}\otimes L(-j)w_{(2)}+\Biggl(\sum_{i=0}
^{1-j}{1-j\choose i}z^iL(-j-i)\Biggr)w_{(1)}\otimes w_{(2)}
\Biggr)\nn
&&=(L'_{P(z)}(k)\lambda)
(w_{(1)}\otimes L(-j)w_{(2)})\nn
&&\quad +(L'_{P(z)}(k)\lambda)
\Biggl(\Biggl(\sum_{i=0}
^{1-j}{1-j\choose i}z^iL(-j-i)\Biggr)w_{(1)}\otimes w_{(2)}\Biggr)\nn
&&=\lambda
\Biggl(w_{(1)}\otimes L(-k)L(-j)w_{(2)}+\Biggl(\sum_{l=0}
^{1-k}{1-k\choose l}z^lL(-k-l)\Biggr)w_{(1)}\otimes L(-j)w_{(2)}\Biggr)
\nn
&&\quad +\lambda
\Biggl(\Biggl(\sum_{i=0}
^{1-j}{1-j\choose i}z^iL(-j-i)\Biggr)w_{(1)}\otimes L(-k)w_{(2)}\nn
&&\quad\quad\quad\quad +\Biggl(\sum_{l=0}
^{1-k}{1-k\choose l}z^lL(-k-l)\Biggr)\Biggl(\sum_{i=0}
^{1-j}{1-j\choose i}z^iL(-j-i)\Biggr)w_{(1)}\otimes w_{(2)}\Biggr).\nn
&&
\end{eqnarray}
{}From formula (\ref{kj}) we obtain
\begin{eqnarray}\label{kj-comm}
\lefteqn{([L'_{P(z)}(j), L'_{P(z)}(k)]\lambda)(w_{(1)}\otimes w_{(2)})}\nn
&&=\lambda
\Biggl(w_{(1)}\otimes [L(-k), L(-j)]w_{(2)}\nn
&&\quad\quad\quad +\Biggl(\sum_{l=0}^{1-k}\sum_{i=0}^{1-j}
{1-k\choose l}{1-j\choose i}z^{l+i}
[L(-k-l), L(-j-i)]\Biggr)w_{(1)}\otimes w_{(2)}\Biggr)\nn
&&=\lambda
\Biggl(w_{(1)}\otimes (j-k)L(-k-j)w_{(2)}\nn
&&\quad\quad\quad +\Biggl(\sum_{l=0}^{1-k}\sum_{i=0}^{1-j}
{1-k\choose l}{1-j\choose i}z^{l+i}(j+i-k-l)
L(-k-l-j-i)\Biggr)w_{(1)}\otimes w_{(2)}\Biggr)\nn
\end{eqnarray}

Taking $j=1$ and $k=-1, 0$ in (\ref{kj-comm}), 
we obtain 
\begin{eqnarray*}
\lefteqn{([L'_{P(z)}(1), L'_{P(z)}(k)]\lambda)(w_{(1)}\otimes w_{(2)})}\nn
&&=\lambda
\Biggl(w_{(1)}\otimes (1-k)L(-k-1)w_{(2)}\nn
&&\quad\quad\quad +\Biggl(\sum_{l=0}^{1-k}
{1-k\choose l}z^{l}(1-k-l)
L(-k-l-1)\Biggr)w_{(1)}\otimes w_{(2)}\Biggr)\nn
&&=\lambda
\Biggl(w_{(1)}\otimes (1-k)L(-k-1)w_{(2)}\nn
&&\quad\quad\quad +\Biggl(\sum_{l=0}^{1-(k+1)}
{1-k\choose l}z^{l}(1-k-l)
L(-k-l-1)\Biggr)w_{(1)}\otimes w_{(2)}\Biggr)\nn
&&=(1-k)\lambda
\Biggl(w_{(1)}\otimes L(-1-k)w_{(2)}\nn
&&\quad\quad\quad\quad\quad\quad\quad +\Biggl(\sum_{l=0}^{1-(k+1)}
{1-(k+1)\choose l}z^{l}
L(-(k+1)-l)\Biggr)w_{(1)}\otimes w_{(2)}\Biggr)\nn
&&=((1-k)L'_{P(z)}(1+k)\lambda)(w_{(1)}\otimes w_{(2)}),
\end{eqnarray*}
proving the commutator formula
\[
[L'_{P(z)}(1), L'_{P(z)}(k)]=(1-k)L'_{P(z)}(1+k)
\]
for $k=-1, 0$. 

Taking $j=0$ and $k=-1$ in (\ref{kj-comm}), 
we obtain 
\begin{eqnarray*}
\lefteqn{([L'_{P(z)}(0), L'_{P(z)}(-1)]\lambda)(w_{(1)}\otimes w_{(2)})}\nn
&&=\lambda
\Biggl(w_{(1)}\otimes L(1)w_{(2)} +\Biggl(\sum_{l=0}^{2}\sum_{i=0}^{1}
{2\choose l}z^{l+i}(i+1-l)
L(1-l-i)\Biggr)w_{(1)}\otimes w_{(2)}\Biggr)\nn
&&=\lambda(w_{(1)}\otimes L(1)w_{(2)} +
(L(1)+2zL(0)+z^{2}L(-1))w_{(1)}\otimes w_{(2)})\nn
&&=\lambda\Biggl(w_{(1)}\otimes L(1)w_{(2)} 
+\Biggl(\sum_{m=0}^{1-(-1)}{2\choose m}z^{m}L(-(-1)-m)\Biggr)
w_{(1)}\otimes w_{(2)}\Biggr)\nn
&&=(L'_{P(z)}(-1)\lambda)(w_{(1)}\otimes w_{(2)}),
\end{eqnarray*}
proving that
\[
[L'_{P(z)}(0), L'_{P(z)}(-1)]=L'_{P(z)}(-1).
\]
The three commutator formulas we have proved show that 
$L'_{P(z)}(-1)$,
$L'_{P(z)}(0)$ and $L'_{P(z)}(1)$ indeed realize the actions of $L_{-1}$,
$L_0$ and $L_1$ in ${\mathfrak s}{\mathfrak l}(2)$.
\epfv

The commutator formulas corresponding to (\ref{sl2-1})--(\ref{sl2-3})
(recall Definition \ref{moduleMobius})) also need to be proved in the
M\"obius case:

\begin{propo}\label{pz-l-y-comm}
Let $V$ be a M\"obius vertex algebra and let $W_{1}$ and $W_{2}$ be
generalized $V$-modules.  Then for $v\in V$, we have the following
commutator formulas:
\begin{eqnarray}
{[L(-1), Y'_{P(z)}(v, x)]}&=&Y'_{P(z)}(L(-1)v, x),\label{pz-sl-2-pz-y--2}\\
{[L(0), Y'_{P(z)}(v, x)]}&=&Y'_{P(z)}(L(0)v, x)+xY'_{P(z)}(L(-1)v, x),
\label{pz-sl-2-pz-y--1}\\
{[L(1), Y'_{P(z)}(v, x)]}&=&Y'_{P(z)}(L(1)v, x)
+2xY'_{P(z)}(L(0)v, x)+x^{2}Y'_{P(z)}(L(-1)v, x),\label{pz-sl-2-pz-y}\nn
&&
\end{eqnarray}
where for brevity we write $L'_{P(z)}(j)$ acting on $(W_1\otimes
W_2)^*$ as $L(j)$.
\end{propo}
\pf Let $\lambda\in (W_{1}\otimes W_{2})^{*}$, $w_{(1)}\in W_{1}$ and
$w_{(2)}\in W_{2}$. Using (\ref{LP'(j)}), (\ref{Y'def}), the
commutator formulas for $L(j)$ and $Y_{1}(v, x_{0})$ for $j=-1, 0, 1$
and $v\in V$ (recall Definition \ref{moduleMobius}), and the
commutator formulas for $L(j)$ and $Y_{2}^{o}(v, x)$ for $j=-1, 0, 1$
and $v\in V$ (recall Lemma \ref{sl2opposite}), we obtain, for $j=-1,
0, 1$,
\begin{eqnarray}\label{pz-sl-2-pz-y-1}
\lefteqn{([L(j), Y'_{P(z)}(v, x)]\lambda)(w_{(1)}\otimes w_{(2)})}\nn
&&=(L(j)Y'_{P(z)}(v, x)\lambda)(w_{(1)}\otimes w_{(2)})
-(Y'_{P(z)}(v, x)L(j)\lambda)(w_{(1)}\otimes w_{(2)})\nn
&&=(Y'_{P(z)}(v, x)\lambda)(w_{(1)}\otimes L(-j)w_{(2)})\nn
&&\quad +\sum_{i=0}^{1-j}{1-j\choose i}
(Y'_{P(z)}(v, x)\lambda)(z^iL(-j-i)w_{(1)}\otimes w_{(2)})\nn
&&\quad -(L(j)\lambda)(w_{(1)}\otimes Y_{2}^{o}(v, x)w_{(2)})\nn
&&\quad -\res_{x_{0}}z^{-1}\delta\left(\frac{x^{-1}-x_{0}}{z}\right)
(L(j)\lambda)(Y_{1}(e^{xL(1)}(-x^{-2})^{L(0)}v, x_{0})w_{(1)}
\otimes w_{(2)})\nn
&&=\lambda(w_{(1)}\otimes Y^{o}_{2}(v, x)L(-j)w_{(2)})\nn
&&\quad +\res_{x_{0}}z^{-1}\delta\left(\frac{x^{-1}-x_{0}}{z}\right)
\lambda(Y_{1}(e^{xL(1)}(-x^{-2})^{L(0)}v, x_{0})w_{(1)}
\otimes L(-j)w_{(2)})\nn
&&\quad +\sum_{i=0}^{1-j}{1-j\choose i}\lambda(z^iL(-j-i)w_{(1)}
\otimes Y_{2}^{o}(v, x)w_{(2)})\nn
&&\quad +\sum_{i=0}^{1-j}{1-j\choose i}
\res_{x_{0}}z^{-1}\delta\left(\frac{x^{-1}-x_{0}}{z}\right)\cdot\nn
&&\quad\quad\quad\quad\quad\quad\quad\quad\quad\quad\quad\quad \cdot
\lambda(Y_{1}(e^{xL(1)}(-x^{-2})^{L(0)}v, x_{0})
z^iL(-j-i)w_{(1)}\otimes w_{(2)})\nn
&&\quad -\lambda(w_{(1)}\otimes L(-j)Y_{2}^{o}(v, x)w_{(2)})\nn
&&\quad -\sum_{i=0}^{1-j}{1-j\choose i}
\lambda(z^iL(-j-i)w_{(1)}\otimes Y_{2}^{o}(v, x)w_{(2)})\nn
&&\quad -\res_{x_{0}}z^{-1}\delta\left(\frac{x^{-1}-x_{0}}{z}\right)
\lambda(Y_{1}(e^{xL(1)}(-x^{-2})^{L(0)}v, x_{0})w_{(1)}
\otimes L(-j)w_{(2)})\nn
&&\quad -\sum_{i=0}^{1-j}{1-j\choose i}
\res_{x_{0}}z^{-1}\delta\left(\frac{x^{-1}-x_{0}}{z}\right)\cdot\nn
&&\quad\quad\quad\quad\quad\quad\quad\quad\quad\quad\quad\quad \cdot
\lambda(z^iL(-j-i)Y_{1}(e^{xL(1)}(-x^{-2})^{L(0)}v, x_{0})w_{(1)}
\otimes w_{(2)})\nn
&&=\lambda(w_{(1)}\otimes [Y^{o}_{2}(v, x), L(-j)]w_{(2)})\nn
&&\quad -\sum_{i=0}^{1-j}{1-j\choose i}
\res_{x_{0}}z^{-1}\delta\left(\frac{x^{-1}-x_{0}}{z}\right)\cdot\nn
&&\quad\quad\quad\quad\quad\quad\quad\quad\quad\quad\quad\quad \cdot
\lambda([z^iL(-j-i), Y_{1}(e^{xL(1)}(-x^{-2})^{L(0)}v, x_{0})]
w_{(1)}\otimes w_{(2)})\nn
&&=\sum_{k=0}^{j+1}{j+1\choose k}x^{j+1-k}
\lambda(w_{(1)}\otimes Y^{o}_{2}(L(k-1)v, x)w_{(2)})\nn
&&\quad -\res_{x_{0}}z^{-1}\delta\left(\frac{x^{-1}-x_{0}}{z}\right)
\sum_{i=0}^{1-j}\sum_{k=0}^{-j-i+1}{1-j\choose i}
{-j-i+1\choose k}z^ix_{0}^{-j-i+1-k}
\cdot\nn
&&\quad\quad\quad\quad\quad\quad\quad\quad\quad\quad\quad\quad \cdot
\lambda(Y_{1}(L(k-1)e^{xL(1)}(-x^{-2})^{L(0)}v, x_{0})
w_{(1)}\otimes w_{(2)}).\nn
\end{eqnarray}
Using (\ref{2termdeltarelation}) and
(\ref{deltafunctionsubstitutionformula}) when necessary, we see that
the second term on the right-hand side of (\ref{pz-sl-2-pz-y-1}) is
equal to the following expressions for $j=1, 0$ and $-1$, respectively:
\begin{eqnarray}\label{pz-sl-2-pz-y-2}
-\res_{x_{0}}z^{-1}\delta\left(\frac{x^{-1}-x_{0}}{z}\right)
\lambda(Y_{1}(L(-1)e^{xL(1)}(-x^{-2})^{L(0)}v, x_{0})
w_{(1)}\otimes w_{(2)}),
\end{eqnarray}
\begin{eqnarray}\label{pz-sl-2-pz-y-3}
\lefteqn{-\res_{x_{0}}z^{-1}\delta\left(\frac{x^{-1}-x_{0}}{z}\right)
\sum_{i=0}^{1}\sum_{k=0}^{-i+1}{1\choose i}
{-j-i+1\choose k}z^ix_{0}^{-j-i+1-k}
\cdot}\nn
&&\quad\quad\quad\quad\quad\quad\quad\cdot
\lambda(Y_{1}(L(k-1)e^{xL(1)}(-x^{-2})^{L(0)}v, x_{0})
w_{(1)}\otimes w_{(2)})\nn
&&=-\res_{x_{0}}z^{-1}\delta\left(\frac{x^{-1}-x_{0}}{z}\right)
x_{0}
\lambda(Y_{1}(L(-1)e^{xL(1)}(-x^{-2})^{L(0)}v, x_{0})
w_{(1)}\otimes w_{(2)})\nn
&&\quad -\res_{x_{0}}z^{-1}\delta\left(\frac{x^{-1}-x_{0}}{z}\right)
\lambda(Y_{1}(L(0)e^{xL(1)}(-x^{-2})^{L(0)}v, x_{0})
w_{(1)}\otimes w_{(2)})\nn
&&\quad -\res_{x_{0}}z^{-1}\delta\left(\frac{x^{-1}-x_{0}}{z}\right)
z\lambda(Y_{1}(L(-1)e^{xL(1)}(-x^{-2})^{L(0)}v, x_{0})
w_{(1)}\otimes w_{(2)})\nn
&&=-\res_{x_{0}}x\delta\left(\frac{z+x_{0}}{x^{-1}}\right)
(x_{0}+z)
\lambda(Y_{1}(L(-1)e^{xL(1)}(-x^{-2})^{L(0)}v, x_{0})
w_{(1)}\otimes w_{(2)})\nn
&&\quad -\res_{x_{0}}z^{-1}\delta\left(\frac{x^{-1}-x_{0}}{z}\right)
\lambda(Y_{1}(L(0)e^{xL(1)}(-x^{-2})^{L(0)}v, x_{0})
w_{(1)}\otimes w_{(2)})\nn
&&=-\res_{x_{0}}z^{-1}\delta\left(\frac{x^{-1}-x_{0}}{z}\right)x^{-1}
\lambda(Y_{1}(L(-1)e^{xL(1)}(-x^{-2})^{L(0)}v, x_{0})
w_{(1)}\otimes w_{(2)})\nn
&&\quad -\res_{x_{0}}z^{-1}\delta\left(\frac{x^{-1}-x_{0}}{z}\right)
\lambda(Y_{1}(L(0)e^{xL(1)}(-x^{-2})^{L(0)}v, x_{0})
w_{(1)}\otimes w_{(2)})\nn
\end{eqnarray}
and
\begin{eqnarray}\label{pz-sl-2-pz-y-4}
\lefteqn{-\res_{x_{0}}z^{-1}\delta\left(\frac{x^{-1}-x_{0}}{z}\right)
\sum_{i=0}^{2}\sum_{k=0}^{2-i}{2\choose i}
{2-i\choose k}z^ix_{0}^{2-i-k}
\cdot}\nn
&&\quad\quad\quad\quad\quad\quad\quad\quad\quad\quad\quad\quad \cdot
\lambda(Y_{1}(L(k-1)e^{xL(1)}(-x^{-2})^{L(0)}v, x_{0})
w_{(1)}\otimes w_{(2)})\nn
&&=-\res_{x_{0}}z^{-1}\delta\left(\frac{x^{-1}-x_{0}}{z}\right)
x_{0}^{2}
\lambda(Y_{1}(L(-1)e^{xL(1)}(-x^{-2})^{L(0)}v, x_{0})
w_{(1)}\otimes w_{(2)})\nn
&&\quad -\res_{x_{0}}z^{-1}\delta\left(\frac{x^{-1}-x_{0}}{z}\right)
2x_{0}
\lambda(Y_{1}(L(0)e^{xL(1)}(-x^{-2})^{L(0)}v, x_{0})
w_{(1)}\otimes w_{(2)})\nn
&&\quad -\res_{x_{0}}z^{-1}\delta\left(\frac{x^{-1}-x_{0}}{z}\right)
\lambda(Y_{1}(L(1)e^{xL(1)}(-x^{-2})^{L(0)}v, x_{0})
w_{(1)}\otimes w_{(2)})\nn
&&\quad -\res_{x_{0}}z^{-1}\delta\left(\frac{x^{-1}-x_{0}}{z}\right)
2zx_{0}
\lambda(Y_{1}(L(-1)e^{xL(1)}(-x^{-2})^{L(0)}v, x_{0})
w_{(1)}\otimes w_{(2)})\nn
&&\quad -\res_{x_{0}}z^{-1}\delta\left(\frac{x^{-1}-x_{0}}{z}\right)
2z\lambda(Y_{1}(L(0)e^{xL(1)}(-x^{-2})^{L(0)}v, x_{0})
w_{(1)}\otimes w_{(2)})\nn
&&\quad -\res_{x_{0}}z^{-1}\delta\left(\frac{x^{-1}-x_{0}}{z}\right)
z^2\lambda(Y_{1}(L(-1)e^{xL(1)}(-x^{-2})^{L(0)}v, x_{0})
w_{(1)}\otimes w_{(2)})\nn
&&=-\res_{x_{0}}z^{-1}\delta\left(\frac{x^{-1}-x_{0}}{z}\right)
x^{-2}
\lambda(Y_{1}(L(-1)e^{xL(1)}(-x^{-2})^{L(0)}v, x_{0})
w_{(1)}\otimes w_{(2)})\nn
&&\quad -\res_{x_{0}}z^{-1}\delta\left(\frac{x^{-1}-x_{0}}{z}\right)
2x^{-1}
\lambda(Y_{1}(L(0)e^{xL(1)}(-x^{-2})^{L(0)}v, x_{0})
w_{(1)}\otimes w_{(2)})\nn
&&\quad -\res_{x_{0}}z^{-1}\delta\left(\frac{x^{-1}-x_{0}}{z}\right)
\lambda(Y_{1}(L(1)e^{xL(1)}(-x^{-2})^{L(0)}v, x_{0})
w_{(1)}\otimes w_{(2)}).\nn
\end{eqnarray}
Using (\ref{log:SL2-3}) and (\ref{log:xLx^}), we see that
(\ref{pz-sl-2-pz-y-2}), (\ref{pz-sl-2-pz-y-3}) and
(\ref{pz-sl-2-pz-y-4}) are respectively equal to
\begin{eqnarray}\label{pz-sl-2-pz-y-5}
\lefteqn{-\res_{x_{0}}z^{-1}\delta\left(\frac{x^{-1}-x_{0}}{z}\right)
\cdot}\nn
&&\quad\quad
\cdot\lambda(Y_{1}(e^{xL(1)}(x^{2}L(1)-2xL(0)+L(-1))(-x^{-2})^{L(0)}v, x_{0})
w_{(1)}\otimes w_{(2)})\nn
&&=\res_{x_{0}}z^{-1}\delta\left(\frac{x^{-1}-x_{0}}{z}\right)\cdot\nn
&&\quad\quad
\cdot
\lambda(Y_{1}(e^{xL(1)}(-x^{-2})^{L(0)}(L(1)+2xL(0)+x^{2}L(-1))v, x_{0})
w_{(1)}\otimes w_{(2)}),\nn
\end{eqnarray}
\begin{eqnarray}\label{pz-sl-2-pz-y-6}
\lefteqn{-\res_{x_{0}}z^{-1}\delta\left(\frac{x^{-1}-x_{0}}{z}\right)x^{-1}
\cdot}\nn
&&\quad\quad\cdot
\lambda(Y_{1}(e^{xL(1)}(x^{2}L(1)-2xL(0)+L(-1))(-x^{-2})^{L(0)}v, x_{0})
w_{(1)}\otimes w_{(2)})\nn
&&\quad -\res_{x_{0}}z^{-1}\delta\left(\frac{x^{-1}-x_{0}}{z}\right)
\lambda(Y_{1}(e^{xL(1)}(-xL(1)+L(0))(-x^{-2})^{L(0)}v, x_{0})
w_{(1)}\otimes w_{(2)})\nn
&&=\res_{x_{0}}z^{-1}\delta\left(\frac{x^{-1}-x_{0}}{z}\right)x^{-1}\cdot\nn
&&\quad\quad\cdot
\lambda(Y_{1}(e^{xL(1)}(-x^{-2})^{L(0)}(L(1)+2xL(0)+x^{2}L(-1))v, x_{0})
w_{(1)}\otimes w_{(2)})\nn
&&\quad +\res_{x_{0}}z^{-1}\delta\left(\frac{x^{-1}-x_{0}}{z}\right)
\lambda(Y_{1}(e^{xL(1)}(-x^{-2})^{L(0)}(-x^{-1}L(1)-L(0))v, x_{0})
w_{(1)}\otimes w_{(2)})\nn
&&=\res_{x_{0}}z^{-1}\delta\left(\frac{x^{-1}-x_{0}}{z}\right)
\lambda(Y_{1}(e^{xL(1)}(-x^{-2})^{L(0)}(L(0)+xL(-1))v, x_{0})
w_{(1)}\otimes w_{(2)})\nn
\end{eqnarray}
and
\begin{eqnarray}\label{pz-sl-2-pz-y-7}
\lefteqn {-\res_{x_{0}}z^{-1}\delta\left(\frac{x^{-1}-x_{0}}{z}\right)
x^{-2}\cdot}\nn
&&\quad\quad\cdot
\lambda(Y_{1}(e^{xL(1)}(x^{2}L(1)-2xL(0)+L(-1))(-x^{-2})^{L(0)}v, x_{0})
w_{(1)}\otimes w_{(2)})\nn
&&\quad -\res_{x_{0}}z^{-1}\delta\left(\frac{x^{-1}-x_{0}}{z}\right)
2x^{-1}
\lambda(Y_{1}(e^{xL(1)}(-xL(1)+L(0))(-x^{-2})^{L(0)}v, x_{0})
w_{(1)}\otimes w_{(2)})\nn
&&\quad -\res_{x_{0}}z^{-1}\delta\left(\frac{x^{-1}-x_{0}}{z}\right)
\lambda(Y_{1}(e^{xL(1)}L(1)(-x^{-2})^{L(0)}v, x_{0})
w_{(1)}\otimes w_{(2)})\nn
&&=\res_{x_{0}}z^{-1}\delta\left(\frac{x^{-1}-x_{0}}{z}\right)
x^{-2}\cdot\nn
&&\quad\quad
\cdot
\lambda(Y_{1}(e^{xL(1)}(-x^{-2})^{L(0)}
(L(1)+2xL(0)+x^{2}L(-1))v, x_{0})
w_{(1)}\otimes w_{(2)})\nn
&&\quad +\res_{x_{0}}z^{-1}\delta\left(\frac{x^{-1}-x_{0}}{z}\right)
2x^{-1}
\lambda(Y_{1}(e^{xL(1)}(-x^{-2})^{L(0)}(-x^{-1}L(1)-L(0))v, x_{0})
w_{(1)}\otimes w_{(2)})\nn
&&\quad +\res_{x_{0}}z^{-1}\delta\left(\frac{x^{-1}-x_{0}}{z}\right)
\lambda(Y_{1}(e^{xL(1)}(-x^{-2})^{L(0)}x^{-2}L(1)v, x_{0})
w_{(1)}\otimes w_{(2)})\nn
&&=\res_{x_{0}}z^{-1}\delta\left(\frac{x^{-1}-x_{0}}{z}\right)
\lambda(Y_{1}(e^{xL(1)}(-x^{-2})^{L(0)}
L(-1)v, x_{0})
w_{(1)}\otimes w_{(2)}).
\end{eqnarray}
The right-hand sides of (\ref{pz-sl-2-pz-y-5}),
(\ref{pz-sl-2-pz-y-6}) and   (\ref{pz-sl-2-pz-y-7}) can be written as 
\[
\sum_{k=0}^{j+1}{j+1\choose k}x^{j+1-k}
\res_{x_{0}}z^{-1}\delta\left(\frac{x^{-1}-x_{0}}{z}\right)
\lambda(Y_{1}(e^{xL(1)}(-x^{-2})^{L(0)}
L(k-1)v, x_{0})
w_{(1)}\otimes w_{(2)}),
\]
for $j=1, 0, -1$, respectively. Thus the right-hand side of 
(\ref{pz-sl-2-pz-y-1}) is equal to 
\begin{eqnarray}\label{pz-sl-2-pz-y-8}
\lefteqn{\sum_{k=0}^{j+1}{j+1\choose k}x^{j+1-k}
\lambda(w_{(1)}\otimes Y^{o}_{2}(L(k-1)v, x)w_{(2)})}\nn
&&\quad +\sum_{k=0}^{j+1}{j+1\choose k}x^{j+1-k}
\res_{x_{0}}z^{-1}\delta\left(\frac{x^{-1}-x_{0}}{z}\right)\cdot\nn
&&\quad\quad\quad\quad\quad\quad\quad\quad\quad\quad\quad\quad\quad\quad
\cdot
\lambda(Y_{1}(e^{xL(1)}(-x^{-2})^{L(0)}
L(k-1)v, x_{0})
w_{(1)}\otimes w_{(2)})\nn
&&=\sum_{k=0}^{j+1}{j+1\choose k}x^{j+1-k}
(Y'_{P(z)}(L(k-1)v, x)\lambda)(w_{(1)}\otimes w_{(2)}),
\end{eqnarray}
proving the proposition.
\epfv

We have seen in (\ref{tauw}), (\ref{deltaY3'}) and (\ref{tausubW3'})
that for a generalized $V$-module $(W,Y_W)$, the space $V \otimes
{\C}((t))$, and in particular, the space $V \otimes \iota_{+}{\mathbb
C}[t,t^{- 1}, (z^{-1}-t)^{-1}]$, acts naturally on $W$ via the action
$\tau_W$, in view of (\ref{set:wtvn}) and Assumption \ref{assum};
recall that $v \otimes t^n$ ($v \in V$, $n \in {\mathbb Z}$) acts as
the component $v_n$ of $Y_W(v,x)$, and that more generally,
\begin{equation}\label{tau-w-comp}
\tau_W \left(v \otimes \sum_{n > N} a_nt^n\right) = \sum_{n > N} a_nv_n
\end{equation}
for $a_n \in {\mathbb C}$.  For generalized $V$-modules $W_1$, $W_2$
and $W_2$, we shall next relate the $P(z)$-intertwining maps of type
${W_3\choose W_1\, W_2}$ to certain linear maps {}from $W'_{3}$ to
$(W_1\otimes W_2)^{*}$ intertwining the actions of $V \otimes
\iota_{+}{\mathbb C}[t,t^{- 1}, (z^{-1}-t)^{-1}]$ and of ${\mathfrak
s} {\mathfrak l}(2)$ on $W'_{3}$ and on $(W_1\otimes W_2)^{*}$ (see
Proposition \ref{pz} and Notation \ref{scriptN} below).  For this, as
is suggested by Lemma \ref{4.36} and Proposition \ref{tau-a-comp}, we
need to consider $\tilde{A}$-compatibility for linear maps {}from $W_3$
to $(W_1\otimes W_2)^{*}$:

\begin{defi}
{\rm We call a map $J \in {\rm Hom}(W_3,(W_1\otimes W_2)^{*})$ 
{\it $\tilde{A}$-compatible} if
\begin{equation}\label{JAtildecompat}
J((W_{3})^{(\beta)}) \subset ((W_{1}\otimes W_{2})^{*})^{(\beta)}
\end{equation}
for $\beta \in {\tilde A}$.
}
\end{defi}

As in the discussion preceding Lemma \ref{4.36}, we see that an
element $\lambda$ of $(W_{1}\otimes W_{2}\otimes W_{3})^{*}$ amounts
exactly to a linear map
\[
J_{\lambda}: W_3 \rightarrow (W_{1}\otimes W_{2})^{*}.
\]
If $\lambda$ is $\tilde{A}$-compatible (see that discussion), then for
$w_{(1)}\in W_{1}^{(\beta)}$, $w_{(2)}\in W_{2}^{(\gamma)}$ and
$w_{(3)}\in W_{3}^{(\delta)}$ such that $\beta +\gamma +\delta \ne 0$,
\[
J_{\lambda}(w_{(3)})(w_{(1)}\otimes w_{(2)})=\lambda(w_{(1)}\otimes
w_{(2)}\otimes w_{(3)}) = 0,
\]
so that
\[
J_{\lambda}(w_{(3)}) \in ((W_{1}\otimes W_{2})^{*})^{(\delta)},
\]
and so $J_{\lambda}$ is $\tilde{A}$-compatible.  Similarly, if
$J_{\lambda}$ is $\tilde{A}$-compatible, then so is $\lambda$.  Thus
using Lemma \ref{4.36} we have:

\begin{lemma}\label{IlambdatoJlambda}
The linear functional $\lambda\in (W_{1}\otimes W_{2}\otimes
W_{3})^{*}$ is $\tilde{A}$-compatible if and only if $J_{\lambda}$ is
$\tilde{A}$-compatible.  The map given by $\lambda\mapsto J_{\lambda}$
is the unique linear isomorphism {}from the space of
$\tilde{A}$-compatible elements of $(W_{1}\otimes W_{2}\otimes
W_{3})^{*}$ to the space of $\tilde{A}$-compatible linear maps {}from
$W_3$ to $(W_{1}\otimes W_{2})^{*}$ such that
\[
J_{\lambda}(w_{(3)})(w_{(1)}\otimes w_{(2)})
=\lambda(w_{(1)}\otimes w_{(2)}\otimes w_{(3)})
\]
for $w_{(1)}\in W_{1}$, $w_{(2)}\in W_{2}$ and $w_{(3)}\in W_{3}$.  In
particular, the correspondence $I_{\lambda} \mapsto J_{\lambda}$
defines a (unique) linear isomorphism {}from the space of
$\tilde{A}$-compatible linear maps
\[
I = I_{\lambda}:W_{1}\otimes W_{2} \rightarrow \overline{W_{3}'}
\]
to the space of $\tilde{A}$-compatible linear maps
\[
J = J_{\lambda}: W_3 \rightarrow (W_{1}\otimes W_{2})^{*}
\]
such that
\[
\langle w_{(3)}, I(w_{(1)}\otimes w_{(2)})\rangle
=J(w_{(3)})(w_{(1)}\otimes w_{(2)})
\]
for $w_{(1)}\in W_{1}$, $w_{(2)}\in W_{2}$ and $w_{(3)}\in W_{3}$.
\epf
\end{lemma}

\begin{rema}\label{alternateformoflemma}{\rm
{}From Lemma \ref{IlambdatoJlambda} (with $W_3$ replaced by $W'_{3}$) we
have a canonical isomorphism {}from the space of
$\tilde{A}$-compatible linear maps
\[
I:W_{1}\otimes W_{2} \rightarrow \overline{W_3}
\]
to the space of $\tilde{A}$-compatible linear maps
\[
J:W'_{3} \rightarrow (W_{1}\otimes W_{2})^{*},
\]
determined by:
\begin{equation}\label{IcorrespondstoJ}
\langle w'_{(3)}, I(w_{(1)}\otimes w_{(2)})\rangle
=J(w'_{(3)})(w_{(1)}\otimes w_{(2)})
\end{equation}
for $w_{(1)}\in W_{1}$, $w_{(2)}\in W_{2}$ and $w'_{(3)}\in W'_{3}$,
or equivalently,
\begin{equation}\label{IcorrespondstoJalternateform}
w'_{(3)}\circ I = J(w'_{(3)})
\end{equation}
for $w'_{(3)} \in W'_{3}$.
}
\end{rema}

We introduce another notion, corresponding to the lower truncation
condition (\ref{im:ltc}) for $P(z)$-intertwining maps:

\begin{defi}\label{gradingrestrictedmapJ}{\rm
We call a map $J\in \hom(W_3, (W_1\otimes W_2)^{*})$ {\em grading
restricted} if for $n\in {\mathbb C}$, $w_{(1)}\in W_1$ and
$w_{(2)}\in W_2$,
\[
J((W_3)_{[n-m]})(w_{(1)}\otimes w_{(2)})=0\;\;\mbox{ for }\;m\in
{\mathbb N}\;\mbox{ sufficiently large.}
\]
}
\end{defi}

\begin{rema}\label{Jcompatimpliesgradingrestr}
{\rm
If $J \in {\rm Hom}(W_3,(W_1\otimes W_2)^{*})$ is
$\tilde{A}$-compatible, then $J$ is also grading restricted, as we see
using (\ref{set:dmltc}).
}
\end{rema}

\begin{rema}{\rm
Under the natural isomorphism given in Remark
\ref{alternateformoflemma} (see (\ref{IcorrespondstoJ})) in the
$\tilde{A}$-compatible setting, the map $J : W'_{3} \rightarrow
(W_{1}\otimes W_{2})^{*}$ is grading restricted (recall Definition
\ref{gradingrestrictedmapJ}) if and only if the map $I:W_{1}\otimes
W_{2} \rightarrow \overline{W_{3}}$ satisfies the lower truncation
condition (\ref{im:ltc}).  But notice also that in this
$\tilde{A}$-compatible setting, we have seen that both $I$ and $J$
automatically have these properties.
}
\end{rema}

Using the above together with Remarks \ref{I-intw} and \ref{I-intw2},
we now have the following result, generalizing Proposition 13.1 
in \cite{tensor3}:

\begin{propo}\label{pz}
Let $W_1$, $W_2$ and $W_3$ be generalized $V$-modules.  Under the
natural isomorphism described in Remark \ref{alternateformoflemma}
between the space of $\tilde{A}$-compatible linear maps
\[
I:W_{1}\otimes W_{2} \rightarrow \overline{W_{3}}
\]
and the space of $\tilde{A}$-compatible linear maps
\[
J:W'_{3} \rightarrow (W_{1}\otimes W_{2})^{*}
\]
determined by (\ref{IcorrespondstoJ}), the $P(z)$-intertwining maps
$I$ of type ${W_3\choose W_1\, W_2}$ correspond exactly to the
(grading restricted) $\tilde{A}$-compatible maps $J$ that intertwine
the actions of both
\[
V \otimes \iota_{+}{\mathbb C}[t,t^{- 1}, (z^{-1}-t)^{-1}]
\]
and ${\mathfrak s} {\mathfrak l}(2)$ on $W'_{3}$ and on $(W_1\otimes
W_2)^{*}$.
\end{propo}
\pf 
In view of (\ref{IcorrespondstoJalternateform}), Remark
\ref{I-intw} asserts that (\ref{im:def'}), or equivalently,
(\ref{im:def}), is equivalent to the condition
\begin{equation}\label{j-tau}
J\left(\tau_{W'_3}\left(x_0^{-1}\delta\left(\frac{x^{-1}_1-z}{x_0}\right)
Y_{t}(v, x_1)\right)w'_{(3)}\right)
=\tau_{P(z)}\left(x_0^{-1}\delta\left(\frac{x^{-1}_1-z}{x_0}\right)
Y_{t}(v, x_1)\right)J(w'_{(3)}),
\end{equation}
that is, the condition that $J$ intertwines the actions of $V \otimes
\iota_{+}{\mathbb C}[t,t^{- 1}, (z^{-1}-t)^{-1}]$ on $W'_{3}$ and on
$(W_1\otimes W_2)^{*}$ (recall (\ref{3.18-1})--(\ref{3.19-1})).
Similarly, Remark \ref{I-intw2} asserts that (\ref{im:Lj}) is
equivalent to the condition
\begin{equation}\label{j-lj}
J(L'(j)w'_{(3)}) = L'_{P(z)}(j)J(w'_{(3)})
\end{equation}
for $j=-1$, $0$, $1$, that is, the condition that $J$ intertwines the
actions of ${\mathfrak s} {\mathfrak l}(2)$ on $W'_{3}$ and on
$(W_1\otimes W_2)^{*}$.
\epfv

\begin{nota}\label{scriptN}
{\rm Given generalized $V$-modules $W_1$, $W_2$ and $W_3$, we shall
write ${\cal N}[P(z)]_{W'_3}^{(W_1 \otimes W_2)^{*}}$, or ${\cal
N}_{W'_3}^{(W_1 \otimes W_2)^{*}}$ if there is no ambiguity, for the
space of (grading restricted) $\tilde{A}$-compatible linear maps
\[
J:W'_{3} \rightarrow (W_{1}\otimes W_{2})^{*}
\]
that intertwine the actions of both
\[
V \otimes \iota_{+}{\mathbb C}[t,t^{- 1}, (z^{-1}-t)^{-1}]
\]
and ${\mathfrak s} {\mathfrak l}(2)$ on $W'_{3}$ and on $(W_1\otimes
W_2)^{*}$.
Note that Proposition \ref{pz} gives a natural linear isomorphism
\begin{eqnarray*}
{\cal M}[P(z)]^{W_3}_{W_1 W_2} = {\cal M}^{W_3}_{W_1 W_2} &
\stackrel{\sim}{\longrightarrow} & {\cal N}_{W'_3}^{(W_1 \otimes
W_2)^{*}}\nno\\ I & \mapsto & J
\end{eqnarray*}
(recall {}from Definition \ref{im:imdef} the notations for the space of
$P(z)$-intertwining maps).  Let us use the symbol ``prime'' to denote
this isomorphism in both directions:
\begin{eqnarray*}
{\cal M}^{W_3}_{W_1 W_2} & \stackrel{\sim}{\longrightarrow} & {\cal
N}_{W'_3}^{(W_1 \otimes W_2)^{*}}\nno\\
I & \mapsto & I'\nno\\
J' & \leftarrow\!\!\!{\scriptstyle |} & J,
\end{eqnarray*}
so that in particular,
\[
I'' = I \;\;\mbox{ and }\;\; J'' = J
\]
for $I \in {\cal M}^{W_3}_{W_1 W_2}$ and $J \in {\cal N}_{W'_3}^{(W_1
\otimes W_2)^{*}}$, and the relation between $I$ and $I'$ is
determined by
\[
\langle w'_{(3)}, I(w_{(1)}\otimes w_{(2)})\rangle
=I'(w'_{(3)})(w_{(1)}\otimes w_{(2)})
\]
for $w_{(1)}\in W_{1}$, $w_{(2)}\in W_{2}$ and $w'_{(3)}\in W'_{3}$,
or equivalently,
\[
w'_{(3)}\circ I = I'(w'_{(3)}).
\]
}
\end{nota}

\begin{rema}{\rm
Combining Proposition \ref{pz} with Proposition \ref{im:correspond},
we see that for any integer $p$, we also have a natural linear
isomorphism
\[
{\cal N}_{W'_3}^{(W_1 \otimes W_2)^{*}} \stackrel{\sim}{\longrightarrow}
{\cal V}^{W_3}_{W_1 W_2}
\]
{}from ${\cal N}_{W'_3}^{(W_1 \otimes W_2)^{*}}$ to the space of
logarithmic intertwining operators of type ${W_3\choose W_1\,W_2}$.
In particular, given a logarithmic intertwining operator ${\cal Y}$ of
type ${W_3\choose W_1\,W_2}$, the map
\[
I'_{{\cal Y},p}: W'_3\to (W_1\otimes W_2)^*
\]
defined by
\[
I'_{{\cal Y},p}(w'_{(3)})(w_{(1)}\otimes w_{(2)})=\langle
w'_{(3)},{\cal Y}(w_{(1)},e^{l_p(z)})w_{(2)}\rangle_{W_3}
\]
is $\tilde{A}$-compatible and intertwines both actions on both spaces.
}
\end{rema}

Recall that we have formulated the notions of $P(z)$-product and
$P(z)$-tensor product using $P(z)$-intertwining maps (Definitions
\ref{pz-product} and \ref{pz-tp}).  Since we now know that
$P(z)$-intertwining maps can be interpreted as in Proposition \ref{pz}
(and Notation \ref{scriptN}), we can easily reformulate the notions of
$P(z)$-product and $P(z)$-tensor product correspondingly:
\begin{propo}\label{productusingI'}
Let ${\cal C}_1$ be either of the categories ${\cal M}_{sg}$ or ${\cal
GM}_{sg}$, as in Definition \ref{pz-product}.  For $W_1, W_2\in
\ob{\cal C}_1$, a $P(z)$-product $(W_3;I_3)$ of $W_1$ and $W_2$
(recall Definition \ref{pz-product}) amounts to an object $(W_3,Y_3)$
of ${\cal C}_1$ equipped with a map $I'_3 \in {\cal N}_{W'_3}^{(W_1
\otimes W_2)^{*}}$, that is, equipped with an $\tilde{A}$-compatible
map
\[
I'_3:W'_{3} \rightarrow (W_{1}\otimes W_{2})^{*}
\]
that intertwines the two actions of $V \otimes \iota_{+}{\mathbb
C}[t,t^{- 1}, (z^{-1}-t)^{-1}]$ and of ${\mathfrak s} {\mathfrak
l}(2)$.  The map $I'_3$ corresponds to the $P(z)$-intertwining map
\[
I_3:W_{1}\otimes W_{2} \rightarrow \overline{W_{3}}
\]
as above:
\[
I'_3(w'_{(3)}) = w'_{(3)}\circ I_3
\]
for $w'_{(3)} \in W'_{3}$ (recall
\ref{IcorrespondstoJalternateform})).  Denoting this structure by
$(W_3,Y_3;I'_3)$ or simply by $(W_3;I'_3)$, let $(W_4;I'_4)$ be
another such structure.  Then a morphism of $P(z)$-products {}from $W_3$
to $W_4$ amounts to a module map $\eta: W_3 \to W_4$ such that the
diagram
\begin{center}
\begin{picture}(100,60)
\put(-2,0){$W'_4$}
\put(13,4){\vector(1,0){104}}
\put(119,0){$W'_3$}
\put(38,50){$(W_1\otimes W_2)^*$}
\put(13,12){\vector(3,2){50}}
\put(118,12){\vector(-3,2){50}}
\put(65,8){$\eta'$}
\put(23,27){$I'_4$}
\put(98,27){$I'_3$}
\end{picture}
\end{center}
commutes, where $\eta'$ is the natural map given by (\ref{fprime}).
\end{propo}
\pf All we need to check is that the diagram in Definition
\ref{pz-product} commutes if and only if the diagram above commutes.
But this follows {}from the definitions and the fact that for
$\overline{w_{(3)}} \in \overline{W_{3}}$ and $w'_{(4)} \in W'_{4}$,
\[
\langle \eta'(w'_{(4)}),\overline{w_{(3)}} \rangle=\langle
w'_{(4)},\overline{\eta}(\overline{w_{(3)}})\rangle,
\]
which in turn follows {}from (\ref{fprime}).  \epfv

\begin{corol}\label{tensorproductusingI'}
Let ${\cal C}$ be a full subcategory of either ${\cal M}_{sg}$ or
${\cal GM}_{sg}$, as in Definition \ref{pz-tp}.  For $W_1, W_2\in
\ob{\cal C}$, a $P(z)$-tensor product $(W_0; I_0)$ of $W_1$ and $W_2$
in ${\cal C}$, if it exists, amounts to an object $W_0 =
W_1\boxtimes_{P(z)} W_2$ of ${\cal C}$ and a structure $(W_0 =
W_1\boxtimes_{P(z)} W_2; I'_0)$ as in Proposition
\ref{productusingI'}, with
\[
I'_0: (W_1\boxtimes_{P(z)} W_2)' \longrightarrow (W_1\otimes W_2)^*
\]
in ${\cal N}_{(W_1\boxtimes_{P(z)} W_2)'}^{(W_1 \otimes W_2)^{*}}$,
such that for any such pair $(W; I')$ $(W\in \ob \mathcal{C})$, with
\[
I': W' \longrightarrow (W_1\otimes W_2)^*
\]
in ${\cal N}_{W'}^{(W_1 \otimes W_2)^{*}}$, there is a unique module
map
\[
\chi: W' \longrightarrow (W_1\boxtimes_{P(z)} W_2)'
\]
such that the diagram
\begin{center}
\begin{picture}(100,60)
\put(-2,0){$W'$}
\put(13,4){\vector(1,0){104}}
\put(119,0){$(W_1\boxtimes_{P(z)} W_2)'$}
\put(38,50){$(W_1\otimes W_2)^*$}
\put(13,12){\vector(3,2){50}}
\put(118,12){\vector(-3,2){50}}
\put(65,8){$\chi$}
\put(23,27){$I'$}
\put(98,27){$I'_0$}
\end{picture}
\end{center}
commutes.  Here $\chi = \eta'$, where $\eta$ is a correspondingly
unique module map
\[
\eta: W_1\boxtimes_{P(z)} W_2 \longrightarrow W.
\]
Also, the map $I_0'$, which is $\tilde{A}$-compatible and which
intertwines the two actions of $V \otimes \iota_{+}{\mathbb C}[t,t^{-
1}, (z^{-1}-t)^{-1}]$ and of ${\mathfrak s} {\mathfrak l}(2)$, is
related to the $P(z)$-intertwining map
\[
I_0 = \boxtimes_{P(z)}: W_1\otimes W_2 \longrightarrow 
\overline{W_1\boxtimes_{P(z)} W_2}
\]
by
\[
I_0'(w') = w' \circ \boxtimes_{P(z)}
\]
for $w' \in (W_1\boxtimes_{P(z)} W_2)'$, that is,
\[
I_0'(w')(w_{(1)}\otimes w_{(2)}) = \langle w',w_{(1)}\boxtimes_{P(z)}
w_{(2)} \rangle
\]
for $w_{(1)}\in W_{1}$ and $w_{(2)}\in W_{2}$, using the notation 
(\ref{boxtensorofelements}).  
\epf
\end{corol}

\begin{rema}\label{motivationofbackslash}{\rm
{}From Corollary \ref{tensorproductusingI'} we see that it is natural to
try to construct $W_1\boxtimes_{P(z)} W_2$, when it exists, as the
contragredient of a suitable natural substructure of $(W_1\otimes
W_2)^*$.  We shall now proceed to do this.  Under suitable assumptions, 
we shall in fact construct
a module-like structure
\[
W_1\hboxtr_{P(z)} W_2 \subset (W_1\otimes W_2)^*
\]
for $W_1,W_1 \in \ob \mathcal{C}$, and we will show that 
$W_1\hboxtr_{P(z)} W_2$ is an object of $\mathcal{C}$ if and only if 
$W_1\boxtimes_{P(z)} W_2$ exists in $\mathcal{C}$, in which case we
will have 
\[
W_1\boxtimes_{P(z)} W_2 = (W_1\hboxtr_{P(z)} W_2)'
\]
(observe the notation
\[
\boxtimes = \hboxtr',
\]
as in the special cases studied in \cite{tensor1}--\cite{tensor3}).
It is important to keep in mind that the space $W_1\hboxtr_{P(z)} W_2$
will depend on the category $\mathcal{C}$.}
\end{rema}

We formalize certain of the properties of the category $\mathcal{C}$ that we 
have been using, and some new ones, as follows:

\begin{assum}\label{assum-c}
Throughout the remainder of this work, we shall assume that
$\mathcal{C}$ is a full subcategory of the category $\mathcal{M}_{sg}$
or $\mathcal{G}\mathcal{M}_{sg}$ closed under the contragredient
functor (recall Notation \ref{MGM}; for now, we are not assuming that
$V\in \ob \mathcal{C}$). We shall also assume that $\mathcal{C}$ is
closed under taking finite direct sums.
\end{assum}

\begin{defi}
{\rm For $W_{1}, W_{2}\in \ob \mathcal{C}$, define the subset 
\[
W_{1}\hboxtr_{P(z)}W_{2}\subset (W_{1}\otimes W_{2})^{*}
\]
of $(W_{1}\otimes W_{2})^{*}$ to be the union of the images
\[
I'(W')\subset (W_{1}\otimes W_{2})^{*}
\]
as $(W; I)$ ranges through all the $P(z)$-products of $W_{1}$ and
$W_{2}$ with $W\in \ob \mathcal{C}$. Equivalently,
$W_{1}\hboxtr_{P(z)}W_{2}$ is the union of the images $I'(W')$ as $W$
(or $W'$) ranges through $\ob \mathcal{C}$ and $I'$ ranges through
$\mathcal{N}_{W'}^{(W_{1}\otimes W_{2})^{*}}$---the space of
$\tilde{A}$-compatible linear maps
\[
W'\to (W_{1}\otimes W_{2})^{*}
\]
intertwining the actions of both 
\[
V\otimes \iota_{+}\C[t, t^{-1}, (z^{-1}-t)^{-1}]
\]
and $\mathfrak{s}\mathfrak{l}(2)$ on both spaces.}
\end{defi}

\begin{rema}
{\rm Since $\mathcal{C}$ is closed under direct sums (Assumption 
\ref{assum-c}), it is clear that $W_{1}\hboxtr_{P(z)}W_{2}$ is
in fact a linear subspace of $(W_{1}\otimes W_{2})^{*}$, and in particular,
it can be defined alternatively as the sum of all the images $I'(W')$:
\begin{equation}\label{hboxtr-sum}
W_{1}\hboxtr_{P(z)}W_{2}=\sum I'(W') = \bigcup I'(W')\subset 
(W_{1}\otimes W_{2})^{*},
\end{equation}
where the sum and union both range over $W\in \ob
\mathcal{C}$, $I\in \mathcal{M}_{W_{1}W_{2}}^{W}$.}
\end{rema}

For any generalized $V$-modules $W_{1}$ and $W_{2}$,
using the operator $L'_{P(z)}(0)$ (recall (\ref{LP'(j)}))
on $(W_{1}\otimes W_{2})^{*}$ we define
the generalized $L'_{P(z)}(0)$-eigenspaces 
$((W_{1}\otimes W_{2})^{*})_{[n]}$ for $n\in \C$ in the usual way:
\begin{equation}
((W_{1}\otimes W_{2})^{*})_{[n]}=\{w\in (W_{1}\otimes W_{2})^{*}\;|\;
(L'_{P(z)}(0)-n)^{m}w=0 \;{\rm for}\; m\in \N \;
\mbox{\rm sufficiently large}\}.
\end{equation}
Then we have the (proper) subspace 
\begin{equation}
\coprod_{n\in \C}((W_{1}\otimes W_{2})^{*})_{[n]}\subset 
(W_{1}\otimes W_{2})^{*}.
\end{equation}
We also define the ordinary $L'_{P(z)}(0)$-eigenspaces
$((W_{1}\otimes W_{2})^{*})_{(n)}$ in the usual way:
\begin{equation}
((W_{1}\otimes W_{2})^{*})_{(n)}=\{w\in (W_{1}\otimes W_{2})^{*}\;|\;
L'_{P(z)}(0)w=nw \}.
\end{equation}
Then we have the (proper) subspace
\begin{equation}
\coprod_{n\in \C}((W_{1}\otimes W_{2})^{*})_{(n)}\subset 
(W_{1}\otimes W_{2})^{*}.
\end{equation}

\begin{propo}\label{im:abc}
Let $W_{1}, W_{2}\in \ob \mathcal{C}$. 

(a)
The elements of $W_1\hboxtr_{P(z)} W_2$ are exactly the linear functionals
on $W_{1}\otimes W_{2}$ of the form $w'\circ
I(\cdot\otimes \cdot)$ for some $P(z)$-intertwining map $I$ of type
${W\choose W_1\,W_2}$ and some $w'\in W'$, $W\in\ob{\cal C}$.

(b) Let $(W; I)$ be any $P(z)$-product of $W_{1}$ and $W_{2}$, with 
$W$ any generalized $V$-module. Then for $n\in \C$,
\[
I'(W'_{[n]}) \subset  ((W_{1}\otimes W_{2})^{*})_{[n]}
\]
and 
\[
I'(W'_{(n)})\subset ((W_{1}\otimes W_{2})^{*})_{(n)}.
\]

(c) The structure $(W_1\hboxtr_{P(z)} W_2,Y'_{P(z)})$ (recall
(\ref{y'-p-z})) satisfies all the axioms in the definition of
(strongly $\tilde{A}$-graded) generalized $V$-module except perhaps
for the two grading conditions (\ref{set:dmltc}) and
(\ref{set:dmfin}).

(d) Suppose that the objects of the category $\mathcal{C}$ consist
only of (strongly $\tilde{A}$-graded) {\em ordinary}, as opposed to
{\em generalized}, $V$-modules. Then the structure $(W_1\hboxtr_{P(z)}
W_2,Y'_{P(z)})$ satisfies all the axioms in the definition of
(strongly $\tilde{A}$-graded ordinary) $V$-module except perhaps for
(\ref{set:dmltc}) and (\ref{set:dmfin}).
\end{propo}
\pf 
Part (a) is clear {}from the definition of $W_1\hboxtr_{P(z)} W_2$, and 
(b) follows {}from (\ref{j-lj}) with $j=0$.

For (c), let $(W; I)$ be any any $P(z)$-product of $W_{1}$ and $W_{2}$, with 
$W$ any generalized $V$-module. Then $(I'(W'), Y'_{P(z)})$ satisfies all
the conditions in the definition of (strongly $\tilde{A}$-graded) generalized 
$V$-module since $I'$ is $\tilde{A}$-compatible and intertwines the actions of 
$V\otimes {\mathbb C}[t,t^{-1}]$ and of ${\mathfrak s}{\mathfrak
l}(2)$; the $\C$-grading follows {}from Part (b). Since $W_1\hboxtr_{P(z)} W_2$
is the sum of these structures $I'(W')$ over $W\in \ob \mathcal{C}$ (recall
(\ref{hboxtr-sum})), we see that $(W_1\hboxtr_{P(z)} W_2, Y'_{P(z)})$ satisfies
all the conditions in the definition of generalized module except perhaps for 
(\ref{set:dmltc}) and (\ref{set:dmfin}).

Finally, Part (d) is proved by the same argument as for (c). In fact, for 
$(W; I)$ any $P(z)$-product of possibly generalized $V$-modules $W_{1}$ and $W_{2}$,
with $W$ any ordinary $V$-module, $(I'(W'), Y'_{P(z)})$ satisfies all the conditions 
in the definition of (strongly $\tilde{A}$-graded) ordinary $V$-module; the 
$\C$-grading (this time, by ordinary $L'_{P(z)}(0)$-eigenspaces) 
again follows {}from Part (b).
\epfv

We now have the following generalization of Proposition 13.7 in
\cite{tensor3}, characterizing $W_{1}\boxtimes_{P(z)}W_{2}$, including its
existence, in terms of $W_1\hboxtr_{P(z)} W_2$:

\begin{propo}\label{tensor1-13.7}
Let $W_{1}, W_{2}\in \ob \mathcal{C}$. 
If $(W_1\hboxtr_{P(z)} W_2, Y'_{P(z)})$ is an object of ${\cal C}$,
denote by $(W_1\boxtimes_{P(z)} W_2, Y_{P(z)})$ its contragredient
module. Then the $P(z)$-tensor product of $W_{1}$ and $W_{2}$ 
in ${\cal C}$ exists and is
$(W_1\boxtimes_{P(z)} W_2, Y_{P(z)}; i')$, where $i$ is the natural
inclusion {}from $W_1\hboxtr_{P(z)} W_2$ to $(W_1\otimes W_2)^*$ (recall
Notation \ref{scriptN}).
Conversely, let us assume that $\mathcal{C}$ is closed under quotients. 
If the $P(z)$-tensor product of $W_1$ and $W_2$ in ${\cal
C}$ exists, then $(W_1\hboxtr_{P(z)} W_2, Y'_{P(z)})$ is an object of
${\cal C}$.
\end{propo}
\pf 
Suppose that $(W_1\hboxtr_{P(z)} W_2, Y'_{P(z)})$ is an object of
${\cal C}$ and take $(W_1\boxtimes_{P(z)} W_2, Y_{P(z)})$ and 
the map $i$ as indicated. Then 
\[
i\in \mathcal{N}_{W_1\hboxtr_{P(z)} W_2}^{W_{1}\otimes W_{2})^{*}},
\]
and
\[
i'\in \mathcal{M}_{W_{1}W_{2}}^{W_1\boxtimes_{P(z)} W_2}.
\]
In the notation of Corollary \ref{tensorproductusingI'}, we take
$I_{0}=i'$, $I_{0}'=i$. For any pair $(W; I')$ as in Corollary
\ref{tensorproductusingI'}, we have $I'(W')\subset W_1\hboxtr_{P(z)}
W_2$ (which is the union of all such images), so that there is
certainly a unique module map
\[
\chi: W'\to W_1\hboxtr_{P(z)} W_2
\]
such that 
\[
i\circ \chi=I',
\]
namely, $I'$ itself, viewed as a module map. Thus by Corollary
\ref{tensorproductusingI'}, $W_1\boxtimes_{P(z)} W_2$ exists as
indicated.

Conversely, if the $P(z)$-tensor product of $W_1$ and $W_2$ in ${\cal
C}$ exists and is $(W_0; I_0)$, then for any $P(z)$-product $(W; I)$
with $W \in \ob \mathcal{C}$, we have a unique module map $\chi: W'\to
W'_0$ as in Corollary \ref{tensorproductusingI'} such that
$I'=I'_0\circ \chi$, so that $I'(W')\subset I'_0(W'_0)$, proving that
$W_1\hboxtr_{P(z)} W_2\subset I'_0(W'_0)$. On the other hand, $(W_0;
I_0)$ is itself a $P(z)$-product, so that $I'_0(W'_0)\subset
W_1\hboxtr_{P(z)} W_2$. Thus $W_1\hboxtr_{P(z)} W_2=I'_0(W'_0)$, and
so $W_1\hboxtr_{P(z)} W_2$ is a generalized $V$-module and is the
image of the module map
\[
I_{0}': W_{0}'\to W_1\hboxtr_{P(z)} W_2.
\]
Since $\mathcal{C}$ is closed under quotients by assumption, we have that 
$W_1\hboxtr_{P(z)} W_2 \in \ob \mathcal{C}$.
\epfv

\begin{rema}
{\rm Suppose that $W_1\hboxtr_{P(z)} W_2$ is an object of
$\mathcal{C}$.  {}From Corollary \ref{tensorproductusingI'} and
Proposition \ref{tensor1-13.7} we see that
\begin{equation}\label{boxpair}
\langle\lambda, w_{(1)}\boxtimes_{P(z)}w_{(2)}\rangle
\lbar_{W_1\boxtimes_{P(z)} W_2}=
\lambda(w_{(1)}\otimes w_{(2)})
\end{equation}
for $\lambda\in W_1\hboxtr_{P(z)} W_2\subset (W_1\otimes W_2)^*$,
$w_{(1)}\in W_1$ and $w_{(2)}\in W_2$.}
\end{rema}

Our next goal is to present a crucial alternative description of the
subspace $W_1\hboxtr_{P(z)} W_2$ of $(W_1\otimes W_2)^*$. The main
ingredient of this description will be the ``$P(z)$-compatibility
condition,'' which was a cornerstone of the development of tensor
product theory in the special cases treated in
\cite{tensor1}--\cite{tensor3} and \cite{tensor4}.

Assume now that $W_{1}$ and $W_{2}$ are arbitrary generalized $V$-modules. 
Let $(W; I)$ ($W$ a generalized $V$-module) 
be a $P(z)$-product of $W_{1}$ and $W_{2}$ and let 
$w'\in W'$. Then {}from (\ref{j-tau}), Proposition \ref{productusingI'},
(\ref{tau-w-comp}), (\ref{3.7}) and (\ref{y'-p-z}), we have, for all
$v\in V$, 
\begin{eqnarray}\label{5.18-p}
\lefteqn{\tau_{P(z)}\left(x_0^{-1}\delta\bigg(\frac{x^{-1}_1-z}{x_0}
\bigg)
Y_{t}(v, x_1)\right)I'(w')}\nno\\
&&=I'\left(\tau_{W'}\left(x_0^{-1}\delta\bigg(\frac{x^{-1}_1-z}{x_0}
\bigg)
Y_{t}(v, x_1)\right)w'\right)\nno\\
&&=I'\left(x_0^{-1}\delta\bigg(\frac{x^{-1}_1-z}{x_0}
\bigg)Y_{W'}(v, x_1)w'\right)\nno\\
&&=x_0^{-1}\delta\bigg(\frac{x^{-1}_1-z}{x_0}
\bigg)I'(Y_{W'}(v,
x_1)w')\nno\\
&&=x_0^{-1}\delta\bigg(\frac{x^{-1}_1-z}{x_0}
\bigg)I'(\tau_{W'}(Y_{t}(v,
x_1))w')\nno\\
&&=x_0^{-1}\delta\bigg(\frac{x^{-1}_1-z}{x_0}
\bigg)\tau_{P(z)}(Y_{t}(v,
x_1))I'(w')\nn
&&=x_0^{-1}\delta\bigg(\frac{x^{-1}_1-z}{x_0}
\bigg)
Y'_{P(z)}(v, x_1)I'(w').
\end{eqnarray}
That is, 
$I'(w')$ satisfies the following nontrivial and subtle condition on
$\lambda \in (W_1\otimes W_2)^{*}$:

\begin{description}
\item{\bf The $P(z)$-compatibility condition}

(a) The {\em $P(z)$-lower truncation condition}: For all $v\in V$, the formal
Laurent series $Y'_{P(z)}(v, x)\lambda$ involves only finitely many
negative powers of $x$.

(b) The following formula holds:
\begin{eqnarray}\label{cpb}
\lefteqn{\tau_{P(z)}\bigg(x_0^{-1}\delta\bigg(\frac{x^{-1}_1-z}{x_0}
\bigg)
Y_{t}(v, x_1)\bigg)\lambda}\nno\\
&&=x_0^{-1}\delta\bigg(\frac{x^{-1}_1-z}{x_0}\bigg)
Y'_{P(z)}(v, x_1)\lambda  \;\;\mbox{ for all }\;v\in V.
\end{eqnarray}
(Note that the two sides of (\ref{cpb}) are not {\it a priori} equal
for general $\lambda\in (W_1\otimes W_2)^{*}$. Note also that Condition 
(a) insures that the right-hand side in Condition (b) is 
well defined.)
\end{description}

\begin{nota}
{\rm Note that the set of elements of $(W_1\otimes W_2)^*$ satisfying
either  the full $P(z)$-compatibility
condition or Part (a) of this condition forms a subspace. 
We shall denote the space of elements of $(W_1\otimes W_2)^*$ satisfying
the $P(z)$-compatibility
condition by
\[
\comp_{P(z)}((W_1\otimes W_2)^*).
\]}
\end{nota}

Recall {}from the comments preceding Definition
\ref{linearactioncompatible} that for each $\beta\in \tilde{A}$ we
have the subspace $((W_1\otimes W_2)^*)^{(\beta)}$ of $(W_1\otimes
W_2)^*$.  The sum of these subspaces is of course direct:
\[
\sum_{\beta\in \tilde{A}}((W_1\otimes W_2)^*)^{(\beta)}
=\coprod_{\beta\in \tilde{A}}((W_1\otimes W_2)^*)^{(\beta)}.
\]
Each space $((W_1\otimes W_2)^*)^{(\beta)}$ is $L'_{P(z)}(0)$-stable
(recall Proposition \ref{tau-a-comp} and Remark
\ref{L'jpreservesbetaspace}), so that we may consider the subspaces
\[
\coprod_{n\in \C}((W_1\otimes W_2)^*)_{[n]}^{(\beta)} \subset
((W_1\otimes W_2)^*)^{(\beta)}
\]
and 
\[
\coprod_{n\in \C}((W_1\otimes W_2)^*)_{(n)}^{(\beta)} \subset
((W_1\otimes W_2)^*)^{(\beta)}
\]
(recall Remark \ref{generalizedeigenspacedecomp}).  We now define the
two subspaces
\begin{equation}\label{W1W2_[C]^Atilde}
((W_1\otimes W_2)^*)_{[{\mathbb C}]}^{( \tilde A )}=
\coprod_{n\in
\C}\coprod_{\beta\in \tilde{A}}((W_1\otimes W_2)^*)_{[n]}^{(\beta)}
\subset (W_1\otimes W_2)^*
\end{equation}
and 
\begin{equation}\label{W1W2_(C)^Atilde}
((W_1\otimes W_2)^*)_{({\mathbb C})}^{( \tilde A )}=
\coprod_{n\in
\C}\coprod_{\beta\in \tilde{A}}((W_1\otimes W_2)^*)_{(n)}^{(\beta)}
\subset (W_1\otimes W_2)^*.
\end{equation}

\begin{rema}\label{singleanddoublegraded}
{\rm Any $L'_{P(z)}(0)$-stable subspace of $((W_1\otimes
W_2)^*)_{[{\mathbb C}]}^{( \tilde A )}$ is graded by generalized
eigenspaces (again recall Remark \ref{generalizedeigenspacedecomp}),
and if such a subspace is also $\tilde A$-graded, then it is doubly
graded; similarly for subspaces of $((W_1\otimes W_2)^*)_{({\mathbb
C})}^{( \tilde A )}$.}
\end{rema}

We have:

\begin{lemma}\label{a-tilde-comp}
Suppose that $\lambda\in ((W_1\otimes W_2)^*)_{[{\mathbb C}]}^{(\tilde
A )}$ satisfies the $P(z)$-compatibility condition. Then every
$\tilde{A}$-homogeneous component of $\lambda$ also satisfies this
condition.
\end{lemma}
\pf
When $v\in V$ is $\tilde{A}$-homogeneous, 
\[
\tau_{P(z)}\bigg(x_0^{-1}\delta\bigg(\frac{x^{-1}_1-z}{x_0} \bigg)
Y_{t}(v, x_1)\bigg)\;\;\mbox{ and }\;\;
x_0^{-1}\delta\bigg(\frac{x^{-1}_1-z}{x_0}\bigg) Y'_{P(z)}(v, x_1)
\]
are both $\tilde{A}$-homogeneous as operators, in the obvious sense.
By comparing the $\tilde{A}$-homogeneous components of both sides of
(\ref{cpb}), we see that the $\tilde{A}$-homogeneous components of
$\lambda$ also satisfy the $P(z)$-compatibility condition.  \epfv

\begin{rema}\label{stableundercomponentops}
{\rm Both the spaces $((W_1\otimes W_2)^*)_{[{\mathbb C}]}^{( \tilde A
)}$ and $((W_1\otimes W_2)^*)_{({\mathbb C})}^{( \tilde A )}$ are
stable under the component operators $\tau_{P(z)}(v\otimes t^m)$ of
the operators $Y'_{P(z)}(v,x)$ for $v\in V$, $m\in {\mathbb Z}$, and
under the operators $L'_{P(z)}(-1)$, $L'_{P(z)}(0)$ and
$L'_{P(z)}(1)$.  For the $\tilde A$-grading, this follows {}from
Proposition \ref{tau-a-comp} and Remark \ref{L'jpreservesbetaspace},
and for the $\C$-gradings, we simply follow the proof of Proposition
\ref{gweight}, using Propositions \ref{id-dev} and \ref{pz-comm}
together with (\ref{pz-sl-2-pz-y--1}).}
\end{rema}

Again let $(W; I)$ ($W$ a generalized $V$-module) 
be a $P(z)$-product of $W_{1}$ and $W_{2}$ and let 
$w'\in W'$. 
Since $I'$ in particular intertwines the actions of
$V\otimes{\mathbb C}[t, t^{-1}]$ and of $\mathfrak{s}\mathfrak{l}(2)$,
and is $\tilde{A}$-compatible, $I'(W')$ is a generalized $V$-module,
as we have seen in the proof of Proposition \ref{im:abc}.
Therefore, for each $w'\in W'$, $I'(w')$ also satisfies
the following condition on $\lambda \in (W_1\otimes W_2)^*$:
\begin{description}
\item{\bf The $P(z)$-local grading restriction condition}

(a) The {\em $P(z)$-grading condition}: $\lambda$ is a (finite) sum of
generalized eigenvectors  for the operator
$L'_{P(z)}(0)$ on $(W_1\otimes W_2)^*$ that 
are also homogeneous with respect to $\tilde A$, that is, 
\[
\lambda\in ((W_1\otimes W_2)^*)_{[{\mathbb C}]}^{( \tilde A )}.
\]
\label{homo}

(b) Let $W_{\lambda}$ be the smallest doubly graded (or equivalently,
$\tilde A$-graded; recall Remark \ref{singleanddoublegraded}) subspace
of $((W_1\otimes W_2)^*)_{[ {\mathbb C} ]}^{( \tilde A )}$ containing
$\lambda$ and stable under the component operators
$\tau_{P(z)}(v\otimes t^m)$ of the operators $Y'_{P(z)}(v,x)$ for
$v\in V$, $m\in {\mathbb Z}$, and under the operators $L'_{P(z)}(-1)$,
$L'_{P(z)}(0)$ and $L'_{P(z)}(1)$.  (In view of Remark
\ref{stableundercomponentops}, $W_{\lambda}$ indeed exists.)  Then
$W_{\lambda}$ has the properties
\begin{eqnarray}
&\dim(W_{\lambda})^{(\beta)}_{[n]}<\infty,&\label{lgrc1}\\
&(W_{\lambda})^{(\beta)}_{[n+k]}=0\;\;\mbox{ for }\;k\in {\mathbb Z}
\;\mbox{ sufficiently negative},&\label{lgrc2}
\end{eqnarray}
for any $n\in {\mathbb C}$ and $\beta\in \tilde A$, where as usual the
subscripts denote the ${\mathbb C}$-grading and the superscripts
denote the $\tilde A$-grading.
\end{description}

In the case that $W$ is an (ordinary) $V$-module and $w'\in W'$,
$I'(w')$ also satisfies the following $L(0)$-semisimple version of
this condition on $\lambda \in (W_1\otimes W_2)^*$:

\begin{description}
\item{\bf The $L(0)$-semisimple $P(z)$-local grading restriction condition}

(a) The {\em $L(0)$-semisimple 
$P(z)$-grading condition}: $\lambda$ is a (finite) sum of
eigenvectors  for the operator
$L'_{P(z)}(0)$ on $(W_1\otimes W_2)^*$ that 
are also homogeneous with respect to $\tilde A$, that is,
\[
\lambda\in ((W_1\otimes W_2)^*)_{({\mathbb C})}^{( \tilde A )}.
\]
\label{semi-homo}

(b) Consider $W_\lambda$ as above, which in this case is in fact the
smallest doubly graded (or equivalently, $\tilde A$-graded) subspace
of $((W_1\otimes W_2)^*)_{({\mathbb C})}^{( \tilde A )}$ containing
$\lambda$ and stable under the component operators
$\tau_{P(z)}(v\otimes t^m)$ of the operators $Y'_{P(z)}(v,x)$ for
$v\in V$, $m\in {\mathbb Z}$, and under the operators $L'_{P(z)}(-1)$,
$L'_{P(z)}(0)$ and $L'_{P(z)}(1)$.  Then $W_\lambda$ has the
properties
\begin{eqnarray}
&\dim(W_\lambda)^{(\beta)}_{(n)}<\infty,&\label{semi-lgrc1}\\
&(W_\lambda)^{(\beta)}_{(n+k)}=0\;\;\mbox{ for }\;k\in {\mathbb Z}
\;\mbox{ sufficiently negative},&\label{semi-lgrc2}
\end{eqnarray}
for any $n\in {\mathbb C}$ and $\beta\in \tilde A$, where  the
subscripts denote the ${\mathbb C}$-grading and the superscripts
denote the $\tilde A$-grading.

\end{description}

\begin{nota}
{\rm Note that the set of elements of $(W_1\otimes W_2)^*$ satisfying
either of these two $P(z)$-local grading restriction conditions, or
either of the Part (a)'s in these conditions, forms a subspace.  We
shall denote the space of elements of $(W_1\otimes W_2)^*$ satisfying
the $P(z)$-local grading restriction condition and the
$L(0)$-semisimple $P(z)$-local grading restriction condition by
\[
\lgr_{[\C]; P(z)}((W_1\otimes W_2)^*)
\]
and 
\[
\lgr_{(\C); P(z)}((W_1\otimes W_2)^*),
\]
respectively.}
\end{nota}

The following theorems are among the most important in this work.
Note that even in the finitely reductive case studied in
\cite{tensor3}, they are stronger and more general than (the last
assertion of) Theorem 13.9 in \cite{tensor3}.  The proofs of these two
theorems will be given in the next section.

\begin{theo}\label{comp=>jcb}
Let $\lambda$ be an element of $(W_{1}\otimes W_{2})^{*}$ satisfying
the $P(z)$-compatibility condition. Then when acting on $\lambda$, the
Jacobi identity for $Y'_{P(z)}$ holds, that is,
\begin{eqnarray}
\lefteqn{x_{0}^{-1}\delta
\left({\displaystyle\frac{x_{1}-x_{2}}{x_{0}}}\right)Y'_{P(z)}(u, x_{1})
Y'_{P(z)}(v, x_{2})\lambda}\nno\\
&&\hspace{2ex}-x_{0}^{-1} \delta
\left({\displaystyle\frac{x_{2}-x_{1}}{-x_{0}}}\right)Y'_{P(z)}(v, x_{2})
Y'_{P(z)}(u, x_{1})\lambda\nonumber \\
&&=x_{2}^{-1} \delta
\left({\displaystyle\frac{x_{1}-x_{0}}{x_{2}}}\right)Y'_{P(z)}(Y(u, x_{0})v,
x_{2})\lambda\label{cjcb}
\end{eqnarray}
for $u, v\in V$.
\end{theo}

\begin{theo}\label{stable}
The subspace $\comp_{P(z)}((W_1\otimes W_2)^*)$ of $(W_{1}\otimes
W_{2})^{*}$ is stable under the operators $\tau_{P(z)}(v\otimes
t^{n})$ for $v\in V$ and $n\in {\mathbb Z}$, and in the M\"obius case,
also under the operators $L'_{P(z)}(-1)$, $L'_{P(z)}(0)$ and
$L'_{P(z)}(1)$; similarly for the subspaces $\lgr_{[\C];
P(z)}((W_1\otimes W_2)^*$ and $\lgr_{(\C); P(z)}((W_1\otimes W_2)^*$.
\end{theo}

\begin{rema}{\rm
The converse of Theorem \ref{comp=>jcb} is not true. One can see this
in the tensor product theory of the ``trivial'' case where $V$ is a
vertex operator algebra associated with a finite-dimensional unital
commutative associative algebra $(A,\cdot,1)$ with derivation $D=0$
(cf.\ Remark \ref{va>cva}). In this case, $(V,Y,{\bf 1},\omega)=
(A,\cdot,1,0)$ and the Jacobi identity for a $V$-module $W$ reduces to
\begin{eqnarray*}
&{\dps x^{-1}_0\delta \bigg( {x_1-x_2\over x_0}\bigg) u\cdot (v\cdot
w) - x^{-1}_0\delta \bigg( {x_2-x_1\over -x_0}\bigg) v\cdot (u\cdot
w) }&\nno\\
&{\dps = x^{-1}_2\delta \bigg( {x_1-x_0\over x_2}\bigg)
(u\cdot v)\cdot w}
\end{eqnarray*}
for $u,v\in A$ and $w\in W$, where we also use ``$\cdot$'' to denote
the action of $A$ on its modules.  In particular, a $V$-module is just
a finite-dimensional module for the associative algebra $A$. Given
$V$-modules $W_1$ and $W_2$, the action $Y'_{P(z)}$ given in
(\ref{Y'def}) now becomes
\[
(Y'_{P(z)}(v,x)\lambda)(w_{(1)}\otimes w_{(2)})=
\lambda(w_{(1)}\otimes v\cdot w_{(2)}).
\]
{}From this it is clear that (\ref{cjcb}) holds for any element
$\lambda\in (W_1\otimes W_2)^*$. However, the $P(z)$-compatibility
condition (\ref{cpb}) in this case reduces to
\[
\lambda(v\cdot w_{(1)}\otimes w_{(2)})=
\lambda(w_{(1)}\otimes v\cdot w_{(2)})
\]
for all $v\in A$, $w_{(1)}\in W_1$ and $w_{(2)}\in W_2$, which is not
necessarily true for every $\lambda$.  This example is discussed
further in Remark 2.20 of \cite{HLLZ}, which treats a range of issues
related to the compatibility condition, intertwining operators, and
tensor product theory.}
\end{rema}

We now generalize the notion of ``weak module'' for a vertex operator 
algebra to our M\"obius or
conformal vertex algebra $V$:

\begin{defi}
{\rm A {\it weak module for $V$} (or {\it weak $V$-module}) is a
vector space $W$ equipped with a vertex operator map $Y_{W}: V\otimes
W\to W[[x,x^{-1}]]$ satisfying (only) the axioms (\ref{ltc-w}),
(\ref{m-1left}), (\ref{m-Jacobi}) and (\ref{L-1}) in Definition
\ref{cvamodule} (note that there is no grading given on $W$) and in
case $V$ is M\"obius, also the existence of a representation of
$\mathfrak{s}\mathfrak{l}(2)$ on $W$, as in Definition
\ref{moduleMobius}, satisfying the conditions
(\ref{sl2-1})--(\ref{sl2-3}).  }
\end{defi}

Then
we have:

\begin{theo}\label{wk-mod}
The space $\comp_{P(z)}((W_1\otimes W_2)^*)$, equipped with the vertex
operator map $Y'_{P(z)}$ and, in case $V$ is M\"{o}bius, also equipped
with the operators $L_{P(z)}'(-1)$, $L_{P(z)}'(0)$ and $L_{P(z)}'(1)$,
is a weak $V$-module; similarly for the spaces
\begin{equation}\label{COMPintLGR[]}
(\comp_{P(z)}((W_1\otimes W_2)^*))\cap 
(\lgr_{[\C]; P(z)}((W_1\otimes W_2)^*))
\end{equation}
and 
\begin{equation}\label{COMPintLGR()}
(\comp_{P(z)}((W_1\otimes W_2)^*))\cap 
(\lgr_{(\C); P(z)}((W_1\otimes W_2)^*)).
\end{equation}
\end{theo}
\pf By Theorem \ref{stable}, $Y'_{P(z)}$ is a map {}from the tensor
product of $V$ with any of these three subspaces to the space of
formal Laurent series with elements of the subspace as coefficients.
By Proposition \ref{id-dev} and Theorem \ref{comp=>jcb} and, in the
case that $V$ is a M\"{o}bius vertex algebra, also by Propositions
\ref{sl-2} and \ref{pz-l-y-comm}, we see that all the axioms for weak
$V$-module are satisfied.  \epfv

We also have:

\begin{theo}\label{generation}
Let 
\[
\lambda\in 
\comp_{P(z)}((W_1\otimes W_2)^*)\cap \lgr_{[\C]; P(z)}((W_1\otimes W_2)^*).
\] 
Then $W_{\lambda}$ (recall Part (b) of the $P(z)$-local grading
restriction condition) equipped with the vertex operator map
$Y_{P(z)}'$ and, in case $V$ is M\"{o}bius, also equipped with the
operators $L'_{P(z)}(-1)$, $L'_{P(z)}(0)$ and $L'_{P(z)}(1)$, is a
(strongly-graded) generalized $V$-module. If in addition
\[
\lambda\in 
\comp_{P(z)}((W_1\otimes W_2)^*)\cap \lgr_{(\C); P(z)}((W_1\otimes W_2)^*),
\] 
that is, $\lambda$ is a sum of eigenvectors of $L_{P(z)}'(0)$, then
$W_{\lambda}$ ($\subset ((W_{1}\otimes
W_{2})^{*})_{(\C)}^{(\tilde{A})}$) is a (strongly-graded) $V$-module.
\end{theo}
\pf Decompose $\lambda$ as
\[
\lambda = \sum_{\beta\in \tilde{A}}\lambda^{(\beta)}
\]
(finite sum), where $\lambda^{(\beta)} \in ((W_1\otimes
W_2)^*)^{(\beta)}$.  By Lemma \ref{a-tilde-comp}, each
$\lambda^{(\beta)}$ satisfies the $P(z)$-compatibility condition.
Also, each $\lambda^{(\beta)}$ satisfies the $P(z)$-grading condition
(and in the semisimple case, the $L(0)$-semisimple $P(z)$-grading
condition), and each $W_{{\lambda}^{(\beta)}}$ is simply the smallest
subspace containing $\lambda^{(\beta)}$ and stable under the operators
listed above (without the $\tilde A$-gradedness condition).  Moreover,
each $W_{{\lambda}^{(\beta)}} \subset W_{\lambda}$ and in fact
\[
W_{\lambda} = \sum_{\beta\in \tilde{A}}W_{{\lambda}^{(\beta)}}.
\] 
Thus each $\lambda^{(\beta)}$ lies in the space (\ref{COMPintLGR[]})
(or (\ref{COMPintLGR()})).  (Note that we have reduced Theorem
\ref{generation} to the $\tilde A$-homogeneous case.)  By Theorem
\ref{wk-mod}, each $W_{{\lambda}^{(\beta)}}$ is a weak submodule of
the weak module (\ref{COMPintLGR[]}) (or (\ref{COMPintLGR()})), and
hence is a (strongly-graded) generalized module (or module).  Thus
$W_{\lambda}$ has the same properties.
\epfv

Now we can give an alternative description of $W_1\hboxtr_{P(z)} W_2$
by characterizing the elements of $W_1\hboxtr_{P(z)} W_2$ using the
$P(z)$-compatibility condition and the $P(z)$-local grading
restriction conditions, generalizing Theorem 13.10 in
\cite{tensor3}. This description will be crucial in later sections,
especially in the construction of the associativity isomorphisms.

\begin{theo}\label{characterizationofbackslash}
Suppose that for every element 
\[
\lambda
\in \comp_{P(z)}((W_1\otimes W_2)^*)\cap \lgr_{[\C]; P(z)}((W_1\otimes W_2)^*)
\]
the (strongly-graded) generalized module $W_{\lambda}$ given in
Theorem \ref{generation} is an object of ${\cal C}$ (this of course
holds in particular if ${\cal C}$ is ${\cal GM}_{sg}$). Then
\[
W_1\hboxtr_{P(z)}W_2=\comp_{P(z)}((W_1\otimes W_2)^*)
\cap \lgr_{[\C]; P(z)}((W_1\otimes W_2)^*).
\]
Suppose that ${\cal C}$ is a category of strongly-graded $V$-modules
(that is, ${\cal C}\subset {\cal M}_{sg}$) and that for every element
\[
\lambda
\in \comp_{P(z)}((W_1\otimes W_2)^*)\cap \lgr_{(\C); P(z)}((W_1\otimes W_2)^*)
\]
the (strongly-graded) $V$-module $W_{\lambda}$ given in Theorem
\ref{generation} is an object of ${\cal C}$ (which of course holds in
particular if ${\cal C}$ is ${\cal M}_{sg}$).  Then
\[
W_1\hboxtr_{P(z)}W_2=\comp_{P(z)}((W_1\otimes W_2)^*)
\cap \lgr_{(\C); P(z)}((W_1\otimes W_2)^*).
\]
\end{theo}
\pf 
We have seen that
\[
W_1\hboxtr_{P(z)}W_2\subset 
\comp_{P(z)}((W_1\otimes W_2)^*)\cap \lgr_{[\C]; P(z)}((W_1\otimes W_2)^*)
\] 
and, in case ${\cal C}\subset {\cal M}_{sg}$,
\[
W_1\hboxtr_{P(z)}W_2\subset 
\comp_{P(z)}((W_1\otimes W_2)^*)\cap \lgr_{(\C); P(z)}((W_1\otimes W_2)^*).
\] 
On the other hand, by the assumptions, every element $\lambda$ of
\[
\comp_{P(z)}((W_1\otimes W_2)^*)\cap \lgr_{[\C]; P(z)}((W_1\otimes W_2)^*)
\] 
and, in case ${\cal C}\subset {\cal M}_{sg}$, every element $\lambda$
of
\[
\comp_{P(z)}((W_1\otimes W_2)^*)\cap \lgr_{(\C); P(z)}((W_1\otimes W_2)^*),
\] 
is contained in some object of ${\cal C}$, namely, $W_{\lambda}$, and
for any such (generalized) module, the inclusion map into $(W_1\otimes
W_2)^*$ satisfies the intertwining conditions in Proposition \ref{pz}.
Thus $\lambda$ lies in $W_1\hboxtr_{P(z)}W_2$, proving the desired
inclusion.  \epfv

\subsection{Constructions of $Q(z)$-tensor products}

We now give the construction of $Q(z)$-tensor products.  It is
analogous to that of $P(z)$-tensor products, and the formulations,
results and proofs in this subsection largely parallel those in
Subsection 5.2.  As usual, $z\in \C^{\times}$.

Given generalized $V$-modules $W_1$ and $W_2$, we shall be
constructing an action of the space $V \otimes \iota_{+}{\mathbb
C}[t,t^{- 1},(z+t)^{-1}]$ on the space $(W_1 \otimes W_2)^*$.

Let $I$ be a $Q(z)$-intertwining map of type ${W_3\choose W_1\, W_2}$,
as in Definition \ref{im:qimdef}. Consider the contragredient
generalized $V$-module $(W_{3}', Y_{3}')$, recall the opposite vertex
operator (\ref{yo}) and formula (\ref{y'}), and recall why the
ingredients of formula (\ref{imq:def}) are well defined. For $v\in V$,
$w_{(1)}\in W_{1}$, $w_{(2)}\in W_{2}$ and $w'_{(3)}\in W'_3$,
applying $w'_{(3)}$ to (\ref{imq:def}) we obtain
\begin{eqnarray}\label{imq:def'}
\lefteqn{\left\langle z^{-1}\delta\bigg(\frac{x_1-x_0}{z}\bigg)
Y_{3}'(v, x_0)w'_{(3)}, I(w_{(1)}\otimes w_{(2)})\right\rangle}\nno\\
&&=\left\langle w'_{(3)},x_{0}^{-1}\delta\bigg(\frac{x_1-z}{x_0}\bigg)
I(Y_1^{o}(v, x_1)w_{(1)}\otimes
w_{(2)})\right\rangle\nno\\
&&\quad -\left\langle w'_{(3)}, x^{-1}_0\delta\bigg(\frac{z-x_1}{-x_0}\bigg)
I(w_{(1)}\otimes Y_2(v, x_1)w_{(2)})\right\rangle.
\end{eqnarray}
We shall use this to motivate our action.

As we discussed in Subsection 5.1 (see (\ref{3.12}) and
(\ref{3.13})), in the left-hand side of (\ref{imq:def'}), the
coefficients of
\begin{equation}\label{qdeltaY3'}
z^{-1}\delta\bigg(\frac{x_1-x_0}{z}\bigg) Y_{3}'(v, x_0)
\end{equation}
in powers of $x_0$ and $x_1$, for all $v\in V$, span
\begin{equation}\label{qtausubW3'}
\tau_{W'_3}(V
\otimes \iota_{+}{\mathbb C}[t,t^{- 1},(z+t)^{-1}])
\end{equation}
(recall (\ref{tauw}) and (\ref{3.7})).  We now define a linear action
of $V \otimes \iota_{+}{\mathbb C}[t,t^{- 1},(z+t)^{-1}]$ on $(W_1
\otimes W_2)^*$, that is, a linear map $$\tau_{Q(z)}: V\otimes
\iota_{+}{\mathbb C}[t, t^{-1}, (z+t)^{-1}]\to{\rm End}\;(W_1\otimes
W_2)^{*}.$$ Recall the notations $T_{-z}^{+}$ and $T_{-z}^{o}$ {}from
Subsection 5.1 ((\ref{Tpm-z}) and (\ref{To-z})).

\begin{defi}
{\rm We define the linear action $\tau_{Q(z)}$ of 
\[
V
\otimes \iota_{+}{\mathbb C}[t,t^{- 1},(z+t)^{-1}]
\]
on $(W_1 \otimes
W_2)^*$
by
\begin{equation}\label{(5.1)}
(\tau_{Q(z)}(\xi)\lambda)(w_{(1)}\otimes w_{(2)})
=\lambda(\tau_{W_1}(T_{-z}^{o}\xi)w_{(1)}\otimes w_{(2)})
-\lambda(w_{(1)}\otimes \tau_{W_2}(T_{-z}^{+}\xi)w_{(2)})
\end{equation}
for $\xi\in V\otimes \iota_{+}{\mathbb C}[t, t^{-1}, (z+t)^{-1}]$,
$\lambda\in (W_1\otimes W_2)^{*}$, $w_{(1)}\in W_1$, $w_{(2)}\in W_2$,
and denote by $Y'_{Q(z)}$ the action of $V\otimes{\mathbb C}[t,t^{-1}]$
on $(W_1\otimes W_2)^*$ thus defined, that is,
\begin{equation}\label{y'-q-z}
Y'_{Q(z)}(v, x)=\tau_{Q(z)}(Y_{t}(v, x))
\end{equation}
for $v\in V$.}
\end{defi}

Using Lemma \ref{lemma5.2}, (\ref{3.7}) and (\ref{tauw-yto}),
we see
that (\ref{(5.1)}) can be written using generating
functions as
\begin{eqnarray}\label{5.2}
\lefteqn{\left(\tau_{Q(z)}
\left(z^{-1}\delta\left(\frac{x_1-x_0}{z}\right) Y_{t}(v,
x_0)\right)\lambda\right)(w_{(1)}\otimes w_{(2)})}\nno\\
&&=x^{-1}_0\delta\left(\frac{x_1-z}{x_0}\right) \lambda(Y_1^{o}(v,
x_1)w_{(1)}\otimes w_{(2)})\nno\\
&&\hspace{2em}-x_0^{-1}\delta\left(\frac{z-x_1}{-x_0}\right)
\lambda(w_{(1)}\otimes Y_2(v, x_1)w_{(2)})
\end{eqnarray}
for $v\in V$, $\lambda\in (W_1\otimes W_2)^{*}$, $w_{(1)}\in W_1$,
$w_{(2)}\in W_2$; compare this with (\ref{imq:def'}).  The generating
function form of the action $Y'_{Q(z)}$ can be obtained by taking
$\res_{x_1}$ of both sides of (\ref{5.2}):
\begin{eqnarray}\label{Y'qdef}
\lefteqn{(Y'_{Q(z)}(v,x_0)\lambda)(w_{(1)} \otimes
w_{(2)})}\nno\\
&&= \res_{x_1} x^{-1}_0 \delta\left(\frac{x_1-z}{x_0}\right)
\lambda(Y^o_1(v,x_1)w_{(1)} \otimes w_{(2)})\nno\\
&&\quad-
\res_{x_1}x^{-1}_0 \delta \left(\frac{z-x_1}{-x_0}\right)\lambda(w_{(1)}
\otimes Y_2(v,x_1)w_{(2)})\nno\\
&&= \lambda(Y^o_1(v,x_0 + z)w_{(1)}
\otimes w_{(2)}) \nno\\
&&\quad -
\res_{x_1}x^{-1}_0 \delta \left(\frac{z-x_1}{-x_0}\right)\lambda(w_{(1)}
\otimes Y_2(v,x_1)w_{(2)}).
\end{eqnarray}

\begin{rema}\label{I-intw-q}{\rm
Using the actions $\tau_{W'_3}$ and $\tau_{Q(z)}$, we can write
(\ref{imq:def'}) as
\[
\left(z^{-1}\delta\left(\frac{x_1-x_0}{z}\right)
Y_{3}'(v, x_0)w'_{(3)}\right)\circ I=
\tau_{Q(z)}\left(z^{-1}\delta\left(\frac{x_1-x_0}{z}\right)
Y_t(v, x_0)\right)(w'_{(3)}\circ I)
\]
or equivalently, as
\[
\left(\tau_{W'_3}\left(z^{-1}\delta\left(\frac{x_1-x_0}{z}\right)
Y_t(v, x_0)\right)w'_{(3)}\right)\circ I=
\tau_{Q(z)}\left(z^{-1}\delta\left(\frac{x_1-x_0}{z}\right)
Y_t(v, x_0)\right)(w'_{(3)}\circ I).
\]
}
\end{rema}

Recall the $\tilde{A}$-grading on $(W_1\otimes W_2)^{*}$ and the $A$-grading 
on $V
\otimes \iota_{+}{\mathbb C}[t,t^{- 1},(z^{-1}-t)^{-1}]$. 
Similarly, we also have an $A$-grading 
on $V
\otimes \iota_{+}{\mathbb C}[t,t^{- 1},(z+t)^{-1}]$. 
Definition \ref{linearactioncompatible} also 
applies to a linear action of 
$V \otimes \iota_{+}{\mathbb C}[t,t^{- 1}, (z+t)^{-1}]$
on $(W_1 \otimes W_2)^*$.
{}From 
(\ref{(5.1)}) or (\ref{5.2}), we have:

\begin{propo}\label{tau-q-a-comp}
The action $\tau_{Q(z)}$ is $\tilde{A}$-compatible. \epf
\end{propo}

We also have:

\begin{propo}\label{5.1}
The action $Y'_{Q(z)}$ has the property
\begin{equation}\label{Q-id}
Y'_{Q(z)}({\bf 1}, x)=1
\end{equation}
and the $L(-1)$-derivative property
\begin{equation}\label{QL-1}
\frac{d}{dx}Y'_{Q(z)}(v, x)=Y'_{Q(z)}(L(-1)v, x)
\end{equation}
for $v\in V$. 
\end{propo}
\pf
{From} (\ref{Y'qdef}), (\ref{yo}) and (\ref{3termdeltarelation}), 
\begin{eqnarray}
(Y'_{Q(z)}({\bf 1},x)\lambda)(w_{(1)} \otimes
w_{(2)}) 
&=& \res_{x_1}x^{-1} \delta\left(\frac{x_1-
z}{x}\right)\lambda(w_{(1)} \otimes w_{(2)})\nno\\
&&-\res_{x_1}x^{-1}
\delta\left(\frac{z-x_1}{-x}\right)\lambda(w_{(1)} \otimes w_{(2)})\nno\\
&=&
\res_{x_1}x^{-1}_1 \delta\left(\frac{z+x}{x_1}\right)\lambda(w_{(1)}
\otimes w_{(2)})\nno\\
&=& \lambda(w_{(1)} \otimes
w_{(2)}),
\end{eqnarray}
proving (\ref{Q-id}).  To prove the $L(-1)$-derivative property,
observe that {from} (\ref{Y'qdef}),
\begin{eqnarray}\label{5.8}
\lefteqn{\left(\left(\frac{d}{dx}
Y'_{Q(z)}(v,x)\right)\lambda\right)(w_{(1)} \otimes w_{(2)})}\nno\\
&&= \frac{d}{dx}
\lambda(Y^o_1(v,x + z)w_{(1)} \otimes w_{(2)})\nno\\
&&\quad
+\res_{x_1}\left(\frac{d}{dx} z^{-1}\delta\left(\frac{-x + x_1}{z}\right)
\right)\lambda(w_{(1)}
\otimes Y_2(v,x_1)w_{(2)}).
\end{eqnarray} 
But for any formal Laurent series $f(x)$, we have
\begin{equation}
\frac{d}{dx}f\left(\frac{-x+x_{1}}{z}\right)
=-\frac{d}{dx_{1}}f\left(\frac{-x+x_{1}}{z}\right)
\end{equation}
and if $f(x)$ involves only finitely many negative powers of $x$,
\begin{equation}
\res_{x_{1}}\left(\frac{d}{dx_{1}}z^{-1}\delta\left(
\frac{-x+x_{1}}{z}\right)\right)f(x_{1})=
-\res_{x_{1}}z^{-1}\delta\left(
\frac{-x+x_{1}}{z}\right)\frac{d}{dx_{1}}f(x_{1})
\end{equation}
(since the residue of a derivative is $0$).
We also have the $L(-1)$-derivative property (\ref{yo-l-1}) 
for $Y^o$.
Thus the right-hand side of (\ref{5.8}) equals
\begin{eqnarray}
\lefteqn{\lambda(Y^o_1(L(-
1)v,x+z)w_{(1)} \otimes w_{(2)})}\nno\\
&&\quad +
\res_{x_1}z^{-1}\delta\left(
\frac{-x+x_{1}}{z}\right)\frac{d}{dx_1}\lambda(w_{(1)} \otimes
Y_2(v,x_1)w_{(2)})\nno\\
&&=  \lambda(Y^o_1(L(-1)v,x+z)w_{(1)}
\otimes w_{(2)})\nno\\
&&\quad + \res_{x_1}z^{-1}\delta\left(
\frac{-x+x_{1}}{z}\right)\lambda(w_{(1)}
\otimes Y_2(L(-1)v,x_1)w_{(2)})\nno\\
&&=  (Y'_{Q(z)}(L(-1)v,x)\lambda)(w_{(1)}
\otimes w_{(2)}),
\end{eqnarray}
proving (\ref{QL-1}).  \epfv

\begin{propo}\label{qz-comm}
The action $Y'_{Q(z)}$ satisfies the commutator formula for vertex
operators: On $(W_1\otimes W_2)^{*}$,
\begin{eqnarray}\label{commu-q-z}
\lefteqn{[Y'_{Q(z)}(v_1, x_1), Y'_{Q(z)}(v_2, x_2)]}\nn
&&=\res_{x_0}x_2^{-1}\delta\left(\frac{x_1-x_0}{x_2}\right)
Y'_{Q(z)}(Y(v_1, x_0)v_2, x_2)
\end{eqnarray}
for $v_1, v_2\in V$.
\end{propo}
\pf As usual, the reader should note the well-definedness of each
expression and the justifiability of each use of a $\delta$-function
property in the argument that follows.  This argument is the same as
the proof of Proposition 5.2 of \cite{tensor1}, given in Section 8 of
\cite{tensor2}.  Let $\lambda \in (W_{1}\otimes W_{2})^{*}$, $v_{1},
v_{2}\in V$, $w_{(1)}\in W_{1}$ and $w_{(2)}\in W_{2}$. By
(\ref{Y'qdef}),
\bea\label{8.1}
\lefteqn{(Y'_{Q(z)}(v_{1}, x_{1})Y'_{Q(z)}(v_{2}, x_{2})
\lambda)(w_{(1)}\otimes w_{(2)})}\nno\\
&&=\mbox{\rm Res}_{y_{1}}x_{1}^{-1}\delta
\left(\frac{y_{1}-z}{x_{1}}\right)(Y'_{Q(z)}(v_{2}, x_{2})
\lambda)(Y_{1}^{o}(v_{1}, y_{1})w_{(1)}\otimes w_{(2)})\nno\\
&&\quad -\mbox{\rm Res}_{y_{1}}x_{1}^{-1}\delta
\left(\frac{z-y_{1}}{-x_{1}}\right)(Y'_{Q(z)}(v_{2}, x_{2})
\lambda)(w_{(1)}\otimes Y_{2}(v_{1}, y_{1})w_{(2)})\nno\\
&&=\mbox{\rm Res}_{y_{1}}\mbox{\rm Res}_{y_{2}}x_{1}^{-1}\delta
\left(\frac{y_{1}-z}{x_{1}}\right)x_{2}^{-1}\delta
\left(\frac{y_{2}-z}{x_{2}}\right)\cdot \nno\\
&&\hspace{4em}\cdot \lambda(Y_{1}^{o}(v_{2}, y_{2})Y_{1}^{o}
(v_{1}, y_{1})w_{(1)}\otimes w_{(2)})
\nno\\
&&\quad -\mbox{\rm Res}_{y_{1}}\mbox{\rm Res}_{y_{2}}x_{1}^{-1}\delta
\left(\frac{y_{1}-z}{x_{1}}\right)x_{2}^{-1}\delta
\left(\frac{z-y_{2}}{-x_{2}}\right)\cdot \nno\\
&&\hspace{4em}\cdot \lambda(Y_{1}^{o}(v_{1}, y_{1})w_{(1)}
\otimes Y_{2}(v_{2}, y_{2})w_{(2)})\nno\\
&&\quad -\mbox{\rm Res}_{y_{1}}\mbox{\rm Res}_{y_{2}}x_{1}^{-1}\delta
\left(\frac{z-y_{1}}{-x_{1}}\right)x_{2}^{-1}\delta
\left(\frac{y_{2}-z}{x_{2}}\right)\cdot \nno\\
&&\hspace{4em} \cdot \lambda(Y_{1}^{o}(v_{2}, y_{2})w_{(1)}
\otimes Y_{2}(v_{1}, y_{1})w_{(2)})\nno\\
&&\quad +\mbox{\rm Res}_{y_{1}}\mbox{\rm Res}_{y_{2}}x_{1}^{-1}\delta
\left(\frac{z-y_{1}}{-x_{1}}\right)x_{2}^{-1}\delta
\left(\frac{z-y_{2}}{-x_{2}}\right)\cdot \nno\\
&&\hspace{4em}\cdot \lambda(w_{(1)}
\otimes Y_{2}(v_{2}, y_{2})Y_{2}(v_{1}, y_{1})w_{(2)}).
\eea
Transposing the subscripts $1$ and $2$ of the symbols $v$, $x$ and $y$, 
we have
\bea\label{8.2}
\lefteqn{(Y'_{Q(z)}(v_{2}, x_{2})Y'_{Q(z)}(v_{1}, x_{1})\lambda)(w_{(1)}
\otimes w_{(2)})}\nno\\
&&=\mbox{\rm Res}_{y_{2}}\mbox{\rm Res}_{y_{1}}x_{2}^{-1}
\delta\left(\frac{y_{2}-z}{x_{2}}\right)x_{1}^{-1}
\delta\left(\frac{y_{1}-z}{x_{1}}\right)\cdot \nno\\
&&\hspace{4em}\cdot \lambda(Y_{1}^{o}
(v_{1}, y_{1})Y_{1}^{o}(v_{2}, y_{2})w_{(1)}\otimes w_{(2)})\nno\\
&&\quad -\mbox{\rm Res}_{y_{2}}\mbox{\rm Res}_{y_{1}}x_{2}^{-1}
\delta\left(\frac{y_{2}-z}{x_{2}}\right)x_{1}^{-1}
\delta\left(\frac{z-y_{1}}{-x_{1}}\right)\cdot \nno\\
&&\hspace{4em}\cdot \lambda(Y_{1}^{o}(v_{2}, y_{2})w_{(1)}
\otimes Y_{2}(v_{1}, y_{1})w_{(2)})\nno\\
&&\quad -\mbox{\rm Res}_{y_{2}}\mbox{\rm Res}_{y_{1}}x_{2}^{-1}
\delta\left(\frac{z-y_{2}}{-x_{2}}\right)x_{1}^{-1}\delta
\left(\frac{y_{1}-z}{x_{1}}\right)\cdot \nno\\
&&\hspace{4em}\cdot \lambda(Y_{1}^{o}(v_{1}, y_{1})w_{(1)}
\otimes Y_{2}(v_{2}, y_{2})w_{(2)})\nno\\
&&\quad +\mbox{\rm Res}_{y_{2}}\mbox{\rm Res}_{y_{1}}x_{2}^{-1}\delta
\left(\frac{z-y_{2}}{-x_{2}}\right)x_{1}^{-1}\delta\left
(\frac{z-y_{1}}{-x_{1}}\right)\cdot \nno\\
&&\hspace{4em}\cdot \lambda(w_{(1)}\otimes Y_{2}(v_{1}, y_{1})
Y_{2}(v_{2}, y_{2})w_{(2)}).
\eea
Formulas (\ref{8.1}) and (\ref{8.2}) give
\bea\label{8.3}
\lefteqn{([Y'_{Q(z)}(v_{1}, x_{1}), Y'_{Q(z)}(v_{2}, x_{2})]\lambda)
(w_{(1)}\otimes w_{(2)})}\nno\\
&&=\mbox{\rm Res}_{y_{2}}\mbox{\rm Res}_{y_{1}}x_{1}^{-1}
\delta\left(\frac{y_{1}-z}{x_{1}}\right)x_{2}^{-1}\delta
\left(\frac{y_{2}-z}{x_{2}}\right)\cdot \nno\\
&&\hspace{4em}\cdot \lambda([Y_{1}^{o}(v_{2}, y_{2}), Y_{1}^{o}(v_{1}, y_{1})]w_{(1)}
\otimes w_{(2)})\nno\\
&&\quad -\mbox{\rm Res}_{y_{2}}\mbox{\rm Res}_{y_{1}}x_{1}^{-1}
\delta\left(\frac{z-y_{1}}{-x_{1}}\right)x_{2}^{-1}
\delta\left(\frac{z-y_{2}}{-x_{2}}\right)\cdot \nno\\
&&\hspace{4em}\cdot \lambda(w_{(1)}\otimes [Y_{2}(v_{1}, y_{1}), Y_{2}(v_{2}, y_{2})]
w_{(2)})\nno\\
&&=\mbox{\rm Res}_{y_{2}}\mbox{\rm Res}_{y_{1}}x_{1}^{-1}
\delta\left(\frac{y_{1}-z}{x_{1}}\right)x_{2}^{-1}
\delta\left(\frac{y_{2}-z}{x_{2}}\right)\cdot \nno\\
&&\hspace{4em}\cdot \lambda\left(\mbox{\rm Res}_{x_{0}}y_{2}^{-1}
\delta\left(\frac{y_{1}-x_{0}}{y_{2}}\right)
Y^{o}_{1}(Y(v_{1}, x_{0})v_{2}, y_{2})w_{(1)}\otimes w_{(2)}\right) \nno\\
&&\quad -\mbox{\rm Res}_{y_{2}}\mbox{\rm Res}_{y_{1}}x_{1}^{-1}
\delta\left(\frac{z-y_{1}}{-x_{1}}\right)x_{2}^{-1}
\delta\left(\frac{z-y_{2}}{-x_{2}}\right)\cdot \nno\\
&&\hspace{4em}\cdot \lambda\left(
w_{(1)}\otimes \mbox{\rm Res}_{x_{0}}y_{2}^{-1}
\delta\left(\frac{y_{1}-x_{0}}{y_{2}}\right)Y_{2}(Y(v_{1}, x_{0})v_{2}, 
y_{2})w_{(2)}\right)\nno\\
&&=\mbox{\rm Res}_{x_{0}}\mbox{\rm Res}_{y_{2}}
\mbox{\rm Res}_{y_{1}}x_{1}^{-1}\delta\left(\frac{y_{1}-z}{x_{1}}\right)
x_{2}^{-1}\delta\left(\frac{y_{2}-z}{x_{2}}\right)y_{2}^{-1}
\delta\left(\frac{y_{1}-x_{0}}{y_{2}}\right)\cdot \nno\\
&&\hspace{4em}\cdot \lambda(
Y^{o}_{1}(Y(v_{1}, x_{0})v_{2}, y_{2})w_{(1)}\otimes w_{(2)})\nno\\
&&\quad -\mbox{\rm Res}_{x_{0}}\mbox{\rm Res}_{y_{2}}\mbox{\rm Res}_{y_{1}}
x_{1}^{-1}\delta\left(\frac{z-y_{1}}{-x_{1}}\right)x_{2}^{-1}
\delta\left(\frac{z-y_{2}}{-x_{2}}\right)y_{2}^{-1}\delta
\left(\frac{y_{1}-x_{0}}{y_{2}}\right)\cdot \nno\\
&&\hspace{4em}\cdot \lambda(w_{(1)}\otimes Y_{2}(Y(v_{1}, x_{0})v_{2}, y_{2})w_{(2)}).
\eea
But
\bea
\lefteqn{x_{1}^{-1}\delta\left(\frac{y_{1}-z}{x_{1}}\right)x_{2}^{-1}
\delta\left(\frac{y_{2}-z}{x_{2}}\right)y_{2}^{-1}
\delta\left(\frac{y_{1}-x_{0}}{y_{2}}\right)}\nno\\
&&=y_{1}^{-1}\delta\left(\frac{x_{1}+z}{y_{1}}\right)y_{2}^{-1}
\delta\left(\frac{x_{2}+z}{y_{2}}\right)(x_{2}+z)^{-1}
\delta\left(\frac{(x_{1}+z)-x_{0}}{x_{2}+z}\right)\nno\\
&&=y_{1}^{-1}\delta\left(\frac{x_{1}+z}{y_{1}}\right)y_{2}^{-1}
\delta\left(\frac{x_{2}+z}{y_{2}}\right)x_{2}^{-1}
\delta\left(\frac{x_{1}-x_{0}}{x_{2}}\right)\nno\\
&&=y_{1}^{-1}\delta\left(\frac{x_{1}+z}{y_{1}}\right)x_{2}^{-1}
\delta\left(\frac{y_{2}-z}{x_{2}}\right)x_{2}^{-1}
\delta\left(\frac{x_{1}-x_{0}}{x_{2}}\right)
\eea
and
\bea
\lefteqn{x_{1}^{-1}\delta\left(\frac{z-y_{1}}{-x_{1}}\right)x_{2}^{-1}
\delta\left(\frac{z-y_{2}}{-x_{2}}\right)y_{2}^{-1}
\delta\left(\frac{y_{1}-x_{0}}{y_{2}}\right)}\nno\\
&&=z^{-1}\delta\left(\frac{-x_{1}+y_{1}}{z}\right)z^{-1}
\delta\left(\frac{-x_{2}+y_{2}}{z}\right)y_{2}^{-1}
\delta\left(\frac{y_{1}-x_{0}}{y_{2}}\right)\nno\\
&&=\left({\displaystyle \sum_{m, n\in {\Z}}}
\frac{(-x_{1}+y_{1})^{m}}{z^{m+1}}
\frac{(-x_{2}+y_{2})^{n}}{z^{n+1}}\right) 
y_{2}^{-1}\delta\left(\frac{y_{1}-x_{0}}{y_{2}}\right)\nno\\
&&=\left({\displaystyle \sum_{m, n\in {\Z}}}(-x_{2}+y_{2})^{-1}
\left(\frac{-x_{1}+y_{1}}{-x_{2}+y_{2}}\right)^{m}\frac{(-x_{2}+y_{2})^{m+n+1}}
{z^{m+n+2}} \right)
y_{2}^{-1}\delta\left(\frac{y_{1}-x_{0}}{y_{2}}\right)\nno\\
&&=\left({\displaystyle \sum_{m, k\in {\Z}}}(-x_{2}+y_{2})^{-1}
\left(\frac{-x_{1}+y_{1}}{-x_{2}+y_{2}}\right)^{m}
z^{-1}\left(\frac{-x_{2}+y_{2}}{z}\right)^{k}\right) 
y_{2}^{-1}\delta\left(\frac{y_{1}-x_{0}}{y_{2}}\right)\nno\\
&&=(-x_{2}+y_{2})^{-1}\delta\left(\frac{-x_{1}+y_{1}}{-x_{2}+y_{2}}\right)
z^{-1}\delta\left(\frac{-x_{2}+y_{2}}{z}\right)
y_{2}^{-1}\delta\left(\frac{y_{1}-x_{0}}
{y_{2}}\right)\nno\\
&&=(-x_{2})^{-1}\delta\left(\frac{x_{1}-(y_{1}-y_{2})}{x_{2}}\right)
z^{-1}\delta\left(\frac{-x_{2}+y_{2}}{z}\right)
y_{1}^{-1}\delta\left(\frac{y_{2}+x_{0}}
{y_{1}}\right)\nno\\
&&=x_{2}^{-1}\delta\left(\frac{x_{1}-x_{0}}{x_{2}}\right)x_{2}^{-1}
\delta\left(\frac{z-y_{2}}{-x_{2}}\right)y_{1}^{-1}\delta
\left(\frac{y_{2}+x_{0}}{y_{1}}\right).
\eea
Thus (\ref{8.3}) becomes
\bea
\lefteqn{([Y'_{Q(z)}(v_{1}, x_{1}), Y'_{Q(z)}(v_{2}, x_{2})]\lambda)(w_{(1)}
\otimes w_{(2)})}\nno\\
&&=\mbox{\rm Res}_{x_{0}}\mbox{\rm Res}_{y_{2}}
\mbox{\rm Res}_{y_{1}}y_{1}^{-1}\delta\left(\frac{x_{1}+z}{y_{1}}\right)
x_{2}^{-1}\delta\left(\frac{y_{2}-z}{x_{2}}\right)x_{2}^{-1}
\delta\left(\frac{x_{1}-x_{0}}{x_{2}}\right)\cdot \nno\\
&&\hspace{4em}\cdot \lambda(
Y^{o}_{1}(Y(v_{1}, x_{0})v_{2}, y_{2})w_{(1)}\otimes w_{(2)})\nno\\
&&\quad -\mbox{\rm Res}_{x_{0}}
\mbox{\rm Res}_{y_{2}}\mbox{\rm Res}_{y_{1}}x_{2}^{-1}
\delta\left(\frac{x_{1}-x_{0}}{x_{2}}\right)x_{2}^{-1}
\delta\left(\frac{z-y_{2}}{-x_{2}}\right)y_{1}^{-1}
\delta\left(\frac{y_{2}+x_{0}}{y_{1}}\right)\cdot \nno\\
&&\hspace{4em}\cdot \lambda(w_{(1)}\otimes Y_{2}(Y(v_{1}, x_{0})v_{2},
y_{2})w_{(2)})\nno\\
&&=\mbox{\rm Res}_{x_{0}}x_{2}^{-1}\delta\left(\frac{x_{1}-x_{0}}
{x_{2}}\right)\cdot\nno\\
&&\hspace{2em}\cdot \biggl(\mbox{\rm Res}_{y_{2}}x_{2}^{-1}\delta\left(\frac{y_{2}-z}
{x_{2}}\right)\lambda(
Y^{o}_{1}(Y(v_{1}, x_{0})v_{2}, y_{2})w_{(1)}\otimes w_{(2)})\nno\\
&&\hspace{4em}-\mbox{\rm Res}_{y_{2}}x_{2}^{-1}\delta\left(\frac{z-y_{2}}
{-x_{2}}\right)\lambda(w_{(1)}\otimes Y_{2}(Y(v_{1}, x_{0})v_{2}, y_{2})w_{(2)})\biggr)
\nno\\
&&=\mbox{\rm Res}_{x_{0}}x_{2}^{-1}\delta\left(\frac{x_{1}-x_{0}}{x_{2}}\right)
(Y'_{Q(z)}(Y(v_{1}, x_{0})v_{2}, x_{2})\lambda)(w_{(1)}\otimes w_{(2)}).
\eea
Since $\lambda$, $w_{(1)}$ and $w_{(2)}$ are arbitrary, 
this  equality gives the commutator formula (\ref{commu-q-z})
for $Y'_{Q(z)}$. \epfv

When $V$ is in fact a conformal vertex algebra, we write
\begin{equation}\label{13.11-qz}
Y'_{Q(z)}(\omega, x)=\sum_{n\in {\mathbb Z}}L'_{Q(z)}(n)x^{-n-2}.
\end{equation}
Then {}from the last two propositions we see that the coefficient operators
of $Y'_{Q(z)}(\omega, x)$ satisfy the Virasoro algebra commutator
relations:
\begin{equation}\label{5.14}
[L'_{Q(z)}(m), L'_{Q(z)}(n)]
=(m-n)L'_{Q(z)}(m+n)+{\displaystyle\frac1{12}}
(m^3-m)\delta_{m+n,0}c.
\end{equation}
Moreover, in this case, by setting $v=\omega$ in (\ref{Y'qdef}) and
taking $\res_{x_0}x_0^{j+1}$ for $j=-1,0,1$, we see that
\begin{eqnarray}\label{LQ'(j)}
\lefteqn{(L'_{Q(z)}(j)\lambda)(w_{(1)}\otimes w_{(2)})}\nno\\
&&=\res_{x_1}(x_1-z)^{j+1}\lambda(Y^o_1(\omega,x_1)w_{(1)} \otimes
w_{(2)})\nno\\ 
&&\quad- \res_{x_1}
(-z+x_1)^{j+1}\lambda(w_{(1)} \otimes Y_2(\omega,x_1)w_{(2)})\nn
&&=\sum_{i=0}^{j+1}{j+1\choose i}(-z)^{i}\lambda(L(i-j)w_{(1)}
\otimes w_{(2)})\nn
&&\quad -\sum_{i=0}^{j+1}{j+1\choose i}(-z)^{i}
\lambda(w_{(1)}
\otimes L(j-i)w_{(2)})
\end{eqnarray}
for $j=-1,0,1$. If $V$ is just a M\"obius vertex algebra, we
define the actions $L'_{Q(z)}(j)$ on $(W_1\otimes W_2)^*$ by
the right-hand side of 
(\ref{LQ'(j)}) for $j=-1, 0$ and $1$. 

\begin{rema}\label{I-q-intw2}{\rm
In view of the action $L'_{Q(z)}(j)$, the ${\mathfrak s}{\mathfrak
l}(2)$-bracket relations (\ref{imq:Lj}) for a $Q(z)$-intertwining map
can be written as
\begin{equation}\label{I-q-intw2f}
(L'(j)w'_{(3)})\circ I=L'_{Q(z)}(j)(w'_{(3)}\circ I)
\end{equation}
for $w'_{(3)} \in W'_3$ and $j=-1$, $0$, and $1$.  }
\end{rema}

\begin{rema}\label{L'qjpreservesbetaspace}
{\rm  We have 
\[
L'_{Q(z)}(j)((W_{1}\otimes W_{2})^{*})^{(\beta)}\subset 
((W_{1}\otimes W_{2})^{*})^{(\beta)}
\]
for $j=-1, 0, 1$ and $\beta\in \tilde{A}$ (cf. Proposition \ref{tau-q-a-comp}).
}
\end{rema}

In the case that $V$ is a conformal vertex algebra, $L'_{Q(z)}(-1)$,
$L'_{Q(z)}(0)$ and $L'_{Q(z)}(1)$ realize the actions of $L_{-1}$,
$L_0$ and $L_1$ in ${\mathfrak s}{\mathfrak l}(2)$ (cf.\ (\ref{L_*}))
on $(W_1\otimes W_2)^*$.  In the case that $V$ is just a M\"{o}bius
vertex algebra, we now state this fact as a proposition.  This
proposition is needed in the proof of Theorem \ref{q-wk-mod} and
therefore also for Theorems \ref{q-generation} and
\ref{q-characterizationofbackslash}, but neither this proposition nor
any of these three theorems are needed anywhere else in this work, so
we omit the proof of this proposition.  Of course, however, the proof
is straightforward, as is the case with all the ${\mathfrak
s}{\mathfrak l}(2)$ formulas.

\begin{propo}\label{q-sl-2}
Let $V$ be a M\"{o}bius vertex algebra and let $W_{1}$ and $W_{2}$ be
generalized $V$-modules.  Then the operators $L'_{Q(z)}(-1)$,
$L'_{Q(z)}(0)$ and $L'_{Q(z)}(1)$ realize the actions of $L_{-1}$,
$L_0$ and $L_1$ in ${\mathfrak s}{\mathfrak l}(2)$ on $(W_1\otimes
W_2)^*$.
\end{propo}

We also have:

\begin{propo}\label{qz-l-y-comm}
Let $V$ be a M\"obius vertex algebra and let $W_{1}$ and $W_{2}$ be
generalized $V$-modules.  Then for $v\in V$,
\begin{eqnarray}
{[L(-1), Y'_{Q(z)}(v, x)]}&=&Y'_{Q(z)}(L(-1)v, x),\label{qz-sl-2-qz-y-1}\\
{[L(0), Y'_{Q(z)}(v, x)]}&=&Y'_{Q(z)}(L(0)v, x)+xY'_{Q(z)}(L(-1)v, x),
\label{qz-sl-2-qz-y-2}\\
{[L(1), Y'_{Q(z)}(v, x)]}&=&Y'_{Q(z)}(L(1)v, x)
+2xY'_{Q(z)}(L(0)v, x)+x^{2}Y'_{Q(z)}(L(-1)v, x),\label{qz-sl-2-qz-y-3}\nn
&&
\end{eqnarray}
where for brevity we write $L'_{Q(z)}(j)$ acting on 
$(W_{1}\otimes W_{2})^{*}$ as $L(j)$. 
\end{propo}
\pf We prove only (\ref{qz-sl-2-qz-y-2}) since it is needed for Remark
\ref{q-stableundercomponentops} and in Section 6.  We omit the proofs
of (\ref{qz-sl-2-qz-y-1}) and (\ref{qz-sl-2-qz-y-3}) for the same
reasons as above; they are used only for Theorems
\ref{q-wk-mod}--\ref{q-characterizationofbackslash}.

Let $\lambda\in (W_{1}\otimes W_{2})^{*}$, $w_{(1)}\in W_{1}$ and
$w_{(2)}\in W_{2}$. Using (\ref{LQ'(j)}), (\ref{Y'qdef}), the
commutator formulas for $L(j)$ and $Y_{1}(v, x_{0})$ for $j=-1, 0, 1$
and $v\in V$ (recall Definition \ref{moduleMobius}), and the
commutator formulas for $L(j)$ and $Y_{2}^{o}(v, x)$ for $j=-1, 0, 1$
and $v\in V$ (recall Lemma \ref{sl2opposite}), we obtain
\begin{eqnarray*}
\lefteqn{([L(0), Y'_{Q(z)}(v, x)]\lambda)(w_{(1)}\otimes w_{(2)})}\nn
&&=(Y'_{Q(z)}(v, x)\lambda)(L(0)w_{(1)}\otimes w_{(2)})\nn
&&\quad -z(Y'_{Q(z)}(v, x)\lambda)(L(1)w_{(1)}\otimes w_{(2)})\nn
&&\quad -
(Y'_{Q(z)}(v, x)\lambda)(w_{(1)}\otimes L(0)w_{(2)})\nn
&&\quad +z
(Y'_{Q(z)}(v, x)\lambda)(w_{(1)}\otimes L(-1)w_{(2)})\nn
&&\quad -(L(0)\lambda)(Y^o_1(v, x + z)w_{(1)}
\otimes w_{(2)})\nn
&&\quad +\res_{x_1}x^{-1} \delta \left(\frac{z-x_1}{-x}\right)
(L(0)\lambda)(w_{(1)}
\otimes Y_2(v,x_1)w_{(2)})\nn
&&=\lambda(Y^o_1(v, x + z)L(0)w_{(1)}\otimes w_{(2)})\nn
&&\quad-\res_{x_1}x^{-1} \delta \left(\frac{z-x_1}{-x}\right)
\lambda(L(0)w_{(1)}\otimes Y_2(v,x_1)w_{(2)})\nn
&&\quad -z\lambda(Y^o_1(v, x + z)L(1)w_{(1)}\otimes w_{(2)})\nn
&&\quad +z\res_{x_1}x^{-1} \delta \left(\frac{z-x_1}{-x}\right)
\lambda(L(1)w_{(1)}\otimes Y_2(v,x_1)w_{(2)})\nn
&&\quad -
\lambda(Y^o_1(v, x + z)w_{(1)}\otimes L(0)w_{(2)})\nn
&&\quad +\res_{x_1}x^{-1} \delta \left(\frac{z-x_1}{-x}\right)
\lambda(w_{(1)}\otimes Y_2(v,x_1)L(0)w_{(2)})\nn
&&\quad +z
\lambda(Y^o_1(v, x + z)w_{(1)}\otimes L(-1)w_{(2)})\nn
&&\quad -z\res_{x_1}x^{-1} \delta \left(\frac{z-x_1}{-x}\right)
\lambda(w_{(1)}\otimes Y_2(v,x_1)L(-1)w_{(2)})\nn
&&\quad -\lambda(L(0)Y^o_1(v, x + z)w_{(1)}
\otimes w_{(2)})\nn
&&\quad +z\lambda(L(1)Y^o_1(v, x + z)w_{(1)}
\otimes w_{(2)})\nn
&&\quad +\lambda(Y^o_1(v, x + z)w_{(1)}
\otimes L(0)w_{(2)})\nn
&&\quad -z\lambda(Y^o_1(v, x + z)w_{(1)}
\otimes L(-1)w_{(2)})\nn
&&\quad +\res_{x_1}x^{-1} \delta \left(\frac{z-x_1}{-x}\right)
\lambda(L(0)w_{(1)}
\otimes Y_2(v,x_1)w_{(2)})\nn
&&\quad -z\res_{x_1}x^{-1} \delta \left(\frac{z-x_1}{-x}\right)
\lambda(L(1)w_{(1)}
\otimes Y_2(v,x_1)w_{(2)})\nn
&&\quad -\res_{x_1}x^{-1} \delta \left(\frac{z-x_1}{-x}\right)
\lambda(w_{(1)}
\otimes L(0)Y_2(v,x_1)w_{(2)})\nn
&&\quad +z\res_{x_1}x^{-1} \delta \left(\frac{z-x_1}{-x}\right)
\lambda(w_{(1)}
\otimes L(-1)Y_2(v,x_1)w_{(2)})\nn
&&=\lambda([Y^o_1(v, x + z), L(0)]w_{(1)}\otimes w_{(2)})\nn
&&\quad -z\lambda([Y^o_1(v, x + z), L(1)]w_{(1)}\otimes w_{(2)})\nn
&&\quad +z\res_{x_1}x^{-1} \delta \left(\frac{z-x_1}{-x}\right)
\lambda(w_{(1)}\otimes [L(-1), Y_2(v,x_1)]w_{(2)})\nn
&&\quad -\res_{x_1}x^{-1} \delta \left(\frac{z-x_1}{-x}\right)
\lambda(w_{(1)}
\otimes [L(0), Y_2(v,x_1)]w_{(2)})\nn
&&=\lambda(Y^o_1((L(0)+(x+z)L(-1)v, x + z)w_{(1)}\otimes w_{(2)})\nn
&&\quad -z\lambda(Y^o_1(L(-1)v, x + z)w_{(1)}\otimes w_{(2)})\nn
&&\quad +z\res_{x_1}x^{-1} \delta \left(\frac{z-x_1}{-x}\right)
\lambda(w_{(1)}\otimes Y_2(L(-1)v,x_1)w_{(2)})\nn
&&\quad -\res_{x_1}x^{-1} \delta \left(\frac{z-x_1}{-x}\right)
\lambda(w_{(1)}
\otimes Y_2((L(0)+x_{1}L(-1))v,x_1)w_{(2)})\nn
&&=\lambda(Y^o_1((L(0)+xL(-1)v, x + z)w_{(1)}\otimes w_{(2)})\nn
&&\quad -\res_{x_1}x^{-1} \delta \left(\frac{z-x_1}{-x}\right)
\lambda(w_{(1)}
\otimes Y_2((L(0)+(x-z)L(-1))v,x_1)w_{(2)})\nn
&&=\lambda(Y^o_1((L(0)+xL(-1)v, x + z)w_{(1)}\otimes w_{(2)})\nn
&&\quad -\res_{x_1}x^{-1} \delta \left(\frac{z-x_1}{-x}\right)
\lambda(w_{(1)}
\otimes Y_2((L(0)+xL(-1))v,x_1)w_{(2)})\nn
&&=(Y'_{Q(z)}((L(0)+xL(-1))v, x)\lambda)(w_{(1)}\otimes w_{(2)})\nn
&&=(Y'_{Q(z)}(L(0)v, x)\lambda)(w_{(1)}\otimes w_{(2)})
+(xY'_{Q(z)}(L(-1)v, x)\lambda)(w_{(1)}\otimes w_{(2)}),
\end{eqnarray*}
proving (\ref{qz-sl-2-qz-y-2}).
\epfv

Let $W_3$ also be an object of ${\cal C}$. Note that $V\otimes
\iota_{+}{\mathbb C}[t, t^{-1}, (z+t)^{-1}]$ acts on $W'_3$ in the
natural way.  The following result provides further motivation for the
definition of our action (\ref{5.2}) on $(W_1\otimes W_2)^{*}$; recall
the discussion preceding Proposition \ref{pz}:

\begin{propo}\label{qz}
Let $W_{1}$, $W_{2}$ and $W_{3}$ be generalized $V$-modules.
Under the natural isomorphism described in Remark \ref{alternateformoflemma}
between the space of $\tilde{A}$-compatible linear maps
\[
I:W_{1}\otimes W_{2} \rightarrow \overline{W_{3}}
\]
and the space of $\tilde{A}$-compatible linear maps
\[
J:W'_{3} \rightarrow (W_{1}\otimes W_{2})^{*}
\]
determined by (\ref{IcorrespondstoJ}), the $Q(z)$-intertwining maps
$I$ of type ${W_3\choose W_1\, W_2}$ correspond exactly to the
(grading restricted) $\tilde{A}$-compatible maps $J$ that intertwine
the actions of both
\[
V \otimes \iota_{+}{\mathbb C}[t, t^{-1}, (z+t)^{-1}]
\]
and ${\mathfrak s} {\mathfrak l}(2)$ on $W'_{3}$ and on $(W_1\otimes
W_2)^{*}$.
\end{propo}
\pf 
In view of (\ref{IcorrespondstoJalternateform}), Remark
\ref{I-intw-q} asserts that (\ref{imq:def'}), or equivalently,
(\ref{imq:def}), is equivalent to the condition
\begin{equation}\label{q-j-tau}
J\left(\tau_{W'_3}\left(z^{-1}\delta\left(\frac{x_1-x_0}{z}\right)
Y_{t}(v, x_0)\right)w'_{(3)}\right)
=\tau_{Q(z)}\left(z^{-1}\delta\left(\frac{x_1-x_0}{z}\right)
Y_{t}(v, x_0)\right)J(w'_{(3)}),
\end{equation}
that is, the condition that $J$ intertwines the actions of $V \otimes
\iota_{+}{\mathbb C}[t,t^{- 1}, (z+t)^{-1}]$ on $W'_{3}$ and on
$(W_1\otimes W_2)^{*}$ (recall (\ref{3.12})--(\ref{3.13})).
Similarly, Remark \ref{I-q-intw2} asserts that (\ref{imq:Lj}) is
equivalent to the condition
\begin{equation}\label{q-j-lj}
J(L'(j)w'_{(3)}) = L'_{Q(z)}(j)J(w'_{(3)})
\end{equation}
for $j=-1$, $0$, $1$, that is, the condition that $J$ intertwines the
actions of ${\mathfrak s} {\mathfrak l}(2)$ on $W'_{3}$ and on
$(W_1\otimes W_2)^{*}$.
\epfv

\begin{nota}\label{qscriptN}
{\rm Given generalized $V$-modules $W_1$, $W_2$ and $W_3$, we shall
write ${\cal N}[Q(z)]_{W'_3}^{(W_1 \otimes W_2)^{*}}$ for the
space of (grading restricted) $\tilde{A}$-compatible linear maps
\[
J:W'_{3} \rightarrow (W_{1}\otimes W_{2})^{*}
\]
that intertwine the actions of both
\[
V \otimes \iota_{+}{\mathbb C}[t,t^{- 1}, (z+t)^{-1}]
\]
and ${\mathfrak s} {\mathfrak l}(2)$ on $W'_{3}$ and on $(W_1\otimes
W_2)^{*}$.
Note that Proposition \ref{qz} gives a natural linear isomorphism
\begin{eqnarray*}
{\cal M}[Q(z)]^{W_3}_{W_1 W_2} &
\stackrel{\sim}{\longrightarrow} & {\cal N}[Q(z)]_{W'_3}^{(W_1 \otimes
W_2)^{*}}\nno\\ I & \mapsto & J
\end{eqnarray*}
(recall {}from Definition \ref{im:qimdef} the notation for the space of
$Q(z)$-intertwining maps).  As in Notation \ref{scriptN}, we still use
the symbol ``prime'' to denote this isomorphism in both directions:
\begin{eqnarray*}
{\cal M}[Q(z)]^{W_3}_{W_1 W_2} & \stackrel{\sim}{\longrightarrow} & {\cal
N}[Q(z)]_{W'_3}^{(W_1 \otimes W_2)^{*}}\nno\\
I & \mapsto & I'\nno\\
J' & \leftarrow\!\!\!{\scriptstyle |} & J,
\end{eqnarray*}
so that in particular,
\[
I'' = I \;\;\mbox{ and }\;\; J'' = J
\]
for $I \in {\cal M}[Q(z)]^{W_3}_{W_1 W_2}$ and 
$J \in {\cal N}[Q(z)]_{W'_3}^{(W_1
\otimes W_2)^{*}}$, and the relation between $I$ and $I'$ is
determined by
\[
\langle w'_{(3)}, I(w_{(1)}\otimes w_{(2)})\rangle
=I'(w'_{(3)})(w_{(1)}\otimes w_{(2)})
\]
for $w_{(1)}\in W_{1}$, $w_{(2)}\in W_{2}$ and $w'_{(3)}\in W'_{3}$,
or equivalently,
\[
w'_{(3)}\circ I = I'(w'_{(3)}).
\]
}
\end{nota}

\begin{rema}
{\rm Combining Proposition \ref{qz} with Proposition \ref{Q-cor}, we see
that for any integer $p$, we also have a natural linear
isomorphism
\[
{\cal N}[Q(z)]_{W'_3}^{(W_1 \otimes W_2)^{*}} \stackrel{\sim}{\longrightarrow}
{\cal V}^{W'_1}_{W'_3 W_2}
\]
{}from ${\cal N}[Q(z)]_{W'_3}^{(W_1 \otimes W_2)^{*}}$ to the space of
logarithmic intertwining operators of type ${W'_1\choose W'_3\,W_2}$.
In particular, given any such logarithmic intertwining
operator ${\cal Y}$ and integer $p$, the map 
\[
(I^{Q(z)}_{{\cal Y}, p})':
W'_3\to (W_1\otimes W_2)^{*}
\]
defined by
\[
(I^{Q(z)}_{{\cal Y}, p})'(w'_{(3)})(w_{(1)}\otimes w_{(2)})=\bra w_{(1)},
{\cal Y}(w'_{(3)}, e^{l_{p}(z)})w_{(2)}\ket_{W'_1}
\]
is $\tilde{A}$-compatible and intertwines both actions on both spaces.}
\end{rema}

We have formulated the notions of $Q(z)$-product and $Q(z)$-tensor
product using $Q(z)$-intertwining maps (Definitions \ref{qz-product}
and \ref{qz-tp}).  Now that we know that $Q(z)$-intertwining maps can
be interpreted as in Proposition \ref{qz} (and Notation
\ref{qscriptN}), we can reformulate the notions of $Q(z)$-product and
$Q(z)$-tensor product correspondingly (the proof of the next result is
the same as that of Proposition \ref{productusingI'}):

\begin{propo}\label{q-productusingI'}
Let ${\cal C}_1$ be either of the categories ${\cal M}_{sg}$ or ${\cal
GM}_{sg}$, as in Definition \ref{qz-product}.  For $W_1, W_2\in
\ob{\cal C}_1$, a $Q(z)$-product $(W_3;I_3)$ of $W_1$ and $W_2$
(recall Definition \ref{qz-product}) amounts to an object $(W_3,Y_3)$
of ${\cal C}_1$ equipped with a map $I'_3 \in {\cal N}[Q(z)]_{W'_3}^{(W_1
\otimes W_2)^{*}}$, that is, equipped with an $\tilde{A}$-compatible
map
\[
I'_3:W'_{3} \rightarrow (W_{1}\otimes W_{2})^{*}
\]
that intertwines the two actions of $V \otimes \iota_{+}{\mathbb
C}[t,t^{- 1}, (z+t)^{-1}]$ and of ${\mathfrak s} {\mathfrak
l}(2)$.  The map $I'_3$ corresponds to the $Q(z)$-intertwining map
\[
I_3:W_{1}\otimes W_{2} \rightarrow \overline{W_{3}}
\]
as above:
\[
I'_3(w'_{(3)}) = w'_{(3)}\circ I_3
\]
for $w'_{(3)} \in W'_{3}$ (recall
\ref{IcorrespondstoJalternateform})).  Denoting this structure by
$(W_3,Y_3;I'_3)$ or simply by $(W_3;I'_3)$, let $(W_4;I'_4)$ be
another such structure.  Then a morphism of $Q(z)$-products {}from $W_3$
to $W_4$ amounts to a module map $\eta: W_3 \to W_4$ such that the
diagram
\begin{center}
\begin{picture}(100,60)
\put(-2,0){$W'_4$}
\put(13,4){\vector(1,0){104}}
\put(119,0){$W'_3$}
\put(38,50){$(W_1\otimes W_2)^*$}
\put(13,12){\vector(3,2){50}}
\put(118,12){\vector(-3,2){50}}
\put(65,8){$\eta'$}
\put(23,27){$I'_4$}
\put(98,27){$I'_3$}
\end{picture}
\end{center}
commutes, where $\eta'$ is the natural map given by (\ref{fprime}).
\epf
\end{propo}

\begin{corol}\label{q-tensorproductusingI'}
Let ${\cal C}$ be a full subcategory of either ${\cal M}_{sg}$ or
${\cal GM}_{sg}$, as in Definition \ref{qz-tp}.  For $W_1, W_2\in
\ob{\cal C}$, a $Q(z)$-tensor product $(W_0; I_0)$ of $W_1$ and $W_2$
in ${\cal C}$, if it exists, amounts to an object $W_0 =
W_1\boxtimes_{Q(z)} W_2$ of ${\cal C}$ and a structure $(W_0 =
W_1\boxtimes_{Q(z)} W_2; I'_0)$ as in Proposition
\ref{q-productusingI'}, with
\[
I'_0: (W_1\boxtimes_{Q(z)} W_2)' \longrightarrow (W_1\otimes W_2)^*
\]
in ${\cal N}[Q(z)]_{(W_1\boxtimes_{Q(z)} W_2)'}^{(W_1 \otimes W_2)^{*}}$,
such that for any such pair $(W; I')$ $(W\in \ob \mathcal{C})$, with
\[
I': W' \longrightarrow (W_1\otimes W_2)^*
\]
in ${\cal N}[Q(z)]_{W'}^{(W_1 \otimes W_2)^{*}}$, there is a unique module
map
\[
\chi: W' \longrightarrow (W_1\boxtimes_{Q(z)} W_2)'
\]
such that the diagram
\begin{center}
\begin{picture}(100,60)
\put(-2,0){$W'$}
\put(13,4){\vector(1,0){104}}
\put(119,0){$(W_1\boxtimes_{Q(z)} W_2)'$}
\put(38,50){$(W_1\otimes W_2)^*$}
\put(13,12){\vector(3,2){50}}
\put(118,12){\vector(-3,2){50}}
\put(65,8){$\chi$}
\put(23,27){$I'$}
\put(98,27){$I'_0$}
\end{picture}
\end{center}
commutes.  Here $\chi = \eta'$, where $\eta$ is a correspondingly
unique module map
\[
\eta: W_1\boxtimes_{Q(z)} W_2 \longrightarrow W.
\]
Also, the map $I_0'$, which is $\tilde{A}$-compatible and which
intertwines the two actions of $V \otimes \iota_{+}{\mathbb C}[t,t^{-
1}, (z+t)^{-1}]$ and of ${\mathfrak s} {\mathfrak l}(2)$, is
related to the $Q(z)$-intertwining map
\[
I_0 = \boxtimes_{Q(z)}: W_1\otimes W_2 \longrightarrow 
\overline{W_1\boxtimes_{Q(z)} W_2}
\]
by
\[
I_0'(w') = w' \circ \boxtimes_{Q(z)}
\]
for $w' \in (W_1\boxtimes_{Q(z)} W_2)'$, that is,
\[
I_0'(w')(w_{(1)}\otimes w_{(2)}) = \langle w',w_{(1)}\boxtimes_{Q(z)}
w_{(2)} \rangle
\]
for $w_{(1)}\in W_{1}$ and $w_{(2)}\in W_{2}$, using the notation 
(\ref{q-boxtensorofelements}).
\epf
\end{corol}

\begin{defi}
{\rm For $W_{1}, W_{2}\in \ob \mathcal{C}$, define the subset 
\[
W_{1}\hboxtr_{Q(z)}W_{2}\subset (W_{1}\otimes W_{2})^{*}
\]
of $(W_{1}\otimes W_{2})^{*}$ to be the union of the images
\[
I'(W')\subset (W_{1}\otimes W_{2})^{*}
\]
as $(W; I)$ ranges through all the $Q(z)$-products of $W_{1}$ and $W_{2}$ with 
$W\in \ob
\mathcal{C}$. Equivalently, $W_{1}\hboxtr_{P(z)}W_{2}$ is the union 
of the images 
$I'(W')$ as $W$ (or $W'$) ranges through $\ob
\mathcal{C}$ and $I'$ ranges through 
$\mathcal{N}[Q(z)]_{W'}^{(W_{1}\otimes W_{2})^{*}}$---the space of 
$\tilde{A}$-compatible linear maps
\[
W'\to (W_{1}\otimes W_{2})^{*}
\]
intertwining the actions of both 
\[
V\otimes \iota_{+}\C[t, t^{-1}, (z+t)^{-1}]
\]
and $\mathfrak{s}\mathfrak{l}(2)$ on both spaces.}
\end{defi}

\begin{rema}
{\rm Since $\mathcal{C}$ is closed under direct sums (Assumption 
\ref{assum-c}), it is clear that $W_{1}\hboxtr_{Q(z)}W_{2}$ is
in fact a linear subspace of $(W_{1}\otimes W_{2})^{*}$, and in particular,
it can be defined alternatively as the sum of all the images $I'(W')$:
\begin{equation}\label{q-hboxtr-sum}
W_{1}\hboxtr_{Q(z)}W_{2}=\sum I'(W') = \bigcup I'(W')\subset 
(W_{1}\otimes W_{2})^{*},
\end{equation}
where the sum and union both range over $W\in \ob
\mathcal{C}$, $I\in \mathcal{M}[Q(z)]_{W_{1}W_{2}}^{W}$.}
\end{rema}

For any generalized $V$-modules $W_{1}$ and $W_{2}$,
using the operator $L'_{Q(z)}(0)$ (recall (\ref{LQ'(j)}))
on $(W_{1}\otimes W_{2})^{*}$ we define
the generalized $L'_{Q(z)}(0)$-eigenspaces 
$((W_{1}\otimes W_{2})^{*})_{[n]; Q(z)}$ for $n\in \C$ in the usual way:
\begin{equation}
((W_{1}\otimes W_{2})^{*})_{[n]; Q(z)}=\{w\in (W_{1}\otimes W_{2})^{*}\;|\;
(L'_{Q(z)}(0)-n)^{m}w=0 \;{\rm for}\; m\in \N \;
\mbox{\rm sufficiently large}\}.
\end{equation}
Then we have the (proper) subspace 
\begin{equation}
\coprod_{n\in \C}((W_{1}\otimes W_{2})^{*})_{[n]; Q(z)}\subset 
(W_{1}\otimes W_{2})^{*}.
\end{equation}
We also define the ordinary $L'_{Q(z)}(0)$-eigenspaces
$((W_{1}\otimes W_{2})^{*})_{(n); Q(z)}$ in the usual way:
\begin{equation}
((W_{1}\otimes W_{2})^{*})_{(n); Q(z)}=\{w\in (W_{1}\otimes W_{2})^{*}\;|\;
L'_{P(z)}(0)w=nw \}.
\end{equation}
Then we have the (proper) subspace
\begin{equation}
\coprod_{n\in \C}((W_{1}\otimes W_{2})^{*})_{(n); Q(z)}\subset 
(W_{1}\otimes W_{2})^{*}.
\end{equation}

Just as in Proposition \ref{im:abc}, we have:

\begin{propo}\label{im-q:abc}
Let $W_{1}, W_{2}\in \ob \mathcal{C}$. 

(a)
The elements of $W_1\hboxtr_{Q(z)} W_2$ are exactly the linear functionals
on $W_{1}\otimes W_{2}$ of the form $w'\circ
I(\cdot\otimes \cdot)$ for some $Q(z)$-intertwining map $I$ of type
${W\choose W_1\,W_2}$ and some $w'\in W'$, $W\in\ob{\cal C}$.

(b) Let $(W; I)$ be any $Q(z)$-product of $W_{1}$ and $W_{2}$, with 
$W$ any generalized $V$-module. Then for $n\in \C$,
\[
I'(W'_{[n]}) \subset  ((W_{1}\otimes W_{2})^{*})_{[n]; Q(z)}
\]
and 
\[
I'(W'_{(n)})\subset ((W_{1}\otimes W_{2})^{*})_{(n); Q(z)}.
\]

(c) The structure $(W_1\hboxtr_{Q(z)} W_2,Y'_{Q(z)})$ (recall
(\ref{y'-q-z})) satisfies all the axioms in the definition of
(strongly $\tilde{A}$-graded) generalized $V$-module except perhaps
for the two grading conditions (\ref{set:dmltc}) and
(\ref{set:dmfin}).

(d) Suppose that the objects of the category $\mathcal{C}$ consist
only of (strongly $\tilde{A}$-graded) {\em ordinary}, as opposed to
{\em generalized}, $V$-modules.  Then the structure
$(W_1\hboxtr_{Q(z)} W_2,Y'_{Q(z)})$ satisfies all the axioms in the
definition of (strongly $\tilde{A}$-graded ordinary) $V$-module except
perhaps for (\ref{set:dmltc}) and (\ref{set:dmfin}).
\end{propo}
\pf 
Part (a) is clear {}from the definition of $W_1\hboxtr_{Q(z)} W_2$, and 
(b) follows {}from (\ref{q-j-lj}) with $j=0$.

To prove (c), let $(W; I)$ be any any $Q(z)$-product of $W_{1}$ and
$W_{2}$, with $W$ any generalized $V$-module. Then $(I'(W'),
Y'_{Q(z)})$ satisfies all the conditions in the definition of
(strongly $\tilde{A}$-graded) generalized $V$-module since $I'$ is
$\tilde{A}$-compatible and intertwines the actions of $V\otimes
{\mathbb C}[t,t^{-1}]$ and of ${\mathfrak s}{\mathfrak l}(2)$; the
$\C$-grading follows {}from Part (b). Since $W_1\hboxtr_{Q(z)} W_2$ is
the sum of these structures $I'(W')$ over $W\in \ob \mathcal{C}$
(recall (\ref{hboxtr-sum})), $(W_1\hboxtr_{Q(z)} W_2, Y'_{Q(z)})$
satisfies all the conditions in the definition of generalized module
except perhaps for (\ref{set:dmltc}) and (\ref{set:dmfin}).

Part (d) is proved by the same argument as for (c): For $(W; I)$ any
$Q(z)$-product of possibly generalized $V$-modules $W_{1}$ and
$W_{2}$, with $W$ any ordinary $V$-module, $(I'(W'), Y'_{Q(z)})$
satisfies all the conditions in the definition of (strongly
$\tilde{A}$-graded) ordinary $V$-module; the $\C$-grading (by ordinary
$L'_{Q(z)}(0)$-eigenspaces) again follows {}from Part (b).  \epfv

We now have the following generalization of Proposition 5.8 in
\cite{tensor1}, characterizing $W_{1}\boxtimes_{Q(z)}W_{2}$, including its
existence, in terms of $W_1\hboxtr_{Q(z)} W_2$; 
the proof is the same as that of Proposition 
\ref{tensor1-13.7}:

\begin{propo}\label{tensor1-5.7}
Let $W_{1}, W_{2}\in \ob \mathcal{C}$. 
If $(W_1\hboxtr_{Q(z)} W_2, Y'_{Q(z)})$ is an object of ${\cal C}$,
denote by $(W_1\boxtimes_{Q(z)} W_2, Y_{Q(z)})$ its contragredient
module. Then the $Q(z)$-tensor product of $W_{1}$ and $W_{2}$ 
in ${\cal C}$ exists and is
$(W_1\boxtimes_{Q(z)} W_2, Y_{Q(z)}; i')$, where $i$ is the natural
inclusion {}from $W_1\hboxtr_{Q(z)} W_2$ to $(W_1\otimes W_2)^*$ (recall
Notation \ref{qscriptN}).
Conversely, let us assume that $\mathcal{C}$ is closed under quotients. 
If the $Q(z)$-tensor product of $W_1$ and $W_2$ in ${\cal
C}$ exists, then $(W_1\hboxtr_{Q(z)} W_2, Y'_{Q(z)})$ is an object of
${\cal C}$.\epf
\end{propo}

\begin{rema}
{\rm Suppose that $W_1\hboxtr_{Q(z)} W_2$ is an object of
$\mathcal{C}$.  {}From Corollary \ref{q-tensorproductusingI'} and
Proposition \ref{tensor1-5.7} we see that
\begin{equation}\label{boxpair-q}
\langle\lambda, w_{(1)}\boxtimes_{Q(z)}w_{(2)}\rangle
\lbar_{W_1\boxtimes_{Q(z)} W_2}=
\lambda(w_{(1)}\otimes w_{(2)})
\end{equation}
for $\lambda\in W_1\hboxtr_{Q(z)} W_2\subset (W_1\otimes W_2)^*$,
$w_{(1)}\in W_1$ and $w_{(2)}\in W_2$.}
\end{rema}

As in the $P(z)$-case, our next goal is 
to present an alternative description of the
subspace $W_1\hboxtr_{Q(z)} W_2$ of $(W_1\otimes W_2)^*$. The main
ingredient of this description will be the ``$Q(z)$-compatibility
condition,'' as was the case in \cite{tensor1}--\cite{tensor2}.

Take $W_{1}$ and $W_{2}$ to be arbitrary generalized $V$-modules.  Let
$(W,I)$ ($W$ a generalized $V$-module) be a $Q(z)$-product of $W_{1}$
and $W_{2}$ and let $w'\in W'$. Then {}from (\ref{q-j-tau}),
Proposition \ref{q-productusingI'}, (\ref{tau-w-comp}), (\ref{3.7})
and (\ref{y'-q-z}), we have, for all $v\in V$,
\begin{eqnarray}\label{5.18}
\lefteqn{\tau_{Q(z)}\left(z^{-1}\delta\left(\frac{x_1-x_0}{z}\right)
Y_{t}(v, x_0)\right)I'(w')}\nno\\
&&=I'\left(\tau_{W'}\left(z^{-1}\delta\left(\frac{x_1-x_0}{z}\right)
Y_{t}(v, x_0)\right)w'\right)\nno\\
&&=I'\left(z^{-1}\delta\left(\frac{x_1-x_0}{z}\right)
Y_{W'}(v, x_0)w'\right)\nno\\
&&=z^{-1}\delta\left(\frac{x_1-x_0}{z}\right)I'(Y_{W'}(v,
x_0)w')\nno\\
&&=z^{-1}\delta\left(\frac{x_1-x_0}{z}\right)I'(\tau_{W'}(Y_{t}(v,
x_0))w')\nno\\
&&=z^{-1}\delta\left(\frac{x_1-x_0}{z}\right)\tau_{Q(z)}(Y_{t}(v,
x_0))I'(w')\nn
&&=z^{-1}\delta\left(\frac{x_1-x_0}{z}\right)
Y'_{Q(z)}(v, x_0)I'(w').
\end{eqnarray}
That is, 
$I'(w')$ satisfies the following nontrivial and subtle condition on
$\lambda \in (W_1\otimes W_2)^{*}$:

\begin{description}
\item{\bf The $Q(z)$-compatibility condition}

(a) The {\em $Q(z)$-lower truncation condition}: For all $v\in V$, the formal
Laurent series $Y'_{Q(z)}(v, x)\lambda$ involves only finitely many
negative powers of $x$.

(b) The following formula holds:
\begin{eqnarray}\label{cpb-q}
\lefteqn{\tau_{Q(z)}\left(z^{-1}\delta\left(\frac{x_1-x_0}{z}\right)
Y_{t}(v, x_0)\right)\lambda}\nno\\
&&=z^{-1}\delta\left(\frac{x_1-x_0}{z}\right)
Y'_{Q(z)}(v, x_0)\lambda  \;\;\mbox{ for all }\;v\in V.
\end{eqnarray}
(Note that the two sides of (\ref{cpb-q}) are not {\it a priori} equal
for general $\lambda\in (W_1\otimes W_2)^{*}$. Note also that Condition 
(a) insures that the right-hand side in Condition (b) is 
well defined.)
\end{description}

\begin{nota}
{\rm Note that the set of elements of $(W_1\otimes W_2)^*$ satisfying
either  the full $Q(z)$-compatibility
condition or Part (a) of this condition forms a subspace. 
We shall denote the space of elements of $(W_1\otimes W_2)^*$ satisfying
the $Q(z)$-compatibility
condition by
\[
\comp_{Q(z)}((W_1\otimes W_2)^*).
\]}
\end{nota}

We know that each space $((W_1\otimes W_2)^*)^{(\beta)}$ is 
$L'_{Q(z)}(0)$-stable
(recall Proposition \ref{tau-q-a-comp} and Remark
\ref{L'qjpreservesbetaspace}), so that we may consider the subspaces
\[
\coprod_{n\in \C}((W_1\otimes W_2)^*)_{[n];Q(z)}^{(\beta)} \subset
((W_1\otimes W_2)^*)^{(\beta)}
\]
and 
\[
\coprod_{n\in \C}((W_1\otimes W_2)^*)_{(n);Q(z)}^{(\beta)} \subset
((W_1\otimes W_2)^*)^{(\beta)}
\]
(recall Remark \ref{generalizedeigenspacedecomp}).  We define the
two subspaces
\begin{equation}\label{W1W2_[C];q^Atilde}
((W_1\otimes W_2)^*)_{[{\mathbb C}]; Q(z)}^{( \tilde A )}=
\coprod_{n\in
\C}\coprod_{\beta\in \tilde{A}}((W_1\otimes W_2)^*)_{[n];Q(z)}^{(\beta)}
\subset (W_1\otimes W_2)^*
\end{equation}
and 
\begin{equation}\label{W1W2_(C);q^Atilde}
((W_1\otimes W_2)^*)_{({\mathbb C});Q(z)}^{( \tilde A )}=
\coprod_{n\in
\C}\coprod_{\beta\in \tilde{A}}((W_1\otimes W_2)^*)_{(n);Q(z)}^{(\beta)}
\subset (W_1\otimes W_2)^*.
\end{equation}

\begin{rema}\label{q-singleanddoublegraded}
{\rm Any $L'_{Q(z)}(0)$-stable subspace of $((W_1\otimes
W_2)^*)_{[{\mathbb C}];Q(z)}^{( \tilde A )}$ is graded by generalized
eigenspaces (again recall Remark \ref{generalizedeigenspacedecomp}),
and if such a subspace is also $\tilde A$-graded, then it is doubly
graded; similarly for subspaces of $((W_1\otimes W_2)^*)_{({\mathbb
C});Q(z)}^{( \tilde A )}$.}
\end{rema}

We have:

\begin{lemma}\label{q-a-tilde-comp}
Suppose that $\lambda\in ((W_1\otimes W_2)^*)_{[{\mathbb C}];Q(z)}^{(\tilde
A )}$ satisfies the $Q(z)$-compatibility condition. Then every
$\tilde{A}$-homogeneous component of $\lambda$ also satisfies this
condition.
\end{lemma}
\pf
When $v\in V$ is $\tilde{A}$-homogeneous, 
\[
\tau_{Q(z)}\bigg(z^{-1}\delta\bigg(\frac{x_1-x_{0}}{z} \bigg)
Y_{t}(v, x_0)\bigg)\;\;\mbox{ and }\;\;
z^{-1}\delta\bigg(\frac{x_1-x_{0}}{z} \bigg)Y'_{Q(z)}(v, x_0)
\]
are both $\tilde{A}$-homogeneous as operators.  By comparing the
$\tilde{A}$-homogeneous components of both sides of (\ref{cpb-q}), we
see that the $\tilde{A}$-homogeneous components of $\lambda$ also
satisfy the $Q(z)$-compatibility condition.  \epfv

\begin{rema}\label{q-stableundercomponentops}
{\rm Just as in Remark \ref{stableundercomponentops}, note that both
the spaces $((W_1\otimes W_2)^*)_{[{\mathbb C}];Q(z)}^{( \tilde A )}$
and $((W_1\otimes W_2)^*)_{({\mathbb C});Q(z)}^{( \tilde A )}$ are
stable under the component operators $\tau_{Q(z)}(v\otimes t^m)$ of
the operators $Y'_{Q(z)}(v,x)$ for $v\in V$, $m\in {\mathbb Z}$, and
under the operators $L'_{Q(z)}(-1)$, $L'_{Q(z)}(0)$ and
$L'_{Q(z)}(1)$; this uses Proposition \ref{tau-q-a-comp}, Remark
\ref{L'qjpreservesbetaspace}, Propositions \ref{5.1} and
\ref{qz-comm}, and (\ref{qz-sl-2-qz-y-2}).}
\end{rema}

Again let $(W; I)$ ($W$ a generalized $V$-module) be a $Q(z)$-product
of $W_{1}$ and $W_{2}$ and let $w'\in W'$.  Since $I'$ in particular
intertwines the actions of $V\otimes{\mathbb C}[t, t^{-1}]$ and of
$\mathfrak{s}\mathfrak{l}(2)$, and is $\tilde{A}$-compatible, $I'(W')$
is a generalized $V$-module (recall the proof of Proposition
\ref{im-q:abc}).  Thus for every $w'\in W'$, $I'(w')$ also
satisfies the following condition on $\lambda \in (W_1\otimes W_2)^*$:
\begin{description}
\item{\bf The $Q(z)$-local grading restriction condition}

(a) The {\em $Q(z)$-grading condition}: $\lambda$ is a (finite) sum of
generalized eigenvectors  for the operator
$L'_{Q(z)}(0)$ on $(W_1\otimes W_2)^*$ that 
are also homogeneous with respect to $\tilde A$, that is, 
\[
\lambda\in ((W_1\otimes W_2)^*)_{[{\mathbb C}]; Q(z)}^{( \tilde A )}.
\]
\label{q-homo}

(b) Let $W_{\lambda;Q(z)}$ be the smallest doubly graded (or
equivalently, $\tilde A$-graded; recall Remark
\ref{q-singleanddoublegraded}) subspace of $((W_1\otimes
W_2)^*)_{[{\mathbb C}];Q(z)}^{(\tilde A )}$ containing $\lambda$ and
stable under the component operators $\tau_{Q(z)}(v\otimes t^m)$ of
the operators $Y'_{Q(z)}(v,x)$ for $v\in V$, $m\in {\mathbb Z}$, and
under the operators $L'_{Q(z)}(-1)$, $L'_{Q(z)}(0)$ and
$L'_{Q(z)}(1)$.  (In view of Remark \ref{q-stableundercomponentops},
$W_{\lambda;Q(z)}$ indeed exists.)  Then $W_{\lambda; Q(z)}$ has the
properties
\begin{eqnarray}
&\dim(W_{\lambda})^{(\beta)}_{[n];Q(z)}<\infty,&\label{q-lgrc1}\\
&(W_{\lambda})^{(\beta)}_{[n+k];Q(z)}=0\;\;\mbox{ for }\;k\in {\mathbb Z}
\;\mbox{ sufficiently negative},&\label{q-lgrc2}
\end{eqnarray}
for any $n\in {\mathbb C}$ and $\beta\in \tilde A$, where as usual the
subscripts denote the ${\mathbb C}$-grading and the superscripts
denote the $\tilde A$-grading.
\end{description}

In the case that $W$ is an (ordinary) $V$-module and $w'\in W'$,
$I'(w')$ also satisfies the following $L(0)$-semisimple version of
this condition on $\lambda \in (W_1\otimes W_2)^*$:

\begin{description}
\item{\bf The $L(0)$-semisimple $Q(z)$-local grading restriction condition}

(a) The {\em $L(0)$-semisimple 
$Q(z)$-grading condition}: $\lambda$ is a (finite) sum of
eigenvectors  for the operator
$L'_{Q(z)}(0)$ on $(W_1\otimes W_2)^*$ that 
are also homogeneous with respect to $\tilde A$, that is,
\[
\lambda\in ((W_1\otimes W_2)^*)_{({\mathbb C});Q(z)}^{( \tilde A )}.
\]
\label{q-semi-homo}

(b) Consider $W_{\lambda;Q(z)}$ as above, which in this case is in fact the
smallest doubly graded (or equivalently, $\tilde A$-graded) subspace
of $((W_1\otimes W_2)^*)_{({\mathbb C});Q(z)}^{( \tilde A )}$ containing
$\lambda$ and stable under the component operators
$\tau_{Q(z)}(v\otimes t^m)$ of the operators $Y'_{Q(z)}(v,x)$ for
$v\in V$, $m\in {\mathbb Z}$, and under the operators $L'_{Q(z)}(-1)$,
$L'_{Q(z)}(0)$ and $L'_{Q(z)}(1)$.  Then $W_{\lambda;Q(z)}$ has the
properties
\begin{eqnarray}
&\dim(W_{\lambda;Q(z)})^{(\beta)}_{(n);Q(z)}<\infty,&\label{q-semi-lgrc1}\\
&(W_\lambda)^{(\beta)}_{(n+k);Q(z)}=0\;\;\mbox{ for }\;k\in {\mathbb Z}
\;\mbox{ sufficiently negative},&\label{q-semi-lgrc2}
\end{eqnarray}
for any $n\in {\mathbb C}$ and $\beta\in \tilde A$, where  the
subscripts denote the ${\mathbb C}$-grading and the superscripts
denote the $\tilde A$-grading.
\end{description}

\begin{nota}
{\rm Note that the set of elements of $(W_1\otimes W_2)^*$ satisfying
either of these two $Q(z)$-local grading restriction conditions, or
either of the Part (a)'s in these conditions, forms a subspace.  We
shall denote the space of elements of $(W_1\otimes W_2)^*$ satisfying
the $Q(z)$-local grading restriction condition and the
$L(0)$-semisimple $Q(z)$-local grading restriction condition by
\[
\lgr_{[\C]; Q(z)}((W_1\otimes W_2)^*)
\]
and 
\[
\lgr_{(\C); Q(z)}((W_1\otimes W_2)^*),
\]
respectively.}
\end{nota}

We have the following important theorems generalizing 
the corresponding results stated in \cite{tensor1} and 
proved in \cite{tensor2}. The proofs of these theorems
will be given in
the next section.

\begin{theo}\label{6.1}
Let $\lambda$ be an element of $(W_1\otimes W_2)^{*}$ satisfying the
$Q(z)$-compatibility condition. Then when acting on $\lambda$, the Jacobi
identity for $Y'_{Q(z)}$ holds, that is,
\begin{eqnarray}
\lefteqn{x_0^{-1}\delta
\left({\displaystyle\frac{x_1-x_2}{x_0}}\right)Y'_{Q(z)}(u, x_1)
Y'_{Q(z)}(v, x_2)\lambda}\nno\\
&&\hspace{2ex}-x_0^{-1} \delta
\left({\displaystyle\frac{x_2-x_1}{-x_0}}\right)Y'_{Q(z)}(v, x_2)
Y'_{Q(z)}(u, x_1)\lambda\nonumber \\
&&=x_2^{-1} \delta
\left({\displaystyle\frac{x_1-x_0}{x_2}}\right)Y'_{Q(z)}(Y(u, x_0)v,
x_2)\lambda
\end{eqnarray}
for $u, v\in V$.
\end{theo}

\begin{theo}\label{6.2}
The subspace $\comp_{Q(z)}((W_1\otimes W_2)^*)$ of $(W_1\otimes W_2)^{*}$
is stable under the operators
$\tau_{Q(z)}(v\otimes t^{n})$ for $v\in V$ and $n\in {\mathbb Z}$, 
and in the M\"obius case,
also under the operators $L'_{Q(z)}(-1)$, $L'_{Q(z)}(0)$ and
$L'_{Q(z)}(1)$;
similarly for the  subspaces $\lgr_{[\C];
Q(z)}((W_1\otimes W_2)^*$ and $\lgr_{(\C); Q(z)}((W_1\otimes W_2)^*$.
\end{theo}

We have:

\begin{theo}\label{q-wk-mod}
The space $\comp_{Q(z)}((W_1\otimes W_2)^*)$, equipped with the vertex
operator map $Y'_{Q(z)}$ and, in case $V$ is M\"{o}bius, also equipped
with the operators $L_{Q(z)}'(-1)$, $L_{Q(z)}'(0)$ and $L_{Q(z)}'(1)$,
is a weak $V$-module; similarly for the spaces
\[
(\comp_{Q(z)}((W_1\otimes W_2)^*))\cap 
(\lgr_{[\C]; Q(z)}((W_1\otimes W_2)^*))
\]
and 
\[
(\comp_{Q(z)}((W_1\otimes W_2)^*))\cap 
(\lgr_{(\C); Q(z)}((W_1\otimes W_2)^*)).
\]
\end{theo}
\pf By Theorem \ref{6.2}, $Y'_{Q(z)}$ is a map {}from the tensor product
of $V$ with any of these three subspaces to the space of formal
Laurent series with elements of the subspace as coefficients.  By
Proposition \ref{5.1} and Theorem \ref{6.1} and, in the case that $V$
is M\"{o}bius, also by Propositions \ref{q-sl-2} and
\ref{qz-l-y-comm}, we see that all the axioms for weak $V$-module are
satisfied.  \epfv

Moreover, we have the following consequence of Theorem \ref{q-wk-mod}
and Lemma \ref{q-a-tilde-comp}, just as in Theorem \ref{generation}:

\begin{theo}\label{q-generation}
Let 
\[
\lambda\in 
\comp_{Q(z)}((W_1\otimes W_2)^*)\cap \lgr_{[\C]; Q(z)}((W_1\otimes W_2)^*).
\] 
Then $W_{\lambda;Q(z)}$ (recall Part (b) of the $Q(z)$-local grading
restriction condition) equipped with the vertex operator map
$Y_{Q(z)}'$ and, in case $V$ is M\"{o}bius, also equipped with the
operators $L'_{Q(z)}(-1)$, $L'_{Q(z)}(0)$ and $L'_{Q(z)}(1)$, is a
(strongly-graded) generalized $V$-module. If in addition
\[
\lambda\in 
\comp_{Q(z)}((W_1\otimes W_2)^*)\cap \lgr_{(\C); Q(z)}((W_1\otimes W_2)^*),
\] 
that is, $\lambda$ is a sum of eigenvectors of $L_{Q(z)}'(0)$, then
$W_{\lambda;Q(z)}$ ($\subset ((W_{1}\otimes
W_{2})^{*})_{(\C);Q(z)}^{(\tilde{A})}$) is a (strongly-graded) $V$-module.
\epfv
\end{theo}

Finally, as in Theorem \ref{characterizationofbackslash}, we can give
an alternative description of $W_1\hboxtr_{Q(z)} W_2$ by
characterizing the elements of $W_1\hboxtr_{Q(z)} W_2$ using the
$Q(z)$-compatibility condition and the $Q(z)$-local grading
restriction conditions, generalizing Theorem 6.3 in
\cite{tensor1}. The proof of the following theorem is the same as that
of Theorem \ref{characterizationofbackslash}.

\begin{theo}\label{q-characterizationofbackslash}
Suppose that for every element 
\[
\lambda
\in \comp_{Q(z)}((W_1\otimes W_2)^*)\cap \lgr_{[\C]; Q(z)}((W_1\otimes W_2)^*)
\]
the (strongly-graded) generalized module $W_{\lambda;Q(z)}$ given in
Theorem \ref{q-generation} is an object of ${\cal C}$ (this of course
holds in particular if ${\cal C}$ is ${\cal GM}_{sg}$). Then
\[
W_1\hboxtr_{Q(z)}W_2=\comp_{Q(z)}((W_1\otimes W_2)^*)
\cap \lgr_{[\C]; Q(z)}((W_1\otimes W_2)^*).
\]
Suppose that ${\cal C}$ is a category of strongly-graded $V$-modules
(that is, ${\cal C}\subset {\cal M}_{sg}$) and that for every element
\[
\lambda
\in \comp_{Q(z)}((W_1\otimes W_2)^*)\cap \lgr_{(\C); Q(z)}((W_1\otimes W_2)^*)
\]
the (strongly-graded) $V$-module $W_{\lambda;Q(z)}$ given in Theorem
\ref{q-generation} is an object of ${\cal C}$ (which of course holds in
particular if ${\cal C}$ is ${\cal M}_{sg}$).  Then
\[
W_1\hboxtr_{Q(z)}W_2=\comp_{Q(z)}((W_1\otimes W_2)^*)
\cap \lgr_{(\C); Q(z)}((W_1\otimes W_2)^*). \quad\quad\quad\quad
 \square
\]
\end{theo}

\newpage

\setcounter{equation}{0}
\setcounter{rema}{0}

\section{Proof of the theorems used in the constructions}

The primary goal of this section is to prove Theorems \ref{comp=>jcb},
\ref{stable}, \ref{6.1} and \ref{6.2}.  In Subsection 6.1 we prove
Theorems \ref{comp=>jcb} and \ref{stable}, and in Subsection 6.2,
Theorems \ref{6.1} and \ref{6.2}.  The proofs in Subsection 6.1 are
new, even for the category of (ordinary) modules for a vertex operator
algebra satisfying the finiteness and reductivity conditions treated
in \cite{tensor1}--\cite{tensor3}. In \cite{tensor1}--\cite{tensor3},
for a vertex operator algebra satisfying these conditions, Theorems
\ref{6.1} and \ref{6.2}, in the $Q(z)$ case, were proved first, and
then Theorems \ref{comp=>jcb} and \ref{stable}, in the $P(z)$ case,
were proved using results {}from the $Q(z^{-1})$ case and relations
between $P(z)$-tensor products and $Q(z^{-1})$-tensor products. In
Subsection 6.1, we prove Theorems \ref{comp=>jcb} and \ref{stable}
directly, without using any results {}from the $Q(z)$ case.  As usual,
the reader should observe the justifiability of each step in the
arguments (the well-definedness of the formal series, etc.); again as
usual, this is sometimes quite subtle.

\subsection{Proofs of Theorems \ref{comp=>jcb} and
\ref{stable}}

We first prove a formula for vertex operators that will be needed in
the proofs of both Theorem \ref{comp=>jcb} and Theorem \ref{stable}.

\begin{lemma}
For $u, v \in V$, we have 
\begin{eqnarray}\label{comp=>jcb-9}
&{\displaystyle x_{2}^{-1}\delta\left(\frac{x_{1}-x_{0}}{x_{2}}\right)
Y(e^{x_{1}L(1)}
(-x_{1}^{-2})^{L(0)}u, -x_{0}x_{1}^{-1}x_{2}^{-1})e^{x_{2}L(1)}
(-x_{2}^{-2})^{L(0)}}v &\nn
&{\displaystyle =x_{2}^{-1}\delta\left(\frac{x_{1}-x_{0}}{x_{2}}\right)
e^{x_{2}L(1)}(-x_{2}^{-2})^{L(0)}
Y(u, x_{0})}v. &
\end{eqnarray}
\end{lemma}
\pf Using (\ref{log:p2}), (\ref{log:p3}), (\ref{xe^Lx}) and
(\ref{deltafunctionsubstitutionformula}), we have
\begin{eqnarray*}
\lefteqn{x_{2}^{-1}\delta\left(\frac{x_{1}-x_{0}}{x_{2}}\right)
Y(e^{x_{1}L(1)}
(-x_{1}^{-2})^{L(0)}u, -x_{0}x_{1}^{-1}x_{2}^{-1})e^{x_{2}L(1)}
(-x_{2}^{-2})^{L(0)}}\nn
&&=x_{2}^{-1}\delta\left(\frac{x_{1}-x_{0}}{x_{2}}\right)e^{x_{2}L(1)}
Y(e^{-x_{2}(1-x_{0}x_{1}^{-1})L(1)}(1-x_{0}x_{1}^{-1})^{-2L(0)}\cdot\nn
&&\quad\quad\quad\quad\quad \cdot
e^{x_{1}L(1)}
(-x_{1}^{-2})^{L(0)}u, -x_{0}x_{1}^{-1}x_{2}^{-1}(1-x_{0} x_{1}^{-1})^{-1})
(-x_{2}^{-2})^{L(0)}\nn
&&=x_{2}^{-1}\delta\left(\frac{x_{2}^{-1}-x_{0}x^{-1}_{1}
x_{2}^{-1}}{x^{-1}_{1}}\right)
e^{x_{2}L(1)}(-x_{2}^{-2})^{L(0)}\cdot\nn
&&\quad\quad\cdot Y((-x_{2}^{2})^{L(0)}
e^{-x_{2}(1-x_{0}x_{1}^{-1})L(1)}(1-x_{0}x_{1}^{-1})^{-2L(0)}\cdot\nn
&&\quad\quad\quad\quad\quad\quad\quad\quad\quad\quad \cdot
e^{x_{1}L(1)}
(-x_{1}^{-2})^{L(0)}u, x_{0}x_{1}^{-1}
(x_{2}^{-1}-x_{0} x_{1}^{-1}x_{2}^{-1})^{-1})\nn
&&=x_{2}^{-1}\delta\left(\frac{x_{2}^{-1}-x_{0}x^{-1}_{1}
x_{2}^{-1}}{x^{-1}_{1}}\right)e^{x_{2}L(1)}(-x_{2}^{-2})^{L(0)}\cdot\nn
&&\quad\quad\cdot
Y(e^{x_{2}^{-1}(1-x_{0}x_{1}^{-1})L(1)}
(-x_{2}^{2})^{L(0)}
(1-x_{0}x_{1}^{-1})^{-2L(0)}
e^{x_{1}L(1)}
(-x_{1}^{-2})^{L(0)}u, x_{0})\nn
&&=x_{2}^{-1}\delta\left(\frac{x_{1}-x_{0}}{x_{2}}\right)
e^{x_{2}L(1)}(-x_{2}^{-2})^{L(0)}\cdot\nn
&&\quad\quad \cdot
Y(e^{x_{2}^{-1}(1-x_{0}x_{1}^{-1})L(1)}
(-(x_{2}^{-1}(1-x_{0}x_{1}^{-1}))^{-2})^{L(0)}
e^{x_{1}L(1)}
(-x_{1}^{-2})^{L(0)}u, x_{0})\nn
&&=x_{2}^{-1}\delta\left(\frac{x_{1}-x_{0}}{x_{2}}\right)
e^{x_{2}L(1)}(-x_{2}^{-2})^{L(0)}\cdot\nn
&&\quad\quad \cdot
Y(e^{x_{2}^{-1}(1-x_{0}x_{1}^{-1})L(1)}
e^{-x_{1}x_{2}^{-2}(1-x_{0}x_{1}^{-1})^{2}L(1)}
(-(x_{2}^{-1}(1-x_{0}x_{1}^{-1}))^{-2})^{L(0)}
(-x_{1}^{-2})^{L(0)}u, x_{0})\nn
&&=x_{2}^{-1}\delta\left(\frac{x_{1}-x_{0}}{x_{2}}\right)
e^{x_{2}L(1)}(-x_{2}^{-2})^{L(0)}\cdot\nn
&&\quad\quad \cdot
Y(e^{x_{2}^{-1}(1-x_{0}x_{1}^{-1})
(1-x_{2}^{-1}(x_{1}-x_{0}))L(1)}
(x_{2}^{-1}(x_{1}-x_{0}))^{-2L(0)}
u, x_{0})\nn
&&=x_{2}^{-1}\delta\left(\frac{x_{1}-x_{0}}{x_{2}}\right)
e^{x_{2}L(1)}(-x_{2}^{-2})^{L(0)}
Y(u, x_{0}).
\end{eqnarray*}
\epfv

\noindent {\it Proof of Theorem \ref{comp=>jcb}} Let $\lambda$ be an
element of $(W_{1}\otimes W_{2})^{*}$ satisfying the
$P(z)$-compatibility condition, that is, satisfying (a) the
$P(z)$-lower truncation condition---for all $v\in V$, the formal
Laurent series $Y'_{P(z)}(v, x)\lambda$ involves only finitely many
negative powers of $x$, and (b) formula (\ref{cpb}) for all $v\in V$.

For $u, v\in V$, $w_{(1)}\in W_{1}$ and $w_{(2)}\in W_{2}$, by 
definition, 
\begin{eqnarray}\label{comp=>jcb-1}
\lefteqn{\left(x_{0}^{-1}\delta\left(\frac{x_{1}-x_{2}}{x_{0}}\right)
Y'_{P(z)}(u, x_{1})Y'_{P(z)}(v, x_{2})\lambda\right)(w_{(1)}
\otimes w_{(2)})}\nn
&&=x_{0}^{-1}\delta\left(\frac{x_{1}-x_{2}}{x_{0}}\right)
\Biggl((Y'_{P(z)}(v, x_{2})\lambda)(w_{(1)}\otimes 
Y_{2}^{o}(u, x_{1})w_{(2)})\nn
&&\quad\quad +\res_{y_{1}}z^{-1}
\delta\left(\frac{x_{1}^{-1}-y_{1}}{z}\right)(Y'_{P(z)}(v, x_{2})\lambda)
(Y_{1}(e^{x_{1}L(1)}(-x_{1}^{-2})^{L(0)}u, y_{1})w_{(1)}\otimes w_{(2)})
\Biggr)\nn
&&=x_{0}^{-1}\delta\left(\frac{x_{1}-x_{2}}{x_{0}}\right)
\Biggl(\lambda(w_{(1)}\otimes 
Y_{2}^{o}(v, x_{2})Y_{2}^{o}(u, x_{1})w_{(2)})\nn
&&\quad\quad +\res_{y_{2}}z^{-1}
\delta\left(\frac{x_{2}^{-1}-y_{2}}{z}\right)\lambda(
Y_{1}(e^{x_{2}L(1)}(-x_{2}^{-2})^{L(0)}v, y_{2})w_{(1)}\otimes 
Y_{2}^{o}(u, x_{1})w_{(2)})\nn
&&\quad\quad +\res_{y_{1}}z^{-1}
\delta\left(\frac{x_{1}^{-1}-y_{1}}{z}\right)(Y'_{P(z)}(v, x_{2})\lambda)
(Y_{1}(e^{x_{1}L(1)}(-x_{1}^{-2})^{L(0)}u, y_{1})w_{(1)}\otimes w_{(2)})
\Biggr)\nn
&&=x_{0}^{-1}\delta\left(\frac{x_{1}-x_{2}}{x_{0}}\right)
\lambda(w_{(1)}\otimes 
Y_{2}^{o}(v, x_{2})Y_{2}^{o}(u, x_{1})w_{(2)})\nn
&&\quad\quad +x_{0}^{-1}\delta\left(\frac{x_{1}-x_{2}}{x_{0}}\right)
\res_{y_{2}}z^{-1}
\delta\left(\frac{x_{2}^{-1}-y_{2}}{z}\right)\cdot\nn
&&\quad\quad\quad\quad\quad \cdot\lambda(
Y_{1}(e^{x_{2}L(1)}(-x_{2}^{-2})^{L(0)}v, y_{2})w_{(1)}\otimes 
Y_{2}^{o}(u, x_{1})w_{(2)})\nn
&&\quad\quad +x_{0}^{-1}\delta\left(\frac{x_{1}-x_{2}}{x_{0}}\right)
\res_{y_{1}}z^{-1}
\delta\left(\frac{x_{1}^{-1}-y_{1}}{z}\right)\cdot\nn
&&\quad\quad\quad \quad\quad\cdot(Y'_{P(z)}(v, x_{2})\lambda)
(Y_{1}(e^{x_{1}L(1)}(-x_{1}^{-2})^{L(0)}u, y_{1})w_{(1)}\otimes w_{(2)}).
\end{eqnarray}

Using (\ref{2termdeltarelation}) and (\ref{cpb}), we see that the
third term on the right-hand side of (\ref{comp=>jcb-1}) is equal to
\begin{eqnarray}\label{comp=>jcb-2}
\lefteqn{x_{0}^{-1}\delta\left(\frac{x^{-1}_{2}-x^{-1}_{1}}
{x_{0}x_{1}^{-1}x_{2}^{-1}}\right)
\res_{y_{1}}z^{-1}
\delta\left(\frac{x_{1}^{-1}-y_{1}}{z}\right)\cdot}
\nn
&&\quad\quad\quad\cdot (Y'_{P(z)}(v, x_{2})\lambda)
(Y_{1}(e^{x_{1}L(1)}(-x_{1}^{-2})^{L(0)}u, y_{1})w_{(1)}\otimes w_{(2)})\nn
&&=\res_{y_{1}}x_{1}^{-1}x_{2}^{-1}(x_{0}x_{1}^{-1}x_{2}^{-1})^{-1}
\delta\left(\frac{x^{-1}_{2}-y_{1}-z}
{x_{0}x_{1}^{-1}x_{2}^{-1}}\right)
z^{-1}
\delta\left(\frac{x_{1}^{-1}-y_{1}}{z}\right)\cdot
\nn
&&\quad\quad\quad\cdot (Y'_{P(z)}(v, x_{2})\lambda)
(Y_{1}(e^{x_{1}L(1)}(-x_{1}^{-2})^{L(0)}u, y_{1})w_{(1)}\otimes w_{(2)})\nn
&&=\res_{y_{1}}x_{1}^{-1}x_{2}^{-1}(x_{0}x_{1}^{-1}x_{2}^{-1}+y_{1})^{-1}
\delta\left(\frac{x^{-1}_{2}-z}
{x_{0}x_{1}^{-1}x_{2}^{-1}+y_{1}}\right)z^{-1}
\delta\left(\frac{x_{1}^{-1}-y_{1}}{z}\right)\cdot
\nn
&&\quad\quad\quad\cdot (Y'_{P(z)}(v, x_{2})\lambda)
(Y_{1}(e^{x_{1}L(1)}(-x_{1}^{-2})^{L(0)}u, y_{1})w_{(1)}\otimes w_{(2)})\nn
&&=\res_{y_{1}}x_{1}^{-1}x_{2}^{-1}z^{-1}
\delta\left(\frac{x_{1}^{-1}-y_{1}}{z}\right)\cdot
\nn
&&\quad\quad\quad\cdot 
\left(\tau_{P(z)}\left((x_{0}x_{1}^{-1}x_{2}^{-1}+y_{1})^{-1}
\delta\left(\frac{x^{-1}_{2}-z}
{x_{0}x_{1}^{-1}x_{2}^{-1}+y_{1}}\right)Y_{t}(v, x_{2})\right)
\lambda\right)\nn
&&\quad\quad\quad\quad\quad\quad(Y_{1}(e^{x_{1}L(1)}
(-x_{1}^{-2})^{L(0)}u, y_{1})w_{(1)}\otimes w_{(2)})\nn
&&=\res_{y_{1}}x_{1}^{-1}x_{2}^{-1}
z^{-1}
\delta\left(\frac{x_{1}^{-1}-y_{1}}{z}\right)\cdot
\nn
&&\quad\quad\cdot \Biggl(z^{-1}
\delta\left(\frac{x_{2}^{-1}-x_{0}x_{1}^{-1}x_{2}^{-1}-y_{1}}{z}\right)
\cdot\nn
&&\quad\quad\quad\; \cdot\lambda(Y_{1}(e^{x_{2}L(1)}
(-x_{2}^{-2})^{L(0)}v, x_{0}x_{1}^{-1}x_{2}^{-1}+y_{1})Y_{1}(e^{x_{1}L(1)}
(-x_{1}^{-2})^{L(0)}u, y_{1})w_{(1)}\otimes w_{(2)})\nn
&&\quad\quad\;\; +
(x_{0}x_{1}^{-1}x_{2}^{-1}+y_{1})^{-1}
\delta\left(\frac{z-x_{2}^{-1}}{-x_{0}x_{1}^{-1}x_{2}^{-1}-y_{1}}\right)
\cdot\nn
&&\quad\quad\quad \;\cdot\lambda(Y_{1}(e^{x_{1}L(1)}
(-x_{1}^{-2})^{L(0)}u, y_{1})w_{(1)}\otimes Y_{2}^{o}(v, x_{2})
w_{(2)})\Biggr)\nn
&&=\res_{y_{1}}x_{2}^{-1}
\delta\left(\frac{z+y_{1}}{x_{1}^{-1}}\right) (z+y_{1})^{-1}
\delta\left(\frac{x_{2}^{-1}-x_{0}x_{1}^{-1}x_{2}^{-1}}{z+y_{1}}\right)
\cdot\nn
&&\quad\quad\quad\; \cdot\lambda(Y_{1}(e^{x_{2}L(1)}
(-x_{2}^{-2})^{L(0)}v, x_{0}x_{1}^{-1}x_{2}^{-1}+y_{1})Y_{1}(e^{x_{1}L(1)}
(-x_{1}^{-2})^{L(0)}u, y_{1})w_{(1)}\otimes w_{(2)})\nn
&&\quad\quad +\res_{y_{1}}x_{2}^{-1}
\delta\left(\frac{z+y_{1}}{x_{1}^{-1}}\right)
(x_{0}x_{1}^{-1}x_{2}^{-1})^{-1}
\delta\left(\frac{z+y_{1}-x_{2}^{-1}}{-x_{0}x_{1}^{-1}x_{2}^{-1}}\right)
\cdot\nn
&&\quad\quad\quad\; \cdot\lambda(Y_{1}(e^{x_{1}L(1)}
(-x_{1}^{-2})^{L(0)}u, y_{1})w_{(1)}\otimes Y_{2}^{o}(v, x_{2})
w_{(2)}).
\end{eqnarray}
By (\ref{deltafunctionsubstitutionformula}) and
(\ref{2termdeltarelation}), the right-hand side of (\ref{comp=>jcb-2})
is equal to
\begin{eqnarray}\label{comp=>jcb-3}
\lefteqn{\res_{y_{1}}x_{2}^{-1}
\delta\left(\frac{z+y_{1}}{x_{1}^{-1}}\right)x_{1}
\delta\left(\frac{x_{2}^{-1}-x_{0}x_{1}^{-1}x_{2}^{-1}}{x_{1}^{-1}}\right)
\cdot}\nn
&&\quad\quad \cdot\lambda(Y_{1}(e^{x_{2}L(1)}
(-x_{2}^{-2})^{L(0)}v, x_{0}x_{1}^{-1}x_{2}^{-1}+y_{1})Y_{1}(e^{x_{1}L(1)}
(-x_{1}^{-2})^{L(0)}u, y_{1})w_{(1)}\otimes w_{(2)})\nn
&&\quad +\res_{y_{1}}x_{2}^{-1}
\delta\left(\frac{z+y_{1}}{x_{1}^{-1}}\right)
(x_{0}x_{1}^{-1}x_{2}^{-1})^{-1}
\delta\left(\frac{x_{1}^{-1}-x_{2}^{-1}}{-x_{0}x_{1}^{-1}x_{2}^{-1}}\right)
\cdot\nn
&&\quad\quad \cdot\lambda(Y_{1}(e^{x_{1}L(1)}
(-x_{1}^{-2})^{L(0)}u, y_{1})w_{(1)}\otimes Y_{2}^{o}(v, x_{2})
w_{(2)})\nn
&&=\res_{y_{1}}z^{-1}
\delta\left(\frac{x_{1}^{-1}-y_{1}}{z}\right)x_{2}^{-1}
\delta\left(\frac{x_{1}-x_{0}}{x_{2}}\right)
\cdot\nn
&&\quad\quad \cdot\lambda(Y_{1}(e^{x_{2}L(1)}
(-x_{2}^{-2})^{L(0)}v, x_{0}x_{1}^{-1}x_{2}^{-1}+y_{1})Y_{1}(e^{x_{1}L(1)}
(-x_{1}^{-2})^{L(0)}u, y_{1})w_{(1)}\otimes w_{(2)})\nn
&&\quad +\res_{y_{1}}z^{-1}
\delta\left(\frac{x_{1}^{-1}-y_{1}}{z}\right)
x_{0}^{-1}
\delta\left(\frac{x_{2}-x_{1}}{-x_{0}}\right)
\cdot\nn
&&\quad\quad \cdot\lambda(Y_{1}(e^{x_{1}L(1)}
(-x_{1}^{-2})^{L(0)}u, y_{1})w_{(1)}\otimes Y_{2}^{o}(v, x_{2})
w_{(2)}).
\end{eqnarray}

Since 
\[
\res_{y_{2}}y_{2}^{-1}
\delta\left(\frac{x_{0}x_{1}^{-1}x_{2}^{-1}+y_{1}}{y_{2}}\right)=1,
\]
the first term on the right-hand side of (\ref{comp=>jcb-3}) 
can be written as
\begin{eqnarray}\label{comp=>jcb-4}
\lefteqn{\res_{y_{1}}z^{-1}
\delta\left(\frac{x_{1}^{-1}-y_{1}}{z}\right)
x_{2}^{-1}
\delta\left(\frac{x_{1}-x_{0}}{x_{2}}\right)
\res_{y_{2}}y_{2}^{-1}
\delta\left(\frac{x_{0}x_{1}^{-1}x_{2}^{-1}+y_{1}}{y_{2}}\right)
\cdot }\nn
&&\quad\quad\quad \cdot\lambda(Y_{1}(e^{x_{2}L(1)}
(-x_{2}^{-2})^{L(0)}v, x_{0}x_{1}^{-1}x_{2}^{-1}+y_{1})Y_{1}(e^{x_{1}L(1)}
(-x_{1}^{-2})^{L(0)}u, y_{1})w_{(1)}\otimes w_{(2)})\nn
&&=x_{2}^{-1}
\delta\left(\frac{x_{1}-x_{0}}{x_{2}}\right)\res_{y_{1}}
\res_{y_{2}}z^{-1}
\delta\left(\frac{x_{1}^{-1}-y_{1}}{z}\right)
(x_{0}x_{1}^{-1}x_{2}^{-1})^{-1}
\delta\left(\frac{y_{2}-y_{1}}{x_{0}x_{1}^{-1}x_{2}^{-1}}\right)
\cdot \nn
&&\quad\quad\quad \cdot\lambda(Y_{1}(e^{x_{2}L(1)}
(-x_{2}^{-2})^{L(0)}v, y_{2})Y_{1}(e^{x_{1}L(1)}
(-x_{1}^{-2})^{L(0)}u, y_{1})w_{(1)}\otimes w_{(2)}),
\end{eqnarray}
where we have also used (\ref{deltafunctionsubstitutionformula}) and
(\ref{2termdeltarelation}).  Again using (\ref{2termdeltarelation})
and (\ref{deltafunctionsubstitutionformula}), we see that the
right-hand side of (\ref{comp=>jcb-4}) is also equal to
\begin{eqnarray}\label{comp=>jcb-5}
\lefteqn{x_{2}^{-1}
\delta\left(\frac{x_{2}^{-1}-x_{0}x_{1}^{-1}x_{2}^{-1}}
{x_{1}^{-1}}\right)\res_{y_{1}}
\res_{y_{2}}z^{-1}
\delta\left(\frac{x_{1}^{-1}-y_{1}}{z}\right)
y_{2}^{-1}
\delta\left(\frac{x_{0}x_{1}^{-1}x_{2}^{-1}+y_{1}}{y_{2}}\right)
\cdot} \nn
&&\quad\quad\quad\quad \cdot\lambda(Y_{1}(e^{x_{2}L(1)}
(-x_{2}^{-2})^{L(0)}v, y_{2})Y_{1}(e^{x_{1}L(1)}
(-x_{1}^{-2})^{L(0)}u, y_{1})w_{(1)}\otimes w_{(2)})\nn
&&=x_{2}^{-1}
\delta\left(\frac{x_{2}^{-1}-x_{0}x_{1}^{-1}x_{2}^{-1}}
{x_{1}^{-1}}\right)\res_{y_{1}}
\res_{y_{2}}z^{-1}
\delta\left(\frac{x_{2}^{-1}-x_{0}x_{1}^{-1}x_{2}^{-1}-y_{1}}{z}\right)
\cdot\nn
&&\quad\quad\quad\quad \cdot y_{2}^{-1}
\delta\left(\frac{x_{0}x_{1}^{-1}x_{2}^{-1}+y_{1}}{y_{2}}\right)
\cdot \nn
&&\quad\quad\quad\quad \cdot\lambda(Y_{1}(e^{x_{2}L(1)}
(-x_{2}^{-2})^{L(0)}v, y_{2})Y_{1}(e^{x_{1}L(1)}
(-x_{1}^{-2})^{L(0)}u, y_{1})w_{(1)}\otimes w_{(2)})\nn
&&=x_{2}^{-1}
\delta\left(\frac{x_{1}-x_{0}}
{x_{2}}\right)\res_{y_{1}}
\res_{y_{2}}z^{-1}
\delta\left(\frac{x_{2}^{-1}-y_{2}}{z}\right)
y_{2}^{-1}
\delta\left(\frac{x_{0}x_{1}^{-1}x_{2}^{-1}+y_{1}}{y_{2}}\right)
\cdot \nn
&&\quad\quad\quad\quad \cdot\lambda(Y_{1}(e^{x_{2}L(1)}
(-x_{2}^{-2})^{L(0)}v, y_{2})Y_{1}(e^{x_{1}L(1)}
(-x_{1}^{-2})^{L(0)}u, y_{1})w_{(1)}\otimes w_{(2)})\nn
&&=x_{2}^{-1}
\delta\left(\frac{x_{1}-x_{0}}
{x_{2}}\right)\res_{y_{1}}
\res_{y_{2}}z^{-1}
\delta\left(\frac{x_{2}^{-1}-y_{2}}{z}\right)
(x_{0}x_{1}^{-1}x_{2}^{-1})^{-1}
\delta\left(\frac{y_{2}-y_{1}}{x_{0}x_{1}^{-1}x_{2}^{-1}}\right)
\cdot \nn
&&\quad\quad\quad\quad \cdot\lambda(Y_{1}(e^{x_{2}L(1)}
(-x_{2}^{-2})^{L(0)}v, y_{2})Y_{1}(e^{x_{1}L(1)}
(-x_{1}^{-2})^{L(0)}u, y_{1})w_{(1)}\otimes w_{(2)}).
\end{eqnarray}
That is, in the middle delta-function expression in the right-hand
side of (\ref{comp=>jcb-4}), we may replace $x_1$ by $x_2$ and $y_1$
by $y_2$.

{}From (\ref{comp=>jcb-1})--(\ref{comp=>jcb-5}), we obtain
\begin{eqnarray}\label{comp=>jcb-6}
\lefteqn{\left(x_{0}^{-1}\delta\left(\frac{x_{1}-x_{2}}{x_{0}}\right)
Y'_{P(z)}(u, x_{1})Y'_{P(z)}(v, x_{2})\lambda\right)(w_{(1)}
\otimes w_{(2)})}\nn
&&=x_{0}^{-1}\delta\left(\frac{x_{1}-x_{2}}{x_{0}}\right)
\lambda(w_{(1)}\otimes 
Y_{2}^{o}(v, x_{2})Y_{2}^{o}(u, x_{1})w_{(2)})\nn
&&\quad +x_{0}^{-1}\delta\left(\frac{x_{1}-x_{2}}{x_{0}}\right)
\res_{y_{2}}z^{-1}
\delta\left(\frac{x_{2}^{-1}-y_{2}}{z}\right)\cdot \nn
&&\quad\quad\quad\quad \cdot\lambda(
Y_{1}(e^{x_{2}L(1)}(-x_{2}^{-2})^{L(0)}v, y_{2})w_{(1)}\otimes 
Y_{2}^{o}(u, x_{1})w_{(2)})\nn
&&\quad +x_{2}^{-1}
\delta\left(\frac{x_{1}-x_{0}}{x_{2}}\right)\res_{y_{1}}
\res_{y_{2}}z^{-1}
\delta\left(\frac{x_{2}^{-1}-y_{2}}{z}\right)
(x_{0}x_{1}^{-1}x_{2}^{-1})^{-1}  
\delta\left(\frac{y_{2}-y_{1}}{x_{0}x_{1}^{-1}x_{2}^{-1}}\right)
\cdot \nn
&&\quad\quad\quad\quad \cdot\lambda(Y_{1}(e^{x_{2}L(1)}
(-x_{2}^{-2})^{L(0)}v, y_{2})Y_{1}(e^{x_{1}L(1)}
(-x_{1}^{-2})^{L(0)}u, y_{1})w_{(1)}\otimes w_{(2)})\nn
&&\quad +x_{0}^{-1}
\delta\left(\frac{x_{2}-x_{1}}{-x_{0}}\right)\res_{y_{1}}
z^{-1}
\delta\left(\frac{x_{1}^{-1}-y_{1}}{z}\right)
\cdot\nn
&&\quad\quad\quad\quad \cdot\lambda(Y_{1}(e^{x_{1}L(1)}
(-x_{1}^{-2})^{L(0)}u, y_{1})w_{(1)}\otimes Y_{2}^{o}(v, x_{2})
w_{(2)})\Biggr).
\end{eqnarray}
{}From (\ref{comp=>jcb-5}) and (\ref{comp=>jcb-6}), replacing $u, v,
x_{1}, x_{2}, x_{0}$ by $v, u, x_{2}, x_{1}, -x_{0}$, respectively,
and also using (\ref{2termdeltarelation}), we find that
\begin{eqnarray}\label{comp=>jcb-7}
\lefteqn{\left(-x_{0}^{-1}\delta\left(\frac{x_{2}-x_{1}}{-x_{0}}\right)
Y'_{P(z)}(v, x_{2})Y'_{P(z)}(u, x_{1})\lambda\right)(w_{(1)}
\otimes w_{(2)})}\nn
&&=-x_{0}^{-1}\delta\left(\frac{x_{2}-x_{1}}{-x_{0}}\right)
\lambda(w_{(1)}\otimes 
Y_{2}^{o}(u, x_{1})Y_{2}^{o}(v, x_{2})w_{(2)})\nn
&&\quad -x_{0}^{-1}\delta\left(\frac{x_{2}-x_{1}}{-x_{0}}\right)
\res_{y_{1}}z^{-1}
\delta\left(\frac{x_{1}^{-1}-y_{1}}{z}\right)\cdot \nn
&&\quad\quad\quad\quad \cdot\lambda(
Y_{1}(e^{x_{1}L(1)}(-x_{1}^{-2})^{L(0)}u, y_{1})w_{(1)}\otimes 
Y_{2}^{o}(v, x_{2})w_{(2)})\nn
&&\quad -x_{2}^{-1}
\delta\left(\frac{x_{1}-x_{0}}{x_{2}}\right)\res_{y_{1}}
\res_{y_{2}}z^{-1}
\delta\left(\frac{x_{2}^{-1}-y_{2}}{z}\right)
(x_{0}x_{1}^{-1}x_{2}^{-1})^{-1}  
\delta\left(\frac{y_{1}-y_{2}}{-x_{0}x_{1}^{-1}x_{2}^{-1}}\right)
\cdot \nn
&&\quad\quad\quad\quad \cdot\lambda(Y_{1}(e^{x_{1}L(1)}
(-x_{1}^{-2})^{L(0)}u, y_{1})Y_{1}(e^{x_{2}L(1)}
(-x_{2}^{-2})^{L(0)}v, y_{2})w_{(1)}\otimes w_{(2)})\nn
&&\quad -x_{0}^{-1}
\delta\left(\frac{x_{1}-x_{2}}{x_{0}}\right)\res_{y_{2}}
z^{-1}
\delta\left(\frac{x_{2}^{-1}-y_{2}}{z}\right)
\cdot\nn
&&\quad\quad\quad\quad \cdot\lambda(Y_{1}(e^{x_{2}L(1)}
(-x_{2}^{-2})^{L(0)}v, y_{2})w_{(1)}\otimes Y_{2}^{o}(u, x_{1})
w_{(2)})\Biggr).
\end{eqnarray}

Using (\ref{comp=>jcb-6}), (\ref{comp=>jcb-7}), the Jacobi identity,
the opposite Jacobi identity (\ref{op-jac-id}) and
(\ref{2termdeltarelation}), we obtain
\begin{eqnarray}\label{comp=>jcb-8}
\lefteqn{\Biggl(x_{0}^{-1}\delta\left(\frac{x_{1}-x_{2}}{x_{0}}\right)
Y'_{P(z)}(u, x_{1})Y'_{P(z)}(v, x_{2})\lambda}\nn
&&\quad
-x_{0}^{-1}\delta\left(\frac{x_{2}-x_{1}}{-x_{0}}\right)
Y'_{P(z)}(v, x_{2})Y'_{P(z)}(u, x_{1})\lambda\Biggr)(w_{(1)}
\otimes w_{(2)})\nn
&&=\lambda\Biggl(w_{(1)}\otimes 
\Biggl(x_{0}^{-1}\delta\left(\frac{x_{1}-x_{2}}{x_{0}}\right)
Y_{2}^{o}(v, x_{2})Y_{2}^{o}(u, x_{1})\nn
&&\quad\quad\quad\quad \quad\quad \quad
-x_{0}^{-1}\delta\left(\frac{x_{2}-x_{1}}{-x_{0}}\right)
Y_{2}^{o}(u, x_{1})Y_{2}^{o}(v, x_{2})\Biggr) w_{(2)}\Biggr)\nn
&&\quad +x_{2}^{-1}
\delta\left(\frac{x_{1}-x_{0}}{x_{2}}\right)\res_{y_{1}}
\res_{y_{2}}z^{-1}
\delta\left(\frac{x_{2}^{-1}-y_{2}}{z}\right)
\cdot \nn
&&\quad\quad \cdot\lambda\Biggl(
\Biggl((x_{0}x_{1}^{-1}x_{2}^{-1})^{-1}
\delta\left(\frac{y_{2}-y_{1}}{x_{0}x_{1}^{-1}x_{2}^{-1}}\right)\cdot\nn
&&\quad\quad\quad\quad\quad\quad\quad \cdot Y_{1}(e^{x_{2}L(1)}
(-x_{2}^{-2})^{L(0)}v, y_{2})Y_{1}(e^{x_{1}L(1)}
(-x_{1}^{-2})^{L(0)}u, y_{1})\nn
&&\quad\quad\quad\quad\;\;-(x_{0}x_{1}^{-1}x_{2}^{-1})^{-1}
\delta\left(\frac{y_{1}-y_{2}}{-x_{0}x_{1}^{-1}x_{2}^{-1}}\right)\cdot\nn
&&\quad\quad\quad\quad\quad\quad \quad\cdot
Y_{1}(e^{x_{1}L(1)}
(-x_{1}^{-2})^{L(0)}u, y_{1})Y_{1}(e^{x_{2}L(1)}
(-x_{2}^{-2})^{L(0)}v, y_{2})\Biggr)w_{(1)}\otimes w_{(2)}\Biggr)\nn
&&=\lambda\Biggl(w_{(1)}\otimes 
\Biggl(x_{2}^{-1}\delta\left(\frac{x_{1}-x_{0}}{x_{2}}\right)
Y_{2}^{o}(Y(u, x_{0})v, x_{2})\Biggr) w_{(2)}\Biggr)\nn
&&\quad +x_{2}^{-1}
\delta\left(\frac{x_{1}-x_{0}}{x_{2}}\right)\res_{y_{1}}
\res_{y_{2}}z^{-1}
\delta\left(\frac{x_{2}^{-1}-y_{2}}{z}\right)
\cdot \nn
&&\quad\quad \cdot\lambda\Biggl(
\Biggl(y_{2}^{-1}
\delta\left(\frac{y_{1}+x_{0}x_{1}^{-1}x_{2}^{-1}}{y_{2}}\right)\cdot\nn
&&\quad\quad\quad\quad \cdot Y_{1}(Y(e^{x_{1}L(1)}
(-x_{1}^{-2})^{L(0)}u, -x_{0}x_{1}^{-1}x_{2}^{-1})e^{x_{2}L(1)}
(-x_{2}^{-2})^{L(0)}v, y_{2})\Biggr)w_{(1)}\otimes w_{(2)}\Biggr)\nn
&&=x_{2}^{-1}\delta\left(\frac{x_{1}-x_{0}}{x_{2}}\right)
\lambda(w_{(1)}\otimes Y_{2}^{o}(Y(u, x_{0})v, x_{2}) w_{(2)})\nn
&&\quad +x_{2}^{-1}
\delta\left(\frac{x_{1}-x_{0}}{x_{2}}\right)\res_{y_{1}}
\res_{y_{2}}z^{-1}
\delta\left(\frac{x_{2}^{-1}-y_{2}}{z}\right)
\cdot \nn
&&\quad\quad \cdot\lambda\Biggl(
\Biggl(y_{1}^{-1}
\delta\left(\frac{y_{2}-x_{0}x_{1}^{-1}x_{2}^{-1}}{y_{1}}\right)\cdot\nn
&&\quad\quad\quad\quad \cdot Y_{1}(Y(e^{x_{1}L(1)}
(-x_{1}^{-2})^{L(0)}u, -x_{0}x_{1}^{-1}x_{2}^{-1})e^{x_{2}L(1)}
(-x_{2}^{-2})^{L(0)}v, y_{2})\Biggr)w_{(1)}\otimes w_{(2)}\Biggr)\nn
&&=x_{2}^{-1}\delta\left(\frac{x_{1}-x_{0}}{x_{2}}\right)
\lambda(w_{(1)}\otimes Y_{2}^{o}(Y(u, x_{0})v, x_{2}) w_{(2)})\nn
&&\quad +x_{2}^{-1}
\delta\left(\frac{x_{1}-x_{0}}{x_{2}}\right)
\res_{y_{2}}z^{-1}
\delta\left(\frac{x_{2}^{-1}-y_{2}}{z}\right)
\cdot \nn
&&\quad\quad \cdot\lambda(
Y_{1}(Y(e^{x_{1}L(1)}
(-x_{1}^{-2})^{L(0)}u, -x_{0}x_{1}^{-1}x_{2}^{-1})e^{x_{2}L(1)}
(-x_{2}^{-2})^{L(0)}v, y_{2})w_{(1)}\otimes w_{(2)})\nn
\end{eqnarray}

Finally, {}from (\ref{comp=>jcb-9}) we see that the right-hand side of
(\ref{comp=>jcb-8}) becomes
\begin{eqnarray}\label{comp=>jcb-10}
\lefteqn{x_{2}^{-1}\delta\left(\frac{x_{1}-x_{0}}{x_{2}}\right)
\lambda(w_{(1)}\otimes Y_{2}^{o}(Y(u, x_{0})v, x_{2}) w_{(2)})}\nn
&&\quad +x_{2}^{-1}
\delta\left(\frac{x_{1}-x_{0}}{x_{2}}\right)
\res_{y_{2}}z^{-1}
\delta\left(\frac{x_{2}^{-1}-y_{2}}{z}\right)
\cdot \nn
&&\quad\quad \cdot\lambda(
Y_{1}(e^{x_{2}L(1)}(-x_{2}^{-2})^{L(0)}
Y(u, x_{0})
v, y_{2})w_{(1)}\otimes w_{(2)})\nn
&&=\left(x_{2}^{-1}\delta\left(\frac{x_{1}-x_{0}}{x_{2}}\right)
Y'_{P(z)}(Y(u, x_{0})v, x_{2})\lambda\right)(w_{(1)}\otimes w_{(2)}),
\end{eqnarray}
and we have proved the Jacobi identity and hence Theorem
\ref{comp=>jcb}.
\epfv

\noindent {\it Proof of Theorem \ref{stable}} Let $\lambda$ be an
element of $(W_{1}\otimes W_{2})^{*}$ satisfying the
$P(z)$-compatibility condition. We first want to prove that the
coefficient of each power of $x$ in $Y'_{P(z)}(u, x_0)Y'_{P(z)}(v,
x)\lambda$ is a formal Laurent series involving only finitely many
negative powers of $x_0$ and that
\begin{eqnarray}\label{stable-1}
\lefteqn{\tau_{P(z)}\left(x_{0}^{-1}\delta\left(\frac{x_{1}^{-1}-z}
{x_{0}}\right)
Y_t(u, x_1)\right)
Y'_{P(z)}(v, x)\lambda}\nno\\
&&=x_{0}^{-1}\delta\left(\frac{x_{1}^{-1}-z}
{x_{0}}\right)
Y'_{P(z)}(u, x_1)Y'_{P(z)}(v,
x) \lambda
\end{eqnarray}
for all $u, v\in V$.  Using the commutator formula (Proposition
\ref{pz-comm}) for $Y'_{P(z)}$, we have
\begin{eqnarray}\label{stable-2}
\lefteqn{Y'_{P(z)}(u, x_0)Y'_{P(z)}(v, x)\lambda}\nno\\
&&=Y'_{P(z)}(v, x)Y'_{P(z)}(u, x_0)\lambda\nno\\
&&\quad -\res_{y}x^{-1}_0\delta\left(\frac{x-y}{x_0}\right)
Y'_{P(z)}(Y(v, y)u, x_0)\lambda.
\end{eqnarray}
Each coefficient in $x$ of the right-hand side of (\ref{stable-2}) is
a formal Laurent series involving only finitely many negative powers
of $x_0$ since $\lambda$ satisfies the $P(z)$-lower truncation
condition.  Thus the coefficients in $x$ of $Y'_{P(z)}(v, x)\lambda$
satisfy the $P(z)$-lower truncation condition.

By (\ref{taudef}) and (\ref{Y'def}), we have
\begin{eqnarray}\label{stable-3}
\lefteqn{\left(\tau_{P(z)}\left(x_{0}^{-1}\delta\left(\frac{x_{1}^{-1}-z}
{x_{0}}\right)
Y_t(u, x_1)\right)
Y'_{P(z)}(v, x) \lambda\right)(w_{(1)}\otimes w_{(2)})}\nno\\
&&=z^{-1}\delta\left(\frac{x_1^{-1}-x_{0}}{z}\right)(Y'_{P(z)}(v, x)
\lambda)(Y_1(e^{x_{1}L(1)}(-x_{1}^{-2})^{L(0)}u, x_0)w_{(1)}
\otimes w_{(2)})\nno\\
&&\quad +x^{-1}_0\delta\left(\frac{z-x_1^{-1}}{-x_0}\right)(Y'_{P(z)}(v, x)
\lambda)(w_{(1)}\otimes Y_2^{o}(u, x_1)w_{(2)})\nno\\
&&=z^{-1}\delta\left(\frac{x_1^{-1}-x_{0}}{z}\right)
\Biggl(\lambda(Y_1(e^{x_{1}L(1)}(-x_{1}^{-2})^{L(0)}u, x_0)w_{(1)}
\otimes Y_{2}^{o}(v, x)w_{(2)})\nno\\
&&\quad \quad +\res_{x_2}
z^{-1}\delta\left(\frac{x^{-1}-x_{2}}{z}\right)\cdot\nn
&&\quad \quad\quad \quad\cdot \lambda(Y_1(e^{xL(1)}(-x^{-2})^{L(0)}v, x_2)
Y_1(e^{x_{1}L(1)}(-x_{1}^{-2})^{L(0)}u,
x_0)w_{(1)}\otimes w_{(2)})\Biggr)\nno\\
&&\quad +x^{-1}_0\delta\left(\frac{z-x_1^{-1}}{-x_0}\right)
\Biggl(\lambda(w_{(1)}\otimes Y_2^o(v, x)Y_2^{o}(u, x_1)w_{(2)})\nno\\
&&\quad \quad +\res_{x_2} z^{-1}\delta\left(\frac{x^{-1}-x_{2}}{z}\right)
\lambda(Y_1(e^{xL(1)}(-x^{-2})^{L(0)}v, x_2)w_{(1)}\otimes
Y_2^{o}(u, x_1)w_{(2)})\Biggr).\nn
&&
\end{eqnarray}
Now the distributive law applies, giving us four terms. Then using the
opposite commutator formula for $Y^{o}_{2}$ (recall (\ref{op-jac-id}))
and the commutator formula for $Y_{2}$, and (\ref{taudef}), we can
write the right-hand side of (\ref{stable-3}) as
\begin{eqnarray}\label{stable-4}
\lefteqn{z^{-1}\delta\left(\frac{x_1^{-1}-x_{0}}{z}\right)\lambda
(Y_1(e^{x_{1}L(1)}(-x_{1}^{-2})^{L(0)}u, x_0)w_{(1)}
\otimes Y_{2}^{o}(v, x)w_{(2)})}\nno\\
&&\quad  +z^{-1}\delta\left(\frac{x_1^{-1}-x_{0}}{z}\right)\res_{x_2}
z^{-1}\delta\left(\frac{x^{-1}-x_{2}}{z}\right)\cdot\nn
&&\quad \quad\cdot \lambda(Y_1(e^{x_{1}L(1)}(-x_{1}^{-2})^{L(0)}u,
x_0)Y_1(e^{xL(1)}(-x^{-2})^{L(0)}v, x_2)w_{(1)}\otimes w_{(2)})\nno\\
&&\quad  +z^{-1}\delta\left(\frac{x_1^{-1}-x_{0}}{z}\right)\res_{x_2}
z^{-1}\delta\left(\frac{x^{-1}-x_{2}}{z}\right)\res_{x_{3}}
x_{0}^{-1}\delta\left(\frac{x_{2}-x_{3}}{x_{0}}\right)\cdot\nn
&&\quad \quad \cdot 
\lambda(Y_1(Y(e^{xL(1)}(-x^{-2})^{L(0)}v, x_3)
e^{x_{1}L(1)}(-x_{1}^{-2})^{L(0)}u,
x_0)w_{(1)}\otimes w_{(2)})\nn
&&\quad +x^{-1}_0\delta\left(\frac{z-x_1^{-1}}{-x_0}\right)
\lambda(w_{(1)}\otimes Y_2^{o}(u, x_1)Y_2^o(v, x)w_{(2)})\nn
&&\quad -x^{-1}_0\delta\left(\frac{z-x_1^{-1}}{-x_0}\right)
\res_{x_{3}}x^{-1}_1\delta\left(\frac{x-x_3}{x_1}\right)
\lambda(w_{(1)}\otimes Y_2^{o}(Y(v, x_{3})u, x_{1})w_{(2)})\nn
&&\quad  +x^{-1}_0\delta\left(\frac{z-x_1^{-1}}{-x_0}\right)
\res_{x_2} z^{-1}\delta\left(\frac{x^{-1}-x_{2}}{z}\right)\cdot\nn
&&\quad \quad \cdot 
\lambda(Y_1(e^{xL(1)}(-x^{-2})^{L(0)}v, x_2)w_{(1)}\otimes
Y_2^{o}(u, x_{1})w_{(2)})\nn
&&=\left(\tau_{P(z)}\left(x_{0}^{-1}\delta\left(\frac{x_{1}^{-1}-z}
{x_{0}}\right)
Y_t(u, x_1)\right)\lambda\right)(w_{(1)}\otimes
Y_2^{o}(v, x)w_{(2)})\nn
&&\quad +\res_{x_2}
z^{-1}\delta\left(\frac{x^{-1}-x_{2}}{z}\right)\cdot\nn
&&\quad \quad\cdot \left(\tau_{P(z)}
\left(x_{0}^{-1}\delta\left(\frac{x_{1}^{-1}-z}
{x_{0}}\right)
Y_t(u, x_1)\right)\lambda\right)
(Y_1(e^{xL(1)}(-x^{-2})^{L(0)}v, x_2)w_{(1)}\otimes w_{(2)})\nn
&&\quad  +\res_{x_2}
z^{-1}\delta\left(\frac{x^{-1}-x_{2}}{z}\right)\res_{x_{3}}
x_{0}^{-1}\delta\left(\frac{x_{2}-x_{3}}{x_{0}}\right)
z^{-1}\delta\left(\frac{x_1^{-1}-x_{0}}{z}\right)\cdot\nn
&&\quad \quad \cdot 
\lambda(Y_1(Y(e^{xL(1)}(-x^{-2})^{L(0)}v, x_3)
e^{x_{1}L(1)}(-x_{1}^{-2})^{L(0)}u,
x_0)w_{(1)}\otimes w_{(2)})\nn
&&\quad -\res_{x_{3}}x^{-1}_1\delta\left(\frac{x-x_3}{x_1}\right)
x^{-1}_0\delta\left(\frac{z-x_1^{-1}}{-x_0}\right)
\lambda(w_{(1)}\otimes Y_2^{o}(Y(v, x_{3})u, x_{1})w_{(2)}).
\end{eqnarray}

Since $\lambda$ satisfies the $P(z)$-compatibility condition
(\ref{cpb}), by (\ref{Y'def}) the sum of the first two terms of
(\ref{stable-4}) is equal to
\begin{eqnarray}\label{stable-5}
\lefteqn{\left(Y'_{P(z)}(v, x)
\tau_{P(z)}\left(x_{0}^{-1}\delta\left(\frac{x_{1}^{-1}-z}
{x_{0}}\right)
Y_t(u, x_1)\right)\lambda\right)(w_{(1)}\otimes
w_{(2)})}\nn
&&=x_{0}^{-1}\delta\left(\frac{x_{1}^{-1}-z}
{x_{0}}\right)\left(Y'_{P(z)}(v, x)
Y'_{P(z)}(u, x_1)\lambda\right)(w_{(1)}\otimes
w_{(2)}).
\end{eqnarray}

Changing the dummy variable $x_{3}$ to $-x_{3}x^{-1}x_{1}^{-1}$ where
we use $x_{3}$ to denote the new dummy variable, using
(\ref{deltafunctionsubstitutionformula}), (\ref{2termdeltarelation})
and (\ref{comp=>jcb-9}), and then evaluating $\res_{x_{2}}$, we see
that the third term of (\ref{stable-4}) is equal to
\begin{eqnarray}\label{stable-6}
\lefteqn{-\res_{x_2}
\res_{x_{3}}
z^{-1}\delta\left(\frac{x^{-1}-x_{2}}{z}\right)x^{-1}x_{1}^{-1}
x_{0}^{-1}\delta\left(\frac{x_{2}+x_{3}x^{-1}x_{1}^{-1}}
{x_{0}}\right)z^{-1}\delta\left(\frac{x_1^{-1}-x_{0}}{z}\right)\cdot}\nn
&&\quad \quad \cdot 
\lambda(Y_1(Y(e^{xL(1)}(-x^{-2})^{L(0)}v, -x_3x^{-1}x_{1}^{-1})
e^{x_{1}L(1)}(-x_{1}^{-2})^{L(0)}u,
x_0)w_{(1)}\otimes w_{(2)})\nn
&&=-\res_{x_2}
z^{-1}\delta\left(\frac{x^{-1}-x_{2}}{z}\right)\res_{x_{3}}x^{-1}x_{1}^{-1}
x_{2}^{-1}\delta\left(\frac{x_{0}-x_{3}x^{-1}x_{1}^{-1}}
{x_{2}}\right)z^{-1}\delta\left(\frac{x_1^{-1}-x_{0}}{z}\right)\cdot\nn
&&\quad \quad \cdot 
\lambda(Y_1(Y(e^{xL(1)}(-x^{-2})^{L(0)}v, -x_3x^{-1}x_{1}^{-1})
e^{x_{1}L(1)}(-x_{1}^{-2})^{L(0)}u,
x_0)w_{(1)}\otimes w_{(2)})\nn
&&=-\res_{x_{3}}\res_{x_2}
z^{-1}\delta\left(\frac{x^{-1}-x_{0}+x_{3}x^{-1}x_{1}^{-1}}{z}\right)
\cdot\nn
&&\quad\quad
 \cdot x^{-1}x_{1}^{-1}
x_{2}^{-1}\delta\left(\frac{x_{0}-x_{3}x^{-1}x_{1}^{-1}}
{x_{2}}\right)z^{-1}\delta\left(\frac{x_1^{-1}-x_{0}}{z}\right)\cdot\nn
&&\quad \quad \cdot 
\lambda(Y_1(Y(e^{xL(1)}(-x^{-2})^{L(0)}v, -x_3x^{-1}x_{1}^{-1})
e^{x_{1}L(1)}(-x_{1}^{-2})^{L(0)}u,
x_0)w_{(1)}\otimes w_{(2)})\nn
&&=-\res_{x_{3}}\res_{x_2}
(z+x_{0})^{-1}\delta\left(\frac{x^{-1}+x_{3}x^{-1}x_{1}^{-1}}{z+x_{0}}\right)
\cdot\nn
&&\quad\quad\cdot x^{-1}x_{1}^{-1}
x_{2}^{-1}\delta\left(\frac{x_{0}-x_{3}x^{-1}x_{1}^{-1}}
{x_{2}}\right)x_{1}\delta\left(\frac{z+x_{0}}{x_1^{-1}}\right)\cdot\nn
&&\quad \quad \cdot 
\lambda(Y_1(Y(e^{xL(1)}(-x^{-2})^{L(0)}v, -x_3x^{-1}x_{1}^{-1})
e^{x_{1}L(1)}(-x_{1}^{-2})^{L(0)}u,
x_0)w_{(1)}\otimes w_{(2)})\nn
&&=-\res_{x_{3}}\res_{x_2}
x_{1}\delta\left(\frac{x^{-1}+x_{3}x^{-1}x_{1}^{-1}}
{x_{1}^{-1}}\right)
\cdot\nn
&&\quad\quad\cdot x^{-1}x_{1}^{-1}
x_{2}^{-1}\delta\left(\frac{x_{0}-x_{3}x^{-1}x_{1}^{-1}}
{x_{2}}\right)x_{1}\delta\left(\frac{z+x_{0}}{x_1^{-1}}\right)\cdot\nn
&&\quad \quad \cdot 
\lambda(Y_1(Y(e^{xL(1)}(-x^{-2})^{L(0)}v, -x_3x^{-1}x_{1}^{-1})
e^{x_{1}L(1)}(-x_{1}^{-2})^{L(0)}u,
x_0)w_{(1)}\otimes w_{(2)})\nn
&&=-\res_{x_{3}}\res_{x_2}
x_{1}^{-1}\delta\left(\frac{x-x_{3}}
{x_{1}}\right)
\cdot\nn
&&\quad\quad\cdot 
x_{2}^{-1}\delta\left(\frac{x_{0}-x_{3}x^{-1}x_{1}^{-1}}
{x_{2}}\right)z^{-1}\delta\left(\frac{x_1^{-1}-x_{0}}{z}\right)\cdot\nn
&&\quad \quad \cdot 
\lambda(Y_1(Y(e^{xL(1)}(-x^{-2})^{L(0)}v, -x_3x^{-1}x_{1}^{-1})
e^{x_{1}L(1)}(-x_{1}^{-2})^{L(0)}u,
x_0)w_{(1)}\otimes w_{(2)})\nn
&&=-\res_{x_{3}}\res_{x_2}
x_{1}^{-1}\delta\left(\frac{x-x_{3}}
{x_{1}}\right)
x_{2}^{-1}\delta\left(\frac{x_{0}-x_{3}x^{-1}x_{1}^{-1}}
{x_{2}}\right)z^{-1}\delta\left(\frac{x_1^{-1}-x_{0}}{z}\right)\cdot\nn
&&\quad \quad \cdot 
\lambda(Y_1(e^{x_{1}L(1)}(-x_{1}^{-2})^{L(0)}
Y(v, x_{3})u,
x_0)w_{(1)}\otimes w_{(2)})\nn
&&=-\res_{x_{3}}
x_{1}^{-1}\delta\left(\frac{x-x_{3}}
{x_{1}}\right)z^{-1}\delta\left(\frac{x_1^{-1}-x_{0}}{z}\right)\cdot\nn
&&\quad \quad \cdot 
\lambda(Y_1(e^{x_{1}L(1)}(-x_{1}^{-2})^{L(0)}
Y(v, x_{3})u,
x_0)w_{(1)}\otimes w_{(2)}).
\end{eqnarray}

{}From (\ref{stable-6}), (\ref{taudef}) and (\ref{cpb}), the sum of
the last two terms of (\ref{stable-4}) becomes
\begin{eqnarray}\label{stable-7}
\lefteqn{-\res_{x_{3}}
x_{1}^{-1}\delta\left(\frac{x-x_{3}}
{x_{1}}\right)z^{-1}\delta\left(\frac{x_1^{-1}-x_{0}}{z}\right)\cdot}\nn
&&\quad \quad \cdot
\lambda(Y_1(e^{x_{1}L(1)}(-x_{1}^{-2})^{L(0)}
Y(v, x_{3})u,
x_0)w_{(1)}\otimes w_{(2)})\nn
&&\quad -\res_{x_{3}}x^{-1}_1\delta\left(\frac{x-x_3}{x_1}\right)
x^{-1}_0\delta\left(\frac{z-x_1^{-1}}{-x_0}\right)
\lambda(w_{(1)}\otimes Y_2^{o}(Y(v, x_{3})u, x_{1})w_{(2)})\nn
&&=-\res_{x_{3}}
x_{1}^{-1}\delta\left(\frac{x-x_{3}}{x_{1}}\right)\left(\tau_{P(z)}\left(
x_{0}^{-1}\delta\left(\frac{x_1^{-1}-z}{x_{0}}\right)
Y_{t}(Y(v, x_{3})u, x_{1})\right)\lambda\right)(w_{(1)}\otimes w_{(2)})\nn
&&=-x_{0}^{-1}\delta\left(\frac{x_1^{-1}-z}{x_{0}}\right)\res_{x_{3}}
x_{1}^{-1}\delta\left(\frac{x-x_{3}}{x_{1}}\right)
\left(Y'_{P(z)}(Y(v, x_{3})u, x_{1})
\lambda\right)(w_{(1)}\otimes w_{(2)}).
\end{eqnarray}
Using (\ref{stable-5}), (\ref{stable-7}) and the commutator formula
for $Y'_{P(z)}$, we now see that the right-hand side of
(\ref{stable-4}) is equal to
\begin{eqnarray}\label{stable-8}
\lefteqn{x_{0}^{-1}\delta\left(\frac{x_{1}^{-1}-z}
{x_{0}}\right)\left(Y'_{P(z)}(v, x)
Y'_{P(z)}(u, x_1)\lambda\right)(w_{(1)}\otimes
w_{(2)})}\nn
&&-x_{0}^{-1}\delta\left(\frac{x_1^{-1}-z}{x_{0}}\right)\res_{x_{3}}
x_{1}^{-1}\delta\left(\frac{x-x_{3}}{x_{1}}\right)
\left(Y'_{P(z)}(Y(v, x_{3})u, x_{1})
\lambda\right)(w_{(1)}\otimes w_{(2)})\nn
&&=x_{0}^{-1}\delta\left(\frac{x_{1}^{-1}-z}
{x_{0}}\right)\left(
Y'_{P(z)}(u, x_1)Y'_{P(z)}(v, x)\lambda\right)(w_{(1)}\otimes
w_{(2)}).
\end{eqnarray}
The formulas (\ref{stable-3}), (\ref{stable-4}) and (\ref{stable-8})
together prove (\ref{stable-1}), as desired.  For the M\"obius case,
the corresponding verification for $L'_{P(z)}(-1)$, $L'_{P(z)}(0)$ and
$L'_{P(z)}(1)$ is straightforward, as usual, and we omit this
verification.  The first half of Theorem \ref{stable} holds.

For the second half of Theorem \ref{stable}, suppose that $\lambda\in
(W_{1}\otimes W_{2})^{*}$ satisfies either the $P(z)$-local grading
restriction condition or the $L(0)$-semisimple condition.  Assume
without loss of generality that $\lambda$ is doubly homogeneous.  {}From
Remark \ref{stableundercomponentops}, we see that for $v \in V$ doubly
homogeneous, $m \in \Z$ and $j=-1,0,1$, the elements
$\tau_{P(z)}(v\otimes t^{m})\lambda$ and $L'_{P(z)}(j)\lambda$ are
also doubly homogeneous.  Each such element $\mu$ lies in
$W_{\lambda}$, and so $W_{\mu} \subset W_{\lambda}$.  Thus $\mu$
satisfies the $P(z)$-local grading restriction condition (or the
$L(0)$-semisimple condition), and the second half of Theorem
\ref{stable} follows.  \epfv

\subsection{Proofs of Theorems \ref{6.1} and \ref{6.2}}

In this subsection, we follow \cite{tensor2}; the arguments given
there carry over to our more general situation with very little
change.  We first prove certain formulas that will be useful later.

Let $\lambda\in (W_{1}\otimes W_{2})^{*}$, $w_{(1)}\in W_{1}$ and
$w_{(2)}\in W_{2}$.  {}From (\ref{LQ'(j)}) we have
\begin{eqnarray}\label{9.1}
\lefteqn{(L'_{Q(z)}(0)\lambda)(w_{(1)} \otimes w_{(2)})}\nno\\
&&= \lambda((L(0) -
zL(1))w_{(1)} \otimes w_{(2)})- \lambda(w_{(1)} \otimes (L(0) -
zL(-1))w_{(2)}),\;\;\;\;\;
\end{eqnarray}
where (as usual) we have used the same notations $L(0), L(-1), L(1)$
to denote operators on both $W_1$ and $W_2$.  For convenience we write
$L(-1)=L'_{Q(z)}(-1)$, $L(0) = L'_{Q(z)}(0)$ and $L(1)=L'_{Q(z)}(1)$
in the rest of this section.  There will be no confusion since the
operators act on different spaces.

\begin{lemma}\label{1-y1zL(0)}
For $\lambda\in (W_{1}\otimes W_{2})^{*}$, $w_{(1)}\in W_{1}$
and $w_{(2)}\in W_{2}$,
we have
\begin{eqnarray}\label{9.2}
\lefteqn{\biggl(\biggl(1-\frac{y_1}{z}\biggr)^{L(0)} \lambda\biggr)(w_{(1)}
\otimes w_{(2)})}\nno\\
&&= \lambda\biggl(\biggl(1-\frac{y_1}{z}\biggr)^{L(0)-zL(1)}w_{(1)}
\otimes \biggl(1-\frac{y_1}{z}\biggr)^{-(L(0)-z L(-
1))}w_{(2)}\biggr).
\end{eqnarray}
\end{lemma}
 \pf
{}From (\ref{9.1}),
\begin{eqnarray}
\lefteqn{\biggl(\biggl(1-\frac{y_1}{z}\biggr)^{L(0)} \lambda)(w_{(1)}
\otimes w_{(2)})}\nno\\
&&= (e^{L(0)\log(1-\frac{y_1}{z})} \lambda)(w_{(1)}
\otimes w_{(2)})\nno\\
&&= \lambda(e^{(L(0)-z L(1))\log(1-\frac{y_1}{z})}w_{(1)}
\otimes e^{-(L(0)-z L(-1))\log(1-\frac{y_1}{z})}w_{(2)})\nno\\
&&= \lambda\biggl(\biggl(1-\frac{y_1}{z}\biggr)^{L(0)-zL(1)}w_{(1)}
\otimes \biggr(1-\frac{y_1}{z}\biggr)^{-(L(0)-z L(-
1))}w_{(2)}\biggr). \hspace{1.5em}\square
\end{eqnarray}

\begin{lemma}\label{Y'Q(z)L(0)}
For $v\in V$,
\begin{equation}\label{9.4}
Y'_{Q(z)}(v,x) = \biggl(1-
\frac{y_1}{z}\biggr)^{L(0)}Y'_{Q(z)}
\biggl(\biggl(1-\frac{y_1}{z}\biggr)^{-L(0)}
v,\frac{x}{1-y_1/z}\biggr)\biggr(1-\frac{y_1}{z}\biggr)^{-
L(0)}.
\end{equation}
This formula also holds for the vertex operators associated with any
generalized $V$-module.
\end{lemma}
\pf
The identity (\ref{9.4}) will follow {}from the formula
\begin{equation}\label{9.5}
e^{yL(0)}Y'_{Q(z)}(v, x)e^{-yL(0)}=Y'_{Q(z)}(e^{yL(0)}
v, e^{y}x).
\end{equation}
To prove this, assume without loss of generality that $\mbox{wt}\
v=h\in {\mathbb Z}$, and use the $L(-1)$-derivative property
(\ref{QL-1}) and the commutator formulas (\ref{commu-q-z}) and
(\ref{qz-sl-2-qz-y-2}) to get
\begin{equation}\label{9.6}
[L(0), Y'_{Q(z)}(v, x)]=\left(x\frac{d}{dx}+h\right)Y'_{Q(z)}(v, x).
\end{equation}
Formula (\ref{9.5}) now follows {}from (an easier version of) the proof
of (\ref{710}).
\epf

\begin{lemma}\label{L(0)L(-1)formula}
Let $L(-1)$ and $L(0)$ be any operators satisfying the commutator
relation
\begin{equation}
[L(0), L(-1)]=L(-1).
\end{equation}
 Then
\begin{equation}\label{9.8}
\biggr(1-\frac{y_1}{x}\biggr)^{L(0)-xL(-1)}=e^{y_1L(-1)}
\biggr(1-\frac{y_1}{x}\biggr)^{L(0)}.
\end{equation}
\end{lemma}
\pf
We first prove that the derivative with respect to $y$ of
$$(1-y)^{L(0)-xL(-1)}
(1-y)^{-L(0)}
e^{-xyL(-1)}$$ is $0$.  Write $A=(1-y)^{L(0)-xL(-1)}$,
$B=(1-y)^{-L(0)}$,
$C=e^{-xyL(-1)}$.  Then
\begin{eqnarray}\label{9.9}
\frac{d}{dy}(ABC)
&=&-A(1-y)^{-1}(L(0)-xL(-1)))BC\nno\\
&&+A(1-y)^{-1}L(0)BC\nno\\
&&-xABL(-1)C.
\end{eqnarray}
Using (\ref{log:SL2-2}) we have
\begin{eqnarray}\label{9.10}
BL(-1)&=&(1-y)^{-L(0)}L(-1)\nno\\
&=&e^{(-\log(1-y))L(0)}L(-1)\nno\\
&=&L(-1)e^{(-\log(1-y))L(0)}e^{-\log(1-y)}\nno\\
&=&L(-1)(1-y)^{-L(0)}
(1-y)^{-1}\nno\\
&=&(1-y)^{-1}L(-1)B,
\end{eqnarray}
and substituting (\ref{9.10}) into (\ref{9.9}) gives
\begin{eqnarray*}
\frac{d}{dy}(ABC)&=&-A(1-y)^{-1}L(0)BC
+xA(1-y)^{-1}L(-1)BC
\nno\\
&&+A(1-y)^{-1}L(0)BC
-xA(1-y)^{-1}L(-1)BC\nno\\
&=&0.
\end{eqnarray*}
Thus $ABC$ is constant in $y$, and since $ABC\lbar_{y=0}=1$, we have
$ABC=1$, which is equivalent to (\ref{9.8}).\epfv

\noindent {\it Proof of Theorem \ref{6.1}} As always, the reader
should again observe the justifiability of each formal step in the
argument.

Let $\lambda$ be an element of  $(W_{1}\otimes W_{2})^{*}$
satisfying the $Q(z)$-compatibility condition, that is,
(a) the $Q(z)$-lower truncation condition---for all $v\in V$,
$Y'_{Q(z)}(v,x)\lambda = \tau_{Q(z)}(Y_t(v,x))\lambda$
involves only finitely many negative powers of $x$,
and (b)
\begin{eqnarray}\label{10.3}
\lefteqn{
\tau_{Q(z)}\left(z^{-1}\delta\left(\frac{x_{1}-x_0}{z}\right)Y_t(v,x_0)
\right)\lambda}\nno\\
&&=z^{-1}\delta\left(\frac{x_{1}-x_0}{z}\right)
Y'_{Q(z)}(v,x_0)\lambda\;\;\;
\mbox{\rm for all}\ v\in V.
\end{eqnarray}
By (\ref{5.2}) and (\ref{Y'qdef}), (\ref{10.3}) is equivalent to
\begin{eqnarray}
\lefteqn{x^{-1}_0 \delta\left(\frac{x_1-
z}{x_0}\right)\lambda(Y^o_1(v,x_1)w_{(1)} \otimes w_{(2)})}\nno\\
&&\quad- x^{-1}_0 \delta\left(\frac{z-x_1}{-x_0}\right)\lambda(w_{(1)}
\otimes
Y_2(v,x_1)w_{(2)})\nno\\
&&= z^{-1}\delta\left(\frac{x_{1}-x_0}{z}\right)\biggl(\res_{y_{1}}x^{-1}_0
\delta\left(\frac{y_1-
z}{x_0}\right)\lambda(Y^o_1(v,y_1)w_{(1)} \otimes
w_{(2)})\nno\\
&&\quad- \res_{y_1}x^{-1}_0 \delta\left(\frac{z-y_1}{-
x_0}\right)\lambda(w_{(1)} \otimes Y_2(v,y_1)w_{(2)})\biggr)
\end{eqnarray}
for all $v\in V$, $w_{(1)}\in W_{1}$ and $w_{(2)}\in W_{2}$.  It is
important to note that on the right-hand side the distributive law is
not valid since the two individual products are not defined. One
critical feature of the argument that follows is that we must rewrite
expressions to allow the application of distributivity.

By (\ref{Y'qdef}), we have
\begin{eqnarray}\label{10.5}
\lefteqn{\left(x^{-
1}_0\delta\left(\frac{x_1-x_2}{x_0}\right)
Y'_{Q(z)}(v_1,x_1)Y'_{Q(z)}(v_2,x_2)
\lambda\right)(w_{(1)} \otimes w_{(2)})}\nno\\
&&= x^{-1}_0 \delta\left(\frac{x_1-
x_2}{x_0}\right)(Y'_{Q(z)}(v_1,x_1)Y'_{Q(z)}(v_2,x_2) \lambda)(w_{(1)} \otimes
w_{(2)})\nno\\
&&= x^{-1}_0 \delta\left(\frac{x_1-x_2}{x_0}\right)\biggl(\res_{y_1}x^{-
1}_1 \delta\left(\frac{y_1-z}{x_1}\right)\cdot\nno\\
&&\hspace{6em}\cdot (Y'_{Q(z)}(v_2,x_2)
\lambda)(Y^o_1(v_1,y_1)w_{(1)} \otimes w_{(2)})\nno\\
 &&\quad -\res_{y_1} x^{-
1}_1 \delta \left(\frac{z - y_1}{-x_1}\right)(Y'_{Q(z)}(v_2,x_2)
\lambda)(w_{(1)} \otimes Y_2(v_1,y_1)w_{(2)})\biggr)\nno\\
&&= x^{-1}_0
\delta\left(\frac{x_1-x_2}{x_0}\right)\biggl(\res_{y_1}x^{-1}_1
\delta\left(\frac{y_1-
z}{x_1}\right)\res_{y_2}x^{-1}_2 \delta \left(\frac{y_2-
z}{x_2}\right)\cdot\nno\\
&&\hspace{6em}\cdot \lambda(Y^o_1(v_2,y_2)Y^o_1(v_1,y_1)w_{(1)} \otimes
w_{(2)})\nno\\
&&\quad -\res_{y_1}x^{-1}_1 \delta\left(\frac{y_1-
z}{x_1}\right)\res_{y_2}x^{-1}_2 \delta\left(\frac{z-y_2}{-
x_2}\right)\cdot \nno\\
&&\hspace{6em}\cdot \lambda(Y^o_1(v_1,y_1)w_{(1)} \otimes Y_2(v_2,y_2)w_{(2)})\nno\\
&&\quad -
\res_{y_1}x^{-1}_1 \delta\left(\frac{z-y_1}{-x_1}\right)(Y'_{Q(z)}(v_2,x_2)
\lambda)(w_{(1)} \otimes Y_2(v_1,y_1)w_{(2)})\biggr).\;\;\;\;\;\;
\end{eqnarray}
{}From the properties of the formal $\delta$-function,
we see that the right-hand side of (\ref{10.5})
is equal to
\begin{eqnarray}\label{10.6}
\lefteqn{x^{-1}_0
\delta\left(\frac{x_1-x_2}{x_0}\right)\biggl(\res_{y_1}y^{-1}_1
\delta\left(\frac{x_1+z}{y_1}\right)\res_{y_2}y^{-1}_2
\delta\left(\frac{x_2+z}{y_2}\right)\cdot} \nno\\
&&\hspace{6em}\cdot \lambda(Y^o_1(v_2,x_2+z)Y^o_1(v_1,x_1+z)w_{(
1)} \otimes w_{(2)})\nno\\
&&\quad -\res_{y_1}y^{-1}_1
\delta\left(\frac{x_1+z}{y_1}\right)\res_{y_2}x^{-1}_2
\delta\left(\frac{z-y_2}{-
x_2}\right)\cdot \nno\\
&&\hspace{6em}\cdot \lambda(Y^o_1(v_1,x_1 +z)w_{(1)} \otimes
Y_2(v_2,y_2)w_{(2)})\nno\\
&&\quad - \res_{y_1}x^{-1}_1 \delta\left(\frac{z-y_1}{-
x_1}\right)(Y'_{Q(z)}(v_2,x_2) \lambda)(w_{(1)} \otimes
Y_2(v_1,y_1)w_{(2)})\biggr)\nno\\
&&= x^{-1}_0 \delta \left(\frac{x_1-
x_2}{x_0}\right)\biggl(\lambda(Y^o_1(v_2,x_2+z)
Y^o_1(v_1,x_1+z)w_{(1)}
\otimes w_{(2)})\nno\\
&&\quad - \res_{y_2}x^{-1}_2 \delta\left(\frac{z-y_2}{-
x_2}\right)\lambda(Y^o_1(v_1,x_1+z)w_{(1)} \otimes
Y_2(v_2,y_2)w_{(2)})\nno\\
&&\quad - \res_{y_1}x^{-1}_1 \delta\left(\frac{z-y_1}{-
x_1}\right)(Y'_{Q(z)}(v_2,x_2) \lambda)(w_{(1)} \otimes
Y_2(v_1,y_1)w_{(2)})\biggr).\;\;\;\;\;\;\;\;\;
\end{eqnarray}
{}From the $L(-1)$-derivative property for $Y_{1}$, the
$L(-1)$-derivative property (\ref{yo-l-1}) for $Y_{1}^{o}$ and the
commutator formulas for $L(-1)$, $Y_{1}(\cdot, x)$ and for $L(1)$,
$Y^{o}_{1}(\cdot, x)$ (recall Lemma \ref{sl2opposite}), we obtain
\begin{eqnarray}\label{10.7}
Y_{1}(v, x+z)&=&Y_{1}(e^{zL(-1)}v, x)\nno\\
&=&e^{zL(-1)}Y_{1}(v, x)e^{-zL(-1)}\nno\\
&=&\sum_{n\ge
0}\frac{z^{n}}{n!}\frac{d^{n}}{dx^{n}}Y_{1}(v, x)
\end{eqnarray}
and
\begin{eqnarray}\label{10.8}
Y^{o}_{1}(v, x+z)&=&Y^{o}_{1}(e^{zL(-1)}v, x)\nno\\
&=&e^{-zL(1)}Y^{o}_{1}(v, x)e^{zL(1)}\nno\\
&=&\sum_{n\ge
0}\frac{z^{n}}{n!}\frac{d^{n}}{dx^{n}}Y^{o}_{1}(v, x).
\end{eqnarray}
(Note that all these expressions are in fact defined.)
Using (\ref{10.8}), we see that the right-hand side of (\ref{10.6})
can be written as
\begin{eqnarray}\label{10.9}
\lefteqn{x^{-1}_0 \delta \left(\frac{x_1-
x_2}{x_0}\right)\biggl(\lambda(Y^o_1(e^{zL(-1)}v_2,x_2)Y^o_1(e^{zL(-1)}v_1,x_1)w_{(1)}
\otimes w_{(2)})}\nno\\
&&\quad - \res_{y_2}x^{-1}_2 \delta\left(\frac{z-y_2}{-
x_2}\right)\lambda(Y^o_1(e^{zL(-1)}v_1,x_1)w_{(1)} \otimes
Y_2(v_2,y_2)w_{(2)})\nno\\
&&\quad - \res_{y_1}x^{-1}_1 \delta\left(\frac{z-y_1}{-
x_1}\right)(Y'_{Q(z)}(v_2,x_2) \lambda)(w_{(1)} \otimes
Y_2(v_1,y_1)w_{(2)})\biggr)\nno\\
&&=x^{-1}_0 \delta \left(\frac{x_1-
x_2}{x_0}\right)\biggl(\lambda(e^{-zL(1)}Y^o_1(v_2,x_2)Y^o_1(v_1,x_1)e^{zL(1)}
w_{(1)}
\otimes w_{(2)})\nno\\
&&\quad - \res_{y_2}x^{-1}_2 \delta\left(\frac{z-y_2}{-
x_2}\right)\lambda(e^{-zL(1)}Y^o_1(v_1,x_1)e^{zL(1)}w_{(1)} \otimes
Y_2(v_2,y_2)w_{(2)})\nno\\
&&\quad - \res_{y_1}x^{-1}_1 \delta\left(\frac{z-y_1}{-
x_1}\right)(Y'_{Q(z)}(v_2,x_2) \lambda)(w_{(1)} \otimes
Y_2(v_1,y_1)w_{(2)})\biggr).
\end{eqnarray}
Note that it is easier to verify the well-definedness of the terms on
the right-hand side of (\ref{10.6}) than that of the terms in (\ref{10.9}),
though
(\ref{10.9}) is sometimes easier to use if it is known that every term is
well defined.  Below we shall write expressions like those on the right-hand
side of (\ref{10.6}) in whichever way suits our needs.
The distributive law applies to the right-hand side of (\ref{10.6})
(or (\ref{10.9})) since
all three of the following expressions are defined:
$$x^{-1}_0
\delta\left(\frac{x_1-x_2}{x_0}\right)\lambda(Y^o_1(v_2,x_2+z)
Y^o_1(v_1,x_1+z)w_{(1)} \otimes w_{(2)}),$$ $$x^{-1}_0
\delta\left(\frac{x_1-x_2}{x_0}\right)\res_{y_2}x^{-1}_2
\delta\left(\frac{z-y_2}{-x_2}\right)\lambda(Y^o_1(v_1,x_1+z)w_{(1)} \otimes
Y_2(v_2,y_2)w_{(2)}),$$
$$x^{-1}_0
\delta\left(\frac{x_1-x_2}{x_0}\right)\res_{y_1}x^{-1}_1 \delta\left(\frac{z-
y_1}{-x_1}\right)(Y'_{Q(z)}(v_2,x_2) \lambda)(w_{(1)} \otimes
Y_2(v_1,y_1)w_{(2)}).$$

Now we examine the last expression in (\ref{10.9}).  Rewriting the
formal $\delta$-functions
$x_{0}^{-1}\delta\left(\frac{x_{1}-x_{2}}{x_{0}}\right)$ and
$x_{1}^{-1}\delta\left(\frac{z-y_1}{-x_1}\right)$, and using Lemma
\ref{Y'Q(z)L(0)} and (\ref{deltafunctionsubstitutionformula}), we
have:
\begin{eqnarray}\label{10.10}
\lefteqn{x^{-1}_0\delta\left(\frac{x_1-x_2}{x_0}\right)\res_{y_{1}}x^{-1}_1
\delta\left(\frac{z-y_1}{-x_1}\right)\cdot}\nno\\
&&\hspace{4em}\cdot (Y'_{Q(z)}(v_2, x_2)\lambda)(w_{(1)}\otimes
 Y_2(v_1,
y_1)w_{(2)})\nno\\
&&=\biggl(\frac{x_1}{z}\biggr)^{-1}\biggl(\frac{x_0}{x_1/z}\biggr)^{-
1}\delta\left(\frac{z+(\frac{x_2}{-x_1/z})}{\frac{x_0}{x_1/z}}\right)
\res_{y_1}x_1^{-1}\delta\left(\frac{1-y_1/z}{-x_1/z}\right)\cdot \nno\\
&&\hspace{4em}\cdot\biggl(\biggl(1-\frac{y_1}{z}\biggr)^{L(0)}Y'_{Q(z)}
\biggl(\biggl(1-\frac{y_1}{z}\biggr)^{-
L(0)}v_2, \frac{x_2}{1-y_{1}/z}\biggr)\cdot \nno\\
&&\hspace{5em}\cdot \biggl(1-\frac{y_1}{z}\biggr)^{-L(0)}
\lambda\biggr)(w_{(1)}\otimes Y_2(v_1, y_1)w_{(2)})\nno\\
&&=\biggl(\frac{x_1}{z}\biggr)^{-1}\biggl(\frac{x_0}{x_1/z}\biggr)^{-
1}\delta\left(\frac{z+(\frac{x_2}{-x_1/z})}{\frac{x_0}{x_1/z}}\right)
\res_{y_1}x_1^{-1}\delta\left(\frac{1-y_1/z}{-x_1/z}\right)\cdot \nno\\
&&\hspace{4em}\cdot\biggl(\biggl(1-\frac{y_1}{z}\biggr)^{L(0)}Y'_{Q(z)}
\biggl(\biggl(1-\frac{y_1}{z}\biggr)^{-
L(0)}v_2, \frac{x_2}{-x_{1}/z}\biggr)\cdot \nno\\
&&\hspace{5em}\cdot \biggl(1-\frac{y_1}{z}\biggr)^{-L(0)}
\lambda\biggr)(w_{(1)}\otimes Y_2(v_1, y_1)w_{(2)}).
\end{eqnarray}
By Lemma \ref{1-y1zL(0)} and (\ref{10.3}), the right-hand side of
 (\ref{10.10}) is equal to
\begin{eqnarray}\label{10.11}
\lefteqn{\biggl(\frac{x_1}{z}\biggr)^{-1}\biggl(\frac{x_0}{x_1/z}\biggr)^{-
1}\delta\left(\frac{z+(\frac{x_2}{-x_1/z})}{\frac{x_0}{x_1/z}}\right)
\res_{y_1}x_1^{-1}\delta\left(\frac{1-y_1/z}{-x_1/z}\right)\cdot} \nno\\
&&\hspace{1em}\cdot \biggl(Y'_{Q(z)}\biggl(\biggl(1-\frac{y_1}{z}\biggr)^{-L(0)}
v_2, \frac{x_2}{-x_1/z}\biggr)\biggl(1-
\frac{y_1}{z}\biggr)^{-L(0)} \lambda\biggr)\cdot \nno\\
&&\hspace{3em}\cdot \biggl(\biggl(1-\frac{y_1}{z}\biggr)^{L(0)-
zL(1)}w_{(1)}\otimes\biggl(1-\frac{y_1}{z}\biggr)^{-
(L(0)-zL(-1))}Y_2(v_1, y_1)w_{(2)}\biggr)\nno\\
&&=\biggl(\frac{x_1}{z}\biggr)^{-1}
\res_{y_1}x_1^{-1}\delta\left(\frac{1-y_1/z}{-x_1/z}\right)\cdot \nno\\
&&\hspace{2em} \cdot\Biggl(\tau_{Q(z)}\Biggl(\biggl(\frac{x_0}{x_1/z}
\biggr)^{-
1}\delta\left(\frac{z+(\frac{x_2}{-x_1/z})}{\frac{x_0}{x_1/z}}\right)
\cdot \nno\\
&&\hspace{3em}\cdot Y_{t}\biggl(\biggl(1-\frac{y_1}{z}\biggr)^{-L(0)}v_2,
\frac{x_2}{-x_1/z}\biggr)\Biggr)\biggl(1-
\frac{y_1}{z}\biggr)^{-L(0)} \lambda\Biggr)\cdot \nno\\
&&\hspace{4em}\cdot \biggl(\biggl(1-\frac{y_1}{z}\biggr)^{L(0)-
zL(1)}w_{(1)}\otimes\biggl(1-\frac{y_1}{z}\biggr)^{-
(L(0)-zL(-1))}Y_2(v_1, y_1)w_{(2)}\biggr).
\end{eqnarray}
By (\ref{5.2}), the right-hand side of (\ref{10.11}) becomes
\begin{eqnarray}\label{10.12}
\lefteqn{\biggl(\frac{x_1}{z}\biggr)^{-1}
\res_{y_1}x_1^{-1}\delta\left(\frac{1-y_1/z}{-x_1/z}\right)
\biggl(\biggl(\frac{x_2}{-x_1/z}\biggr)^{-
1}\delta\left(\frac{\frac{x_0}{x_1/z}-z}{\frac{x_2}{-x_1/z}}\right)\cdot}
\nno\\
&&\hspace{1em}\cdot \biggl(\biggl(1-\frac{y_1}{z}\biggr)^{-L(0)}\lambda\biggr)
\biggl(Y_1^o\biggl(\biggl(1-\frac{y_1}{z}\biggr)^{-L(0)}v_2,\frac{x_0}
{x_{1}/z}\biggr)\cdot \nno\\
&&\hspace{2em}\cdot \biggl(1-\frac{y_1}{z}\biggr)^{L(0)-zL(1)}
w_{(1)}\otimes \biggl(1-\frac{y_1}{z}\biggr)^{-(L(0)-zL(1))}Y_2(v_1,
y_1)w_{(2)}\biggr)\nno\\
&&\quad -\biggl(\frac{x_2}{-x_1/z}\biggr)^{-
1}\delta\left(\frac{z-\frac{x_0}{x_1/z}}{-\frac{x_2}{-x_1/z}}\right)
\biggl(\biggl(1-\frac{y_1}{z}\biggr)^{-L(0)}\lambda\biggr)
\biggl(\biggl(1-\frac{y_1}{z}\biggr)^{L(0)-zL(1)} \cdot\nno\\
&&\hspace{2em}\cdot w_{(1)} \otimes
 Y_2\biggl(\biggl(1-\frac{y_1}{z}\biggr)^{-L(0)}v_2,
\frac{x_0}{x_1/z}\biggr)\cdot\nno\\
&&\hspace{3em}\cdot
\biggl(1-\frac{y_1}{z}\biggr)^{-(L(0)-zL(-1))}Y_2(v_1,
y_1)w_{(2)}\biggr)\biggr).\;\;\;\;\;\;\;
\end{eqnarray}
Using Lemma \ref{1-y1zL(0)} again but with $1-\frac{y_{1}}{z}$
replaced by $(1-\frac{y_{1}}{z})^{-1}$, rewriting formal
$\delta$-functions and then using the distributive law, we see that
(\ref{10.12}) is equal to
\begin{eqnarray}\label{10.13}
\lefteqn{\res_{y_1}x_1^{-1}\delta\left(\frac{1-y_1/z}{-x_1/z}\right)
\biggl(-x_2^{-1}\delta\left(\frac{x_0-x_1}{-x_2}\right)\cdot} \nno\\
&&\hspace{4em}\cdot \lambda\biggl(\biggl(1-\frac{y_1}{z}\biggr)^{-(L(0)-zL(1))}
Y_1^o\biggl(\biggl(1-\frac{y_1}{z}\biggr)^{-L(0)}v_2,\frac{x_0}{-(1-
y_1/z)}\biggr)\cdot \nno\\
&&\hspace{5em}\cdot \biggl(1-\frac{y_1}{z}\biggr)^{L(0)-zL(1)}
w_{(1)}\otimes Y_2(v_1, y_1)w_{(2)}\biggr)\nno\\
&&\quad +x_2^{-1}\delta\left(\frac{x_1-x_0}{x_2}\right)
\lambda\biggl(w_{(1)}\otimes\biggl(1-\frac{y_1}{z}\biggr)^{L(0)-zL(-1)}\cdot \nno\\
&&\hspace{4em}\cdot Y_2\biggl(\biggl(1-\frac{y_1}{z}\biggr)^{-L(0)}v_2,
\frac{x_0}{-(1-
y_1/z)}\biggr)\cdot \nno\\
&&\hspace{5em}\cdot \biggl(1-\frac{y_1}{z}\biggr)^{-
(L(0)-zL(-1))}Y_2(v_1, y_1)w_{(2)})\biggr)\nno\\
&&=-\res_{y_1}x_1^{-1}\delta\left(\frac{1-y_1/z}{-x_1/z}\right)
x_2^{-1}\delta\left(\frac{x_0-x_1}{-x_2}\right)\cdot \nno\\
&&\hspace{4em}\cdot \lambda\biggl(\biggl(1-\frac{y_1}{z}\biggr)^{-(L(0)-zL(1))}
Y_1^o\biggl(\biggl(1-\frac{y_1}{z}\biggr)^{-L(0)}v_2,\frac{x_0}{-(1-
y_1/z)}\biggr)\cdot \nno\\
&&\hspace{5em}\cdot \biggl(1-\frac{y_1}{z}\biggr)^{L(0)-zL(1)}
w_{(1)}\otimes Y_2(v_1, y_1)w_{(2)}\biggr)\nno\\
&&\quad +\res_{y_1}x_1^{-1}\delta\left(\frac{1-y_1/z}{-x_1/z}\right)
x_2^{-1}\delta\left(\frac{x_1-x_0}{x_2}\right)
\lambda\biggl(w_{(1)}\otimes\nno\\
&&\hspace{4em}\otimes\biggl(1-\frac{y_1}{z}\biggr)^{L(0)-zL(-1)}
Y_2\biggl(\biggl(1-\frac{y_1}{z}\biggr)^{-L(0)}v_2,
\frac{x_0}{-(1-
y_1/z)}\biggr)\cdot \nno\\
&&\hspace{5em}\cdot \biggl(1-\frac{y_1}{z}\biggr)^{-
(L(0)-zL(-1))}Y_2(v_1, y_1)w_{(2)}\biggr).
\end{eqnarray}

But by Lemmas \ref{L(0)L(-1)formula} and \ref{Y'Q(z)L(0)},
\begin{eqnarray}\label{10.14}
\lefteqn{\biggl(1-\frac{y_1}{z}\biggr)^{L(0)-zL(-1)} Y_2
\biggl(\biggl(1-\frac{y_1}{z}\biggr)^{-L(0)}v_2, \frac{x_0}{-(1-
y_1/z)}\biggr)\biggl(1-\frac{y_1}{z}\biggr)^{-(L(0)-zL(-
1))}}\nno\\
&&=e^{y_{1}L(-1)}\biggl(1-\frac{y_1}{z}\biggr)^{L(0)}Y_2
\biggl(\biggl(1-\frac{y_1}{z}\biggr)^{-L(0)}v_2, \frac{x_0}{-(1-
y_1/z)}\biggr)\cdot\hspace{4em}\nno\\
&&\hspace{10em}\cdot\biggl(1-\frac{y_1}{z}\biggr)^{-L(0)}e^{-y_{1}L(-1)}\nno\\
&&=e^{y_{1}L(-1)}Y_2(v_2, -x_0) e^{-y_{1}L(-1)}\nno\\
&&=Y_2(v_2, -x_0+ y_1)\nno\\
&&=Y_2(v_2, -x_0-(z-y_1)+z)\nno\\
&&=Y_2(e^{zL(-1)}v_2, -x_0-
(z-y_1)).\hspace{17em}
\end{eqnarray}
We similarly have, using Lemmas \ref{L(0)L(-1)formula} and
\ref{Y'Q(z)L(0)} for $Y'_1(v_2, x)$ and then using (\ref{y'}) and
Theorem \ref{set:W'},
\begin{eqnarray}\label{10.15}
\lefteqn{\biggl(1-\frac{y_1}{z}\biggr)^{-(L(0)-zL(1))}
Y_1^o\biggl(\biggl(1-\frac{y_1}{z}\biggr)^{-L(0)}v_2,\frac{x_0}{-(1-
y_1/z)}\biggr)\biggl(1-\frac{y_1}{z}\biggr)^{L(0)-zL(1)}}\nno\\
&&=e^{-y_{1}L(1)}Y_{1}^{o}(v_{2}, -x_{0})e^{y_{1}L(1)}\nno\\
&&=Y_{1}^{o}(v_{2}, -x_{0}+y_{1})\nno\\
&&=Y_{1}^{o}(v_{2}, -x_{0}-(z-y_{1})+z)\nno\\
&& =Y_1^o(e^{zL(-1)}v_2, -x_0-
(z-y_1)).\hspace{17em}
\end{eqnarray}
Substituting (\ref{10.14}) and (\ref{10.15}) into the
right-hand side of (\ref{10.13}) and then
combining  with  (\ref{10.10})--(\ref{10.13}),
we obtain
\begin{eqnarray}\label{10.16}
\lefteqn{x^{-1}_0\delta\left(\frac{x_1-x_2}{x_0}\right)\res_{y_1}
x^{-1}_1\delta\left(\frac{z-y_1}{-x_1}\right)(Y'_{Q(z)}(v_2, x_2) \lambda)
(w_{(1)}\otimes Y_2(v_1, y_1)w_{(2)})}\nno\\
&&=-\res_{y_1} x^{-1}_1\delta\left(\frac{z-y_1}{-x_1}\right)
x^{-1}_2\delta\left(\frac{x_0-x_1}{-x_2}\right)\cdot \nno\\
&&\hspace{4em}\cdot \lambda(Y_1^o\biggl(v_2,
-x_0-(z-y_1)+z\biggr)
w_{(1)}\otimes Y_2(v_1, y_1)w_{(2)})\nno\\
&&\quad +\res_{y_1} x^{-1}_1\delta\left(\frac{z-y_1}{-x_1}\right)
x^{-1}_2\delta\left(\frac{x_1-x_0}{x_2}\right)\cdot \nno\\
&&\hspace{4em}\cdot \lambda(w_{(1)}\otimes Y_2(v_2, -x_0+y_1)
Y_2(v_1, y_1)w_{(2)}).\hspace{8em}
\end{eqnarray}
(We choose the form of the expression {}from (\ref{10.15})
in anticipation of the
next step.)

By (\ref{deltafunctionsubstitutionformula}) and (\ref{10.7}), the
right-hand side of (\ref{10.16}) is equal to
\begin{eqnarray}\label{10.17}
\lefteqn{-x^{-1}_2\delta\left(\frac{-x_0+x_1}{x_2}\right)\res_{y_1}
x^{-1}_1\delta\left(\frac{z-y_1}{-x_1}\right)\cdot} \nno\\
&&\hspace{4em}\cdot \lambda(Y_1^o(v_2, -x_0+x_1+z)w_{(1)}\otimes Y_2(v_1,
y_1)w_{(2)})\nno\\
&&\quad +x^{-1}_2\delta\left(\frac{x_1-x_0}{x_2}\right)\res_{y_1}
x^{-1}_1\delta\left(\frac{z-y_1}{-x_1}\right)\cdot \nno\\
&&\hspace{4em}\cdot \lambda(w_{(1)}\otimes Y_2(v_2, -x_0+y_1)
Y_2(v_1, y_1)w_{(2)})\nno\\
&&=-x^{-1}_2\delta\left(\frac{-x_0+x_1}{x_2}\right)\res_{y_1}
x^{-1}_1\delta\left(\frac{z-y_1}{-x_1}\right)\cdot \nno\\
&&\hspace{4em}\cdot \lambda(Y_1^o(e^{zL(-1)}v_2, x_{2})w_{(1)}\otimes Y_2(v_1,
y_1)w_{(2)})\nno\\
&&\quad +x^{-1}_2\delta\left(\frac{x_1-x_0}{x_2}\right)\res_{y_1}
x^{-1}_1\delta\left(\frac{z-y_1}{-x_1}\right)\cdot \nno\\
&&\hspace{4em}\cdot \lambda(w_{(1)}\otimes Y_2(v_2, -x_0+y_1)
Y_2(v_1, y_1)w_{(2)}).
\end{eqnarray}
Since
$$\res_{y_{2}}y_{2}^{-1}\delta\left(\frac{-x_{0}+y_{1}}{y_{2}}\right)=1,$$
the right-hand side of (\ref{10.17}) can be written as
\begin{eqnarray}\label{10.18}
&&-x^{-1}_2\delta\left(\frac{-x_0+x_1}{x_2}\right)\res_{y_1}
x^{-1}_1\delta\left(\frac{z-y_1}{-x_1}\right)\cdot \nno\\
&&\hspace{4em}\cdot \lambda(Y_1^o(e^{zL(-1)}v_2, x_{2})w_{(1)}\otimes Y_2(v_1,
y_1)w_{(2)})\nno\\
&&\quad +x^{-1}_2\delta\left(\frac{x_1-x_0}{x_2}\right)\res_{y_1}
x^{-1}_1\delta\left(\frac{z-y_1}{-x_1}\right)\res_{y_2}
y^{-1}_2\delta\left(\frac{-x_0+y_1}{y_2}\right)\cdot \nno\\
&&\hspace{4em}\cdot \lambda(w_{(1)}\otimes Y_2(v_2, -x_0+y_1)
Y_2(v_1, y_1)w_{(2)}).
\end{eqnarray}

By (\ref{2termdeltarelation}) and 
(\ref{deltafunctionsubstitutionformula}), (\ref{10.18}) becomes
\begin{eqnarray}\label{10.19}
\lefteqn{x^{-1}_0\delta\left(\frac{x_2-x_1}{-x_0}\right)
\res_{y_1} x^{-1}_1\delta\left(\frac{z-y_1}{-x_1}\right)\cdot} \nno\\
&&\hspace{4em}\cdot \lambda(Y_1^o(e^{zL(-1)}v_2, x_2)w_{(1)}\otimes Y_2(v_1,
y_1)w_{(2)})\nno\\
&&\quad +x^{-1}_2\delta\left(\frac{x_1-x_0}{x_2}\right)\res_{y_1}
x^{-1}_1\delta\left(\frac{z-y_1}{-x_1}\right)\res_{y_2}
y^{-1}_2\delta\left(\frac{-x_0+y_1}{y_2}\right)\cdot \nno\\
&&\hspace{4em}\cdot \lambda(w_{(1)}\otimes Y_2(v_2, y_2)Y_2(v_1, y_1)w_{(2)})\nno\\
&&=x^{-1}_0\delta\left(\frac{x_2-x_1}{-x_0}\right)\res_{y_1}
x^{-1}_1\delta\left(\frac{z-y_1}{-x_1}\right)\cdot \nno\\
&&\hspace{4em}\cdot \lambda(Y_1^o(e^{zL(-1)}v_2, x_2)w_{(1)}\otimes Y_2(v_1,
y_1)w_{(2)})\nno\\
&&\quad +x^{-1}_2\delta\left(\frac{x_1-x_0}{x_2}\right)\res_{y_1}
\res_{y_2} x^{-1}_1\delta\left(\frac{z-y_1}{-x_1}\right)
y^{-1}_2\delta\left(\frac{-x_0+y_1}{y_2}\right)\cdot \nno\\
&&\hspace{4em}\cdot \lambda(w_{(1)}\otimes Y_2(v_2, y_2)Y_2(v_1, y_1)w_{(2)})\nno\\
&&=x^{-1}_0\delta\left(\frac{x_2-x_1}{-x_0}\right)\res_{y_1}
x^{-1}_1\delta\left(\frac{z-y_1}{-x_1}\right)\cdot \nno\\
&&\hspace{4em}\cdot \lambda(Y_1^o(e^{zL(-1)}v_2, x_2)w_{(1)}\otimes Y_2(v_1,
y_1)w_{(2)})\nno\\
&&\quad -x^{-1}_2\delta\left(\frac{x_1-x_0}{x_2}\right)\res_{y_1} \res_{y_2}
\cdot\nno\\
&&\hspace{3em}\cdot
(x_2+x_0)^{-1}\delta\left(\frac{z-y_1}{-x_2-x_0}\right)
x^{-1}_0\delta\left(\frac{y_2-y_1}{-x_0}\right)\cdot \nno\\
&&\hspace{4em}\cdot \lambda(w_{(1)}\otimes Y_2(v_2, y_2)Y_2(v_1, y_1)w_{(2)})\nno\\
&&=x^{-1}_0\delta\left(\frac{x_2-x_1}{-x_0}\right)\res_{y_1}
x^{-1}_1\delta\left(\frac{z-y_1}{-x_1}\right)\cdot \nno\\
&&\hspace{4em}\cdot \lambda(Y_1^o(e^{zL(-1)}v_2, x_2)w_{(1)}\otimes Y_2(v_1,
y_1)w_{(2)})\nno\\
&&\quad -x^{-1}_2\delta\left(\frac{x_1-x_0}{x_2}\right)\res_{y_1} \res_{y_2}
x^{-1}_2\delta\left(\frac{z-y_2}{-x_2}\right)x^{-1}_0\delta\left(\frac{y_2-
y_1}{-x_0}\right)\cdot \nno\\
&&\hspace{4em}\cdot \lambda(w_{(1)}\otimes Y_2(v_2, y_2)Y_2(v_1, y_1)w_{(2)}).
\end{eqnarray}
Substituting  (\ref{10.16})--(\ref{10.19}) into (\ref{10.9}) we obtain
\begin{eqnarray}\label{10.20}
\lefteqn{\left(x^{-1}_0\delta\left(\frac{x_1-x_2}{x_0}\right)Y'_{Q(z)}(v_1, x_1)
Y'_{Q(z)}(v_2,
x_2) \lambda\right)(w_{(1)}\otimes w_{(2)})}\nno\\
&&=x^{-1}_0\delta\left(\frac{x_1-x_2}{x_0}\right)\lambda(Y_1^o(e^{zL(-
1)}v_2, x_2)Y_1^o(e^{zL(-
1)}v_1, x_1)w_{(1)}\otimes w_{(2)})\nno\\
&&\quad -x^{-1}_0\delta\left(\frac{x_1-x_2}{x_0}\right)
\res_{y_2} x^{-1}_2\delta\left(\frac{z-y_2}{-x_2}\right)\cdot\nno\\
&&\hspace{6em}\cdot\lambda(Y_1^o(e^{zL(-1)}v_1, x_1)w_{(1)}\otimes Y_2(v_2,
x_2)w_{(2)})\nno\\
&&\quad -x^{-1}_0\delta\left(\frac{x_2-x_1}{-x_0}\right)\res_{y_1}
x^{-1}_1\delta\left(\frac{z-y_1}{-x_1}\right)\cdot\nno\\
&&\hspace{6em}\cdot\lambda(Y_1^o(e^{zL(-1)}v_2, x_2)w_{(1)}\otimes Y_2(v_1,
y_1)w_{(2)})\nno\\
&&\quad +x^{-1}_2\delta\left(\frac{x_1-x_0}{x_2}\right)\res_{y_1}
\res_{y_2} x^{-1}_2\delta\left(\frac{z-y_2}{-x_2}\right)
x^{-1}_0\delta\left(\frac{y_2-y_1}{-x_0}\right)\cdot\nno\\
&&\hspace{6em}\cdot\lambda(w_{(1)}\otimes Y_2(v_2, y_2)Y_2(v_1, y_1)w_{(2)}).
\end{eqnarray}

Now consider the result of the calculation {}from (\ref{10.16})
to (\ref{10.19})
except for the last two steps in (\ref{10.19}). Reversing
the subscripts 1 and 2 of the symbols $v$, $x$ and $y$ and replacing 
$x_0$ by $-x_0$ in this result and then using (\ref{2termdeltarelation}), 
we have
\begin{eqnarray}\label{10.21}
\lefteqn{x^{-1}_0\delta\left(\frac{x_2-x_1}{-x_0}\right)\res_{y_2}
x^{-1}_2\delta\left(\frac{z-y_2}{-x_2}\right)\cdot }\nno\\
&&\hspace{6em}\cdot (Y'_{Q(z)}(v_1, x_1) \lambda)
(w_{(1)}\otimes Y_2(v_2, y_2)w_{(2)})\nno\\
&& =x^{-1}_0\delta\left(\frac{x_1-x_2}{x_0}\right)\res_{y_2}
x^{-1}_2\delta\left(\frac{z-y_2}{-x_2}\right)\cdot \nno\\
&& \hspace{6em}\cdot \lambda(Y_1^o(e^{zL(-1)}v_1, x_1)w_{(1)}\otimes Y_2(v_2,
x_2)w_{(2)})\nno\\
&&\quad -x^{-1}_2\delta\left(\frac{x_1-x_0}{x_2}\right)\res_{y_1}\res_{y_2}
x^{-1}_2\delta\left(\frac{z-y_2}{-x_2}\right)
x^{-1}_0\delta\left(\frac{y_1-y_2}{x_0}\right)\cdot \nno\\
&&\hspace{6em}\cdot \lambda(w_{(1)}\otimes Y_2(v_1, y_1)Y_2(v_2, y_2)w_{(2)}).
\end{eqnarray}
{}From (\ref{10.5})--(\ref{10.9}),
again reversing the subscripts 1 and 2 of the symbols
$v$, $x$ and $y$
and replacing $x_0$ by $-x_0$, and (\ref{10.21}), we have
 \begin{eqnarray}\label{10.22}
\lefteqn{\left(-x^{-1}_0\delta\left(\frac{x_2-x_1}{-x_0}\right)Y'_{Q(z)}(v_2,
x_2)Y'_{Q(z)}(v_1, x_1)\cdot \lambda\right) (w_{(1)}\otimes w_{(2)})}\nno\\
&&=-x^{-1}_0\delta\left(\frac{x_2-x_1}{-
x_0}\right)\lambda(Y_1^o(e^{zL(-1)}v_1, x_1)Y_1^o(e^{zL(-
1)}v_2, x_2)w_{(1)}\otimes w_{(2)})\nno\\
&&\quad +x^{-1}_0\delta\left(\frac{x_2-x_1}{-x_0}\right)\res_{y_1}
x^{-1}_1\delta\left(\frac{z-y_1}{-x_1}\right)\cdot \nno\\
&&\hspace{6em}\cdot \lambda(Y_1^o(e^{zL(-1)}v_2,
x_2)w_{(1)}\otimes Y_2(v_1, y_1)w_{(2)})\nno\\
&&\quad +x^{-1}_0\delta\left(\frac{x_1-x_2}{x_0}\right)\res_{y_2}
x^{-1}_2\delta\left(\frac{z-y_2}{-x_2}\right)\cdot \nno\\
&&\hspace{6em}\cdot \lambda(Y_1^o(e^{zL(-1)}v_1, x_1)w_{(1)}\otimes Y_2(v_2,
x_2)w_{(2)})\nno\\
&&-x^{-1}_2\delta\left(\frac{x_1-x_0}{x_2}\right)\res_{y_1}\res_{y_2}
x^{-1}_2\delta\left(\frac{z-y_2}{-x_2}\right) x^{-1}_0\delta\left(\frac{y_1-
y_2}{x_0}\right)\cdot \nno\\
&&\hspace{6em}\cdot \lambda(w_{(1)}\otimes Y_2(v_1, y_1)Y_2(v_2, y_2)w_{(2)}).
\end{eqnarray}

The formulas (\ref{10.20}) and (\ref{10.22}) give:
\begin{eqnarray}\label{10.23}
\lefteqn{\biggl(\biggl(x^{-1}_0\delta\left(\frac{x_1-x_2}{x_0}\right)
Y'_{Q(z)}(v_1, x_1)
Y'_{Q(z)}(v_2,
x_2)}\nno\\
&&\quad  -x^{-1}_0\delta\left(\frac{x_2-x_1}{-x_0}\right)
Y'_{Q(z)}(v_2, x_2)Y'_{Q(z)}(v_1,
x_1)\biggr) \lambda\biggr)(w_{(1)}\otimes w_{(2)})\nno\\
&&=\lambda\biggl(\biggl(x^{-1}_0\delta\left(\frac{x_1-x_2}{x_0}\right)
Y_1^o(e^{zL(-
1)}v_2, x_2)Y_1^o(e^{zL(-1)}v_1, x_1)\nno\\
&&\hspace{2em}-x^{-1}_0\delta\left(\frac{x_2-x_1}{-x_0}\right)
Y_1^o(e^{zL(-1)}v_1,
x_1)Y_1^o(e^{zL(-1)}v_2, x_2)\biggr)w_{(1)}\otimes w_{(2)}\biggr)\nno\\
&&\quad -x^{-1}_2\delta\left(\frac{x_1-x_0}{x_2}\right)\res_{y_1}\res_{y_2}
x^{-1}_2\delta\left(\frac{z-y_2}{-x_2}\right)\cdot \nno\\
&&\hspace{4em}\cdot \lambda\biggl(w_{(1)}\otimes \biggl(x^{-1}_0
\delta\left(\frac{y_1-y_2}{x_0}\right)
Y_2(v_1, y_1)Y_2(v_2, y_2)\nno\\
&&\hspace{2em}-x^{-1}_0\delta\left(\frac{y_2-y_1}{-x_0}\right)
Y_2(v_2, y_2)Y_2(v_1,
y_1)\biggr)w_{(2)}\biggr).
\end{eqnarray}
{}From the Jacobi identities for $Y_{1}^{o}$ and $Y_{2}$ and 
(\ref{log:p1}), the right-hand
side of (\ref{10.23}) is equal to
\begin{eqnarray}\label{10.24}
\lefteqn{x^{-1}_2\delta\left(\frac{x_1-x_0}{x_2}\right)\lambda(Y_1^o(Y(e^{zL(-
1)}v_1, x_0)e^{zL(-1)}v_2, x_2)w_{(1)}\otimes w_{(2)})}\nno\\
&&\quad -x^{-1}_2\delta\left(\frac{x_1-x_0}{x_2}\right)\res_{y_1}\res_{y_2}
x^{-1}_2\delta\left(\frac{z-y_2}{-x_2}\right) y^{-1}_2\delta\left(\frac{y_1-
x_0}{y_2}\right)\cdot \nno\\
&&\hspace{6em}\cdot \lambda(w_{(1)}\otimes Y_2(Y(v_1, x_0)v_2, y_2)w_{(2)})
\nno\\
&&=x^{-1}_2\delta\left(\frac{x_1-x_0}{x_2}\right)\lambda(Y_1^o(e^{zL(-
1)}Y(v_1, x_0)v_2, x_2)w_{(1)}\otimes w_{(2)})\nno\\
&&\quad -x^{-1}_2\delta\left(\frac{x_1-x_0}{x_2}\right)\res_{y_1}\res_{y_2}
x^{-1}_2\delta\left(\frac{z-y_2}{-x_2}\right)
y^{-1}_2\delta\left(\frac{y_1-
x_0}{y_2}\right)\cdot \nno\\
&&\hspace{6em}\cdot \lambda(w_{(1)}\otimes Y_2(Y(v_1, x_0)v_2, y_2)w_{(2)}).
\end{eqnarray}
Using (\ref{10.8}), evaluating
$\res_{y_{1}}$ and then using the definition of $Y'_{Q(z)}$ (recall
(\ref{Y'qdef})),
we finally see that the
right-hand side of (\ref{10.24}) is equal to
\begin{eqnarray}\label{10.25}
\lefteqn{x^{-1}_2\delta\left(\frac{x_1-x_0}{x_2}\right)\lambda(Y_1^o(Y(v_1,
x_0)v_2, x_2+z)w_{(1)}\otimes w_{(2)})}\nno\\
&&\quad -x^{-1}_2\delta\left(\frac{x_1-x_0}{x_2}\right)\res_{y_2}
x^{-1}_2\delta\left(\frac{z-y_2}{-x_2}\right)\cdot \nno\\
&&\hspace{6em}\cdot \lambda(w_{(1)}\otimes Y_2(Y(v_1, x_0)v_2, y_2)w_{(2)})\nno\\
&&=x^{-1}_2\delta\left(\frac{x_1-x_0}{x_2}\right)(Y'_{Q(z)}(Y(v_1,
x_0)v_2, x_2) \lambda)(w_{(1)}\otimes w_{(2)}),
\end{eqnarray}
proving Theorem \ref{6.1}.
\epfv

\noindent {\it Proof of Theorem \ref{6.2}} Let $\lambda$ be an element
of $(W_{1}\otimes W_{2})^{*}$ satisfying the $Q(z)$-compatibility
condition. We first want to prove that each coefficient in $x$ of
$Y'_{Q(z)}(u, x_0)Y'_{Q(z)}(v, x)\lambda$ is a formal Laurent series
involving only finitely many negative powers of $x_0$ and that
\begin{eqnarray}\label{11.1}
\lefteqn{\tau_{Q(z)}\left(z^{-1}\delta\left(\frac{x_{1}-x_0}{z}\right)
Y_t(u, x_0)\right)
Y'_{Q(z)}(v, x)\lambda}\nno\\
&&=z^{-1}\delta\left(\frac{x_{1}-x_0}{z}\right)
Y'_{Q(z)}(u, x_0)Y'_{Q(z)}(v,
x) \lambda
\end{eqnarray}
for all $u, v\in V$.
Using the commutator formula for $Y'_{Q(z)}$, we have
\begin{eqnarray}\label{11.2}
\lefteqn{Y'_{Q(z)}(u, x_0)Y'_{Q(z)}(v, x)\lambda}\nno\\
&&=Y'_{Q(z)}(v, x)Y'_{Q(z)}(u, x_0)\lambda\nno\\
&&\quad -\res_{y}x^{-1}_0\delta\left(\frac{x-y}{x_0}\right)
Y'_{Q(z)}(Y(v, y)u, x_0)\lambda.
\end{eqnarray}
Each coefficient in $x$ of the right-hand side of (\ref{11.2}) is a
formal Laurent series involving only finitely many negative powers of
$x_0$ since $\lambda$ satisfies the $Q(z)$-lower truncation condition.
Thus the coefficients in $x$ of $Y'_{Q(z)}(v, x)\lambda$ satisfy the
$Q(z)$-lower truncation condition.

By (\ref{5.2}) and (\ref{Y'qdef}), we have
\begin{eqnarray}\label{11.3}
\lefteqn{\left(\tau_{Q(z)}\left(z^{-1}\delta\left(\frac{x_{1}-x_0}{z}\right)
Y_t(u, x_0)\right)
Y'_{Q(z)}(v, x) \lambda\right)(w_{(1)}\otimes w_{(2)})}\nno\\
&&=x^{-1}_0\delta\left(\frac{x_1-z}{x_0}\right)(Y'_{Q(z)}(v, x)
\lambda)(Y_1^o(u, x_1)w_{(1)}\otimes w_{(2)})\nno\\
&&\quad -x^{-1}_0\delta\left(\frac{z-x_1}{-x_0}\right)(Y'_{Q(z)}(v, x)
\lambda)(w_{(1)}\otimes Y_2(u, x_1)w_{(2)})\nno\\
&&=x^{-1}_0\delta\left(\frac{x_1-z}{x_0}\right)\biggl(\lambda(Y_1^o(v,
x+z)Y_1^o(u, x_1)w_{(1)}\otimes w_{(2)})\nno\\
&&\quad \quad -\res_{x_2}x^{-1}\delta\left(\frac{z-x_2}{-x}\right)
\lambda(Y_1^o(u,
x_1)w_{(1)}\otimes Y_2(v, x_2)w_{(2)})\biggr)\nno\\
&&\quad -x^{-1}_0\delta\left(\frac{z-x_1}{-x_0}\right)\biggl(\lambda(Y_1^o(e^{zL(-
1)}v, x)w_{(1)}\otimes Y_2(u, x_1)w_{(2)})\nno\\
&&\quad \quad-\res_{x_2} x^{-1}\delta\left(\frac{z-x_2}{-x}\right)
\lambda(w_{(1)}\otimes
Y_2(v, x_2)Y_2(u, x_1)w_{(2)})\biggr).\;\;\;\;\;
\end{eqnarray}
Now the distributive law applies, giving us four terms. Inserting
\[
\res_{x_{4}}x_{4}^{-1}\delta\left(\frac{x+z}{x_{4}}\right)=1
\]
into the first of these terms and correspondingly replacing $x+z$ by $x_{4}$
in $Y^{o}_{1}(v, x+z)$, we can apply the commutator formula for $Y^{o}_{1}$
in the usual way. Also using the commutator formula for $Y_{2}$,
(\ref{2termdeltarelation}) and 
(\ref{deltafunctionsubstitutionformula}), we write the
 right-hand side of (\ref{11.3}) as
\begin{eqnarray}\label{11.4}
\lefteqn{x^{-1}_0\delta\left(\frac{x_1-z}{x_0}\right)\lambda(Y_1^o(u, x_1)Y_1^o(v,
x+z)w_{(1)}\otimes w_{(2)})}\nno\\
&&\quad -x^{-1}_0\delta\left(\frac{x_1-z}{x_0}\right)\res_{x_4}\res_{x_3}
x^{-1}_1\delta\left(\frac{x_{4}-x_3}{x_1}\right)
x^{-1}_4\delta\left(\frac{x+z}{x_4}\right)\cdot \nno\\
&&\hspace{6em}\cdot \lambda(Y_1^o(Y(v, x_3)u, x_1)w_{(1)}\otimes
w_{(2)})\nno\\
&&\quad -x^{-1}_0\delta\left(\frac{x_1-z}{x_0}\right)\res_{x_2}
x^{-1}\delta\left(\frac{z-x_2}{-x}\right)\lambda(Y_1^o(u, x_1)w_{(1)}\otimes
Y_2(v, x_2)w_{(2)})\nno\\
&&\quad -x^{-1}_0\delta\left(\frac{z-x_1}{-x_0}\right)\lambda(Y_1^o(e^{zL(-
1)}v, x)w_{(1)}\otimes Y_2(u, x_1)w_{(2)})\nno\\
&&\quad +x^{-1}_0\delta\left(\frac{z-x_1}{-x_0}\right)\res_{x_2}
x^{-1}\delta\left(\frac{z-x_2}{-x}\right)\lambda(w_{(1)}\otimes
Y_2(u, x_1)Y_2(v, x_2)w_{(2)})\nno\\
&&\quad +x^{-1}_0\delta\left(\frac{z-x_1}{-x_0}\right)\res_{x_2}
x^{-1}\delta\left(\frac{z-x_2}{-x}\right)\res_{x_3}
x^{-1}_1\delta\left(\frac{x_{2}-x_3}{x_1}\right)\cdot \nno\\
&&\hspace{6em}\cdot \lambda(w_{(1)}\otimes Y_2(Y(v, x_3)u, x_1)w_{(2)})\nno\\
&&=\Biggl(x^{-1}_0\delta\left(\frac{x_1-z}{x_0}\right)\lambda(Y_1^o(u,
x_1)Y_1^o(e^{zL(-1)}v, x)w_{(1)}\otimes w_{(2)})\nno\\
&&\quad \quad-x^{-1}_0\delta\left(\frac{z-x_1}{-x_0}\right)\lambda(Y_1^o(e^{zL(-
1)}v, x)w_{(1)}\otimes Y_2(u, x_1)w_{(2)})\Biggr)\nno\\
&&\quad -\res_{x_2} x^{-1}\delta\left(\frac{z-x_2}{-x}\right)
\Biggl(x^{-1}_0\delta\left(\frac{x_1-z}{x_0}\right)\cdot\nno\\
&&\hspace{6em}\cdot\lambda(Y_1^o(u,
x_1)w_{(1)}\otimes Y_2(v, x_2)w_{(2)})\nno\\
&&\quad \quad -x^{-1}_0\delta\left(\frac{z-x_1}{-x_0}\right)
\lambda(w_{(1)}\otimes
Y_2(u, x_1)Y_2(v, x_2)w_{(2)})\Biggr)\nno\\
&&\quad -\Biggl( \res_{x_4}\res_{x_3}(x_0+z)^{-1}\delta\left(\frac{x+z-x_3}{x_0+
z}\right) \cdot\nno\\
&&\hspace{5em}\cdot x^{-1}_4\delta\left(\frac{x+z}{x_4}\right)
x^{-1}_0\delta\left(\frac{x_1-z}{x_0}\right)\cdot \nno\\
&&\hspace{6em}\cdot \lambda(Y_1^o(Y(v, x_3)u, x_1)w_{(1)}\otimes w_{(2)})\nno\\
&&\quad \quad-\res_{x_2}\res_{x_3}x^{-1}\delta\left(\frac{z-(x_1+x_3)}{-
x}\right) x^{-1}_2\delta\left(\frac{x_1+x_3}{x_2}\right)\cdot \nno\\
&&\hspace{6em}\cdot x^{-1}_0
\delta\left(\frac{z-
x_1}{-x_0}\right)\lambda(w_{(1)}\otimes Y_2(Y(v, x_3)u, x_1)w_{(2)})\Biggr).
\end{eqnarray}

Using (\ref{5.2}) and (\ref{2termdeltarelation}) and evaluating
suitable residues, we see that the right-hand side of (\ref{11.4}) is
equal to
\begin{eqnarray}\label{11.5}
\lefteqn{\left(\tau_{Q(z)}\left(z^{-1}\delta\left(\frac{x_{1}-x_0}{z}\right)Y_t(u,
x_0)\right)\lambda\right)(Y_1^o(e^{zL(-1)}v, x)w_{(1)}\otimes w_{(2)})}\nno\\
&&\quad -\res_{x_2} x^{-1}\delta\left(\frac{z-x_2}{-x}\right)\cdot \nno\\
&&\hspace{4em}\cdot
\left(\tau_{Q(z)}\left(z^{-1}\delta\left(\frac{x_{1}-x_0}{z}\right)
Y_t(u, x_0)\right)\lambda\right)(w_{(1)}\otimes Y_2(v, x_2)w_{(2)})\nno\\
&&\quad -\biggl(\res_{x_3}x^{-1}_0
\delta\left(\frac{x-x_3}{x_0}\right)\cdot\nno\\
&&\hspace{4em}\cdot x^{-1}_0\delta\left(\frac{x_1-z}{x_0}\right)
\lambda(Y_1^o(Y(v, x_3)u,
x_1)w_{(1)}\otimes w_{(2)})\nno\\
&&\quad -\res_{x_3}x^{-1}\delta
\left(\frac{x_0+x_{3}}{x}\right)\cdot \nno\\
&&\hspace{4em}\cdot x^{-1}_0\delta\left(\frac{z-x_1}{-x_0}\right)
\lambda(w_{(1)}\otimes Y_2(Y(v, x_3)u, x_1)w_{(2)})\biggr).
\end{eqnarray}
Using (\ref{5.2}) and (\ref{Y'qdef}), we find that the right-hand side of
(\ref{11.5}) becomes
\begin{eqnarray}\label{11.6}
\lefteqn{\left(Y'_{Q(z)}(v, x)\tau_{Q(z)}\left(z^{-1}
\delta\left(\frac{x_{1}-x_0}{z}\right)
Y_t(u, x_0) \right)\lambda\right)(w_{(1)}\otimes w_{(2)})}\nno\\
&&\quad -\biggl( \res_{x_3} x^{-1}_0\delta\left(\frac{x-x_3}{x_0}\right)\cdot \nno\\
&&\hspace{3em}\cdot \tau_{Q(z)}\biggl(z^{-1}\delta\left(\frac{x_{1}-x_0}{z}
\right)
Y_t(Y(v, x_3)u,
x_0)\biggr) \lambda\biggr)(w_{(1)}\otimes w_{(2)}).\;\;\;\;\;
\end{eqnarray}
By the compatibility condition for $\lambda$ and the commutator formula for
$Y'_{Q(z)}$, the right-hand side of (\ref{11.6}) is
equal to
\begin{eqnarray}
\lefteqn{z^{-1}\delta\left(\frac{x_{1}-x_0}{z}\right)
(Y'_{Q(z)}(v, x)Y'_{Q(z)}(u,
x_0) \lambda)(w_{(1)}\otimes w_{(2)})}\nno\\
&&\quad -z^{-1}\delta\left(\frac{x_{1}-x_0}{z}\right)\biggl( \res_{x_3}
x^{-1}_0\delta\left(\frac{x-x_3}{x_0}\right)\cdot \nno\\
&&\hspace{6em}\cdot Y'_{Q(z)}(Y(v, x_3)u, x_0)
\lambda\biggr)(w_{(1)}\otimes w_{(2)})\nno\\
&&= z^{-1}\delta\left(\frac{x_{1}-x_0}{z}\right)
\biggl(\biggl(Y'_{Q(z)}(v, x)Y'_{Q(z)}(u, x_0)\nno\\
&&\quad -\res_{x_3} x^{-1}_0\delta\left(\frac{x-x_3}{x_0}\right)
Y'_{Q(z)}(Y(v, x_3)u,
x_0)\biggr) \lambda\biggr)(w_{(1)}\otimes w_{(2)})\nno\\
&&=z^{-1}\delta\left(\frac{x_{1}-x_0}{z}\right)
(Y'_{Q(z)}(u, x_0)Y'_{Q(z)}(v,
x) \lambda)(w_{(1)}\otimes w_{(2)}).
\end{eqnarray}
This proves (\ref{11.1}).  In the M\"obius case, the three operators
are handled in the usual way.  The first part of Theorem \ref{6.2} is
established.

The proof of the second half of Theorem \ref{6.2} is exactly like that
for Theorem \ref{stable}.
\epfv

\newpage

\setcounter{equation}{0}
\setcounter{rema}{0}

\section{The convergence condition}\label{convsec}

Now that we have constructed tensor product modules and functors, our
next goal is to construct natural associativity isomorphisms.  More
precisely, under suitable conditions, we shall construct a natural
isomorphism between two functors {}from ${\cal C}\times{\cal
C}\times{\cal C}$ to ${\cal C}$, one given by
\[
(W_1, W_2, W_3)\mapsto (W_1\boxtimes_{P(z_1-z_{2})}
W_2)\boxtimes_{P(z_2)}W_3,
\]
and the other by
\[
(W_1, W_2, W_3)\mapsto
W_1\boxtimes_{P(z_1)} (W_2\boxtimes_{P(z_2)}W_3),
\]
where $W_1$, $W_2$ and $W_3$ are objects of ${\cal C}$ and $z_1$ and
$z_2$ are suitable complex numbers.  This will give us natural module
isomorphisms
\[
\alpha_{P(z_{1}), P(z_{2})}^{P(z_{1}-z_{2}), P(z_{2})}:
(W_1\boxtimes_{P(z_1-z_{2})} W_2)\boxtimes_{P(z_2)}W_3 \to
W_1\boxtimes_{P(z_1)} (W_2\boxtimes_{P(z_2)}W_3),
\]
which we will call the ``associativity isomorphisms.''  We have seen
that geometric data plays a crucial role in the tensor product itself,
and we will see that it continues to be a crucial ingredient in the
construction of the associativity isomorphisms.

We will mainly follow the ideas developed in \cite{tensor4}, and it
will be natural for us to work only in the case where all tensor
products involved are of type $P(z)$, for various nonzero complex
numbers $z$ (recall Remark \ref{motivate-Mobius}).

As we have stated in Sections 4 and 5, in the remainder of this work,
in particular in this section, Assumptions \ref{assum} and
\ref{assum-c} hold.

In this section we study one of the prerequisites for the existence of
the associativity isomorphisms. As we have discussed Subsection 1.4,
in order to construct the associativity isomorphisms between tensor
products of three objects, we must have that the intertwining maps
involved are ``composable,'' which means that certain convergence
conditions have to be satisfied.

More precisely, we first need to consider the composition of a
$P(z_1)$-intertwining map and a $P(z_2)$-intertwining map for suitable
nonzero complex numbers $z_1$ and $z_2$.  Geometrically, these
compositions correspond to sewing operations (see \cite{H1}) of
Riemann surfaces with punctures and local coordinates.  Compositions
(that is, products and iterates) of maps of this type have been
defined in \cite{tensor4} for intertwining maps among ordinary
modules.  The same definitions carry over to the greater generality of
this work:

Recall {}from Definition \ref{im:imdef} the space
\[
\mathcal{M}[P(z)]_{W_{1}W_{2}}^{W_{3}}
\]
of $P(z)$-intertwining maps of type ${W_3}\choose {W_1W_{2}}$ for 
$z\in \C^{\times}$ and $W_{1}, W_{2}, W_{3}$ objects of $\mathcal{C}$.
Let $W_1$, $W_2$, $W_3$, $W_4$ and $M_1$ be objects of ${\cal C}$.
Let $z_1,z_2 \in {\mathbb C}^{\times}$, 
$I_1\in \mathcal{M}[P(z_{1})]_{W_{1}M_{1}}^{W_{4}}$ and 
$I_2\in \mathcal{M}[P(z_{2})]_{W_{2}W_{3}}^{M_{1}}$.
If for any $w_{(1)}\in W_1$,
$w_{(2)}\in W_2$, $w_{(3)}\in W_3$ and $w'_{(4)}\in W'_4$, the series
\begin{equation}\label{convp}
\sum_{n\in {\mathbb C}}\langle w'_{(4)}, I_1(w_{(1)}\otimes
\pi_n(I_2(w_{(2)}\otimes w_{(3)})))\rangle_{W_4}
\end{equation}
(recall the notation $\pi_n$ {}from (\ref{pi_n}) and Definition
\ref{Wbardef} and note that
$\pi_n(I_2(w_{(2)}\otimes w_{(3)}))\in M_1$) is absolutely 
convergent, then the sums of these series give a linear map
\[
W_1\otimes W_2\otimes W_3 \to (W'_4)^{*}.
\]
Recalling the arguments in Lemmas \ref{4.36} and
\ref{IlambdatoJlambda}, we see that the image of this map is actually
in $\overline{W_{4}}$, so that we obtain a linear map
\[
W_1\otimes W_2\otimes W_3 \to \overline{W_{4}}.
\]
Analogously, let $W_1$, $W_2$, $W_3$, $W_4$ and $M_2$ be objects of
${\cal C}$. let $z_2,z_0 \in {\mathbb C}^{\times}$, 
$I^1\in \mathcal{M}[P(z_{2})]_{M_{2}W_{3}}^{W_{4}}$ and 
$I^2\in \mathcal{M}[P(z_{0})]_{W_{1}W_{2}}^{M_{2}}$.
If for any
$w_{(1)}\in W_1$, $w_{(2)}\in W_2$, $w_{(3)}\in W_3$ and $w'_{(4)}\in
W'_4$, the series
\begin{equation}\label{convi}
\sum_{n\in {\mathbb C}}\langle w'_{(4)}, I^1(\pi_n
(I^2(w_{(1)}\otimes w_{(2)}))\otimes w_{(3)})\rangle_{W_4}
\end{equation}
is absolutely convergent, then the sums of these series also 
give a linear map
\[
W_1\otimes W_2\otimes W_3 \to \overline{W_{4}}.
\]

\begin{defi}\label{productanditerateexisting}{\rm Let $W_1$, $W_2$, $W_3$, 
$W_4$ and $M_1$ be objects of ${\cal C}$.  Let $z_1,z_2 \in {\mathbb
C}^{\times}$, $I_1\in \mathcal{M}[P(z_{1})]_{W_{1}M_{1}}^{W_{4}}$ and
$I_2\in \mathcal{M}[P(z_{2})]_{W_{2}W_{3}}^{M_{1}}$.  We say that {\it
the product of $I_1$ and $I_2$ exists} if for any $w_{(1)}\in W_1$,
$w_{(2)}\in W_2$, $w_{(3)}\in W_3$ and $w'_{(4)}\in W'_4$, the series
(\ref{convp}) is absolutely convergent. In this case, we denote the
sum (\ref{convp}) by
\begin{equation}\label{I-prod}
\langle w'_{(4)}, I_1(w_{(1)}\otimes
I_2(w_{(2)}\otimes w_{(3)}))\rangle.
\end{equation}
We call
the map
\[
W_1\otimes W_2\otimes W_3 \to \overline{W}_4,
\]
defined by (\ref{I-prod})
the {\it product} of $I_1$ and $I_2$ and denote
it by 
\[
I_1 \circ (1_{W_{1}}\otimes I_{2}).
\]
In particular, we have
\[
\langle w'_{(4)}, I_1(w_{(1)}\otimes
I_2(w_{(2)}\otimes w_{(3)}))\rangle=\langle w'_{(4)},
(I_1 \circ (1_{W_{1}}\otimes I_{2}))(w_{(1)}\otimes w_{(2)}\otimes
w_{(3)})\rangle.
\]
Analogously, let $W_1$, $W_2$, $W_3$, $W_4$ and $M_2$ be objects of
${\cal C}$, and let $z_2,z_0 \in {\mathbb C}^{\times}$, 
$I^1\in \mathcal{M}[P(z_{2})]_{M_{2}W_{3}}^{W_{4}}$ and 
$I^2\in \mathcal{M}[P(z_{0})]_{W_{1}W_{2}}^{M_{2}}$.
We say that {\it the iterate of $I^1$ and $I^2$ exists} if for any
$w_{(1)}\in W_1$, $w_{(2)}\in W_2$, $w_{(3)}\in W_3$ and $w'_{(4)}\in
W'_4$, the series (\ref{convi}) is absolutely convergent. In this
case, we denote the sum (\ref{convi}) by
\begin{equation}\label{I-iter}
\langle w'_{(4)},
I^1(I^2(w_{(1)}\otimes w_{(2)})\otimes w_{(3)})\rangle
\end{equation}
and we call the map
\[
W_1\otimes W_2\otimes W_3\to
\overline{W}_4
\]
defined by (\ref{I-iter})
the {\it iterate} of
$I^1$ and $I^2$ and denote
it by 
\[
I^1\circ
(I^2\otimes 1_{W_3}).
\]
In particular, we have
\[
\langle w'_{(4)},
I^1(I^2(w_{(1)}\otimes w_{(2)})\otimes w_{(3)})\rangle
=\langle w'_{(4)}, (I^1\circ
(I^2\otimes 1_{W_3}))(w_{(1)}\otimes w_{(2)}\otimes w_{(3)})\rangle.
\]  }
\end{defi}

\begin{rema}\label{grad-comp-prod-iter}
{\rm Note that {}from the grading compatibility condition
(\ref{grad-comp}) for $P(z)$-intertwining maps, the product and the
iterate defined above, when they exist, also satisfy the following
{\it grading compatibility conditions}: With the notation as in
Definition \ref{productanditerateexisting}, suppose that $w_{(1)}\in
W_1^{(\beta)}$, $w_{(2)}\in W_2^{(\gamma)}$ and $w_{(3)}\in
W_3^{(\delta)}$, where $\beta, \gamma, \delta \in \tilde A$.  Then
\[
(I_1 \circ (1_{W_{1}}\otimes I_{2}))(w_{(1)}\otimes 
w_{(2)}\otimes w_{(3)})\in \overline{W_{4}^{(\beta+\gamma+\delta)}}
\]
if the product of $I_1$ and $I_2$ exists, and 
\[
(I^1\circ
(I^2\otimes 1_{W_3}))(w_{(1)}\otimes
w_{(2)}\otimes w_{(3)})\in \overline{W_{4}^{(\beta+\gamma+\delta)}}
\]
if the iterate of $I^1$ and $I^2$ exists.}
\end{rema}

\begin{propo}\label{convergence}
The following two conditions are equivalent:
\begin{enumerate}
\item Let $W_1$, $W_2$, $W_3$, $W_4$ and $M_1$ be arbitrary
objects of ${\cal
C}$ and let $z_1$ and $z_2$ be arbitrary nonzero complex numbers satisfying
\[
|z_1|>|z_2|>0.
\] 
Then for any $I_1\in \mathcal{M}[P(z_{1})]_{W_{1}M_{1}}^{W_{4}}$ and
$I_2\in \mathcal{M}[P(z_{2})]_{W_{2}W_{3}}^{M_{1}}$, the product of
$I_1$ and $I_2$ exists.

\item Let $W_1$, $W_2$, $W_3$, $W_4$ and $M_2$ be arbitrary
objects of ${\cal
C}$ and let  $z_0$ and $z_2$ be arbitrary nonzero complex numbers satisfying
\[
|z_2|>|z_0|>0.
\]
Then for any $I^1\in \mathcal{M}[P(z_{2})]_{M_{2}W_{3}}^{W_{4}}$ and
$I^2\in \mathcal{M}[P(z_{0})]_{W_{1}W_{2}}^{M_{2}}$, the iterate of
$I^1$ and $I^2$ exists.
\end{enumerate}
\end{propo}
\pf We shall use the isomorphism $\Omega_0$ given by (\ref{Omega_r}))
and its inverse $\Omega_{-1}$ (recall Proposition \ref{log:omega}) to
prove this result.  Suppose that Condition 1 holds.  Let $z_0$ and
$z_2$ be any nonzero complex numbers. For any intertwining maps $I^1$
and $I^2$ as in the statement of Condition 2, let 
${\cal Y}^1=\Y_{I^{1}, 0}$ and
${\cal Y}^2=\Y_{I^{2}, 0}$ be the logarithmic intertwining 
operators corresponding
to $I^1$ and $I^2$, respectively, according to Proposition
\ref{im:correspond}.  We need to prove that when $|z_2|>|z_0|>0$, the
series (\ref{convi}), which can now be written as
\begin{equation}\label{4itm}
\sum_{n\in \C}\langle w'_{(4)}, 
{\cal Y}^1(\pi_{n}({\cal Y}^2(w_{(1)}, x_0)w_{(2)}),
x_2)w_{(3)} \rangle_{W_4}\lbar_{x_0=z_0, \;x_2=z_2}
\end{equation}
(recall the ``substitution'' notation {}from (\ref{im:f(z)}), where we
choose $p=0$ again), is absolutely convergent for any $w_{(1)}\in
W_1$, $w_{(2)}\in W_2$, $w_{(3)}\in W_3$ and $w'_{(4)}\in W'_4$.

Using the linear isomorphism $\Omega_{-1}: {\cal V}_{W_3 M_2}^{W_4}
\to {\cal V}_{M_2 W_3}^{W_4}$ (see (\ref{Omega_r})),
\[
\Omega_{-1}({\cal Y})(w, x)w_{(3)}=e^{xL(-1)}{\cal Y}(w_{(3)}, e^{-\pi
i}x)w,
\]
for ${\cal Y}\in {\cal V}_{W_3\,
M_2}^{W_4}$, $w\in M_2$ and 
$w_{(3)}\in W_3$, and its inverse $\Omega_0:{\cal V}_{M_2 W_3}^{W_4}\to
{\cal V}_{W_3 M_2}^{W_4}$, we have
\begin{eqnarray}\label{nosub}
\lefteqn{\langle w'_{(4)}, {\cal Y}^1({\cal Y}^2(w_{(1)},
x_0)w_{(2)}, x_2)w_{(3)} \rangle_{W_4}}\nno\\
&&=\langle w'_{(4)}, \Omega_{-1}(\Omega_0({\cal Y}^1))({\cal
Y}^2(w_{(1)}, x_0)w_{(2)}, x_2)w_{(3)} \rangle_{W_4} \nno\\
&&=\langle w'_{(4)}, e^{x_2L(-1)}\Omega_0({\cal Y}^1) (w_{(3)},
e^{-\pi i}x_2) {\cal Y}^2(w_{(1)}, x_0)w_{(2)}\rangle_{W_4}
\nno\\
&&=\langle e^{x_2L'(1)}w'_{(4)}, \Omega_0({\cal Y}^1) (w_{(3)},
e^{-\pi i}x_2) {\cal Y}^2(w_{(1)}, x_0)w_{(2)}\rangle_{W_4}
\end{eqnarray}
for $w_{(1)}\in W_1$, $w_{(2)}\in W_2$, $w_{(3)}\in W_3$ and
$w'_{(4)}\in W'_4$. Hence for $n\in \C$, 
\begin{eqnarray*}
\lefteqn{\langle w'_{(4)}, {\cal Y}^1(\pi_{n}({\cal Y}^2(w_{(1)},
x_0)w_{(2)}), x_2)w_{(3)} \rangle_{W_4}\lbar_{x_0= z_0,\;
x_2=z_2}}\nno\\
&&=\langle e^{x_2L'(1)}w'_{(4)}, \Omega_0({\cal Y}^1) (w_{(3)},
e^{-\pi i}x_2) \pi_{n}({\cal Y}^2(w_{(1)},
x_0)w_{(2)})\rangle_{W_4}\lbar_{x_0= z_0,\; x_2=z_2} \nno\\
&&=\langle e^{z_2L'(1)}w'_{(4)}, \Omega_0({\cal Y}^1)
(w_{(3)}, x_2) \pi_{n}({\cal Y}^2(w_{(1)},
x_0)w_{(2)})\rangle_{W_4}\lbar_{x_0= z_0,\; x_2=e^{-\pi i}
z_2},\nno\\ &&
\end{eqnarray*}
that is,
\begin{eqnarray}\label{i2p}
\lefteqn{\langle w'_{(4)}, {\cal Y}^1({\cal Y}^2(w_{(1)},
x_0)w_{(2)}, x_2)w_{(3)} \rangle_{W_4}\lbar_{x_0= z_0,\;
x_2=z_2}}\nno\\
&&=\langle e^{z_2L'(1)}w'_{(4)}, \Omega_0({\cal Y}^1)
(w_{(3)}, x_2) {\cal Y}^2(w_{(1)},
x_0)w_{(2)}\rangle_{W_4}\lbar_{x_0= z_0,\; x_2=e^{-\pi i}
z_2}.\nno\\ &&
\end{eqnarray}
Since the last expression is equal to the product of a
$P(-z_2)$-intertwining map and a $P(z_0)$-intertwining map evaluated
at $w_{(3)}\otimes w_{(1)}\otimes w_{(2)}\in W_3\otimes W_1\otimes
W_2$ and paired with $e^{z_2L'(1)}w'_{(4)}\in W'_4$, it converges
absolutely when $|-z_2|>|z_0|>0$, or equivalently, when
$|z_2|>|z_0|>0$.

Conversely, suppose that Condition 2 holds, and let $z_1$ and $z_2$ be
any nonzero complex numbers. For any intertwining maps $I_1$ and $I_2$
as in the statement of Condition 1, let ${\cal Y}_1=\Y_{I_{1}, 0}$ and
${\cal Y}_2=\Y_{I_{2}, 0}$ be the logarithmic intertwining operators
corresponding to $I_1$ and $I_2$, respectively.  We need to prove that
when $|z_1|>|z_2|>0$, the series (\ref{convp}), which can now be
written as
\begin{equation}\label{4prm}
\sum_{n\in \C}\langle w'_{(4)}, {\cal Y}_1(w_{(1)}, x_1) 
\pi_{n}({\cal Y}_2(w_{(2)},
x_2)w_{(3)})\rangle_{W_4}\lbar_{x_1=z_1, \; x_2=z_2},
\end{equation}
is absolutely convergent for any $w_{(1)}\in W_1$, $w_{(2)}\in W_2$,
$w_{(3)}\in W_3$ and $w'_{(4)}\in W'_4$.

Using the linear isomorphism $\Omega_0: {\cal V}_{M_1 W_1}^{W_4}
\to {\cal V}_{W_1 M_1}^{W_4}$,
\[
\Omega_0({\cal Y})(w_{(1)}, x)w=e^{xL(-1)}{\cal Y}(w, e^{\pi
i}x)w_{(1)},
\]
for  ${\cal Y}\in {\cal V}_{M_1\,
W_1}^{W_4}$, $w_{(1)}\in W_1$ and $w\in M_1$, 
and its inverse $\Omega_{-1}:{\cal V}_{W_1 M_1}^{W_4}\to
{\cal V}_{M_1 W_1}^{W_4}$, we have
\begin{eqnarray}\label{nosub2}
\lefteqn{\langle w'_{(4)}, {\cal Y}_1 (w_{(1)},
x_1) {\cal Y}_2(w_{(2)}, x_2)w_{(3)}\rangle_{W_4}}\nno\\
&&=\langle w'_{(4)}, \Omega_0(\Omega_{-1}({\cal Y}_1)) (w_{(1)},
x_1) {\cal Y}_2(w_{(2)}, x_2)w_{(3)}\rangle_{W_4}=\nno\\
&&=\langle w'_{(4)}, e^{x_1L(-1)}\Omega_{-1}({\cal Y}_1)({\cal
Y}_2(w_{(2)}, x_2)w_{(3)}, e^{\pi i}x_1)w_{(1)} \rangle_{W_4} \nno\\
&&=\langle e^{x_1L'(1)}w'_{(4)}, \Omega_{-1}({\cal Y}_1)({\cal
Y}_2(w_{(2)},
x_2)w_{(3)}, e^{\pi i}x_1)w_{(1)} \rangle_{W_4}
\end{eqnarray}
for $w_{(1)}\in W_1$, $w_{(2)}\in W_2$, $w_{(3)}\in W_3$ and
$w'_{(4)}\in W'_4$. Hence for $n\in \C$,
\begin{eqnarray*}
\lefteqn{\langle w'_{(4)}, {\cal Y}_1 (w_{(1)},
x_1) \pi_{n}({\cal Y}_2(w_{(2)}, x_2)w_{(3)})\rangle_{W_4}\lbar_{x_1= z_1,\;
x_2=z_2}}\nno\\
&&=\langle e^{x_1L'(1)}w'_{(4)}, \Omega_{-1}({\cal Y}_1)(
\pi_{n}({\cal
Y}_2(w_{(2)},
x_2)w_{(3)}), e^{\pi i}x_1)w_{(1)} \rangle_{W_4}\lbar_{x_1= z_1,\;
x_2=z_2}\nno\\
&&=\langle e^{z_1L'(1)}w'_{(4)}, \Omega_{-1}({\cal Y}_1)(
\pi_{n}({\cal
Y}_2(w_{(2)},
x_2)w_{(3)}), x_1)w_{(1)} \rangle_{W_4}\lbar_{x_1= e^{\pi i}z_1,\;
x_2=z_2},\nno\\ &&
\end{eqnarray*}
that is,
\begin{eqnarray}\label{p2i}
\lefteqn{\langle w'_{(4)}, {\cal Y}_1 (w_{(1)},
x_1) {\cal Y}_2(w_{(2)}, x_2)w_{(3)}\rangle_{W_4}\lbar_{x_1= z_1,\;
x_2=z_2}}\nno\\
&&=\langle e^{z_1L'(1)}w'_{(4)}, \Omega_{-1}({\cal Y}_1)({\cal
Y}_2(w_{(2)},
x_2)w_{(3)}, x_1)w_{(1)} \rangle_{W_4}\lbar_{x_1= e^{\pi i}z_1,\;
x_2=z_2}.\nno\\ &&
\end{eqnarray}
Since the last expression is equal to the iterate of a
$P(-z_1)$-intertwining map and a $P(z_2)$-intertwining map
evaluated at $w_{(2)}\otimes w_{(3)}\otimes w_{(1)}\in W_2\otimes
W_3\otimes W_1$ and paired with $e^{z_1L'(1)}w'_{(4)}\in W'_4$,
it converges absolutely when $|-z_1|>|z_2|>0$, or equivalently,
when $|z_1|>|z_2|>0$.
\epfv

For convenience, we shall use the notations 
\[
\langle w'_{(4)}, 
{\cal Y}^1({\cal Y}^2(w_{(1)}, x_0)w_{(2)},
x_2)w_{(3)} \rangle_{W_4}\lbar_{x_0=z_0, \;x_2=z_2}
\]
and 
\[
\langle w'_{(4)}, {\cal Y}_1(w_{(1)}, x_1) 
{\cal Y}_2(w_{(2)},
x_2)w_{(3)}\rangle_{W_4}\lbar_{x_1=z_1, \; x_2=z_2},
\]
or even more simply, the notations
\begin{equation}\label{iterateabbreviation}
\langle w'_{(4)}, 
{\cal Y}^1({\cal Y}^2(w_{(1)}, z_0)w_{(2)},
z_2)w_{(3)} \rangle_{W_4}
\end{equation}
and 
\begin{equation}\label{productabbreviation}
\langle w'_{(4)}, {\cal Y}_1(w_{(1)}, z_1) 
{\cal Y}_2(w_{(2)},
z_2)w_{(3)}\rangle_{W_4}
\end{equation}
to denote (\ref{4itm}) and (\ref{4prm}), respectively.
We shall also use similar notations to denote 
series obtained {}from products and iterates of more than 
two intertwining operators. 

\begin{defi}\label{conv-conditions}
{\rm We call either of the two conditions in Proposition
\ref{convergence} the {\it convergence condition for intertwining maps
in the category ${\cal C}$}.}
\end{defi}

We need the following concept concerning unique expansion of an
analytic function in terms of powers of $z$ and $\log z$  
(recall our choice of the branch of $\log z$ in
(\ref{branch1}) and thus the branch of $z^{\alpha}$, $\alpha\in \C$):

\begin{defi}
{\rm We call a subset ${\cal S}$ of ${\mathbb C}\times{\mathbb C}$ a
{\em unique expansion set} if the absolute convergence to $0$ on some
nonempty open subset of ${\mathbb C}^{\times}$ of any series
$\sum_{(\alpha,\beta)\in{\cal S}} a_{\alpha,\beta}z^\alpha(\log
z)^\beta$, $a_{\alpha,\beta} \in {\mathbb C}$, implies that
$a_{\alpha,\beta}=0$ for all $(\alpha,\beta)\in{\cal S}$. }
\end{defi}

\begin{rema}
{\rm It is easy to show, using the Laplace transform, 
that ${\mathbb Z}\times\{ 0,1,\dots,N\}$ is a
unique expansion set for any $N\in{\mathbb N}$; this is a special case
of the next proposition, in which we also use a Laplace transform argument.  
On the other hand, it is known that
${\mathbb C}\times\{0\}$ is {\em not} a unique expansion
set\footnote{We thank A.~Eremenko for informing us of this result.}.}
\end{rema}

\begin{lemma}\label{po-ser-an}
Let $D$ be a subset of $\R$ and let $\sum_{\alpha\in
D}a_{\alpha}z^{\alpha}$ ($a_{\alpha} \in {\mathbb C}$)
be absolutely convergent on a (nonempty) open
subset of $\C^{\times}$. Then $\sum_{\alpha\in D}a_{\alpha}\alpha
z^{\alpha}$ is absolutely and uniformly convergent near any $z$ in the
open subset. In particular, the sum $\sum_{\alpha\in
D}a_{\alpha}z^{\alpha}$ as a function of $z$ is analytic in the sense that 
it is analytic at $z$ when $z$ is in the open subset of $\C^{\times}$ and 
$\arg z>0$, and that it can be analytically extended to an analytic 
function in a neighborhood of $z$
when $z$ is in the intersection of the open subset and the 
positive real line.
\end{lemma}
\pf
We need only prove that 
$\sum_{\alpha\in D}a_{\alpha}\alpha z^{\alpha}$
is absolutely and uniformly convergent near any 
$z$ in the open subset. Note that since the original series is 
absolutely convergent on an open subset of $\C^{\times}$, 
$\sum_{\alpha\in D, \; \alpha\ge 0}a_{\alpha}z^{\alpha}$
and 
$\sum_{\alpha\in D, \; \alpha< 0}a_{\alpha}z^{\alpha}$
are also absolutely 
convergent on the set. For any fixed $z_{0}$ in the 
set, we can always find $z_{1}$ and $z_{2}$ in the set
such that $|z_{1}|<|z_{0}|<|z_{2}|$ and both
$\sum_{\alpha\in D, \; \alpha\ge 0}a_{\alpha}z_{2}^{\alpha}$
and 
$\sum_{\alpha\in D, \; \alpha< 0}a_{\alpha}z_{1}^{\alpha}$
are absolutely convergent. Let $r_{1}$ and $r_{2}$ be numbers 
such that $|z_{1}|<r_{1}<|z_{0}|<r_{2}<|z_{2}|$. 
Since 
\[
\lim_{\alpha\to \infty}~^{\alpha}\!\!\!\sqrt{\alpha}=1,
\]
we can find $M>0$ such that 
\[
^{\alpha}\!\!\!\sqrt{\alpha}<\min\left(\frac{|z_{2}|}{r_{2}}, 
\frac{r_{1}}{|z_{1}|}\right)
\]
when  $\alpha>M$. 
But when $r_{1}<|z|<r_{2}$, 
\[
^{\alpha}\!\!\!\sqrt{\alpha}<\min\left(\frac{|z_{2}|}{r_{2}}, 
\frac{r_{1}}{|z_{1}|}\right)< \min\left(\frac{|z_{2}|}{|z|}, 
\frac{|z|}{|z_{1}|}\right)
\]
for $\alpha>M$, so for $z$ in the open subset and satisfying
$r_{1}<|z|<r_{2}$, we have
\begin{eqnarray}\label{po-ser-an-1}
\sum_{\alpha\in D, \; \alpha>M}|a_{\alpha}\alpha z^{\alpha}|
&=&\sum_{\alpha\in D, \; \alpha>M}|a_{\alpha}|
|^{\alpha}\!\!\!\sqrt{\alpha}z|^{\alpha}\nn
&\leq&\sum_{\alpha\in D, \; \alpha>M}|a_{\alpha}z_{2}^{\alpha}|
\end{eqnarray}
and 
\begin{eqnarray}\label{po-ser-an-2}
\sum_{\alpha\in D, \; \alpha< -M}|a_{\alpha}\alpha z^{\alpha}|
&=&\sum_{\alpha\in D, \; \alpha<-M}|a_{\alpha}|
\left|\frac{^{-\alpha}\!\!\!\sqrt{-\alpha}}{z}\right|^{-\alpha}\nn
&\leq&\sum_{\alpha\in D, \; \alpha<-M}|a_{\alpha}z_{1}^{\alpha}|.
\end{eqnarray}
On the other hand, we have
\begin{eqnarray}\label{po-ser-an-3}
\sum_{\alpha\in D, \; 0\le \alpha\le M}|a_{\alpha}\alpha z^{\alpha}|
&\le& M\sum_{\alpha\in D, \; 0\le \alpha\le M}|a_{\alpha} z^{\alpha}|\nn
&\le & M\sum_{\alpha\in D, \; 0\le \alpha\le M}|a_{\alpha} z_{2}^{\alpha}|
\end{eqnarray}
and 
\begin{eqnarray}\label{po-ser-an-4}
\sum_{\alpha\in D, \; 0>\alpha\ge  -M}|a_{\alpha}\alpha z^{\alpha}|
&\le &M\sum_{\alpha\in D, \; 0>\alpha\ge  -M}|a_{\alpha} z^{\alpha}|\nn
&\le & M\sum_{\alpha\in D, \; 0>\alpha\ge  -M}|a_{\alpha} z_{1}^{\alpha}|.
\end{eqnarray}
{}From (\ref{po-ser-an-1})--(\ref{po-ser-an-4}), we see that 
$\sum_{\alpha\in D}a_{\alpha}\alpha z^{\alpha}$ is absolutely and uniformly
convergent in the neighborhood of $z_{0}$ consisting of $z$ in the 
open subset satisfying $r_{1}<|z|<r_{2}$. 
\epfv

\begin{propo}\label{discrete-exp-set}
For any discrete subset $D$ of $\mathbb{R}$, $D\times\{0,1,\dots,N\}$ is a
unique expansion set.
\end{propo}
\pf
Assume that 
\[
\sum_{\alpha\in D}\sum_{\beta=0}^{N}a_{\alpha, \beta}
z^\alpha (\log z)^\beta
\]
is absolutely convergent to $0$ on some nonempty open subset of
${\mathbb C}^{\times}$. We can assume that the open subset 
does not intersect the unit circle; otherwise, we can delete the 
intersection with the unit circle to obtain the open subset we want.
Then for $\beta=0, \dots, N$, $\sum_{\alpha\in
D}a_{\alpha, \beta} z^\alpha$ is absolutely convergent.  Let $z_{0}$
be in the open subset, so that we have
\[
\sum_{\alpha\in D}\sum_{\beta=0}^{N}
a_{\alpha, \beta}z_{0}^\alpha (\log z_{0})^\beta
\]
is absolutely convergent to $0$.
Then for any $z$ such that $|z|=|z_{0}|$, even if it is not
in the open subset,
\begin{eqnarray*}
\sum_{\alpha\in D}\sum_{\beta=0}^{N}
|a_{\alpha, \beta}z^\alpha(\log z)^\beta|
&=&
\sum_{\beta=0}^{N}\left(\sum_{\alpha\in D}
|a_{\alpha, \beta}
z_{0}^\alpha|\left|\frac{z}{z_{0}}\right|^\alpha\right)
|(\log z)^\beta|
\nn
&\le&\sum_{\beta=0}^{N}\left(\sum_{\alpha\in D}
|a_{\alpha, \beta}z_{0}^\alpha|\right) |(\log z)^\beta|.
\end{eqnarray*}
So for such $z$, 
\[
\sum_{\alpha\in D}\sum_{\beta=0}^{N}
a_{\alpha, \beta}z^\alpha(\log z)^\beta
\]
is also absolutely convergent. Since 
$|z|=|z_{0}|\ne 1$, $z\ne 1$, and so we conclude that
for $\beta=0, \dots, N$,
$\sum_{\alpha\in D}a_{\alpha, \beta}
z^\alpha$ is absolutely convergent. 
By Lemma \ref{po-ser-an},
$\sum_{\alpha\in D}a_{\alpha, \beta}
z^\alpha$ represents an analytic function at $z$. Thus we see that 
\[
\sum_{\alpha\in D}\sum_{\beta=0}^{N}a_{\alpha, \beta}
z^\alpha (\log z)^\beta=\sum_{\beta=0}^{N}
\left(\sum_{\alpha\in D}a_{\alpha, \beta}
z^\alpha\right) (\log z)^\beta
\]
is also analytic at $z$ for $z\in \C$ such that 
$|z|=|z_{0}|$ for some $z_{0}$ in the open subset. 
So we can use analytic extension. 
Since 
\[
\sum_{\alpha\in D}\sum_{\beta=0}^{N}a_{\alpha, \beta}
z^\alpha (\log z)^\beta=0
\]
in the nonempty open subset, its analytic extension is also 
$0$. Using analytic extension,  we conclude that for $z$ as above,
\[
\sum_{\alpha\in D}\sum_{\beta=0}^{N}
a_{\alpha, \beta}z^\alpha e^{2p\pi \alpha i}(\log z+2p\pi i)^\beta
\]
is absolutely convergent to $0$ for $p\in \Z$. 

Now writing $z=re^{i\theta}$ with $r>0$ and $0\le \theta<2\pi$, we obtain
that 
\[
\sum_{\alpha\in D}\sum_{\beta=0}^{N}\sum_{k=0}^{\beta}{\beta\choose k}
a_{\alpha, \beta}
r^\alpha (\log r)^{\beta-k} i^{k}(\theta+2p\pi)^k e^{i\alpha (\theta+2p\pi)}
\]
is absolutely convergent to $0$ for $0\le \theta<2 \pi$ and $p\in \Z$,
or equivalently, that
\begin{equation}\label{discrete-exp-set-1}
\sum_{\alpha\in D}\sum_{\beta=0}^{N}\sum_{k=0}^{\beta}{\beta\choose k}
a_{\alpha, \beta}
r^\alpha (\log r)^{\beta-k} i^{k}\theta^k e^{i\alpha \theta}
\end{equation}
is absolutely convergent to $0$ for $\theta\in \R$. 
In particular, for $M>0$, we have 
\begin{equation}\label{discrete-exp-set-1.5}
\int_{0}^{M}e^{-s\theta}\left(\sum_{\alpha\in D}\sum_{\beta=0}^{N}
\sum_{k=0}^{\beta}{\beta\choose k}
a_{\alpha, \beta}
r^\alpha (\log r)^{\beta-k} i^{k}\theta^k e^{i\alpha \theta}\right)d\theta=0
\end{equation}
for $s\in \C$.

Since when $0\le \theta \le M$,
\begin{eqnarray*}
\lefteqn{\sum_{\alpha\in D}\sum_{\beta=0}^{N}\left|\sum_{k=0}^{\beta}{\beta\choose k}
a_{\alpha, \beta}
r^\alpha (\log r)^{\beta-k} i^{k}\theta^k e^{i\alpha \theta}\right|}\nn
&&\le \sum_{\alpha\in D}\sum_{\beta=0}^{N}\sum_{k=0}^{\beta}{\beta\choose k}
|a_{\alpha, \beta}|
r^\alpha |(\log r)^{\beta-k}| M^k \nn
&&=\sum_{\beta=0}^{N}\sum_{k=0}^{\beta}{\beta\choose k}
\left(\sum_{\alpha\in D}|a_{\alpha, \beta}|
r^\alpha\right) |(\log r)^{\beta-k}| M^k \nn
&&<\infty,
\end{eqnarray*}
(\ref{discrete-exp-set-1}) is in fact 
uniformly convergent with respect to $\theta$ when $0\le \theta \le M$.
So the left-hand side of (\ref{discrete-exp-set-1.5}) 
can be calculated term by term and, when $s\ne i\alpha$ for $\alpha\in D$, we obtain
\begin{eqnarray}\label{discrete-exp-set-1.7}
\sum_{\alpha\in D}\sum_{\beta=0}^{N}\sum_{k=0}^{\beta}
\Biggl({\beta\choose k}\frac{a_{\alpha, \beta}
r^\alpha (\log r)^{\beta-k} i^{k} k!}{(s-i\alpha)^{k+1}}
-\sum_{j=0}^{k}
{\beta\choose k}
\frac{a_{\alpha, \beta}
r^\alpha (\log r)^{\beta-k} i^{k}k!M^{k-j}e^{-(s-i\alpha)M}}{(k-j)!
(s-i\alpha)^{j+1}}\Biggr)
=0.\nn
\end{eqnarray}
(To compute $\int \theta^k e^{(i\alpha-s)\theta}d\theta$, one can for
example simply extract the coefficient of $u^k$, $u$ a formal
variable, {}from $\int e^{(i\alpha-s-u)\theta}d\theta$.)
Now assume that $\Re{(s)}>0$. Then
\begin{eqnarray*}
\lefteqn{\left|\sum_{\alpha\in D}\sum_{\beta=0}^{N}\sum_{k=0}^{\beta}
{\beta\choose k}\frac{a_{\alpha, \beta}
r^\alpha (\log r)^{\beta-k}i^{k}k!}{(s-i\alpha)^{k+1}}\right|}\nn
&&\le \sum_{\beta=0}^{N}\sum_{k=0}^{\beta}
{\beta\choose k}\left(\sum_{\alpha\in D}\frac{|a_{\alpha, \beta}|
r^\alpha }{|s-i\alpha|^{k+1}}\right)|(\log r)^{\beta-k}|k!\nn
&&\le \sum_{\beta=0}^{N}\sum_{k=0}^{\beta}
{\beta\choose k}\left(\sum_{\alpha\in D}|a_{\alpha, \beta}|
r^\alpha \right)\frac{|(\log r)^{\beta-k}|k!}{(\Re{(s)})^{k+1}},
\end{eqnarray*}
so that
\[
\sum_{\alpha\in D}\sum_{\beta=0}^{N}\sum_{k=0}^{\beta}
{\beta\choose k}\frac{a_{\alpha, \beta}
r^\alpha (\log r)^{\beta-k}i^{k}k!}{(s-i\alpha)^{k+1}}
\] 
is absolutely convergent. Thus (\ref{discrete-exp-set-1.7})
can be rewritten as 
\begin{eqnarray}\label{discrete-exp-set-2}
\lefteqn{\sum_{\alpha\in D}\sum_{\beta=0}^{N}\sum_{k=0}^{\beta}
{\beta\choose k}\frac{a_{\alpha, \beta}
r^\alpha (\log r)^{\beta-k}i^{k}k!}{(s-i\alpha)^{k+1}}}\nn
&&\quad 
-\sum_{\alpha\in D}\sum_{\beta=0}^{N}\sum_{k=0}^{\beta}\sum_{j=0}^{k}
{\beta\choose k}
\frac{a_{\alpha, \beta}
r^\alpha (\log r)^{\beta-k} i^{k}k!M^{k-j}e^{-(s-i\alpha)M}}{(k-j)!
(s-i\alpha)^{j+1}}\nn
&&=0.\nn
\end{eqnarray}

We claim that the limit of the second term on the 
left-hand side of 
(\ref{discrete-exp-set-2}) when $M$ goes to $\infty$ is $0$.
In fact, since for any $j\in \N$,
\begin{eqnarray*}
\left|\sum_{\alpha\in D}
\frac{a_{\alpha, \beta}
r^\alpha e^{i\alpha M}}{(s-i\alpha)^{j+1}}\right|
&\le& \sum_{\alpha\in D}
\frac{|a_{\alpha, \beta}|
r^\alpha}{|s-i\alpha|^{j+1}}\nn
&\le& \frac{1}{(\Re{(s)})^{j+1}}\sum_{\alpha\in D}
|a_{\alpha, \beta}| r^\alpha,
\end{eqnarray*}
we see that \[
\sum_{\alpha\in D}
\frac{a_{\alpha, \beta}
r^\alpha e^{i\alpha M}}{(s-i\alpha)^{j+1}}
\]
is absolutely convergent. Then we have
\begin{eqnarray*}
\lefteqn{\lim_{M\to \infty}\left|
\sum_{\alpha\in D}\sum_{\beta=0}^{N}\sum_{k=0}^{\beta}\sum_{j=0}^{k}
{\beta\choose k}
\frac{a_{\alpha, \beta}
r^\alpha (\log r)^{\beta-k} i^k k!M^{k-j}e^{-(s-i\alpha)M}}{(k-j)!(s-i\alpha)^{j+1}}
\right|}\nn
&&=\lim_{M\to \infty}|e^{-sM}|\left|
\sum_{\beta=0}^{N}\sum_{k=0}^{\beta}\sum_{j=0}^{k}
{\beta\choose k}\left(\sum_{\alpha\in D}
\frac{a_{\alpha, \beta}
r^\alpha e^{i\alpha M}}{(s-i\alpha)^{j+1}}\right)
\frac{(\log r)^{\beta-k} i^k k!M^{k-j}}{(k-j)!}
\right|\nn
&&= 0.
\end{eqnarray*}
Taking limits $\lim_{M\to \infty}$ on both sides 
of (\ref{discrete-exp-set-2}), we obtain
\begin{equation}\label{discrete-exp-set-3}
\sum_{\alpha\in D}\sum_{\beta=0}^{N}\sum_{k=0}^{\beta}
{\beta\choose k}\frac{a_{\alpha, \beta}
r^\alpha (\log r)^{\beta-k} i^k k!}{(s-i\alpha)^{k+1}}=0.
\end{equation}

Fix $\alpha_{0}\in D$. 
Since $D$ is discrete, we can find $\delta>0$ such that 
$|\alpha-\alpha_{0}|>\delta$ for $\alpha \in D\setminus 
\{\alpha_{0}\}$.
When $|s-i\alpha_{0}|<\frac{\delta}{2}$, we have 
$|s-i\alpha|>\frac{\delta}{2}$
for $\alpha \in D\setminus \{\alpha_{0}\}$.
Thus when $|s-i\alpha_{0}|<\frac{\delta}{2}$, even if $\Re{(s)}\le 0$,
\begin{eqnarray*}
\lefteqn{\sum_{\alpha\in D\setminus \{\alpha_{0}\}}\sum_{\beta=0}^{N}
\left|\sum_{k=0}^{\beta}
{\beta\choose k}
\frac{a_{\alpha, \beta}r^\alpha (\log r)^{\beta-k}i^{k}k!}{(s-i\alpha)^{k+1}}\right|}\nn
&&\le \sum_{\alpha\in D\setminus \{\alpha_{0}\}}\sum_{\beta=0}^{N}
\sum_{k=0}^{\beta}
{\beta\choose k}
\frac{|a_{\alpha, \beta}|r^\alpha |(\log r)^{\beta-k}|k!}{|s-i\alpha|^{k+1} }\nn
&&\le \sum_{\alpha\in D\setminus \{\alpha_{0}\}}\sum_{\beta=0}^{N}
\sum_{k=0}^{\beta}
{\beta\choose k}
\frac{2^{k+1}|a_{\alpha, \beta}|r^\alpha |(\log r)^{\beta-k}|k!}{\delta^{k+1}}\nn
&&=\sum_{\beta=0}^{N}
\sum_{k=0}^{\beta}
{\beta\choose k}\frac{2^{k+1}|(\log r)^{\beta-k}|k!}{\delta^{k+1}}
\left(\sum_{\alpha\in D\setminus \{\alpha_{0}\}}
|a_{\alpha, \beta}|r^\alpha\right)\nn
&&<\infty.
\end{eqnarray*}
This shows that for such $s$ (even if $\Re{(s)}\le 0$),
\begin{equation}\label{discrete-exp-set-4}
\sum_{\alpha\in D\setminus \{\alpha_{0}\}}\sum_{\beta=0}^{N}
\sum_{k=0}^{\beta}
{\beta\choose k}
\frac{a_{\alpha, \beta}r^\alpha(\log r)^{\beta-k}i^{k}k!}{(s-i\alpha)^{k+1}}
\end{equation}
is absolutely and uniformly convergent. 
In particular, taking $s=i\alpha_{0}$, we see that
\[
\sum_{\alpha\in D\setminus \{\alpha_{0}\}}\sum_{\beta=0}^{N}
\sum_{k=0}^{\beta}
{\beta\choose k}
\frac{a_{\alpha, \beta}r^\alpha(\log r)^{\beta-k}i^{k}k!}{(i\alpha_{0}-i\alpha)^{k+1}}
\]
is absolutely convergent. Moreover, since 
(\ref{discrete-exp-set-4}) is uniformly convergent,
we have
\begin{eqnarray*}
\lefteqn{\lim_{s\to i\alpha_{0}}
\sum_{\alpha\in D\setminus \{\alpha_{0}\}}\sum_{\beta=0}^{N}
\sum_{k=0}^{\beta}
{\beta\choose k}
\frac{a_{\alpha, \beta}r^\alpha(\log r)^{\beta-k}i^{k}k!}{(s-i\alpha)^{k+1}}}\nn
&&=\sum_{\alpha\in D\setminus \{\alpha_{0}\}}\sum_{\beta=0}^{N}
\sum_{k=0}^{\beta}
{\beta\choose k}
\frac{a_{\alpha, \beta}r^\alpha(\log r)^{\beta-k}i^{k}k!}{(i\alpha_{0}-i\alpha)^{k+1}}.
\end{eqnarray*}
Thus by (\ref{discrete-exp-set-3}),
\begin{eqnarray*}
\lefteqn{\lim_{s\to i\alpha_{0}}\sum_{\beta=0}^{N}\sum_{k=0}^{\beta}
{\beta\choose k}\frac{a_{\alpha_{0}, \beta}
r^{\alpha_{0}} (\log r)^{\beta-k}i^{k}k!}{(s-i\alpha_{0})^{k+1}}}\nn
&&=\lim_{s\to i\alpha_{0}}
\sum_{\alpha\in D}\sum_{\beta=0}^{N}
\sum_{k=0}^{\beta}
{\beta\choose k}
\frac{a_{\alpha, \beta}r^\alpha(\log r)^{\beta-k}i^{k}k!}{(s-i\alpha)^{k+1}}\nn
&&\quad -\lim_{s\to i\alpha_{0}}
\sum_{\alpha\in D\setminus \{\alpha_{0}\}}\sum_{\beta=0}^{N}
\sum_{k=0}^{\beta}
{\beta\choose k}
\frac{a_{\alpha, \beta}r^\alpha(\log r)^{\beta-k}i^{k}k!}{(s-i\alpha)^{k+1}}\nn
&&=-
\sum_{\alpha\in D\setminus \{\alpha_{0}\}}\sum_{\beta=0}^{N}
\sum_{k=0}^{\beta}
{\beta\choose k}
\frac{a_{\alpha, \beta}r^\alpha
(\log r)^{\beta-k}i^{k}k!}{(i\alpha_{0}-i\alpha)^{k+1}}.
\end{eqnarray*}
In particular, this limit exists and is finite 
for $r$ in a nonempty open subset of the positive 
real line.  Thus for such $r$, the coefficient of 
$(s-i\alpha_{0})^{-(k+1)}$ in
\[
\sum_{k=0}^{N}\sum_{\beta=k}^{N}
{\beta\choose k}\frac{a_{\alpha_{0}, \beta}
r^{\alpha_{0}} (\log r)^{\beta-k}i^{k}k!}{(s-i\alpha_{0})^{k+1}}
\]
is 0 for $k=0, \dots, N$, and in particular for $k=0$, and since the
powers of $\log r$ are linearly independent, we have $a_{\alpha_{0},
\beta}=0$ for $\beta=0, \dots, N$, as desired.  \epfv

We will also need the following:
\begin{propo}\label{double-conv<=>iterate-conv}
Let $D$ be a subset of $\R$ and $N$ a nonnegative integer. Then the series 
\begin{equation}\label{double-series}
\sum_{\alpha\in D}\sum_{\beta=0}^{N}a_{\alpha, \beta}
z^\alpha (\log z)^\beta
\end{equation}
is absolutely convergent on some (nonempty) open subset of $\C^{\times}$
if and only if the series
\begin{equation}\label{iterate-series}
\sum_{\alpha\in D}\left(\sum_{\beta=0}^{N}a_{\alpha, \beta}
(\log z)^\beta\right) z^\alpha
\end{equation}
and the corresponding series of first and higher derivatives with respect to 
$z$,
viewed as series whose terms are the expressions
\[
\left(\sum_{\beta=0}^{N}a_{\alpha, \beta}
(\log z)^\beta\right) z^\alpha
\]
and their derivatives with respect to $z$, are absolutely convergent
on the same open subset.  The series of derivatives of
(\ref{iterate-series}) have the same format as (\ref{iterate-series}),
except that for the $n$-th derivative, the outer sum is over the set
$D-n$ and the inner sum has new coefficients in $\C$.
\end{propo}
\pf The last assertion is clear.

We know that the absolute convergence of (\ref{double-series}) and
its (higher) derivatives implies the absolute convergence of
(\ref{iterate-series}) and its derivatives.  If (\ref{double-series})
is absolutely convergent, then so are its derivatives, as we see by using
Lemma \ref{po-ser-an}.  Thus (\ref{iterate-series}) and its
derivatives are absolutely convergent.

Conversely, assume that (\ref{iterate-series}) and its derivatives are
absolutely convergent.  We need to show that (\ref{double-series}) is
absolutely convergent at any $z_{0}$ in the open subset.  We consider
the series
\begin{equation}\label{iterate-series-1}
\sum_{\alpha\in D, \;\alpha\ge 0}\left(\sum_{\beta=0}^{N}a_{\alpha, \beta}
z_{2}^\beta\right) z_{1}^\alpha
\end{equation}
of functions
\[
\left(\sum_{\beta=0}^{N}a_{\alpha, \beta}
z_{2}^\beta\right) z_{1}^\alpha
\]
in two variables $z_{1}$ and $z_{2}$. Since $z_{0}$ is in the open
subset, we can find a smaller open subset inside the original one such
that for $z$ in this smaller one, $|z_{0}|<|z|$ and $|\log
z_{0}|<|\log z|$. We know that the series (\ref{iterate-series-1}) is
absolutely convergent when $z_{1}=z$, $z_{2}=\log z$ and $z$ is in the
original open subset.  For any $z_{1}$ and $z_{2}$ satisfying
$0<|z_{1}|<|z|$ and $z_{2}=\log z$ where $z$ is in the smaller open
subset,
\[
\sum_{\alpha\in D, \;\alpha\ge 0}\left|\sum_{\beta=0}^{N}a_{\alpha, \beta}
z_{2}^\beta\right| |z_{1}^\alpha|
\le \sum_{\alpha\in D, \;\alpha\ge 0}\left|\sum_{\beta=0}^{N}a_{\alpha, \beta}
(\log z)^\beta\right| |z^\alpha|
\]
is convergent. So in this case (\ref{iterate-series-1}) is absolutely
convergent. Since for any fixed $z_{2}=\log z$ 
where $z$ is in the smaller open subset, the numbers $z_{1}$ satisfying 
$0<|z_{1}|<|z|$ form an open subset, we can apply Lemma \ref{po-ser-an}
to obtain that
\begin{equation}\label{iterate-series-2}
\sum_{\alpha\in D, \;\alpha\ge 0}\frac{\partial}{\partial z_{1}}
\left(\left(\sum_{\beta=0}^{N}a_{\alpha, \beta}
z_{2}^\beta\right)  z_{1}^{\alpha}\right)
\end{equation}
is also absolutely convergent for any $z_{1}$ and $z_{2}$ satisfying
$0<|z_{1}|<|z|$ and $z_{2}=\log z$ where $z$ is in the smaller open
subset.

Also, by assumption, 
\begin{equation}\label{iterate-series-3}
\sum_{\alpha\in D, \; \alpha\ge 0}\left(\frac{\partial}{\partial z_{1}}
+\frac{1}{z_{1}}
\frac{\partial}{\partial z_{2}}\right)\left(
\left(\sum_{\beta=0}^{N}a_{\alpha, \beta}
z_{2}^{\beta}\right) z_{1}^\alpha\right)
\end{equation}
is absolutely convergent when $z_{1}=z$ and $z_{2}=\log z$ when $z$ is
in the original open subset.  For $z$ in the smaller open subset
and any $z_{1}$ and $z_{2}$ satisfying $0<|z_{1}|<|z|$ and $z_{2}=\log z$,
\begin{eqnarray*}
\lefteqn{\sum_{\alpha\in D, \; \alpha\ge 0}\left|
\left(\frac{\partial}{\partial z_{1}}
+\frac{1}{z_{1}}
\frac{\partial}{\partial z_{2}}\right)\left(
\left(\sum_{\beta=0}^{N}a_{\alpha, \beta}
z_{2}^{\beta}\right) z_{1}^\alpha\right)\right|}\nn
&&=\sum_{\alpha\in D, \; \alpha\ge 0}\left|
\left(\left(\sum_{\beta=0}^{N}a_{\alpha, \beta}
z_{2}^{\beta}\right) \alpha z_{1}^{\alpha-1}\right)
+\left(\left(\sum_{\beta=0}^{N}a_{\alpha, \beta}\beta
z_{2}^{\beta-1}\right) z_{1}^{\alpha-1}\right)\right|\nn
&&=\sum_{\alpha\in D, \; \alpha\ge 0}\left|
\left(\left(\sum_{\beta=0}^{N}a_{\alpha, \beta}
z_{2}^{\beta}\right)\alpha 
+\sum_{\beta=0}^{N}a_{\alpha, \beta}\beta
z_{2}^{\beta-1} \right)\right|\left|z_{1}^{\alpha-1}\right|\nn
&&\le |z {z_1}^{-1}|
\sum_{\alpha\in D, \; \alpha\ge 0}\left|
\left(\left(\sum_{\beta=0}^{N}a_{\alpha, \beta}
(\log z)^{\beta}\right)\alpha 
+\sum_{\beta=0}^{N}a_{\alpha, \beta}\beta
(\log z)^{\beta-1} \right)\right|\left|z^{\alpha-1}\right|\nn
&&=|z {z_1}^{-1}|
\sum_{\alpha\in D, \; \alpha\ge 0}\left|
\left(\left(\sum_{\beta=0}^{N}a_{\alpha, \beta}
(\log z)^{\beta}\right)\alpha z^{\alpha-1}
+\left(\sum_{\beta=0}^{N}a_{\alpha, \beta}\beta
(\log z)^{\beta-1}\right)z^{\alpha-1} \right)\right|\nn
&&=|z {z_1}^{-1}|
\sum_{\alpha\in D, \; \alpha\ge 0}\left|\frac{\partial}{\partial z}
\left(\left(\sum_{\beta=0}^{N}a_{\alpha, \beta}
(\log z)^{\beta}\right) 
z^{\alpha}\right)\right|
\end{eqnarray*}
(where we keep in mind that $\alpha - 1$ could be negative)
is convergent, so that (\ref{iterate-series-3}) is absolutely convergent
for such $z_{1}$ and $z_{2}$.
Thus, subtracting, we see that
\begin{equation}\label{iterate-series-4}
\sum_{\alpha\in D, \;\alpha\ge 0}
\frac{\partial}{\partial z_{2}}\left(
\left(\sum_{\beta=0}^{N}a_{\alpha, \beta}
z_{2}^{\beta}\right) z_{1}^\alpha\right)=
\sum_{\alpha\in D, \;\alpha\ge 0}
\left(\frac{\partial}{\partial z_{2}}
\left(\sum_{\beta=0}^{N}a_{\alpha, \beta}
z_{2}^{\beta}\right)\right) z_{1}^\alpha
\end{equation}
is also absolutely convergent 
for such $z_{1}$ and $z_{2}$. By Lemma \ref{po-ser-an},
\[
\sum_{\alpha\in D, \;\alpha\ge 0}
\frac{\partial}{\partial z_{1}}
\frac{\partial}{\partial z_{2}}\left(
\left(\sum_{\beta=0}^{N}a_{\alpha, \beta}
z_{2}^{\beta}\right) z_{1}^\alpha\right)
\]
is absolutely convergent 
for such $z_{1}$ and $z_{2}$.

Since $\frac{\partial}{\partial z_{1}}
+\frac{1}{z_{1}}
\frac{\partial}{\partial z_{2}}$ and $\frac{\partial}{\partial z_{2}}$
commute with each other,  we have 
\begin{eqnarray}\label{iterate-series-4.5}
\lefteqn{\sum_{\alpha\in D, \;\alpha\ge 0}
\left(\frac{\partial}{\partial z_{1}}
+\frac{1}{z_{1}}
\frac{\partial}{\partial z_{2}}\right)
\frac{\partial}{\partial z_{2}}\left(
\left(\sum_{\beta=0}^{N}a_{\alpha, \beta}
z_{2}^{\beta}\right) z_{1}^\alpha\right)}\nn
&&=\sum_{\alpha\in D, \;\alpha\ge 0}\frac{\partial}{\partial z_{2}}
\left(\frac{\partial}{\partial z_{1}}
+\frac{1}{z_{1}}
\frac{\partial}{\partial z_{2}}\right)
\left(\left(\sum_{\beta=0}^{N}a_{\alpha, \beta}
z_{2}^{\beta}\right) z_{1}^{\alpha}\right)\nn
&&=z_{1}\sum_{\alpha\in D, \;\alpha\ge 0}
\left(\frac{\partial}{\partial z_{1}}
+\frac{1}{z_{1}}
\frac{\partial}{\partial z_{2}}\right)^{2}
\left(\left(\sum_{\beta=0}^{N}a_{\alpha, \beta}
z_{2}^{\beta}\right) z_{1}^{\alpha}\right)\nn
&&\quad-z_{1}\sum_{\alpha\in D, \;\alpha\ge 0}\frac{\partial}{\partial z_{1}}
\left(\frac{\partial}{\partial z_{1}}
+\frac{1}{z_{1}}
\frac{\partial}{\partial z_{2}}\right)
\left(\left(\sum_{\beta=0}^{N}a_{\alpha, \beta}
z_{2}^{\beta}\right) z_{1}^{\alpha}\right).
\end{eqnarray}
By assumption, the first term on the right-hand side of 
(\ref{iterate-series-4.5})
is absolutely convergent when $z_{1}=z$ and $z_{2}=\log z$
and $z$ is in the original open subset, and then, by the same argument
as above,  is also absolutely 
convergent for $z_{1}$ and $z_{2}$ satisfying $0<|z_{1}|<|z|$,
$z_{2}=\log z$ and $z$ in the smaller open subset. 
By Lemma \ref{po-ser-an} and the absolute convergence of 
(\ref{iterate-series-3}) for such $z_{1}$ and $z_{2}$, the second 
term on the right-hand side of 
(\ref{iterate-series-4.5})
is also absolutely convergent for $z_{1}$ and $z_{2}$ 
satisfying $0<|z_{1}|<|z|$,
$z_{2}=\log z$ and $z$ in the smaller open subset. 
So the left-hand side of (\ref{iterate-series-4.5}) is absolutely 
convergent for such $z_{1}$ and $z_{2}$. Thus 
\begin{eqnarray*}
\lefteqn{\sum_{\alpha\in D, \;\alpha\ge 0}
\left(\left(\frac{\partial}{\partial z_{2}}\right)^{2}
\left(\sum_{\beta=0}^{N}a_{\alpha, \beta}
z_{2}^{\beta}\right)\right) z_{1}^\alpha}\nn
&&=\sum_{\alpha\in D, \;\alpha\ge 0}
\left(\frac{\partial}{\partial z_{2}}\right)^{2}\left(
\left(\sum_{\beta=0}^{N}a_{\alpha, \beta}
z_{2}^{\beta}\right) z_{1}^\alpha\right)\nn
&&=z_{1}\sum_{\alpha\in D, \;\alpha\ge 0}
\left(\frac{\partial}{\partial z_{1}}
+\frac{1}{z_{1}}
\frac{\partial}{\partial z_{2}}\right)
\frac{\partial}{\partial z_{2}}\left(
\left(\sum_{\beta=0}^{N}a_{\alpha, \beta}
z_{2}^{\beta}\right) z_{1}^\alpha\right)\nn
&&\quad -z_{1}\sum_{\alpha\in D, \;\alpha\ge 0}
\frac{\partial}{\partial z_{1}}
\frac{\partial}{\partial z_{2}}\left(
\left(\sum_{\beta=0}^{N}a_{\alpha, \beta}
z_{2}^{\beta}\right) z_{1}^\alpha\right)
\end{eqnarray*}
is absolutely convergent for $z$ in the smaller open subset
and any $z_{1}$ and $z_{2}$ satisfying $0<|z_{1}|<|z|$ and 
$z_{2}=\log z$.

Repeating these arguments, we obtain that 
\begin{equation}\label{iterate-series-5}
\sum_{\alpha\in D, \;\alpha\ge 0}
\left(\left(\frac{\partial}{\partial z_{2}}\right)^{k}
\left(\sum_{\beta=0}^{N}a_{\alpha, \beta}
z_{2}^{\beta}\right)\right) z_{1}^\alpha
\end{equation}
is absolutely convergent for such $z_{1}$ and $z_{2}$ and for $k\in
\N$.  Taking $k=N$, we see that
\[
\sum_{\alpha\in D, \; \alpha\ge 0}
a_{\alpha, N}z_{1}^\alpha
\]
is absolutely convergent for such $z_{1}$.
Continuing this process with $k=N-1, \dots, 0$ we obtain that 
\[
\sum_{\alpha\in D, \; \alpha\ge 0}
a_{\alpha, \beta}z_{1}^\alpha
\]
is absolutely convergent for such $z_{1}$ and each $\beta=0,
\dots, N$. Since $0<|z_{0}|<|z|$, we see that in the case $z_{1}=z_{0}$,
\[
\sum_{\alpha\in D, \;\alpha\ge 0}
a_{\alpha, \beta}z_{0}^\alpha
\]
is absolutely convergent for $\beta=0, \dots, N$.  

We also need to prove the absolute convergence of
\[
\sum_{\alpha\in D, \;\alpha< 0}
a_{\alpha, \beta}z_{0}^\alpha
\]
for $\beta=0, \dots, N$. The proof is completely analogous to the 
proof above except that 
we take a smaller open subset such that for $z$ in this smaller 
one, $|z_{0}|> |z|>0$ and $|\log z_{0}|>|\log z|$ instead of 
$|z_{0}|< |z|$ and $|\log z_{0}|<|\log z|$. 
Thus
\[
\sum_{\alpha\in D}
a_{\alpha, \beta}z_{0}^\alpha
\]
is absolutely convergent for $\beta=0, \dots, N$.  This absolute
convergence is equivalent to the absolute convergence of
(\ref{double-series}).  \epfv

\begin{assum}\label{assum-exp-set}
Throughout the remainder of this work, we shall assume that ${\cal C}$
satisfies the condition that for any object $W$ of $\mathcal{C}$, 
the set $\{(n, i)\in \C\times \N\;|\; 
(L(0)-n)^{i}W_{(n)}\ne 0\}$ is a unique expansion set. 
We also assume that for any objects $W_{1}$, $W_{2}$ and
$W_{3}$ of $\mathcal{C}$, any logarithmic intertwining operator $\Y$
of type ${W_{3}\choose W_{1}W_{2}}$, and any $w_{(1)}\in W_{1}$ and
$w_{(2)}\in W_{2}$, the powers of $x$ and $\log x$ occurring in
\[\Y(w_{(1)}, x)w_{(2)}
\]
form a unique expansion set of the form $D\times \{0, \dots, N\}$
where $D$ is a set of real numbers.
\end{assum}

In practice, the following result guarantees that ``virtually all the 
interesting examples'' satisfy this assumption:

\begin{propo}
A sufficient condition for Assumption \ref{assum-exp-set} to hold is
that for any object of ${\cal C}$, the (generalized) weights form a
discrete set of real numbers and in addition there exists $K\in
\Z_{+}$ such that $(L(0)-L(0)_{s})^{K}=0$ on the module.
\end{propo}
\pf
This result follows immediately {}from Proposition \ref{discrete-exp-set}
and Proposition \ref{log:logwt}; see Remark \ref{log:compM}.
\epfv

Recall the projections $\pi_{p}$ for $p\in \C$ {}from (\ref{pi_n})
and Definition \ref{Wbardef}, and recall the notations
(\ref{iterateabbreviation}), (\ref{productabbreviation}). 
Using Assumption \ref{assum-exp-set} and Proposition
\ref{double-conv<=>iterate-conv},
we shall prove the following result:

\begin{propo}\label{prod=0=>comp=0}
Assume the convergence condition for intertwining maps
in $\mathcal{C}$ (recall Definition \ref{conv-conditions}).  Let $z_1$,
$z_2$ be two nonzero complex numbers satisfying 
\[
|z_1|>|z_2|>0,
\]
and let $I_1\in \mathcal{M}[P(z_{1})]_{W_{1}M_{1}}^{W_{4}}$ and $I_2\in
\mathcal{M}[P(z_{2})]_{W_{2}W_{3}}^{M_{1}}$.  Let $w_{(1)}\in W_1$, $w_{(3)}\in
W_3$ and $w'_{(4)}\in W'_4$ be homogeneous elements with respect to 
the (generalized) weight gradings. Suppose that for all homogeneous
$w_{(2)}\in W_2$,
\[
\langle w'_{(4)}, I_1(w_{(1)}\otimes I_2(w_{(2)}\otimes w_{(3)}))\rangle=0.
\]
Then
\[
\langle w'_{(4)},  I_1({w_{(1)}}\otimes \pi_p I_2({w_{(2)}}\otimes 
w_{(3)}))\rangle=0
\]
for all $p\in \C$ and all $w_{(2)}\in W_2$. In particular, 
\[
\langle w'_{(4)},  \pi_{p} I_1({w_{(1)}}\otimes \pi_q I_2({w_{(2)}}\otimes 
w_{(3)}))\rangle=0
\]
for all $p, q\in \C$ and all $w_{(2)}\in W_2$. 
\end{propo}
\pf Recall the correspondence between $P(z)$-intertwining maps and
logarithmic intertwining operators of the same type 
(Proposition \ref{im:correspond}),
and the 
notation $\mathcal{Y}_{I, p}$, $p\in \Z$, for the logarithmic intertwining
operators corresponding to a $P(z)$-intertwining map $I$ ((\ref{YIp}) and 
(\ref{recover})). 

For a new formal variable $y$, by (\ref{log:ck1}) and the
$L(-1)$-derivative property (\ref{log:L(-1)dev}),
\begin{eqnarray*}
\lefteqn{\langle w'_{(4)}, {\cal Y}_{I_1,0}(w_{(1)}, x_1){\cal
Y}_{I_2,0}(w_{(2)},x_2+y)w_{(3)}\rangle}\\
&&=\langle w'_{(4)}, {\cal
Y}_{I_1,0}(w_{(1)},x_1){\cal
Y}_{I_2,0}(e^{yL(-1)}w_{(2)},x_2)w_{(3)}\rangle\\
&&=\sum_{i\in {\mathbb N}}\frac{y^i}{i!}\langle w'_{(4)}, {\cal
Y}_{I_1,0}(w_{(1)},x_1){\cal
Y}_{I_2,0}(L(-1)^iw_{(2)},x_2)w_{(3)}\rangle.
\end{eqnarray*}
Hence,  for any
$\epsilon\in {\mathbb C}$ such that $|\epsilon|<|z_2|$ and
$0<|z_2+\epsilon|<|z_1|$, we have
\begin{eqnarray*}
\lefteqn{\langle w'_{(4)}, \mathcal{Y}_{I_1,0}(w_{(1)},
z_1)\mathcal{Y}_{I_2,0}(w_{(2)},z_2+\epsilon)w_{(3)}\rangle}\\
&&=\sum_{i\in {\mathbb N}}\frac{\epsilon^i}{i!}\langle w'_{(4)},
\mathcal{Y}_{I_1,0}(w_{(1)},z_1)
\mathcal{Y}_{I_2,0}(L(-1)^iw_{(2)},z_2)w_{(3)}\rangle\\
&&=\sum_{i\in {\mathbb N}}\frac{\epsilon^i}{i!}0=0.
\end{eqnarray*}

By suitably repeating this, we find that
\begin{equation}\label{w2z}
\langle w'_{(4)}, \mathcal{Y}_{I_1,0}(w_{(1)}, z_1)
\mathcal{Y}_{I_2,0}(w_{(2)},z)w_{(3)}\rangle=0
\end{equation}
for all homogeneous $w_{(2)}\in W_2$ and all $z\in {\mathbb C}$ with
$0<|z|<|z_1|$. Now by Assumption \ref{assum-exp-set},
Proposition \ref{log:logwt}(b) and the meaning of the absolutely 
convergent series on the
left-hand side of (\ref{w2z}), for fixed $z_{1}\ne 0$,
this series on
the left-hand side is of the form (\ref{iterate-series}), with $D$ a
unique expansion set.  By the
$L(-1)$-derivative property for logarithmic intertwining operators,
the higher-derivative series of the
left-hand side of (\ref{w2z}) are
absolutely convergent. So we can apply Propostion
\ref{double-conv<=>iterate-conv} to obtain that the double series obtained 
{}from the left-hand side of (\ref{w2z}) by taking the terms to be monomials in $z$ and
$\log z$ is also absolutely convergent to 0 for $0<|z|<|z_1|$.
By the definition of unique expansion 
set, we see that all of the coefficients of the monomials in $z$ and
$\log z$ of this double
series must be zero. Hence we get
\[
\langle w'_{(4)}, \mathcal{Y}_{I_1,0}(w_{(1)}, z_1)
({w_{(2)}}^{\mathcal{Y}_{I_2,0}}_{n;\,k}
w_{(3)})\rangle=0
\]
for any homogeneous $w_{(2)}\in W_2$, $n\in {\mathbb C}$ and $k\in \N$. 
Thus
\[
\langle w'_{(4)}, I_1(w_{(1)}\otimes 
\pi_p I_2({w_{(2)}}\otimes  w_{(3)}))
\rangle=0
\]
for any homogeneous $w_{(2)}\in W_2$ and $p\in {\mathbb C}$, in view of Proposition
\ref{log:logwt}(b) and Proposition \ref{im:correspond}, and this
remains true for any $w_{(2)}\in W_2$.  The last statement is clear.
\epf

\begin{corol}\label{prospan}
Assume the convergence condition for intertwining maps
in $\mathcal{C}$ (recall Definition \ref{conv-conditions}).  Let $z_1$,
$z_2$ be two nonzero complex numbers satisfying 
\[
|z_1|>|z_2|>0.
\]
Suppose that the $P(z_2)$-tensor product of $W_2$ and $W_3$ and the
$P(z_1)$-tensor product of $W_1$ and $W_2 \boxtimes_{P(z_2)} W_3$ both
exist (recall Definition \ref{pz-tp}).
Then $W_1\boxtimes_{P(z_1)} (W_2\boxtimes_{P(z_2)}
W_3)$ is spanned (as a vector space) by all the elements of the form
\[
\pi_{n}(w_{(1)}\boxtimes_{P(z_1)}(w_{(2)}\boxtimes_{P(z_2)}w_{(3)}))
\] 
where $w_{(1)}\in W_1$, $w_{(2)}\in W_2$ and $w_{(3)}\in W_3$ are 
homogeneous with respect to the (generalized) weight gradings 
and $n\in \mathbb{C}$ (recall the notation
(\ref{boxtensorofelements})). 
\end{corol}
\pf Let $w'_{(4)}\in (W_1\boxtimes_{P(z_1)} (W_2\boxtimes_{P(z_2)}
W_3))'$ be homogeneous such that
\[
\langle w'_{(4)}, w_{(1)}\boxtimes_{P(z_1)} (w_{(2)}\boxtimes_{P(z_2)}
w_{(3)})\rangle=0
\]
for all homogeneous 
$w_{(1)}\in W_1$, $w_{(2)}\in W_2$ and $w_{(3)}\in W_3$. {}From
Proposition \ref{prod=0=>comp=0} we see that
\[
\langle w'_{(4)}, \pi_p({w_{(1)}}\boxtimes_{P(z_1)}
\pi_q({w_{(2)}}\boxtimes_{P(z_2)} w_{(3)}))\rangle=0
\]
for all $p, q\in {\mathbb C}$ and all 
$w_{(1)}\in W_1$, $w_{(2)}\in W_2$ and
$w_{(3)}\in W_3$. Since by Proposition \ref{span}, the set
\[
\{ \pi_p({w_{(1)}}\boxtimes_{P(z_1)} \pi_q({w_{(2)}}\boxtimes_{P(z_2)}
w_{(3)})) |\;p,q\in {\mathbb C}, w_{(1)}\in W_1, w_{(2)}\in W_2,
w_{(3)}\in W_3\}
\]
spans the space $W_1\boxtimes_{P(z_1)} (W_2\boxtimes_{P(z_2)}
W_3)$, we must have $w'_{(4)}=0$, and the result follows.
\epfv

Analogously, by similar proofs we have:

\begin{propo}\label{iter=0=>comp=0}
Assume the convergence condition for intertwining maps in ${\cal
C}$. Let $z_0$, $z_2$ be two nonzero complex numbers satisfying
\[
|z_2|>|z_0|>0,
\] 
and let $I^1\in \mathcal{M}[P(z_{2})]_{M_{2}W_{3}}^{W_{4}}$ and $I^2\in
\mathcal{M}[P(z_{0})]_{W_{1}W_{2}}^{M_{2}}$. Let $w'_{(4)}\in W'_4$, $w_{(2)}\in
W_2$ and $w_{(3)}\in W_3$ be homogeneous with respect to the (generalized)
weight gradings. Suppose that for all homogeneous $w_{(1)}\in W_1$,
\[
\langle w'_{(4)}, I^1(I^2(w_{(1)}\otimes w_{(2)}) \otimes w_{(3)})\rangle=0.
\]
Then
\[
\langle w'_{(4)}, I^1(\pi_p I^2(w_{(1)}\otimes w_{(2)})\otimes w_{(3)})
\rangle=0
\]
for all $p\in \C$ and all $w_{(1)}\in W_1$. In particular, 
\[
\langle w'_{(4)}, \pi_{p} I^1(\pi_q I^2(w_{(1)}\otimes w_{(2)})\otimes w_{(3)})
\rangle=0
\]
for all $p, q\in \C$ and all $w_{(1)}\in W_1$. \epf
\end{propo}

\begin{corol}\label{iterspan}
Assume the convergence condition for intertwining maps in ${\cal
C}$. Let $z_0$, $z_2$ be two nonzero complex numbers satisfying
\[
|z_2|>|z_0|>0.
\] 
Suppose that the $P(z_0)$-tensor product of $W_1$ and $W_2$ and the
$P(z_2)$-tensor product of $W_1 \boxtimes_{P(z_0)} W_2$ and $W_3$ both
exist. Then $(W_1\boxtimes_{P(z_0)} W_2)\boxtimes_{P(z_2)}
W_3$ is spanned by all the elements of the form
\[
\pi_{n}((w_{(1)}\boxtimes_{P(z_0)} w_{(2)})\boxtimes_{P(z_2)} w_{(3)})
\]
where $w_{(1)}\in W_1$, $w_{(2)}\in W_2$ and $w_{(3)}\in W_3$ are
homogeneous and 
$n\in \C$.\epf
\end{corol}

\begin{rema}
{\rm One can generalize the convergence condition for two intertwining
maps and the results above to products and iterates of any number of
intertwining maps. The convergence conditions for three
intertwining maps and the spanning properties in the case of four
generalized modules will be needed in Section 12 in the proof of the
commutativity of the pentagon diagram, and we will discuss these
conditions and properties in Section 12. }
\end{rema}

\begin{rema}\label{weakly-abs-conv}
{\rm The convergence studied in this section can easily be formulated
as special cases of the following general notion: Let $W$ be a
(complex) vector space and let $\langle \cdot,\cdot \rangle:
W^{*}\times W\to \C$ be the pairing between the dual space $W^{*}$ and
$W$.  Consider the weak topology on $W^{*}$ defined by this pairing,
so that $W^{*}$ becomes a Hausdorff locally convex topological vector
space.  Let $\sum_{n\in \C}w^{*}_{n}$ be a formal series (indexed by
$\C$) in $W^{*}$. We say that $\sum_{n\in \C}w^{*}_{n}$ is {\it weakly
absolutely convergent} if for all $w\in W$, the formal series
\begin{equation}\label{sum-w}
\sum_{n\in \C}\langle w^{*}_{n}, w\rangle
\end{equation}
(again indexed by $\C$) of complex numbers is absolutely
convergent. Note that if $\sum_{n\in \C}w^{*}_{n}$ is weakly
absolutely convergent, then (\ref{sum-w}) as $w$ ranges through $W$
defines a (unique) element of $W^{*}$, and the formal series is in
fact convergent to this element in the weak topology. This element is
the {\it sum} of the series and is denoted using the same notation
$\sum_{n\in \C}w^{*}_{n}$. In this section, the convergence that we
have been discussing amounts to the weak absolute convergence of
formal series in $(W')^{*}$ for an object $W$ of $\mathcal{C}$, and
this kind of convergence will again be used in Section 8.  In Section
9, we will use this notion for $(W_{1}\otimes W_{2})^{*}$ where
$W_{1}$ and $W_{2}$ are generalized $V$-modules, and in Section 12, we
will be using more general cases.}
\end{rema}

\newpage

\setcounter{equation}{0}
\setcounter{rema}{0}

\section{$P(z_1,z_2)$-intertwining maps}

In this section we first prove some natural identities satisfied by
products and iterates of logarithmic intertwining operators and of
intertwining maps.  These identities were first proved in
\cite{tensor4} for intertwining operators and intertwining maps among
ordinary modules.  We also prove a list of identities relating
products of formal delta functions, as was done in \cite{tensor4}.
Using all these identities as motivation, we define
``$P(z_1,z_2)$-intertwining maps'' and study their basic properties,
by analogy with the relevant parts of the study of $P(z)$-intertwining
maps in Sections 4 and 5.  The notion of $P(z_1,z_2)$-intertwining map
is new; the treatment in this section is different {}from that in
\cite{tensor4}, even for the case of ordinary intertwining operators.

It is possible to define ``tensor products of three modules'' (as
opposed to iterated tensor products) based on the theory in this
section, and $P(z_1,z_2)$-intertwining maps play the same role for
``tensor products of three modules'' that $P(z)$-intertwining maps
play for tensor products of two modules.  But since we do not need
such ``tensor products of three modules'' in this work, we will not
formally introduce and study them.

Recall the Jacobi identity (\ref{log:jacobi}) in the definition of the
notion of logarithmic intertwining operator associated with
generalized modules $(W_1,Y_1)$, $(W_2,Y_2)$ and $(W_3,Y_3)$ for a
M\"obius (or conformal) vertex algebra $V$.  Suppose that we also have
generalized modules $(W_4,Y_4)$, $(M_1,Y_{M_1})$ and $(M_2,Y_{M_2})$.
Then {}from (\ref{log:jacobi}) we see that a product of logarithmic
intertwining operators of types ${W_4\choose W_1 M_1}$ and
${M_1\choose W_2 W_3}$ satisfies an identity analogous to
(\ref{log:jacobi}), as does an iterate of logarithmic intertwining
operators of types ${W_4\choose M_2 W_3}$ and ${M_2\choose W_1 W_2}$:

Let ${\cal Y}_1$ and ${\cal Y}_2$ be logarithmic intertwining
operators of types ${W_4\choose W_1 M_1}$ and ${M_1\choose W_2 W_3}$,
respectively.  Then for $v\in V$, $w_{(1)}\in W_1$, $w_{(2)}\in W_2$
and $w_{(3)}\in W_3$, the product of ${\cal Y}_1$ and ${\cal Y}_2$
satisfies the identity
\begin{eqnarray}
\lefteqn{\dlt{x_1}{x_0}{-y_1}\dlt{x_2}{x_0}{-y_2}Y_4(v,x_0){\cal
Y}_1(w_{(1)},y_1){\cal Y}_2(w_{(2)},y_2)w_{(3)}}\nno \\
&&= \dlt{y_1}{x_0}{-x_1}\dlt{x_2}{x_0}{-y_2}{\cal
Y}_1(Y_1(v,x_1)w_{(1)},y_1){\cal Y}_2(w_{(2)},y_2)w_{(3)}\nno \\
&& \quad+\dlt{x_1}{-y_1}{+x_0}\dlt{x_2}{x_0}{-y_2}{\cal
Y}_1(w_{(1)},y_1)Y_{M_1}(v,x_0){\cal Y}_2(w_{(2)},y_2)w_{(3)}\nno \\
&&= \dlt{y_1}{x_0}{-x_1}\dlt{x_2}{x_0}{-y_2}{\cal
Y}_1(Y_1(v,x_1)w_{(1)},y_1){\cal Y}_2(w_{(2)},y_2)w_{(3)}\nno \\
&& \quad+\dlt{x_1}{-y_1}{+x_0}\dlt{y_2}{x_0}{-x_2}{\cal
Y}_1(w_{(1)},y_1){\cal Y}_2(Y_2(v,x_2)w_{(2)},y_2)w_{(3)}\nno \\
&& \quad+\dlt{x_1}{-y_1}{+x_0}\dlt{x_2}{-y_2}{+x_0}{\cal
Y}_1(w_{(1)},y_1){\cal Y}_2(w_{(2)},y_2)Y_3(v,x_0)w_{(3)}.\nno \\
&& \label{Y12}
\end{eqnarray}
In addition, let ${\cal Y}^1$ and ${\cal Y}^2$ be logarithmic
intertwining operators of types ${W_4\choose M_2 W_3}$ and
${M_2\choose W_1 W_2}$, respectively.  Then for $v\in V$, $w_{(1)}\in
W_1$, $w_{(2)}\in W_2$ and $w_{(3)}\in W_3$, the iterate of ${\cal
Y}^1$ and ${\cal Y}^2$ satisfies the identity
\begin{eqnarray}
\lefteqn{\dlt{x_2}{x_0}{-y_2}\dlt{x_1}{x_2}{-y_0}Y_4(v,x_0){\cal
Y}^1({\cal Y}^2(w_{(1)},y_0)w_{(2)},y_2)w_{(3)}}\nno \\
&&= \dlt{y_2}{x_0}{-x_2}\dlt{x_1}{x_2}{-y_0}{\cal
Y}^1(Y_{M_2}(v,x_2){\cal Y}^2(w_{(1)},y_0)w_{(2)},y_2)w_{(3)}\nno
\\
&& \quad+\dlt{x_2}{-y_2}{+x_0}\dlt{x_1}{x_2}{-y_0}{\cal Y}^1({\cal
Y}^2(w_{(1)},y_0)w_{(2)},y_2)Y_3(v,x_0)w_{(3)}\nno \\
&&= \dlt{y_2}{x_0}{-x_2}\dlt{y_0}{x_2}{-x_1}{\cal Y}^1({\cal
Y}^2(Y_1(v,x_1)w_{(1)}),y_0)w_{(2)},y_2)w_{(3)}\nno \\
&& \quad+\dlt{y_2}{x_0}{-x_2}\dlt{x_1}{-y_0}{+x_2}{\cal Y}^1({\cal
Y}^2(w_{(1)},y_1)Y_2(v,x_2)w_{(2)},y_2)w_{(3)}\nno \\
&& \quad+\dlt{x_2}{-y_2}{+x_0}\dlt{x_1}{x_2}{-y_0}{\cal Y}^1({\cal
Y}^2(w_{(1)},y_0)w_{(2)},y_2)Y_3(v,x_0)w_{(3)}.\nno \\ &&
\label{Y34}
\end{eqnarray}

Under natural hypotheses motivated by Section 7, we will need to
specialize the formal variables $y_1$, $y_2$ and $y_0$ to complex
numbers $z_1$, $z_2$ and $z_0$, respectively, in (\ref{Y12}) and
(\ref{Y34}), when $|z_1|>|z_2|>0$ and $|z_2|>|z_0|>0$.  For this,
following \cite{tensor4}, we will need the next lemma, on products of
formal delta functions, with certain of the variables being complex
variables in suitable domains.  Our formulation and proof here are
different {}from those in \cite{tensor4}.  In addition to justifying the
specializations just indicated, this lemma will give us the natural
relation between the specialized expressions (\ref{F12}) and
(\ref{F34}) below.

\begin{lemma}\label{deltalemma}
Let $z_1$ and $z_2$ be complex numbers and set $z_0=z_1-z_2$.  Then
the left-hand sides of the following expressions converge absolutely
in the indicated domains, in the sense that the coefficient of each
monomial in the formal variables $x_0$, $x_1$ and $x_2$ is an
absolutely convergent series in the two variables related by the
inequalities, and the following identities hold:
\begin{eqnarray}
\dlt{x_1}{x_0}{-z_1}\dlt{x_2}{x_0}{-z_2}&=&
\dlt{x_2}{x_0}{-z_2}\dlt{x_1}{x_2}{-z_0}\nno\\
&&\mbox{\rm for arbitrary } z_1, z_2;\label{l1}\\
\dlt{z_1}{x_0}{-x_1}\dlt{x_2}{x_0}{-z_2} & = &
\dlt{x_0}{z_1}{+x_1}\dlt{x_2}{z_0}{+x_1}\nno \\
&& \mbox{\rm if }|z_1|>|z_2|;\label{l2-1}\\
\dlt{z_2}{x_0}{-x_2}\dlt{z_0}{x_2}{-x_1}&=&
\dlt{x_0}{z_1}{+x_1}\dlt{x_2}{z_0}{+x_1}\nno \\
&& \mbox{\rm if }|z_2|>|z_0|>0;\label{l2-2}\\
\dlt{x_1}{-z_1}{+x_0}\dlt{z_2}{x_0}{-x_2}&=&
\dlt{x_0}{z_2}{+x_2}\dlt{x_1}{-z_0}{+x_2}\nno \\
&& \mbox{\rm if }|z_1|>|z_2|>0;\label{l3}\\
\dlt{x_2}{-z_2}{+x_0}\dlt{x_1}{x_2}{-z_0}&=&
\dlt{x_1}{-z_1}{+x_0}\dlt{x_2}{-z_2}{+x_0}\nno \\
&& \mbox{\rm if }|z_2|>|z_0|.\label{l4}
\end{eqnarray}
(Note that the first identity does not require a restricted domain for
$z_1$, $z_2$ and $z_0$, while the others need certain conditions among
the complex numbers $z_i$ in order for the expressions on the
left-hand sides to be well defined, that is, absolutely convergent.
None of the five expressions on the right-hand sides require
restricted domains for absolute convergence.)
\end{lemma}

\pf In this proof we will use additional formal variables $y_0$,
$y_1$, $y_2$, and repeatedly use Remark 2.3.25 in \cite{LL} about
delta function substitution.

First, we have
\begin{eqnarray*}
\dlt{x_1}{x_0}{-y_1}\dlt{x_2}{x_0}{-y_2}&=&
\dlt{x_1}{x_0}{-y_1}\dlt{x_0}{x_2}{+y_2}\\
&=& \dlt{x_1}{x_2+y_2}{-y_1}\dlt{x_0}{x_2}{+y_2}\\
&=& \dlt{x_2}{x_0}{-y_2}\dlt{x_1}{x_2}{-(y_1-y_2)}.
\end{eqnarray*}
(Note that the notation $(x_0+y_2-y_1)^{n}$ is unambiguous: it is the
power series expansion in nonnegative powers of $y_1$ and $y_2$.)
Since it is clear that the left-hand side of this identity lies in
\[
{\mathbb C}[y_1, y_2]((x_0^{-1}))[[x_1, x_1^{-1}, x_2, x_2^{-1}]],
\]
one
can substitute any complex numbers $z_1$, $z_2$ for $y_1$, $y_2$,
respectively, and get the identity (\ref{l1}).

For (\ref{l2-1}), we have
\begin{eqnarray*}
\dlt{y_1}{x_0}{-x_1}\dlt{x_2}{x_0}{-y_2}&=&
\dlt{x_0}{y_1}{+x_1}\dlt{x_2}{x_0}{-y_2}\\
& = & \dlt{x_0}{y_1}{+x_1}\dlt{x_2}{y_1+x_1}{-y_2}\\
& = & \dlt{x_0}{y_1}{+x_1}\dlt{x_2}{(y_1-y_2)}{+x_1},
\end{eqnarray*}
and the right-hand side and hence the left-hand side lies in 
\[{\mathbb C}[y_1, y_1^{-1}, (y_1-y_2),
(y_1-y_2)^{-1}][[x_0, x_0^{-1}, x_1,x_2, x_2^{-1}]].\]
Thus if
$|z_1|>|z_2|>0$, so that the binomial expansion of $(z_1-z_2)^{n}$
converges for all $n$, we can substitute $z_1$, $z_2$ for $y_1$,
$y_2$ and obtain (\ref{l2-1}). On the other hand,
\begin{eqnarray*}
\dlt{y_2}{x_0}{-x_2}\dlt{y_0}{x_2}{-x_1}&=&
\dlt{x_0}{y_2}{+x_2}\dlt{x_2}{y_0}{+x_1}\\
&=&\dlt{x_0}{y_2}{+y_0+x_1}\dlt{x_2}{y_0}{+x_1}.
\end{eqnarray*}
It is clear {}from the right-hand side that both sides lie in 
\[{\mathbb
C}[y_0, y_0^{-1}, (y_2+y_0), (y_2+y_0)^{-1}][[x_0, x_0^{-1}, x_1, x_2,
x_2^{-1}]],\]
so if $|z_2|>|z_0|>0$ we can substitute $z_2$, $z_0$
for $y_2$, $y_0$ and obtain (\ref{l2-2}).

To prove (\ref{l3}), we see that
\begin{eqnarray*}
\dlt{x_1}{-y_1}{+x_0}\dlt{y_2}{x_0}{-x_2}&=&
\dlt{x_1}{-y_1}{+x_0}\dlt{x_0}{y_2}{+x_2}\\ &=&
\dlt{x_1}{-y_1}{+y_2+x_2}\dlt{x_0}{y_2}{+x_2},
\end{eqnarray*}
and the right-hand side and hence both sides lie in 
\[{\mathbb C}[y_2, y_2^{-1}, (y_1-y_2), (y_1-y_2)^{-1}][[x_0,
x_0^{-1}, x_1, x_1^{-1}, x_2]],\]
so that when $|z_1|>|z_2|>0$ we can
substitute $z_1$, $z_2$ for $y_1$, $y_2$ and obtain (\ref{l3}).
Finally, we have
\begin{eqnarray*}
\dlt{x_2}{-y_2}{+x_0}\dlt{x_1}{x_2}{-y_0}&=&
\dlt{x_2}{-y_2}{+x_0}\dlt{x_1}{-y_2+x_0}{-y_0}\\
&=&\dlt{x_2}{-y_2}{+x_0}\dlt{x_1}{-(y_2+y_0)}{+x_0},
\end{eqnarray*}
and {}from the right-hand side we see that both sides lie in 
\[{\mathbb
C}[(y_2+y_0), (y_2+y_0)^{-1}, y_2, y_2^{-1}][[x_0, x_1, x_1^{-1}, x_2,
x_2^{-1}]],\]
so that when $|z_2|>|z_0|$ we can substitute $z_2$,
$z_0$ for $y_2$, $y_0$ and obtain the identity (\ref{l4}).  \epfv

If we assume that the convergence condition for intertwining maps in
${\cal C}$ holds and that our generalized modules are objects of
${\cal C}$, in the setting of Section 7, then after pairing with an
element $w'_{(4)} \in W'_4$, we can specialize the formal variables
$y_1$, $y_2$ to complex numbers $z_1$, $z_2$ in (\ref{Y12}) whenever
$|z_1|>|z_2|>0$, and we can specialize $y_2$, $y_0$ to complex numbers
$z_2$, $z_0$ in (\ref{Y34}) whenever $|z_2|>|z_0|>0$, using Lemma
\ref{deltalemma}:

\begin{propo}\label{compositeJacobiforproductsanditerates}
Assume that the convergence condition for intertwining maps in ${\cal
C}$ holds and that the generalized modules entering into (\ref{Y12})
and (\ref{Y34}) are objects of ${\cal C}$.  Continuing to use the
notation of (\ref{Y12}) and (\ref{Y34}), also let $w'_{(4)} \in W'_4$.
Let $z_1$, $z_2$ be complex numbers satisfying $|z_1|>|z_2|>0$.  Then
for a $P(z_1)$-intertwining map $I_1$ of type ${W_4\choose W_1 M_1}$
and a $P(z_2)$-intertwining map $I_2$ of type ${M_1\choose W_2 W_3}$,
the following expressions are absolutely convergent, and the following
formula for the product $I_{1}\circ (1_{W_{1}}\otimes I_{2})$ of $I_1$
and $I_2$ holds:
\begin{eqnarray}\label{F12}
\lefteqn{\Bigg\langle w'_{(4)},
\dlt{x_1}{x_0}{-z_1}\dlt{x_2}{x_0}{-z_2}Y_4(v,x_0)
(I_{1}\circ (1_{W_{1}}\otimes I_{2}))
(w_{(1)},w_{(2)},w_{(3)})\Bigg\rangle}\nno \\
&&=\Bigg\langle w'_{(4)},\dlt{z_1}{x_0}{-x_1}\dlt{x_2}{x_0}{-z_2}\cdot\nno\\
&&\hspace{5cm}\cdot(I_{1}\circ (1_{W_{1}}\otimes I_{2}))
(Y_1(v,x_1)w_{(1)},w_{(2)},
w_{(3)})\Bigg\rangle\nno\\
&&\quad+\Bigg\langle w'_{(4)},\dlt{x_1}{-z_1}{+x_0}\dlt{z_2}{x_0}{-x_2}\cdot\nno\\
&&\hspace{5cm}\cdot(I_{1}\circ (1_{W_{1}}\otimes I_{2}))
(w_{(1)},Y_2(v,x_2)w_{(2)},
w_{(3)})\Bigg\rangle\nno\\
&&\quad +\Bigg\langle w'_{(4)},\dlt{x_1}{-z_1}{+x_0}\dlt{x_2}{-z_2}{+x_0}\cdot\nno\\
&&\hspace{5cm}\cdot(I_{1}\circ (1_{W_{1}}\otimes I_{2}))
(w_{(1)},w_{(2)},Y_3(v,x_0)
w_{(3)})\Bigg\rangle . \nno\\
\end{eqnarray}
Moreover, let $z_2$, $z_0$ be complex numbers satisfying
$|z_2|>|z_0|>0$. Then for a $P(z_2)$-intertwining map $I^1$ of type
${W_4\choose M_2 W_3}$ and a $P(z_0)$-intertwining map $I^2$ of type
${M_2\choose W_1 W_2}$, the following expressions are absolutely
convergent, and the following formula for the iterate $I^{1}\circ
(I^{2}\otimes 1_{W_{3}})$ of $I^1$ and $I^2$ holds:
\begin{eqnarray}\label{F34}
\lefteqn{\Bigg\langle w'_{(4)},
\dlt{x_2}{x_0}{-z_2}\dlt{x_1}{x_2}{-z_0}Y_4(v,x_0)
(I^{1}\circ (I^{2}\otimes 1_{W_{3}}))
(w_{(1)},w_{(2)},w_{(3)})\Bigg\rangle}\nno\\
&&=\Bigg\langle w'_{(4)},\dlt{z_2}{x_0}{-x_2}\dlt{z_0}{x_2}{-x_1}\cdot\nno\\
&&\hspace{5cm}\cdot(I^{1}\circ (I^{2}\otimes 1_{W_{3}}))
(Y_1(v,x_1)w_{(1)},w_{(2)},
w_{(3)})\Bigg\rangle\nno\\
&&\quad +\Bigg\langle w'_{(4)},\dlt{z_2}{x_0}{-x_2}\dlt{x_1}{-z_0}{+x_2}\cdot\nno\\
&&\hspace{5cm}\cdot(I^{1}\circ (I^{2}\otimes 1_{W_{3}}))
(w_{(1)},Y_2(v,x_2)w_{(2)},
w_{(3)})\Bigg\rangle\nno \\
&&\quad+\Bigg\langle w'_{(4)},\dlt{x_2}{-z_2}{+x_0}\dlt{x_1}{x_2}{-z_0}\cdot\nno\\
&&\hspace{5cm}\cdot(I^{1}\circ (I^{2}\otimes 1_{W_{3}}))
(w_{(1)},w_{(2)},Y_3(v,x_0)
w_{(3)})\Bigg\rangle . \nno\\
\end{eqnarray}
\end{propo}

\pf When $y_1$ and $y_2$ are specialized to $z_1$ and $z_2$,
respectively, the product of the two delta-function expressions on the
left-hand side of (\ref{Y12}) and the three products of pairs of
delta-function expressions on the right-hand side of (\ref{Y12}) all
converge absolutely in the domain $|z_1|>|z_2|>0$, by Lemma
\ref{deltalemma}; note that for the last of the three products of
pairs of delta-function expressions on the right-hand side of
(\ref{Y12}), the convergence is immediate.  Analogously, {}from Lemma
\ref{deltalemma} we see that the corresponding statements also hold
for (\ref{Y34}), when $y_0$ and $y_2$ are specialized to $z_0$ and
$z_2$, respectively, in the domain $|z_1|>|z_2|>0$.  Recalling the
notations (\ref{4itm}), (\ref{4prm}), (\ref{iterateabbreviation}) and
(\ref{productabbreviation}), we see that the result follows {}from the
convergence condition.  \epfv

Considering the ${\mathfrak s}{\mathfrak l}(2)$-action instead of the
$V$-action, by (\ref{log:L(j)b}) we have
\begin{eqnarray}
\lefteqn{L(j)\mathcal{Y}_1(w_{(1)},y_1)\mathcal{Y}_2(w_{(2)},y_2)w_{(3)}}
\nno\\
&&=\sum_{i=0}^{j+1}{j+1\choose i}y_1^i\mathcal{Y}_1(L(j-i)w_{(1)},y_1)
\mathcal{Y}_2(w_{(2)},y_2)w_{(3)}\nno\\ &&\quad+\mathcal{Y}_1(w_{(1)},y_1)L(j)
\mathcal{Y}_2(w_{(2)},y_2)w_{(3)}\nno\\
&&=\sum_{i=0}^{j+1}{j+1\choose i}y_1^i\mathcal{Y}_1(L(j-i)w_{(1)},y_1)
\mathcal{Y}_2(w_{(2)},y_2)w_{(3)}\nno\\ &&\quad+\mathcal{Y}_1(w_{(1)},y_1)
\sum_{k=0}^{j+1}{j+1\choose k}y_2^k{\cal
Y}_2(L(j-k)w_{(2)},y_2)w_{(3)}\nno\\ &&\quad+{\cal
Y}_1(w_{(1)},y_1)\mathcal{Y}_2(w_{(2)},y_2)L(j)w_{(3)}
\end{eqnarray}
for $j=-1, 0$ and $1$.  In the setting of Proposition
\ref{compositeJacobiforproductsanditerates}, if $|z_1|>|z_2|>0$ we can
substitute $z_1$, $z_2$ for $y_1$, $y_2$, respectively, and we obtain,
setting $z_0 = z_1 - z_2$,
\begin{eqnarray}\label{zz:sl2p}
\lefteqn{\langle w'_{(4)},L(j)(I_{1}\circ (1_{W_{1}}\otimes I_{2}))
(w_{(1)},w_{(2)},w_{(3)})\rangle}\nno\\
&&=\Bigg\langle w'_{(4)},\sum_{i=0}^{j+1}{j+1\choose i}(z_2+z_0)^i
(I_{1}\circ (1_{W_{1}}\otimes I_{2}))
(L(j-i)w_{(1)},w_{(2)},w_{(3)})\Bigg\rangle\nno\\
&&\quad+\Bigg\langle w'_{(4)},\sum_{k=0}^{j+1}{j+1\choose k}z_2^k
(I_{1}\circ (1_{W_{1}}\otimes I_{2}))
(w_{(1)},L(j-k)w_{(2)},w_{(3)})\Bigg\rangle\nno\\
&&\quad+\langle w'_{(4)},
(I_{1}\circ (1_{W_{1}}\otimes I_{2}))
(w_{(1)},w_{(2)},L(j)w_{(3)})\rangle
\end{eqnarray}
for $j=-1, 0$ and $1$.

On the other hand, by (\ref{log:L(j)b}) we also have
\begin{eqnarray}\label{zz:sl2i0}
\lefteqn{L(j){\cal Y}^1({\cal Y}^2(w_{(1)},y_0)w_{(2)},y_2)w_{(3)}}
\nno\\
&&=\sum_{i=0}^{j+1}{j+1\choose i}y_2^i{\cal Y}^1(L(j-i){\cal Y}^2
(w_{(1)}, y_0)w_{(2)},y_2)w_{(3)}\nno\\ &&\quad+{\cal Y}^1({\cal
Y}^2(w_{(1)},y_0)w_{(2)},y_2)L(j)w_{(3)} \nno\\
&&=\sum_{i=0}^{j+1}{j+1\choose i}y_2^i{\cal Y}^1
\Bigg(\sum_{k=0}^{j-i+1}{j-i+1\choose k}y_0^k{\cal
Y}^2(L(j-i-k)w_{(1)},y_0)w_{(2)} ,y_2\Bigg)w_{(3)}\nno\\
&&\quad+\sum_{i=0}^{j+1}{j+1\choose i}y_2^i{\cal Y}^1( {\cal
Y}^2(w_{(1)}, y_0)L(j-i)w_{(2)} ,y_2)w_{(3)}\nno\\ &&\quad+{\cal
Y}^1({\cal Y}^2(w_{(1)},y_0)w_{(2)},y_2)L(j)w_{(3)}
\end{eqnarray}
for $j=-1, 0$ and $1$, The first term of the right-hand side is
\begin{eqnarray*}
\lefteqn{\sum_{i=0}^{j+1}{j+1\choose i}y_2^i\sum_{k=0}^{j-i+1}
{j-i+1\choose k}y_0^k{\cal Y}^1({\cal Y}^2(L(j-i-k)w_{(1)},y_0)w_{(2)}
,y_2)w_{(3)}}\nno\\
&&=\sum_{t=0}^{j+1}\sum_{k=0}^{j+1-t+k}
{j+1\choose t-k}{j+1-t+k\choose k}
y_2^{t-k}y_0^k{\cal Y}^1({\cal Y}^2(L(j-t)w_{(1)},y_0)w_{(2)}
,y_2)w_{(3)}\nno\\
&&=\sum_{t=0}^{j+1}\sum_{k=0}^t{j+1\choose t-k}{j+1-t+k\choose k}
y_2^{t-k}y_0^k{\cal Y}^1({\cal Y}^2(L(j-t)w_{(1)},y_0)w_{(2)}
,y_2)w_{(3)}\nno\\
&&=\sum_{t=0}^{j+1}{j+1\choose t}(y_2+y_0)^t{\cal Y}^1({\cal
Y}^2(L(j-t)w_{(1)},y_0)w_{(2)} ,y_2)w_{(3)},
\end{eqnarray*}
where we have used the identity $\displaystyle{j+1\choose t-k}{j+1-t+k
\choose k}={j+1\choose t}{t\choose k}$ in the last step.  Thus in the
setting of Proposition \ref{compositeJacobiforproductsanditerates}, if
$|z_2|>|z_0|>0$ we can substitute $z_2$, $z_0$ for $y_2$, $y_0$,
respectively, in (\ref{zz:sl2i0}), and we obtain
\begin{eqnarray}\label{zz:sl2i}
\lefteqn{\langle w'_{(4)},L(j)(I^{1}\circ (I^{2}\otimes 1_{W_3}))
(w_{(1)},w_{(2)},w_{(3)})\rangle}\nno\\
&&=\Bigg\langle w'_{(4)},\sum_{t=0}^{j+1}{j+1\choose t}(z_2+z_0)^t
(I^{1}\circ (I^{2}\otimes 1_{W_3}))
(L(j-t)w_{(1)},w_{(2)},w_{(3)})\Bigg\rangle\nno\\
&&\quad+\Bigg\langle w'_{(4)},\sum_{i=0}^{j+1}{j+1\choose i}z_2^i
(I^{1}\circ (I^{2}\otimes 1_{W_3}))
(w_{(1)},L(j-i)w_{(2)},w_{(3)})\Bigg\rangle\nno\\
&&\quad+\langle w'_{(4)},(I^{1}\circ (I^{2}\otimes 1_{W_3}))
(w_{(1)},w_{(2)},L(j)w_{(3)})\rangle
\end{eqnarray}
for $j=-1, 0$ and $1$.

Of course, in case $V$ is a conformal vertex algebra, these formulas
follow {}from the earlier computation for the $V$-action (Proposition
\ref{compositeJacobiforproductsanditerates}), by setting $v=\omega$
and taking $\res_{x_1}\res_{x_2}\res_{x_0}x_0^{j+1}$, $j=-1, 0, 1$.

Lemma \ref{deltalemma}, Proposition
\ref{compositeJacobiforproductsanditerates}, (\ref{zz:sl2p}),
(\ref{zz:sl2i}) and Remark \ref{grad-comp-prod-iter} motivate the
following definition, which is analogous to the definition of the
notion of $P(z)$-intertwining map (Definition \ref{im:imdef}):

\begin{defi}
{\rm Let $z_0, z_1, z_2\in {\mathbb C}^{\times }$ with $z_0=z_1-z_2$
(so that in particular $z_1\neq z_2$, $z_0\neq z_1$ and $z_0\neq
-z_2$). Let $(W_1,Y_1)$, $(W_2,Y_2)$, $(W_3,Y_3)$ and $(W_4,Y_4)$ be
generalized modules for a M\"obius (or conformal) vertex algebra $V$.
A {\it $P(z_1,z_2)$-intertwining map} is a linear map $F:\, W_1\otimes
W_2\otimes W_3\to \overline{W}_4$ such that the following conditions
are satisfied: the {\it grading compatibility condition}: For $\beta,
\gamma, \delta \in \tilde A$ and $w_{(1)}\in W_1^{(\beta)}$,
$w_{(2)}\in W_2^{(\gamma)}$, $w_{(3)}\in W_3^{(\delta)}$,
\begin{equation} \label{grad-comp-F}
F(w_{(1)}\otimes 
w_{(2)}\otimes w_{(3)})\in \overline{W_{4}^{(\beta+\gamma+\delta)}};
\end{equation}
the {\em lower truncation condition}: for any elements $w_{(1)}\in
W_1$, $w_{(2)}\in W_2$ and $w_{(3)}\in W_3$, and any $n\in {\mathbb
C}$,
\begin{equation}\label{zz:ltc}
\pi_{n-m}F(w_{(1)}\otimes w_{(2)}\otimes w_{(3)})=0\;\;\mbox{ for }\;m\in
{\mathbb N}\;\mbox{ sufficiently large}
\end{equation}
(which follows {}from (\ref{grad-comp-F}), in view of the grading
restriction condition (\ref{set:dmltc})); the {\em composite Jacobi
identity}:
\begin{eqnarray}\label{zz:Y}
\lefteqn{\dlt{x_1}{x_0}{-z_1}\dlt{x_2}{x_0}{-z_2}Y_4(v,x_0)F(w_{(1)}\otimes
w_{(2)}\otimes w_{(3)})}\nno \\
&& =\dlt{x_0}{z_1}{+x_1}\dlt{x_2}{z_0}{+x_1}F(Y_1(v,x_1)w_{(1)}\otimes
w_{(2)}\otimes w_{(3)})\nno \\
&& \quad+\dlt{x_0}{z_2}{+x_2}\dlt{x_1}{-z_0}{+x_2}F(w_{(1)}\otimes
Y_2(v,x_2)w_{(2)}\otimes w_{(3)})\nno \\
&& \quad+\dlt{x_1}{-z_1}{+x_0}\dlt{x_2}{-z_2}{+x_0}F(w_{(1)}\otimes
w_{(2)}\otimes Y_3(v,x_0)w_{(3)})\nno \\
\end{eqnarray}
for $v\in V$, $w_{(1)}\in W_1$, $w_{(2)}\in W_2$ and $w_{(3)}\in W_3$
(note that all the expressions in the right-hand side of (\ref{zz:Y})
are well defined, that none of the products of delta-function
expressions require restricted domains, and that the left-hand side of
(\ref{zz:Y}) is meaningful because any infinite linear combination of
$v_n$ ($n \in {\mathbb Z}$) of the form $\sum_{n<N}a_nv_n$ ($a_n\in
{\mathbb C}$) acts in a well-defined way on any $F(w_{(1)}\otimes
w_{(2)}\otimes w_{(3)})$, in view of (\ref{zz:ltc})); and the {\em
${\mathfrak s}{\mathfrak l}(2)$-bracket relations}: for any
$w_{(1)}\in W_1$, $w_{(2)}\in W_2$ and $w_{(3)}\in W_3$,
\begin{eqnarray}\label{zz:L}
\lefteqn{L(j)F(w_{(1)}\otimes w_{(2)}\otimes w_{(3)})}\nno\\
&&=\sum_{i=0}^{j+1}{j+1\choose i}z_{1}^{i}F
(L(j-i)w_{(1)}\otimes w_{(2)}\otimes w_{(3)})\nno\\
&&\quad+\sum_{k=0}^{j+1}{j+1\choose k}z_2^kF
(w_{(1)}\otimes L(j-k)w_{(2)}\otimes w_{(3)})\nno\\
&&\quad+F(w_{(1)}\otimes w_{(2)}\otimes L(j)w_{(3)})
\end{eqnarray}
for $j=-1, 0$ and $1$ (again, in case $V$ is a conformal
vertex algebra, this follows {}from (\ref{zz:Y}) by setting $v=\omega$
and taking $\res_{x_1}\res_{x_2}\res_{x_0}x_0^{j+1}$).  }
\end{defi}

We emphasize that every term in (\ref{zz:Y}) and (\ref{zz:L}) in this
definition is purely algebraic; that is, no convergence is involved.

{}From Lemma \ref{deltalemma}, Proposition
\ref{compositeJacobiforproductsanditerates}, (\ref{zz:sl2p}),
(\ref{zz:sl2i}) and Remark \ref{grad-comp-prod-iter}, we have the
following:

\begin{propo}\label{productanditerateareintwmaps}
In the setting of Proposition
\ref{compositeJacobiforproductsanditerates}, for intertwining maps
$I_1$, $I_2$, $I^1$ and $I^2$ as indicated, when $|z_1|>|z_2|>0$,
$I_1\circ (1_{W_1}\otimes I_2)$ is a $P(z_1,z_2)$-intertwining map and
when $|z_2|>|z_0|>0$, $I^1\circ (I^2\otimes 1_{W_3})$ is a
$P(z_2+z_0,z_2)$-intertwining map. \epf
\end{propo}

Now we consider $P(z_1,z_2)$-intertwining maps {}from a ``dual''
viewpoint, and we use this to motivate an analogue $\tau_{P(z_{1},
z_{2})}$ of the action $\tau_{P(z)}$ introduced in Subsection 5.2.
Fix any $w'_{(4)}\in W'_4$. Then (\ref{zz:Y}) implies:
\begin{eqnarray}
\lefteqn{\Bigg\langle
w'_{(4)},\dlt{x_1}{x_0}{-z_1}\dlt{x_2}{x_0}{-z_2}Y_4(v,x_0)F(w_{(1)}
\otimes w_{(2)}\otimes w_{(3)})\Bigg\rangle} \nno \\
&&=\Bigg\langle
w'_{(4)},\dlt{x_0}{z_1}{+x_1}\dlt{x_2}{z_0}{+x_1}F(Y_1(v,x_1)w_{(1)}
\otimes w_{(2)}\otimes w_{(3)})\Bigg\rangle \nno \\
&& \quad +\Bigg\langle
w'_{(4)},\dlt{x_0}{z_2}{+x_2}\dlt{x_1}{-z_0}{+x_2}F(w_{(1)}\otimes
Y_2(v,x_2)w_{(2)}\otimes w_{(3)})\Bigg\rangle \nno \\
&& \quad +\Bigg\langle
w'_{(4)},\dlt{x_1}{-z_1}{+x_0}\dlt{x_2}{-z_2}{+x_0}F(w_{(1)}\otimes
w_{(2)}\otimes Y_3(v,x_0)w_{(3)})\Bigg\rangle .\nno\\
&& \label{cmpF}
\end{eqnarray}
The left-hand side can be written as
\begin{eqnarray*}
\Bigg\langle
\dlt{x_1}{x_0}{-z_1}\dlt{x_2}{x_0}{-z_2}Y'_4(e^{x_0L(1)}(-x_0^2)^{-L(0)}v,
x_0^{-1})w'_{(4)},
F(w_{(1)}\otimes w_{(2)}\otimes w_{(3)})\Bigg\rangle,
\end{eqnarray*}
and so by replacing $v$ by $(-x_0^2)^{L(0)}e^{-x_0L(1)}v$ and then
replacing $x_0$ by $x_0^{-1}$ in both sides of (\ref{cmpF}) we
see that
\begin{eqnarray}\label{taumot}
\lefteqn{\Bigg\langle \dlt{x_1}{x^{-1}_0}{-z_1}\dlt{x_2}{x^{-1}_0}{-z_2}
Y'_4(v,x_0)w'_{(4)},F(w_{(1)}\otimes w_{(2)}\otimes
w_{(3)})\Bigg\rangle} \nno\\ 
&&=\Bigg\langle
w'_{(4)},\dlti{x_0}{z_1}{+x_1}\dlt{x_2}{z_0}{+x_1}\cdot\nno\\
&&\qquad\qquad
F(Y_1((-x_0^{-2})^{L(0)}e^{-x_0^{-1}L(1)}v,x_1)w_{(1)}\otimes
w_{(2)}\otimes w_{(3)})\Bigg\rangle \nno\\ 
&&+\Bigg\langle
w'_{(4)},\dlti{x_0}{z_2}{+x_2}\dlt{x_1}{-z_0}{+x_2}\cdot\nno\\
&&\qquad\qquad F(w_{(1)}\otimes
Y_2((-x_0^{-2})^{L(0)}e^{-x_0^{-1}L(1)}v,x_2)w_{(2)}\otimes
w_{(3)})\Bigg\rangle \nno\\ 
&&+\Bigg\langle
w'_{(4)},\dlt{x_1}{-z_1}{+x^{-1}_0}\dlt{x_2}{-z_2}{+x^{-1}_0}\cdot\nno\\
&&\qquad\qquad F(w_{(1)}\otimes w_{(2)}\otimes
Y_3((-x_0^{-2})^{L(0)}e^{-x_0^{-1}L(1)}v,x_0^{-1})w_{(3)})\Bigg\rangle.
\end{eqnarray}

Arguing just as in (\ref{y-t-delta})--(\ref{3.19-1}),
we note that in the left-hand side of (\ref{taumot}), the coefficients of
\[
\dlt{x_1}{x^{-1}_0}{-z_1}\dlt{x_2}{x^{-1}_0}{-z_2}Y'_4(v,x_0)
\]
in powers of $x_0$, $x_1$ and $x_2$, for all $v\in V$, span
\[
\tau_{W'_4}(V\otimes \iota_{+}{\mathbb C}[t,t^{-1},(z_1^{-1}-t)^{-1},
(z_2^{-1}-t)^{-1}])
\]
(recall the notation $\tau_W$ {}from (\ref{tauW}), (\ref{tauw}),
(\ref{3.7}) and the notation $\iota_{\pm}$ {}from (\ref{iota+-})).
By analogy with the case of $P(z)$-intertwining maps, we shall
define an action of $V\otimes \iota_{+}{\mathbb C}[t,t^{-1},(z_1^{-1}
-t)^{-1},(z_2^{-1}-t)^{-1}]$ on $(W_1\otimes W_2\otimes W_3)^{*}$. We
shall need the following analogue of Lemma \ref{tauP}, where we use the 
notations $Y_{t}$, $T_{z}$ and $o$ introduced in Subsection 5.1,
and where we recall that $z_{0}=z_{1}-z_{2}$:

\begin{lemma}\label{tauzzlm}
We have
\begin{eqnarray}
\lefteqn{o\bigg(\dlt{x_1}{x_0^{-1}}{-z_1}\dlt{x_2}
{x_0^{-1}}{-z_2}Y_{t}(v,x_0)\bigg)}\nno\\
&&=\dlt{x_1}{x_0^{-1}}{-z_1}\dlt{x_2}
{x_0^{-1}}{-z_2}Y^o_{t}(v,x_0),\label{zztr1}\\
\lefteqn{(\iota_+\circ\iota_-^{-1}\circ o)\bigg(\dlt{x_1}{x_0^{-1}}{-z_1}
\dlt{x_2}{x_0^{-1}}{-z_2}Y_{t}(v,x_0)\bigg)}\nno\\
&&=\dlt{x_1}{-z_1}{+x^{-1}_0}\dlt{x_2}{-z_2}{+x^{-1}_0}Y^o_t(v,x_0),
\label{zztr2}\\
\lefteqn{(\iota_+\circ T_{z_1}\circ\iota_-^{-1}\circ o)\bigg(\dlt{x_1}
{x_0^{-1}}{-z_1}\dlt{x_2}{x_0^{-1}}{-z_2}Y_{t}(v,x_0)\bigg)}\nno\\
&&=\dlti{x_0}{z_1}{+x_1}\dlt{x_2}{z_0}{+x_1}Y_t((-x_0^{-2})^{L(0)}
e^{-x_0^{-1}L(1)}v,x_1),\label{zztr3}\\
\lefteqn{(\iota_+\circ T_{z_2}\circ\iota_-^{-1}\circ o)\bigg(\dlt{x_1}
{x_0^{-1}}{-z_1}\dlt{x_2}{x_0^{-1}}{-z_2}Y_{t}(v,x_0)\bigg)}\nno\\
&&=\dlti{x_0}{z_2}{+x_2}\dlt{x_1}{-z_0}{+x_2}Y_t((-x_0^{-2})^{L(0)}
e^{-x_0^{-1}L(1)}v,x_2).\label{zztr4}
\end{eqnarray}
\end{lemma}
\pf The identity (\ref{zztr1})  immediately follows {}from (\ref{3.38}), and
(\ref{zztr2}) follows {}from (\ref{zztr1}), as in the proof of
(\ref{ztr2}). For (\ref{zztr3}), note that by (\ref{op-y-t-2}),
the coefficient of $x_1^{-m-1}x_2^{-n-1}$ in the right-hand side of
(\ref{zztr1}) is
\begin{eqnarray*}
\lefteqn{(x_0^{-1}-z_1)^m(x_0^{-1}-z_2)^n\bigg(e^{x_0L(1)}(-x_0^{-2})^{L(0)}v\otimes
x_0\delta\bigg(\frac t{x_0^{-1}}\bigg)\bigg)}\nn
&&=(t-z_1)^m(t-z_2)^n\bigg(e^{x_0L(1)}(-x_0^{-2})^{L(0)}v\otimes
x_0\delta\bigg(\frac t{x_0^{-1}}\bigg)\bigg).
\end{eqnarray*}
Acted on by $\iota_+\circ T_{z_1}\circ\iota_-^{-1}$, this becomes
\begin{eqnarray*}
\lefteqn{t^m(z_0+t)^n\bigg(e^{x_0L(1)}(-x_0^{-2})^{L(0)}v\otimes
x_0\delta\bigg(\frac {z_1+t}{x_0^{-1}}\bigg)\bigg)}\nn
&&=x_0\delta\bigg(\frac {z_1+t}{x_0^{-1}}\bigg)(z_0+t)^n\bigg(
e^{x_0L(1)}(-x_0^{-2})^{L(0)}v\otimes t^m\bigg)\nn
&&=x_0\delta\bigg(\frac {z_1+t}{x_0^{-1}}\bigg)(z_0+t)^n\bigg(
(-x_0^{-2})^{L(0)}e^{-x_0^{-1}L(1)}v\otimes t^m\bigg),
\end{eqnarray*}
by formula (5.3.1) in \cite{FHL}, and using (\ref{3.5}), we see that 
this is the coefficient of $x_1^{-m-1}x_2^{-n-1}$ in the right-hand side
of (\ref{zztr3}). The analogous
identity (\ref{zztr4}) is proved similarly. \epfv

Our analogue of Definition \ref{deftau} is:

\begin{defi}\label{tauzzdef}{\rm
Let $z_1, z_2\in {\mathbb C}^{\times}$, $z_1\neq z_2$. We define
a linear action
$\tau_{P(z_1,z_2)}$ of the space
\begin{equation}\label{thez1z2space}
V\otimes \iota _{+}{\mathbb C}[t,t^{-1},(z_1^{-1}-t)^{-1},(z_2^{-1}
-t)^{-1}]
\end{equation}
on $(W_1\otimes W_2\otimes W_3)^{*}$  by
\begin{eqnarray}\label{tauzzdef0}
\lefteqn{(\tau_{P(z_1,z_2)}(\xi)\lambda)(w_{(1)}\otimes w_{(2)}\otimes
w_{(3)})}\nno\\
&&=\lambda(\tau_{W_1}((\iota_+\circ T_{z_1}\circ\iota_-^{-1}\circ
o)\xi)w_{(1)}\otimes w_{(2)}\otimes w_{(3)})\nno\\
&&\quad +\lambda(w_{(1)}\otimes\tau_{W_2}((\iota_+\circ T_{z_2}\circ
\iota_-^{-1}\circ o)\xi)w_{(2)}\otimes w_{(3)})\nno\\
&&\quad +\lambda(w_{(1)}\otimes w_{(2)}\otimes\tau_{W_3}((\iota_+\circ
\iota_-^{-1}\circ o)\xi)w_{(3)})
\end{eqnarray}
for $\xi\in V\otimes \iota _{+}{\mathbb C}[t,t^{-1},(z_1^{-1}-t)^{-1},
(z_2^{-1}-t)^{-1}]$, $\lambda \in (W_1\otimes W_2\otimes W_3)^{*}$,
$w_{(1)}\in W_1$, $w_{(2)}\in W_2$ and $w_{(3)}\in W_3$. (The fact that
the right-hand side is in fact well defined follows immediately 
{}from the generating function reformulation of (\ref{tauzzdef0}) given
in (\ref{tauzzgf}) below.)
Denote by $Y'_{P(z_1,z_2)}$ the action of $V\otimes{\mathbb
C}[t,t^{-1}]$ on $(W_1\otimes W_2\otimes W_3)^*$ thus defined, that
is,
\begin{equation}\label{y'-zz}
Y'_{P(z_1,z_2)}(v,x)=\tau_{P(z_1,z_2)}(Y_t(v,x)).
\end{equation}
}
\end{defi}

By Lemma \ref{tauzzlm}, (\ref{3.7}) and (\ref{tauw-yto}), 
we see that (\ref{tauzzdef0}) can be written in terms of
generating functions as
\begin{eqnarray}\label{tauzzgf}
\lefteqn{\left(\tau_{P(z_1,z_2)}\left(\dlt{x_1}{x_0^{-1}}{-z_1}\dlt{x_2}
{x_0^{-1}}{-z_2}Y_{t}(v,x_0)\right)
\lambda \right)(w_{(1)}\otimes w_{(2)}\otimes
w_{(3)})}\nno\\
&=& \dlti{x_0}{z_1}{+x_1}\dlt{x_2}{z_0}{+x_1}\cdot\nno\\
&& \qquad\qquad \lambda (Y_1((-x_0^{-2})^{L(0)}e^{-x_0^{-1}L(1)}v,x_1)
w_{(1)}\otimes w_{(2)}\otimes w_{(3)})\nno\\
&& +\dlti{x_0}{z_2}{+x_2}\dlt{x_1}{-z_0}{+x_2}\cdot\nno\\
&& \qquad\qquad \lambda (w_{(1)}\otimes
Y_2((-x_0^{-2})^{L(0)}e^{-x_0^{-1}L(1)}v,x_2)w_{(2)}\otimes w_{(3)})\nno\\
&& +\dlt{x_1}{-z_1}{+x^{-1}_0}\dlt{x_2}{-z_2}{+x^{-1}_0}\lambda
(w_{(1)}\otimes w_{(2)}\otimes Y_3^o(v,x_0)w_{(3)})\nno\\
\end{eqnarray}
for $v\in V$, $\lambda \in (W_1\otimes W_2\otimes W_3)^{*}$,
$w_{(1)}\in W_1$, $w_{(2)}\in W_2$ and $w_{(3)}\in W_3$; the expansion 
coefficients in $x_{0}$, $x_{1}$ and $x_{2}$ of the left-hand side span 
the space of elements in the left-hand side of (\ref{tauzzdef0}). 
Compare this with the motivating formula (\ref{taumot}).  The
generating function form (\ref{y'-zz}) 
of the action $Y'_{P(z_1,z_2)}$ (\ref{y'-zz}) can be
obtained by taking $\res_{x_1}\res_{x_2}$ of both sides of
(\ref{tauzzgf}).

\begin{rema}{\rm
The action $\tau
_{P(z_1,z_2)}$ of $V\otimes \iota _{+}{\mathbb
C}[t,t^{-1},(z_1^{-1}-t)^{-1},(z_2^{-1}-t)^{-1}]$ on $(W_1 \otimes W_2
\otimes W_3)^*$,  defined for all $z_1,z_2\in {\mathbb C}^{\times
}$ with $z_1\neq z_2$,  coincides with the action $\tau
^{(1)}_{P(z_1,z_2)}$ when $|z_1|>|z_2|>0$, and coincides with the
action $\tau ^{(2)}_{P(z_1,z_2)}$ when $|z_2|>|z_1-z_2|>0$, where 
$\tau^{(1)}_{P(z_{1}, z_{2})}$ and $\tau^{(2)}_{P(z_{1}, z_{2})}$
are the two actions defined in Section 14 of \cite{tensor4}. The action 
$\tau_{P(z_{1}, z_{2})}$ and the related notion of 
$P(z_{1}, z_{2})$-intertwining map extend the 
corresponding considerations
in \cite{tensor4} in a natural way. }
\end{rema}

\begin{rema}\label{F-intw}
{\rm (cf. Remark \ref{I-intw})
Using the action $\tau_{P(z_1,z_2)}$, we can write
the equality (\ref{taumot}) as
\begin{eqnarray}\label{intw}
\lefteqn{\left(\dlt{x_1}{x^{-1}_0}{-z_1}\dlt{x_2}{x^{-1}_0}{-z_2}
Y'_4(v,x_0)w'_{(4)}\right)\circ F }\nno\\
&&=\tau_{P(z_1,z_2)}\left(\dlt{x_1}{x_0^{-1}}{-z_1}\dlt{x_2}
{x_0^{-1}}{-z_2}Y_{t}(v,x_0)\right)(w'_{(4)}\circ F).
\end{eqnarray}
Furthermore, using the action of 
$\iota_{+}{\mathbb C}[t,t^{-1}, (z_1^{-1}-t)^{-1},
(z_2^{-1}-t)^{-1}]$ on $W_{4}'$ (recall (\ref{tauW}), 
(\ref{tauw}) and (\ref{3.7})), we can also write
(\ref{intw}) as
\begin{eqnarray}\label{zz:Psi}
\lefteqn{\left(\tau_{W_{4}'}\left(\dlt{x_1}{x^{-1}_0}{-z_1}\dlt{x_2}{x^{-1}_0}
{-z_2} Y_{t}(v,x_0)\right) w'_{(4)}\right)\circ F}\nn
&&=\tau_{P(z_1,z_2)}\left(\dlt{x_1}{x_0^{-1}}{-z_1}\dlt{x_2}
{x_0^{-1}}{-z_2}Y_{t}(v,x_0)\right)(w'_{(4)}\circ F).
\end{eqnarray}}
\end{rema}

As in Section 5, we need to consider gradings by $A$ and $\tilde{A}$.

The space $W_1 \otimes W_2 \otimes W_3$ is naturally
$\tilde{A}$-graded, and this gives us naturally-defined subspaces
$((W_{1}\otimes W_{2}\otimes W_{3})^{*})^{(\beta)}$ for $\beta \in
\tilde{A}$, as in the discussion after Remark \ref{I-intw}.

The space (\ref{thez1z2space}) is naturally $A$-graded, {}from the
$A$-grading on $V$: For $\alpha\in A$,
\begin{equation}
(V\otimes \iota _{+}{\mathbb C}[t,t^{-1},(z_1^{-1}-t)^{-1},(z_2^{-1}
-t)^{-1}])^{(\alpha)}
=V^{(\alpha)}\otimes\iota _{+}{\mathbb C}
[t,t^{-1},(z_1^{-1}-t)^{-1},(z_2^{-1}
-t)^{-1}].
\end{equation}

\begin{defi}\label{3-mod-actioncompatible}
{\rm We call a linear action $\tau$ of 
\[
V\otimes \iota _{+}{\mathbb C}[t,t^{-1},(z_1^{-1}-t)^{-1},(z_2^{-1}
-t)^{-1}]
\]
on $(W_1 \otimes W_2 \otimes W_{3})^*$
{\it $\tilde{A}$-compatible} if 
for $\alpha\in A$, $\beta\in \tilde{A}$,
\[
\xi\in 
(V\otimes \iota _{+}{\mathbb C}[t,t^{-1},(z_1^{-1}-t)^{-1},(z_2^{-1}
-t)^{-1}])^{(\alpha)}
\]
and $\lambda\in ((W_1 \otimes W_2 \otimes W_{3})^{*})^{(\beta)}$,
\[
\tau(\xi)\lambda\in ((W_{1}\otimes W_{2} \otimes W_{3})^{*})^{(\alpha+\beta)}.
\]}
\end{defi}

{}From (\ref{tauzzdef0}) or (\ref{tauzzgf}), we have:

\begin{propo}\label{tauzz-a-comp}
The action  $\tau_{P(z_1,z_2)}$ is $\tilde{A}$-compatible. \epf
\end{propo}

Again as in Section 5, when $V$ is a conformal vertex algebra, we write
\[
Y'_{P(z_1,z_2)}(\omega ,x)=\sum _{n\in {\mathbb Z}}L'_{P(z_1,z_2)}(n)x^{-n-2}.
\]
In this case, by setting $v=\omega$ in (\ref{tauzzgf}) and taking
$\res_{x_0}{x_0}^{j+1}$ for $j=-1, 0, 1$, we see that
\begin{eqnarray}\label{LP'(j)F}
\lefteqn{(L'_{P(z_1,z_2)}(j)\lambda)(w_{(1)}\otimes w_{(2)}\otimes
w_{(3)})}\nno\\
&&=\lambda\Bigg(\Bigg(\sum_{i=0}^{1-j}{1-j\choose i}z_1^iL(-j-i)\Bigg)
w_{(1)}\otimes
w_{(2)}\otimes w_{(3)} \nno\\
&&\quad\quad\quad\quad +\sum_{i=0}^{1-j}{1-j\choose i}z_2^iw_{(1)}\otimes
L(-j-i)w_{(2)}\otimes w_{(3)}\nno\\
&&\quad\quad\quad\quad +w_{(1)}\otimes w_{(2)}\otimes L(-j)w_{(3)}\Bigg).
\end{eqnarray}
If $V$ is  a M\"obius vertex algebra, we
define the actions $L'_{P(z_1,z_2)}(j)$ on $(W_1\otimes W_2\otimes
W_3)^*$ by (\ref{LP'(j)F}) for $j=-1, 0$ and $1$.  Using these
notations, the ${\mathfrak s}{\mathfrak l}(2)$-bracket relations (\ref{zz:L})
for a $P(z_1, z_2)$-intertwining map $F$ can be written as
\begin{equation}\label{LwF=LwF}
(L'(j)w'_{(4)})\circ F= L'_{P(z_1,z_2)}(j)(w'_{(4)}\circ F)
\end{equation}
for $w'_{(4)}\in W'_4$, $j=-1, 0, 1$ (cf. Remarks \ref{I-intw2} and
\ref{F-intw}).  We have
\[
L'_{P(z_1,z_2)}(j)((W_{1}\otimes W_{2} \otimes W_{3})^{*})^{(\beta)}\subset 
((W_{1}\otimes W_{2}\otimes W_{3})^{*})^{(\beta)}
\]
for $j=-1, 0, 1$ and $\beta\in \tilde{A}$ (cf. Remark
\ref{L'jpreservesbetaspace} and Proposition \ref{tauzz-a-comp}).

For the natural analogue of Proposition \ref{pz} (see Proposition
\ref{zzcor} below), we shall use the following analogues of the
relevant notions in Sections 4 and 5:  A map
\[
F\in \hom(W_1\otimes W_2\otimes W_3, (W_{4}')^{*})
\]
is {\it $\tilde{A}$-compatible} if
\[
F\in \hom(W_1\otimes W_2\otimes W_3, \overline{W_{4}})
\]
and if $F$ satisfies the natural analogue of the condition in
(\ref{IAtildecompat}), as in (\ref{grad-comp-F}).  A map
\[
G \in \hom(W_4',(W_1\otimes W_2\otimes W_3)^{*})
\]
is {\it $\tilde{A}$-compatible} if $G$ satisfies the analogue of 
(\ref{JAtildecompat}).  Then just as in Lemma \ref{IlambdatoJlambda}
and Remark \ref{alternateformoflemma}:

\begin{rema}\label{Atildecompatcorrespondence}{\rm
We have a canonical isomorphism {}from the space of
$\tilde{A}$-compatible linear maps
\[
F:W_{1}\otimes W_{2}\otimes W_{3} \rightarrow \overline{W_4}
\]
to the space of $\tilde{A}$-compatible linear maps
\[
G:W'_{4} \rightarrow (W_{1}\otimes W_{2} \otimes W_{3})^{*},
\]
determined by:
\begin{equation}
\langle w'_{(4)}, F(w_{(1)}\otimes w_{(2)}\otimes w_{(3)})\rangle
=G(w'_{(4)})(w_{(1)}\otimes w_{(2)}\otimes w_{(3)})
\end{equation}
for $w_{(1)}\in W_{1}$, $w_{(2)}\in W_{2}$, $w_{(3)}\in W_{3}$  
and $w'_{(4)}\in W'_{4}$,
or equivalently,
\begin{equation}\label{wF=Gw}
w'_{(4)}\circ F = G(w'_{(4)})
\end{equation}
for $w'_{(4)} \in W'_{4}$.
}
\end{rema}

We also have the natural analogues of Definition
\ref{gradingrestrictedmapJ} and Remark
\ref{Jcompatimpliesgradingrestr}:

\begin{defi}{\rm
A map $G\in \hom(W'_4, (W_1\otimes W_2\otimes W_3)^{*})$ is {\em
grading restricted} if for $n\in {\mathbb C}$, $w_{(1)}\in W_1$,
$w_{(2)}\in W_2$ and $w_{(3)}\in W_3$,
\[
G((W'_4)_{[n-m]})(w_{(1)}\otimes w_{(2)}\otimes w_{(3)})=0\;\;\mbox{ for
}\;m\in {\mathbb N}\;\mbox{ sufficiently large.}
\]
}
\end{defi}

\begin{rema}
{\rm
If $G \in {\rm Hom}(W'_4,(W_1\otimes W_2 \otimes W_3)^{*})$ is
$\tilde{A}$-compatible, then $G$ is also grading restricted.
}
\end{rema}

As in Proposition \ref{pz} we now have:

\begin{propo}\label{zzcor}
Let $z_1,z_2\in {\mathbb C}^{\times }$, $z_1\neq z_2$. Let $W_1$,
$W_2$, $W_3$ and $W_4$ be generalized $V$-modules.  Then under the
canonical isomorphism described in Remark
\ref{Atildecompatcorrespondence}, the $P(z_1,z_2)$-intertwining maps $F$
correspond exactly to the (grading restricted) $\tilde{A}$-compatible
maps $G$ that intertwine the actions of
\[
V\otimes \iota _{+}{\mathbb C}[t,t^{-1},(z_1^{-1}-t)^{-1},(z_2^{-1}-t)^{-1}]
\]
and of $L'(j)$ and $L'_{P(z_1,z_2)}(j)$, $j=-1,0,1$, on $W_4'$ and on
$(W_1\otimes W_2\otimes W_3)^*$.
\end{propo}

\pf By (\ref{wF=Gw}), Remark \ref{F-intw} asserts that (\ref{taumot}),
or equivalently, (\ref{zz:Y}), is equivalent to the condition
\begin{eqnarray}\label{Gtau=tauG}
\lefteqn{G\left(\tau_{W_{4}'}\left(\dlt{x_1}{x^{-1}_0}{-z_1}\dlt{x_2}{x^{-1}_0}
{-z_2} Y_{t}(v,x_0)\right) w'_{(4)}\right)}\nn
&&=\tau_{P(z_1,z_2)}\left(\dlt{x_1}{x_0^{-1}}{-z_1}\dlt{x_2}
{x_0^{-1}}{-z_2}Y_{t}(v,x_0)\right)G(w'_{(4)}),
\end{eqnarray}
that is, the condition that $G$ intertwines the actions of 
\[
V\otimes \iota _{+}{\mathbb C}[t,t^{-1},(z_1^{-1}-t)^{-1},(z_2^{-1}-t)^{-1}]
\]
on $W_4'$ and on $(W_1\otimes W_2\otimes W_3)^*$.  Analogously, {}from 
(\ref{LwF=LwF}) we see that (\ref{zz:L}) is equivalent to the
condition
\begin{equation}
G(L'(j)w'_{(4)})= L'_{P(z_1,z_2)}(j)G(w'_{(4)})
\end{equation}
for $j=-1, 0, 1$, that is, the condition that $G$ intertwines the
actions of $L'(j)$ and $L'_{P(z_1,z_2)}(j)$.
\epfv

Let $W_{1}$, $W_{2}$ and $W_{3}$ be generalized $V$-modules.  By
analogy with (\ref{W1W2_[C]^Atilde}) and (\ref{W1W2_(C)^Atilde}), we
have the spaces
\begin{equation}\label{3-mod-2-gradings}
((W_1\otimes W_2\otimes W_3)^{*})_{[\C]}^{(\tilde{A})}
=\coprod_{n\in \C}\coprod_{\beta\in \tilde{A}}
((W_1\otimes W_2\otimes W_3)^{*})_{[n]}^{(\beta)}\subset
(W_1\otimes W_2\otimes W_3)^{*}
\end{equation}
and 
\begin{equation}\label{3-mod-2-s-gradings}
((W_1\otimes W_2\otimes W_3)^{*})_{(\C)}^{(\tilde{A})}
=\coprod_{n\in \C}\coprod_{\beta\in \tilde{A}}
((W_1\otimes W_2\otimes W_3)^{*})_{(n)}^{(\beta)}\subset
(W_1\otimes W_2\otimes W_3)^{*},
\end{equation}
defined by means of the operator $L'_{P(z_1,z_2)}(0)$.

Again by analogy with the situation in Section 5, consider the
following conditions for elements $\lambda \in (W_1\otimes W_2\otimes
W_3)^*$:

\begin{description}
\item{\bf The $P(z_1, z_2)$-compatibility condition}

(a) The \emph{$P(z_1, z_2)$-lower truncation condition}: For all $v\in
V$, the formal Laurent series $Y'_{P(z_1,z_2)}(v,x)\lambda$ involves
only finitely many negative powers of $x$.

(b) The following formula holds for all $v\in V$:
\begin{eqnarray}\label{zz:cpb}
\lefteqn{\tau_{P(z_1,z_2)}\bigg(\dlt{x_1}{x_0^{-1}}{-z_1}\dlt{x_2}{x_0^{-1}}
{-z_2}Y_{t}(v,x_0)\bigg)\lambda}\nno\\
&&=\dlt{x_1}{x_0^{-1}}{-z_1}\dlt{x_2}{x_0^{-1}}{-z_2}Y'_{P(z_1,z_2)}
(v,x_0)\lambda.
\end{eqnarray}
(Note that the two sides of (\ref{zz:cpb}) are not {\it a priori}
equal for general $\lambda\in (W_1\otimes W_2\otimes W_3)^{*}$.
Condition (a) implies that the right-hand side in Condition (b) is
well defined.)

\item{\bf The $P(z_1,z_2)$-local grading restriction condition}

(a) The {\em $P(z_1,z_2)$-grading condition}: There exists a doubly
graded subspace of the space (\ref{3-mod-2-gradings}) containing
$\lambda$ and stable under the component operators
$\tau_{P(z_1,z_2)}(v\otimes t^{m})$ of the operators
$Y'_{P(z_1,z_2)}(v,x)$ for $v\in V$, $m\in {\mathbb Z}$, and under the
operators $L'_{P(z_1,z_2)}(-1)$, $L'_{P(z_1,z_2)}(0)$ and
$L'_{P(z_1,z_2)}(1)$. In particular, $\lambda$ is a (finite) sum of
generalized eigenvectors for $L'_{P(z_1,z_2)}(0)$ that are also
homogeneous with respect to $\tilde A$.

(b) Let $W_{\lambda; P(z_{1}, z_{2})}$ be the smallest doubly graded
(or equivalently, $\tilde A$-graded) subspace of the space
(\ref{3-mod-2-gradings}) containing $\lambda$ and stable under the
component operators $\tau_{P(z_1,z_2)}(v\otimes t^{m})$ of the
operators $Y'_{P(z_1,z_2)}(v,x)$ for $v\in V$, $m\in {\mathbb Z}$, and
under the operators $L'_{P(z_1,z_2)}(-1)$, $L'_{P(z_1,z_2)}(0)$ and
$L'_{P(z_1,z_2)}(1)$ (the existence being guaranteed by Condition
(a)).  Then $W_{\lambda; P(z_{1}, z_{2})}$ has the properties
\begin{eqnarray}
&\dim(W_{\lambda; P(z_{1}, z_{2})})^{(\beta)}_{[n]}<\infty,&\\
&(W_{\lambda; P(z_{1}, z_{2})})^{(\beta)}_{[n+k]}=0\;\;
\mbox{ for }\;k\in {\mathbb Z}
\;\mbox{ sufficiently negative,}&
\end{eqnarray}
for any $n\in {\mathbb C}$ and $\beta\in \tilde A$, where the
subscripts denote the ${\mathbb C}$-grading by (generalized)
$L'_{P(z_1,z_2)}(0)$-eigenvalues and the superscripts denote the
$\tilde A$-grading.

\item{\bf The $L(0)$-semisimple $P(z_{1}, z_{2})$-local 
grading restriction condition}

(a) The {\em $L(0)$-semisimple $P(z_1,z_2)$-grading condition}: There
exists a doubly graded subspace of the space
(\ref{3-mod-2-s-gradings}) containing $\lambda$ and stable under the
component operators $\tau_{P(z_1,z_2)}(v\otimes t^{m})$ of the
operators $Y'_{P(z_1,z_2)}(v,x)$ for $v\in V$, $m\in {\mathbb Z}$, and
under the operators $L'_{P(z_1,z_2)}(-1)$, $L'_{P(z_1,z_2)}(0)$ and
$L'_{P(z_1,z_2)}(1)$. In particular, $\lambda$ is a (finite) sum of
eigenvectors for $L'_{P(z_1,z_2)}(0)$ that are also homogeneous with
respect to $\tilde A$.

(b) Consider $W_{\lambda; P(z_{1}, z_{2})}$ as above, which in this
case is in fact the smallest doubly graded subspace of the space
(\ref{3-mod-2-s-gradings}) containing $\lambda$ and stable under the
component operators $\tau_{P(z_1,z_2)}(v\otimes t^{m})$ of the
operators $Y'_{P(z_1,z_2)}(v,x)$ for $v\in V$, $m\in {\mathbb Z}$, and
under the operators $L'_{P(z_1,z_2)}(-1)$, $L'_{P(z_1,z_2)}(0)$ and
$L'_{P(z_1,z_2)}(1)$.  Then $W_{\lambda; P(z_{1}, z_{2})}$ has the
properties
\begin{eqnarray}
&\dim(W_{\lambda; P(z_{1}, z_{2})})^{(\beta)}_{(n)}
<\infty,&\label{zz-semi-lgrc1}\\
&(W_{\lambda; P(z_{1}, z_{2})})^{(\beta)}_{(n+k)}=0\;\;
\mbox{ for }\;k\in {\mathbb Z}
\;\mbox{ sufficiently negative},&\label{zz-semi-lgrc2}
\end{eqnarray}
for any $n\in {\mathbb C}$ and $\beta\in \tilde A$, where the
subscripts denote the ${\mathbb C}$-grading by
$L'_{P(z_1,z_2)}(0)$-eigenvalues and the superscripts denote the
$\tilde A$-grading.

\end{description}

Then we have the following, by analogy with the comments preceding the
statement of the $P(z)$-compatibility condition (recall
(\ref{5.18-p})) and the $P(z)$-local grading restriction conditions:

\begin{propo}\label{8.12}
Suppose that $G\in \hom (W'_4, (W_1\otimes W_2\otimes W_3)^*)$
corresponds to a $P(z_1, z_2)$-intertwining map as in Proposition
\ref{zzcor}. Then for any $w'_{(4)}\in W'_4$, $G(w'_{(4)})$ satisfies
the $P(z_1, z_2)$-compatibility condition
and the $P(z_1,z_2)$-local grading restriction condition.
If $W_{4}$ is an ordinary $V$-module, then $G(w'_{(4)})$
satisfies the $L(0)$-semisimple $P(z_1,z_2)$-local grading 
restriction condition.
\end{propo}
\pf For any $w'_{(4)}\in W'_4$, the fact that $G(w'_{(4)})$ satisfies
the $P(z_1, z_2)$-compatibility condition follows {}from
(\ref{Gtau=tauG}), just as in (\ref{5.18-p}).  Since $G$ in particular
intertwines the actions of $V\otimes{\mathbb C}[t, t^{-1}]$ and of the
$L(j)$-operators and is $\tilde A$-compatible, $G(W'_4)$ is a
generalized $V$-module and thus $G(w'_{(4)})$ satisfies the
$P(z_1,z_2)$-local grading restriction condition, and if $W_{4}$ is an
ordinary $V$-module, then $G(W'_4)$ must also be an ordinary
$V$-module and thus $G(w'_{(4)})$ satisfies the $L(0)$-semisimple
$P(z_1,z_2)$-local grading restriction condition, just as in the
comments preceding the statement of the $P(z)$-local grading
restriction conditions. \epf

\begin{rema}\label{consequenceofPz1z2compat}
{\rm In the next section we will use the following: Assume the
$P(z_1,z_2)$-compatibility condition.  By (\ref{l1}) (a ``purely
algebraic'' identity, involving no convergence issues), (\ref{zz:cpb})
can be written as:
\begin{eqnarray}\label{alternatecompat}
\lefteqn{\tau_{P(z_1,z_2)}\bigg(\dlt{x_1}{x_2}{-z_0}\dlt{x_2}{x_0^{-1}}
{-z_2}Y_{t}(v,x_0)\bigg)\lambda}\nn
&&=\dlt{x_1}{x_2}{-z_0}\biggl(\dlt{x_2}{x_0^{-1}}{-z_2}Y'_{P(z_1,z_2)}
(v,x_0)\lambda\biggr);
\end{eqnarray}
the right-hand side is defined because of the $P(z_1, z_2)$-lower
truncation condition.  We can take ${\rm Res}_{x_{1}}$ to obtain
\begin{eqnarray}\label{resofconsequence}
\tau_{P(z_1,z_2)}\bigg(\dlt{x_2}{x_0^{-1}}{-z_2}Y_{t}(v,x_0)\bigg)\lambda
=\dlt{x_2}{x_0^{-1}}{-z_2}Y'_{P(z_1,z_2)} (v,x_0)\lambda,
\end{eqnarray}
which is reminiscent of the $P(z)$-compatibility condition
(\ref{cpb}) for $z=z_2$.  Now we can multiply both sides by
$\displaystyle\dlt{x_1}{x_2}{-z_0}$, giving
\begin{eqnarray*}
\lefteqn{\dlt{x_1}{x_2}{-z_0}\tau_{P(z_1,z_2)}\bigg(\dlt{x_2}{x_0^{-1}}
{-z_2}Y_{t}(v,x_0)\bigg)\lambda}\nn
&&=\dlt{x_1}{x_2}{-z_0}\dlt{x_2}{x_0^{-1}}{-z_2}Y'_{P(z_1,z_2)}
(v,x_0)\lambda,
\end{eqnarray*}
and thus by (\ref{zz:cpb}) and (\ref{alternatecompat}),
\begin{eqnarray}\label{consequenceofPz1z2compatformula}
\lefteqn{\tau_{P(z_1,z_2)}\bigg(\dlt{x_1}{x_0^{-1}}{-z_1}\dlt{x_2}{x_0^{-1}}
{-z_2}Y_{t}(v,x_0)\bigg)\lambda}\nn
&&=\dlt{x_1}{x_2}{-z_0}\tau_{P(z_1,z_2)}\bigg(\dlt{x_2}{x_0^{-1}}
{-z_2}Y_{t}(v,x_0)\bigg)\lambda.
\end{eqnarray}}
\end{rema}

\newpage

\setcounter{equation}{0}
\setcounter{rema}{0}

\section{The expansion condition}\label{extsec}

In Section 7 we studied the conditions necessary for products and
iterates of certain intertwining maps to exist. Assuming that these
conditions are satisfied, in this section we study the condition for
the product of two suitable intertwining maps to be expressible as the
iterate of some suitable intertwining maps, and vice versa.  To do
this, following the idea in \cite{tensor4}, we first study certain
properties satisfied by the product, and separately, the iterate of
two intertwining maps.  Then we prove that the condition that the
product satisfies its properties and the condition that the iterate
satifies its properties are equivalent to each other.  In this work,
we introduce the term ``expansion condition'' for either of these
equivalent conditions.  We show that if a product or iterate of two
intertwining maps satisfies the expansion condition, then it can be
expressed in the other form.  These results are generalizations to the
logarithmic setting of the corresponding results in the finitely
reductive case obtained in \cite{tensor4}.

Unless otherwise stated, in this section $z_1$ and $z_2$ will be
distinct nonzero complex numbers, and $z_0=z_1-z_2$.

For objects $W_1$, $W_2$, $W_3$, $W_4$, $M_1$ and $M_2$ of ${\cal C}$,
let $I_1$, $I_2$, $I^1$ and $I^2$ be $P(z_1)$-, $P(z_2)$-, $P(z_2)$-
and $P(z_0)$-intertwining maps of types ${W_4 \choose W_1\, M_1}$,
${M_1 \choose W_2\, W_3}$, ${W_4 \choose M_2\, W_3}$ and ${M_2\choose
W_1\, W_2}$, respectively.  Then under the assumption of the
convergence condition for intertwining maps in ${\cal C}$ (recall
Proposition \ref{convergence} and Definition \ref{conv-conditions}),
when $|z_1|>|z_2|>|z_0|>0$, both the product $I_1\circ (1_{W_1}\otimes
I_2)$ and the iterate $I^1\circ (I^2\otimes 1_{W_3})$ exist and are
$P(z_1,z_2)$-intertwining maps, by Proposition
\ref{productanditerateareintwmaps}.  In this section we will consider
when such a product can be expressed as such an iterate and vice
versa.

To compare these two types of maps, we shall study some conditions
specific to each type.  For this, let $W_1$, $W_2$ and $W_3$ be
arbitrary generalized $V$-modules.  We start with the following:

\begin{defi}\label{mudef}{\rm
Let $\lambda \in (W_1\otimes W_2\otimes W_3)^*$.  For $w_{(1)}\in
W_{1}$, we define the {\it evaluation of $\lambda$ at $w_{(1)}$} to be
the element $\mu^{(1)}_{\lambda, w_{(1)}}$ of $(W_2\otimes W_3)^{*}$
given by
\[
\mu^{(1)}_{\lambda, w_{(1)}}(w_{(2)}\otimes w_{(3)})
=\lambda(w_{(1)}\otimes w_{(2)}\otimes w_{(3)})
\]
for $w_{(2)}\in W_2$ and $w_{(3)}\in W_3$.  For $w_{(3)}\in W_3$, we
define the {\it evaluation of $\lambda$ at $w_{(3)}$} to be the
element $\mu^{(2)}_{\lambda, w_{(3)}}$ of $(W_1\otimes W_2)^{*}$ given
by
\[
\mu^{(2)}_{\lambda, w_{(3)}}(w_{(1)}\otimes w_{(2)})
=\lambda(w_{(1)}\otimes w_{(2)}\otimes w_{(3)})
\]
for $w_{(1)}\in W_1$ and $w_{(2)}\in W_2$.
}
\end{defi}

\begin{rema}
{\rm Given $\lambda\in (W_1\otimes W_2\otimes W_3)^{*}$, $w_{(1)}\in
W_{1}$ and $w_{(3)}\in W_3$, it is natural to ask whether the
evaluations $\mu^{(1)}_{\lambda, w_{(1)}}\in (W_{2}\otimes W_{3})^{*}$
and $\mu^{(2)}_{\lambda, w_{(3)}}\in (W_{1}\otimes W_{2})^{*}$ of
$\lambda$ satisfy the $P(z)$-compatibility condition (recall
(\ref{cpb})) for some suitable nonzero complex numbers $z$, under
suitable conditions.  In fact, the formal computations underlying the
next lemma suggest that when $\lambda$ satisfies the $P(z_{1},
z_{2})$-compatibility condition (recall (\ref{zz:cpb})), these
evaluations of $\lambda$ ``almost'' satisfy the
$P(z_{2})$-compatibility condition and the $P(z_{0})$-compatibility
condition, respectively.  However, even when $\lambda$ does satisfy
the $P(z_{1}, z_{2})$-compatibility condition, in general these
evaluations of $\lambda$ do not even satisfy the $P(z_{2})$-lower
truncation condition (Part (a) of the $P(z_{2})$-compatibility
condition) or the $P(z_{0})$-lower truncation condition (Part (a) of
the $P(z_{0})$-compatibility condition).  In particular, when
$\lambda$ in (\ref{cpb}) is replaced by $\mu^{(1)}_{\lambda, w_{(1)}}$
and $z=z_{2}$, the right-hand side of (\ref{cpb}) might not exist in
the usual algebraic sense, and similarly, when $\lambda$ in
(\ref{cpb}) is replaced by $\mu^{(2)}_{\lambda, w_{(3)}}$ and
$z=z_{0}$, the right-hand side of (\ref{cpb}) might not exist
algebraically, and for this reason $\mu^{(1)}_{\lambda, w_{(1)}}$ and
$\mu^{(2)}_{\lambda, w_{(3)}}$ do not in general satisfy the
$P(z_{2})$-compatibility condition or the $P(z_{0})$-compatibility
condition.  But the next result, which generalizes Lemma 14.3 in
\cite{tensor4}, asserts that if $\lambda$ satisfies the $P(z_1,
z_2)$-compatibility condition, then in both cases, under the natural
assumptions on the complex numbers, the right-hand side of (\ref{cpb})
exists {\it analytically} and (\ref{cpb}) holds {\it analytically}, in
the sense of weak absolute convergence, as discussed in Remark
\ref{weakly-abs-conv}.}
\end{rema}

\begin{lemma}\label{mulemma}
Assume that $\lambda\in (W_1\otimes W_2\otimes W_3)^{*}$ satisfies the
$P(z_1, z_2)$-compatibility condition (recall (\ref{zz:cpb})).  If
$|z_2|>|z_0|$ $(>0)$, then for any $v\in V$ and $w_{(1)}\in W_1$,
$w_{(2)}\in W_2$ and $w_{(3)}\in W_3$, the coefficients of the
monomials in $x$ and $x_{1}$ in
\[
x^{-1}_1 \delta\bigg(\frac{x^{-1}-z_0}{x_1}\bigg)
\biggl(Y'_{P(z_0)}(v, x) \mu^{(2)}_{\lambda, w_{(3)}}\biggr)
(w_{(1)} \otimes w_{(2)})
\]
are absolutely convergent and we have
\begin{eqnarray}\label{mu12}
\lefteqn{\biggl(\tau_{P(z_0)}\biggl( x^{-1}_1
\delta\bigg(\frac{x^{-1}-z_0}{x_1}\bigg) Y_{t}(v, x)\biggr)
\mu^{(2)}_{\lambda, w_{(3)}}\biggr)(w_{(1)} \otimes w_{(2)})}\nno\\
&&=x^{-1}_1 \delta\bigg(\frac{x^{-1}-z_0}{x_1}\bigg)
\biggl(Y'_{P(z_0)}(v, x) \mu^{(2)}_{\lambda, w_{(3)}}\biggr)
(w_{(1)} \otimes w_{(2)}).
\end{eqnarray}
Analogously, if $|z_1|>|z_2|$ $(>0)$, then for any $v\in V$ and any
$w_{(j)}\in W_j$, the coefficients of the monomials in $x$ and $x_{1}$
in
\[
x^{-1}_1 \delta\bigg(\frac{x^{-1}-z_2}{x_1}\bigg)\biggl(Y'_{P(z_2)}(v, x)
\mu^{(1)}_{\lambda, w_{(1)}}\biggr)(w_{(2)} \otimes w_{(3)})
\]
are absolutely convergent and we have
\begin{eqnarray}\label{mu23}
\lefteqn{\biggl(\tau_{P(z_2)}\biggl( x^{-1}_1
\delta\bigg(\frac{x^{-1}-z_2}{x_1}\bigg) Y_{t}(v, x)\biggr)
\mu^{(1)}_{\lambda, w_{(1)}}\biggr)(w_{(2)} \otimes w_{(3)})}\nno\\
&&=x^{-1}_1 \delta\bigg(\frac{x^{-1}-z_2}{x_1}\bigg)
\biggl(Y'_{P(z_2)}(v, x) \mu^{(1)}_{\lambda, w_{(1)}}\biggr)
(w_{(2)} \otimes w_{(3)}).
\end{eqnarray}
\end{lemma}
\pf First, for our distinct nonzero complex numbers $z_1$ and $z_2$
with $z_0=z_1-z_2$, by definition of the action $\tau _{P(z_1,z_2)}$
(\ref{tauzzgf}) we have
\begin{eqnarray}\label{lm:1}
\lefteqn{\dlti{x_0}{z_1}{+x_1}\dlt{x_2}{z_0}{+x_1} \lambda
(Y_1((-x_0^{-2})^{L(0)}e^{-x_0^{-1}L(1)}v,x_1)w_{(1)}\otimes
w_{(2)}\otimes w_{(3)})}\nno\\
&&+ \dlti{x_0}{z_2}{+x_2}\dlt{x_1}{-z_0}{+x_2}\lambda (w_{(1)}\otimes
Y_2((-x_0^{-2})^{L(0)}e^{-x_0^{-1}L(1)}v,x_2)w_{(2)}\otimes
w_{(3)})\nno\\
&&= \left(\tau _{P(z_1,z_2)}\left(\dlt{x_1}{x_0^{-1}}{-z_1}\dlt{x_2}{x_0^{-1}}
{-z_2}Y_{t}(v,x_0)\right)\lambda \right)(w_{(1)}\otimes w_{(2)}\otimes
w_{(3)})\nno\\
&&\quad-\dlt{x_1}{-z_1}{+x^{-1}_0}\dlt{x_2}{-z_2}{+x^{-1}_0}\lambda
(w_{(1)}\otimes w_{(2)}\otimes Y_3^o(v,x_0)w_{(3)})
\end{eqnarray}
for any $v\in V$, $w_{(1)}\in W_1$, $w_{(2)}\in W_2$ and $w_{(3)}\in W_3$.
Replacing $v$ by
\[
(-x_0^2)^{L(0)}e^{-x_0L(1)}e^{x_2^{-1}L(1)}(-x_2^2)^{L(0)}v,
\]
using formula (5.3.1) in \cite{FHL}, and then
taking $\res_{x_0^{-1}}$ we get:
\begin{eqnarray}\label{lm4}
\lefteqn{\dlt{x_2}{z_0}{+x_1}\lambda
(Y_1(e^{x_2^{-1}L(1)}(-x_2^2)^{L(0)}v,x_1)w_{(1)}\otimes
w_{(2)}\otimes w_{(3)})}\nno\\
&& +\dlt{x_1}{-z_0}{+x_2}\lambda (w_{(1)}\otimes
Y_2(e^{x_2^{-1}L(1)}(-x_2^2)^{L(0)}v,x_2)w_{(2)}\otimes
w_{(3)})\nno\\
&&= {\rm
Res}_{x_0^{-1}}\Bigg(\tau_{P(z_1,z_2)}\Bigg(\dlt{x_1}{x_0^{-1}}{-z_1}
\dlt{x_2}{x_0^{-1}}{-z_2}\cdot\nno\\
&&\qquad\cdot
Y_{t}((-x_0^2)^{L(0)}e^{-x_0L(1)}e^{x_2^{-1}L(1)}(-x_2^2)^{L(0)}v,x_0)
\Bigg)\lambda
\Bigg)(w_{(1)}\otimes w_{(2)}\otimes w_{(3)}))\nno\\
&&\quad -{\rm
Res}_{x_0^{-1}}\dlt{x_1}{-z_1}{+x^{-1}_0}\dlt{x_2}{-z_2}{+x^{-1}_0}\cdot\nno\\
&&\qquad \cdot\lambda (w_{(1)}\otimes w_{(2)}\otimes
Y_3(e^{x_2^{-1}L(1)}(-x_2^2)^{L(0)}v,x_0^{-1})w_{(3)}).
\end{eqnarray}

By (\ref{taudef}), the left-hand side of (\ref{lm4}) is equal to
\begin{eqnarray}\label{lefthandside}
\left(\tau_{P(z_0)}\left(\dlt{x_1}{x_2}{-z_0}Y_t(v,x_2^{-1})\right)\mu^{(2)}_{\lambda,
w_{(3)}}\right)(w_{(1)}\otimes w_{(2)}).
\end{eqnarray}
Taking $\res_{x_1}$ gives
\begin{eqnarray}\label{resoflefthandside}
(\tau_{P(z_0)}(Y_t(v,x_2^{-1}))\mu^{(2)}_{\lambda,
w_{(3)}})(w_{(1)}\otimes
w_{(2)})=(Y'_{P(z_0)}(v,x_2^{-1})\mu^{(2)}_{\lambda,
w_{(3)}})(w_{(1)}\otimes w_{(2)}).
\end{eqnarray}

By the $P(z_1, z_2)$-compatibility condition and formula
(\ref{consequenceofPz1z2compatformula}) in Remark
\ref{consequenceofPz1z2compat}, the first term on the right-hand side
of (\ref{lm4}) equals
\begin{eqnarray}\label{RHSexpression}
\lefteqn{\dlt{x_1}{x_2}{-z_0}{\rm Res}_{x_0^{-1}}\Bigg(\tau_{P(z_1,z_2)}\Bigg(
\dlt{x_2}{x_0^{-1}}{-z_2}\cdot}\nno\\
&&\qquad\cdot
Y_{t}((-x_0^2)^{L(0)}e^{-x_0L(1)}e^{x_2^{-1}L(1)}(-x_2^2)^{L(0)}v,x_0)
\Bigg)\lambda \Bigg)(w_{(1)}\otimes w_{(2)}\otimes w_{(3)})).
\end{eqnarray}

Now suppose that $|z_2|>|z_0|$ $(>0)$.  Then formula (\ref{l4})
holds. {}From this and (\ref{RHSexpression}), the right-hand side of
(\ref{lm4}) becomes
\begin{eqnarray*}
\lefteqn{\dlt{x_1}{x_2}{-z_0}{\rm
Res}_{x_0^{-1}}\Biggl(\tau_{P(z_1,z_2)}\Biggl(
\dlt{x_2}{x_0^{-1}}{-z_2}\cdot}\\
&&\qquad\cdot
Y_t((-x_0^2)^{L(0)}e^{-x_0L(1)}e^{x_2^{-1}L(1)}(-x_2^2)^{L(0)}v,x_0)
\Biggr)\lambda\Biggr)(w_{(1)}\otimes w_{(2)}\otimes w_{(3)})\\
&&-\dlt{x_1}{x_2}{-z_0}{\rm
Res}_{x_0^{-1}}\dlt{x_2}{-z_2}{+x_0^{-1}}\cdot\\
&&\qquad \cdot\lambda (w_{(1)}\otimes w_{(2)}\otimes
Y_3(e^{x_2^{-1}L(1)}(-x_2^2)^{L(0)}v,x_0^{-1})w_{(3)}).
\end{eqnarray*}

By taking $\res_{x_1}$ of this expression, we erase the two factors
$\displaystyle\dlt{x_1}{x_2}{-z_0}$, leaving an expression which,
while not lower truncated in $x_2^{-1}$, can still be multiplied by
$\displaystyle\dlt{x_1}{x_2}{-z_0}$ (when $|z_2|>|z_0|>0$), in the
sense of absolute convergence, yielding this expression again.  That
is, let $X$ be either side of (\ref{lm4}).  Then
\[
X=\dlt{x_1}{x_2}{-z_0}\res_{x_1}X,
\]
in this sense of convergence.  Applying this to (\ref{lefthandside})
and (\ref{resoflefthandside}) gives
\begin{eqnarray}\label{lmu12}
\lefteqn{\tau_{P(z_0)}\left(\dlt{x_1}{x_2}{-z_0}Y_t(v,x_2^{-1})\right)
\mu^{(2)}_{\lambda,w_{(3)}}(w_{(1)}\otimes w_{(2)})}\nno\\
&&\qquad
=\dlt{x_1}{x_2}{-z_0}\biggl(Y'_{P(z_0)}(v,x_2^{-1})\mu^{(2)}_{\lambda,
w_{(3)}}\biggr)(w_{(1)}\otimes w_{(2)}),
\end{eqnarray}
proving (\ref{mu12}) (with $x$ in (\ref{mu12}) replaced by
$x_2^{-1}$).

Analogously, for our distinct nonzero complex numbers $z_1$ and $z_2$
with $z_0=z_1-z_2$, we can also write the definition of
$\tau_{P(z_1,z_2)}$ as
\begin{eqnarray}\label{lm:2}
\lefteqn{\dlti{x_0}{z_2}{+x_2}\dlt{x_1}{-z_0}{+x_2} \lambda (w_{(1)}\otimes
Y_2((-x_0^{-2})^{L(0)}e^{-x_0^{-1}L(1)}v,x_2)w_{(2)}\otimes
w_{(3)})}\nno\\
&&+\dlt{x_1}{-z_1}{+x^{-1}_0}\dlt{x_2}{-z_2}{+x^{-1}_0}\lambda
(w_{(1)}\otimes w_{(2)}\otimes Y_3^o(v,x_0)w_{(3)})\nno\\
&&=\left(\tau_{P(z_1,z_2)}\left(\dlt{x_1}{x_0^{-1}}{-z_1}\dlt{x_2}{x_0^{-1}}
{-z_2}Y_{t}(v,x_0)\right)\lambda \right)(w_{(1)}\otimes w_{(2)}\otimes w_{(3)})\nno\\
&  &\quad- \dlti{x_0}{z_1}{+x_1}\dlt{x_2}{z_0}{+x_1}\lambda
(Y_1((-x_0^{-2})^{L(0)}e^{-x_0^{-1}L(1)}v,x_1)w_{(1)}\otimes
w_{(2)}\otimes w_{(3)})\nno\\
\end{eqnarray}
for $v\in V$, $w_{(1)}\in W_1$, $w_{(2)}\in W_2$ and $w_{(3)}\in W_3$.
Taking $\res_{x_1}$ we get
\begin{eqnarray}\label{lm7}
\lefteqn{\dlti{x_0}{z_2}{+x_2}\lambda (w_{(1)}\otimes
Y_2((-x_0^{-2})^{L(0)}e^{-x_0^{-1}L(1)}v,x_2)w_{(2)}\otimes
w_{(3)})}\nno\\
&  &+ \dlt{x_2}{-z_2}{+x^{-1}_0}\lambda (w_{(1)}\otimes
w_{(2)}\otimes Y^o_3(v,x_0)w_{(3)})\nno\\
& &=\left(\tau
_{P(z_1,z_2)}\left(\dlt{x_2}{x_0^{-1}}{-z_2}Y_{t}(v,x_0)\right)\lambda
\right)(w_{(1)}\otimes w_{(2)}\otimes w_{(3)})\nno\\
& &\quad-\res_{x_1} \dlti{x_0}{z_1}{+x_1}\dlt{x_2}{z_0}{+x_1}\cdot\nno\\
&&\qquad\qquad\cdot\lambda
(Y_1((-x_0^{-2})^{L(0)}e^{-x_0^{-1}L(1)}v,x_1)w_{(1)}\otimes
w_{(2)}\otimes w_{(3)}).
\end{eqnarray}

By the definition of $\tau_{P(z_2)}$ (\ref{taudef}) and formula
(5.3.1) in \cite{FHL}, the left-hand side of (\ref{lm7}) is equal to
\[
\left(\tau_{P(z_2)}\left(\dlt{x_2}{x_0^{-1}}{-z_2}Y_t(v,x_0)\right)\mu^{(1)}_{\lambda,
w_{(1)}}\right)(w_{(2)}\otimes w_{(3)}),
\]
and taking ${\rm Res}_{x_2}$ gives
\[
(\tau_{P(z_2)}(Y_t(v,x_0))\mu^{(1)}_{\lambda,
w_{(1)}})(w_{(2)}\otimes w_{(3)})\\
=(Y'_{P(z_2)}(v,x_0)\mu^{(1)}_{\lambda, w_{(1)}})(w_{(2)}\otimes
w_{(3)}).
\]

Now suppose that $|z_1|>|z_2|>0$.  Then (\ref{l2-1}) holds, by
(\ref{resofconsequence}), which follows from the
$P(z_1,z_2)$-compatibility condition, the right-hand side of
(\ref{lm7}) becomes
\begin{eqnarray*}
\lefteqn{\dlt{x_2}{x_0^{-1}}{-z_2}(Y'_{P(z_1,z_2)}(v,x_0)
\lambda )(w_{(1)}\otimes w_{(2)}\otimes w_{(3)})}\\
&&-\dlt{x_2}{x_0^{-1}}{-z_2}\res_{x_1}\dlti{x_0}{z_1}{+x_1}\cdot\\
&&\qquad\cdot\lambda (Y_1((-x_0^{-2})^{L(0)}e^{-x_0^{-1}L(1)}v,x_1)
w_{(1)}\otimes w_{(2)}\otimes w_{(3)}).
\end{eqnarray*}
Just as in the proof of (\ref{mu12}), we take $\res_{x_2}$ and then
multiply by $\displaystyle\dlt{x_2}{x_0^{-1}}{-z_2}$, yielding the
same expression, and we obtain
\begin{eqnarray}\label{lmu23}
\lefteqn{\biggl(\tau_{P(z_2)}\biggl(\dlt{x_2}{x_0^{-1}}{-z_2} Y_{t}(v,x_0)\biggr)
\mu^{(1)}_{\lambda, w_{(1)}}\biggr)
(w_{(2)}\otimes w_{(3)})}\nno\\
&&\qquad
=\dlt{x_2}{x_0^{-1}}{-z_2}\biggl(Y'_{P(z_2)}(v, x_0)
\mu^{(1)}_{\lambda,w_{(1)}}\biggr)(w_{(2)}\otimes w_{(3)}),
\end{eqnarray}
proving (\ref{mu23}).
\epf

\begin{rema}{\rm
As we discussed above, Lemma \ref{mulemma} says that under the
appropriate conditions, $\mu^{(2)}_{\lambda, w_{(3)}}$ and
$\mu^{(1)}_{\lambda, w_{(1)}}$ satisfy natural analytic analogues of
the $P(z_0)$-compatibility condition and the $P(z_2)$-compatibility
condition, respectively.  Note that from the proof of (\ref{mu12}),
$(Y'_{P(z_0)}(v,x_2^{-1})\mu^{(2)}_{\lambda, w_{(3)}})(w_{(1)}\otimes
w_{(2)})$ ``behaves qualitatively'' like
$\displaystyle\delta\bigg(\frac{z_2}{-x_2}\bigg)$ and so (\ref{lmu12})
``behaves qualitatively'' like
\[
\dlt{x_1}{x_2}{-z_0}\delta\bigg(\frac{z_2}{-x_2}\bigg),
\]
suggesting the expected convergence when $|z_2|>|z_0|$.  Analogously,
from the proof of (\ref{mu23}), $Y'_{P(z_2)}(v,x_0)\mu^{(1)}_{\lambda,
w_{(1)}}(w_{(2)}\otimes w_{(3)})$ ``behaves qualitatively'' like
$\displaystyle\delta\bigg(\frac{z_1}{x_0^{-1}}\bigg)$ and so
(\ref{lmu23}) ``behaves qualitatively'' like
\[
\dlt{x_2}{x_0^{-1}}{-z_2}\delta\bigg(\frac{z_1}{x_0^{-1}}\bigg),
\]
again suggesting the expected convergence, this time when
$|z_1|>|z_2|$.  (By the $P(z_1 ,z_2)$-lower truncation condition,
${\rm Res}_{x_1}$ of (\ref{lefthandside}) is upper-truncated in $x_2$,
independently of $w_{(1)}$, $w_{(2)}$ and $w_{(3)}$, and ${\rm
Res}_{x_2}$ of the first term on the right-hand side of (\ref{lm7}),
$Y'_{P(z_1,z_2)}(v,x_0)\lambda)(w_{(1)}\otimes w_{(2)}\otimes
w_{(3)})$, is lower-truncated in $x_0$, independently of the
$w_{(j)}$.)}
\end{rema}

\begin{rema}
{\rm Given $\lambda\in (W_1\otimes W_2\otimes W_3)^{*}$, $w_{(1)}\in
W_{1}$ and $w_{(3)}\in W_3$, it is also natural to ask whether, under
suitable conditions, the evaluations $\mu^{(1)}_{\lambda, w_{(1)}}$
and $\mu^{(2)}_{\lambda, w_{(3)}}$ of $\lambda$ satisfy the
$P(z_{2})$-local grading restriction condition and $P(z_{0})$-local
grading restriction condition, respectively.  In general, even for
$\lambda$ obtained {}from a product or an iterate of intertwining
maps, these conditions are not satisfied by $\mu^{(1)}_{\lambda,
w_{(1)}}$ or $\mu^{(2)}_{\lambda, w_{(3)}}$, but as we will see below,
for such $\lambda$, $\mu^{(1)}_{\lambda, w_{(1)}}$ and
$\mu^{(2)}_{\lambda, w_{(3)}}$ satisfy certain generalizations of
these conditions.  These generalizations motivate the next four
important conditions on $\lambda\in (W_1\otimes W_2\otimes W_3)^*$.
An element $\lambda\in (W_1\otimes W_2\otimes W_3)^*$ satisfying one
of these conditions means essentially that either $\mu^{(1)}_{\lambda,
w_{(1)}}\in (W_2\otimes W_3)^{*}$ or $\mu^{(2)}_{\lambda, w_{(3)}}\in
(W_1\otimes W_2)^{*}$ is the sum of a weakly convergent series, in the
sense of Remark \ref{weakly-abs-conv}, rather than just a finite sum,
of $P(z)$-generalized weight vectors or of ordinary weight vectors
satisfying the $P(z)$-local grading restriction condition or the
$L(0)$-semisimple $P(z)$-local grading restriction condition in
Section 5, and that all of the summands lie in the same subspace whose
grading is restricted.  We shall typically use these conditions for
$z=z_2$ when we consider $\mu^{(1)}_{\lambda, w_{(1)}}$ and for
$z=z_{0}$ when we consider $\mu^{(2)}_{\lambda, w_{(3)}}$.  Recall the
spaces (\ref{W1W2_[C]^Atilde}) and (\ref{W1W2_(C)^Atilde}).}
\end{rema}

\begin{description}
\item{\bf The $P^{(1)}(z)$-local grading restriction condition}

(a) The {\em $P^{(1)}(z)$-grading condition}: For any $w_{(1)}\in
W_1$, there exists a formal series $\sum_{n\in \C}\lambda^{(1)}_{n}$
with
\[
\lambda^{(1)}_{n}\in \coprod_{\beta\in \tilde{A}}
((W_{2}\otimes W_{3})^{*})_{[n]}^{(\beta)}
\]
for $n\in \C$ such that for $w_{(2)}\in W_{2}$ and $w_{(3)}\in W_{3}$,
\[
\{(n, i)\in \C\times \N\;|\;
((L_{P(z)}'(0)-n)^{i}\lambda_{n}^{(1)})(w_{(2)}\otimes w_{(3)})\ne 0\}
\]
is a unique expansion set and
\[
\sum_{n\in \C}\lambda^{(1)}_{n}(w_{(2)}\otimes w_{(3)})
\]
is absolutely convergent to $\mu^{(1)}_{\lambda,
w_{(1)}}(w_{(2)}\otimes w_{(3)})$.

(b) For any $w_{(1)}\in W_1$, let $W^{(1)}_{\lambda, w_{(1)}}$ be the
smallest doubly graded (or equivalently, $\tilde A$-graded; recall
Remark \ref{singleanddoublegraded}) subspace of $((W_2\otimes
W_3)^{*})_{[\C]}^{(\tilde{A})}$ containing all the terms $\lambda^{(1)}_{n}$
in the formal series in (a) and stable under the component operators
$\tau_{P(z)}(v\otimes t^{m})$ of the operators $Y'_{P(z)}(v, x)$ for
$v\in V$, $m\in {\mathbb Z}$, and under the operators $L'_{P(z)}(-1)$,
$L'_{P(z)}(0)$ and $L'_{P(z)}(1)$.  (In view of Remark
\ref{stableundercomponentops}, $W^{(1)}_{\lambda, w_{(1)}}$ indeed
exists, just as in the case of the $P(z)$-local grading restriction
condition.)  Then $W^{(1)}_{\lambda, w_{(1)}}$ has the properties
\begin{eqnarray*}
&\dim(W^{(1)}_{\lambda, w_{(1)}})^{(\beta)}_{[n]}<\infty,&\\
&(W^{(1)}_{\lambda, w_{(1)}})^{(\beta)}_{[n+k]}=0\;\;\mbox{ for
}\;k\in {\mathbb Z} \;\mbox{ sufficiently negative}&
\end{eqnarray*}
for any $n\in {\mathbb C}$ and $\beta\in \tilde A$, where the
subscripts denote the ${\mathbb C}$-grading by
$L'_{P(z)}(0)$-(generalized) eigenvalues and the superscripts denote
the $\tilde A$-grading.
\end{description}

\begin{description}
\item{\bf The $L(0)$-semisimple 
$P^{(1)}(z)$-local grading restriction condition}

(a) The {\em $L(0)$-semisimple $P^{(1)}(z)$-grading condition}: For
any $w_{(1)}\in W_1$, there exists  a
formal series $\sum_{n\in \C}\lambda^{(1)}_{n}$ 
with
\[
\lambda^{(1)}_{n}\in \coprod_{\beta\in \tilde{A}}
((W_{2}\otimes W_{3})^{*})_{(n)}^{(\beta)}
\]
for $n\in \C$ such that for $w_{(2)}\in W_{2}$ and $w_{(3)}\in W_{3}$,
\[
\{(n, 0)\in \C\times \N\;|\;\lambda_{n}^{(1)}(w_{(2)}\otimes w_{(3)})\ne 0\}
\]
is a unique expansion set and
\[
\sum_{n\in \C}\lambda^{(1)}_{n}(w_{(2)}\otimes w_{(3)})
\]
is absolutely convergent to $\mu^{(1)}_{\lambda,
w_{(1)}}(w_{(2)}\otimes w_{(3)})$.  (Note that such an element
$\lambda$ also satisfies the $P^{(1)}(z)$-grading condition above with
the same elements $\lambda^{(1)}_{n}$.)

(b) For any $w_{(1)}\in W_1$, consider the space $W^{(1)}_{\lambda,
w_{(1)}}$ as above, which in this case is in fact the smallest doubly
graded (or equivalently, $\tilde A$-graded) subspace of $((W_2\otimes
W_3)^{*})_{(\C)}^{(\tilde{A})}$ containing all the terms $\lambda^{(1)}_{n}$
in the formal series in (a) and stable under the component operators
$\tau_{P(z)}(v\otimes t^{m})$ of the operators $Y'_{P(z)}(v, x)$ for
$v\in V$, $m\in {\mathbb Z}$, and under the operators $L'_{P(z)}(-1)$,
$L'_{P(z)}(0)$ and $L'_{P(z)}(1)$. Then $W^{(1)}_{\lambda, w_{(1)}}$
has the properties
\begin{eqnarray*}
&\dim(W^{(1)}_{\lambda, w_{(1)}})^{(\beta)}_{(n)}<\infty,&\\
&(W^{(1)}_{\lambda, w_{(1)}})^{(\beta)}_{(n+k)}=0\;\;\mbox{ for
}\;k\in {\mathbb Z} \;\mbox{ sufficiently negative}&
\end{eqnarray*}
for any $n\in {\mathbb C}$ and $\beta\in \tilde A$, where the
subscripts denote the ${\mathbb C}$-grading by
$L'_{P(z)}(0)$-eigenvalues and the superscripts denote the $\tilde
A$-grading.
\end{description}

\begin{description}
\item{\bf The $P^{(2)}(z)$-local grading restriction condition}

(a) The {\em $P^{(2)}(z)$-grading condition}: For any $w_{(3)}\in
W_3$,  there exists  a formal series $\sum_{n\in \C}\lambda^{(2)}_{n}$ 
with
\[
\lambda^{(2)}_{n}\in \coprod_{\beta\in \tilde{A}}
((W_{1}\otimes W_{2})^{*})_{[n]}^{(\beta)}
\]
for $n\in \C$ such that for $w_{(1)}\in W_{1}$ and $w_{(2)}\in W_{2}$,
\[
\{(n, i)\in \C\times \N\;|\;
((L_{P(z)}'(0)-n)^{i}\lambda_{n}^{(2)})(w_{(1)}\otimes w_{(2)})\ne 0\}
\]
is a unique expansion set and 
\[
\sum_{n\in \C}\lambda^{(2)}_{n}(w_{(1)}\otimes w_{(2)})
\]
is absolutely convergent to $\mu^{(2)}_{\lambda,
w_{(3)}}(w_{(1)}\otimes w_{(2)})$.

(b) For any $w_{(3)}\in W_3$, let $W^{(2)}_{\lambda, w_{(3)}}$ be the
smallest doubly graded (or equivalently, $\tilde A$-graded) subspace
of $((W_1\otimes W_2)^{*})_{[\C]}^{(\tilde{A})}$ containing all the
terms $\lambda^{(2)}_{n}$ in the formal series in (a) and stable under the
component operators $\tau_{P(z)}(v\otimes t^{m})$ of the operators
$Y'_{P(z)}(v, x)$ for $v\in V$, $m\in {\mathbb Z}$, and under the
operators $L'_{P(z)}(-1)$, $L'_{P(z)}(0)$ and $L'_{P(z)}(1)$.  (As
above, $W^{(2)}_{\lambda, w_{(3)}}$ indeed exists.)  Then
$W^{(2)}_{\lambda, w_{(3)}}$ has the properties
\begin{eqnarray*}
&\dim(W^{(2)}_{\lambda, w_{(3)}})^{(\beta)}_{[n]}<\infty,&\\
&(W^{(2)}_{\lambda, w_{(3)}})^{(\beta)}_{[n+k]}=0\;\;\mbox{ for
}\;k\in {\mathbb Z} \;\mbox{ sufficiently negative}&
\end{eqnarray*}
for any $n\in {\mathbb C}$ and $\beta\in \tilde A$, where the
subscripts denote the ${\mathbb C}$-grading by
$L'_{P(z)}(0)$-(generalized) eigenvalues and the superscripts denote
the $\tilde A$-grading.
\end{description}

\begin{description}
\item{\bf The $L(0)$-semisimple 
$P^{(2)}(z)$-local grading restriction condition}

(a) The {\em $L(0)$-semisimple $P^{(2)}(z)$-grading condition}: For
any $w_{(3)}\in W_3$, there exists a formal series $\sum_{n\in
\C}\lambda^{(2)}_{n}$ with
\[
\lambda^{(2)}_{n}\in \coprod_{\beta\in \tilde{A}}
((W_{1}\otimes W_{2})^{*})_{(n)}^{(\beta)}
\]
for $n\in \C$ such that for $w_{(1)}\in W_{1}$ and $w_{(2)}\in W_{2}$,
\[
\{(n, 0)\in \C\times \N\;|\;\lambda_{n}^{(1)}(w_{(1)}\otimes w_{(2)})\ne 0\}
\]
is a unique expansion set and 
\[
\sum_{n\in \C}\lambda^{(2)}_{n}(w_{(1)}\otimes w_{(2)})
\]
is absolutely convergent to $\mu^{(2)}_{\lambda,
w_{(3)}}(w_{(1)}\otimes w_{(2)})$.  (Note that such an element
$\lambda$ also satisfies the $P^{(2)}(z)$-grading condition above with
the same elements $\lambda^{(2)}_{n}$.)

(b) For any $w_{(3)}\in W_3$, consider the space $W^{(2)}_{\lambda,
w_{(3)}}$ as above, which in this case is in fact the smallest doubly
graded (or equivalently, $\tilde A$-graded) subspace of $((W_1\otimes
W_2)^{*})_{(\C)}^{(\tilde{A})}$ containing all the terms
$\lambda^{(2)}_{n}$ in the formal series in (a) and stable under the
component operators $\tau_{P(z)}(v\otimes t^{m})$ of the operators
$Y'_{P(z)}(v, x)$ for $v\in V$, $m\in {\mathbb Z}$, and under the
operators $L'_{P(z)}(-1)$, $L'_{P(z)}(0)$ and $L'_{P(z)}(1)$. Then
$W^{(2)}_{\lambda, w_{(3)}}$ has the properties
\begin{eqnarray*}
&\dim(W^{(2)}_{\lambda, w_{(3)}})^{(\beta)}_{(n)}<\infty,&\\
&(W^{(2)}_{\lambda, w_{(3)}})^{(\beta)}_{(n+k)}=0\;\;\mbox{ for
}\;k\in {\mathbb Z} \;\mbox{ sufficiently negative}&
\end{eqnarray*}
for any $n\in {\mathbb C}$ and $\beta\in \tilde A$, where the
subscripts denote the ${\mathbb C}$-grading by
$L'_{P(z)}(0)$-eigenvalues and the superscripts denote the $\tilde
A$-grading.
\end{description}

\begin{rema}\label{part-a}
{\rm Part (a) of each of these conditions can be reformulated using
the language of weak absolute convergence, in the sense of Remark
\ref{weakly-abs-conv}.  For example, the $P^{(1)}(z)$-grading
condition includes the assertion that for any $w_{(1)}\in W_1$,
$\mu^{(1)}_{\lambda, w_{(1)}}$ is the sum of a weakly absolutely
convergent series $\sum_{n\in \C}\lambda^{(1)}_{n}$ with
$\lambda^{(1)}_{n}\in \coprod_{\beta\in \tilde{A}} ((W_2\otimes
W_3)^{*})_{[n]}^{(\beta)}$.}
\end{rema}

In the rest of this section, we shall focus on the case that the
convergence condition for intertwining maps in $\mathcal{C}$ holds and
that generalized $V$-modules that we start with are objects of
$\mathcal{C}$. Recall the categories $\mathcal{M}_{sg}$ and
$\mathcal{GM}_{sg}$ from Notation $\ref{MGM}$, and recall Assumptions
\ref{assum}, \ref{assum-c} and \ref{assum-exp-set} on the category
$\mathcal{C}$.

Let $W_{1}$, $W_{2}$, $W_{3}$ and $W_{4}$ be generalized $V$-modules.
Given an $\tilde{A}$-compatible map $F: W_{1}\otimes W_{2}\otimes
W_{3}\to \overline{W}_{4}$ as in Remark
\ref{Atildecompatcorrespondence}, there is a canonical
$\tilde{A}$-compatible map $G$ from $W'_{4}$ to $(W_{1}\otimes
W_{2}\otimes W_{3})^{*}$ corresponding to $F$ under the indicated
canonical isomorphism between the spaces of such maps. We shall denote
the map $G$ corresponding to $F$ by $F'$.  Assume the convergence
condition for intertwining maps in $\mathcal{C}$ and that all
generalized $V$-modules considered are objects of $\mathcal{C}$.  Let
$I_{1}$, $I_{2}$, $I^1$ and $I^2$ be $P(z_1)$-, $P(z_2)$-, $P(z_2)$-
and $P(z_0)$-intertwining maps of types ${W_4}\choose {W_1M_1}$,
${M_1}\choose {W_2W_3}$, ${W_4}\choose {M_2W_3}$ and ${M_2}\choose
{W_1W_2}$, respectively. Then the maps $(I_1\circ (1_{W_1}\otimes
I_2))': W_{4}'\to (W_{1}\otimes W_{2}\otimes W_{3})^{*}$ and
$(I^1\circ (I^2\otimes 1_{W_3}))': W_{4}'\to (W_{1}\otimes
W_{2}\otimes W_{3})^{*}$ are well defined.

\begin{propo}\label{9.7}
Assume the convergence condition for intertwining maps in
$\mathcal{C}$ and that all generalized $V$-modules considered are
objects of $\mathcal{C}$.  Let $W_{1}$, $W_{2}$,
$W_{3}$, $W_{4}$, $M_{1}$ and $M_{2}$ be objects of $\mathcal{C}$ 
and let $I_{1}$, $I_{2}$, $I^1$ and $I^2$ be
$P(z_1)$-, $P(z_2)$-, $P(z_2)$- and $P(z_0)$-intertwining maps of
types ${W_4}\choose {W_1M_1}$, ${M_1}\choose {W_2W_3}$, ${W_4}\choose
{M_2W_3}$ and ${M_2}\choose {W_1W_2}$, respectively.  Let $w'_{(4)}\in
W'_4$.  If $|z_1|>|z_2|>0$, then $(I_1\circ (1_{W_1}\otimes
I_2))'(w'_{(4)})$ satisfies the $P^{(1)}(z_2)$-local grading
restriction condition (or the $L(0)$-semisimple $P^{(1)}(z_2)$-local
grading restriction condition when $\mathcal{C}$ is in
$\mathcal{M}_{sg}$), and if $|z_2|>|z_0|>0$, then $(I^1\circ
(I^2\otimes 1_{W_3}))'(w'_{(4)})$ satisfies the $P^{(2)}(z_0)$-local
grading restriction condition (or the $L(0)$-semisimple
$P^{(2)}(z_2)$-local grading restriction condition when $\mathcal{C}$
is in $\mathcal{M}_{sg}$).  Moreover, $W^{(1)}_{(I_1\circ
(1_{W_1}\otimes I_2))'(w'_{(4)}), w_{(1)}}$ and $W^{(1)}_{(I^1\circ
(I^2\otimes 1_{W_3}))'(w'_{(4)}), w_{(3)}}$, equipped with the vertex
operator maps given by $Y'_{P(z_{2})}$ and $Y'_{P(z_{0})}$,
respectively, and the operators $L'_{P(z_{2})}(j)$ and
$L'_{P(z_{0})}(j)$, respectively, for $j=-1, 0, 1$, are objects of
$\mathcal{C}$; in particular, they are generalized $V$-modules.  In
particular, for any $w_{(1)}\in W_{1}$ and any $w_{(3)}\in W_{3}$, if
we let $\sum_{n\in \C}\lambda_{n}^{(1)}$ be weakly absolutely
convergent to $\mu^{(1)}_{(I_1\circ (1_{W_1}\otimes I_2))'(w'_{(4)}),
w_{(1)}}$ as given by the $P^{(1)}(z_{2})$-grading condition and let
$\sum_{n\in \C}\lambda_{n}^{(2)}$ be weakly absolutely convergent to
$\mu^{(2)}_{(I^1\circ (I^2\otimes 1_{W_3}))'(w'_{(4)}), w_{(3)}}$ as
given by the $P^{(2)}(z_{0})$-grading condition, then for any $n\in
\C$, the generalized $V$-submodule $W_{\lambda_{n}^{(1)}}$ (recall the
notation used in the $P(z)$-local grading restriction condition in
Section 5) of $W^{(1)}_{(I_1\circ (1_{W_1}\otimes I_2))'(w'_{(4)}),
w_{(1)}}$ and the generalized $V$-submodule $W_{\lambda_{n}^{(2)}}$ of
$W^{(1)}_{(I^1\circ (I^2\otimes 1_{W_3}))'(w'_{(4)}), w_{(3)}}$ are
objects of $\mathcal{C}$.
\end{propo}
\pf 
Let $w_{(1)}\in W_{1}$. 
For $n\in \C$, let $m_{(1), n}'\in M_{1}^{*}$ be defined by 
\[
m_{(1), n}'(m_{(1)})=\langle w'_{(4)}, I_1(w_{(1)}\otimes  
\pi_{n}(m_{(1)}))\rangle
\]
for $m_{(1)}\in M_{1}$. Since by definition,
for $m_{(1)}\in (M_{1})_{[m]}$,
$m_{(1), n}'(m_{(1)})=0$ when  $m\ne n$, 
we see that $m_{(1), n}'\in (M'_{1})_{[n]}$. 
Since $w'_{(4)}$ is a finite sum of $\tilde{A}$-homogeneous 
elements and $I_{1}$ presevers the $\tilde{A}$-gradings, 
$m_{(1), n}'$ is also a finite sum of $\tilde{A}$-homogeneous 
elements. 

Let $\lambda^{(1)}_{n}=m_{(1), n}'\circ I_{2}
\in (W_{2}\otimes W_{3})^{*}$. 
{}From Notation \ref{scriptN}, we have
$\lambda^{(1)}_{n}= I_{2}'(m_{(1), n}')$. Since $m_{(1), n}'\in (M'_{1})_{[n]}$,
by Proposition \ref{im:abc}(b), $\lambda^{(1)}_{n}= I_{2}'(m_{(1), n}')
\in ((W_{2}\otimes W_{3})^{*})_{[n]}$ for $n\in \C$. 
In addition, since $I_{2}'$ preserves the $\tilde{A}$-gradings and 
$m_{(1), n}'$ is a finite sum of $\tilde{A}$-homogeneous 
elements, 
$\lambda^{(1)}_{n}$ is also a finite sum of $\tilde{A}$-homogeneous elements. 
By Assumption \ref{assum-exp-set}, the set 
$\{(n, i)\in \C\times \N\;|\;(M'_{1})_{(n)}\ne 0,\;
(L(0)-n)^{i}(M'_{1})_{(n)}\ne 0\}$ is a unique expansion set.
In particular, its  subsets 
$\{(n, i)\in \C\times \N\;|\; 
((L(0)-n)^{i} \lambda^{(1)}_{n})(w_{(2)}\otimes w_{(3)})\ne 0\}$
for $w_{(2)}\in W_{2}$ and $w_{(3)}\in W_{3}$ are unique expansion sets.

When $|z_1|>|z_2|>0$, the product of $I_1$ and $I_2$ exists.
For any $w'_{(4)}\in W'_4$ and $w_{(1)}\in W_1$, we have
\begin{eqnarray*}
\mu^{(1)}_{(I_1\circ (1_{W_1}\otimes I_2))'(w'_{(4)}),
w_{(1)}}(w_{(2)}\otimes w_{(3)})
&=&\langle w'_{(4)}, I_1(w_{(1)}\otimes I_2(w_{(2)}\otimes w_{(3)}))
\rangle\nno\\
&=&\sum_{n\in \C}\langle w'_{(4)}, 
I_1(w_{(1)}\otimes \pi_{n}(I_2(w_{(2)}\otimes w_{(3)})))
\rangle\nno\\
&=&\sum_{n\in \C}m_{(1), n}'(I_2(w_{(2)}\otimes w_{(3)}))\nno\\
&=&\sum_{n\in \C}\lambda^{(1)}_{n}(w_{(2)}\otimes w_{(3)})
\end{eqnarray*}
for $w_{(2)}\in W_{2}$ and $w_{(3)}\in W_{3}$, proving 
that $(I_1\circ (1_{W_1}\otimes I_2))'(w'_{(4)})$ satisfies the
$P^{(1)}(z_2)$-grading condition. 

By Proposition \ref{im:abc}(b), the map $I_{2}'$ preserves
generalized weights. By definition, we know that $I_{2}'$ also
preserves the $\tilde{A}$-grading. Since $M'_{1}$ is an object of
$\mathcal{C}$, the image under $I_{2}'$ of the generalized
$V$-submodule of $M'_{1}$ generated by $m_{(1), n}'$ for $n\in \C$
satisfies the two grading restriction conditions (\ref{lgrc1}) and
(\ref{lgrc2}).  Since $W^{(1)}_{(I_1\circ (1_{W_1}\otimes
I_2))'(w'_{(4)}), w_{(1)}}$ is the image of this generalized
$V$-submodule under $I'_{2}$, Part (b) of the $P^{(1)}(z_2)$-local
grading restriction condition holds.

Let $w_{(3)}\in W_{3}$. Analogously,
for $n\in \C$, let $m_{(2), n}'\in M_{2}^{*}$ be defined by 
\[
m_{(2), n}'(m_{(2)})=\langle w'_{(4)}, I^1(\pi_{n}(m_{(2)})
\otimes w_{(3)})\rangle
\]
for $m_{(2)}\in M_{2}$. Then $m_{(2), n}'\in (M'_{2})_{[n]}$ and 
is a finite sum of $\tilde{A}$-homogeneous 
elements. 
Let $\lambda^{(2)}_{n}=m_{(2), n}'\circ I^{2}
\in (W_{1}\otimes W_{2})^{*}$. Then 
$\lambda^{(2)}_{n}= (I^{2})'(m_{(2), n}')
\in ((W_{2}\otimes W_{3})^{*})_{[n]}$ for $n\in \C$ and is 
a finite sum of $\tilde{A}$-homogeneous elements. 
By Assumption \ref{assum-exp-set}, the set 
$\{(n, i)\in \C\times \N\;|\;(M'_{2})_{(n)}\ne 0,\;
(L(0)-n)^{i}(M'_{2})_{(n)}\ne 0\}$ is a unique expansion set.
In particular, its  subsets 
$\{(n, i)\in \C\times \N\;|\; 
(L(0)-n)^{i} \lambda^{(2)}_{n}(w_{(1)}\otimes w_{(2)})\ne 0\}$ 
are  unique expansion sets.

When $|z_2|>|z_0|>0$, the iterate of $I^1$ and $I^2$ exists.
For any $w'_{(4)}\in W'_4$ and $w_{(1)}\in W_1$, we have
\begin{eqnarray*}
\mu^{(2)}_{(I^1\circ (I^{2}\otimes 1_{W_2}))'(w'_{(4)}),
w_{(3)}}(w_{(1)}\otimes w_{(2)})
&=&\langle w'_{(4)}, I^1(I^2(w_{(1)}\otimes w_{(2)})\otimes w_{(3)})
\rangle\nno\\
&=&\sum_{n\in \C}\langle w'_{(4)}, 
I^1(\pi_{n}(I^2(w_{(1)}\otimes w_{(2)}))\otimes w_{(3)})
\rangle\nno\\
&=&\sum_{n\in \C}m_{(2), n}'(I^2(w_{(1)}\otimes w_{(2)}))\nno\\
&=&\sum_{n\in \C}\lambda^{(2)}_{n}(w_{(1)}\otimes w_{(2)})
\end{eqnarray*}
for $w_{(1)}\in W_{1}$ and $w_{(2)}\in W_{2}$, proving 
that $(I^1\circ  (I_2)\otimes 1_{W_3}))'(w'_{(4)})$ satisfies the
$P^{(2)}(z_0)$-grading condition. 

By Proposition \ref{im:abc}(b), the map $(I^{2})'$ preserves
generalized weights and by definition, we know that $(I_{2})'$ also
preserves the $\tilde{A}$-grading. Since $M'_{2}$ is an object of
$\mathcal{C}$, the image under $(I^{2})'$ of the generalized
$V$-submodule of $M_{2}'$ generated by $m_{(2), n}'$ for $n\in \C$
satisfies the two grading restriction conditions (\ref{lgrc1}) and
(\ref{lgrc2}).  Since $W^{(2)}_{(I_1\circ (1_{W_1}\otimes
I_2))'(w'_{(4)}), w_{(1)}}$ is the image of this generalized
$V$-submodule under $(I^{2})'$, Part (b) of the $P^{(2)}(z_0)$-local
grading restriction condition holds.

By Proposition \ref{pz},
$I_{2}'$ intertwinines the actions of vertex operators
and elements of $\mathfrak{sl}(2)$. Since
$\lambda^{(1)}_{n}= I_{2}'(m_{(1), n}')$, $\lambda^{(1)}_{n}$ is in fact 
in the image $I_{2}'(M_{1}')$ 
of the homomorphism $I'_{2}$ of generalized $V$-modules from 
$M_{1}'$ to $I_{2}'(M_{1}')$. Since $M_{1}$ is an object of $\mathcal{C}$,
$M_{1}'$ is also an object of $\mathcal{C}$. Thus $I_{2}'(M_{1}')$
is an object of $\mathcal{C}$. So $W^{(1)}_{(I_1\circ (1_{W_1}\otimes
I_2))'(w'_{(4)}), w_{(1)}}$
is an object of $\mathcal{C}$. 

Similarly, we can show that $W^{(1)}_{(I^1\circ
(I^2\otimes 1_{W_3}))'(w'_{(4)}), w_{(3)}}\subset 
(I^{2})'(M_{2}')$   is an object of $\mathcal{C}$. 
\epf

\begin{theo}\label{9.7-1}
Assume the convergence condition for intertwining maps in
$\mathcal{C}$ and $|z_1|>|z_2|>|z_{0}|>0$.  Let $W_{1}$, $W_{2}$,
$W_{3}$, $W_{4}$, $M_{1}$ and $M_{2}$ be objects of $\mathcal{C}$ 
and let $I_{1}$, $I_{2}$,
$I^1$ and $I^2$ be $P(z_1)$-, $P(z_2)$-, $P(z_2)$- and
$P(z_0)$-intertwining maps of types ${W_4}\choose {W_1M_1}$,
${M_1}\choose {W_2W_3}$, ${W_4}\choose {M_2W_3}$ and ${M_2}\choose
{W_1W_2}$, respectively.  Let $w'_{(4)}\in W'_4$. If $(I_1\circ
(1_{W_1}\otimes I_2))'(w'_{(4)})$ satisfies the $P^{(2)}(z_0)$-grading 
condition such that in particular, for any $w_{(3)}\in
W_{3}$, $\sum_{n\in \C}\lambda_{n}^{(2)}$ is weakly absolutely
convergent to $\mu^{(2)}_{(I_1\circ (1_{W_1}\otimes I_2))'(w'_{(4)}),
w_{(3)}}$ and, in addition,  $\lambda_{n}^{(2)}$, $n\in \C$,
satisfy the $P(z_{0})$-lower truncation 
condition (Part (a) of the $P(z_{0})$-compatibility condition
in Section 5), 
then $W^{(2)}_{(I_1\circ (1_{W_1}\otimes I_2))'(w'_{(4)}),
w_{(3)}}$, equipped with the vertex operator map given by
$Y'_{P(z_{0})}$ and the operators $L'_{P(z_{0})}(j)$ for $j=-1, 0, 1$,
is a lower-truncated generalized $V$-module. In particular,
if $(I_1\circ
(1_{W_1}\otimes I_2))'(w'_{(4)})$ satisfies the $P(z_{0})$-local 
grading restriction condition (or the $L(0)$-semisimple
$P^{(2)}(z_0)$-local grading restriction condition when $\mathcal{C}$
is in $\mathcal{M}_{sg}$), then 
$W^{(2)}_{(I_1\circ (1_{W_1}\otimes I_2))'(w'_{(4)}),
w_{(3)}}$ is
an object of $\mathcal{GM}_{sg}$ (or $\mathcal{M}_{sg}$ when
$\mathcal{C}$ is in $\mathcal{M}_{sg}$). Analogously, if $(I^1\circ
(I^2\otimes 1_{W_3}))'(w'_{(4)})$ satisfies the $P^{(1)}(z_2)$-grading 
condition such that in particular, for any $w_{(1)}\in
W_{1}$, $\sum_{n\in \C}\lambda_{n}^{(1)}$ is weakly absolutely
convergent to $\mu^{(1)}_{(I^1\circ (I^2\otimes 1_{W_3}))'(w'_{(4)}),
w_{(3)}}$ and, in addition,  $\lambda_{n}^{(1)}$, $n\in \C$,
satisfy the $P(z_{2})$-lower truncation 
condition, then $W^{(1)}_{(I^1\circ (I^2\otimes 1_{W_3}))'(w'_{(4)}),
w_{(1)}}$, equipped with the vertex operator map given by
$Y'_{P(z_{2})}$ and the operators $L'_{P(z_{2})}(j)$ for $j=-1, 0, 1$,
is a generalized $V$-module. In particular, 
if $(I^1\circ
(I^2\otimes 1_{W_3}))'(w'_{(4)})$ satisfies the $P^{(1)}(z_2)$-local
grading restriction condition (or the $L(0)$-semisimple
$P^{(1)}(z_2)$-local grading restriction condition when $\mathcal{C}$
is in $\mathcal{M}_{sg}$) is an object of $\mathcal{GM}_{sg}$ 
(or $\mathcal{M}_{sg}$ when $\mathcal{C}$ is
in $\mathcal{M}_{sg}$).
\end{theo}
\pf
By Theorem \ref{wk-mod}, to prove $W^{(2)}_{(I_1\circ (1_{W_1}\otimes
I_2))'(w'_{(4)}), w_{(3)}}$ is a generalized $V$-module, 
we need only prove that its elements satisfies the
$P(z_{0})$-compatibility condition. This is equivalent to proving that
$\lambda_{n}^{(2)}$, $n\in \C$, satisfy the $P(z_{0})$-compatibility
condition. By assumption, the $P(z_{0})$-lower truncation condition is
clearly satisfied by $\lambda_{n}^{(2)}$, $n\in \C$. We need only prove
that they satisfy Part (b) of the $P(z_{0})$-compatibility
condition. The proof below is a generalization of the proof of (14.51)
in \cite{tensor4}.

Let $\mathcal{Y}_{1}=\mathcal{Y}_{I_{1}, 0}$ and 
$\mathcal{Y}_{2}=\mathcal{Y}_{I_{2}, 0}$.
For $z\in \C^{\times}$, Let $I_{1}^{z}$ and $I_{2}^{z}$ 
be the $P(z_{0}+zz_{2})$- and 
$P(zz_{2})$-intertwining maps $I_{\mathcal{Y}_{1}, -1}$ and 
$I_{\mathcal{Y}_{2}, -1}$.
If $0<|z|< \frac{|z_{0}|}{2|z_{2}|}$, we have 
$|z_{0}+zz_{2}|>|zz_{2}|>0$. Thus the product 
$I_{1, z}\circ (1_{W_{2}}\otimes I_{2, z})$ exists. We now assume that 
$z$ is in the region  $0<|z|< \frac{|z_{0}|}{2|z_{2}|}$.

By assumption, for $w_{(1)}\in W_{1}$ and $w_{(2)}\in W_{2}$,
the series $\sum_{n\in {\Bbb C}}\lambda^{(2)}_{n}(w_{(1)}\otimes w_{(2)})$ 
converges  absolutely to 
$\mu^{(2)}_{(I_1\circ (1_{W_1}\otimes
I_2))'(w'_{(4)}), w_{(3)}}(w_{(1)}\otimes w_{(2)})$. 
Note that $\lambda^{(2)}_{n}$ for $n\in \C$ depend on $w'_{(4)}\in W'_{4}$ 
and $w_{(3)}\in W_{3}$. Since in this proof, we need to use this dependence,
we shall write $\lambda^{(2)}_{n}$ as 
$\lambda^{(2)}_{n}(w'_{(4)}, w_{(3)})$.

The case $j=0$, $z=z_{0}$ and 
$\lambda=\mu^{(2)}_{(I_1\circ (1_{W_1}\otimes
I_2))'(w'_{(4)}), w_{(3)}}$
of (\ref{LP'(j)})
gives 
\begin{eqnarray}\label{9.7-1-0}
\lefteqn{(L'_{P(z_{0})}(0)\mu^{(2)}_{(I_1\circ (1_{W_1}\otimes
I_2))'(w'_{(4)}), w_{(3)}})(w_{(1)}\otimes w_{(2)})}\nno\\
&&=\mu^{(2)}_{(I_1\circ (1_{W_1}\otimes
I_2))'(w'_{(4)}), w_{(3)}}(w_{(1)}\otimes L(0)w_{(2)}
+(L(0)+z_{0}L(-1))w_{(1)}\otimes w_{(2)}).
\end{eqnarray}
Let $x$ be a formal variable. Then we have
\begin{eqnarray}\label{14.43}
\lefteqn{((1-x)^{-L'_{P(z_{0})}(0)}
\mu^{(2)}_{(I_1\circ (1_{W_1}\otimes
I_2))'(w'_{(4)}), w_{(3)}})(w_{(1)}\otimes w_{(2)})}\nno\\
&&=(e^{-\log(1-x)L'_{P(z_{0})}(0)}
\mu^{(2)}_{(I_1\circ (1_{W_1}\otimes
I_2))'(w'_{(4)}), w_{(3)}})(w_{(1)}\otimes w_{(2)})\nno\\
&&=\mu^{(2)}_{(I_1\circ (1_{W_1}\otimes
I_2))'(w'_{(4)}), w_{(3)}}(e^{\log(1-x)(-z_{0}L(-1)-L(0))}w_{(1)}\otimes 
e^{-\log(1-x)L(0)}w_{(2)})\nno\\
&&=\mu^{(2)}_{(I_1\circ (1_{W_1}\otimes
I_2))'(w'_{(4)}), w_{(3)}}((1-x)^{-z_{0}L(-1)-L(0)}w_{(1)}\otimes 
(1-x)^{-L(0)}w_{(2)})\nno\\
&&=((I_1\circ (1_{W_1}\otimes
I_2))'(w'_{(4)}))((1-x)^{-z_{0}L(-1)-L(0)}
w_{(1)}\otimes 
(1-x)^{-L(0)}w_{(2)}\otimes w_{(3)})\nno\\
&&=\langle w'_{(4)}, {\cal Y}_{1}((1-x)^{-z_{0}L(-1)-L(0)}
w_{(1)}, x_{1}){\cal Y}_{2}(
(1-x)^{-L(0)}w_{(2)}, x_{2})
w_{(3)}\rangle_{W_{4}}\lbar_{x_{1}=z_{1}, x_{2}=z_{2}}.\nno\\
&&
\end{eqnarray}
Using (\ref{log:L(j)b2}), we see that
the right-hand side of (\ref{14.43}) is equal to
\begin{eqnarray}\label{14.44}
\lefteqn{\langle w'_{(4)}, (1-x)^{-(L(0)-z_{2}L(-1))}{\cal Y}_{1}(
w_{(1)}, x_{1})\cdot}\nno\\
&&\hspace{3em}\cdot{\cal Y}_{2}(
w_{(2)}, x_{2})(1-x)^{L(0)-z_{2}L(-1)}
w_{(3)}\rangle_{W_{4}}\lbar_{x_{1}=z_{1}, x_{2}=z_{2}}.
\end{eqnarray}
By Lemma 9.3 in \cite{tensor2} which gives the formula
$$(1-x)^{L(0)-z_{2}L(-1)}=e^{z_{2}xL(-1)}(1-x)^{L(0)},$$ 
(\ref{14.44}) is equal to
\begin{eqnarray}\label{14.45}
\lefteqn{\langle w'_{(4)}, (1-x)^{-L(0)}e^{-z_{2}xL(-1)}{\cal Y}_{1}(
w_{(1)}, x_{1})\cdot}\nno\\
&&\hspace{3em}\cdot {\cal Y}_{2}(
w_{(2)}, x_{2})e^{z_{2}xL(-1)}(1-x)^{L(0)}
w_{(3)}\rangle_{W_{4}}\lbar_{x_{1}=z_{1}, x_{2}=z_{2}}\nno\\
&&=\langle (1-x)^{-L(0)}w'_{(4)}, {\cal Y}_{1}(
w_{(1)}, x_{1}-z_{2}x)\cdot\nno\\
&&\hspace{3em}\cdot{\cal Y}_{2}(
w_{(2)}, x_{2}-z_{2}x)(1-x)^{L(0)}
w_{(3)}\rangle_{W_{4}}\lbar_{x_{1}=z_{1}, x_{2}=z_{2}}\nno\\
&&=\langle (1-x)^{-L(0)}w'_{(4)}, {\cal Y}_{1}(
w_{(1)}, x_{1})\cdot\nno\\
&&\hspace{3em}\cdot{\cal Y}_{2}(
w_{(2)}, x_{2})(1-x)^{L(0)}
w_{(3)}\rangle_{W_{4}}\lbar_{x_{1}=z_{1}-z_{2}x, x_{2}=z_{2}-z_{2}x}\nno\\
&&=\langle (1-x)^{-L(0)}w'_{(4)}, {\cal Y}_{1}(
w_{(1)}, x_{1})\cdot\nno\\
&&\hspace{3em}\cdot{\cal Y}_{2}(
w_{(2)}, x_{2})(1-x)^{L(0)}
w_{(3)}\rangle_{W_{4}}\lbar_{x_{1}=z_{1}-z_{2}+z_{2}(1-x), x_{2}=z_{2}(1-x)}.
\end{eqnarray}

The calculations above show that the left-hand side of (\ref{14.43}) 
is equal
to the right-hand side of (\ref{14.45}).
Since $0<|z|<\frac{|z_{0}|}{2|z_{2}|}$,
the series obtained by writing the right-hand side of (\ref{14.45})
as a series in $\log(1-x)$, substituting $y$ for $\log(1-x)$,
and then substituting 
$\log z$ for $y$ is equal to the absolutely convergent series
\begin{equation}\label{9.7-1-1}
\langle e^{-(\log z)L'(0)}w'_{(4)}, {\cal Y}_{1}(
w_{(1)}, x_{1}){\cal Y}_{2}(
w_{(2)}, x_{2})e^{(\log z)L(0)}
w_{(3)}\rangle_{W_{4}}\lbar_{x_{1}=z_{0}+zz_{2}, x_{2}=zz_{2}}.
\end{equation}
So the series obtained by writing the right-hand side of (\ref{14.43})
as a series in $\log(1-x)$, substituting $y$ for $\log(1-x)$,
and then substituting 
$\log z$ for $y$ is also absolutely convergent and, 
its sum is equal to (\ref{9.7-1-1}).
We now calculate 
this series.

Since 
\[
\sum_{n\in \C}(\lambda^{(2)}_{n}(w'_{(4)}, w_{(3)}))(w_{(1)}\otimes 
w_{(2)})
\]
is absolutely convergent to 
\[
\mu^{(2)}_{(I_1\circ (1_{W_1}\otimes
I_2))'(w'_{(4)}), w_{(3)}}(w_{(1)}\otimes w_{(2)}),
\]
using (\ref{9.7-1-0}), we have 
\begin{eqnarray}\label{9.7-1-2}
\lefteqn{((1-x)^{-L'_{P(z_{0})}(0)}
\mu^{(2)}_{(I_1\circ (1_{W_1}\otimes
I_2))'(w'_{(4)}), w_{(3)}})(w_{(1)}\otimes w_{(2)})}\nno\\
&&=(e^{-\log (1-x) L'_{P(z_{0})}(0)}
\mu^{(2)}_{(I_1\circ (1_{W_1}\otimes
I_2))'(w'_{(4)}), w_{(3)}})(w_{(1)}\otimes w_{(2)})\nno\\
&&=\mu^{(2)}_{(I_1\circ (1_{W_1}\otimes
I_2))'(w'_{(4)}), w_{(3)}}(e^{\log (1-x)(-z_{0}L(-1)-L(0))}w_{(1)}\otimes 
e^{-\log (1-x)L(0)}w_{(2)})\nno\\
&&=\mu^{(2)}_{(I_1\circ (1_{W_1}\otimes
I_2))'(w'_{(4)}), w_{(3)}}((1-x)^{-z_{0}L(-1)-L(0)}w_{(1)}\otimes 
(1-x)^{-L(0)}w_{(2)})\nno\\
&&=\sum_{n\in \C}(\lambda^{(2)}_{n}(w'_{(4)}, w_{(3)}))
(e^{\log (1-x)(-z_{0}L(-1)-L(0))}w_{(1)}\otimes 
e^{-\log (1-x)L(0)}w_{(2)})\nno\\
&&=\sum_{n\in \C}(e^{-\log (1-x)L'_{P(z_{0})}(0)}\lambda^{(2)}_{n}(w'_{(4)}, w_{(3)}))
(w_{(1)}\otimes w_{(2)})\nno\\
&&=\sum_{n\in \C}((1-x)^{-L'_{P(z_{0})}(0)}\lambda^{(2)}_{n}(w'_{(4)}, w_{(3)}))
(w_{(1)}\otimes w_{(2)}).
\end{eqnarray}
Since $\lambda_{n}^{(2)}\in \coprod_{\beta\in \tilde{A}}(W_{1}\otimes 
W_{2})^{*})_{[n]}^{(\beta)}$, 
$L'_{P(z_{0})}(0)_{s}\lambda_{n}^{(2)}=n\lambda_{n}^{(2)}$
and there exist $K_{n}\in \Z_{+}$ such that 
$(L'_{P(z_{0})}(0)-L'_{P(z_{0})}(0)_{s})^{K_{n}}\lambda_{n}^{(2)}=0$. 
Now the calculations in 
Remark \ref{set:L(0)s} still work for the space 
$W^{(2)}_{(I_1\circ (1_{W_1}\otimes
I_2))'(w'_{(4)}), w_{(3)}}$. In particular, 
$L'_{P(z_{0})}(0)-L'_{P(z_{0})}(0)_{s}$ commutes with 
$L'_{P(z_{0})}(0)$ and therefore also with $L'_{P(z_{0})}(0)_{s}$. 
Using these facts, we obtain 
\begin{eqnarray}\label{9.7-1-3}
\lefteqn{(1-x)^{-L'_{P(z_{0})}(0)}\lambda^{(2)}_{n}(w'_{(4)}, w_{(3)})}\nno\\
&&=
e^{-\log(1-x)L'_{P(z_{0})}(0)}\lambda^{(2)}_{n}(w'_{(4)}, w_{(3)})\nno\\
&&=
e^{-\log (1-x)(L'_{P(z_{0})}(0)-L'_{P(z_{0})}(0)_{s})}
e^{-\log (1-x)L'_{P(z_{0})}(0)_{s}}\lambda^{(2)}_{n}(w'_{(4)}, w_{(3)})\nno\\
&&=e^{-n\log (1-x)}\sum_{i=0}^{K_{n}-1}\frac{(-\log (1-x))^{i}}{i!}
(L'_{P(z_{0})}(0)-L'_{P(z_{0})}(0)_{s})^{i}\lambda^{(2)}_{n}(w'_{(4)}, w_{(3)})\nno\\
&&=(1-x)^{-n}\sum_{i=0}^{K_{n}-1}\frac{(-\log (1-x))^{i}}{i!}
(L'_{P(z_{0})}(0)-n)^{i}\lambda^{(2)}_{n}(w'_{(4)}, w_{(3)}).
\end{eqnarray}
{}From  (\ref{9.7-1-2}) and (\ref{9.7-1-3}), we obtain
\begin{eqnarray}\label{9.7-1-4}
\lefteqn{\sum_{n\in \C}((1-x)^{-L'_{P(z_{0})}(0)}\lambda^{(2)}_{n}(w'_{(4)}, w_{(3)}))
(w_{(1)}\otimes w_{(2)})}\nno\\
&&=\sum_{n\in \C}(1-x)^{-n}\sum_{i=0}^{K_{n}-1}\frac{(-\log(1-x))^{i}}{i!}
((L'_{P(z_{0})}(0)-n)^{i}\lambda^{(2)}_{n}(w'_{(4)}, w_{(3)}))
(w_{(1)}\otimes w_{(2)}).\nno\\
&&
\end{eqnarray}
Writing the left-hand side of (\ref{9.7-1-4}) as a series 
in $\log(1-x)$, substituting $y$ for $\log(1-x)$ 
and then substituting 
$\log z$ for $y$, we obtain 
\begin{equation}\label{9.7-1-5}
\sum_{n\in \C}e^{-n\log z}\sum_{i=0}^{K_{n}-1}\frac{(-\log z)^{i}}{i!}
((L'_{P(z_{0})}(0)-n)^{i}\lambda^{(2)}_{n}(w'_{(4)}, w_{(3)}))
(w_{(1)}\otimes w_{(2)}).
\end{equation}
{}From the discussion above, we see that (\ref{9.7-1-5}) is absolutely 
convergent to (\ref{9.7-1-1}). 

But by definition, (\ref{9.7-1-1}) is equal to 
\[
(\mu^{(2)}_{(I^{z}_1\circ (1_{W_1}\otimes
I^{z}_2))'(e^{-(\log z)L'(0)}w'_{(4)}), e^{(\log z)L(0)}w_{(3)}})
(w_{(1)}\otimes w_{(2)}).
\]
So we obtain
\begin{eqnarray*}
\lefteqn{(\mu^{(2)}_{(I^{z}_1\circ (1_{W_1}\otimes
I^{z}_2))'(e^{-(\log z)L'(0)}w'_{(4)}), e^{(\log z)L(0)}w_{(3)}})
(w_{(1)}\otimes w_{(2)})}\nno\\
&&=\sum_{n\in \C}e^{-n\log z}\sum_{i=0}^{K_{n}-1}\frac{(-\log z)^{i}}{i!}
((L'_{P(z_{0})}(0)-n)^{i}\lambda^{(2)}_{n}(w'_{(4)}, w_{(3)}))
(w_{(1)}\otimes w_{(2)}),
\end{eqnarray*}
or equivalently,
\begin{eqnarray}\label{9.7-1-6}
\lefteqn{(\mu^{(2)}_{(I^{z}_1\circ (1_{W_1}\otimes
I^{z}_2))'(w'_{(4)}), w_{(3)}})
(w_{(1)}\otimes w_{(2)})}\nno\\
&&=\sum_{n\in \C}e^{-n\log z}\sum_{i=0}^{K_{n}-1}\frac{(-\log z)^{i}}{i!}
\cdot\nno\\
&&\quad\quad\quad\quad\cdot 
((L'_{P(z_{0})}(0)-n)^{i}\lambda^{(2)}_{n}(e^{(\log z)L'(0)}w'_{(4)}, 
e^{-(\log z)L(0)}w_{(3)}))
(w_{(1)}\otimes w_{(2)}).
\end{eqnarray}
When $|z_{0}+zz_{2}|>|zz_{2}|>|z_{0}|>0$, since 
$(I^{z}_1\circ (1_{W_1}\otimes
I^{z}_2))'(w'_{(4)})$ satisfies the $P(z_{1}, z_{2})$-compatibility
condition, by Lemma \ref{mulemma},
the coefficients of 
\[
x^{-1}_1 \delta\left(\frac{x^{-1}-z_{0})}{x_1}\right)
Y'_{P(z_{0})}(v, x)(\mu^{(2)}_{(I^{z}_1\circ (1_{W_1}\otimes
I^{z}_2))'(w'_{(4)}), w_{(3)}})
(w_{(1)}\otimes w_{(2)})
\]
in monomials of $x_{1}$ and $x$ are absolutely convergent and 
we have
\begin{eqnarray}\label{mu12-1}
\lefteqn{\biggl(\tau_{P(z_0)}\biggl( x^{-1}_1
\delta\bigg(\frac{x^{-1}-z_0}{x_1}\bigg) Y_{t}(v, x)\biggr)
\mu^{(2)}_{(I^{z}_1\circ (1_{W_1}\otimes
I^{z}_2))'(w'_{(4)}), w_{(3)}}\biggr)(w_{(1)} \otimes w_{(2)})}\nno\\
&&=x^{-1}_1 \delta\bigg(\frac{x^{-1}-z_0}{x_1}\bigg)
\biggl(Y'_{P(z_0)}(v, x) \mu^{(2)}_{(I^{z}_1\circ (1_{W_1}\otimes
I^{z}_2))'(w'_{(4)}), w_{(3)}}\biggr)
(w_{(1)} \otimes w_{(2)}).
\end{eqnarray}
{}From the definitions of $I_{1}^{z}$, $I_{2}^{z}$ and 
$\mu^{(2)}_{(I^{z}_1\circ (1_{W_1}\otimes
I^{z}_2))'(w'_{(4)}), w_{(3)}}$, we see that the coefficients of 
both sides of (\ref{mu12-1}) in monomials of $x_{1}$ and $x$
are absolutely convergent series in $z$. From (\ref{9.7-1-6})
and (\ref{mu12-1}), we obtain 
\begin{eqnarray*}
\lefteqn{\Biggl(\tau_{P(z_0)}\biggl( x^{-1}_1
\delta\biggl(\frac{x^{-1}-z_0}{x_1}\bigg) Y_{t}(v, x)\biggr)
\sum_{n\in \C}e^{-n\log z}\sum_{i=0}^{K_{n}-1}\frac{(-\log z)^{i}}{i!}
\cdot}\nno\\
&&\quad\quad\quad\quad\cdot 
((L'_{P(z_{0})}(0)-n)^{i}\lambda^{(2)}_{n}(e^{(\log z)L'(0)}w'_{(4)}, 
e^{-(\log z)L(0)}w_{(3)}))\Biggr)
(w_{(1)}\otimes w_{(2)})\nno\\
&&=x^{-1}_1 \delta\bigg(\frac{x^{-1}-z_0}{x_1}\bigg)
\Biggl(Y'_{P(z_0)}(v, x) \sum_{n\in \C}e^{-n\log z}\sum_{i=0}^{K_{n}-1}\frac{(-\log z)^{i}}{i!}
\cdot\nno\\
&&\quad\quad\quad\quad\cdot 
((L'_{P(z_{0})}(0)-n)^{i}\lambda^{(2)}_{n}(e^{(\log z)L'(0)}w'_{(4)}, 
e^{-(\log z)L(0)}w_{(3)}))\Biggr)
(w_{(1)}\otimes w_{(2)}),
\end{eqnarray*}
or equivalent, by replacing $e^{(\log z)L'(0)}w'_{(4)}$
and $e^{-(\log z)L(0)}w_{(3)}$ with $w'_{(4)}$ and $w_{(3)}$,
respectively,
\begin{eqnarray}\label{mu12-2}
\lefteqn{\Biggl(\tau_{P(z_0)}\biggl( x^{-1}_1
\delta\biggl(\frac{x^{-1}-z_0}{x_1}\bigg) Y_{t}(v, x)\biggr)\cdot}\nno\\
&&\quad\quad\quad\quad\cdot 
\sum_{n\in \C}e^{-n\log z}\sum_{i=0}^{K_{n}-1}\frac{(-\log z)^{i}}{i!}
((L'_{P(z_{0})}(0)-n)^{i}\lambda^{(2)}_{n}(w'_{(4)}, 
w_{(3)}))\Biggr)
(w_{(1)}\otimes w_{(2)})\nno\\
&&=x^{-1}_1 \delta\bigg(\frac{x^{-1}-z_0}{x_1}\bigg)
\Biggl(Y'_{P(z_0)}(v, x) \cdot\nno\\
&&\quad\quad\quad\quad\cdot 
\sum_{n\in \C}e^{-n\log z}\sum_{i=0}^{K_{n}-1}\frac{(-\log z)^{i}}{i!}
((L'_{P(z_{0})}(0)-n)^{i}\lambda^{(2)}_{n}(w'_{(4)}, 
w_{(3)}))\Biggr)
(w_{(1)}\otimes w_{(2)}).\nno\\&&
\end{eqnarray}
Note that (\ref{mu12-2})
holds in the region
$|z_{0}+zz_{2}|>|zz_{2}|>|z_{0}|>0$ and the powers of 
$z$ and $\log z$ in the nonzero terms of both sides of (\ref{mu12-2}) 
form unique expansion sets, by the $P(z_{0})$-grading condition.
Thus, by the definition of unique expansion sets,
we see that the expansion coefficients of the left- and the right-hand
sides of (\ref{mu12-2}) are equal. In particular, 
\begin{eqnarray*}
\lefteqn{\Biggl(\tau_{P(z_0)}\biggl( x^{-1}_1
\delta\biggl(\frac{x^{-1}-z_0}{x_1}\bigg) Y_{t}(v, x)\biggr)
(\lambda^{(2)}_{n}(w'_{(4)}, 
w_{(3)}))\Biggr)
(w_{(1)}\otimes w_{(2)})}\nno\\
&&=x^{-1}_1 \delta\Bigg(\frac{x^{-1}-z_0}{x_1}\bigg)
(Y'_{P(z_0)}(v, x)
(\lambda^{(2)}_{n}(w'_{(4)}, 
w_{(3)})))
(w_{(1)}\otimes w_{(2)})
\end{eqnarray*}
which is Part (b) of the compatibility condition 
for $\lambda^{(2)}_{n}(w'_{(4)}, 
w_{(3)})$. 

Analogously, we can prove that under the stated conditions,
$W^{(1)}_{(I^1\circ
(I^2\otimes 1_{W_3}))'(w'_{(4)}), w_{(1)}}$ is a generalized 
$V$-module.
\epfv

\begin{rema}\label{bar-boxbackslash}
{\rm Let $W_1$, $W_2$ and $W_3$ be objects of $\mathcal{C}$.  If
$\lambda\in (W_1\otimes W_2\otimes W_3)^*$ satisfies the
$P^{(2)}(z)$-local grading restriction condition for suitable $z\in
\C^{\times}$, then as we discussed in Remark \ref{part-a}, Part (a) of
the condition includes the assertion that for any $w_{(3)}\in W_3$,
$\mu^{(2)}_{\lambda, w_{(3)}}=\lambda\circ(\cdot\otimes w_{(3)})$ is
the sum of a weakly absolutely convergent series of $P(z)$-generalized
weight vectors in $(W_1\otimes W_2)^*$. If $\lambda$ is obtained from
the product or iterate of intertwining maps, then by Proposition
\ref{9.7} and Theorem \ref{9.7-1}, $W_{\lambda, w_{(3)}}$ is a
strongly-graded generalized $V$-module and therefore each of such
$P(z)$-generalized weight vectors mentioned above generates a
generalized $V$-submodule of $W_{\lambda, w_{(3)}}$.  If, in addition,
these generalized $V$-modules are objects of $\mathcal{C}$, then these
$P(z)$-generalized weight vectors in $(W_1\otimes W_2)^*$ actually
belong to $W_1\hboxtr_{P(z_0)} W_2$. Thus, in this case,
$\mu^{(2)}_{\lambda, w_{(3)}}$ can be identified with an element of
$\overline{W_1\hboxtr_{P(z_0)} W_2}$ for any $w_{(3)}\in W_{3}$.  If
we further assume that $W_1\boxtimes_{P(z_0)} W_2$ exists in
$\mathcal{C}$, then $W_1\hboxtr_{P(z_0)} W_2$ is an object of
$\mathcal{C}$ and in particular, $W_{\lambda, w_{(3)}}$ is an object
of $\mathcal{C}$.}
\end{rema}

Let $W_1$, $W_2$ and $W_3$ be objects of $\mathcal{C}$ and let $G\in
\hom (W'_4, (W_1\otimes W_2\otimes W_3)^*)$ be such that $G(w'_{(4)})$
satisfies the $P^{(2)}(z_0)$-local grading restriction condition. For
$w_{(3)}\in w_{3}$, let $\sum_{n\in \C}\lambda_{n}^{(2)}$ be the
series weakly absolutely convergent to $\mu^{(2)}_{G(w'_{(4)}),
w_{(3)}}$. If the generalized $V$-submodule of $W^{(2)}_{G(w'_{(4)}),
w_{(3)}}$ generated by each $\lambda_{n}^{(2)}$ is an object of
$\mathcal{C}$, then by Remark \ref{bar-boxbackslash},
$\mu^{(2)}_{G(w'_{(4)}), w_{(3)}}$ can be identified with an element
of $\overline{W_1\hboxtr_{P(z_0)} W_2}$ for any $w_{(3)}\in W_{3}$.
If $W_1\boxtimes_{P(z_0)} W_2$ exists in $\mathcal{C}$, then we have a
map $\tilde G$ {}from $W_1\boxtimes_{P(z_0)} W_2$ to $(W'_4\otimes
W_3)^*$ defined by
\begin{equation}\label{tildeG}
\tilde G(w)(w'_{(4)}\otimes w_{(3)})=\langle
w,\mu^{(2)}_{G(w'_{(4)}),w_{(3)}}\rangle_{W_1\hboxtr_{P(z_0)} W_2}
\end{equation}
for $w\in W_1\boxtimes_{P(z_0)} W_2$.
Moreover, we have, as in \cite{tensor4}:

\begin{lemma}
Assume the convergence condition for
intertwining maps in $\mathcal{C}$ and that all generalized
$V$-modules considered are objects of $\mathcal{C}$. 
Let $I_{1}$ and $I_{2}$ be intertwining maps in Lemma \ref{9.7}
and let $G=(I_1\circ (1_{W_1}\otimes I_2))'$. 
We also assume that
$W_1\boxtimes_{P(z_0)} W_2$ exists in $\mathcal{C}$.
For $w'_{(4)}\in W'_{4}$
and $w_{(3)}\in w_{3}$, let 
$\sum_{n\in \C}\lambda_{n}^{(2)}$ be the series weakly absolutely 
convergent to $\mu^{(2)}_{G(w'_{(4)}), w_{(3)}}$.
If $G(w'_{(4)})$ satisfies the
$P^{(2)}(z_0)$-local grading restriction condition (respectively,
the $L(0)$-semisimple
$P^{(2)}(z_0)$-local grading restriction condition when $\mathcal{C}$ 
is in $\mathcal{M}_{sg}$) and the 
generalized $V$-submodule $W_{\lambda_{n}^{(2)}}$
of $W^{(2)}_{G(w'_{(4)}), w_{(3)}}$ generated by each $\lambda_{n}^{(2)}$ is 
an object of $\mathcal{C}$, 
then for $w\in W_1\boxtimes_{P(z_0)} W_2$, 
$\tilde G(w)\in\hom(W_1\boxtimes_{P(z_0)}
W_2, (W'_4\otimes W_3)^*)$ intertwines the two actions
$\tau_{W_1\hboxtr_{P(z_0)} W_2}$ and $\tau_{Q(z_2)}$ of $V\otimes
\iota_+{\mathbb C}[t,t^{-1}, (z_2+t)^{-1}]$ on $W_1\hboxtr_{P(z_0)} W_2$
and $(W'_4\otimes W_3)^*$.
\end{lemma}
\pf By (\ref{tildeG}) and the definition of $\tau_{Q(z)}$-action we
need to show that
\begin{eqnarray*}
\lefteqn{z^{-1}_2\delta\left(\frac{x_1-x_0}{z_2}\right)
\tilde{G}(Y_{P(z_0)}(v, x_0)w)}\nno\\
&&=\tau_{Q(z_2)}\Biggl(
z^{-1}_2\delta\left(\frac{x_1-x_0}{z_2}\right) Y_{t}(v,
x_0)\Biggr)\tilde{G}(w)
\end{eqnarray*}
for any $w\in W_1\boxtimes_{P(z_0)} W_2$,
or equivalently, we need to show that 
\begin{eqnarray}\label{as:need0}
\lefteqn{\left\langle z^{-1}_2\delta\left(\frac{x_1-x_0}{z_2}\right)
Y_{P(z_0)}(v, x_0)w, \mu^{(2)}_{G(w_{(4)}'),
w_{(3)}} \right\rangle_{W_1\shboxtr_{P(z_0)}W_2}}\nno\\
&&=\tau_{Q(z_2)}\Biggl(
z^{-1}_2\delta\left(\frac{x_1-x_0}{z_2}\right) Y_{t}(v,
x_0)\Biggr)(\langle w, \mu^{(2)}_{G(w_{(4)}'),
w_{(3)}} \rangle_{W_1\shboxtr_{P(z_0)}W_2})
\end{eqnarray}
for any $w\in W_1\boxtimes_{P(z_0)} W_2$, $w_{(4)}'\in W_{4}'$,
$w_{(3)}\in W_{3}$.

By Proposition \ref{span}, 
we need only prove this for 
$w=\pi_{n}(w_{(1)}\boxtimes_{P(z_0)} w_{(2)})$ for any 
$n\in \C$ and $w_{(1)}\in W_{1}$ and $w_{(2)}\in W_{2}$. To prove 
(\ref{as:need0}) for such $w$, it is enough to prove that when we substitute 
$w=w_{(1)}\boxtimes_{P(z_0)} w_{(2)}$ into both sides of
(\ref{as:need0}), they are formal series in $x_{1}$ and $x_{0}$
whose coefficients are absolutely convergent series and the formal series
in $x_{1}$ and $x_{0}$ obtained by taking sums of their coefficients
of monomials in $x_{1}$ and $x_{0}$  are equal. 

Since $G(w'_{(4)})$ satisfies the
$P^{(2)}(z_0)$-local grading restriction condition and 
the generalized $V$-submodule $W_{\lambda_{n}^{(2)}}$
of $W^{(2)}_{G(w'_{(4)}), w_{(3)}}$ generated by $\lambda_{n}^{(2)}$ is
an object of $\mathcal{C}$, 
by Remark \ref{bar-boxbackslash},  $\mu^{(2)}_{G(w'_{(4)}),
w_{(3)}}\in \overline{W_{1}\hboxtr_{P(z_{0})} W_{2}}$.
So there exist 
$\lambda_{n}\in (W_{1}\hboxtr_{P(z_{0})} W_{2})_{[n]}\subset
(W_{1}\otimes W_{2})^{*}$
for $n\in \C$
such that $\mu^{(2)}_{G(w'_{(4)}),
w_{(3)}}(w_{(1)}\otimes w_{(2)})$ is the sum of the 
absolutely convergent series $\sum_{n\in \C}
\lambda_{n}(w_{(1)}\otimes w_{(2)})$ and as elements of 
$\overline{W_{1}\hboxtr_{P(z_{0})} W_{2}}$, 
$\mu^{(2)}_{G(w'_{(4)}),
w_{(3)}}=\sum_{n\in \C}
\lambda_{n}$. 
By (\ref{boxpair}), we have 
\[
\lambda_{n}(w_{(1)}\otimes w_{(2)})=\langle \lambda_{n}, 
w_{(1)}\boxtimes_{P(z_{0})} w_{(2)}\rangle.
\]
Taking sum over $\C$ of both sides and using the fact that 
$\sum_{n\in \C}
\lambda_{n}(w_{(1)}\otimes w_{(2)})$ is absolutely convergent to 
$\mu^{(2)}_{G(w'_{(4)}),
w_{(3)}}(w_{(1)}\otimes w_{(2)})$ and as elements of 
$\overline{W_{1}\hboxtr_{P(z_{0})} W_{2}}$, 
$\mu^{(2)}_{G(w'_{(4)}),
w_{(3)}}=\sum_{n\in \C}
\lambda_{n}$,
we obtain
\[
\mu^{(2)}_{G(w'_{(4)}),
w_{(3)}}(w_{(1)}\otimes w_{(2)})=\langle \mu^{(2)}_{G(w'_{(4)}),
w_{(3)}}, 
w_{(1)}\boxtimes_{P(z_{0})} w_{(2)}\rangle,
\]
where the right-hand side is defined to be the sum 
of the absolutely convergent series
\[
\sum_{n\in \C}\langle \mu^{(2)}_{G(w'_{(4)}),
w_{(3)}}, 
\pi_{n}(w_{(1)}\boxtimes_{P(z_{0})} w_{(2)})\rangle.
\]
Repeat the same auguments for $Y'_{P(z_{0})}(v, x_{0})\mu^{(2)}_{G(w'_{(4)}),
w_{(3)}}$ instead of $\mu^{(2)}_{G(w'_{(4)}),
w_{(3)}}$, we obtain
\begin{equation}\label{*-to-box}
(Y'_{P(z_{0})}(v, x_{0})\mu^{(2)}_{G(w'_{(4)}),
w_{(3)}})(w_{(1)}\otimes w_{(2)})=\langle Y'_{P(z_{0})}(v, x_{0})
\mu^{(2)}_{G(w'_{(4)}),
w_{(3)}}, 
w_{(1)}\boxtimes_{P(z_{0})} w_{(2)}\rangle.
\end{equation}

The equality (\ref{*-to-box}) shows that in particular, 
we can take $w=w_{(1)}\boxtimes_{P(z_0)} w_{(2)}$ on 
both sides of (\ref{as:need0}). 
Now taking  $w=w_{(1)}\boxtimes_{P(z_0)} w_{(2)}$ in 
the left-hand side of (\ref{as:need0}) and using (\ref{y'}) and
(\ref{yo}) (and, as we have done above, 
dropping the subscripts for the
inner product), we see that the left-hand side of (\ref{as:need0})
gives
\begin{eqnarray*}
\lefteqn{\left\langle z^{-1}_2\delta\left(\frac{x_1-x_0}{z_2}\right)
Y_{P(z_0)}(v, x_0)(w_{(1)}\boxtimes_{P(z_0)} w_{(2)}), 
\mu^{(2)}_{G(w'_{(4)}),w_{(3)}} \right\rangle}\nno\\
&&=z^{-1}_2\delta\left(\frac{x_1-x_0}{z_2}\right)\langle
Y'^o_{P(z_0)}(v, x_0)\mu^{(2)}_{G(w'_{(4)}),
w_{(3)}}, w_{(1)}\boxtimes_{P(z_0)}
w_{(2)}\rangle\nno\\
&&=z^{-1}_2\delta\left(\frac{x_1-x_0}{z_2}\right)
(Y'^o_{P(z_0)}(v, x_0)\mu^{(2)}_{G(w'_{(4)}),
w_{(3)}})(w_{(1)}\boxtimes_{P(z_0)}
w_{(2)})\nno\\
&&=z^{-1}_2\delta\left(\frac{x_1-x_0}{z_2}\right)
\big(Y'_{P(z_0)}(e^{x_0L(1)}(-x_0^{-2})^{L(0)}v,
x_0^{-1})\mu^{(2)}_{G(w'_{(4)}),
w_{(3)}}\big)(w_{(1)}\otimes w_{(2)}).
\end{eqnarray*}
By the definition of $Y'_{P(z_0)}$ (see (\ref{Y'def})) and the
definition of $\mu^{(2)}_{G(w'_{(4)}),w_{(3)}}$, the right-hand side gives
\begin{eqnarray}\label{as:l}
\lefteqn{z^{-1}_2\delta\left(\frac{x_1-x_0}{z_2}\right)
\mu^{(2)}_{G(w'_{(4)}),w_{(3)}}
(w_{(1)}\otimes Y_2(v,x_0)w_{(2)})}\nno\\
&&+z^{-1}_2\delta\left(\frac{x_1-x_0}{z_2}\right)\res_{x_2}z_0^{-1}
\delta\bigg(\frac{x_0-x_2}{z_0}\bigg)\mu^{(2)}_{G(w'_{(4)}),w_{(3)}}
(Y_1(v, x_2)w_{(1)}\otimes w_{(2)}))\nno\\
&&=z^{-1}_2\delta\left(\frac{x_1-x_0}{z_2}\right)\langle w'_{(4)},
I_1(w_{(1)},z_1)I_2(Y_2(v,x_0)w_{(2)},z_2)w_{(3)}\rangle\nno\\
&&\quad +z^{-1}_2\delta\left(\frac{x_1-x_0}{z_2}\right)\res_{x_2}z_0^{-1}
\delta\bigg(\frac{x_0-x_2}{z_0}\bigg)\langle w'_{(4)},
I_1(Y_1(v,x_2)w_{(1)},z_1)I_2(w_{(2)},z_2)w_{(3)}\rangle.\nn
\end{eqnarray}

Taking 
$w=w_{(1)}\boxtimes_{P(z_0)} w_{(2)}$ in the right-hand side of 
(\ref{as:need0}) and using the definition 
(\ref{5.2}) of
$\tau_{Q(z_2)}$, (\ref{boxpair})
and the  definition of $\mu^{(2)}_{G(w_{(4)}',w_{(3)}}$, we obtain
\begin{eqnarray*}
\lefteqn{\bigg(\tau_{Q(z_2)}\biggl(z^{-1}_2\delta\left(\frac{x_1-x_0}
{z_2}\right) Y_{t}(v,x_0)\biggr)(\langle w_{(1)}\boxtimes_{P(z_{0})}
w_{(2)}, \mu^{(2)}_{G(\cdot), \cdot} \rangle)\bigg)(w'_{(4)}\otimes
w_{(3)})}\\
&&=x_0^{-1}\delta\bigg(\frac{x_1-z_2}{x_0}\bigg)\langle w,
\mu^{(2)}_{G(Y'^o_4(v,x_1)w'_{(4)}),w_{(3)}}
\rangle\\
&&\qquad-x_0^{-1}\delta\bigg(\frac{z_2-x_1}{-x_0}\bigg)\langle w,
\mu^{(2)}_{G(w'_{(4)}),Y_3(v,x_1)w_{(3)}}
\rangle\\
&&=x_0^{-1}\delta\bigg(\frac{x_1-z_2}{x_0}\bigg)\langle w_{(1)}\otimes
w_{(2)},\mu^{(2)}_{G(Y'^o_4(v,x_1)w'_{(4)}),
w_{(3)}}\rangle\\
&&\qquad-x_0^{-1}\delta\bigg(\frac{z_2-x_1}{-x_0}\bigg)\langle w_{(1)}
\otimes w_{(2)},\mu^{(2)}_{G(w'_{(4)}),
Y_3(v,x_1)w_{(3)}}\rangle\\
&&=x_0^{-1}\delta\bigg(\frac{x_1-z_2}{x_0}\bigg)\langle
Y'^o_4(v,x_1)w'_{(4)},I_1(w_{(1)},z_1)I_2(w_{(2)},z_2) w_{(3)}
\rangle\\
&&\qquad-x_0^{-1}\delta\bigg(\frac{z_2-x_1}{-x_0}\bigg)\langle
w'_{(4)},I_1(w_{(1)},z_1)I_2(w_{(2)},z_2)Y_3(v,x_1)w_{(3)}
\rangle\\
&&=x_0^{-1}\delta\bigg(\frac{x_1-z_2}{x_0}\bigg)\langle
w'_{(4)},Y_4(v,x_1)I_1(w_{(1)},z_1)I_2(w_{(2)},z_2) w_{(3)}
\rangle\\
&&\qquad-x_0^{-1}\delta\bigg(\frac{z_2-x_1}{-x_0}\bigg)\langle
w'_{(4)},I_1(w_{(1)},z_1)I_2(w_{(2)},z_2)Y_3(v,x_1)w_{(3)}\rangle.
\end{eqnarray*}

Now using the formula obtained by taking $\res_{x_1}$ of (\ref{F12})
then replacing $x_0$ by $x_1$ and $x_2$ by $x_0$, we see that this
is equal to
\begin{eqnarray}\label{as:r}
&&z_2^{-1}\delta\bigg(\frac{x_1-x_0}{z_2}\bigg)\langle
w'_{(4)},I_1(w_{(1)},z_1)I_2(Y_2(v,x_0)w_{(2)},z_2)w_{(3)}
\rangle\nno\\
&&\quad+x_0^{-1}\delta\bigg(\frac{x_1-z_2}{x_0}\bigg)\res_{x_2}
z_1^{-1}\delta\bigg(\frac{x_1-x_2}{z_1}\bigg)\langle w'_{(4)},
I_1(Y_1(v,x_2)w_{(1)},z_1)I_2(w_{(2)},z_2)w_{(3)}\rangle.\nn
\end{eqnarray}

We now need to show that the right-hand side of (\ref{as:l}) equals
(\ref{as:r}). But their first terms are the same, while their second
terms are equal due to the following delta function identities 
obtained using (\ref{l2-2}) and (\ref{l2-1}):
\begin{eqnarray*}
z^{-1}_2\delta\left(\frac{x_1-x_0}{z_2}\right)
z_0^{-1}\delta\bigg(\frac{x_0-x_2}{z_0}\bigg)
&=&x^{-1}_1\delta\left(\frac{z_1+x_2}{x_1}\right)
x_0^{-1}\delta\bigg(\frac{z_0+x_2}{x_0}\bigg)
\\
&=&z_1^{-1}\delta\bigg(\frac{x_1-x_2}{z_1}\bigg)
x_0^{-1}\delta\bigg(\frac{x_1-z_2}{x_0}\bigg)
.
\end{eqnarray*}
\epfv

Under the assumptions of this lemma, by Proposition
\ref{pz} there is a $P(z_2)$-intertwining map $I$ of type ${W_4\choose
W_1\boxtimes_{P(z_0)} W_2\,\,W_3}$ such that
\[
\tilde G(w)(w'_{(4)}\otimes w_3)=\langle w'_{(4)},I(w,z_2)
w_{(3)}\rangle.
\]
Hence by (\ref{tildeG}) we have
\[
\langle w,\mu^{(2)}_{G(w'_{(4)}),w_{(3)}}\rangle_{W_1\hboxtr_{P(z_0)}
W_2}= \langle w'_{(4)},I(w,z_2) w_{(3)}\rangle
\]
for any $w\in W_1\boxtimes_{P(z_0)} W_2$. In particular, since the
convergence condition for intertwining maps 
in ${\cal C}$ is satisfied, we take
$w=\pi_{n}(w_{(1)}\boxtimes_{P(z_0)} w_{(2)})$, sum over $n\in \C$ and get
\[
\langle w'_{(4)},I_1(w_{(1)},z_1)I_2(w_{(2)},z_2)w_{(3)}\rangle
=\langle w'_{(4)},I(w_{(1)}\boxtimes_{P(z_0)} w_{(2)},z_2)
w_{(3)}\rangle.
\]

We have now proved the first half of the following theorem,
which states in particular that under natural conditions we 
have discussed above,
this product of intertwining maps can be written 
as an iterate of certain intertwining maps and the 
intermediate module can be taken as $W_{1}\boxtimes_{P(z_0)} W_{2}$; the other
half can be proved similarly:

\begin{theo}\label{lgr=>asso}
Assume the convergence condition for
intertwining maps in $\mathcal{C}$ and that all generalized
$V$-modules considered are objects of $\mathcal{C}$. 
We also assume that
$W_1\boxtimes_{P(z_0)} W_2$ exists in $\mathcal{C}$.
Let $I_1$ and $I_2$ be $P(z_1)$- and $P(z_2)$-intertwining maps of
type ${W_4 \choose W_1\, M_1}$ and ${M_1 \choose W_2\, W_3}$,
respectively.  Suppose that $(I_1\circ (1_{W_1}\otimes
I_2))'(w'_{(4)})$ satisfies the $P^{(2)}(z_0)$-local grading
restriction condition (or the $L(0)$-semisimple $P^{(2)}(z_0)$-local
grading restriction condition when $\mathcal{C}$ is in
$\mathcal{M}_{sg}$). For $w'_{(4)}\in W'_{4}$
and $w_{(3)}\in w_{3}$, let 
$\sum_{n\in \C}\lambda_{n}^{(2)}$ be the series weakly absolutely 
convergent to $\mu^{(2)}_{(I_1\circ (1_{W_1}\otimes
I_2))'(w'_{(4)}), w_{(3)}}$. Suppose also that 
the 
generalized $V$-submodule $W_{\lambda_{n}^{(2)}}$
of $W^{(2)}_{(I_1\circ (1_{W_1}\otimes
I_2))'(w'_{(4)}), w_{(3)}}$ generated by $\lambda_{n}^{(2)}$
for any $n\in \C$, $w'_{(4)}\in W'_4$ and $w_{(3)}\in W_{3}$  is 
an object of $\mathcal{C}$.  Then there is a
$P(z_2)$-intertwining map $I$ of type ${W_4\choose
W_1\boxtimes_{P(z_0)} W_2\,\,W_3}$ such that
\[
\langle w'_{(4)},I_1(w_{(1)} \otimes I_2(w_{(2)} \otimes w_{(3)}))\rangle
=\langle w'_{(4)}, I((w_{(1)}\boxtimes_{P(z_0)} w_{(2)})\otimes 
w_{(3)})\rangle
\]
for any $w_{(1)}\in W_1$, $w_{(2)}\in W_2$, $w_{(3)}\in W_3$ and
$w'_{(4)}\in W'_4$.  Analogously, let $I^1$ and $I^2$ be $P(z_2)$- and
$P(z_0)$-intertwining maps of type ${W_4 \choose M_2\, W_3}$ and ${M_2
\choose W_1\, W_2}$, respectively.  Suppose that $(I^1\circ (I^2
\otimes 1_{W_3}))'(w'_{(4)})$ satisfies the $P^{(1)}(z_2)$-local grading
restriction condition (or the $L(0)$-semisimple $P^{(1)}(z_2)$-local
grading restriction condition when $\mathcal{C}$ is in
$\mathcal{M}_{sg}$). For $w'_{(4)}\in W'_{4}$
and $w_{(1)}\in w_{1}$, let 
$\sum_{n\in \C}\lambda_{n}^{(1)}$ be the series weakly absolutely 
convergent to $\mu^{(1)}_{(I^1\circ (I^2
\otimes 1_{W_3}))'(w'_{(4)}), w_{(1)}}$. Suppose also that 
the 
generalized $V$-submodule $W_{\lambda_{n}^{(1)}}$
of $W^{(1)}_{(I^1\circ (I^2
\otimes 1_{W_3}))'(w'_{(4)}), w_{(1)}}$ generated by $\lambda_{n}^{(1)}$
for any $n\in \C$, $w'_{(4)}\in W'_4$ and $w_{(1)}\in W_{3}$  is 
an object of $\mathcal{C}$.  Then there is a
$P(z_1)$-intertwining map $I$ of type ${W_4\choose
W_1\,\,W_2\boxtimes_{P(z_2)} W_3}$ such that
\[
\langle w'_{(4)}, I^1(I^2(w_{(1)}\otimes w_{(2)})\otimes w_{(3)})\rangle
=\langle w'_{(4)}, I(w_{(1)}\otimes (w_{(2)}\boxtimes_{P(z_2)}w_{(3)}))
\rangle
\]
for any $w_{(1)}\in W_1$, $w_{(2)}\in W_2$, $w_{(3)}\in W_3$ and
$w'_{(4)}\in W'_4$. \epf
\end{theo}

Recall that in Section \ref{convsec} we have proved two formulas
((\ref{nosub}) and (\ref{nosub2})) on writing products of intertwining
operators satisfying certain conditions in terms of iterates and vice
versa. Now we will need a refinement of these two formulas.  The
two formulas in the 
following lemma should be compared to the intuitive guides
``$(12)3=3(12)=3(21)$'' and ``$1(23)=(23)1=(32)1$,'' respectively:

\begin{lemma}
\begin{enumerate}

\item For any intertwining operators ${\cal Y}^1$ and ${\cal Y}^2$ of
types ${W_4}\choose {M_2 W_3}$ and ${M_2}\choose {W_1 W_2}$,
respectively, and any nonzero complex numbers $z_1$, $z_2$, when
$|z_1|>|z_1-z_2|>0$ and $|z_2|>|z_1-z_2|>0$ we have
\begin{eqnarray}\label{(12)3-3(21)}
\lefteqn{\langle w'_{(4)}, {\cal Y}^1({\cal
Y}^2(w_{(1)},x_0)w_{(2)},x_2)w_{(3)}\rangle_{W_4}\lbar_{x_0=z_1-z_2,\;
x_2=z_2}}\nno\\
&&=\langle e^{z_1L'(1)}w'_{(4)}, \Omega_0({\cal Y}^1) (w_{(3)},
y_1)\Omega_{-1}({\cal Y}^2)(w_{(2)}, y_2)
w_{(1)}\rangle_{W_4}\lbar_{y_1=e^{\pi i}z_1,\; y_2=e^{\pi
i}(z_1-z_2)}.\nno\\
\end{eqnarray}

\item For any intertwining operators ${\cal Y}_1$ and ${\cal Y}_2$ of
types ${W_4}\choose {W_1 M_1}$ and ${M_1}\choose {W_2 W_3}$,
respectively, and any nonzero complex numbers $z_1$, $z_2$, when
$|z_1|>|z_2|>0$ and $|z_1-z_2|>|z_2|>0$ we have
\begin{eqnarray}\label{1(23)-(32)1}
\lefteqn{\langle w'_{(4)}, {\cal Y}_1(w_{(1)},x_1){\cal
Y}_2(w_{(2)},x_2)w_{(3)}\rangle_{W_4}\lbar_{x_1=z_1,\; x_2=z_2}}\nno\\
&&=\langle e^{z_1L'(1)}w'_{(4)}, \Omega_{-1}({\cal
Y}_1)(\Omega_{-1}({\cal Y}_2)(w_{(3)}, y_0)w_{(2)},
y_2)w_{(1)} \rangle_{W_4}\lbar_{y_0= e^{\pi i}z_2,\; y_2=e^{\pi
i}(z_1-z_2)}.\nno\\
\end{eqnarray}
\end{enumerate}
\end{lemma}
\pf The proofs of these formulas are the following straighforward
calculations, using in particular the formulas (\ref{i2p}) and
(\ref{p2i}), respectively:
\begin{eqnarray*}
\lefteqn{\langle w'_{(4)}, {\cal Y}^1({\cal
Y}^2(w_{(1)},x_0)w_{(2)},x_2)w_{(3)}\rangle_{W_4}\lbar_{x_0=z_1-z_2,\;
x_2=z_2}}\nno\\
&&=\langle e^{z_2L'(1)}w'_{(4)}, \Omega_0({\cal Y}^1) (w_{(3)},
x_2)\Omega_0(\Omega_{-1}({\cal Y}^2))(w_{(1)},
x_0)w_{(2)}\rangle_{W_4}\lbar_{x_0=z_1-z_2,\; x_2=e^{-\pi
i}z_2}\nno\\
&&=\langle e^{z_2L'(1)}w'_{(4)}, \Omega_0({\cal Y}^1) (w_{(3)},
x_2)e^{x_0L(-1)}\Omega_{-1}({\cal Y}^2)(w_{(2)}, e^{\pi i}x_0)
w_{(1)}\rangle_{W_4}\lbar_{x_0=z_1-z_2,\; x_2=e^{-\pi i}z_2}\nno\\
&&=\langle e^{z_2L'(1)}w'_{(4)}, e^{x_0L(-1)}\Omega_0({\cal Y}^1)
(w_{(3)}, x_2-x_0)\Omega_{-1}({\cal Y}^2)(w_{(2)}, e^{\pi i}x_0)
w_{(1)}\rangle_{W_4}\lbar_{x_0=z_1-z_2,\; x_2=e^{-\pi i}z_2}\nno\\
&&=\langle e^{x_0L'(1)}e^{z_2L'(1)}w'_{(4)}, \Omega_0({\cal Y}^1)
(w_{(3)}, x_2-x_0)\Omega_{-1}({\cal Y}^2)(w_{(2)}, e^{\pi i}x_0)
w_{(1)}\rangle_{W_4}\lbar_{x_0=z_1-z_2,\; x_2=e^{-\pi i}z_2}\nno\\
&&=\langle e^{z_1L'(1)}w'_{(4)}, \Omega_0({\cal Y}^1) (w_{(3)},
y_1)\Omega_{-1}({\cal Y}^2)(w_{(2)}, y_2)
w_{(1)}\rangle_{W_4}\lbar_{y_1=e^{\pi i}z_1,\; y_2=e^{\pi
i}(z_1-z_2)},
\end{eqnarray*}
\begin{eqnarray*}
\lefteqn{\langle w'_{(4)}, {\cal Y}_1(w_{(1)},x_1){\cal
Y}_2(w_{(2)},x_2)w_{(3)}\rangle_{W_4}\lbar_{x_1=z_1,\; x_2=z_2}}\nno\\
&&=\langle e^{z_1L'(1)}w'_{(4)}, \Omega_{-1}({\cal Y}_1)
(\Omega_0(\Omega_{-1}({\cal Y}_2))(w_{(2)},
x_2)w_{(3)}, x_1)w_{(1)} \rangle_{W_4}\lbar_{x_1= e^{\pi i}z_1,\;
x_2=z_2}\nno\\
&&=\langle e^{z_1L'(1)}w'_{(4)}, \Omega_{-1}({\cal
Y}_1)(e^{x_2L(-1)}\Omega_{-1}({\cal Y}_2)(w_{(3)}, e^{\pi
i}x_2)w_{(2)}, x_1)w_{(1)} \rangle_{W_4}\lbar_{x_1= e^{\pi i}z_1,\;
x_2=z_2}\nno\\
&&=\langle e^{z_1L'(1)}w'_{(4)}, \Omega_{-1}({\cal
Y}_1)(\Omega_{-1}({\cal Y}_2)(w_{(3)}, e^{\pi
i}x_2)w_{(2)}, x_1+x_2)w_{(1)} \rangle_{W_4}\lbar_{x_1= e^{\pi i}z_1,\;
x_2=z_2}\nno\\
&&=\langle e^{z_1L'(1)}w'_{(4)}, \Omega_{-1}({\cal
Y}_1)(\Omega_{-1}({\cal Y}_2)(w_{(3)}, y_0)w_{(2)},
y_2)w_{(1)} \rangle_{W_4}\lbar_{y_0= e^{\pi i}z_2,\; y_2=e^{\pi
i}(z_1-z_2)}.
\end{eqnarray*}
\epfv

\begin{theo}\label{asso-io}
Assume that the convergence condition for intertwining maps
in ${\cal C}$ holds. Then the
following two conditions are equivalent:
\begin{enumerate}
\item For any objects $W_1$, $W_2$, $W_3$, $W_4$ and $M_1$ of ${\cal
C}$, any nonzero complex numbers $z_1$ and $z_2$ satisfying
$|z_1|>|z_2|>|z_1-z_2|>0$ and any logarithmic 
intertwining operators ${\cal Y}_1$ and
${\cal Y}_2$ of types ${W_4}\choose {W_1 M_1}$ and ${M_1}\choose {W_2
W_3}$, respectively, there exist an object $M_2$ of ${\cal C}$ and 
logarithmic intertwining operators ${\cal Y}^1$ and ${\cal Y}^2$ of types
${W_4}\choose {M_2 W_3}$ and ${M_2}\choose {W_1 W_2}$, respectively,
such that
\begin{eqnarray}\label{extcnd1}
\lefteqn{\langle w'_{(4)}, {\cal Y}_1(w_{(1)},x_1){\cal Y}_2(w_{(2)},
x_2)w_{(3)}\rangle\lbar_{x_1=z_1,\,x_2=z_2}}\nno\\
&&=\langle w'_{(4)}, {\cal Y}^1({\cal Y}^2(w_{(1)},
x_0)w_{(2)},x_2)w_{(3)}\rangle\lbar_{x_0=z_1-z_2,\,x_2=z_2}
\end{eqnarray}
for all $w'_{(4)}\in W'_4$, $w_{(1)}\in W_1$, $w_{(2)}\in W_2$ and
$w_{(3)}\in W_3$.

\item For any objects $W_1$, $W_2$, $W_3$, $W_4$ and $M_2$ of ${\cal
C}$, any nonzero complex numbers $z_1$ and $z_2$ satisfying
$|z_1|>|z_2|>|z_1-z_2|>0$ and any logarithmic
intertwining operators ${\cal Y}^1$ and
${\cal Y}^2$ of types ${W_4}\choose {M_2 W_3}$ and ${M_2}\choose {W_1
W_2}$, respectively, there exist an object $M_1$ of ${\cal C}$ and
logarithmic intertwining operators ${\cal Y}_1$ and ${\cal Y}_2$ of types
${W_4}\choose {W_1 M_1}$ and ${M_1}\choose {W_2 W_3}$, respectively,
such that
\begin{eqnarray}\label{extcnd2}
\lefteqn{\langle w'_{(4)}, {\cal Y}^1({\cal Y}^2(w_{(1)},
x_0)w_{(2)},x_2)w_{(3)}\rangle\lbar_{x_0=z_1-z_2,\,x_2=z_2}}\nno\\
&&=\langle w'_{(4)}, {\cal Y}_1(w_{(1)},x_1){\cal Y}_2(w_{(2)},
x_2)w_{(3)}\rangle\lbar_{x_1=z_1,\,x_2=z_2}
\end{eqnarray}
for all $w'_{(4)}\in W'_4$, $w_{(1)}\in W_1$, $w_{(2)}\in W_2$ and
$w_{(3)}\in W_3$.
\end{enumerate}
\end{theo}
\pf Suppose that Condition 1 holds. Then for any logarithmic intertwining
operators ${\cal Y}^1$ and ${\cal Y}^2$ as in the statement of
Condition 2 and any nonzero complex numbers $z_1$, $z_2$, when
$|z_1|>|z_1-z_2|>0$ and $|z_2|>|z_1-z_2|>0$, by (\ref{(12)3-3(21)}) we
have
\begin{eqnarray}\label{step1}
\lefteqn{\langle w'_{(4)}, {\cal Y}^1({\cal
Y}^2(w_{(1)},x_0)w_{(2)},x_2)w_{(3)}\rangle_{W_4}\lbar_{x_0=z_1-z_2,\;
x_2=z_2}}\nno\\
&&=\langle e^{z_1L'(1)}w'_{(4)}, \Omega_0({\cal Y}^1) (w_{(3)},
y_1)\cdot \nno\\
&&\hspace{4em}\cdot\Omega_{-1}({\cal Y}^2)(w_{(2)}, y_2)
w_{(1)}\rangle_{W_4}\lbar_{y_1=e^{\pi i}z_1,\; y_2=e^{\pi
i}(z_1-z_2)}.
\end{eqnarray}
for any $w_{(1)}\in W_1$, $w_{(2)}\in W_2$, $w_{(3)}\in W_3$ and
$w'_{(4)}\in W'_4$.

Since the last expression is of the same form as the left-hand side of
(\ref{extcnd1}), we have that when $|e^{\pi i}z_1|>|e^{\pi
i}(z_1-z_2)|>|e^{\pi i}z_2|>0$, or equivalently, when
$|z_1|>|z_1-z_2|>|z_2|>0$, there exist object $M_3$ of ${\cal C}$ and
logarithmic intertwining operators ${\cal Y}^3$ and ${\cal Y}^4$ of types
${W_4\choose M_3 W_1}$ and ${M_3\choose W_3 W_2}$, respectively, such
that
\begin{eqnarray}\label{step2}
\lefteqn{\langle e^{z_1L'(1)}w'_{(4)}, \Omega_0({\cal Y}^1)
(w_{(3)}, y_1)\Omega_{-1}({\cal Y}^2)(w_{(2)}, y_2)
w_{(1)}\rangle_{W_4}\lbar_{y_1=e^{\pi i}z_1,\; y_2=e^{\pi
i}(z_1-z_2)}}\nno\\
&&=\langle e^{z_1L'(1)}w'_{(4)}, {\cal Y}^3({\cal
Y}^4(w_{(3)},y_0)w_{(2)},y_2)w_{(1)}\rangle\lbar_{y_0=e^{\pi i}z_2,\;
y_2=e^{\pi i}(z_1-z_2)}.
\end{eqnarray}
By using the fact that ${\cal Y}=\Omega_{-1}(\Omega_0({\cal Y}))$
and comparing the last expression to the right-hand side of
(\ref{1(23)-(32)1}) we get, when $|z_1|>|z_2|>0$ and
$|z_1-z_2|>|z_2|>0$, we have
\begin{eqnarray}\label{step3}
\lefteqn{\langle e^{z_1L'(1)}w'_{(4)}, {\cal Y}^3({\cal
Y}^4(w_{(3)},y_0)w_{(2)},y_2)w_{(1)}\rangle\lbar_{y_0=e^{\pi i}z_2,\;
y_2=e^{\pi i}(z_1-z_2)}}\nno\\
&&=\langle w'_{(4)}, \Omega_0({\cal
Y}^3)(w_{(1)},x_1)\Omega_0({\cal
Y}^4)(w_{(2)},x_2)w_{(3)}\rangle_{W_4}\lbar_{x_1=z_1,\; x_2=z_2}
\end{eqnarray}
for $w_{(1)}\in W_1$, $w_{(2)}\in W_2$, $w_{(3)}\in W_3$ and
$w'_{(4)}\in W'_4$.

By the $L(-1)$-derivative property for intertwining operators, both
sides of (\ref{step1}), (\ref{step2}) and (\ref{step3}) are analytic
(multi-valued) functions of $z_1$ and $z_2$. By these three formulas,
we see that the left-hand side of (\ref{step1}) is equal to a branch
of the analytic extension of the right-hand side of (\ref{step3}) in
the domain $\{(z_1,z_2)\in {\mathbb C}\times{\mathbb C}\;|\;
|z_1|>|z_2|>|z_1-z_2|>0\}$.  This proves Condition 2.

Conversely, suppose that Condition 2 holds. Then for any 
logarithmic intertwining
operators ${\cal Y}_1$ and ${\cal Y}_2$ as in the statement of
Condition 1 and any nonzero complex numbers $z_1$, $z_2$, when
$|z_1|>|z_2|>0$ and $|z_1-z_2|>|z_2|>0$, by (\ref{1(23)-(32)1}) we
have
\begin{eqnarray}\label{step1'}
\lefteqn{\langle w'_{(4)}, {\cal Y}^1({\cal Y}^2(w_{(1)},
x_0)w_{(2)},x_2)w_{(3)}\rangle\lbar_{x_0=z_1-z_2,\,x_2=z_2}}\nno\\
&&=\langle e^{z_1L'(1)}w'_{(4)}, \Omega_{-1}({\cal
Y}_1)(\Omega_{-1}({\cal Y}_2)(w_{(3)}, y_0)w_{(2)},
y_2)w_{(1)} \rangle_{W_4}\lbar_{y_0= e^{\pi i}z_2,\; y_2=e^{\pi
i}(z_1-z_2)}\nno\\
\end{eqnarray}
for $w_{(1)}\in W_1$, $w_{(2)}\in W_2$, $w_{(3)}\in W_3$ and
$w'_{(4)}\in W'_4$.

Since the last expression is in the same form as the left-hand side of
(\ref{extcnd2}), we have that when $|e^{\pi i}z_1|>|e^{\pi
i}(z_1-z_2)|>|e^{\pi i}z_2|>0$, or equivalently, when
$|z_1|>|z_1-z_2|>|z_2|>0$, there exist object $M_4$ of ${\cal C}$ and
logarithmic intertwining operators ${\cal Y}_3$ and ${\cal Y}_4$ of types
${W_4\choose M_4 W_1}$ and ${M_4\choose W_3 W_2}$, respectively, such
that
\begin{eqnarray}\label{step2'}
\lefteqn{\langle e^{z_1L'(1)}w'_{(4)}, \Omega_{-1}({\cal
Y}_1)(\Omega_{-1}({\cal Y}_2)(w_{(3)}, y_0)w_{(2)},
y_2)w_{(1)} \rangle_{W_4}\lbar_{y_0= e^{\pi i}z_2,\; y_2=e^{\pi
i}(z_1-z_2)}}\nno\\
&&=\langle e^{z_1L'(1)}w'_{(4)}, {\cal Y}_3(w_{(3)},y_1)
{\cal Y}_4(w_{(2)},y_2)w_{(1)}\rangle\lbar_{y_1=e^{\pi i}z_1,\;
y_2=e^{\pi i}(z_1-z_2)}
\end{eqnarray}
By using the fact that ${\cal Y}=\Omega_{-1}(\Omega_0({\cal
Y}))=\Omega_0(\Omega_{-1}({\cal Y}))$ and comparing the last
expression to the right-hand side of (\ref{(12)3-3(21)}) we get, when
$|z_1|>|z_2|>0$ and $|z_1-z_2|>|z_2|>0$, we have
\begin{eqnarray}\label{step3'}
\lefteqn{\langle e^{z_1L'(1)}w'_{(4)}, {\cal Y}_3(w_{(3)},y_1){\cal
Y}_4(w_{(2)},y_2)w_{(1)}\rangle\lbar_{y_1=e^{\pi i}z_1,\; y_2=e^{\pi
i}(z_1-z_2)}}\nno\\
&&=\langle w'_{(4)}, \Omega_{-1}({\cal Y}_3)(\Omega_0({\cal
Y}_4)(w_{(3)},x_0)w_{(2)},x_2)w_{(1)}\rangle_{W_4}\lbar_{x_0=z_1-z_2,\;x_2=z_2}
\end{eqnarray}
for $w_{(1)}\in W_1$, $w_{(2)}\in W_2$, $w_{(3)}\in W_3$ and
$w'_{(4)}\in W'_4$.

Since both sides of (\ref{step1'}), (\ref{step2'}) and (\ref{step3'})
are analytic (multi-valued) functions of $z_1$ and $z_2$, by these
three formulas, we see that the left-hand side of (\ref{step1'}) is
equal to a branch of the analytic extension of the right-hand side of
(\ref{step3'}) in the domain $\{(z_1,z_2)\in {\mathbb C}\times{\mathbb
C}\;|\; |z_1|>|z_2|>|z_1-z_2|>0\}$.  This proves Condition 1.
\epf

\begin{theo}\label{expansion}
Assume that the category $\mathcal{C}$ is closed under the $P(z)$-tensor
product operation and the convergence condition for intertwining maps
in ${\cal C}$ holds. Then the
following two conditions are equivalent:
\begin{enumerate}
\item For any objects $W_1$, $W_2$, $W_3$, $W_4$ and $M_1$ of ${\cal
C}$, any nonzero complex numbers $z_1$ and $z_2$ satisfying
$|z_1|>|z_2|>|z_1-z_2|>0$, any $P(z_1)$-intertwining map $I_1$ of type
${W_4}\choose {W_1 M_1}$ and $P(z_2)$-intertwining map $I_2$ of type
${M_1}\choose {W_2W_3}$ and any $w'_{(4)}\in W'_4$, we have that
for any $w'_{(4)}\in W'_4$,
$(I_1\circ (1_{W_1}\otimes I_2))'(w'_{(4)})\in (W_1\otimes W_2\otimes
W_3)^{*}$ satisfies the $P^{(2)}(z_1-z_2)$-local grading restriction
condition (or the $L(0)$-semisimple $P^{(2)}(z_1-z_2)$-local 
grading restriction
condition when $\mathcal{C}$ is in $\mathcal{M}_{sg}$)
and for any $w'_{(4)}\in W'_4$, $w_{(3)}\in W_{3}$ and 
$n\in \C$, the 
generalized $V$-submodule $W_{\lambda_{n}^{(2)}}$
of $W^{(2)}_{(I_1\circ (1_{W_1}\otimes I_2))'(w'_{(4)}), w_{(3)}}$ 
generated by the  term $\lambda_{n}^{(2)}$ of the series 
$\sum_{n\in \C}\lambda_{n}^{(2)}$ 
weakly absolutely convergent to 
$\mu^{(2)}_{(I_1\circ (1_{W_1}\otimes I_2))'(w'_{(4)}), w_{(3)}}$ is 
an object of $\mathcal{C}$.

\item For any objects $W_1$, $W_2$, $W_3$, $W_4$ and $M_2$ of ${\cal
C}$, any nonzero complex numbers $z_1$ and $z_2$ satisfying
$|z_1|>|z_2|>|z_1-z_2|>0$ and any $P(z_2)$-intertwining map $I^1$ of
type ${W_4}\choose {M_2 W_3}$ and $P(z_1-z_2)$-intertwining map $I^2$
of type ${M_2}\choose {W_1W_2}$ and any $w'_{(4)}\in W'_4$, we have
that for any $w'_{(4)}\in W'_4$,
$(I^1\circ (I^2\otimes 1_{W_3}))'(w'_{(4)})\in (W_1\otimes W_2\otimes
W_3)^{*}$ satisfies the $P^{(1)}(z_2)$-local grading restriction
condition (or the $L(0)$-semisimple $P^{(2)}(z_2)$-local 
grading restriction
condition when $\mathcal{C}$ is in $\mathcal{M}_{sg}$) and 
for any $w'_{(4)}\in W'_4$, $w_{(3)}\in W_{3}$ and 
$n\in \C$, the 
generalized $V$-submodule $W_{\lambda_{n}^{(1)}}$
of $W^{(1)}_{(I^1\circ (I^2\otimes 1_{W_3}))'(w'_{(4)}), w_{(1)}}$ 
generated by the term of $\lambda_{n}^{(1)}$ of the series 
$\sum_{n\in \C}\lambda_{n}^{(1)}$ 
weakly absolutely convergent to 
$\mu^{(2)}_{(I^1\circ (I^2\otimes 1_{W_3}))'(w'_{(4)}), w_{(1)}}$ is 
an object of $\mathcal{C}$.
\end{enumerate}
\end{theo}
\pf By Propositions \ref{im:correspond} and \ref{9.7},  
Conditions 1 and  2 in Theorem
\ref{asso-io} imply Conditions 1 and 2, respectively,
in this theorem. By Proposition \ref{im:correspond} and 
Theorem \ref{lgr=>asso}, Conditions 1 and 2 in this theorem imply
Conditions 1 and  2 in Theorem
\ref{asso-io}. Thus this
theorem follows immediately {}from Theorem \ref{asso-io}.  
\epfv

\begin{defi}
{\rm We call either of the two conditions in Theorem
\ref{expansion} the {\it
expansion condition for intertwining maps in the category ${\cal C}$}.}
\end{defi}

\newpage

\setcounter{equation}{0}
\setcounter{rema}{0}

\section{The associativity isomorphisms}

We are now in a position to construct the associativity isomorphisms,
assuming the convergence and expansion conditions. The strategy and
steps in our construction in this section are essentially the same as
those in \cite{tensor4} in the finitely reductive case but, instead of
the corresponding results in \cite{tensor1}, \cite{tensor2},
\cite{tensor3} and \cite{tensor4}, we have to use all the
constructions and results we obtained so far in this work. We remark
that the construction presented here will make easy the proofs of, for
example, the ``coherence property'' needed in braided tensor category
theory.

In the remainder of this work, in additon to Assumptions \ref{assum},
\ref{assum-c} and \ref{assum-exp-set}, we shall also assume that for
some $z\in \C^{\times}$, the category $\mathcal{C}$ is closed under
$P(z)$-tensor products, that is, the $P(z)$-tensor product of $W_{1},
W_{2}\in \ob \mathcal{C}$ exists (in $\mathcal{C}$). For the reader's
convenience, we combine all these assumptions as follows:

\begin{assum}\label{assum-assoc}
Throughout the remainder of this work, we shall assume the following,
unless other assumptions are explicitly made: $A$ is an abelian group
and $\tilde{A}$ is an abelian group containing $A$ as a subgroup; $V$
is a strongly $A$-graded M\"{o}bius or conformal vertex algebra; all
$V$-modules and generalized $V$-modules considered are strongly
$\tilde{A}$-graded; all intertwining operators and logarithmic
intertwining operators considered are grading-compatible; for 
any object $W$ of $\mathcal{C}$, 
the set $\{(n, i)\in \C\times \N\;|\; W_{(n)}\ne 0, \; 
(L(0)-n)^{i}W_{(n)}\ne 0\}$ is a unique expansion set;
for any
objects $W_{1}$, $W_{2}$ and $W_{3}$ of $\mathcal{C}$, any
intertwining operator $\Y$ of type ${W_{3}\choose W_{1}W_{2}}$ and any
$w_{(1)}\in W_{1}$ and $w_{(2)}\in W_{2}$, all the powers of $x$ and
$\log x$ occurring in
\[\Y(w_{(1)}, x)w_{(2)}
\]
form a unique expansion set of the form $D\times \{0, \dots, N\}$
where $D$ is a set of real numbers; $\mathcal{C}$ is a full
subcategory of the category $\mathcal{M}_{sg}$ or
$\mathcal{G}\mathcal{M}_{sg}$ closed under the contragredient functor;
$\mathcal{C}$ is closed under taking finite direct sums; and for some
$z\in \C^{\times}$, $\mathcal{C}$ is closed under $P(z)$-tensor
products.
\end{assum}

\begin{rema}
{\rm {}From Proposition \ref{4.19}, we see that the last part of
Assumption \ref{assum-assoc} is equivalent to the assertion that for
{\it every} $z\in \C^{\times}$, $\mathcal{C}$ is closed under
$P(z)$-tensor product.  Also, by Proposition \ref{tensor1-13.7}, the
last part of Assumption \ref{assum-assoc} is equivalent to the
assumption that for any $W_{1}, W_{2}\in \ob \mathcal{C}$,
$W_{1}\hboxtr_{P(z)}W_{2}$ is an object of $\mathcal{C}$.}
\end{rema}

In this section, we assume further that the convergence condition and
expansion condition for intertwining operators in ${\cal C}$ both
hold.

\begin{theo}
Let $z_1$, $z_2$ be complex numbers satisfying
$|z_1|>|z_2|>|z_1-z_{2}|>0$ (so in particular $z_1\neq 0$, $z_2\neq 0$ and
$z_1\neq z_2$). Let $W_1$, $W_2$, and $W_3$ be
objects of ${\cal C}$. If the convergence condition and expansion
condition for intertwining operators in ${\cal C}$ both hold,
then there exist a unique {\it associativity isomorphism}
\begin{equation}\label{assoc-iso}
\mathcal{A}_{P(z_{1}), P(z_{2})}^{P(z_{1}-z_{2}), P(z_{2})}: 
W_1\boxtimes_{P(z_{1})}
(W_2\boxtimes_{P(z_2)} W_3) \longrightarrow (W_1\boxtimes_{P(z_1-z_2)}
W_2)\boxtimes_{P(z_2)} W_3.
\end{equation}
such that 
\begin{equation}\label{assoc-elt-1}
\overline{\mathcal{A}_{P(z_{1}), P(z_{2})}^{P(z_{1}-z_{2}), P(z_{2})}}
(w_{(1)}\boxtimes_{P(z_1)}
(w_{(2)}\boxtimes_{P(z_2)} w_{(3)})) = (w_{(1)}\boxtimes_{P(z_1-z_2)}
w_{(2)})\boxtimes_{P(z_2)} w_{(3)},
\end{equation}
for all $w_{(1)}\in W_1$, $w_{(2)}\in W_2$ and
$w_{(3)}\in W_3$.
\end{theo}
\pf
As in \cite{tensor4}, we denote the subspace of $(W_1\otimes
W_2\otimes W_3)^{*}$ consisting of all elements satisfying the
$P(z_1,z_2)$-compatibility condition, the $P(z_1, z_2)$-, the
$P^{(1)}(z_2)$- and the $P^{(2)}(z_{1}-z_{2})$-local grading restriction
conditions (or the $L(0)$-semisimple version 
of these conditions when $\mathcal{C}$ is in 
$\mathcal{M}_{sg}$) by $W_{P(z_1, z_2)}$.

We have a linear map $\Psi^{(1)}_{P(z_1, z_2)}$ {}from
\[
W_1\hboxtr_{P(z_1)} (W_2\boxtimes_{P(z_2)} W_3)=(W_1\boxtimes_{P(z_1)}
(W_2\boxtimes_{P(z_2)} W_3))'
\]
to $(W_1\otimes W_2\otimes W_3)^*$ defined by
\begin{eqnarray*}
\lefteqn{\Psi^{(1)}_{P(z_1, z_2)}(\nu)(w_{(1)}\otimes
w_{(2)}\otimes w_{(3)})}\\ 
&&=\langle \nu,
w_{(1)}\boxtimes_{P(z_1)}(w_{(2)}
\boxtimes_{P(z_2)}w_{(3)})\rangle_{W_1\boxtimes_{P(z_1)}(W_2
\boxtimes_{P(z_2)}W_3)}
\end{eqnarray*}
for all $\nu \in W_1\hboxtr_{P(z_1)}(W_2\boxtimes_{P(z_2)}W_3)$,
$w_{(1)}\in W_1$, $w_{(2)}\in W_2$ and $w_{(3)}\in W_3$.

By Proposition \ref{8.12}, the images of
elements of $W_1\hboxtr_{P(z_1)}(W_2\boxtimes_{P(z_2)}W_3)$ under
$\Psi^{(1)}_{P(z_1, z_2)}$ satisfy the $P(z_1, z_2)$-compatibility
condition and the $P(z_1, z_2)$-local grading restriction
condition (or the $L(0)$-semisimple $P(z_1, z_2)$-local grading restriction
condition when $\mathcal{C}$ is in 
$\mathcal{M}_{sg}$). By Proposition \ref{9.7},
they satisfy the
$P^{(1)}(z_2)$-local grading restriction condition
(or the $L(0)$-semisimple $P^{(1)}(z_2)$-local grading restriction
condition when $\mathcal{C}$ is in 
$\mathcal{M}_{sg}$). The expansion 
condition for intertwining operators in $\mathcal{C}$ says that 
they also satisfy the $P^{(2)}(z_{1}-z_{2})$-local grading restriction 
condition (or the $L(0)$-semisimple $P^{(2)}(z_{1}-z_{2})$-local 
grading restriction
condition when $\mathcal{C}$ is in 
$\mathcal{M}_{sg}$). Thus we see that $\Psi^{(1)}_{P(z_1, z_2)}$
in fact maps into $W_{P(z_1, z_2)}$.

On the other hand, we have a linear map $\Psi^{(2)}_{P(z_1, z_2)}$
{}from
\[
(W_1\boxtimes_{P(z_1-z_2)}W_2)\hboxtr_{P(z_2)}W_3
=((W_1\boxtimes_{P(z_1-z_2)}W_2)\boxtimes_{P(z_2)}W_3)'
\]
to $(W_1\otimes W_2\otimes W_3)^{*}$ by
\begin{eqnarray*}
\lefteqn{\Psi^{(2)}_{P(z_1, z_2)}(\nu)(w_{(1)}\otimes w_{(2)}
\otimes w_{(3)})}\\
&&=\langle \nu, (w_{(1)}\boxtimes_{P(z_1-z_2)}w_{(2)})
\boxtimes_{P(z_2)}w_{(3)}\rangle_{(W_1\boxtimes_{P(z_1-z_2)}W_2)
\boxtimes_{P(z_2)}W_3}
\end{eqnarray*}
for all $\nu \in (W_1\boxtimes_{P(z_1z_2)}W_2)\hboxtr_{P(z_2)}W_3$,
$w_{(1)}\in W_1$, $w_{(2)}\in W_2$ and $w_{(3)}\in W_3$.  Similarly,
$\Psi^{(2)}_{P(z_1, z_2)}$ in fact maps into $W_{P(z_1, z_2)}$.

As in \cite{tensor4}, we now show that both of
$\Psi^{(1)}_{P(z_1,z_2)}$ and $\Psi^{(2)}_{P(z_1, z_2)}$ are module
isomorphisms mapping onto $W_{P(z_1, z_2)}$.  We shall do this by
constructing a linear map $(\Psi^{(1)}_{P(z_1, z_2)})^{-1}$ {}from
$W_{P(z_1, z_2)}$ to $W_1\hboxtr_{P(z_1)}(W_2\boxtimes_{P(z_2)}W_3)$
and showing that it is the inverse of $\Psi^{(1)}_{P(z_1, z_2)}$.

For any $w_{(1)}\in W_1$ and any $\lambda\in W_{P(z_1, z_2)}$,
$\mu^{(1)}_{\lambda, w_{(1)}}$ is in
$\overline{W_2\hboxtr_{P(z_2)}W_3}$ by the $P(z_1, z_2)$-compatibility
condition, the $P(z_1, z_2)$-local grading restriction condition and
the $P^{(1)}(z_2)$-local grading restriction condition 
(or the $L(0)$-semisimple version of these two 
conditions when $\mathcal{C}$ is in $\mathcal{M}_{sg}$). We define an
element $(\Psi^{(1)}_{P(z_1, z_2)})^{-1}(\lambda)\in
(W_1\otimes(W_2\boxtimes_{P(z_2)}W_3))^*$ by
\[
(\Psi^{(1)}_{P(z_1, z_2)})^{-1}(\lambda)(w_{(1)}\otimes w)=\bra w,
\mu^{(1)}_{\lambda,
w_{(1)}}\ket_{W_2\shboxtr_{P(z_2)}W_3}
\]
for all $w_{(1)}\in W_1$ and $w\in W_2\boxtimes_{P(z_2)}W_3$. By the
$P(z_1, z_2)$-compatibility condition and the $P(z_1, z_2)$-local
grading restriction condition 
(or the $L(0)$-semisimple $P(z_1, z_2)$-local
grading restriction condition 
when $\mathcal{C}$ is in $\mathcal{M}_{sg}$)
for $\lambda$, $(\Psi^{(1)}_{P(z_1,
z_2)})^{-1}(\lambda)$ is in fact in $W_1\hboxtr_{P(z_1)}
(W_2\boxtimes_{P(z_2)}W_3)$.  Thus $\lambda\mapsto(\Psi^{(1)}_{P(z_1,
z_2)})^{-1}(\lambda)$ defines a map $(\Psi^{(1)}_{P(z_1, z_2)})^{-1}$
{}from $W_{P(z_1, z_2)}$ to $W_1\hboxtr_{P(z_1)}
(W_2\boxtimes_{P(z_2)}W_3)$.

For any $\lambda\in W_{P(z_1, z_2)}$, using the definitions of
$\Psi^{(1)}_{P(z_1, z_2)}$, $(\Psi^{(1)}_{P(z_1, z_2)})^{-1}$ and
$\boxtimes_{P(z_1)}$, we have
\begin{eqnarray*}
\lefteqn{\Psi^{(1)}_{P(z_1, z_2)}
((\Psi^{(1)}_{P(z_1, z_2)})^{-1}(\lambda))(w_{(1)}\otimes
w_{(2)}\otimes w_{(3)})}\\
&&=\sum _{n\in \mathbb{C}}\langle (\Psi^{(1)}_{P(z_1, z_2)})^{-1}(\lambda),
w_{(1)}
\boxtimes_{P(z_1)}
\pi_{n}(w_{(2)}\boxtimes_{P(z_2)} w_{(3)})\rangle_{W_1\boxtimes_{P(z_1)}
(W_2\boxtimes_{P(z_2)}W_3)} \\
&&=\sum _{n\in \mathbb{C}}((\Psi^{(1)}_{P(z_1, z_2)})^{-1}(\lambda))(w_{(1)}
\otimes \pi_{n}(w_{(2)}\boxtimes_{P(z_2)} w_{(3)}))\\
&&=\sum _{n\in \mathbb{C}}\langle
\pi_{n}(w_{(2)}\boxtimes_{P(z_2)} w_{(3)}), \mu^{(1)}_{\lambda, w_{(1)}}
 \rangle_{W_2\shboxtr_{P(z_2)}W_3}\\
&&=\sum _{n\in \mathbb{C}}\langle \pi_{n}(\mu^{(1)}_{\lambda, w_{(1)}}),
w_{(2)}\boxtimes_{P(z_2)}
w_{(3)}\rangle_{W_2\boxtimes_{P(z_2)}W_3}\\
&&=\sum _{n\in \mathbb{C}}(\pi_{n}(\mu^{(1)}_{\lambda, w_{(1)}}))
(w_{(2)}\otimes w_{(3)})\\
&&=\mu^{(1)}_{\lambda, w_{(1)}}
(w_{(2)}\otimes w_{(3)})\\
&&=\lambda(w_{(1)}\otimes w_{(2)}\otimes w_{(3)}),
\end{eqnarray*}
proving
\[
\Psi^{(1)}_{P(z_1, z_2)}(\Psi^{(1)}_{P(z_1, z_2)})^{-1}=1.
\]

Next we want to show that $\Psi^{(1)}_{P(z_1, z_2)}$ is injective.
Let $\nu_1, \nu_2$ be two elements of $W_1\hboxtr_{P(z_1)}
(W_2\boxtimes_{P(z_2)}W_3)$ such that
\[
\Psi^{(1)}_{P(z_1, z_2)}(\nu_1)=\Psi^{(1)}_{P(z_1, z_2)}(\nu_2),
\]
that is,
\begin{eqnarray*}
\lefteqn{\langle \nu_1, w_{(1)}\boxtimes_{P(z_1)}(w_{(2)}
\boxtimes_{P(z_2)}w_{(3)})\rangle_{W_1\boxtimes_{P(z_1)}
(W_2\boxtimes_{P(z_2)}W_3)}}\\
&&=\langle \nu_2, w_{(1)}\boxtimes_{P(z_1)}(w_{(2)}
\boxtimes_{P(z_2)}w_{(3)})\rangle_{W_1\boxtimes_{P(z_1)}
(W_2\boxtimes_{P(z_2)}W_3)}
\end{eqnarray*}
for all $w_{(1)}\in W_1$, $w_{(2)}\in W_2$ and $w_{(3)}\in W_3$.  By
Corollary \ref{prospan} the homogeneous components of all
$w_{(1)}\boxtimes_{P(z_1)}(w_{(2)} \boxtimes_{P(z_2)}w_{(3)})$
linearly span the space $W_1\boxtimes_{P(z_1)} (W_2\boxtimes_{P(z_2)}
W_3)$, hence we must have $\nu_1 =\nu_2$, and the injectivity of
$\Psi^{(1)}_{P(z_1, z_2)}$ follows.

Thus we have an isomorphism:
\[
(\Psi^{(1)}_{P(z_1,z_2)})^{-1} \circ \Psi^{(2)}_{P(z_1, z_2)}:
(W_1\boxtimes_{P(z_1-z_2)}W_2)\hboxtr_{P(z_2)}W_3
\longrightarrow W_1\hboxtr_{P(z_1)} (W_2\boxtimes_{P(z_2)} W_3),
\]
Its contragredient map hence gives a module isomorphism
(\ref{assoc-iso})
such that (\ref{assoc-elt-1}) holds for 
$w_{(1)}\in W_1$, $w_{(2)}\in W_2$ and
$w_{(3)}\in W_3$.

The uniqueness of $\mathcal{A}_{P(z_{1}), P(z_{2})}^{P(z_{1}-z_{2}), P(z_{2})}$ 
follows {}from Corollary \ref{prospan}.
\epfv

\begin{rema}{\rm
We also have the inverse {\it associativity isomorphisms}
\begin{equation}\label{assoc-iso-inv}
\alpha_{P(z_{1}), P(z_{2})}^{P(z_{1}-z_{2}), P(z_{2})}: 
(W_1\boxtimes_{P(z_1-z_2)}
W_2)\boxtimes_{P(z_2)} W_3 \to W_1\boxtimes_{P(z_1)}
(W_2\boxtimes_{P(z_2)} W_3),
\end{equation}
such that 
\begin{equation}\label{assoc-elt-2}
\overline{\alpha_{P(z_{1}), P(z_{2})}^{P(z_{1}-z_{2}), P(z_{2})}}
((w_{(1)}\boxtimes_{P(z_1-z_2)}
w_{(2)})\boxtimes_{P(z_2)} w_{(3)} = w_{(1)}\boxtimes_{P(z_1)}
(w_{(2)}\boxtimes_{P(z_2)} w_{(3)})
\end{equation}
for  $w_{(1)}\in W_1$, $w_{(2)}\in W_2$ and
$w_{(3)}\in W_3$, and (\ref{assoc-elt-2})
determines the associativity isomorphism (\ref{assoc-iso-inv})
uniquely.
}
\end{rema}

\newpage

\setcounter{equation}{0}
\setcounter{rema}{0}

\section{The convergence and extension properties 
and differential equations}

In the construction of the associativity isomorphisms we have needed,
and assumed, the convergence and expansion conditions for intertwining
maps in ${\cal C}$.  In this section we will follow \cite{tensor4} to
give certain sufficient conditions for a category ${\cal C}$ to have
these properties. In Subsection 11.1, we give what we call, as in
\cite{tensor4} the ``convergence and extension properties'' and show
that they imply the convergence and expansion conditions for
intertwining maps in ${\cal C}$. In Subsection 11.2, we show that the
proofs in \cite{diff-eqn} can be adapted to generalize the results of
\cite{diff-eqn} to results in the logarithmic generality. In
particular, we see that two purely algebraic conditions, the
``$C_{1}$-cofiniteness condition'' and the
``quasi-finite-dimensionality condition,'' for all objects of
$\mathcal{C}$ imply the convergence and extension properties and thus
also imply the convergence and expansion conditions for intertwining
maps in ${\cal C}$.

In addition to Assumption \ref{assum-assoc}, we assume the following
in this section:

\begin{assum}
In this section, we assume that all generalized $V$-modules in
$\mathcal{C}$ are $\R$-graded.
\end{assum}

\subsection{The convergence and extension properties}

Given objects $W_1$, $W_2$, $W_3$,
$W_4$, $M_1$ and $M_2$ of the category ${\cal C}$, let ${\cal Y}_1$,
${\cal Y}_2$, ${\cal Y}^1$ and ${\cal Y}^2$ be intertwining operators
of type ${W_4}\choose {W_1M_1}$, ${M_1}\choose {W_2W_3}$,
${W_4}\choose {M_2W_3}$ and ${M_2}\choose {W_1W_2}$,
respectively. For our purpose, it is sufficient to consider only the 
case that $W_1$, $W_2$, $W_3$,
$W_4$, $M_1$ and $M_2$ are
indecomposable generalized $V$-modules. Recall {}from Remark
\ref{congruent} that all of the generalized weights of such a module
are congruent modulo ${\mathbb Z}$ to one another.
Consider the following conditions on the product of
${\cal Y}_1$ and ${\cal Y}_2$ and on the iterate of ${\cal Y}^1$ and
${\cal Y}^2$, respectively:

\begin{description}

\item[Convergence and extension property for products] There exists an
integer $N$ depending only on ${\cal Y}_1$ and ${\cal Y}_2$, and for
any $w_{(1)}\in W_1$, $w_{(2)}\in W_2$, $w_{(3)}\in W_3$, $w'_{(4)}\in
W'_4$, there exist $M\in{\mathbb N}$, $r_{k}, s_{k}\in {\mathbb R}$, $i_{k},
j_{k}\in {\mathbb N}$, $k=1,\dots,M$ and analytic functions
$f_{i_{k}j_{k}}(z)$ on $|z|<1$, $k=1, \dots, M$, satisfying
\begin{equation}\label{c-e-p-1}
\wt w_{(1)}+\wt w_{(2)}+s_{k}>N, \;\;\;k=1, \dots, M,
\end{equation}
such that
\begin{equation}\label{c-e-p-2}
\langle w'_{(4)}, {\cal Y}_1(w_{(1)}, x_2) {\cal Y}_2(w_{(2)},
x_2)w_{(3)}\rangle_{W_4} \lbar_{x_1= z_1, \;x_2=z_2}
\end{equation}
is convergent when $|z_1|>|z_2|>0$ and can be analytically extended to
the multi-valued analytic function
\begin{equation}\label{c-e-p-3}
\sum_{k=1}^{M}z_2^{r_{k}}(z_1-z_2)^{s_{k}}(\log z_2)^{i_{k}}
(\log(z_1-z_2))^{j_{k}}f_{i_{k}j_{k}}\left(\frac{z_1-z_2}{z_2}\right)
\end{equation}
(here $\log (z_{1}-z_{2})$ and $\log z_{2}$ mean the multivalued functions,
not the particular branch we have been using) 
in the region $|z_2|>|z_1-z_2|>0$.

\item[Convergence and extension property without logarithms for
products] When $i_{k}=j_{k}=0$ for $k=1, \dots, M$, we call the
property above the {\it convergence and extension property without
logarithms for products}.

\item[Convergence and extension property for iterates] There exists an
integer $\tilde{N}$ depending only on ${\cal Y}^1$ and ${\cal Y}^2$,
and for any $w_{(1)}\in W_1$, $w_{(2)}\in W_2$, $w_{(3)}\in W_3$,
$w'_{(4)}\in W'_4$, there exist $\tilde M\in{\mathbb N}$, $\tilde{r}_{k},
\tilde{s}_{k}\in {\mathbb R}$, $\tilde{i}_{k}, \tilde{j}_{k}\in {\mathbb
N}$, $k=1,\dots,\tilde M$ and analytic functions
$\tilde{f}_{\tilde{i}_{k}\tilde{j}_{k}}(z)$ on $|z|<1$, $k=1, \dots,
\tilde{M}$, satisfying
\[
\wt w_{(2)}+\wt w_{(3)}+\tilde{s}_{k}>\tilde{N}, \;\;\;k=1, \dots,
\tilde{M},
\]
 such that
\[
\langle w'_{(4)}, {\cal Y}^1({\cal Y}^2(w_{(1)}, x_0)w_{(2)},
x_2)w_{(3)}\rangle_{W_4} \lbar_{x_0=z_1-z_2,\;x_2=z_2}
\]
is convergent when $|z_2|>|z_1-z_2|>0$ and can be analytically
extended to the multi-valued analytic function
\[
\sum_{k=1}^{\tilde{M}} z_1^{\tilde{r}_{k}}z_2^{\tilde{s}_{k}} (\log
z_1)^{\tilde{i}_{k}} (\log
z)^{\tilde{j}_{k}}\tilde{f}_{\tilde{i}_{k}\tilde{j}_{k}}
\left(\frac{z_2}{z_1}\right)
\]
(here $\log z_{1}$ and $\log z_{2}$ mean the multivalued functions,
not the particular branches we have been using) 
in the region $|z_1|>|z_2|>0$.

\item[Convergence and extension property without logarithms 
for iterates] When $i_{k}=j_{k}=0$ for $k=1, \dots, M$, 
we call the property above the {\it convergence and extension 
property without logarithms 
for iterates}.

\end{description}

If for any objects $W_1$, $W_2$, $W_3$, $W_4$ and $M_1$ of ${\cal C}$
and any intertwining operators ${\cal Y}_1$ and ${\cal Y}_2$ of the
types as above, the convergence and extension property, with or without 
logarithms, for products
holds, we say that {\it the (corresponding) convergence and extension property 
for products holds in ${\cal C}$}.  We similarly define the meaning of the
phrase {\it the (corresponding) convergence and extension property 
for iterates holds in ${\cal C}$}.

\begin{rema}\label{power-wt}
{\rm 
If the convergence and extension property for products 
holds, then we can always find $M$,
$r_{k}, s_{k}$, $i_{k}$, $j_{k}$
and $f_{i_{k}j_{k}}(z)$, $k=1, \dots, M$, such that 
\begin{equation}\label{power-wt-1}
r_{k}+s_{k}=\Delta,\;\;\;k=1, \dots, M
\end{equation}
where
\begin{equation}
\Delta=-\wt w_{(1)}-\wt w_{(2)} -\wt w_{(3)}+\wt w'_{(4)}.
\end{equation}
Similarly, if the convergence and extension property for iterates
holds, then
we can always find
$\tilde{M}$, 
$\tilde{r}_{k}, \tilde{s}_{k}$  and 
$\tilde{f}_{\tilde{i}_{k}\tilde{j}_{k}}$, $k=1, \dots, \tilde{M}$,
such that
\begin{equation}
\tilde{r}_{k}+\tilde{s}_{k}=\Delta,\;\;\;k=1, \dots, \tilde{M}.
\end{equation}}
\end{rema}

We also need the following concept: If a generalized $V$-module
$W=\coprod_{n\in {\mathbb C}}W_{[n]}$ satisfies the condition that
$W_{[n]}=0$ for $n$ with real part sufficiently negative, we say that $W$
is a {\it lower-truncated}.

The first result of this section is:

\begin{theo}\label{thm-11.1}
Suppose the following two conditions are satisfied:
\begin{enumerate}
\item Every finitely-generated lower truncated generalized $V$-module
is an object of ${\cal C}$.

\item The convergence and extension property for either products or
iterates holds in ${\cal C}$ (or the 
convergence and extension property without logarithms for either products or
iterates holds in ${\cal C}$ when ${\cal C}$ is in $\mathcal{M}_{sg}$).
\end{enumerate}
Then the convergence  and the expansion
conditions for intertwining maps 
in ${\cal C}$ both hold (recall Proposition
\ref{convergence} and Theorem \ref{expansion}).
\end{theo}
\pf By the convergence and extension property, the convergence
condition for intertwining maps in ${\cal C}$ holds. By Theorem
\ref{expansion}, we need only prove that the first property in
Theorem \ref{expansion} holds, that is, we need only prove that
for any objects $W_1$, $W_2$, $W_3$, $W_4$ and $M_1$ of ${\cal C}$,
any nonzero complex numbers $z_1$ and $z_2$ satisfying
$|z_1|>|z_2|>|z_1-z_2|>0$, any $P(z_1)$-intertwining map $I_1$ of type
${W_4}\choose {W_1 M_1}$ and $P(z_2)$-intertwining map $I_2$ of type
${M_1}\choose {W_2W_3}$ and any $w'_{(4)}\in W'_4$, $(I_1\circ
(1_{W_1}\otimes I_2))'(w'_{(4)})\in (W_1\otimes W_2\otimes W_3)^{*}$
satisfies the $P^{(2)}(z_1-z_2)$-local grading restriction condition.

By assumption and Remark \ref{power-wt}, for
any $w_{(1)}\in W_1$, $w_{(2)}\in W_2$, $w_{(3)}\in W_3$, 
there exist $M\in{\mathbb N}$, $r_{k}, s_{k}\in {\mathbb R}$, $i_{k},
j_{k}\in {\mathbb N}$, $k=1,\dots,M$ and analytic functions
$f_{i_{k}j_{k}}(z)$ on $|z|<1$, $k=1, \dots, M$,
such that (\ref{power-wt-1}) holds and
(\ref{c-e-p-2})
is absolutely convergent when $|z_1|>|z_2|>0$ and 
can be analytically extended to
the multi-valued analytic function
(\ref{c-e-p-3})
in the region $|z_2|>|z_1-z_2|>0$. Then we can always find
$f_{i_{k}j_{k}}(z)$ for $k=1, \dots, M$ such that 
(\ref{c-e-p-2}) is equal to (\ref{c-e-p-3}) when we 
choose the values of $\log z_{2}$ and $\log (z_{1}-z_{2})$ 
to be $|\log (z_{1}-z_{2})|+i\arg z_{2}$
and $|\log (z_{1}-z_{2})|+i\arg (z_{1}-z_{2})$ where 
$0\le \arg z_{2}, \arg (z_{1}-z_{2})<2\pi$ and 
the values of $z_{2}^{r_{k}}$  and $(z_{1}-z_{2})^{s_{k}}$
to be $e^{r_{k}\log z_{2}}$ and $e^{s_{k}\log (z_{1}-z_{2})}$.
Expanding $f_{i_{k}j_{k}}(z)$, $k=1, \dots, M$, we can write 
(\ref{c-e-p-3}) as 
\begin{eqnarray}\label{thm-11.1-1}
\lefteqn{\sum_{k=1}^{M}\sum_{m\in \mathbb{N}}
C_{km}(w'_{(4)}, w_{(1)}, w_{(2)}, w_{(3)})\cdot}\nn
&&\quad\quad
 \cdot e^{(r_{k}-m)\log z_{2}}e^{(s_{k}+m)
\log (z_{1}-z_{2})}(\log z_{2})^{i_{k}}(\log (z_{1}-z_{2}))^{j_{k}}.
\end{eqnarray}
For $n\in \mathbb{C}$ and $k=1, \dots, M$, let 
\begin{equation}\label{thm-11.1-2}
a_{n;k}(w'_{(4)}, w_{(1)}, w_{(2)}, w_{(3)})=
\sum_{-r_{k}+m-1=n}
C_{km}(w'_{(4)}, w_{(1)}, w_{(2)}, w_{(3)}).
\end{equation}
Then (\ref{thm-11.1-1}) can be written as 
\begin{eqnarray}\label{thm-11.1-3}
\lefteqn{\sum_{k=1}^{M}\sum_{n\in \mathbb{C}}a_{n;k}(w'_{(4)}, 
w_{(1)},
w_{(2)}, w_{(3)})\cdot}\nn
&&\quad\quad
 \cdot e^{(\Delta+n+1)\log (z_{1}-z_{2})}e^{(-n-1)\log z_{2}}
(\log z_{2})^{i_{k}}(\log (z_{1}-z_{2}))^{j_{k}}.
\end{eqnarray}
{}From (\ref{power-wt-1}), we see that if  $n\in \mathbb{C}$ satisfies
\begin{equation}\label{thm-11.1-4}
n\ne -r_{k}+m-1=-\Delta+s_{k}+m-1
\end{equation}
for all $m\in \mathbb{N}$, 
then 
\begin{equation}\label{thm-11.1-5}
a_{n; k}(w'_{(4)}, w_{(1)},
w_{(2)}, w_{(3)})=0. 
\end{equation}
Since  $m\ge 0$, by (\ref{c-e-p-1}), (\ref{power-wt-1}) and 
(\ref{thm-11.1-4}),
we see that for $n\in \mathbb{C}$, if
\begin{equation}\label{thm-11.1-6}
n+1+ \wt w'_{(4)}-
 \wt w_{(3)}\le N,
\end{equation}
(\ref{thm-11.1-5}) holds.

Since $|z_{1}|>|z_{2}|>|z_{1}-z_{2}|>0$, we know that
\begin{eqnarray}\label{thm-11.1-7}
\lefteqn{\sum_{k=1}^{M}\sum_{n\in \mathbb{C}}a_{n; k}(w'_{(4)}, 
w_{(1)},
w_{(2)}, w_{(3)})\cdot}\nn
&&\quad
 \cdot e^{(-\Delta+n+1)\log (z_{1}-z_{2})}e^{(-n-1)\log z_{2}}
(\log z_{2})^{i_{k}}(\log (z_{1}-z_{2}))^{j_{k}}
\end{eqnarray}
converges absolutely to (\ref{c-e-p-1}). 
For $w'_{(4)}\in W'_{4}$ and $w_{(3)}\in W_{3}$, 
let $\beta_{n;k}(w'_{(4)}, w_{(3)})\in (W_{1}\otimes W_{2})^{*}$ be defined by
\begin{equation}\label{thm-11.1-8}
\beta_{n;k}(w'_{(4)}, w_{(3)})(w_{(1)}\otimes w_{(2)})
=a_{n;k}(w'_{(4)}, w_{(1)},
w_{(2)}, w_{(3)})
\end{equation}
for all $w_{(1)}\in W_{1}$ and $w_{(2)}\in W_{2}$.
By definition the series 
\[
\sum_{k=1}^{M}\sum_{n\in \mathbb{C}}\beta_{n;k}
(w'_{(4)}, w_{(3)})(w_{(1)}\otimes w_{(2)})
e^{(\Delta+n+1)\log (z_{1}-z_{2})}e^{(-n-1)\log z_{2}}
(\log z_{2})^{i_{k}}(\log (z_{1}-z_{2}))^{j_{k}}
\]
is absolutely 
convergent to 
$\mu^{(2)}_{I_{1}\circ (1_{W_{2}}\otimes I_{2}))'(w'_{(4)}), w_{(3)}}
(w_{(1)}\otimes w_{(2)})$ for $w_{(1)}\in W_{1}$ and $w_{(2)}\in W_{2}$.
To show that $(I_{1}\circ (1_{W_{2}}\otimes I_{2}))'(w'_{(4)})$ satisfies the 
$P^{(2)}(z_{1}-z_{2})$-local grading restriction
condition,
we need to calculate
\[
(v_{1})_{m_{1}}\cdots (v_{r})_{m_{r}}\beta_{n;k}(w'_{(4)}, 
w_{(3)})
\]
and its weight for any $r\in \mathbb{N}$, $v_{1}, \dots, v_{r}\in V,$ 
 $m_{1}, \dots, m_{r}\in \mathbb{Z}$,
$n\in \mathbb{C}$, where $(v_{1})_{m_{1}}, \dots,  (v_{r})_{m_{r}}$, $m_{1}, 
\cdots, m_{r}\in \mathbb{Z}$, are the components of
$Y'_{P(z_{1}-z_{2})}(v_{1}, x)$, $\dots\ $, $Y'_{P(z_{1}-z_{2})}(v_{r}, x)$,
respectively,  on $(W_{1}\otimes W_{2})^{*}$. As in \cite{tensor4},
for convenience, we
instead calculate 
\[
(v^{o}_{1})_{m_{1}}\cdots (v^{o}_{r})_{m_{r}}
\beta_{n;k}(w'_{(4)}, w_{(3)}),
\]
where 
$(v^{o}_{1})_{m_{1}}, \dots,  (v^{o}_{r})_{m_{r}}$, 
$m_{1}, \cdots, m_{r}\in \mathbb{Z}$,  are the components of the opposite
vertex operators
$Y^{\prime\; o}_{P(z_{1}-z_{2})}(v_{1}, x)$, $\dots\ $, 
$Y^{\prime\; o}_{P(z_{1}-z_{2})}(v_{r}, x)$, respectively.

By the definition (\ref{Y'def}) of
$Y'_{P(z_{1}-z_{2})}(v, x)$ and 
$$Y^{\prime\; o}_{P(z_{1}-z_{2})}(v, x)=Y'_{P(z_{1}-z_{2})}
(e^{xL(1)}(-x^{-2})^{L(0)}v, x^{-1}),$$ we have
\begin{eqnarray}\label{16.19}
\lefteqn{(Y^{\prime \; o}_{P(z_{1}-z_{2})}(v, x)
\beta_{n; k}(w'_{(4)}, w_{(3)}))(w_{(1)}\otimes w_{(2)})}\nno\\
&&=(\beta_{n;k}(w'_{(4)}, w_{(3)}))(w_{(1)}\otimes Y_{2}(v, x)
w_{(2)})
\nno\\
&&\quad +\res_{x_{0}}(z_{1}-z_{2})^{-1}\delta\left(\frac{x-x_{0}}{z_{1}-z_{2}}
\right)
(\beta_{n;k}(w'_{(4)}, w_{(3)}))(Y_{1}(v, x_{0})
w_{(1)}\otimes w_{(2)})
\nno\\
&&=a_{n;k}(w'_{(4)}, w_{(1)}, Y_{2}(v, x)
w_{(2)}, w_{(3)})\nno\\
&&\quad +\res_{x_{0}}(z_{1}-z_{2})^{-1}\delta\left(\frac{x-x_{0}}{z_{1}-z_{2}}
\right)
a_{n;k}(w'_{(4)}, Y_{1}(v, x_{0})
w_{(1)}, w_{(2)}, w_{(3)}).
\end{eqnarray}
On the other hand, since (\ref{thm-11.1-7}) is absolutely 
convergent to 
(\ref{c-e-p-1})
when $|z_{1}|>|z_{2}|>|z_{1}-z_{2}|>0$, we have
\begin{eqnarray}\label{16.20}
\lefteqn{\sum_{k=1}^{M}\sum_{n\in {\Bbb C}}a_{n;k}(w'_{(4)}, 
w_{(1)}, Y_{2}(v, x)w_{(2)}, w_{(3)}) }\nno\\
&& \quad \quad\quad\quad\cdot e^{(\Delta+n+1)\log (z_{1}-z_{2})}
e^{(-n-1)\log z_{2}}(\log z_{2})^{i_{k}}(\log (z_{1}-z_{2}))^{j_{k}}\nno\\
&&+ \res_{x_{0}}(z_{1}-z_{2})^{-1}\delta\left(\frac{x-x_{0}}{z_{1}-z_{2}}
\right)\sum_{k=1}^{M}\sum_{n\in {\Bbb C}}a_{n;k}(w'_{(4)}, Y_{1}(v, x_{0})
w_{(1)}, w_{(2)}, w_{(3)})\cdot\nno\\
&&\quad \quad\quad\quad\cdot 
e^{(\Delta+n+1)\log (z_{1}-z_{2})}e^{(-n-1)\log z_{2}}
(\log z_{2})^{i_{k}}(\log (z_{1}-z_{2}))^{j_{k}}\nno\\
&&=\langle w'_{(4)}, {\cal Y}_{1}(w_{(1)}, x_{1})
{\cal Y}_{2}(Y_{2}(v, x)w_{(2)}, x_{2})w_{(3)}\rangle_{W_{4}}
\lbar_{x_{1}
=z_{1}, x_{2}=z_{2}}\nno\\
&&\quad +\res_{x_{0}}(x_{1}-x_{2})^{-1}\delta\left(\frac{x-x_{0}}{x_{1}-x_{2}}
\right)\cdot\nno\\
&&\quad\quad\quad\quad\cdot
\langle w'_{(4)}, {\cal Y}_{1}(Y_{1}(v, x_{0})w_{(1)},
x_{1})
{\cal Y}_{2}(w_{(2)}, x_{2})w_{(3)}\rangle_{W_{4}}
\lbar_{x_{1}
=z_{1}, x_{2}=z_{2}}\nno\\
&&=\langle w'_{(4)}, {\cal Y}_{1}(w_{(1)}, x_{1})
{\cal Y}_{2}(Y_{2}(v, x)w_{(2)}, x_{2})w_{(3)}\rangle_{W_{4}}
\lbar_{x_{1}
=z_{1}, x_{2}=z_{2}}\nno\\
&&\quad +\res_{x_{0}}x_{1}^{-1}\delta\left(\frac{(x+x_{2})-x_{0}}{x_{1}}
\right)\cdot\nno\\
&&\quad\quad\quad\quad\cdot
\langle w'_{(4)}, {\cal Y}_{1}(Y_{1}(v, x_{0})w_{(1)},
x_{1})
{\cal Y}_{2}(w_{(2)}, x_{2})w_{(3)}\rangle_{W_{4}}
\lbar_{x_{1}
=z_{1}, x_{2}=z_{2}}.
\end{eqnarray}

Using the Jacobi identity for the logarithmic 
intertwining operators ${\cal Y}_{1}$ and 
${\cal Y}_{2}$ and the properties of the formal $\delta$-function,  
the right-hand side of (\ref{16.20}) is equal to
\begin{eqnarray}\label{16.21}
\lefteqn{\res_{y}x^{-1}\delta\left(\frac{y-x_{2}}{x}\right)
\langle w'_{(4)}, {\cal Y}_{1}(w_{(1)}, x_{1})Y_{5}(v, y)
{\cal Y}_{2}(w_{(2)}, x_{2})w_{(3)}\rangle_{W_{4}}
\lbar_{x_{1}
=z_{1}, x_{2}=z_{2}}}\nno\\
&&- \res_{y}x^{-1}\delta\left(\frac{x_{2}-y}{-x}\right)
\langle w'_{(4)}, {\cal Y}_{1}(w_{(1)}, x_{1})
{\cal Y}_{2}(w_{(2)}, x_{2})Y_{3}(v, y)w_{(3)}\rangle_{W_{4}}
\lbar_{x_{1}
=z_{1}, x_{2}=z_{2}}\nno\\
&&+\langle w'_{(4)}, Y_{4}(v, x+x_{2}){\cal Y}_{1}(w_{(1)},
x_{1})
{\cal Y}_{2}(w_{(2)}, x_{2})w_{(3)}\rangle_{W_{4}}
\lbar_{x_{1}
=z_{1}, x_{2}=z_{2}}
\nno\\
&&-\langle w'_{(4)}, {\cal Y}_{1}(w_{(1)},
x_{1})Y_{5}(v, x+x_{2})
{\cal Y}_{2}(w_{(2)}, x_{2})w_{(3)}\rangle_{W_{4}}
\lbar_{x_{1}
=z_{1}, x_{2}=z_{2}}
\nno\\
&&=\langle w'_{(4)}, Y_{4}(v, x+x_{2}){\cal Y}_{1}(w_{(1)},
x_{1})
{\cal Y}_{2}(w_{(2)}, x_{2})w_{(3)}\rangle_{W_{4}}
\lbar_{x_{1}
=z_{1}, x_{2}=z_{2}}\nno\\
&&\quad - \res_{y}x^{-1}\delta\left(\frac{x_{2}-y}{-x}\right)
\langle w'_{(4)}, {\cal Y}_{1}(w_{(1)}, x_{1})
{\cal Y}_{2}(w_{(2)}, x_{2})Y_{3}(v, y)w_{(3)}\rangle_{W_{4}}
\lbar_{x_{1}
=z_{1}, x_{2}=z_{2}}\nno\\
&&=\res_{y}x^{-1}\delta\left(\frac{y-x_{2}}{x}\right)
\langle w'_{(4)}, Y_{4}(v, y){\cal Y}_{1}(w_{(1)},
x_{1})
{\cal Y}_{2}(w_{(2)}, x_{2})w_{(3)}\rangle_{W_{4}}
\lbar_{x_{1}
=z_{1}, x_{2}=z_{2}}\nno\\
&&\quad - \res_{y}x^{-1}\delta\left(\frac{x_{2}-y}{-x}\right)
\langle w'_{(4)}, {\cal Y}_{1}(w_{(1)}, x_{1})
{\cal Y}_{2}(w_{(2)}, x_{2})Y_{3}(v, y)w_{(3)}\rangle_{W_{4}}
\lbar_{x_{1}
=z_{1}, x_{2}=z_{2}}.\nno\\
\end{eqnarray}
{}From (\ref{16.20}) and (\ref{16.21}), we obtain
\begin{eqnarray}\label{16.22}
\lefteqn{\sum_{k=1}^{M}\sum_{n\in {\Bbb C}}a_{n;k}(w'_{(4)}, 
w_{(1)}, Y_{2}(v, x)w_{(2)}, w_{(3)}) }\nno\\
&& \quad \quad\quad\quad\cdot e^{(\Delta+n+1)\log (z_{1}-z_{2})}
e^{(-n-1)\log z_{2}}(\log z_{2})^{i_{k}}(\log (z_{1}-z_{2}))^{j_{k}}\nno\\
&&+ \res_{x_{0}}(z_{1}-z_{2})^{-1}\delta\left(\frac{x-x_{0}}{z_{1}-z_{2}}
\right)\sum_{k=1}^{M}\sum_{n\in {\Bbb C}}a_{n;k}(w'_{(4)}, Y_{1}(v, x_{0})
w_{(1)}, w_{(2)}, w_{(3)})\cdot\nno\\
&&\quad \quad\quad\quad\cdot 
e^{(\Delta+n+1)\log (z_{1}-z_{2})}e^{(-n-1)\log z_{2}}
(\log z_{2})^{i_{k}}(\log (z_{1}-z_{2}))^{j_{k}}\nno\\
&&=\res_{y}x^{-1}\delta\left(\frac{y-x_{2}}{x}\right)
\langle w'_{(4)}, Y(v, y){\cal Y}_{1}(w_{(1)}, x_{1})
{\cal Y}_{2}(w_{(2)}, x_{2})w_{(3)}\rangle_{W_{4}}
\lbar_{x_{1}
=z_{1}, x_{2}=z_{2}}\nno\\
&&\quad - \res_{y}x^{-1}\delta\left(\frac{x_{2}-y}{x}\right)
\langle w'_{(4)}, {\cal Y}_{1}(w_{(1)}, x_{1})
{\cal Y}_{2}(w_{(2)}, x_{2})Y(v, y)w_{(3)}\rangle_{W_{4}}
\lbar_{x_{1}
=z_{1}, x_{2}=z_{2}}\nno\\
&&=\sum_{m\in {\Bbb Z}}\sum_{l\ge 0}
(-1)^{l}{{m}\choose {l}}x^{-m-1}x_{2}^{l}\langle w_{(4)}, v_{m-l}
{\cal Y}_{1}(w_{(1)}, x_{1})
{\cal Y}_{2}(w_{(2)}, x_{2})w_{(3)}\rangle_{W_{4}}
\lbar_{x_{1}
=z_{1}, x_{2}=z_{2}}\nno\\
&&\quad -\sum_{m\in {\Bbb Z}}\sum_{ l\ge 0}
(-1)^{l+m}{{m}\choose {l}}x^{-m-1}x_{2}^{m-l}\langle w_{(4)},
{\cal Y}_{1}(w_{(1)}, x_{1})
{\cal Y}_{2}(w_{(2)}, x_{2})v_{l}
w_{(3)}\rangle_{W_{4}}\lbar_{x_{1}
=z_{1}, x_{2}=z_{2}}\nno\\
&&=\sum_{m\in {\Bbb Z}}\sum_{l\ge 0}
(-1)^{l}{{m}\choose {l}}x^{-m-1}\sum_{k=1}^{M}
\sum_{n\in {\Bbb C}}a_{n;k}(v^{*}_{m-l}w'_{(4)},
w_{(1)}, w_{(2)}, w_{(3)})\cdot\nno\\
&&\quad \quad\quad\quad\cdot e^{(\Delta+n+1)
\log (z_{1}-z_{2})}e^{(l-n-1)\log z_{2}}
(\log z_{2})^{i_{k}}(\log (z_{1}-z_{2}))^{j_{k}}\nno\\
&&\quad -\sum_{m\in {\Bbb Z}}\sum_{ l\ge 0}
(-1)^{l+m}{{m}\choose {l}}x^{-m-1}
\sum_{n\in {\Bbb C}}a_{n;k}(w'_{(4)},
w_{(1)}, w_{(2)}, v_{l}w_{(3)})\cdot\nno\\
&&\quad \quad\quad\quad
\cdot e^{(\Delta+n+1)\log (z_{1}-z_{2})}e^{(-m+l-n-1)\log z_{2}}
(\log z_{2})^{i_{k}}(\log (z_{1}-z_{2}))^{j_{k}}.
\end{eqnarray}
The intermediate steps in the equality (\ref{16.22}) hold only when 
$|z_{1}|>|z_{2}|>|z_{1}-z_{2}|>0$. But since
both sides of (\ref{16.22}) can be analytically extended to the 
region $|z_{2}|>|z_{1}-z_{2}|>0$, 
(\ref{16.22}) must hold when 
$|z_{2}|>|z_{1}-z_{2}|>0$. By Assumption \ref{assum-exp-set}, 
the set of $n$'s such that the expansion coefficients of the terms
in (\ref{16.22}) are nonzero together with the set $\{i_{1}, 
\dots, i_{M}\}$ form a unique expansion set. 
Thus the coefficients of both sides of
(\ref{16.22}) in powers of $e^{\log z_{2}}$, $\log z_{2}$ and 
$\log (z_{1}-z_{2})$
are equal, that is,
\begin{eqnarray}\label{16.23}
\lefteqn{a_{n;k}(w'_{(4)}, 
w_{(1)}, Y_{2}(v, x)w_{(2)}, w_{(3)})}
\nno\\
&&+\res_{y}(z_{1}-z_{2})^{-1}\delta\left(\frac{x-y}{z_{1}-z_{2}}
\right) a_{n;k}(w'_{(4)}, Y_{1}(v, y)
w_{(1)}, w_{(2)}, w_{(3)})\nno\\
&&=\sum_{m\in {\Bbb Z}}\sum_{l\ge 0}
(-1)^{l}{{m}\choose {l}}x^{-m-1}
a_{n+l;k}(v^{*}_{m-l}w'_{(4)},
w_{(1)}, w_{(2)}, w_{(3)})(z_{1}-z_{2})^{l}\nno\\
&&\quad -\sum_{m\in {\Bbb Z}}\sum_{ l\ge 0}
(-1)^{l+m}{{m}\choose {l}}x^{-m-1}
a_{n+m-l;k}(w'_{(4)},
w_{(1)}, w_{(2)}, v_{l}w_{(3)})(z_{1}-z_{2})^{m-l}.\nno\\
&&
\end{eqnarray}

By (\ref{16.19}) and (\ref{16.23}), we obtain
\begin{eqnarray}
\lefteqn{(v^{*}_{m}\beta_{n}(w'_{(4)}, w_{(3)}))(w_{(1)}
\otimes w_{(2)})}\nno\\
&&=\sum_{l\ge 0}
(-1)^{l}{{m}\choose {l}}(\beta_{n+l}(v^{*}_{m-l}w'_{(4)}, w_{(3)}))
(w_{(1)}, w_{(2)})(z_{1}-z_{2})^{l}\nno\\
&&\quad -\sum_{ l\ge 0}
(-1)^{l+m}{{m}\choose {l}}(\beta_{n+m-l}(w'_{(4)}, v_{l}w_{(3)}))(w_{(1)},
w_{(2)})(z_{1}-z_{2})^{m-l}.
\end{eqnarray}
By induction, we obtain
\begin{eqnarray}\label{thm-11.1-9}
\lefteqn{((v^{*}_{1})_{m_{1}}\cdots (v^{*}_{r})_{m_{r}}\beta_{n;k}(w'_{(4)}, 
w_{(3)})
(w_{1)}\otimes w_{(2)})=}\nno\\
&&=\sum_{i\in \mathbb{N}}\sum_{\begin{array}{c}
\mbox{\scriptsize $j_{1}>\cdots >j_{i}$}\\
\mbox{\scriptsize $j_{i+1}>\cdots >j_{r}$}
\\ 
\mbox{\scriptsize $\{j_{1}, \dots, j_{r}\}=\{1, \dots, r\}$}\end{array}}
\sum_{l_{1},  \dots, l_{r}\ge 0}(-1)^{l_{1}+\cdots +l_{r}+(m_{j_{i+1}}+1)+
\cdots
+ (m_{j_{r}}+1)}
{{m_{j_{1}}}\choose {l_{1}}}\cdots {{m_{j_{r}}}\choose {l_{r}}}\cdot\nno\\
&&\hspace{2em}\cdot
(\beta_{n+m_{j_{i+1}}\cdots +m_{j_{r}}+l_{1}+\cdots +l_{i}-l_{i+1}-
\cdots -l_{r}; k}\nno\\
&&\hspace{4em}
((v^{*}_{j_{1}})_{m_{j_{1}}-l_{1}}\cdots (v^{*}_{j_{i}})_{m_{j_{i}}-l_{i}}
w'_{(4)}, v_{l_{i+1}}\cdots v_{l_{r}}w_{(3)}))(w_{(1)}\otimes w_{(2)})
\cdot \nno\\
&&\hspace{2em}\cdot (z_{1}-z_{2})^{l_{1}+\cdots +l_{i}+(m_{i+1}-l_{i+1})+\cdots 
+(m_{r}-l_{r})}
\end{eqnarray}
for $m_{1}, \dots, m_{r}\in \mathbb{Z}$,
and any $v_{1}, \dots, 
v_{r}\in V$, and, in particular,
\begin{eqnarray}\label{thm-11.1-10}
\lefteqn{(L'_{P(z_{1}-z_{2})}(0)\beta_{n; k}(w'_{(4)}, 
w_{(3)}))(w_{(1)}
\otimes w_{(2)})}\nno\\
&&=(\beta_{n; k}(L'(0)w'_{(4)}, w_{(3)}))(w_{(1)}
\otimes w_{(2)})\nno\\
&&\quad -(\beta_{n+1; k}(L'(1)w'_{(4)}, w_{(3)}))(w_{(1)}
\otimes w_{(2)})\nno\\
&&\quad -(\beta_{n;k}(w'_{(4)}, L(0)w_{(3)}))(w_{(1)}
\otimes w_{(2)})\nno\\
&&\quad +(\beta_{n+1;k}(w'_{(4)}, L(-1)w_{(3)}))(w_{(1)}
\otimes w_{(2)})\nno\\
&&=(\wt w'_{(4)}-\wt w_{(3)})(\beta_{n; k}(w'_{(4)}, w_{(3)}))(w_{(1)}
\otimes w_{(2)})\nno\\
&&\quad +(\beta_{n; k}((L'(0)-\wt w'_{(4)})w'_{(4)}, w_{(3)}))(w_{(1)}
\otimes w_{(2)})\nno\\
&&\quad -(\beta_{n;k}(w'_{(4)}, (L(0)-\wt w_{(3)})w_{(3)}))(w_{(1)}
\otimes w_{(2)})\nno\\
&&\quad -(\beta_{n+1; k}(L'(1)w'_{(4)}, w_{(3)}))(w_{(1)}
\otimes w_{(2)})\nno\\
&&\quad +(\beta_{n+1;k}(w'_{(4)}, L(-1)w_{(3)}))(w_{(1)}
\otimes w_{(2)}).
\end{eqnarray}

Let $z_{0}=z_{1}-z_{2}$. When $|z_{1}|>|z_{2}|>|z_{0}|>0$, 
\begin{eqnarray}\label{16.27}
\lefteqn{-\sum_{k=1}^{M}\sum_{n\in {\Bbb C}}a_{n; k}(L'(1)w'_{(4)},
w_{(1)}, w_{(2)}, w_{(3)})\cdot}\nno\\
&&\quad\quad\quad\quad\cdot 
e^{((\Delta-1)+n+1)\log (z_{1}-z_{2})}e^{(-n-1)\log z_{2}}
(\log z_{2})^{i_{k}}(\log (z_{1}-z_{2}))^{j_{k}}
\nno\\
&&\quad +\sum_{k=1}^{M}\sum_{n\in {\Bbb C}}a_{n; k}(w'_{(4)},
w_{(1)}, w_{(2)}, L(-1)w_{(3)}) \cdot \nno\\
&&\quad\quad\quad\quad\cdot e^{((\Delta-1)+n+1)\log (z_{1}-z_{2})}
e^{(-n-1)\log z_{2}}
(\log z_{2})^{i_{k}}(\log (z_{1}-z_{2}))^{j_{k}}\nno\\
&&= -
\langle L'(1)w'_{(4)},
{\cal Y}_{1}(
w_{(1)}, x_{1})
{\cal Y}_{2}(w_{(2)}, x_{2})w_{(3)}\rangle_{W_{4}}
\lbar_{x_{1}
=z_{2}+z_{0}, \; x_{2}=z_{2}}\nno\\
&&\quad +
\langle w'_{(4)}, {\cal Y}_{1}(
w_{(1)}, x_{1})
{\cal Y}_{2}(w_{(2)}, x_{2})L(-1)w_{(3)}\rangle_{W_{4}}
\lbar_{x_{1}
=z_{2}+z_{0}, \; x_{2}=z_{2}}.
\end{eqnarray}
By the commutator formula for $L(-1)$ and intertwining operators and the
$L(-1)$-derivative property for intertwining operators, the right-hand side of
(16.27) is equal to
\begin{eqnarray}\label{16.28}
\lefteqn{-
\langle w'_{(4)},
\frac{d}{dx_{1}}({\cal Y}_{1}(
w_{(1)}, x_{1}))
{\cal Y}_{2}(w_{(2)}, x_{2})w_{(3)})\rangle_{W_{4}}
\lbar_{x_{1}
=z_{2}+z_{0}, \; x_{2}=z_{2}}}\nno\\
&&\quad -
\langle w'_{(4)}, {\cal Y}_{1}(
w_{(1)}, x_{1})
\frac{d}{dx_{2}}({\cal Y}_{2}(w_{(2)}, x_{2})w_{(3)})\rangle_{W_{4}}
\lbar_{x_{1}
=z_{2}+z_{0}, \; x_{2}=z_{2}}\nno\\
&&=-\frac{\p}{\p z_{2}}\biggl(
\langle w'_{(4)},
({\cal Y}_{1}(
w_{(1)}, x_{1})
{\cal Y}_{2}(w_{(2)}, x_{2})w_{(3)})\rangle_{W_{4}}
\lbar_{x_{1}
=z_{2}+z_{0}, \; x_{2}=z_{2}}\biggr)\nno\\
&&=-\frac{\p}{\p z_{2}}\Biggl(\sum_{k=1}^{M}
\sum_{n\in {\Bbb C}}a_{n; k}(w'_{(4)},
w_{(1)}, w_{(2)}, w_{(3)}) \cdot \nno\\
&&\quad\quad\quad\quad\quad\quad\quad\quad\cdot
e^{(\Delta+n+1)\log z_{0}}e^{(-n-1)\log z_{2}}
(\log z_{2})^{i_{k}}(\log z_{0})^{j_{k}}
\Biggr)\nno\\
&&=-\sum_{k=1}^{M}\sum_{n\in {\Bbb C}}(-n-1)a_{n}(w'_{(4)},
w_{(1)}, w_{(2)}, w_{(3)}) \cdot \nno\\
&&\quad\quad\quad\quad\quad\quad\quad\quad
\cdot e^{(\Delta+n+1)\log z_{0}}e^{(-n-2)\log z_{2}}
(\log z_{2})^{i_{k}}(\log z_{0})^{j_{k}}\nno\\
&&\quad -\sum_{k=1}^{M}\sum_{n\in {\Bbb C}}a_{n}(w'_{(4)},
w_{(1)}, w_{(2)}, w_{(3)}) \cdot \nno\\
&&\quad\quad\quad\quad\quad\quad\quad\quad
\cdot e^{(\Delta+n+1)\log z_{0}}e^{(-n-1)\log z_{2}}
\left(\frac{\partial}{\partial z_{2}}
(\log z_{2})^{i_{k}}\right)(\log z_{0})^{j_{k}},
\end{eqnarray}
where $\frac{\p}{\p z_{2}}$ is the partial derivation with respect to $z_{2}$ 
acting on functions of $z_{0}$ and $z_{2}$. 

{}From (\ref{16.27}) and (\ref{16.28}), we obtain
\begin{eqnarray}\label{thm-11.1-11}
\lefteqn{-\sum_{k=1}^{M}\sum_{n\in \mathbb{C}}\beta_{n;k}(L'(1)w'_{(4)},
w_{(3)})(w_{(1)}\otimes w_{(2)})\cdot}\nn
&&\quad\quad\quad\quad\quad\quad
 \cdot e^{((\Delta-1)+n+1)\log (z_{1}-z_{2})}e^{(-n-1)\log z_{2}}
(\log z_{2})^{i_{k}}(\log (z_{1}-z_{2}))^{j_{k}}
\nno\\
&&\quad +\sum_{k=1}^{M}\sum_{n\in \mathbb{C}}\beta_{n;k}(w'_{(4)},
L(-1)w_{(3)})(w_{(1)}\otimes w_{(2)}) \cdot\nn
&&\quad\quad\quad\quad\quad\quad
e^{((\Delta-1)+n+1)\log (z_{1}-z_{2})}
e^{(-n-1)\log z_{2}}
(\log z_{2})^{i_{k}}(\log (z_{1}-z_{2}))^{j_{k}}\nno\\
&&=-\sum_{k=1}^{M}\sum_{n\in \mathbb{C}}(-n-1)\beta_{n;k}(w'_{(4)},
w_{(3)})(w_{(1)}\otimes w_{(2)})  \cdot\nn
&&\quad\quad\quad\quad\quad\quad
e^{(\Delta+n+1)\log (z_{1}-z_{2})}e^{(-n-2)\log z_{2}}
(\log z_{2})^{i_{k}}(\log (z_{1}-z_{2}))^{j_{k}}
\nno\\
&&\quad-\sum_{k=1}^{M}\sum_{n\in \mathbb{C}}\beta_{n;k}(w'_{(4)},
w_{(3)})(w_{(1)}\otimes w_{(2)})  \cdot\nn
&&\quad\quad\quad\quad\quad\quad
e^{(\Delta+n+1)\log (z_{1}-z_{2})}e^{(-n-1)\log z_{2}}
\left(\frac{\partial}{\partial z_{2}}
(\log z_{2})^{i_{k}}\right)(\log (z_{1}-z_{2}))^{j_{k}}\nn
\end{eqnarray}
when $|z_{1}|>|z_{2}|>|z_{1}-z_{2}|>0$. 
Since both sides of (\ref{thm-11.1-11}) can be analytically 
extended to the region $|z_{2}|>|z_{1}-z_{2}|>0$, it 
also holds when $|z_{2}|>|z_{1}-z_{2}|>0$. Thus the 
expansion coefficients of both sides of (\ref{thm-11.1-11})
as series in $z_{2}$ and $z_{1}-z_{2}$ are equal. So 
for $k=1, \dots, M$,
\begin{eqnarray}\label{thm-11.1-12}
\lefteqn{\beta_{n+1;k}(L'(1)w'_{(4)},
w_{(3)})(w_{(1)}\otimes w_{(2)})
+\beta_{n+1;k}(w'_{(4)},
L(-1)w_{(3)})(w_{(1)}\otimes w_{(2)})}
\nno\\
&&=-(-n-1)\beta_{n;k}(w'_{(4)},
w_{(3)})(w_{(1)}\otimes w_{(2)}) 
 -d_{n, k}(w'_{(4)}, w_{(1)}, w_{(2)}, w_{(3)}),
\end{eqnarray}
where
\[
d_{n, k}(w'_{(4)}, w_{(1)}, w_{(2)}, w_{(3)})=0
\]
if $i_{k} \ne i_{l}-1$ for $l=1, \dots, M$
and
\[
d_{n, k}(w'_{(4)}, w_{(1)}, w_{(2)}, w_{(3)})=i_{l}\beta_{n;l}(w'_{(4)},
w_{(3)})(w_{(1)}\otimes w_{(2)})
\]
if there exists $l$ such that $i_{k}=i_{l}-1$.
{}From (\ref{thm-11.1-10}) and (\ref{thm-11.1-12}), we obtain
\begin{eqnarray}\label{thm-11.1-12.1}
\lefteqn{((L'_{P(z_{1}-z_{2})}(0)-\wt w'_{(4)}-n-1
+\wt w_{(3)})\beta_{n; k}(w'_{(4)}, 
w_{(3)}))(w_{(1)}
\otimes w_{(2)})}\nno\\
&&=(\beta_{n; k}((L'(0)-\wt w'_{(4)})w'_{(4)}, w_{(3)}))(w_{(1)}
\otimes w_{(2)})\nno\\
&&\quad -(\beta_{n;k}(w'_{(4)}, (L(0)-\wt w_{(3)})w_{(3)}))(w_{(1)}
\otimes w_{(2)})\nno\\
&&\quad -d_{n, k}(w'_{(4)}, w_{(1)}, w_{(2)}, w_{(3)}).
\end{eqnarray}
{}From (\ref{thm-11.1-12.1}), we see
that there exists $N\in \Z_{+}$ such that for all $n\in \C$ and 
$k=1, \dots, M$, 
\begin{equation}\label{thm-11.1-13}
((L'_{P(z_{1}-z_{2})}(0)-\wt w'_{(4)}+\wt w_{(3)}
-n-1)^{N}\beta_{n; k}(w'_{(4)}, 
w_{(3)}))(w_{(1)}
\otimes w_{(2)})=0.
\end{equation}
Thus we obtain the 
following conclusion {}from this calculation:
When
$v_{1}, \dots, v_{r}$, $w_{(3)}$ and $w'_{(4)}$ are homogeneous,
$(v^{*}_{1})_{m_{1}}\cdots 
(v^{*}_{r})_{m_{r}}\beta_{n;k}(w'_{(4)}, 
w_{(3)})$  is also homogeneous of 
(generalized) weight 
\begin{equation}\label{thm-11.1-14}
-(\mbox{\rm wt}\ v_{1}-m_{1}-1)-\cdots -(\mbox{\rm wt}\ v_{r}-m_{r}-1)
+\mbox{\rm wt}\ w'_{(4)}+n+1-\mbox{\rm wt}\ w_{(3)}
\end{equation}
and we have
\begin{eqnarray}\label{thm-11.1-15}
\lefteqn{(v^{*}_{1})_{m_{1}}\cdots (v^{*}_{r})_{m_{r}}\beta_{n;k}(w'_{(4)}, 
w_{(3)})=0}\nno\\
&&\quad\quad \mbox{\rm when}\;\; \Re{(\wt 
(v^{*}_{1})_{m_{1}}\cdots (v^{*}_{r})_{m_{r}}\beta_{n;k}(w'_{(4)}, 
w_{(3)}))}\le N.
\end{eqnarray}
In particular, the (generalized) weight of 
$\beta_{n;k}(w'_{(4)}, 
w_{(3)})$ is $\mbox{\rm wt}\ w'_{(4)}+n+1-\mbox{\rm wt}\ w_{(3)}$.

For $j\in \C$, let 
$\lambda_{j}^{(2)}(w'_{(4)}, 
w_{(3)})\in (W_{1}\otimes W_{2})^{*}$ be defined by 
\begin{eqnarray*}
\lefteqn{\lambda_{j}^{(2)}(w'_{(4)}, 
w_{(3)})(w_{(1)}\otimes w_{(2)})}\nno\\
&&=\sum_{k=1}^{M}
(\beta_{j-\swt w'_{(4)}-1+\swt w_{(3)};k}(w'_{(4)}, 
w_{(3)})(w_{(1)}\otimes w_{(2)})\cdot\nno\\
&&\quad\quad\quad\quad\quad\cdot 
e^{(\swt w_{(1)}+\swt w_{(2)}+j)
\log (z_{1}-z_{2})}e^{(-j+\swt w'_{(4)}-\swt w_{(3)})\log z_{2}}
(\log z_{2})^{i_{k}}(\log (z_{1}-z_{2}))^{j_{k}}
\end{eqnarray*}
for homogeneous $w_{(1)}\in W_{1}$ and $w_{(2)}\in W_{2}$. 
Then $\sum_{j\in \C}\lambda_{j}^{(2)}(w'_{(4)}, 
w_{(3)})(w_{(1)}\otimes w_{(2)})$ is  absolutely convergent  to 
$\mu^{(2)}_{(I_{1}\circ (1_{W_{2}}\otimes I_{2}))'(w'_{(4)}), w_{(3)}}
(w_{(1)}\otimes w_{(2)})$ for $w_{(1)}\in W_{1}$ and $w_{(2)}\in W_{2}$.
Since $\beta_{n;k}(w'_{(4)}, 
w_{(3)})(w_{(1)}\otimes w_{(2)})=a_{n;k}(w'_{(4)}, w_{(1)}, w_{(2)},
w_{(3)})$ and $a_{n;k}(w'_{(4)}, w_{(1)}, w_{(2)},
w_{(3)})=0$ if (\ref{thm-11.1-4}) holds, 
$\lambda_{j}^{(2)}(w'_{(4)}, 
w_{(3)})=0$ if (\ref{thm-11.1-4}) with $n=j-\wt w'_{(4)}-1+\wt w_{(3)}$
holds for $k=1, \dots, M$ and 
$m\in \N$. Since  (\ref{thm-11.1-13}) holds for $n\in \C$ and $k=1, 
\dots, M$, 
\[
((L'_{P(z_{1}-z_{2})}(0)-j)^{N}\lambda_{j}^{(2)}(w'_{(4)}, 
w_{(3)}))(w_{(1)}\otimes w_{(2)})=0
\]
for $j\in \C$. Thus we see that 
the set 
\[
\{(n, i)\in \C\times \N\;|\;(L'_{P(z_{1}-z_{2})}(0)-j)^{i}
\lambda_{j}^{(2)}(w'_{(4)}, 
w_{(3)})(w_{(1)}\otimes w_{(2)})\ne 0\}
\]
is a subset of the set $D\times \{0, \dots, N-1\}$ where 
\[
D=\{\wt w'_{(4)}+1-\wt w_{(3)}-r_{k}+m-1\;|\;k=1, \dots, M,\; m\in \N\}
\]
is clearly a discrete subset of $\R$. By Proposition \ref{discrete-exp-set},
$D\times \{0, \dots, N-1\}$ is a unique expansion set and therefore
its subsets are all unique expansion sets. So
$(I_{1}\circ (1_{W_{2}}\otimes F_{2}))'(w'_{(4)})$ satisfies the 
$P^{(2)}(z_{1}-z_{2})$-grading condition.

Recall the space $W_{\lambda}$ for $\lambda\in (W_{1}\otimes W_{2})^{*}$
in the $P(z)$-local grading restriction condition in Section 5
and the space $W^{(2)}_{\lambda, w_{(3)}}$ for 
$\lambda\in (W_{1}\otimes W_{2}\otimes W_{3})^{*}$ and $w_{(3)}\in W_{3}$
in the $P^{(2)}(z)$-local grading restriction condition in Section 9. 
For fixed $n\in \mathbb{C}$, $w_{(3)}\in W_{3}$ and $w'_{(4)}\in
W'_{4}$, (\ref{thm-11.1-14}) and (\ref{thm-11.1-15}) shows 
that the homogeneous
subspace $(W_{\lambda_{j}^{(2)}(w'_{(4)}, w_{(3)})})_{[l]}$ 
of  a fixed weight
$l\in \mathbb{C}$ of $W_{\lambda_{j}^{(2)}(w'_{(4)}, w_{(3)})}$ 
is $0$ when the 
real part of $l$
is sufficiently small.  In particular, 
$\lambda_{j}^{(2)}(w'_{(4)}, w_{(3)})$, $j\in \C$, satisfy the 
$P(z_{1}-z_{2})$-lower truncation condition. By Theorem \ref{9.7-1},
$W_{\lambda_{j}^{(2)}(w'_{(4)}, w_{(3)})}$ is a generalized $V$-module.
Since $(W_{\lambda_{j}^{(2)}(w'_{(4)}, w_{(3)})})_{[l]}=0$ when the 
real part of $l$
is sufficiently small, it is in fact
is lower truncated. By definition, $W_{\lambda_{j}^{(2)}(w'_{(4)}, w_{(3)})}$ 
is generated by $\lambda_{j}^{(2)}(w'_{(4)}, w_{(3)})$. So
it is a finitely-generated
lower-truncated generalized $V$-module.  By assumption, it is in fact 
in $\ob \mathcal{C}$. Thus $W^{(2)}_{(I_{1}\circ (1_{W_{2}}\otimes I_{2}))'
(w'_{(4)}), w_{(3)}}$ as a sum of these objects of $\mathcal{C}$
is in $W_{1}\hboxtr_{P(z_{1}-z_{2})}W_{2}$. By assumption, 
$W_{1}\hboxtr_{P(z_{1}-z_{2})}W_{2}$ is an object of $\mathcal{C}$.
So $W^{(2)}_{(I_{1}\circ (1_{W_{2}}\otimes I_{2}))'
(w'_{(4)}), w_{(3)}}$ is also an object of $\mathcal{C}$. Since 
$\mathcal{C}$ is a subcategory of $\mathcal{GM}_{sg}$, 
$W^{(2)}_{(I_{1}\circ (1_{W_{2}}\otimes I_{2}))'
(w'_{(4)}), w_{(3)}}$ satisfies the two grading-restriction conditions,
proving that the element
$I_{1}\circ (1_{W_{2}}\otimes I_{2}))'(w'_{(4)})$ of
$(W_{1}\otimes W_{2}\otimes
W_{3})^{*}$ satisfies the $P^{(2)}(z_{1}-z_{2})$-local grading-restriction
condition. 
\epfv

\subsection{Differential equations}

In this subsection, we use differential equations to prove the
convergence and extension property given in the preceding subsection
when $A$ is trivial and objects of $\mathcal{C}$ satisfy natural
considitions.  The results snd the proofs here are in fact the same as
those in \cite{diff-eqn}, except that in this section, we consider
objects of $\mathcal{C}$, not just ordinary $V$-modules, and we
consider logarithmic intertwining operators, not just ordinary
intertwining operators. So in the proofs of the results in this
section, we shall indicate only how the proofs here are different
{}from the corresponding ones in \cite{diff-eqn} and refer the reader
to \cite{diff-eqn} for more details.

Let $V$ be a M\"{o}bius or conformal vertex algebra with $A$ the
trivial group and let $V_{+}=\coprod_{n>0}V_{(n)}$.  Let $W$ be a
generalized $V$-module and let $C_{1}(W)=\{u_{-1}w\;|\; u\in V_{+},\;
w\in W\}$.  If $W/C_{1}(W)$ is finite dimensional, we say that $W$ is
{\it $C_{1}$-cofinite} or satisfies the {\it $C_{1}$-cofiniteness
condition}. If for any real number $N\in {\mathbb R}$,
$\coprod_{\Re{(n)}<N}W_{[n]}$ is finite dimensional, we say that $W$
is {\it quasi-finite dimensional} or satisfies the {\it
quasi-finite-dimensionality condition}.

We have:

\begin{theo}\label{sys}
Let $W_{i}$ for $i=0, \dots, n+1$ be generalized $V$-modules
satisfying the $C_{1}$-cofiniteness condition and the
quasi-finite-dimensionality condition. Then for any $w_{i}\in W_{i}$
for $i=0, \dots, n+1$, there exist
\[
a_{k, \;l}(z_{1}, \dots, z_{n})\in {\mathbb C}[z_{1}^{\pm 1}, \dots,
z_{n}^{\pm 1}, (z_{1}-z_{2})^{-1}, (z_{1}-z_{3})^{-1}, \dots,
(z_{n-1}-z_{n})^{-1}],
\]
for $k=1, \dots, m$ and $l=1, \dots, n,$ such that for any generalized
$V$-modules $\widetilde{W}_{i}$ for $i=1, \dots, n-1$, any
intertwining operators ${\cal Y}_{1}, {\cal Y}_{2}, \dots, {\cal
Y}_{n-1}, {\cal Y}_{n}$, of types ${W'_{0}\choose
W_{1}\widetilde{W}_{1}}$, ${\widetilde{W}_{1}\choose
W_{2}\widetilde{W}_{2}}, \dots, {\widetilde{W}_{n-2}\choose
W_{n-1}\widetilde{W}_{n-1}}$, ${\widetilde{W}_{n-1}\choose
W_{n}W_{n+1}}$, respectively, the series
\begin{equation}
\langle w_{0}, {\cal Y}_{1}(w_{1}, z_{1})\cdots {\cal Y}_{n}(w_{n},
z_{n})w_{n+1}\rangle
\end{equation}
satisfy the expansions of the system of differential equations
\[
\frac{\partial^{m}\varphi}{\partial z_{l}^{m}}+ \sum_{k=1}^{m}a_{k,
\;l}(z_{1}, \dots, z_{n}) \frac{\partial^{m-k}\varphi}{\partial
z_{l}^{m-k}}=0,\;\;\;l=1, \dots, n
\]
in the region $|z_{1}|>\cdots >|z_{n}|>0$. Moreover, there exist
$a_{k, \;l}(z_{1}, \dots, z_{n})$ for $k=1, \dots, m$ and $l=1, \dots,
n$ such that the singular points of the corresponding system are
regular.
\end{theo}
\pf
Proposition 1.1, Corollary 1.2 and Corollary 1.3 in \cite{diff-eqn}
and their proofs still hold when all $V$-modules are 
replaced by strongly-graded generalized $V$-modules. Theorem 
1.4 in \cite{diff-eqn} also still holds, except that in its proof
all $V$-modules are 
replaced by stringly-graded generalized $V$-modules and 
the spaces 
\begin{eqnarray*}
&z_{1}^{\Delta}
\C(\{z_{2}/z_{1}\})[z_{1}^{\pm 1}, z_{2}^{\pm 1}],&\\
&z_{2}^{\Delta}
\mathbb{C}(\{(z_{1}-z_{2})/z_{1}\})[z_{2}^{\pm 1},
(z_{1}-z_{2})^{\pm 1}],&\\
&z_{2}^{\Delta}
\C(\{z_{1}/z_{2}\})[z_{1}^{\pm 1}, z_{2}^{\pm 1}]&
\end{eqnarray*}
are replaced by 
\begin{eqnarray*}
&z_{1}^{\Delta}
\C(\{z_{2}/z_{1}\})[z_{1}^{\pm 1}, z_{2}^{\pm 1}]
[\log z_{2}, \log (z_{1}-z_{2})]],&\\
&z_{2}^{\Delta}
\mathbb{C}(\{(z_{1}-z_{2})/z_{1}\})[z_{2}^{\pm 1},
(z_{1}-z_{2})^{\pm 1}][\log z_{2}, \log (z_{1}-z_{2})],&\\
&z_{2}^{\Delta}
\C(\{z_{1}/z_{2}\})[z_{1}^{\pm 1}, z_{2}^{\pm 1}]
[\log z_{2}, \log (z_{1}-z_{2})],&
\end{eqnarray*}
respectively. 
Theorem 1.6 and its proof also still hold in our setting here.
\epfv

We now have:

\begin{theo}\label{C_1pp}
Suppose that all generalized $V$-modules in ${\cal C}$ satisfy the
$C_{1}$-cofiniteness condition and the quasi-finite-dimensionality
condition. Then the convergence and extension properties for products
and iterates hold in ${\cal C}$. When $\mathcal{C}$ is in
$\mathcal{M}_{sg}$, if every object of $\mathcal{C}$ is a direct sum
of irreducible objects of $\mathcal{C}$ and there are only finitely
many irreducible objects of $\mathcal{C}$ (up to equivalence), then
the convergence and extension properties without logarithms for
products and iterates hold in ${\cal C}$. In addition, for any
$n\in{\mathbb Z}_+$, any $W_i$, $i=0, \dots, n+1$ and $M_i$, $i=1,
\dots, n-1$ of $\ob {\cal C}$, and intertwining operators ${\cal
Y}_{1}, {\cal Y}_{2}, \dots, {\cal Y}_{n-1}, {\cal Y}_{n}$, of types
${W_{0}\choose W_{1}M_{1}}$, ${M_{1}\choose W_{2}M_{2}}, \dots,
{M_{n-2}\choose W_{n-1}M_{n-1}}$, ${M_{n-1}\choose W_{n}W_{n+1}}$,
respectively, and any $w_{(0)}'\in W_{0}'$, $w_{(1)}\in W_{1}, \dots,
w_{(n+1)}\in W_{n+1}$, the series
\begin{equation} 
\langle
w_{(0)}', {\cal Y}_{1}(w_{(1)}, z_{1})\cdots {\cal Y}_{n}(w_{(n)},
z_{n})w_{(n+1)}\rangle
\end{equation}
is absolutely convergent in the region $|z_{1}|>\cdots> |z_{n}|>0$ and
its sum can be analytically extended to a multivalued analytic
function on the region given by $z_{i}\ne 0$, $i=1, \dots, n$,
$z_{i}\ne z_{j}$, $i\ne j$, with regular singular points at $z_{i}=0$,
$i=1, \dots, n$, $z_{i}\ne z_{j}$, $i\ne j$, $i, j=1, \dots, n$.
\end{theo}
\pf The first and the third parts of the theorem follow directly {}from
Theorem \ref{sys} and the theory of differential equations with
regular singular points. The second part has been proved in
\cite{diff-eqn}.  \epfv

\newpage

\setcounter{equation}{0}
\setcounter{rema}{0}

\section{The braided tensor category structure}

In this section, we construct a natural braided tensor category
structure on the category $\mathcal{C}$. The strategy and steps in our
construction in this section are essentially the same as those in
\cite{tensorK}, \cite{tensor5} and \cite{rigidity} in the finitely
reductive case but, instead of the corresponding constructions and
results in \cite{tensor1}, \cite{tensor2}, \cite{tensor3} and
\cite{tensor4}, we of course use all the constructions and results we
have obtained in this work except for those in Section 11. 
The present section is independent of Section 11, which provided
a method for verifying the relevant hypotheses.

We now return to the setting and assumptions of Section 10. 
But in addition to Assumption \ref{assum-assoc},
we also assume the following:

\begin{assum}\label{assum-con}
The convergence and the expansion conditions for intertwining maps in
$\mathcal{C}$ hold.  For objects $W_1$, $W_2$, $W_3$, $W_4$, $W_{5}$,
$M_1$ and $M_{2}$ of $\mathcal{C}$, logarithmic intertwining operators
$\mathcal{Y}_{1}$, $\mathcal{Y}_{2}$ and $\mathcal{Y}_{3}$ of types
${W_5}\choose {W_1M_1}$, ${M_1}\choose {W_2M_{2}}$ and ${M_2}\choose
{W_3W_{4}}$, respectively, $z_{1}, z_{2}, z_{3}\in \C$ satisfying
$|z_{1}|>|z_{2}|>|z_{3}|>0$ and $w_{(1)}\in W_{1}$, $w_{(2)}\in
W_{2}$, $w_{(3)}\in W_{3}$, $w_{(4)}\in W_{4}$ and $w'_{(5)}\in
W'_{5}$, the series
\begin{equation}\label{3-intw-convp}
\sum_{m, n\in {\mathbb C}}\langle w'_{(5)}, \mathcal{Y}_1(w_{(1)}, z_{1})
\pi_{m}(\mathcal{Y}_2(w_{(2)}, z_{2})\pi_{n}(\mathcal{Y}_2(w_{(3)}, z_{3})
w_{(4)}))\rangle_{W_5}
\end{equation}
is absolutely convergent and can be analytically extended to a 
multivalued analytic function on the region given by 
$z_{1}, z_{2}, z_{3}\ne 0$, $z_{1}\ne z_{2}$, $z_{1}\ne z_{3}$
and $z_{2}\ne z_{3}$ with regular singular points at 
$z_{1}=0$, $z_{2}=0$, $z_{3}=0$, $z_{1}=\infty$, $z_{2}=\infty$, 
$z_{3}=\infty$,
$z_{1}= z_{2}$, $z_{1}= z_{3}$
or $z_{2}= z_{3}$. 
\end{assum}

\begin{rema}
{\rm By Theorems \ref{thm-11.1} and \ref{C_1pp}, Assumption
\ref{assum-con} holds if every object of $\mathcal{C}$ satisfies the
$C_{1}$-cofiniteness condition and the quasi-finite dimensionality
condition, or, when $\mathcal{C}$ is in $\mathcal{M}_{sg}$, if every
object of $\mathcal{C}$ is a direct sum of irreducible objects of
$\mathcal{C}$ and there are only finitely many irreducible $C_{1}$-cofinite
objects of
$\mathcal{C}$ (up to equivalence).}
\end{rema}

\subsection{More on tensor products of elements}

Let $W_{1}$, $W_{2}$ and $W_{3}$ be objects of $\mathcal{C}$.
In Section 7, for $z_{1}, z_{2}, z_{3}, z_{4}\in \C$ satisfying 
$|z_{1}|>|z_{2}|>0$ and $|z_{3}|>|z_{4}|>0$, 
we have defined the tensor products
\[
w_{(1)}\boxtimes_{P(z_{1})}(w_{(2)}\boxtimes_{P(z_{2})}w_{(3)})
\in \overline{W_{1}\boxtimes_{P(z_{1})}(W_{2}\boxtimes_{P(z_{2})}W_{3})}
\]
and
\[
(w_{(1)}\boxtimes_{P(z_{4})}w_{(2)})\boxtimes_{P(z_{4})}w_{(3)}
\in \overline{(W_{1}\boxtimes_{P(z_{4})}W_{2})
\boxtimes_{P(z_{4})}W_{3}},
\]
respectively, for $w_{(1)}\in W_{1}$, $w_{(2)}\in W_{2}$ and 
$w_{(3)}\in W_{3}$. In the proof of the commutativity of the 
hexagon diagrams below, we shall also need 
tensor products
of elements $w_{(1)}\in W_{1}$, $w_{(2)}\in W_{2}$ and 
$w_{(3)}\in W_{3}$ in
$\overline{W_{1}\boxtimes_{P(z_{1})}(W_{2}\boxtimes_{P(z_{2})}W_{3})}$
and in
$\overline{(W_{1}\boxtimes_{P(z_{4})}W_{2})
\boxtimes_{P(z_{3})}W_{3}}$ when $z_{1}, z_{2}, z_{3}, z_{4}\in \C^{\times}$ 
satisfy $z_{1}\ne z_{2}$ and $z_{3}\ne z_{4}$ but do not
necessarily satisfy the inequality $|z_{1}|>|z_{2}|>0$ or 
$|z_{3}|>|z_{4}|>0$. Here we first define these elements.

Let 
$\Y_{1}=\Y_{\boxtimes_{P(z_{1})}, 0}$, 
$\Y_{2}=\Y_{\boxtimes_{P(z_{2})}, 0}$,
$\Y_{3}=\Y_{\boxtimes_{P(z_{3})}, 0}$ and 
$\Y_{4}=\Y_{\boxtimes_{P(z_{4})}, 0}$
be intertwining operators of types 
\[
{W_{1}\boxtimes_{P(z_{1})}(W_{2}\boxtimes_{P(z_{2})}W_{3})\choose 
W_{1}\;W_{2}\boxtimes_{P(z_{2})}W_{3}},
\] 
\[
{W_{2}\boxtimes_{P(z_{2})}W_{3}\choose W_{2}\; W_{3}},
\]
\[
{(W_{1}\boxtimes_{P(z_{4})}W_{2})\boxtimes_{P(z_{3})}W_{3}\choose 
W_{1}\boxtimes_{P(z_{4})}W_{2}\;W_{3}}
\] 
and
\[
{W_{1}\boxtimes_{P(z_{4})}W_{2}\choose W_{1}\; W_{2}},
\]
respectively, corresponding to the intertwining maps 
$\boxtimes_{P(z_{1})}$, $\boxtimes_{P(z_{2})}$, $\boxtimes_{P(z_{3})}$
and $\boxtimes_{P(z_{4})}$, respectively (see (\ref{YIp}) and 
(\ref{recover})). 
Then by our assumption, 
\[
\langle w', \Y_{1}(w_{(1)}, \zeta_{1})
\Y_{2}(w_{(2)}, \zeta_{2})w_{(3)}\rangle
\]
and 
\[
\langle \tilde{w}', \Y_{3}(\Y_{4}(w_{(1)}, \zeta_{4})
w_{(2)}, \zeta_{3})w_{(3)}\rangle
\]
are absolutely convergent for 
\[
w'\in (W_{1}\boxtimes_{P(z_{1})}(W_{2}\boxtimes_{P(z_{2})}W_{3}))'
\]
and 
\[
\tilde{w}'\in ((W_{1}\boxtimes_{P(z_{4})}W_{2})\boxtimes_{P(z_{2})}W_{3})',
\]
when $|\zeta_{1}|>|\zeta_{2}|>0$ and when $|\zeta_{3}|>|\zeta_{4}|>0$,
respectively, and can be analytically extended to multivalued analytic
functions in the regions given by $\zeta_{1}, \zeta_{2}\ne 0$ and
$\zeta_{1}\ne \zeta_{2}$ and by $\zeta_{3}, \zeta_{4}\ne 0$ and
$\zeta_{3}\ne -\zeta_{4}$, respectively. If we cut these regions along 
$\zeta_{1}, \zeta_{2}\ge \R_{+}$ and $\zeta_{3}, \zeta_{4}\in \R_{+}$,
respectively, we obtain simply-connected regions and we can choose
single-valued branches of the multivalued analytic functions above.
In particular, we have the branches of these two multivalued 
analytic functions such that their values at points satisfying 
$|\zeta_{1}|>|\zeta_{2}|>0$ and $|\zeta_{3}|>|\zeta_{4}|>0$ are
\[
\langle w', w_{(1)}\boxtimes_{P(\zeta_{1})}(w_{(2)}
\boxtimes_{P(\zeta_{2})} w_{(3)})
\rangle
\]
and 
\[
\langle \tilde{w}', 
(w_{(1)}\boxtimes_{P(\zeta_{4})}w_{(2)})\boxtimes_{P(\zeta_{3})} w_{(3)}
\rangle,
\]
respectively.

Then we immediately have the following result:

\begin{propo}
Let $w_{(1)}\in W_{1}$, $w_{(2)}\in W_{2}$ and $w_{(3)}\in W_{3}$.
Then for any $z_{1}, z_{2}, z_{3}, z_{4}\in \C^{\times}$ 
satisfy $z_{1}\ne z_{2}$ and $z_{3}\ne -z_{4}$, 
there exist unique elements 
\[
w_{(1)}\boxtimes_{P(z_{1})}(w_{(2)}\boxtimes_{P(z_{2})} w_{(3)})
\in \overline{W_{1}\boxtimes_{P(z_{1})}(W_{2}\boxtimes_{P(z_{2})}W_{3})}
\]
and 
\[
(w_{(1)}\boxtimes_{P(z_{4})}w_{(2)})\boxtimes_{P(z_{3})} w_{(3)}
\in \overline{(W_{1}\boxtimes_{P(z_{4})}W_{2})\boxtimes_{P(z_{3})}W_{3}}
\]
such that 
for any 
\[
w'\in (W_{1}\boxtimes_{P(z_{1})}(W_{2}\boxtimes_{P(z_{2})}W_{3}))'
\]
and 
\[
\tilde{w}'\in ((W_{1}\boxtimes_{P(z_{4})}W_{2})\boxtimes_{P(z_{3})}W_{3})',
\]
the numbers
\begin{equation}\label{general-tsr-1}
\langle w', w_{(1)}\boxtimes_{P(z_{1})}(w_{(2)}\boxtimes_{P(z_{2})} w_{(3)})
\rangle
\end{equation}
and 
\begin{equation}\label{general-tsr-2}
\langle \tilde{w}', 
(w_{(1)}\boxtimes_{P(z_{4})}w_{(2)})\boxtimes_{P(z_{3})} w_{(3)}
\rangle
\end{equation}
are values at $(\zeta_{1}, \zeta_{2})=(z_{1}, z_{2})$ and 
$(\zeta_{3}, \zeta_{4})=(z_{3}, z_{4})$, respectively,
of the branches of the multivalued analytic functions above
of $\zeta_{1}$ and $\zeta_{2}$ and of $\zeta_{3}$ and $\zeta_{4}$
above, respectively. \epf
\end{propo}

\begin{rema}
{\rm {}From the definition of 
\[
w_{(1)}\boxtimes_{P(z_{1})}(w_{(2)}\boxtimes_{P(z_{2})} w_{(3)})
\]
and 
\[
(w_{(1)}\boxtimes_{P(z_{4})}w_{(2)})\boxtimes_{P(z_{3})} w_{(3)},
\]
we see that when $|z_{1}|=|z_{2}|$ or $|z_{3}|=|z_{4}|$, 
they are uniquely determined by 
\begin{eqnarray*}
\lefteqn{\langle w', w_{(1)}\boxtimes_{P(\zeta_{1})}(w_{(2)}
\boxtimes_{P(\zeta_{2})} w_{(3)})
\rangle}\nn
&&=\lim_{\zeta_{1}\to z_{1}, \zeta_{2}\to z_{2},
|\zeta_{1}|>|\zeta_{2}|>0}
\langle w', \Y_{1}(w_{(1)}, \zeta_{1})\Y_{2}(w_{(2)},
\zeta_{2}) w_{(3)}
\rangle
\end{eqnarray*}
and 
\begin{eqnarray*}
\lefteqn{\langle \tilde{w}', 
(w_{(1)}\boxtimes_{P(\zeta_{4})}w_{(2)})\boxtimes_{P(\zeta_{3})} w_{(3)}
\rangle}\nn
&&=\lim_{\zeta_{3}\to z_{4}, \zeta_{4}\to z_{4},
|\zeta_{3}|>|\zeta_{4}|>0}\langle \tilde{w}', 
\Y_{3}(\Y_{4}(w_{(1)}, \zeta_{4})w_{(2)})\zeta_{3}) w_{(3)}
\rangle
\end{eqnarray*}
for 
\[
w'\in (W_{1}\boxtimes_{P(z_{1})}(W_{2}\boxtimes_{P(z_{2})}W_{3}))'
\]
and 
\[
\tilde{w}'\in ((W_{1}\boxtimes_{P(z_{4})}W_{2})\boxtimes_{P(z_{3})}W_{3})'.
\]}
\end{rema}

Propositions \ref{prod=0=>comp=0} and 
\ref{iter=0=>comp=0}, Corollaries \ref{prospan} and \ref{iterspan}
and the definitions of tensor products of three elements above
immediately give:

\begin{propo}
For any $z_{1}, z_{2}, z_{3}, z_{4}\in \C^{\times}$ 
satisfying $z_{1}\ne z_{2}$ and $z_{3}\ne -z_{4}$, the elements of the form 
\begin{eqnarray*}
&\pi_{n}(w_{(1)}\boxtimes_{P(z_{1})}(w_{(2)}\boxtimes_{P(z_{2})}
w_{(3)})),&\\
&\pi_{n}((w_{(1)}\boxtimes_{P(z_{4})}w_{(2)})\boxtimes_{P(z_{3})}
w_{(3)})&
\end{eqnarray*}
for $n\in \C$, $w_{(1)}\in W_{1}$, $w_{(2)}\in W_{2}$, 
$w_{(3)}\in W_{3}$ span 
\begin{eqnarray*}
&W_{1}\boxtimes_{P(z_{1})}(W_{2}\boxtimes_{P(z_{2})}
W_{3}),& \\
&(W_{1}\boxtimes_{P(z_{4})}W_{2})\boxtimes_{P(z_{3})}
W_{3},&
\end{eqnarray*}
respectively.\epf
\end{propo}

Next we discuss tensor products of four elements. These
are needed in the proof of the commutativity of 
the pentagon diagram below.

\begin{propo}

\begin{enumerate}

\item Let $W_1$, $W_2$, $W_3$, $W_4$, $W_{5}$, $M_1$ and $M_{2}$ be objects 
of ${\cal C}$, $z_1,z_2, z_{3}$ nonzero complex numbers  satisfying
$|z_1|>|z_2|>|z_{3}|>0$, $I_1$, $I_2$ and 
$I_{3}$ $P(z_1)$-, 
$P(z_2)$- and $P(z_{3})$-intertwining maps of type ${W_5}\choose {W_1M_1}$, 
${M_1}\choose {W_2M_{2}}$ and ${M_2}\choose {W_3W_{4}}$, 
respectively. Then for $w_{(1)}\in W_1$,
$w_{(2)}\in W_2$, $w_{(3)}\in W_3$, $w_{(4)}\in W_4$ and $w'_{(5)}\in W'_5$,
the series
\begin{equation}\label{3-convp}
\sum_{m, n\in {\mathbb C}}\langle w'_{(5)}, I_1(w_{(1)}\otimes
\pi_m(I_2(w_{(2)}\otimes \pi_n(I_3(w_{(3)}\otimes 
w_{(4)}))))\rangle_{W_5}
\end{equation}
is absolutely convergent. 

\item Let $W_1$, $W_2$, $W_3$, $W_4$, $W_{5}$, $M_3$ and $M_{4}$ be
objects of ${\cal C}$, $z_{1}, z_{23}, z_{3}$ nonzero complex numbers
satisfying $|z_3|>|z_{23}|>0$ and $|z_{1}|>|z_{3}|+|z_{23}|>0$, $I_1$,
$I_2$ and $I_{3}$ $P(z_{1})$-, $P(z_3)$- and $P(z_{23})$-intertwining
maps of type ${W_5}\choose {W_1M_3}$, ${M_3}\choose {M_{4}W_4}$ and
${M_4}\choose {W_2W_{3}}$, respectively. Then for $w_{(1)}\in W_1$,
$w_{(2)}\in W_2$, $w_{(3)}\in W_3$, $w_{(4)}\in W_4$ and $w'_{(5)}\in
W'_5$, the series
\begin{equation}\label{3-convip}
\sum_{m, n\in {\mathbb C}}\langle w'_{(5)}, I_1(w_{(1)}\otimes
\pi_m(I_2(\pi_n(I_3(w_{(2)}\otimes w_{(3)}))\otimes 
w_{(4)})))\rangle_{W_5}
\end{equation}
is absolutely convergent. 

\item Let $W_1$, $W_2$, $W_3$, $W_4$, $W_{5}$, $M_5$ and $M_{6}$ be
objects of ${\cal C}$, $z_{3}, z_{13}, z_{23}$ nonzero complex numbers
satisfying $|z_{3}|>|z_{13}|>|z_{23}|>0$, $I_1$, $I_2$ and $I_{3}$
$P(z_{3})$-, $P(z_{13})$- and $P(z_{23})$-intertwining maps of type
${W_5}\choose {M_5W_{4}}$, ${M_5}\choose {W_{1}M_6}$ and ${M_6}\choose
{W_2W_{3}}$, respectively. Then for $w_{(1)}\in W_1$, $w_{(2)}\in
W_2$, $w_{(3)}\in W_3$, $w_{(4)}\in W_4$ and $w'_{(5)}\in W'_5$, the
series
\begin{equation}\label{3-convpi}
\sum_{m, n\in {\mathbb C}}\langle w'_{(5)}, I_1(\pi_m(I_2(w_{(1)}\otimes
\pi_n(I_3(w_{(2)}\otimes w_{(3)}))))\otimes 
w_{(4)})\rangle_{W_5}
\end{equation}
is absolutely convergent. 

\item Let $W_1$, $W_2$, $W_3$, $W_4$, $W_{5}$, $M_7$ and $M_{8}$ be
objects of ${\cal C}$, $z_{3}, z_{23}, z_{12}$ nonzero complex numbers
satisfying $|z_{23}|>|z_{12}|>0$ and $|z_{3}|>|z_{23}|+|z_{12}|>0$,
$I_1$, $I_2$ and $I_{3}$ $P(z_{3})$-, $P(z_{23})$- and
$P(z_{12})$-intertwining maps of type ${W_5}\choose {M_7W_{4}}$,
${M_7}\choose {M_8W_{3}}$ and ${M_8}\choose {W_1W_{2}}$,
respectively. Then for $w_{(1)}\in W_1$, $w_{(2)}\in W_2$, $w_{(3)}\in
W_3$, $w_{(4)}\in W_4$ and $w'_{(5)}\in W'_5$, the series
\begin{equation}\label{3-convi}
\sum_{m, n\in {\mathbb C}}\langle w'_{(5)}, I_1(\pi_m(I_2(\pi_n(I_3(w_{(1)}\otimes
w_{(2)}))\otimes w_{(3)}))\otimes 
w_{(4)})\rangle_{W_5}
\end{equation}
is absolutely convergent. 

\item Let $W_1$, $W_2$, $W_3$, $W_4$, $W_{5}$, $M_9$ and $M_{10}$ be
objects of ${\cal C}$, $z_{12}, z_2, z_{3}$ nonzero complex numbers
satisfying $|z_2|>|z_{12}|+|z_{3}|>0$, $I_1$, $I_2$ and $I_{3}$
$P(z_2)$-, $P(z_{12})$- and $P(z_{3})$-intertwining maps of type
${W_5}\choose {M_9M_{10}}$, ${M_9}\choose {W_1W_{2}}$ and
${M_{10}}\choose {W_3W_{4}}$, respectively. Then for $w_{(1)}\in W_1$,
$w_{(2)}\in W_2$, $w_{(3)}\in W_3$, $w_{(4)}\in W_4$ and $w'_{(5)}\in
W'_5$, the series
\begin{equation}\label{3-convcomp}
\sum_{m, n\in {\mathbb C}}\langle w'_{(5)}, I_1(\pi_m(I_2(w_{(1)}\otimes
w_{(2)}))\otimes \pi_n(I_2(w_{(3)}\otimes 
w_{(4)}))\rangle_{W_5}
\end{equation}
is absolutely convergent. 
\end{enumerate}
\end{propo}
\pf
The absolute convergence of (\ref{3-convp}) follows 
immediately {}from Assumption \ref{assum-con}
and Proposition \ref{im:correspond}.

To prove the absolute convergence of (\ref{3-convip}), let $\Y_{1}$,
$\Y_{2}$ and $\Y_{3}$ be the intertwining operators corresponding to
$I_{1}$, $I_{2}$ and $I_{3}$, using $p=0$ as usual, that is, 
$\Y_{1}=\Y_{I_{1}, 0}$,
$\Y_{2}=\Y_{I_{2}, 0}$ and $\Y_{3}=\Y_{I_{3}, 0}$.  We would like to
prove that for $w_{(1)}\in W_1$, $w_{(2)}\in W_2$, $w_{(3)}\in W_3$,
$w_{(4)}\in W_4$ and $w'_{(5)}\in W'_5$,
\begin{equation}\label{3-convip-1}
\sum_{m, n\in {\mathbb C}}\langle w'_{(5)}, \Y_1(w_{(1)}, z_{1})
\pi_m(\Y_2(\pi_n(\Y_3(w_{(2)}, z_{23}) w_{(3)}), z_{3}) 
w_{(4)}))\rangle_{W_5}
\end{equation}
is absolutely convergent when $|z_3|>|z_{23}|>0$ and 
$|z_{1}|>|z_{3}|+|z_{23}|>0$. By the $L(0)$-conjugation 
property for intertwining 
operators,  (\ref{3-convip-1}) is equal to 
\begin{eqnarray}\label{3-convip-1.5}
\lefteqn{\sum_{m, n\in {\mathbb C}}\langle z_{3}^{L(0)}w'_{(5)}, 
\Y_1(z_{3}^{-L(0)}w_{(1)}, z_{1}z_{3}^{-1})\cdot}\nn
&&\quad\quad\quad
\cdot\pi_m(\Y_2(\pi_n(\Y_3(z_{3}^{-L(0)}w_{(2)}, z_{23}z_{3}^{-1}) 
z_{3}^{-L(0)}w_{(3)}), 1) 
z_{3}^{-L(0)}w_{(4)}))\rangle_{W_5}.\nn
\end{eqnarray}
Since (\ref{3-convip-1}) is equal to (\ref{3-convip-1.5}), we
see that to prove that for $w_{(1)}\in W_1$,
$w_{(2)}\in W_2$, $w_{(3)}\in W_3$, $w_{(4)}\in W_4$ and $w'_{(5)}\in W'_5$,
(\ref{3-convip-1}) is absolutely convergent 
when $|z_3|>|z_{23}|>0$ and 
$|z_{1}|>|z_{3}|+|z_{23}|>0$ is equivalent to prove that for
$w_{(1)}\in W_1$,
$w_{(2)}\in W_2$, $w_{(3)}\in W_3$, $w_{(4)}\in W_4$ and $w'_{(5)}\in W'_5$,
\begin{equation}\label{3-convip-1.7}
\sum_{m, n\in {\mathbb C}}\langle w'_{(5)}, 
\Y_1(w_{(1)}, \zeta_{1})
\pi_m(\Y_2(\pi_n(\Y_3(w_{(2)}, \zeta_{23}) 
w_{(3)}), 1) 
w_{(4)}))\rangle_{W_5}
\end{equation}
is absolutely convergent when $1>|\zeta_{23}|>0$
and $|\zeta_{1}|>1+|\zeta_{23}|>0$.

Using the associativity of intertwining operators, 
we know that there exist an object $M$ and
intertwining operators $\Y_{4}$ and $\Y_{5}$ of types ${M_{3}\choose
W_{2}M}$ and ${M\choose W_{3}W_{4}}$, respectively, such that for
$w'\in M_{3}'$, when $1>|\zeta_{23}|>0$, the series
\[
\sum_{n\in {\mathbb C}}\langle w', 
\Y_2(\pi_n(\Y_3(w_{(2)}, 
\zeta_{23}) w_{(3)}), 1) 
w_{(4)}\rangle_{M_{3}}
=\langle w', 
\Y_2(\Y_3(w_{(2)}, \zeta_{23}) w_{(3)}), 1) 
w_{(4)}\rangle_{M_{3}}
\]
is absolutely convergent and, when $1>|\zeta_{23}|>0$ and
$|\zeta_{23}+1|>1$, its sum is equal to 
\[
\langle w', 
\Y_4(w_{(2)}, \zeta_{23}+1)\Y_{5}(w_{(3)}, 1) 
w_{(4)}\rangle_{M_{3}}.
\]
By Assumption \ref{assum-con}, we know that 
\begin{equation}\label{3-convip-2}
\sum_{m\in {\mathbb C}}\sum_{n\in \C}\langle w'_{(5)}, 
\Y_1(w_{(1)}, \zeta_{1})\pi_{m}(
\Y_4(w_{(2)}, \zeta_{23}+1)
\pi_{n}(\Y_{5}(w_{(3)}, 1) 
w_{(4)}))\rangle_{M_{3}}
\end{equation}
is absolutely convergent and can be analytically extended to 
a multivalued analytic function $f(\zeta_{1}, \zeta_{23})$ 
on the region given by 
$\zeta_{1}, \zeta_{23}+1\ne 0$, 
$\zeta_{1}\ne \zeta_{23}+1$, $\zeta_{1}\ne 1$
and $\zeta_{23}+1\ne 1$ with the only possible 
regular singular points at 
$\zeta_{1}=0$, $\zeta_{23}+1=0$, $\zeta_{1}= \zeta_{23}+1$, $\zeta_{1}= 1$,
$\zeta_{23}=0$, $\zeta_{1}=\infty$ or $\zeta_{23}+1=\infty$.
Since $(\zeta_{1}, \zeta_{23})=(\infty, 0)$ is a regular singular 
point of $f(\zeta_{1}, \zeta_{23})$ and for fixed 
$\zeta_{1}, \zeta_{23}$ satisfying $|\zeta_{1}|>|\zeta_{23}|+1$
and $0<|\zeta_{23}|<1$, there exists a positive real number 
$r<1$ such that $|\zeta_{1}|>r+1$ and $0<|\zeta_{23}|<r$,
$f(\zeta_{1}, \zeta_{23})$ can be expanded as a series 
in powers of $\zeta_{1}$ and $\zeta_{23}$ and in positive 
integral powers of $\log \zeta_{1}$ and $\log \zeta_{23}$
when $|\zeta_{1}|>|\zeta_{23}|+1$
and $0<|\zeta_{23}|<1$.
Thus we see that  
(\ref{3-convip-1.7}) is in fact one value of this expansion of 
$f(\zeta_{1}, \zeta_{23})$, proving that (\ref{3-convip-1.7})
is absolutely convergent when $|\zeta_{1}|>|\zeta_{23}|+1$
and $0<|\zeta_{23}|<1$. 

The proofs of the absolute convergence of
(\ref{3-convpi})--(\ref{3-convcomp}) are similar.  \epfv

The special case that the $P(\cdot)$-intertwining maps considered
are $\boxtimes_{P(\cdot)}$ gives the following:

\begin{corol}\label{t-prod-4-elts}
Let $W_1$, $W_2$, $W_3$, $W_4$ be objects
of ${\cal C}$ and $w_{(1)}\in W_1$,
$w_{(2)}\in W_2$, $w_{(3)}\in W_3$, $w_{(4)}\in W_4$.
Then we have:

\begin{enumerate}

\item For $z_1, z_2, z_{3}\in \C^{\times}$  satisfying
$|z_1|>|z_2|>|z_{3}|>0$ and 
\[
w'\in (W_{1}\boxtimes_{P(z_{1})} 
(W_{2}\boxtimes_{P(z_{2})} (W_{3}\boxtimes_{P(z_{3})} 
W_{4})))'=W_{1}\hboxtr_{P(z_{1})} 
(W_{2}\boxtimes_{P(z_{2})} (W_{3}\boxtimes_{P(z_{3})} 
W_{4})),
\]
the series 
\begin{equation}\label{4-elts-conv-p}
\sum_{m, n\in {\mathbb C}}\langle w', w_{(1)}\boxtimes_{P(z_{1})} 
\pi_{m}(w_{(2)}\boxtimes_{P(z_{2})} \pi_{n}(w_{(3)}\boxtimes_{P(z_{3})} 
w_{(4)}))\rangle
\end{equation}
is absolutely convergent. 

\item For $z_1, z_{23}, z_3, \in \C^{\times}$  satisfying
$|z_3|>|z_{23}|>0$ and $|z_{1}|>|z_{3}|+|z_{23}|>0$ and
\[
w'\in (W_{1}\boxtimes_{P(z_{1})} 
((W_{2}\boxtimes_{P(z_{23})} W_{3})\boxtimes_{P(z_{3})} 
W_{4}))'=W_{1}\hboxtr_{P(z_{1})} 
((W_{2}\boxtimes_{P(z_{23})} W_{3})\boxtimes_{P(z_{3})} 
W_{4}),
\]
the series 
\begin{equation}\label{4-elts-conv-ip}
\sum_{m, n\in {\mathbb C}}\langle w', w_{(1)}\boxtimes_{P(z_{1})} 
\pi_{m}(\pi_{n}(w_{(2)}\boxtimes_{P(z_{2})} w_{(3)})\boxtimes_{P(z_{3})} 
w_{(4)})\rangle
\end{equation}
is absolutely convergent. 

\item For $z_3, z_{13}, z_{23}\in \C^{\times}$  satisfying
$|z_3|>|z_{13}|>|z_{23}|>0$ and 
\[
w'\in ((W_{1}\boxtimes_{P(z_{13})} 
(W_{2}\boxtimes_{P(z_{23})} W_{3}))\boxtimes_{P(z_{3})} 
W_{4})'=(W_{1}\boxtimes_{P(z_{13})} 
(W_{2}\boxtimes_{P(z_{23})} W_{3}))\hboxtr_{P(z_{3})} 
W_{4},
\]
the series 
\begin{equation}\label{4-elts-conv-pi}
\sum_{m, n\in {\mathbb C}}\langle w', \pi_{m}(w_{(1)}\boxtimes_{P(z_{13})} 
\pi_{n}(w_{(2)}\boxtimes_{P(z_{23})} w_{(3)}))\boxtimes_{P(z_{3})} 
w_{(4)}\rangle
\end{equation}
is absolutely convergent. 

\item For $z_3, z_{23}, z_{12}\in \C^{\times}$  satisfying
$|z_{23}|>|z_{12}|>0$ and $|z_{3}|>|z_{23}|+|z_{12}|>0$ and
\[
w'\in (((W_{1}\boxtimes_{P(z_{12})} 
W_{2})\boxtimes_{P(z_{23})} W_{3})\boxtimes_{P(z_{3})} 
W_{4})'=((W_{1}\boxtimes_{P(z_{12})} 
W_{2})\boxtimes_{P(z_{23})} W_{3})\hboxtr_{P(z_{3})} 
W_{4},
\]
the series 
\begin{equation}\label{4-elts-conv-i}
\sum_{m, n\in {\mathbb C}}\langle w', 
\pi_{m}(\pi_{n}(w_{(1)}\boxtimes_{P(z_{12})} 
w_{(2)})\boxtimes_{P(z_{2})} w_{(3)})\boxtimes_{P(z_{3})} 
w_{(4)}\rangle
\end{equation}
is absolutely convergent. 

\item For $z_{12}, z_2, z_{3}\in \C^{\times}$  satisfying
$|z_2|>|z_{12}|+|z_{3}|>0$ and 
\[
w'\in ((W_{1}\boxtimes_{P(z_{12})} 
W_{2})\boxtimes_{P(z_{2})} (W_{3}\boxtimes_{P(z_{3})} 
W_{4})))'=(W_{1}\boxtimes_{P(z_{12})} 
W_{2})\hboxtr_{P(z_{2})} (W_{3}\boxtimes_{P(z_{3})} 
W_{4})),
\]
the series 
\begin{equation}\label{4-elts-conv-comp}
\sum_{m, n\in {\mathbb C}}\langle w', \pi_{m}(w_{(1)}\boxtimes_{P(z_{12})} 
w_{(2)})\boxtimes_{P(z_{2})} \pi_{n}(w_{(3)}\boxtimes_{P(z_{3})} 
w_{(4)})\rangle
\end{equation}
is absolutely convergent. 

\end{enumerate}
\end{corol}

Note that the sums of (\ref{4-elts-conv-p}), (\ref{4-elts-conv-ip}),
(\ref{4-elts-conv-pi}), 
(\ref{4-elts-conv-i}), (\ref{4-elts-conv-comp})
define elements of 
\begin{eqnarray}
&\overline{W_{1}\boxtimes_{P(z_{1})}(W_{2}\boxtimes_{P(z_{2})}
(W_{3}\boxtimes_{P(z_{3})}W_{4}))},&\label{tr-4-mod-p}\\
&\overline{W_{1}\boxtimes_{P(z_{1})}((W_{2}\boxtimes_{P(z_{23})}
W_{3})\boxtimes_{P(z_{1})}W_{4})},&\label{tr-4-mod-ip}\\
&\overline{(W_{1}\boxtimes_{P(z_{13})}(W_{2}\boxtimes_{P(z_{23})}
W_{3}))\boxtimes_{P(z_{3})}W_{4}},&\label{tr-4-mod-pi}\\
&\overline{((W_{1}\boxtimes_{P(z_{12})}W_{2})\boxtimes_{P(z_{23})}
W_{3})\boxtimes_{P(z_{3})}W_{4}},&\label{tr-4-mod-i}\\
&\overline{(W_{1}\boxtimes_{P(z_{12})}W_{2})\boxtimes_{P(z_{2})}
(W_{3}\boxtimes_{P(z_{3})}W_{4})},&\label{tr-4-mod-comp}
\end{eqnarray}
respectively, for $z_{1}, z_{2}, z_{3}, z_{12}, z_{13}, z_{23}$
satisfying the corresponding 
inequalities. 
We shall denote these five elements by 
\begin{eqnarray*}
&w_{(1)}\boxtimes_{P(z_{1})}(w_{(2)}\boxtimes_{P(z_{2})}
(w_{(3)}\boxtimes_{P(z_{3})}w_{(4)})), &\\
&w_{(1)}\boxtimes_{P(z_{1})}((w_{(2)}\boxtimes_{P(z_{23})}
w_{(3)})\boxtimes_{P(z_{1})}w_{(4)}),&\\
&(w_{(1)}\boxtimes_{P(z_{13})}(w_{(2)}\boxtimes_{P(z_{23})}
w_{(3)}))\boxtimes_{P(z_{3})}w_{(4)},&\\
&((w_{(1)}\boxtimes_{P(z_{12})}w_{(2)})\boxtimes_{P(z_{23})}
w_{(3)})\boxtimes_{P(z_{3})}w_{(4)},&\\
&(w_{(1)}\boxtimes_{P(z_{12})}w_{(2)})\boxtimes_{P(z_{2})}
(w_{(3)}\boxtimes_{P(z_{3})}w_{(4)}),&
\end{eqnarray*}
respectively.

\begin{propo}
The elements of the form 
\begin{eqnarray*}
&\pi_{n}(w_{(1)}\boxtimes_{P(z_{1})}(w_{(2)}\boxtimes_{P(z_{2})}
(w_{(3)}\boxtimes_{P(z_{3})}w_{(4)}))),&\\
&\pi_{n}(w_{(1)}\boxtimes_{P(z_{1})}((w_{(2)}\boxtimes_{P(z_{23})}
w_{(3)})\boxtimes_{P(z_{1})}w_{(4)})),&\\
&\pi_{n}((w_{(1)}\boxtimes_{P(z_{13})}(w_{(2)}\boxtimes_{P(z_{23})}
w_{(3)}))\boxtimes_{P(z_{3})}w_{(4)}),&\\
&\pi_{n}(((w_{(1)}\boxtimes_{P(z_{12})}w_{(2)})\boxtimes_{P(z_{23})}
w_{(3)})\boxtimes_{P(z_{3})}w_{(4)}),&\\
&\pi_{n}((w_{(1)}\boxtimes_{P(z_{12})}w_{(2)})\boxtimes_{P(z_{2})}
(w_{(3)}\boxtimes_{P(z_{3})}w_{(4)}))&
\end{eqnarray*}
for $n\in \C$, $w_{(1)}\in W_{1}$, $w_{(2)}\in W_{2}$, 
$w_{(3)}\in W_{3}$, $w_{(4)}\in W_{4}$ span 
\begin{eqnarray*}
&W_{1}\boxtimes_{P(z_{1})}(W_{2}\boxtimes_{P(z_{2})}
(W_{3}\boxtimes_{P(z_{3})}W_{4})),&\\
&W_{(1)}\boxtimes_{P(z_{1})}((W_{2}\boxtimes_{P(z_{23})}
W_{3})\boxtimes_{P(z_{1})}W_{4}),&\\
&(W_{1}\boxtimes_{P(z_{13})}(W_{2}\boxtimes_{P(z_{23})}
W_{3}))\boxtimes_{P(z_{3})}W_{4}),&\\
&((W_{1}\boxtimes_{P(z_{12})}W_{2})\boxtimes_{P(z_{23})}
W_{3})\boxtimes_{P(z_{3})}W_{4}),&\\
&(W_{1}\boxtimes_{P(z_{12})}W_{2})\boxtimes_{P(z_{2})}
(W_{3}\boxtimes_{P(z_{3})}W_{4}))&
\end{eqnarray*}
respectively.
\end{propo}
\pf
The proof is the same as those of Corollaries \ref{prospan}
and \ref{iterspan}. 
\epfv

\begin{propo}
Let $W_1$, $W_2$, $W_3$, $W_4$ be objects
of ${\cal C}$, $w_{(1)}\in W_1$,
$w_{(2)}\in W_2$, $w_{(3)}\in W_3$, $w_{(4)}\in W_4$ and 
$z_{1}, z_{2}, z_{3}\in \C^{\times}$ such that 
$z_{12}=z_{1}-z_{2}\ne 0$, $z_{13}=z_{1}-z_{3}\ne 0$,
$z_{23}=z_{2}-z_{3}\ne 0$.
Then we have:

\begin{enumerate}

\item When
$|z_{1}|>|z_{2}|>|z_{12}|+|z_{3}|>0$, 
we have 
\begin{eqnarray}\label{assoc-4-1}
\lefteqn{\overline{
\mathcal{A}_{P(z_{1}), P(z_{2})}^{P(z_{12}), P(z_{2})}}
(w_{(1)}\boxtimes_{P(z_{1})}(w_{(2)}\boxtimes_{P(z_{2})}
(w_{(3)}\boxtimes_{P(z_{3})}w_{(4)})))}\nn
&&=(w_{(1)}\boxtimes_{P(z_{12})}w_{(2)})\boxtimes_{P(z_{2})}
(w_{(3)}\boxtimes_{P(z_{3})}w_{(4)}),
\end{eqnarray}
where $\overline{
\mathcal{A}_{P(z_{1}), P(z_{2})}^{P(z_{12}), P(z_{2})}}$
is the natural extension of 
$\mathcal{A}_{P(z_{1}), P(z_{2})}^{P(z_{12}), P(z_{2})}$
to (\ref{tr-4-mod-p}).

\item When
$|z_{2}|>|z_{12}|+|z_{3}|>0$, $|z_{3}|>|z_{12}|+|z_{23}|>0$
and $|z_{23}|>|z_{12}|>0$, 
we have 
\begin{eqnarray}\label{assoc-4-2}
\lefteqn{\overline{
\mathcal{A}_{P(z_{2}), P(z_{3})}^{P(z_{23}), P(z_{3})}}
((w_{(1)}\boxtimes_{P(z_{12})}w_{(2)})\boxtimes_{P(z_{2})}
(w_{(3)}\boxtimes_{P(z_{3})}w_{(4)}))}\nn
&&
=((w_{(1)}\boxtimes_{P(z_{12})}w_{(2)})\boxtimes_{P(z_{23})}
w_{(3)})\boxtimes_{P(z_{3})}w_{(4)},
\end{eqnarray}
where $\overline{
\mathcal{A}_{P(z_{2}), P(z_{3})}^{P(z_{23}), P(z_{3})}}$
is the natural extension of 
$\mathcal{A}_{P(z_{2}), P(z_{3})}^{P(z_{23}), P(z_{3})}$
to (\ref{tr-4-mod-comp}).

\item When
$|z_{1}|>|z_{2}|>|z_{3}|>|z_{23}|>0$, $|z_{1}|>|z_{3}|+|z_{23}|>0$,
we have 
\begin{eqnarray}\label{assoc-4-3}
&{\displaystyle \overline{(1_{W_{1}}\boxtimes_{P(z_{1})}
\mathcal{A}_{P(z_{2}), P(z_{3})}^{P(z_{2}-z_{3}), P(z_{3})})}
(w_{(1)}\boxtimes_{P(z_{1})}(w_{(2)}\boxtimes_{P(z_{2})}
(w_{(3)}\boxtimes_{P(z_{3})}w_{(4)})))}&\nn
&{\displaystyle 
=w_{(1)}\boxtimes_{P(z_{1})}((w_{(2)}\boxtimes_{P(z_{23})}
w_{(3)})\boxtimes_{P(z_{3})}w_{(4)}),}&
\end{eqnarray}
where $\overline{(1_{W_{1}}\boxtimes_{P(z_{1})}
\mathcal{A}_{P(z_{2}), P(z_{3})}^{P(z_{2}-z_{3}), P(z_{3})})}$
is the natural extension of 
$1_{W_{1}}\boxtimes_{P(z_{1})}
\mathcal{A}_{P(z_{2}), P(z_{3})}^{P(z_{23}), P(z_{3})}$
to (\ref{tr-4-mod-p}).

\item When
$|z_{3}|>|z_{13}|>|z_{23}|>0$ and $|z_{1}|>|z_{3}|+|z_{23}|>0$,
we have 
\begin{eqnarray}\label{assoc-4-4}
\lefteqn{\overline{
\mathcal{A}_{P(z_{1}), P(z_{3})}^{P(z_{13}), P(z_{3})}}
(w_{(1)}\boxtimes_{P(z_{1})}((w_{(2)}\boxtimes_{P(z_{23})}
w_{(3)})\boxtimes_{P(z_{3})}w_{(4)}))}\nn
&&
=(w_{(1)}\boxtimes_{P(z_{13})}(w_{(2)}\boxtimes_{P(z_{23})}
w_{(3)}))\boxtimes_{P(z_{3})}w_{(4)},
\end{eqnarray}
where $\overline{
\mathcal{A}_{P(z_{1}), P(z_{3})}^{P(z_{13}), P(z_{3})}}$
is the natural extension of 
$\mathcal{A}_{P(z_{1}), P(z_{3})}^{P(z_{13}), P(z_{3})}$
to (\ref{tr-4-mod-ip}).

\item When
$|z_{3}|>|z_{13}|>|z_{23}|>|z_{12}|>0$ and  $|z_{3}|>|z_{12}|+|z_{23}|>0$, 
we have 
\begin{eqnarray}\label{assoc-4-5}
&{\displaystyle \overline{
(\mathcal{A}_{P(z_{13}), P(z_{23})}^{P(z_{12}), P(z_{23})}
\boxtimes_{P(z_{3})}1_{W_{4}})}
((w_{(1)}\boxtimes_{P(z_{13})}(w_{(2)}\boxtimes_{P(z_{23})}
w_{(3)}))\boxtimes_{P(z_{3})}w_{(4)})}&\nn
&{\displaystyle 
=((w_{(1)}\boxtimes_{P(z_{12})}w_{(2)})\boxtimes_{P(z_{23})}
w_{(3)})\boxtimes_{P(z_{3})}w_{(4)},}&
\end{eqnarray}
where $\overline{
(\mathcal{A}_{P(z_{13}), P(z_{23})}^{P(z_{12}), P(z_{23})}
\boxtimes_{P(z_{3})}1_{W_{4}})}$
is the natural extension of 
$(\mathcal{A}_{P(z_{13}), P(z_{23})}^{P(z_{12}), P(z_{23})}
\boxtimes_{P(z_{3})}1_{W_{4}})$
to (\ref{tr-4-mod-pi}).

\end{enumerate}
\end{propo}
\pf
To prove (\ref{assoc-4-1}), 
we note that when $|z_{1}|>|z_{2}|>|z_{12}|+|z_{3}|>0$,
by Proposition \ref{t-prod-4-elts}, for 
\[
w'\in (W_{1}\boxtimes_{P(z_{12})}W_{2})\boxtimes_{P(z_{2})}
(W_{3}\boxtimes_{P(z_{3})}W_{4}),
\]
we have
\begin{eqnarray*}
\lefteqn{\langle w', w_{(1)}\boxtimes_{P(z_{1})}(w_{(2)}\boxtimes_{P(z_{2})}
(w_{(3)}\boxtimes_{P(z_{3})}w_{(4)}))\rangle}\nn
&&=\sum_{m\in \C}
\langle w', w_{(1)}\boxtimes_{P(z_{1})}(w_{(2)}\boxtimes_{P(z_{2})}
\pi_{m}(w_{(3)}\boxtimes_{P(z_{3})}w_{(4)}))\rangle
\end{eqnarray*}
and
\begin{eqnarray*}
\lefteqn{\langle w', (w_{(1)}\boxtimes_{P(z_{12})}w_{(2)})\boxtimes_{P(z_{2})}
(w_{(3)}\boxtimes_{P(z_{3})}w_{(4)})\rangle}\nn
&&=\sum_{m\in \C}\langle w', (w_{(1)}\boxtimes_{P(z_{12})}w_{(2)})
\boxtimes_{P(z_{2})}
\pi_{m}(w_{(3)}\boxtimes_{P(z_{3})}w_{(4)})\rangle.
\end{eqnarray*}
Let $\left(\mathcal{A}_{P(z_{1}), P(z_{2})}^{P(z_{12}), P(z_{2})}\right)'$
be the adjoint of $\mathcal{A}_{P(z_{1}), P(z_{2})}^{P(z_{12}), P(z_{2})}$.
Then
\begin{eqnarray*}
\lefteqn{\left\langle w', 
\overline{\mathcal{A}_{P(z_{1}), P(z_{2})}^{P(z_{12}), P(z_{2})}}
(w_{(1)}\boxtimes_{P(z_{1})}(w_{(2)}\boxtimes_{P(z_{2})}
(w_{(3)}\boxtimes_{P(z_{3})}w_{(4)})))\right\rangle}\nn
&&=\left\langle 
\left(\mathcal{A}_{P(z_{1}), P(z_{2})}^{P(z_{12}), P(z_{2})}\right)'(w'), 
\left(w_{(1)}\boxtimes_{P(z_{1})}(w_{(2)}\boxtimes_{P(z_{2})}
(w_{(3)}\boxtimes_{P(z_{3})}w_{(4)}))\right)\right\rangle\nn
&&=\sum_{m\in \C}\left\langle 
\left(\mathcal{A}_{P(z_{1}), P(z_{2})}^{P(z_{12}), P(z_{2})}\right)'(w'), 
\left(w_{(1)}\boxtimes_{P(z_{1})}(w_{(2)}\boxtimes_{P(z_{2})}
\pi_{m}(w_{(3)}\boxtimes_{P(z_{3})}w_{(4)}))\right)\right\rangle\nn
&&=\sum_{m\in \C}\left\langle w', 
\overline{\mathcal{A}_{P(z_{1}), P(z_{2})}^{P(z_{12}), P(z_{2})}}(
w_{(1)}\boxtimes_{P(z_{1})}(w_{(2)}\boxtimes_{P(z_{2})}
\pi_{m}(w_{(3)}\boxtimes_{P(z_{3})}w_{(4)})))\right\rangle\nn
&&=\sum_{m\in \C}\left\langle w', 
(w_{(1)}\boxtimes_{P(z_{12})}w_{(2)})\boxtimes_{P(z_{2})}
\pi_{m}(w_{(3)}\boxtimes_{P(z_{3})}w_{(4)})\right\rangle\nn
&&=\langle w', (w_{(1)}\boxtimes_{P(z_{12})}w_{(2)})\boxtimes_{P(z_{2})}
(w_{(3)}\boxtimes_{P(z_{3})}w_{(4)})\rangle.
\end{eqnarray*}
Since $w'$ is arbitrary, we obtain (\ref{assoc-4-1}).

The equalities (\ref{assoc-4-2})--(\ref{assoc-4-5})
can be proved similarly.
\epfv

\subsection{The data of the braided tensor category structure}

We choose the tensor product bifunctor of the braided tensor category
that we are constructing to be the bifunctor $\boxtimes_{P(1)}$. 
We shall denote $\boxtimes_{P(1)}$ simply by $\boxtimes$.
We take the unit object to be $V$.  For any $z\in \C^{\times}$ and any
$V$-module $(W, Y_{W})$, we take the {\it left $P(z)$-unit
isomorphism} $l_{W; z}: V\boxtimes_{P(z)}W \to W$ to be the unique
module map {}from $V\boxtimes_{P(z)}W$ to $W$ such that
\[
\overline{l_{W; z}}\circ \boxtimes_{P(z)}=I_{Y_{W}, 0},
\]
where $I_{Y_{W}, 0}=I_{Y_{W}, p}$ for $p\in \mathbb{Z}$
is the unique $P(z)$-intertwining map associated to 
the intertwining operator $Y_{W}$ of type ${W\choose VW}$. 
The existence and uniqueness of $l_{W; z}$ are guaranteed by the 
universal property of the $P(z_{1})$-tensor product 
$\boxtimes_{P(z_{1})}$.
It is
characterized by 
\[
l_{W;
z}(\mathbf{1}\boxtimes_{P(z)} w)=w
\]
for $w\in W$. The {\it right
$P(z)$-unit isomorphism} $r_{W; z}: W\boxtimes_{P(z)}V \to W$ is
the unique module map 
{}from $W\boxtimes_{P(z)}V$ to $W$ such that 
\[
\overline{r_{W; z}}\circ \boxtimes_{P(z)}=I_{\Omega_{0}(Y_{W}), 0},
\]
where $I_{\Omega_{0}(Y_{W}), 0}=I_{\Omega_{0}(Y_{W}), 0}$ for $p\in \mathbb{Z}$ 
is the unique $P(z)$-intertwining map associated to 
the intertwining operator $\Omega_{0}(Y_{W})$ of type ${W\choose WV}$. 
It is
characterized by 
\[
\overline{r_{W; z}}(w\boxtimes_{P(z)}
\mathbf{1})=e^{zL(-1)} w
\]
for $w\in W$.  In particular, we have the
left unit isomorphism $l_{W}=l_{W; 1}: V\boxtimes W \to W$ and the right
unit isomorphism $r_{W}=r_{W; 1}: W\boxtimes V \to W$.

To give the braiding and associativity isomorphisms, we need ``parallel
transport isomorphisms" between $P(z)$-tensor products with different
$z$.  Let $W_{1}$ and $W_{2}$ be objects of $\mathcal{C}$ and $z_{1},
z_{2}\in \C^{\times}$.  Giving a path $\gamma$ in
$\mathbb{C}^{\times}$ {}from $z_{1}$ to $z_{2}$.  Let $\mathcal{Y}$ be
the logarithmic intertwining operator associated to the
$P(z_{2})$-tensor product $W_{1}\boxtimes_{P(z_{2})}W_{2}$ and
$l(z_{1})$ the value of the logarithm of $z_{1}$ determined uniquely
by $\log z_{2}$ and the path $\gamma$. Then we have a
$P(z_{1})$-intertwining map $I$ defined by
\[
I(w_{(1)}\otimes w_{(2)})=\mathcal{Y}(w_{(1)}, e^{l(z_{1})})w_{(2)}
\]
for $w_{(1)}\in W_{1}$ and $w_{(2)}\in W_{2}$. The {\it parallel transport 
isomorphism
$\mathcal{T}_{\gamma}: W_{1}\boxtimes_{P(z_{1})}W_{2}\to
W_{1}\boxtimes_{P(z_{2})}W_{2}$ associated to $\gamma$}
is defined to be 
the unique module map such that 
\[
I=\overline{\mathcal{T}_{\gamma}}\circ \boxtimes_{P(z_{1})}.
\]
where $\overline{\mathcal{T}_{\gamma}}$ is the natural extension of
$\mathcal{T}_{\gamma}$ to the algebraic completion
$\overline{W_{1}\boxtimes_{P(z_{1})} W_{2}}$ of
$W_{1}\boxtimes_{P(z_{1})} W_{2}$.  As in the definition of the left
$P(z)$-unit isomorphism, the existence and uniqueness is guaranteed by
the universal property of the $P(z_{1})$-tensor product
$\boxtimes_{P(z_{1})}$.  The parallel transport isomorphism
$\mathcal{T}_{\gamma}$ is characterized by
\[
\overline{\mathcal{T}_{\gamma}}(w_{(1)}
\boxtimes_{P(z_{1})}w_{(2)})=\mathcal{Y}(w_{(1)}, x)w_{(2)}
|_{\log x=l(z_{1}),\;x^{n}=e^{nl(z_{1})}, \; n\in \C}
\]
for $w_{(1)}\in W_{1}$ and $w_{(2)}\in W_{2}$.  Since the intertwining
map $I$ depends only on the homotopy class of $\gamma$, {}from the
definition, we see that the parallel parallel isomorphism also depends
only on the homotopy class of $\gamma$.

For $z\in \C^{\times}$, let $I$ be the $P(z)$-intertwining map of type
${W_{2}\boxtimes_{P(-z)} W_{1} \choose W_{1}W_{2}}$ defined by
\[
I(w_{1}\otimes w_{2})= e^{zL(-1)}
(w_{2}\boxtimes_{P(-z)}w_{1})
\]
for $w_{1}\in W_{1}$ and $w_{2}\in W_{2}$.  We define a {\it
commutativity isomorphism between the $P(z)$- and $P(-z)$-tensor
products} to be the unique module map
\[
\mathcal{R}_{P(z)}: W_{1}\boxtimes_{P(z)} W_{2}\to 
W_{2}\boxtimes_{P(-z)} W_{1}
\]
such that 
\[
I=\overline{\mathcal{R}_{P(z)}}\circ \boxtimes_{P(z)}.
\]
The existence and uniqueness of $\mathcal{R}_{P(z)}$ is guaranteed by
the universal property of the $P(z)$-tensor product.  By definition,
the commutativity isomorphism $\mathcal{R}_{P(z)}$ is characterized by
\[
\overline{\mathcal{R}_{P(z)}}(w_{(1)}\boxtimes_{P(z)} w_{(2)})=e^{zL(-1)}
(w_{(2)}\boxtimes_{P(-z)} w_{(1)})
\]
for $w_{(1)}\in W_{1}$, $w_{(2)}\in W_{2}$.

Let $\gamma_{1}^{-}$ be a path
{}from $-1$ to $1$  in the
closed upper half plane with $0$ deleted,
$\mathcal{T}_{\gamma_{1}^{-}}$ the corresponding
parallel transport isomorphism.
We define the {\it braiding isomorphism} 
\[
\mathcal{R}: W_{1}\boxtimes W_{2}\to
W_{2}\boxtimes W_{1}
\]
for our
braided tensor category to be 
\[
\mathcal{R}=\mathcal{T}_{\gamma_{1}^{-}}\circ \mathcal{R}_{P(1)}.
\]
This braiding isomorphism $\mathcal{R}$ can also be defined directly
as follows: Let $I$ be the $P(1)$-intertwining map of type
${W_{2}\boxtimes_{P(1)} W_{1} \choose W_{1}W_{2}}$ defined by
\[
I(w_{1}\otimes w_{2})= e^{zL(-1)}
\overline{\mathcal{T}}_{\gamma_{1}^{-}} (w_{2}\boxtimes_{P(1)} w_{1})
\]
for $w_{1}\in W_{1}$ and $w_{2}\in W_{2}$.  Then $\mathcal{R}$ is the
unique module map
\[
\mathcal{R}: W_{1}\boxtimes W_{2}\to 
W_{2}\boxtimes W_{1}
\]
such that 
\[
I=\overline{\mathcal{R}}\circ \boxtimes.
\]
It is characterized by
\[
\overline{\mathcal{R}}(w_{(1)}\boxtimes_{P(1)} w_{(2)})=e^{zL(-1)}
\overline{\mathcal{T}}_{\gamma_{1}^{-}}
(w_{(2)}\boxtimes_{P(1)} w_{(1)})
\]
for $w_{(1)}\in W_{1}$, $w_{(2)}\in W_{2}$.

Let $z_{1}, z_{2}$ be complex numbers satisfying
$|z_{1}|>|z_{2}|>|z_{1}-z_{2}|>0$ and let $W_{1}$, $W_{2}$ and $W_{3}$
$V$-modules. Then the associativity isomorphism (corresponding to the
indicated geometric data)
\[
\alpha_{P(z_{1}), P(z_{2})}^{P(z_{1}-z_{2}), P(z_{2})}:
(W_{1}\boxtimes_{P(z_{1}-z_{2})}W_{2})\boxtimes_{P(z_{2})}W_{3}
\to W_{1}\boxtimes_{P(z_{1})}(W_{2}\boxtimes_{P(z_{2})}W_{3})
\]
and its inverse
\[
\A_{P(z_{1}), P(z_{2})}^{P(z_{1}-z_{2}), P(z_{2})}:
W_{1}\boxtimes_{P(z_{1})}(W_{2}\boxtimes_{P(z_{2})}W_{3})
\to (W_{1}\boxtimes_{P(z_{1}-z_{2})}W_{2})\boxtimes_{P(z_{2})}W_{3}
\]
have been constructed in Section 10 and are determined uniquely by
(\ref{assoc-elt-2}) and (\ref{assoc-elt-1}), respectively.

Let $z_{1}, z_{2}, z_{3}$ and $z_{4}$ be any nonzero complex numbers  
and let $W_{1}$, $W_{2}$ and $W_{3}$ be $V$-modules. We also define
a natural associativity isomorphism
\[
\A_{P(z_{1}), P(z_{2})}^{P(z_{4}), P(z_{3})}:
W_{1}\boxtimes_{P(z_{1})}(W_{2}\boxtimes_{P(z_{2})}W_{3})
\to (W_{1}\boxtimes_{P(z_{4})}W_{2})\boxtimes_{P(z_{3})}W_{3}
\]
using the associativity and parallel transport isomorphisms as
follows: Let $\zeta_{1}$ and $\zeta_{2}$ be nonzero complex numbers
satisfying $|\zeta_{1}|>|\zeta_{2}|>|\zeta_{1}-\zeta_{2}|>0$.  Let
$\gamma_{1}$ and $\gamma_{2}$, be paths {}from $z_{1}$ and $z_{2}$ to
$\zeta_{1}$ and $\zeta_{2}$, respectively, in the complex plane with a
cut along the positive real line, and $\gamma_{3}$ and $\gamma_{4}$ be
paths {}from  $\zeta_{2}$ and
$\zeta_{1}-\zeta_{2}$ to $z_{3}$ and $z_{4}$, 
respectively, also in the complex plane with a
cut along the positive real line.  Then
\[
\A_{P(z_{1}), P(z_{2})}^{P(z_{4}), P(z_{3})}
=\mathcal{T}_{\gamma_{3}}\circ (\mathcal{T}_{\gamma_{4}}
\boxtimes_{P(\zeta_{2})} I_{W_{3}})\circ
\mathcal{A}^{P(\zeta_{1}-\zeta_{2}), P(\zeta_{2})}_{P(\zeta_{1}), P(\zeta_{2})}\circ
(I_{W_{1}} \boxtimes_{P(\zeta_{1})}
\mathcal{T}_{\gamma_{2}})\circ \mathcal{T}_{\gamma_{1}},
\]
that is, $\A_{P(z_{1}), P(z_{2})}^{P(z_{4}), P(z_{3})}$ 
is given by the commutative diagram
\[
\begin{CD}
W_{1}\boxtimes_{P(\zeta_{1})} (W_{2}
\boxtimes_{P(\zeta_{2})} W_{3})
@>\mathcal{A}_{P(\zeta_{1}), P(\zeta_{2})}^{P(\zeta_{1}-z_{2}), P(\zeta_{2})}>>
(W_{1}\boxtimes_{P(\zeta_{1}-\zeta_{2})} W_{2})
\boxtimes_{P(\zeta_{2})} W_{3}\\
@A(I_{W_{1}} \boxtimes_{P(\zeta_{1})}
\mathcal{T}_{\gamma_{2}})\circ \mathcal{T}_{\gamma_{1}}AA
@VV\mathcal{T}_{\gamma_{3}}\circ (\mathcal{T}_{\gamma_{4}}
\boxtimes_{P(\zeta_{2})} I_{W_{3}})V\\
W_{1}\boxtimes_{P(z_{1})} (W_{2}
\boxtimes_{P(z_{2})} W_{3}) @>\mathcal{A}_{P(z_{1}), P(z_{2})}^{P(z_{4}), P(z_{3})}>>
(W_{1}\boxtimes_{P(z_{4})} W_{2})
\boxtimes_{P(z_{3})} W_{3}
\end{CD}
\]
The inverse of $\A_{P(z_{1}), P(z_{2})}^{P(z_{4}), P(z_{3})}$ 
is denoted $\alpha_{P(z_{1}), P(z_{2})}^{P(z_{4}), P(z_{3})}$.

In particular, when $z_{1}=z_{2}=z_{3}=z_{4}=1$, we
call the corresponding associativity isomorphism 
\[
\alpha_{P(1), P(1)}^{P(1), P(1)}:
(W_{1}\boxtimes W_{2})\boxtimes W_{3}\to
W_{1}\boxtimes (W_{2}\boxtimes W_{3})
\]
{\it the associativity isomorphism (for the braided tensor category
structure)} and denote it simply by $\alpha$. Its inverse is denoted
by $\A$.  Because of the importance of this special case, we rewrite
the definition of $\A$ explicitly. Let $r_{1}$ and $r_{2}$ be real
numbers satisfying $r_{1}>r_{2}>r_{1}-r_{2}\ge 0$.  Let $\gamma_{1}$
and $\gamma_{2}$ be paths in $(0, \infty)$ {}from $1$ to $r_{1}$ and
to $r_{2}$, respectively, and let $\gamma_{3}$ and $\gamma_{4}$ be
paths in $(0, \infty)$ {}from $r_{2}$ and {}from $r_{1}-r_{2}$ to $1$,
respectively.  Then
\[
\mathcal{A}=\mathcal{T}_{\gamma_{3}}\circ (\mathcal{T}_{\gamma_{4}}
\boxtimes_{P(r_{2})} I_{W_{3}})\circ
\mathcal{A}^{P(r_{1}-r_{2}), P(r_{2})}_{P(r_{1}), P(r_{2})}\circ
(I_{W_{1}} \boxtimes_{P(r_{1})}
\mathcal{T}_{\gamma_{2}})\circ \mathcal{T}_{\gamma_{1}},
\]
that is, $\mathcal{A}$ is given by the commutative diagram
\[
\begin{CD}
W_{1}\boxtimes_{P(r_{1})} (W_{2}
\boxtimes_{P(r_{2})} W_{3})
@>\mathcal{A}_{P(r_{1}), P(r_{2})}^{P(r_{1}-r_{2}), P(r_{2})}>>
(W_{1}\boxtimes_{P(r_{1}-r_{2})} W_{2})
\boxtimes_{P(r_{2})} W_{3}\\
@A(I_{W_{1}} \boxtimes_{P(r_{1})}
\mathcal{T}_{\gamma_{2}})\circ \mathcal{T}_{\gamma_{1}}AA
@VV\mathcal{T}_{\gamma_{3}}\circ (\mathcal{T}_{\gamma_{4}}
\boxtimes_{P(r_{2})} I_{W_{3}})V\\
W_{1}\boxtimes (W_{2}
\boxtimes W_{3}) @>\mathcal{A}>>
(W_{1}\boxtimes W_{2})
\boxtimes W_{3}.
\end{CD}
\]

\subsection{Actions of the associativity and
commutativity isomorphisms on tensor products of elements}

\begin{propo}
For any $z_{1}, z_{2}\in \C^{\times}$ such that 
$z_{1}\ne z_{2}$ but $|z_{1}|=|z_{2}|=|z_{1}-z_{2}|$,
we have 
\begin{equation}\label{assoc-general-z}
\overline{\A_{P(z_{1}), P(z_{2})}^{P(z_{1}-z_{2}), P(z_{2})}}
(w_{(1)}\boxtimes_{P(z_{1})}(w_{(2)}\boxtimes_{P(z_{2})}w_{(3)}))
=(w_{(1)}\boxtimes_{P(z_{1}-z_{2})}w_{(2)})\boxtimes_{P(z_{2})}w_{(3)}
\end{equation}
for $w_{(1)}\in W_{1}$, $w_{(2)}\in W_{2}$ and $w_{(3)}\in W_{3}$,
where
\[
\overline{\A_{P(z_{1}), P(z_{2})}^{P(z_{1}-z_{2}), P(z_{2})}}:
\overline{W_{1}\boxtimes_{P(z_{1})}(W_{2}\boxtimes_{P(z_{2})}W_{3})}
\to \overline{(W_{1}\boxtimes_{P(z_{1}-z_{2})}W_{2})
\boxtimes_{P(z_{2})}W_{3}}
\]
is the natural extension of
$\A_{P(z_{1}), P(z_{2})}^{P(z_{1}-z_{2}), P(z_{2})}$.
\end{propo}
\pf
We need only prove the case that $w_{(1)}$ and $w_{(2)}$ 
are homogeneous with respect to the generalized-weight grading.
So we now assume that they are homogeneous. 

We can always find $\epsilon_{1}\in \C$
such that 
\begin{eqnarray}
|z_{1}+\epsilon_{1}|&>&|\epsilon_{1}|,\label{assoc-general-z-0}\\
|z_{1}+\epsilon_{1}|&>&|z_{2}|
>|(z_{1}+\epsilon_{1})-z_{2}|>0.\label{assoc-general-z-1}
\end{eqnarray}
Let 
$\Y_{1}=\Y_{\boxtimes_{P(z_{1})}, 0}$ and 
$\Y_{2}=\Y_{\boxtimes_{P(z_{2})}, 0}$
be intertwining operators of types 
\[
{W_{1}\boxtimes_{P(z_{1})}(W_{2}\boxtimes_{P(z_{2})}W_{3})\choose 
W_{1}\;W_{2}\boxtimes_{P(z_{2})}W_{3}}
\]
and 
\[
{W_{2}\boxtimes_{P(z_{2})}W_{3}\choose W_{2}\; W_{3}},
\]
respectively, corresponding to the intertwining maps 
$\boxtimes_{P(z_{1})}$ and $\boxtimes_{P(z_{2})}$, respectively. 
Then the series 
\[
\langle w', \Y_{1}(\pi_{m}(e^{-\epsilon_{1}L(-1)}w_{(1)}), z_{1}+\epsilon_{1})
\Y_{2}(\pi_{n}(e^{-\epsilon L(-1)}w_{(2)}), z_{2})w_{(3)}\rangle
\]
is absolutely convergent for $m, n\in \C$ and 
\[
w'\in 
(W_{1}\boxtimes_{P(z_{1})}(W_{2}\boxtimes_{P(z_{2})}W_{3}))'
\]
 and the sums of these series define elements
\[
\Y_{1}(\pi_{m}(e^{-\epsilon_{1}L(-1)}w_{(1)}), z_{1}+\epsilon_{1})
\Y_{2}(w_{(2)}), z_{2})w_{(3)}
\in \overline{W_{1}\boxtimes_{P(z_{1})}(W_{2}\boxtimes_{P(z_{2})}W_{3})}.
\]

By the definition of the parallel transport isomorphism, 
for any path $\gamma$ {}from $z_{1}+\epsilon_{1}$ to $z_{1}$
 in the complex
plane with a cut along the nonnegative real line, we have
\begin{eqnarray}\label{assoc-general-z-2}
\lefteqn{\overline{\mathcal{T}_{\gamma}}
(\pi_{m}(e^{-\epsilon_{1}L(-1)}w_{(1)}) \boxtimes_{P(z_{1}+\epsilon_{1})}
(w_{(2)}\boxtimes_{P(z_{2})}
w_{(3)}))}\nn
&&=\Y_{1}(\pi_{m}(e^{-\epsilon_{1}L(-1)}w_{(1)}), z_{1}+\epsilon_{1})
\Y_{2}(w_{(2)}), z_{2})w_{(3)}.
\end{eqnarray}
By definition, we know that 
\[
\mathcal{T}_{\gamma}^{-1}=\mathcal{T}_{\gamma^{-1}},
\]
so that (\ref{assoc-general-z-2}) can be written as 
\begin{eqnarray}\label{assoc-general-z-3}
\lefteqn{\overline{ \mathcal{T}_{\gamma_{1}^{-1}}}
(\Y_{1}(\pi_{m}(e^{-\epsilon_{1}L(-1)}w_{(1)}), z_{1}+\epsilon_{1})
\Y_{2}(w_{(2)}), z_{2})w_{(3)})}\nn
&&=\pi_{m}(e^{-\epsilon_{1}L(-1)}w_{(1)})\boxtimes_{P(z_{1}+\epsilon_{1})}
(w_{(2)}\boxtimes_{P(z_{2})}
w_{(3)}).
\end{eqnarray}

Since (\ref{assoc-general-z-1}) holds, by (\ref{assoc-elt-1}), we have
\begin{eqnarray}\label{assoc-general-z-4}
\lefteqn{\overline{\mathcal{A}_{P(z_{1}+\epsilon_{1}), P(z_{2})}
^{P((z_{1}+\epsilon_{1})-z_{1}), P(z_{2})}}
(\pi_{m}(e^{-\epsilon_{1}L(-1)}w_{(1)})\boxtimes_{P(z_{1}+\epsilon_{1})}
(w_{(2)}\boxtimes_{P(z_{2})}
w_{(3)}))}\nn
&&=(\pi_{m}(e^{-\epsilon_{1}L(-1)}w_{(1)})
\boxtimes_{P((z_{1}+\epsilon_{1})-z_{2})}
w_{(2)})\boxtimes_{P(z_{2})}
w_{(3)}.\quad\quad\quad
\end{eqnarray}

Let 
$\Y_{3}=\Y_{\boxtimes_{P(z_{2})}, 0}$ and 
$\Y_{2}=\Y_{\boxtimes_{P(z_{1}-z_{2})}, 0}$
be intertwining operators of types 
\[
{(W_{1}\boxtimes_{P(z_{1}-z_{2})}W_{2})\boxtimes_{P(z_{2})}W_{3}\choose 
(W_{1}\boxtimes_{P(z_{1}-z_{2})}W_{2})\;\;\;\;W_{3}}
\]
and 
\[
{W_{1}\boxtimes_{P(z_{1}-z_{2})}W_{2}\choose W_{1}\;\;\;\; W_{2}},
\]
respectively, corresponding to the intertwining maps 
$\boxtimes_{P(z_{2})}$ and $\boxtimes_{P(z_{1}-z_{2})}$, respectively. 
Then the series 
\[
\langle \tilde{w}', \Y_{3}(\Y_{4}(\pi_{m}(e^{-\epsilon_{1}L(-1)}w_{(1)}), 
(z_{1}+\epsilon_{1})-z_{2})
w_{(2)}), z_{2})w_{(3)}\rangle
\]
is absolutely convergent for $m, n\in \C$ and 
\[
\tilde{w}'\in 
((W_{1}\boxtimes_{P(z_{1}-z_{2})}W_{2})\boxtimes_{P(z_{2})}W_{3})',
\]
 and the sums of these series define elements
\begin{eqnarray*}
&\Y_{3}(\Y_{4}(\pi_{m}(e^{-\epsilon_{1}L(-1)}w_{(1)}), 
(z_{1}+\epsilon_{1})-z_{2})
w_{(2)}), z_{2})w_{(3)}&\nn
& \in \overline{(W_{1}\boxtimes_{P(z_{1}-z_{2})}W_{2})
\boxtimes_{P(z_{2})}W_{3}}.&
\end{eqnarray*}

We can always choose $\gamma$ such that 
the path $\gamma-z_{2}$ {}from  $(z_{1}+\epsilon_{1})
-z_{2}$ to $z_{1}-z_{2}$ is also
in the complex
plane with a cut along the nonnegative real line. Then
by the definition of the parallel transport isomorphism,  we have
\begin{eqnarray}\label{assoc-general-z-5}
\lefteqn{\overline{\mathcal{T}_{\gamma-z_{2}}}
((\pi_{m}(e^{-\epsilon_{1}L(-1)}w_{(1)})\boxtimes_{P((z_{1}+\epsilon_{1})
-z_{2})}w_{(2)})\boxtimes_{P(z_{2})}
w_{(3)})}\nn
&&=\Y_{3}(\Y_{4}(\pi_{m}(e^{-\epsilon_{1}L(-1)}w_{(1)}), (z_{1}+\epsilon_{1})
-z_{2})w_{(2)}), z_{2})w_{(3)}.\nn
\end{eqnarray}

Combining (\ref{assoc-general-z-3})--(\ref{assoc-general-z-5}) and 
using the definition of the associativity isomorphism
$\A_{P(z_{1}), P(z_{2})}^{P(z_{1}-z_{2}), P(z_{2})}$,
we obtain
\begin{eqnarray}\label{assoc-general-z-6}
\lefteqn{\overline{\A_{P(z_{1}), P(z_{2})}
^{P(z_{1}-z_{2}), P(z_{2})}}
(\Y_{1}(\pi_{m}(e^{-\epsilon_{1}L(-1)}w_{(1)}), z_{1}+\epsilon_{1})
\Y_{2}(w_{(2)}), z_{2})w_{(3)})}\nn
&&=\Y_{3}(\Y_{4}(\pi_{m}(e^{-\epsilon_{1}L(-1)}w_{(1)}), (z_{1}+\epsilon_{1})
-z_{2})w_{(2)}), z_{2})w_{(3)}.\nn
\end{eqnarray}

Since (\ref{assoc-general-z-1}) holds, 
both the series 
\[
\sum_{m\in \C}\langle w',
\Y_{1}(\pi_{m}(e^{-\epsilon_{1}L(-1)}w_{(1)}), z_{1}+\epsilon_{1})
\Y_{2}(w_{(2)}), z_{2})w_{(3)}
\rangle
\]
and 
\[
\sum_{m\in \C}\langle \tilde{w}', 
\Y_{3}(\Y_{4}(\pi_{m}(e^{-\epsilon_{1}L(-1)}w_{(1)}), (z_{1}+\epsilon_{1})
-z_{2})
w_{(2)}), z_{2})w_{(3)}\rangle
\]
are absolutely convergent for 
\[
w'\in 
(W_{1}\boxtimes_{P(z_{1})}(W_{2}\boxtimes_{P(z_{2})}W_{3}))'
\]
and 
\[
\tilde{w}'\in 
((W_{1}\boxtimes_{P(z_{1}-z_{2})}W_{2})\boxtimes_{P(z_{2})}W_{3})'.
\]

We know that
\[
\langle \tilde{w}', 
\Y_{1}(w_{(1)}), z_{1}+\epsilon_{1})
\Y_{2}(w_{(2)}), z_{2})w_{(3)}
\rangle
\]
and 
\[
\langle \tilde{w}', 
\Y_{3}(\Y_{4}(w_{(1)}), (z_{1}+\epsilon_{1})
-z_{2})
w_{(2)}), z_{2})w_{(3)}\rangle
\]
is the values of  single-valued analytic 
functions
\[
F(w', w_{(1)}, w_{(2)}, w_{(3)}; \zeta_{1}, \zeta_{2})
\]
and 
\[
G(\tilde{w}', w_{(1)}, w_{(2)}, w_{(3)}; \zeta_{1}, \zeta_{2})
\]
of $\zeta_{1}$ and $\zeta_{2}$
in a neighborhood of the point $(\zeta_{1}, \zeta_{2})=(z_{1}+\epsilon_{1},
z_{2})$ which contains the point $(\zeta_{1}, \zeta_{2})=(z_{1}, z_{2})$. 
Then by the definition of the 
tensor product elements
$w_{(1)}\boxtimes_{P(z_{1})}(w_{(2)}\boxtimes_{P(z_{2})}w_{(3)})$
and $(w_{(1)}\boxtimes_{P(z_{1}-z_{2})}w_{(2)})\boxtimes_{P(z_{2})}w_{(3)}$,
we have 
\begin{equation}\label{assoc-general-z-6.1}
\langle w', w_{(1)}\boxtimes_{P(z_{1})}(w_{(2)}\boxtimes_{P(z_{2})}w_{(3)})
\rangle
=F(w', w_{(1)}, w_{(2)}, w_{(3)}; z_{1}, z_{2})
\end{equation}
and 
\begin{equation}\label{assoc-general-z-6.2}
\langle \tilde{w}', 
(w_{(1)}\boxtimes_{P(z_{1}-z_{2})}w_{(2)})\boxtimes_{P(z_{2})}w_{(3)}
\rangle
=G(\tilde{w}', w_{(1)}, w_{(2)}, w_{(3)}; z_{1}, z_{2}).
\end{equation}
On the other hand, since 
$F(w', w_{(1)}, w_{(2)}, w_{(3)}; \zeta_{1}, \zeta_{2})$
and 
$G(\tilde{w}', w_{(1)}, w_{(2)}, w_{(3)}; \zeta_{1}, \zeta_{2})$
are analytic extensions of matrix elements of products and 
iterates of intertwining maps, properties of these products and 
iterates also hold for these functions if they still make sense. 
In particular, they satisfy the $L(-1)$-derivative property:
\begin{eqnarray}
\frac{\partial}{\partial \zeta_{1}}
F(w', w_{(1)}, w_{(2)}, w_{(3)}; \zeta_{1}, \zeta_{2})
&=&F(w', L(-1)w_{(1)}, w_{(2)}, w_{(3)}; \zeta_{1}, \zeta_{2}),
\label{assoc-general-z-7}\\
\frac{\partial}{\partial \zeta_{1}}
G(\tilde{w}', w_{(1)}, w_{(2)}, w_{(3)}; \zeta_{1}, \zeta_{2})
&=&G(\tilde{w}', L(-1)w_{(1)}, w_{(2)}, w_{(3)}; \zeta_{1}, \zeta_{2}),
\label{assoc-general-z-9}
\end{eqnarray}
{}From the Taylor theorem (which applies since 
(\ref{assoc-general-z-0}) holds) and 
(\ref{assoc-general-z-7})--(\ref{assoc-general-z-9}), we have 
\begin{eqnarray}
F(w', w_{(1)}, w_{(2)}, w_{(3)}; z_{1}, z_{2})
&=&\sum_{m\in \C}F(w', \pi_{m}(e^{-\epsilon_{1}L(-1)})w_{(1)}, 
w_{(2)}, w_{(3)}; 
z_{1}+\epsilon_{1}, z_{2}),\quad\quad\label{assoc-general-z-11}\\
G(\tilde{w}', w_{(1)}, w_{(2)}, w_{(3)}; z_{1}, z_{2})
&=&\sum_{m\in \C}G(\tilde{w}', \pi_{m}(e^{-\epsilon_{1}L(-1)})w_{(1)}, 
w_{(2)}, w_{(3)}; 
z_{1}+\epsilon_{1}, z_{2}).\quad\quad
\label{assoc-general-z-12}
\end{eqnarray}
Thus by the definitions of 
\begin{eqnarray*}
F(w', \pi_{m}(e^{-\epsilon_{1}L(-1)})w_{(1)}, 
w_{(2)}, w_{(3)}; 
z_{1}+\epsilon_{1}, z_{2}),\\
G(\tilde{w}', \pi_{m}(e^{-\epsilon_{1}L(-1)})w_{(1)}, 
w_{(2)}, w_{(3)}; 
z_{1}+\epsilon_{1}, z_{2}),
\end{eqnarray*}
and by (\ref{assoc-general-z-11}),
(\ref{assoc-general-z-12}), (\ref{assoc-general-z-6.1})
and (\ref{assoc-general-z-6.2}), we obtain
\begin{eqnarray}
\lefteqn{\sum_{m\in \C}\langle w',
\Y_{1}(\pi_{m}(e^{-\epsilon_{1}L(-1)}w_{(1)}), z_{1}+\epsilon_{1})
\Y_{2}(w_{(2)}), z_{2})w_{(3)}
\rangle}\nn
&&\quad\quad
=\langle w', w_{(1)}\boxtimes_{P(z_{1})}(w_{(2)}\boxtimes_{P(z_{2})}w_{(3)})
\rangle
\label{assoc-general-z-13}
\end{eqnarray}
and
\begin{eqnarray}
\lefteqn{\sum_{m\in \C}\langle \tilde{w}', 
\Y_{3}(\Y_{4}(\pi_{m}(e^{-\epsilon_{1}L(-1)}w_{(1)}), (z_{1}+\epsilon_{1})
-z_{2}))w_{(2)}), z_{2})w_{(3)}\rangle}\nn
&&\quad\quad\quad
=\langle \tilde{w}', 
(w_{(1)}\boxtimes_{P(z_{1}-z_{2})}w_{(2)})\boxtimes_{P(z_{2})}w_{(3)}
\rangle.
\label{assoc-general-z-14}
\end{eqnarray}
Since $w'$ and $\tilde{w'}'$ are arbitrary, (\ref{assoc-general-z-13}) and
(\ref{assoc-general-z-14}) gives
\begin{eqnarray}
\lefteqn{\sum_{m\in \C}
\Y_{1}(\pi_{m}(e^{-\epsilon_{1}L(-1)}w_{(1)}), z_{1}+\epsilon_{1})
\Y_{2}(w_{(2)}), z_{2})w_{(3)}}\nn
&&\quad\quad
=w_{(1)}\boxtimes_{P(z_{1})}(w_{(2)}\boxtimes_{P(z_{2})}w_{(3)})
\quad\quad\quad\quad\quad\quad\quad\quad\quad
\label{assoc-general-z-15}
\end{eqnarray}
and
\begin{eqnarray}
\lefteqn{\sum_{m\in \C} 
\Y_{3}(\Y_{4}(\pi_{m}(e^{-\epsilon_{1}L(-1)}w_{(1)}), (z_{1}+\epsilon_{1})
-z_{2})w_{(2)}), z_{2})w_{(3)}}\nn
&&\quad\quad\quad
=(w_{(1)}\boxtimes_{P(z_{1}-z_{2})}w_{(2)})\boxtimes_{P(z_{2})}w_{(3)}.
\quad\quad\quad\quad\quad\quad\quad\quad\quad
\label{assoc-general-z-16}
\end{eqnarray}

Taking the sum $\sum_{m\in \C}$ on both sides of 
(\ref{assoc-general-z-6}) and then 
using (\ref{assoc-general-z-15}) and (\ref{assoc-general-z-16}),
we obtain (\ref{assoc-general-z}).
\epfv

We also have:

\begin{propo}
Let $z_{1}, z_{2}$ be nonzero complex numbers such that 
$z_{1}\ne z_{2}$ but $|z_{1}|=|z_{2}|=|z_{1}-z_{2}|$. Let $\gamma$ be 
a path {}from $z_{2}$ to $z_{1}$ in the complex plane with a cut 
along the nonnegative real line. Then
we have 
\begin{equation}\label{commu-1}
\overline{\mathcal{T}_{\gamma}\circ 
(\mathcal{R}_{P(z_{1}-z_{2})}\boxtimes_{P(z_{2})}1_{W_{3}})}
((w_{(1)}\boxtimes_{P(z_{1}-z_{2})}w_{(2)})\boxtimes_{P(z_{2})}w_{(3)})
=((w_{(2)}\boxtimes_{P(z_{2}-z_{1})}w_{(1)})\boxtimes_{P(z_{1})}w_{(3)}
\end{equation}
for $w_{(1)}\in W_{1}$, $w_{(2)}\in W_{2}$ and $w_{(3)}\in W_{3}$.
\end{propo}
\pf
We can find $\epsilon$ such that $|z_{2}|>|\epsilon|$,
\[
|z_{2}+\epsilon|>|z_{1}-z_{2}|>0.
\]
Let 
$\Y_{1}=\Y_{\boxtimes_{P(z_{2})}, 0}$, 
$\tilde{\Y}_{1}=\Y_{\boxtimes_{P(z_{2})}, 0}$ and
$\Y_{2}=\Y_{\boxtimes_{P(z_{1}-z_{2})}, 0}$
be intertwining operators of types 
\[
{(W_{1}\boxtimes_{P(z_{1}-z_{2})}W_{2})\boxtimes_{P(z_{2})}W_{3}\choose 
(W_{1}\boxtimes_{P(z_{1}-z_{2})}W_{2})\;\;\;\;W_{3}},
\]
\[
{(W_{2}\boxtimes_{P(z_{2}-z_{1})}W_{1})\boxtimes_{P(z_{2})}W_{3}\choose 
(W_{2}\boxtimes_{P(z_{2}-z_{1})}W_{1})\;\;\;\;W_{3}},
\]
and 
\[
{W_{1}\boxtimes_{P(z_{1}-z_{2})}W_{2}\choose W_{1}\;\;\;\; W_{2}},
\]
respectively, corresponding to the intertwining maps 
$\boxtimes_{P(z_{2})}$, $\boxtimes_{P(z_{2})}$
and $\boxtimes_{P(z_{1})}$, respectively. 
Note that 
since $|z_{2}+\epsilon|>|z_{1}-z_{2}|>0$,
\begin{eqnarray}\label{commu-1-0-0}
\lefteqn{\Y_{1}(w_{(1)}\boxtimes_{P(z_{1}-z_{2})}w_{(2)}, 
z_{2}+\epsilon)w_{(3)}}\nn
&&=\sum_{m, n\in \C}\sum_{k=1}^{N}
(\pi_{n}(w_{(1)}\boxtimes_{P(z_{1}-z_{2})}w_{(2)}))^{\Y_{1}}_{m;k}
w_{(3)}e^{(-m-1)\log (z_{2}+\epsilon)}(\log (z_{2}+\epsilon))^{k}.
\end{eqnarray}
For $n\in \C$, 
\begin{eqnarray}\label{commu-1-0-1}
&{\displaystyle \sum_{m}\sum_{k=1}^{N}
(\pi_{n}(w_{(1)}\boxtimes_{P(z_{1}-z_{2})}w_{(2)}))^{\Y_{1}}_{m;k}
w_{(3)}e^{(-m-1)\log z_{2}}(\log z_{2})^{k}}&\nn
&{\displaystyle =(\pi_{n}(w_{(1)}\boxtimes_{P(z_{1}-z_{2})}w_{(2)}))
\boxtimes_{P(z_{2})}w_{(3)}}&.
\end{eqnarray}
But by the definition of 
$\mathcal{R}_{P(z_{1}-z_{2})}\boxtimes_{P(z_{2})}
1_{W_{3}}$,
\begin{eqnarray}\label{commu-1-0-2}
\lefteqn{\overline{(\mathcal{R}_{P(z_{1}-z_{2})}
\boxtimes_{P(z_{2})}1_{W_{3}})}
((\pi_{n}(w_{(1)}\boxtimes_{P(z_{1}-z_{2})}w_{(2)}))
\boxtimes_{P(z_{2})}w_{(3)})}\nn
&&\quad =\mathcal{R}_{P(z_{1}-z_{2})}
(\pi_{n}(w_{(1)}\boxtimes_{P(z_{1}-z_{2})}w_{(2)}))
\boxtimes_{P(z_{2})}w_{(3)}\nn
&&\quad =\pi_{n}(\mathcal{R}_{P(z_{1}-z_{2})}
(w_{(1)}\boxtimes_{P(z_{1}-z_{2})}w_{(2)}))
\boxtimes_{P(z_{2})}w_{(3)}\quad\quad
\end{eqnarray}
for $n\in \C$. {}From (\ref{commu-1-0-1}), (\ref{commu-1-0-2}) and 
the definitions of $\Y_{1}$ and $\tilde{\Y}_{1}$,
we obtain
\begin{eqnarray}\label{commu-1-0-3}
\lefteqn{\overline{(\mathcal{R}_{P(z_{1}-z_{2})}
\boxtimes_{P(z_{2})}1_{W_{3}})}
((\pi_{n}(w_{(1)}\boxtimes_{P(z_{1}-z_{2})}w_{(2)}))^{\Y_{1}}_{m;k}
w_{(3)})}\nn
&&=(\pi_{n}(\mathcal{R}_{P(z_{1}-z_{2})}
(w_{(1)}\boxtimes_{P(z_{1}-z_{2})}w_{(2)})))^{\tilde{\Y}_{1}}_{m;k}
w_{(3)}\nn
&&=(\pi_{n}(e^{(z_{1}-z_{2})L(-1)}
(w_{(2)}\boxtimes_{P(z_{2}-z_{1})}w_{(1)})))^{\tilde{\Y}_{1}}_{m;k}
w_{(3)}
\end{eqnarray}
for $m\in \C$ and $k=1, \dots, N$.
Using (\ref{commu-1-0-0}) and (\ref{commu-1-0-3}),
we obtain
\begin{eqnarray}\label{commu-1-1}
\lefteqn{\overline{(\mathcal{R}_{P(z_{1}-z_{2})}\boxtimes_{P(z_{2})}1_{W_{3}})}
(\Y_{1}(w_{(1)}\boxtimes_{P(z_{1}-z_{2})}w_{(2)}, z_{2}+\epsilon)w_{(3)})}\nn
&&=\overline{(\mathcal{R}_{P(z_{1}-z_{2})}\boxtimes_{P(z_{2})}1_{W_{3}})}
\nn
&&\quad\quad\left(\sum_{m, n\in \C}\sum_{k=1}^{N}
(\pi_{n}(w_{(1)}\boxtimes_{P(z_{1}-z_{2})}w_{(2)}))^{\Y_{1}}_{m;k}
w_{(3)}
e^{(-m-1)\log (z_{2}+\epsilon)}(\log (z_{2}+\epsilon))^{k}\right)\nn
&&=\sum_{m, n\in \C}\sum_{k=1}^{N}
\overline{(\mathcal{R}_{P(z_{1}-z_{2})}\boxtimes_{P(z_{2})}1_{W_{3}})}
((\pi_{n}(w_{(1)}\boxtimes_{P(z_{1}-z_{2})}w_{(2)}))^{\Y_{1}}_{m;k}
w_{(3)})\cdot\nn
&&\quad\quad\quad\quad\quad\quad\cdot
e^{(-m-1)\log (z_{2}+\epsilon)}(\log (z_{2}+\epsilon))^{k}\nn
&&=\sum_{m, n\in \C}\sum_{k=1}^{N}
(\pi_{n}(e^{(z_{1}-z_{2})L(-1)}
(w_{(2)}\boxtimes_{P(z_{2}-z_{1})}w_{(1)})))^{\tilde{\Y}_{1}}_{m;k}
w_{(3)}\cdot\nn
&&\quad\quad\quad\quad\quad\quad\cdot
e^{(-m-1)\log (z_{2}+\epsilon)}(\log (z_{2}+\epsilon))^{k}\nn
&&=\tilde{\Y}_{1}(e^{(z_{1}-z_{2})L(-1)}(w_{(2)}\boxtimes_{P(z_{2}-z_{1})}
w_{(1)}), z_{2}+\epsilon)w_{(3)}.
\end{eqnarray}

Let $\Y_{3}=\Y_{\boxtimes_{P(z_{1})}, 0}$
be the intertwining operator of type
\[
{(W_{2}\boxtimes_{P(z_{2}-z_{1})}W_{1})\boxtimes_{P(z_{1})}W_{3}\choose 
(W_{2}\boxtimes_{P(z_{2}-z_{1})}W_{1})\;\;\;\;W_{3}}.
\]
Then by the definition of the parallel transport isomorphism,
for $m\in \C$
we have 
\begin{eqnarray}\label{commu-1-2}
\lefteqn{\overline{\mathcal{T}_{\gamma}}(
\tilde{\Y}_{1}(\pi_{m}(e^{(z_{1}-z_{2})L(-1)}(w_{(2)}\boxtimes_{P(z_{2}-z_{1})}
w_{(1)})), z_{2}+\epsilon)w_{(3)})}\nn
&&=\overline{\mathcal{T}_{\gamma}}(
\tilde{\Y}_{1}(e^{\epsilon L(-1)}
\pi_{m}(e^{(z_{1}-z_{2})L(-1)}(w_{(2)}\boxtimes_{P(z_{2}-z_{1})}
w_{(1)})), z_{2})w_{(3)})\nn
&&=\overline{\mathcal{T}_{\gamma}}(
(e^{\epsilon L(-1)}
\pi_{m}(e^{(z_{1}-z_{2})L(-1)}(w_{(2)}\boxtimes_{P(z_{2}-z_{1})}
w_{(1)})))\boxtimes_{P(z_{2})}w_{(3)})\nn
&&=\Y_{3}(e^{\epsilon L(-1)}
\pi_{m}(e^{(z_{1}-z_{2})L(-1)}(w_{(2)}\boxtimes_{P(z_{2}-z_{1})}
w_{(1)})), z_{2})w_{(3)}\nn
&&=\Y_{3}(
\pi_{m}(e^{(z_{1}-z_{2})L(-1)}(w_{(2)}\boxtimes_{P(z_{2}-z_{1})}
w_{(1)})), z_{2}+\epsilon)w_{(3)}.
\end{eqnarray}
Since $|z_{2}+\epsilon|>|z_{1}-z_{2}|>0$, the sums of both sides
of (\ref{commu-1-2}) for $m\in \C$ are absolutely convergent
(in the 
sense that the series obtained by paring it with elements of 
$((W_{2}\boxtimes_{P(z_{2}-z_{1})}W_{1})\boxtimes_{P(z_{1})}W_{3})'$
are absolutely convergent) and 
we have
\begin{eqnarray}\label{commu-1-3}
\lefteqn{\sum_{m\in \C}\overline{\mathcal{T}_{\gamma}}(
\tilde{\Y}_{1}(\pi_{m}(e^{(z_{1}-z_{2})L(-1)}(w_{(2)}\boxtimes_{P(z_{2}-z_{1})}
w_{(1)})), z_{2}+\epsilon)w_{(3)})}\nn
&&=\overline{\mathcal{T}_{\gamma}}(
\tilde{\Y}_{1}(e^{(z_{1}-z_{2})L(-1)}(w_{(2)}\boxtimes_{P(z_{2}-z_{1})}
w_{(1)}), z_{2}+\epsilon)w_{(3)})
\end{eqnarray}
and
\begin{eqnarray}\label{commu-1-4}
\lefteqn{\sum_{m\in \C}\Y_{3}(
\pi_{m}(e^{(z_{1}-z_{2})L(-1)}(w_{(2)}\boxtimes_{P(z_{2}-z_{1})}
w_{(1)})), z_{2}+\epsilon)w_{(3)}}\nn
&&=\Y_{3}((w_{(2)}\boxtimes_{P(z_{2}-z_{1})}
w_{(1)}), z_{2}+\epsilon+(z_{1}-z_{2}))w_{(3)}\nn
&&=\Y_{3}((w_{(2)}\boxtimes_{P(z_{2}-z_{1})}
w_{(1)}), z_{1}+\epsilon)w_{(3)}.
\end{eqnarray}
{}From (\ref{commu-1-2})--(\ref{commu-1-4}), we obtain
\begin{eqnarray}\label{commu-1-5}
\lefteqn{\overline{\mathcal{T}_{\gamma}}(
\tilde{\Y}_{1}(e^{(z_{1}-z_{2})L(-1)}(w_{(2)}\boxtimes_{P(z_{2}-z_{1})}
w_{(1)}), z_{2}+\epsilon)w_{(3)})}\nn
&&=\Y_{3}((w_{(2)}\boxtimes_{P(z_{2}-z_{1})}
w_{(1)}), z_{1}+\epsilon)w_{(3)}.
\end{eqnarray}
{}From (\ref{commu-1-1}) and (\ref{commu-1-5}), we obtain
\begin{eqnarray*}
\lefteqn{\overline{\mathcal{T}_{\gamma}\circ
(\mathcal{R}_{P(z_{1}-z_{2})}\boxtimes_{P(z_{2})}1_{W_{3}})}
(\Y_{1}(w_{(1)}\boxtimes_{P(z_{1}-z_{2})}w_{(2)}, z_{2}+\epsilon)w_{(3)})}\nn
&&=\overline{\mathcal{T}_{\gamma}}(
\tilde{\Y}_{1}(e^{(z_{1}-z_{2})L(-1)}(w_{(2)}\boxtimes_{P(z_{2}-z_{1})}
w_{(1)}), z_{2}+\epsilon)w_{(3)})\nn
&&=\Y_{3}((w_{(2)}\boxtimes_{P(z_{2}-z_{1})}
w_{(1)}), z_{1}+\epsilon)w_{(3)}.
\end{eqnarray*}
Then for any 
\[
w'\in ((W_{1}\boxtimes_{P(z_{2}-z_{1})}W_{2})\boxtimes_{P(z_{1})}W_{3})',
\]
we have
\begin{eqnarray*}
\lefteqn{\langle w', \overline{\mathcal{T}_{\gamma}\circ
(\mathcal{R}_{P(z_{1}-z_{2})}\boxtimes_{P(z_{2})}1_{W_{3}})}
(\Y_{1}(w_{(1)}\boxtimes_{P(z_{1}-z_{2})}w_{(2)}, z_{2}+\epsilon)w_{(3)})
\rangle}\nn
&&\quad\quad=\langle w', \Y_{3}((w_{(2)}\boxtimes_{P(z_{2}-z_{1})}
w_{(1)}), z_{1}+\epsilon)w_{(3)}\rangle,\quad\quad\quad\quad\quad\quad
\end{eqnarray*}
or equivalently,
\begin{eqnarray}\label{commu-1-7}
\lefteqn{\langle ((\mathcal{R}_{P(z_{1}-z_{2})}\boxtimes_{P(z_{2})}1_{W_{3}})'
\circ \mathcal{T}_{\gamma}')(w'), 
\Y_{1}(w_{(1)}\boxtimes_{P(z_{1}-z_{2})}w_{(2)}, z_{2}+\epsilon)w_{(3)}
\rangle}\nn
&&\quad\quad=\langle w', \Y_{3}((w_{(2)}\boxtimes_{P(z_{2}-z_{1})}
w_{(1)}), z_{1}+\epsilon)w_{(3)}\rangle,\quad\quad\quad\quad\quad\quad
\end{eqnarray}
The left- and right-hand sides of (\ref{commu-1-7})
are values at  and $(\zeta_{1}, \zeta_{2})
=(z_{1}+\epsilon, z_{2}+\epsilon)$ of some single-valued
analytic functions of $\zeta_{1}$ and $\zeta_{2}$
defined in the region 
\[
\{(\zeta_{1}, \zeta_{2})\in \C^{2}\;|\;\zeta_{1}\ne 0,\;
\zeta_{2}\ne 0,\;\zeta_{1}\ne \zeta_{2},\;
0\le \arg \zeta_{1}, \arg \zeta_{2}, 
\arg (\zeta_{1}-\zeta_{2})<2\pi\}.
\]
Also, by the definition of tensor product of three elements
above, the values of these analytic functions at $(\zeta_{1}, \zeta_{2})
=(z_{1}, z_{2})$  are equal to 
\[
\langle ((\mathcal{R}_{P(z_{1}-z_{2})}\boxtimes_{P(z_{2})}1_{W_{3}})'
\circ \mathcal{T}_{\gamma}')(w'), 
(w_{(1)}\boxtimes_{P(z_{1}-z_{2})}w_{(2)})\boxtimes_{P(z_{2})}w_{(3)}
\rangle
\]
and 
\[
\langle w', (w_{(2)}\boxtimes_{P(z_{2}-z_{1})}
w_{(1)})\boxtimes_{P(z_{1})}w_{(3)}\rangle,
\]
respectively.
Thus we can take the limit $\epsilon\to 0$ on both sides of 
(\ref{commu-1-7}) and obtain
\begin{eqnarray}\label{commu-1-8}
\lefteqn{\langle ((\mathcal{R}_{P(z_{1}-z_{2})}\boxtimes_{P(z_{2})}1_{W_{3}})'
\circ \mathcal{T}_{\gamma}')(w'), 
(w_{(1)}\boxtimes_{P(z_{1}-z_{2})}w_{(2)})\boxtimes_{P(z_{2})}w_{(3)}
\rangle}\nn
&&\quad\quad=\langle w', (w_{(2)}\boxtimes_{P(z_{2}-z_{1})}
w_{(1)})\boxtimes_{P(z_{1})}w_{(3)}\rangle.\quad\quad\quad\quad\quad\quad
\end{eqnarray}
Since $w'$ is arbitrary, (\ref{commu-1-8}) is equivalent to 
(\ref{commu-1}).
\epfv

We also prove:

\begin{propo}
Let $z_{1}, z_{2}$ be nonzero complex numbers such that 
$z_{1}\ne z_{2}$ but $|z_{1}|=|z_{2}|=|z_{1}-z_{2}|$. Let $\gamma$ be 
a path {}from $z_{2}$ to $z_{2}-z_{1}$ in the complex plane with a cut 
along the nonnegative real line. Then
we have 
\begin{eqnarray}\label{commu-2}
\lefteqn{\overline{\mathcal{T}_{\gamma}\circ 
(1_{W_{3}}\boxtimes_{P(z_{2})}\mathcal{R}_{P(z_{1})})}
(w_{(2)}\boxtimes_{P(z_{2})}(w_{(3)}\boxtimes_{P(z_{1})}w_{(1)}))}\nn
&&=e^{z_{1}L(-1)}(w_{(2)}\boxtimes_{P(z_{2}-z_{1})}(w_{(3)}
\boxtimes_{P(-z_{1})}w_{(1)}))
\end{eqnarray}
for $w_{(1)}\in W_{1}$, $w_{(2)}\in W_{2}$ and $w_{(3)}\in W_{3}$.
\end{propo}
\pf
We can find $\epsilon$ such that $|z_{2}|>|\epsilon|$,
\[
|z_{2}+\epsilon|, |z_{2}-z_{1}+\epsilon|>|z_{1}|>0.
\]
Let 
$\Y_{1}=\Y_{\boxtimes_{P(z_{2})}, 0}$, $\tilde{\Y}_{1}
=\Y_{\boxtimes_{P(z_{2})}, 0}$ and 
$\Y_{2}=\Y_{\boxtimes_{P(z_{1})}, 0}$
be intertwining operators of types 
\[
{W_{2}\boxtimes_{P(z_{2})}(W_{1}\boxtimes_{P(z_{1})}W_{3})\choose 
W_{2}\;\;\;\;(W_{1}\boxtimes_{P(z_{1})}W_{3})},
\]
\[
{W_{2}\boxtimes_{P(z_{2})}(W_{3}\boxtimes_{P(-z_{1})}W_{1})\choose 
W_{2}\;\;\;\;(W_{3}\boxtimes_{P(-z_{1})}W_{1})}
\]
and 
\[
{W_{1}\boxtimes_{P(z_{1})}W_{3}\choose W_{1}\;\;\;\; W_{3}},
\]
respectively, corresponding to the intertwining maps 
$\boxtimes_{P(z_{2})}$, $\boxtimes_{P(z_{2})}$
and $\boxtimes_{P(z_{1})}$, respectively. 
Since $|z_{2}+\epsilon|>|z_{1}|>0$,
\begin{eqnarray}\label{commu-2-0-0}
\lefteqn{\Y_{1}(w_{(2)}, 
z_{2}+\epsilon)(w_{(1)}\boxtimes_{P(z_{1})}w_{(3)})}\nn
&&=\sum_{m, n\in \C}\sum_{k=1}^{N}
(w_{(2)})^{\Y_{1}}_{m;k}\pi_{n}(w_{(1)}\boxtimes_{P(z_{1})}
w_{(3)})e^{(-m-1)\log (z_{2}+\epsilon)}(\log (z_{2}+\epsilon))^{k}.
\end{eqnarray}
For $n\in \C$, 
\begin{eqnarray}\label{commu-2-0-1}
&{\displaystyle \sum_{m}\sum_{k=1}^{N}
(w_{(2)})^{\Y_{1}}_{m;k}\pi_{n}(w_{(1)}\boxtimes_{P(z_{1})}
w_{(3)})e^{(-m-1)\log z_{2}}(\log z_{2})^{k}}&\nn
&{\displaystyle =w_{(2)}\boxtimes_{P(z_{2})}\pi_{n}(w_{(1)}
\boxtimes_{P(z_{1})}w_{(3)})}&.
\end{eqnarray}
But by the definition of 
$1_{W_{3}}\boxtimes_{P(z_{2})}\mathcal{R}_{P(z_{1})}$,
\begin{eqnarray}\label{commu-2-0-2}
\lefteqn{\overline{(1_{W_{3}}\boxtimes_{P(z_{2})}\mathcal{R}_{P(z_{1})})}
(w_{(2)}\boxtimes_{P(z_{2})}\pi_{n}(w_{(2)}\boxtimes_{P(z_{1})}w_{(3)}))}\nn
&&\quad =w_{(2)}\boxtimes_{P(z_{2})}
\mathcal{R}_{P(z_{1})}
(\pi_{n}(w_{(1)}\boxtimes_{P(z_{1})}w_{(3)}))\nn
&&\quad =w_{(2)}\boxtimes_{P(z_{2})}
\pi_{n}(\mathcal{R}_{P(z_{1})}
(w_{(1)}\boxtimes_{P(z_{1})}w_{(3)}))\quad\quad
\end{eqnarray}
for $n\in \C$. {}From (\ref{commu-2-0-1}), (\ref{commu-2-0-2}) and 
the definitions of $\Y_{1}$ and $\tilde{\Y}_{1}$,
we obtain
\begin{eqnarray}\label{commu-2-0-3}
\lefteqn{\overline{(1_{W_{3}}\boxtimes_{P(z_{2})}\mathcal{R}_{P(z_{1})})}
((w_{(2)})^{\Y_{1}}_{m;k}\pi_{n}(w_{(1)}\boxtimes_{P(z_{1})}
w_{(3)}))}\nn
&&=(w_{(2)})^{\tilde{\Y}_{1}}_{m;k}
\pi_{n}(\mathcal{R}_{P(z_{1})}(w_{(1)}\boxtimes_{P(z_{1})}
w_{(3)}))\nn
&&=(w_{(2)})^{\tilde{\Y}_{1}}_{m;k}
\pi_{n}(e^{z_{1}L(-1)}(w_{(3)}\boxtimes_{P(-z_{1})}
w_{(1)}))
\end{eqnarray}
for $m\in \C$ and $k=1, \dots, N$. 
Using (\ref{commu-2-0-0}) and (\ref{commu-2-0-3}),
we obtain
\begin{eqnarray}\label{commu-2-1}
\lefteqn{\overline{(1_{W_{3}}\boxtimes_{P(z_{2})}\mathcal{R}_{P(z_{1})})}
(\Y_{1}(w_{(2)}, 
z_{2}+\epsilon)(w_{(1)}\boxtimes_{P(z_{1})}w_{(3)}))}\nn
&&=\overline{(1_{W_{3}}\boxtimes_{P(z_{2})}\mathcal{R}_{P(z_{1})})}
\nn
&&\quad\quad\left(\sum_{m, n\in \C}\sum_{k=1}^{N}
(w_{(2)})^{\Y_{1}}_{m;k}\pi_{n}(w_{(1)}\boxtimes_{P(z_{1})}
w_{(3)})e^{(-m-1)\log (z_{2}+\epsilon)}(\log (z_{2}+\epsilon))^{k}\right)\nn
&&=\sum_{m, n\in \C}\sum_{k=1}^{N}
\overline{(1_{W_{3}}\boxtimes_{P(z_{2})}\mathcal{R}_{P(z_{1})})}
((w_{(2)})^{\Y_{1}}_{m;k}\pi_{n}(w_{(1)}\boxtimes_{P(z_{1})}
w_{(3)}))\cdot\nn
&&\quad\quad\quad\quad\quad\quad\cdot
e^{(-m-1)\log (z_{2}+\epsilon)}(\log (z_{2}+\epsilon))^{k}\nn
&&=\sum_{m, n\in \C}\sum_{k=1}^{N}
(w_{(2)})^{\tilde{\Y}_{1}}_{m;k}
\pi_{n}(e^{z_{1}L(-1)}(w_{(3)}\boxtimes_{P(-z_{1})}
w_{(1)}))
e^{(-m-1)\log (z_{2}+\epsilon)}(\log (z_{2}+\epsilon))^{k}\nn
&&=\tilde{\Y}_{1}(w_{(2)}, z_{2}+\epsilon)(e^{z_{1}L(-1)}
(w_{(3)}\boxtimes_{P(-z_{1})}
w_{(1)})).
\end{eqnarray}

Let $\Y_{3}=\Y_{\boxtimes_{P(z_{1})}, 0}$
be the intertwining operator of type
\[
{W_{2}\boxtimes_{P(z_{2}-z_{1})}(W_{1}\boxtimes_{P(-z_{1})}W_{3})\choose 
W_{2}\;\;\;\;W_{1}\boxtimes_{-z_{1})}W_{3}}.
\]
Then by the definition of the parallel transport isomorphism,
for $m\in \C$
we have 
\begin{eqnarray}\label{commu-2-2}
\lefteqn{\overline{\mathcal{T}_{\gamma}}(
\tilde{\Y}_{1}(w_{(2)}, z_{2}+\epsilon)\pi_{m}(e^{z_{1}L(-1)}(w_{(3)}
\boxtimes_{P(-z_{1})}w_{(1)})))}\nn
&&=\overline{\mathcal{T}_{\gamma}}(
e^{\epsilon L(-1)}\tilde{\Y}_{1}(w_{(2)}, z_{2})e^{\epsilon L(-1)}
\pi_{m}(e^{z_{1}L(-1)}(w_{(3)}
\boxtimes_{P(-z_{1})}w_{(1)})))\nn
&&=\overline{\mathcal{T}_{\gamma}}(
e^{\epsilon L(-1)}w_{(2)}\boxtimes_{P(z_{2})}e^{\epsilon L(-1)}
\pi_{m}(e^{z_{1}L(-1)}(w_{(3)}
\boxtimes_{P(-z_{1})}w_{(1)})))\nn
&&=e^{\epsilon L(-1)}\Y_{3}(w_{(2)}, z_{2})e^{\epsilon L(-1)}
\pi_{m}(e^{z_{1} L(-1)}(w_{(3)}
\boxtimes_{P(-z_{1})}w_{(1)}))\nn
&&=\Y_{3}(w_{(2)}, z_{2}+\epsilon)
\pi_{m}(e^{z_{1} L(-1)}(w_{(3)}
\boxtimes_{P(-z_{1})}w_{(1)})).
\end{eqnarray}
Since $|z_{2}+\epsilon|>|z_{1}|>0$, the sums of both
sides of (\ref{commu-2-2}) over $m\in \C$ are absolutely 
convergent. Note also that $|z_{2}-z_{1}+\epsilon|>|z_{1}|>0$.
Thus we have
\begin{eqnarray}\label{commu-2-3}
\lefteqn{\sum_{m\in \C}\overline{\mathcal{T}_{\gamma}}(
\Y_{1}(w_{(2)}, z_{2}+\epsilon)\pi_{m}(e^{z_{1}L(-1)}(w_{(3)}
\boxtimes_{P(-z_{1})}w_{(1)})))}\nn
&&=\overline{\mathcal{T}_{\gamma}}(
\Y_{1}(w_{(2)}, z_{2}+\epsilon)e^{z_{1}L(-1)}(w_{(3)}
\boxtimes_{P(-z_{1})}w_{(1)}))\nn
&&=\overline{\mathcal{T}_{\gamma}}(e^{z_{1}L(-1)}
\Y_{1}(w_{(2)}, z_{2}-z_{1}+\epsilon)(w_{(3)}
\boxtimes_{P(-z_{1})}w_{(1)}))
\end{eqnarray}
and
\begin{eqnarray}\label{commu-2-4}
\lefteqn{\sum_{m\in \C}\Y_{3}(w_{(2)}, z_{2}+\epsilon)
\pi_{m}(e^{z_{1} L(-1)}(w_{(3)}
\boxtimes_{P(-z_{1})}w_{(1)}))}\nn
&&=\Y_{3}(w_{(2)}, z_{2}+\epsilon)
e^{z_{1} L(-1)}(w_{(3)}
\boxtimes_{P(-z_{1})}w_{(1)})\nn
&&=e^{z_{1} L(-1)}\Y_{3}(w_{(2)}, z_{2}-z_{1}+\epsilon)
(w_{(3)}
\boxtimes_{P(-z_{1})}w_{(1)}).
\end{eqnarray}
{}From (\ref{commu-2-2})--(\ref{commu-2-4}), we obtain
\begin{eqnarray}\label{commu-2-5}
\lefteqn{\overline{\mathcal{T}_{\gamma}}(e^{z_{1}L(-1)}
\Y_{1}(w_{(2)}, z_{2}-z_{1}+\epsilon)(w_{(3)}
\boxtimes_{P(-z_{1})}w_{(1)}))}\nn
&&=e^{z_{1} L(-1)}\Y_{3}(w_{(2)}, z_{2}-z_{1}+\epsilon)
(w_{(3)}
\boxtimes_{P(-z_{1})}w_{(1)}).
\end{eqnarray}
{}From (\ref{commu-2-1}) and (\ref{commu-2-5}), we obtain
\begin{eqnarray}\label{commu-2-6}
\lefteqn{\overline{\mathcal{T}_{\gamma}\circ 
(1_{W_{3}}\boxtimes_{P(z_{2})}\mathcal{R}_{P(z_{1})})}
(\Y_{1}(w_{(2)}, z_{2}+\epsilon)(w_{(1)}
\boxtimes_{P(z_{1})}w_{(3)}))}\nn
&&=e^{z_{1} L(-1)}\Y_{3}(w_{(2)}, z_{2}-z_{1}+\epsilon)
(w_{(3)}
\boxtimes_{P(-z_{1})}w_{(1)}).
\end{eqnarray}

For any 
\[
w'\in (W_{2}\boxtimes_{P(z_{2}-z_{1})}(W_{1}\boxtimes_{P(-z_{1})}W_{3}))',
\]
we have
\begin{eqnarray*}
\lefteqn{\langle w', \overline{\mathcal{T}_{\gamma}\circ 
(1_{W_{3}}\boxtimes_{P(z_{2})}\mathcal{R}_{P(z_{1})})}
(\Y_{1}(w_{(2)}, z_{2}+\epsilon)(w_{(1)}
\boxtimes_{P(z_{1})}w_{(3)}))\rangle}\nn
&&=\langle w', e^{z_{1} L(-1)}\Y_{3}(w_{(2)}, z_{2}-z_{1}+\epsilon)
(w_{(3)}
\boxtimes_{P(-z_{1})}w_{(1)})\rangle,
\end{eqnarray*}
or equivalently,
\begin{eqnarray}\label{commu-2-7}
\lefteqn{\langle ((1_{W_{3}}\boxtimes_{P(z_{2})}\mathcal{R}_{P(z_{1})})
\circ \mathcal{T}_{\gamma})(w)', \Y_{1}(w_{(2)}, z_{2}+\epsilon)(w_{(1)}
\boxtimes_{P(z_{1})}w_{(3)})\rangle}\nn
&&=\langle w', e^{z_{1} L(-1)}\Y_{3}(w_{(2)}, z_{2}-z_{1}+\epsilon)
(w_{(3)}
\boxtimes_{P(-z_{1})}w_{(1)})\rangle,
\end{eqnarray}
The left- and right-hand sides of (\ref{commu-2-7})
are values at $(\zeta_{1}, \zeta_{2})
=(z_{1}, z_{2}+\epsilon)$  of  single-valued
analytic functions of $\zeta_{1}$ and $\zeta_{2}$
defined in the region 
\[
\{(\zeta_{1}, \zeta_{2})\in \C^{2}\;|\;\zeta_{1}\ne 0,\;
\zeta_{2}\ne 0,\;\zeta_{1}\ne \zeta_{2},\;
0\le \arg \zeta_{1}, \arg \zeta_{2}, 
\arg (\zeta_{1}-\zeta_{2})<2\pi\}.
\]
Also, by the definition of tensor product of three elements
above, the values of these analytic functions at $(\zeta_{1}, \zeta_{2})
=(z_{1}, z_{2})$ are equal to 
\[
\langle ((1_{W_{3}}\boxtimes_{P(z_{2})}\mathcal{R}_{P(z_{1})})
\circ \mathcal{T}_{\gamma})(w)', w_{(2)}\boxtimes_{P(z_{2})}(w_{(1)}
\boxtimes_{P(z_{1})}w_{(3)})\rangle
\]
and 
\[
\langle w', e^{z_{1} L(-1)}w_{(2)}\boxtimes_{P(z_{2}-z_{1})}
(w_{(3)}
\boxtimes_{P(-z_{1})}w_{(1)})\rangle,
\]
respectively.
Thus we can take the limit $\epsilon\to 0$ on both sides of 
(\ref{commu-2-7}) and obtain
\begin{eqnarray}\label{commu-2-8}
\lefteqn{\langle ((1_{W_{3}}\boxtimes_{P(z_{2})}\mathcal{R}_{P(z_{1})})
\circ \mathcal{T}_{\gamma})(w)', w_{(2)}\boxtimes_{P(z_{2})}(w_{(1)}
\boxtimes_{P(z_{1})}w_{(3)})\rangle}\nn
&&\quad\quad\quad=\langle w', e^{z_{1} L(-1)}w_{(2)}\boxtimes_{P(z_{2}-z_{1})}
(w_{(3)}\boxtimes_{P(-z_{1})}w_{(1)})\rangle.
\end{eqnarray}
Since $w'$ is arbitrary, (\ref{commu-2-8}) is equivalent to 
(\ref{commu-2}).
\epfv

\subsection{The coherence properties}

\begin{theo}
The category $\mathcal{C}$, equipped with the tensor product bifunctor
$\boxtimes$, the unit object $V$, the braiding isomorphisms
$\mathcal{R}$, the associativity isomorphisms $\mathcal{A}$, and the
left and right unit isomorphisms $l$ and $r$, is a braided tensor
category.
\end{theo}
\pf We need only prove the coherence properties. We prove the
commutativity of the pentagon diagram first.  Let $W_{1}$, $W_{2}$,
$W_{3}$ and $W_{4}$ be objects of $\mathcal{C}$ and let $z_{1}, z_{2},
z_{3}\in \R$ such that
\begin{eqnarray}
&|z_{1}|>|z_{2}|>|z_{3}|>
|z_{13}|>|z_{23}|>|z_{12}|>0,&\nn
&|z_{1}|>|z_{3}|+|z_{23}|>0,&\nn
&|z_{2}|>|z_{12}|+|z_{3}|>0,&\nn
&|z_{3}|>|z_{23}|+|z_{12}|>0,&
\end{eqnarray}
where $z_{12}=z_{1}-z_{2}$, $z_{13}=z_{1}-z_{3}$ and $z_{23}=z_{2}-z_{3}$.
For example, we can take $z_{1}=7$, $z_{2}=6$ and $z_{3}=4$.
We first prove the commutativity of the following diagram:

\vspace{3.5em}

\begin{picture}(150,100)(-70,0)

\put(-85,20){\footnotesize $((W_{1}\boxtimes_{P(z_{12})}
W_{2})\boxtimes_{P(z_{23})} W_{3})
\boxtimes_{P(z_{3})} W_{4}$}
\put(123,20){\footnotesize $(W_{1}\boxtimes_{P(z_{13})}
(W_{2}\boxtimes_{P(z_{23})} W_{3}))
\boxtimes_{P(z_{3})} W_{4}$}
\put(-85,68){\footnotesize $(W_{1}\boxtimes_{P(z_{12})}
W_{2})\boxtimes_{P(z_{2})} (W_{3}
\boxtimes_{P(z_{3})} W_{4})$}
\put(123,68){\footnotesize $W_{1}\boxtimes_{P(z_{1})}
((W_{2}\boxtimes_{P(z_{23})} W_{3})
\boxtimes_{P(z_{3})} W_{4}))$}
\put(20,116){\footnotesize $W_{1}\boxtimes_{P(z_{1})}
(W_{2}\boxtimes_{P(z_{2})} (W_{3}
\boxtimes_{P(z_{3})} W_{4})).$}

\put(120,23){\vector(-1,0){25}}
\put(5,60){\vector(0,-1){28}}
\put(212,60){\vector(0,-1){28}}

\put(100,105){\vector(-3,-1){75}}
\put(110,105){\vector(3,-1){75}}
\end{picture}
\begin{equation}\label{pent1}
\end{equation}
For $w_{(1)}\in W_{1}$, $w_{(2)}\in W_{2}$, $w_{(3)}\in W_{3}$ and
$w_{(4)}\in W_{4}$, we consider
$$w_{(1)}\boxtimes_{P(z_{1})} (w_{(2)}\boxtimes_{P(z_{2})}
(w_{(3)}\boxtimes_{P(z_{3})}w_{(4)}))\in
\overline{W_{1}\boxtimes_{P(z_{1})} (W_{2}\boxtimes_{P(z_{2})}
(W_{3}\boxtimes_{P(z_{3})}W_{4}))}.$$
By the characterizations of the associativity isomorphisms,
we see that the compositions of the natural extensions
of the module maps in the two routes in (\ref{pent1}) applied to
this element both give
$$((w_{(1)}\boxtimes_{P(z_{12})} w_{(2)})\boxtimes_{P(z_{23})}
w_{(3)})\boxtimes_{P(z_{3})}w_{(4)}\in
\overline{((W_{1}\boxtimes_{P(z_{12})} W_{2})\boxtimes_{P(z_{23})}
W_{3})\boxtimes_{P(z_{3})}W_{4}}.$$
Since the homogeneous components of
$$w_{(1)}\boxtimes_{P(z_{1})} (w_{(2)}\boxtimes_{P(z_{2})}
(w_{(3)}\boxtimes_{P(z_{3})}w_{(4)}))$$
for $w_{(1)}\in W_{1}$, $w_{(2)}\in W_{2}$, $w_{(3)}\in W_{3}$ and
$w_{(4)}\in W_{4}$ span
$$W_{1}\boxtimes_{P(z_{1})} (W_{2}\boxtimes_{P(z_{2})}
(W_{3}\boxtimes_{P(z_{3})}W_{4})),$$
the diagram (\ref{pent1}) above is commutative.

On the other hand, by the definition of $\A$, the  diagrams
\begin{equation}\label{pent2}
\begin{picture}(60,90)(20,0)
\put(-145,68){\footnotesize $W_{1}\boxtimes_{P(z_{1})}
(W_{2}\boxtimes_{P(z_{2})} (W_{3}
\boxtimes_{P(z_{3})} W_{4}))$}
\put(63,68){\footnotesize $(W_{1}\boxtimes_{P(z_{12})}
W_{2})\boxtimes_{P(z_{2})} (W_{3}
\boxtimes_{P(z_{3})} W_{4})$}
\put(-105,20){\footnotesize $W_{1}\boxtimes
(W_{2}\boxtimes (W_{3}
\boxtimes W_{4}))$}
\put(88,20){\footnotesize $(W_{1}\boxtimes
W_{2})\boxtimes (W_{3}
\boxtimes W_{4})$}

\put(5,23){\vector(1,0){78}}
\put(30,71){\vector(1,0){25}}
\put(-55,60){\vector(0,-1){28}}
\put(152,60){\vector(0,-1){28}}
\end{picture}
\end{equation}
\begin{equation}\label{pent3}
\begin{picture}(60,90)(20,0)
\put(-145,68){\footnotesize $(W_{1}\boxtimes_{P(z_{12})}
W_{2})\boxtimes_{P(z_{2})} (W_{3}
\boxtimes_{P(z_{3})} W_{4})$}
\put(63,68){\footnotesize $((W_{1}\boxtimes_{P(z_{12})}
W_{2})\boxtimes_{P(z_{23})} W_{3})
\boxtimes_{P(z_{3})} W_{4}$}
\put(-105,20){\footnotesize $(W_{1}\boxtimes
W_{2})\boxtimes (W_{3}
\boxtimes W_{4})$}
\put(88,20){\footnotesize $((W_{1}\boxtimes
W_{2})\boxtimes W_{3})
\boxtimes W_{4}$}

\put(5,23){\vector(1,0){78}}
\put(34,71){\vector(1,0){25}}
\put(-55,60){\vector(0,-1){28}}
\put(152,60){\vector(0,-1){28}}
\end{picture}
\end{equation}
\begin{equation}\label{pent4}
\begin{picture}(60,90)(20,0)
\put(-145,68){\footnotesize $W_{1}\boxtimes_{P(z_{1})}
(W_{2}\boxtimes_{P(z_{2})} (W_{3}
\boxtimes_{P(z_{3})} W_{4}))$}
\put(63,68){\footnotesize $W_{1}\boxtimes_{P(z_{1})}
((W_{2}\boxtimes_{P(z_{23})} W_{3})
\boxtimes_{P(z_{3})} W_{4}))$}
\put(-105,20){\footnotesize $W_{1}\boxtimes
(W_{2}\boxtimes (W_{3}
\boxtimes W_{4}$))}
\put(88,20){\footnotesize $W_{1}\boxtimes
((W_{2}\boxtimes W_{3})
\boxtimes W_{4})$}

\put(5,23){\vector(1,0){78}}
\put(34,71){\vector(1,0){25}}
\put(-55,60){\vector(0,-1){28}}
\put(152,60){\vector(0,-1){28}}
\end{picture}
\end{equation}
\begin{equation}\label{pent5}
\begin{picture}(60,80)(20,0)
\put(-145,68){\footnotesize $W_{1}\boxtimes_{P(z_{1})}
((W_{2}\boxtimes_{P(z_{23})} W_{3})
\boxtimes_{P(z_{3})} W_{4}))$}
\put(63,68){\footnotesize $(W_{1}\boxtimes_{P(z_{13})}
(W_{2}\boxtimes_{P(z_{23})} W_{3}))
\boxtimes_{P(z_{3})} W_{4}$}

\put(-105,20){\footnotesize $W_{1}\boxtimes
((W_{2}\boxtimes W_{3})
\boxtimes W_{4})$}
\put(88,20){\footnotesize $(W_{1}\boxtimes
(W_{2}\boxtimes W_{3}))
\boxtimes W_{4}$}

\put(5,23){\vector(1,0){78}}
\put(34,71){\vector(1,0){25}}
\put(-55,60){\vector(0,-1){28}}
\put(152,60){\vector(0,-1){28}}
\end{picture}
\end{equation}
\begin{equation}\label{pent6}
\begin{picture}(60,90)(20,0)
\put(-145,68){\footnotesize $(W_{1}\boxtimes_{P(z_{13})}
(W_{2}\boxtimes_{P(z_{23})} W_{3}))
\boxtimes_{P(z_{3})} W_{4}$}
\put(63,68){\footnotesize $((W_{1}\boxtimes_{P(z_{12})}
W_{2})\boxtimes_{P(z_{23})} W_{3})
\boxtimes_{P(z_{3})} W_{4}$}
\put(-105,20){\footnotesize $(W_{1}\boxtimes
(W_{2}\boxtimes W_{3}))
\boxtimes W_{4}$}
\put(88,20){\footnotesize $((W_{1}\boxtimes
W_{2})\boxtimes W_{3})
\boxtimes W_{4}$}

\put(5,23){\vector(1,0){78}}
\put(34,71){\vector(1,0){25}}
\put(-55,60){\vector(0,-1){28}}
\put(152,60){\vector(0,-1){28}}
\end{picture}
\end{equation}
are all commutative. Combining all the diagrams
(\ref{pent1})--(\ref{pent6}) above, we see that
the pentagon diagram

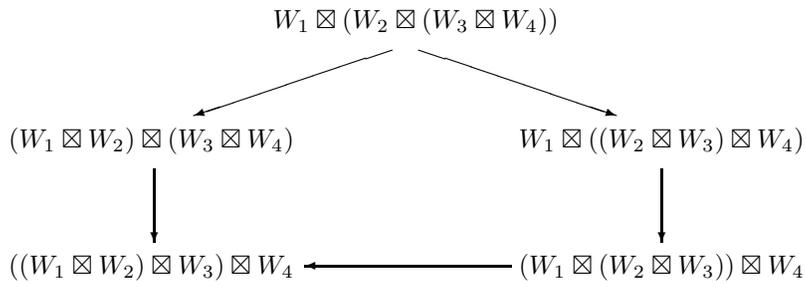
\vspace{3.5em}

\begin{picture}(150,100)(-70,0)

\put(-45,20){\footnotesize $((W_{1}\boxtimes
W_{2})\boxtimes W_{3})
\boxtimes W_{4}$}
\put(148,20){\footnotesize $(W_{1}\boxtimes
(W_{2}\boxtimes W_{3}))
\boxtimes W_{4}$}
\put(-45,68){\footnotesize $(W_{1}\boxtimes
W_{2})\boxtimes (W_{3}
\boxtimes W_{4})$}
\put(148,68){\footnotesize $W_{1}\boxtimes
((W_{2}\boxtimes W_{3})
\boxtimes W_{4})$}
\put(55,113){\footnotesize $W_{1}\boxtimes
(W_{2}\boxtimes (W_{3}
\boxtimes W_{4}))$}

\put(145,23){\vector(-1,0){78}}
\put(10,60){\vector(0,-1){28}}
\put(202,60){\vector(0,-1){28}}

\put(100,105){\vector(-3,-1){75}}
\put(110,105){\vector(3,-1){75}}
\end{picture}

\noindent is also commutative.

Next we prove the commutativity of the hexagon diagrams.  We prove
only the commutativity of the hexagon diagram involving $\mathcal{R}$;
the proof of the commutativity of the other hexagon diagram is the
same.  Let $W_{1}$, $W_{2}$ and $W_{3}$ be objects of $\mathcal{C}$
and let $z_{1}, z_{2}\in \C^{\times}$ satisfying
$|z_{1}|=|z_{2}|=|z_{1}-z_{2}|$ and let $z_{12}=z_{1}-z_{2}$.  We
first prove the commutativity of the following diagram:

\newpage

\begin{picture}(200,175)(-100,0)
\put(50,162){\footnotesize $(W_{1}\boxtimes_{P(z_{12})} W_{2})
\boxtimes_{P(z_{2})} W_{3}$}
\put(-70,66){\footnotesize $(W_{2}\boxtimes_{P(-z_{12})} W_{1})
\boxtimes_{P(z_{2})} W_{3}$}
\put(170,66){\footnotesize $W_{1}\boxtimes_{P(z_{1})} (W_{2}
\boxtimes_{P(z_{2})} W_{3})$}
\put(-70,18){\footnotesize $(W_{2}\boxtimes_{P(-z_{12})} W_{1})
\boxtimes_{P(z_{1})} W_{3}$}
\put(-70,-30){\footnotesize $W_{2}\boxtimes_{P(z_{2})} (W_{1}
\boxtimes_{P(z_{1})} W_{3})$}
\put(170,-30){\footnotesize $(W_{2}\boxtimes_{P(z_{2})} W_{3})
\boxtimes_{P(-z_{1})} W_{1}$}
\put(120,-78){\footnotesize $
W_{2}\boxtimes_{P(-z_{12})} (W_{3}
\boxtimes_{P(-z_{1})} W_{1})$}
\put(50,-126){\footnotesize $W_{2}\boxtimes_{P(z_{2})} (W_{3}
\boxtimes_{P(-z_{1})} W_{1}$)}

\put(-70,114){\footnotesize
$\mathcal{R}_{P(z_{12})}\boxtimes_{P(z_{2})} 1_{W_{3}}$}

\put(200,114){\footnotesize
$\left(\A_{P(z_{1}), P(z_{2})}
^{P(z_{12}), P(z_{2})}\right)^{-1}$}

\put(-100,-5){\footnotesize
$\left(\A_{P(z_{2}), P(z_{1})}
^{P(-z_{12}), P(z_{1})}\right)^{-1}$}

\put(220,-55){\footnotesize
$\left(\A_{P(z_{2}), P(-z_{1})}
^{P(-z_{12}), P(-z_{1})}\right)^{-1}$}

\put(-60,-78){\footnotesize
$1_{W_{2}} \boxtimes_{P(z_{2})} \mathcal{R}_{P(z_{1})}$}

\put(235,18){\footnotesize
$\mathcal{R}_{P(z_{1})}$}

\put(-35,43){\footnotesize
$\mathcal{T}_{\gamma_{1}}$}

\put(170,-105){\footnotesize
$\mathcal{T}_{\gamma_{2}}$}

\put(-10,60){\vector(0,-1){28}}
\put(230,60){\vector(0,-1){78}}

\put(-10,12){\vector(0,-1){28}}

\put(-7,-36){\vector(1,-1){78}}
\put(225,-36){\vector(-1,-1){30}}

\put(75,155){\vector(-1,-1){78}}
\put(145,155){\vector(1,-1){78}}

\put(175,-85){\vector(-1,-1){30}}

\end{picture}
\vskip 1.4in
\begin{equation}\label{hexagon1}
\end{equation}
where $\gamma_{1}$ and $\gamma_{2}$ are paths {}from $z_{2}$ to $z_{1}$
and {}from $z_{2}$ to $-z_{12}$, respectively, in $\C$ with a cut along
the nonnegative real line.

Let $w_{(1)}\in W_{1}$, $w_{(2)}\in W_{2}$ and $w_{(3)}\in W_{3}$.  By
the results proved in the preceding subsection, we see that the images
of the element
\[
(w_{(1)}\boxtimes_{P(z_{12})}
w_{(2)})\boxtimes_{P(z_{2})}w_{(3)}
\]
under the natural extension to
\[
\overline{(W_{1}\boxtimes_{P(z_{12})} W_{2})
\boxtimes_{P(z_{2})} W_{3}}
\]
of the compositions of the maps in both the left and right routes in
(\ref{hexagon1}) {}from 
\[
(W_{1}\boxtimes_{P(z_{12})} W_{2})
\boxtimes_{P(z_{2})} W_{3}
\]
to 
\[
W_{2}\boxtimes_{P(z_{2})} (W_{3}
\boxtimes_{P(-z_{1})} W_{1})
\]
are 
\[
w_{(2)}\boxtimes_{P(z_{2})}
(e^{z_{1}L(-1)}(w_{(1)}\boxtimes_{P(-z_{1})}w_{(3)})).
\]
Since the homogeneous components of 
\[
(w_{(1)}\boxtimes_{P(z_{12})}
w_{(2)})\boxtimes_{P(z_{2})}w_{(3)}
\]
for $w_{(1)}\in W_{1}$, $w_{(2)}\in W_{2}$ and $w_{(3)}\in W_{3}$
span
\[
(W_{1}\boxtimes_{P(z_{12})} W_{2})
\boxtimes_{P(z_{2})} W_{3},
\]
the diagram (\ref{hexagon1}) is commutative.

Now we consider the following diagrams:
\begin{equation}\label{hexagon2}
\begin{CD}
(W_{1}\boxtimes_{P(z_{12})} W_{2})
\boxtimes_{P(z_{2})} W_{3}&@>>>&(W_{1}\boxtimes W_{2})
\boxtimes W_{3}\\
@VVV&&@VVV\\
(W_{2}\boxtimes_{P(-z_{12})} W_{1})
\boxtimes_{P(z_{2})} W_{3}&@>>>&(W_{2}\boxtimes W_{1})
\boxtimes W_{3}
\end{CD}
\end{equation}
\begin{equation}\label{hexagon3}
\begin{CD}
(W_{2}\boxtimes_{P(-z_{12})} W_{1})
\boxtimes_{P(z_{2})} W_{3}&@>>>&(W_{2}\boxtimes W_{1})
\boxtimes W_{3}\\
@VVV&&@VVV\\
(W_{2}\boxtimes_{P(-z_{12})} W_{1})
\boxtimes_{P(z_{1})} W_{3}&@>>>&(W_{2}\boxtimes W_{1})
\boxtimes W_{3}
\end{CD}
\end{equation}
\begin{equation}\label{hexagon4}
\begin{CD}
(W_{2}\boxtimes_{P(-z_{12})} W_{1})
\boxtimes_{P(z_{1})} W_{3}&@>>>&(W_{2}\boxtimes W_{1})
\boxtimes W_{3}\\
@VVV&&@VVV\\
W_{2}\boxtimes_{P(z_{2})} (W_{1}
\boxtimes_{P(z_{1})} W_{3})&@>>>&W_{2}\boxtimes (W_{1}
\boxtimes W_{3})\\
\end{CD}
\end{equation}
\begin{equation}\label{hexagon5}
\begin{CD}
W_{2}\boxtimes_{P(z_{2})} (W_{1}
\boxtimes_{P(z_{1})} W_{3})&@>>>&W_{2}\boxtimes (W_{1}
\boxtimes W_{3})\\
@VVV&&@VVV\\
W_{2}\boxtimes_{P(z_{2})} (W_{3}
\boxtimes_{P(-z_{1})} W_{1})&@>>>&W_{2}\boxtimes (W_{3}
\boxtimes W_{1})
\end{CD}
\end{equation}
\begin{equation}\label{hexagon6}
\begin{CD}
(W_{1}\boxtimes_{P(z_{12})} W_{2})
\boxtimes_{P(z_{2})} W_{3}&@>>>&(W_{1}\boxtimes W_{2})
\boxtimes W_{3}\\
@VVV&&@VVV\\
W_{1}\boxtimes_{P(z_{1})} (W_{2}
\boxtimes_{P(z_{2})} W_{3})&@>>>&W_{1}\boxtimes (W_{1}
\boxtimes W_{3})
\end{CD}
\end{equation}
\begin{equation}\label{hexagon7}
\begin{CD}
W_{1}\boxtimes_{P(z_{1})} (W_{2}
\boxtimes_{P(z_{2})} W_{3})&@>>>&W_{1}\boxtimes (W_{1}
\boxtimes W_{3})\\
@VVV&&@VVV\\
(W_{2}\boxtimes_{P(z_{2})} W_{3})
\boxtimes_{P(-z_{1})} W_{1}&@>>>&(W_{2}\boxtimes W_{3})
\boxtimes W_{1}
\end{CD}
\end{equation}
\begin{equation}\label{hexagon8}
\begin{CD}
(W_{2}\boxtimes_{P(z_{2})} W_{3})
\boxtimes_{P(-z_{1})} W_{1}&@>>>&(W_{2}\boxtimes W_{3})
\boxtimes W_{1}\\
@VVV&&@VVV\\
W_{2}\boxtimes_{P(-z_{12})} (W_{3}
\boxtimes_{P(-z_{1})} W_{1})&@>>>&W_{2}\boxtimes (W_{3}
\boxtimes W_{1})
\end{CD}
\end{equation}
\begin{equation}\label{hexagon9}
\begin{CD}
W_{2}\boxtimes_{P(-z_{12})} (W_{3}
\boxtimes_{P(-z_{1})} W_{1})&@>>>&W_{2}\boxtimes (W_{3}
\boxtimes W_{1})\\
@VVV&&@VVV\\
W_{2}\boxtimes_{P(z_{2})} (W_{3}
\boxtimes_{P(-z_{1})} W_{1})&@>>>&W_{2}\boxtimes (W_{3}
\boxtimes W_{1})
\end{CD}
\end{equation}
The commutativity of the diagrams (\ref{hexagon2}), (\ref{hexagon5})
and (\ref{hexagon7}) follows {}from the definition of the commutativity
isomorphism for the braided tensor category structure and the
naturality of the parallel transport isomorphisms.  The commutativity
of (\ref{hexagon4}), (\ref{hexagon6}) and (\ref{hexagon8}) follows
{}from the definition of the associativity isomorphism for the braided
tensor product structure.  The commutativity of (\ref{hexagon3}) and
(\ref{hexagon9}) follows {}from the facts that compositions of parallel
transport isomorphisms are equal to the parallel transport
isomorphisms associated to the products of the paths and that parallel
transport isomorphisms associated to homotopically equivalent paths
are equal. The commutativity of the hexagon diagram involving 
$\mathcal{R}$ follows {}from (\ref{hexagon1})--(\ref{hexagon9}).

Finally, we prove the commutativity of the triangle diagram for the
unit isomorphisms. Let $z_{1}$ and $z_{2}$ be complex numbers such
that $|z_{1}|>|z_{2}|>|z_{1}-z_{2}|>0$ and let
$z_{12}=z_{1}-z_{2}$. Also let $\gamma$ be a path {}from $z_{2}$ to
$z_{1}$ in $\C$ with a cut along the nonnegative real line. We first
prove the commutativity of the following diagram:
\begin{equation}\label{unit1}
\begin{CD}
(W_{1}\boxtimes_{P(z_{12})} V)
\boxtimes_{P(z_{2})} W_{2}&
@>(\mathcal{A}_{P(z_{1}), P(z_{2})}^{P(z_{12}), P(z_{2})})^{-1}>>
&W_{1}\boxtimes_{P(z_{1})} (V
\boxtimes_{P(z_{2})} W_{2})\\
@Vr_{W_{1};z_{12}}\boxtimes_{P(z_{2})} 1_{W_{2};z_{2}}VV&
&@VV1_{W_{1}}\boxtimes_{P(z_{1})} l_{W_{2}}V\\
W_{1}\boxtimes_{P(z_{2})} W_{2}&@>>\mathcal{T}_{\gamma}>
&W_{1}\boxtimes_{P(z_{1})} W_{2}.
\end{CD}
\end{equation}

Let $w_{(1)}\in W_{1}$ and $w_{(2)}\in W_{2}$. Then we have 
\begin{eqnarray}\label{unit2}
\lefteqn{\overline{(1_{W_{1}}\boxtimes_{P(z_{1})} l_{W_{2};z_{2}})\circ
(\mathcal{A}_{P(z_{1}), P(z_{2})}^{P(z_{1}-z_{2}), P(z_{2})})^{-1}}
((w_{(1)}\boxtimes_{P(z_{12})}\mathbf{1})\boxtimes_{P(z_{2})}w_{(2)})}\nn
&&\quad\quad\quad=\overline{(1_{W_{1}}\boxtimes_{P(z_{1})} l_{W_{2}; z_{2}})}
(w_{(1)}\boxtimes_{P(z_{1})}(\mathbf{1}\boxtimes_{P(z_{2})}w_{(2)})\nn
&&\quad\quad\quad=w_{(1)}\boxtimes_{P(z_{1})}w_{(2)}.
\end{eqnarray}
But 
\begin{eqnarray}\label{unit3}
\lefteqn{\overline{r_{W_{1};z_{12}}\boxtimes_{P(z_{2})} 1_{W_{2}}}
((w_{(1)}\boxtimes_{P(z_{12})}\mathbf{1})\boxtimes_{P(z_{2})}w_{(2)})}\nn
&&=(e^{z_{12}L(-1)}w_{(1)})\boxtimes_{P(z_{2})}w_{(2)}.
\end{eqnarray}
Let $\Y=\Y_{\boxtimes_{P(z_{1})}, 0}$ be the intertwining 
operator of type ${W_{1}\boxtimes_{P(z_{1})}W_{2}\choose
W_{1}\;\;\;\;W_{2}}$ corresponding to the $P(z_{1})$-intertwining 
map $\boxtimes_{P(z_{1})}$. Then by the definition of 
the parallel transport isomorphism and the $L(-1)$-derivative
property for intertwining operators, we have
\begin{eqnarray}\label{unit4}
\overline{\mathcal{T}_{\gamma}}
((e^{z_{12}L(-1)}w_{(1)})\boxtimes_{P(z_{2})}w_{(2)})
&=&\Y(e^{z_{12}L(-1)}w_{(1)}, z_{2})w_{(3)}\nn
&=&\Y(w_{(1)}, z_{1})w_{(3)}\nn
&=&w_{(1)}\boxtimes_{P(z_{1})}w_{(2)}.
\end{eqnarray}
Since the elements 
$(w_{(1)}\boxtimes_{P(z_{12})}\mathbf{1})\boxtimes_{P(z_{2})}w_{(3)}$
for $w_{(1)}\in W_{1}$ and $w_{(2)}\in W_{2}$
span $(W_{1}\boxtimes_{P(z_{12})} V)\boxtimes_{P(z_{2})}W_{3}$,
(\ref{unit2})--(\ref{unit4}) give the commutativity 
of (\ref{unit1}). 

Let $\gamma_{1}$ be a path {}from $z_{1}$ to $1$ in $\C$ with a cut along 
the nonnegative real line.  Let $\gamma_{2}$ be the product of $\gamma$ 
and $\gamma_{1}$. In particular, $\gamma_{2}$ is a path 
{}from $z_{2}$ to $1$ in $\C$ with a cut along 
the nonnegative real line. Also let $\gamma_{12}$ be a path 
{}from $z_{12}=z_{1}-z_{2}$ to $1$ in $\C$ with a cut along 
the nonnegative real line. Then we have the following commutative 
diagrams:
\begin{equation}\label{unit5}
\begin{CD}
(W_{1}\boxtimes V)\boxtimes W_{2}&
@>\mathcal{A}^{-1}>>
&W_{1}\boxtimes (V\boxtimes W_{2})\\
@V\mathcal{T}_{\gamma_{2}}\circ (\mathcal{T}_{\gamma_{12}}
\boxtimes_{P(z_{2})}
1_{W_{2}})VV&
&@VV\mathcal{T}_{\gamma_{1}}\circ (1_{W_{1}}\boxtimes_{P(z_{1})}
\mathcal{T}_{\gamma_{2}})V\\
(W_{1}\boxtimes_{P(z_{12})} V)
\boxtimes_{P(z_{2})} W_{2}&@>>
(\mathcal{A}_{P(z_{1}), P(z_{2})}^{P(z_{12}), P(z_{2})})^{-1}>
&W_{1}\boxtimes_{P(z_{1})} (V
\boxtimes_{P(z_{2})} W_{2}).
\end{CD}
\end{equation}
\begin{equation}\label{unit6}
\begin{CD}
(W_{1}\boxtimes V)
\boxtimes W_{2}&
@>\mathcal{T}_{\gamma_{2}}^{-1}\circ (\mathcal{T}_{\gamma_{12}}^{-1}
\boxtimes 1_{W_{2}})>>
&(W_{1}\boxtimes_{P(z_{12})} V)
\boxtimes_{P(z_{2})} W_{2}\\
@Vr_{W_{1}}\boxtimes 1_{W_{2}}VV&
&@VVr_{W_{1}}\boxtimes_{P(z_{2})} 1_{W_{2}}V\\
W_{1}\boxtimes W_{2}&@>>\mathcal{T}_{\gamma_{2}}^{-1}>
&W_{1}\boxtimes_{P(z_{2})} W_{2}.
\end{CD}
\end{equation}
\begin{equation}\label{unit7}
\begin{CD}
W_{1}\boxtimes_{P(z_{1})} (W_{2}
\boxtimes_{P(z_{2})} W_{2})&
@>\mathcal{T}_{\gamma_{1}}\circ (1_{W_{1}}\boxtimes_{P(z_{1})}
\mathcal{T}_{\gamma_{2}})>>
&W_{1}\boxtimes (V
\boxtimes W_{2})\\
@V1_{W_{1}}\boxtimes_{P(z_{1})} l_{W_{2}}VV&
&@VV1_{W_{1}}\boxtimes l_{W_{2}}V\\
W_{1}\boxtimes_{P(z_{1})} W_{2}&@>>\mathcal{T}_{\gamma_{1}}>
&W_{1}\boxtimes W_{2}.
\end{CD}
\end{equation}
\begin{equation}\label{unit8}
\begin{CD}
W_{1}\boxtimes_{P(z_{2})} W_{2}&
@>\mathcal{T}_{\gamma}>>
&W_{1}\boxtimes_{P(z_{1})} W_{2}\\
@V\mathcal{T}_{\gamma_{2}}VV&
&@VV\mathcal{T}_{\gamma_{1}}V\\
W_{1}\boxtimes W_{2}&&=&
&W_{1}\boxtimes W_{2}.
\end{CD}
\end{equation}
The commutativity of (\ref{unit5}) follows {}from the definition 
of $\mathcal{A}$. The commutativity of (\ref{unit6}) and (\ref{unit7})
follows {}from the definition of the left and right unit isomorphisms 
and the parallel transport isomorphisms. The commutativity of (\ref{unit8})
follows {}from the fact that $\gamma_{2}$ is the product of 
$\gamma$ and $\gamma_{1}$. Combining (\ref{unit1}) and 
(\ref{unit5})--(\ref{unit8}), we obtain the commutativity of
the triangle diagram for the unit isomorphisms.

It is clear {}from the definition that $l_{V}=r_{V}$. 

Thus we have proved that the category $\mathcal{C}$ equipped with the
data given in Subsection 12.2 is a braided tensor category.  
\epfv

\newpage


\bigskip

\noindent {\small \sc Institut des Hautes \'{E}tudes Scientifiques, 
Le Bois-Marie, 35, Route De Chartres, F-91440 Bures-sur-Yvette, 
France}

\noindent {\it and}

\noindent {\small \sc Department of Mathematics, Rutgers University,
Piscataway, NJ 08854 (permanent address)}

\noindent {\em E-mail address}: yzhuang@math.rutgers.edu

\vspace{1em}

\noindent {\small \sc Department of Mathematics, Rutgers University,
Piscataway, NJ 08854}

\noindent {\em E-mail address}: lepowsky@math.rutgers.edu

\vspace{1em}

\noindent {\small \sc Department of Mathematics, Rutgers University,
Piscataway, NJ 08854}

\noindent {\em E-mail address}: linzhang@math.rutgers.edu

\end{document}